\documentclass[a4paper,11pt]{article}
\usepackage[a4paper, bottom=4cm, footskip=1.2cm]{geometry}
\usepackage[utf8]{inputenc}
\usepackage[T1]{fontenc} 
\usepackage[ngerman, english]{babel}
\usepackage{lmodern} 
\usepackage[autostyle]{csquotes}
\usepackage{enumitem}
\usepackage{float}
\usepackage[pdftex,pdfpagelabels,bookmarks,
hyperindex,hyperfigures,hidelinks]{hyperref}

\usepackage{amsthm,amsmath,amssymb,amsfonts} 
\usepackage{mathrsfs}
\usepackage{nicefrac} 
\usepackage{verbatim} 
\usepackage{tikz}
\usetikzlibrary{cd, babel}
\usepackage{extpfeil}

\usepackage{mathtools}
\usepackage{amssymb}
\usepackage{thmtools}
\usepackage{thm-restate}

\makeatletter
\newcommand{\oset}[3][0ex]{%
  \mathrel{\mathop{#3}\limits^{
    \vbox to#1{\kern-1\ex@
    \hbox{$\scriptstyle#2$}\vss}}}}
\makeatother

\usepackage{tabularx}
\usepackage{booktabs}
\usepackage[labelfont=bf,format=plain,justification=raggedright,singlelinecheck=false]{caption}

  \theoremstyle{definition}
  \newtheorem{definition}{Definition}[section]

  \newtheorem{example}[definition]{Example}
  
  \newtheorem*{remark}{Remark}
  
  \theoremstyle{plain}
  \newtheorem{theorem}[definition]{Theorem}
    \newtheorem{conjecture}[definition]{Conjecture}
  \newtheorem{prop}[definition]{Proposition}
  \newtheorem{lemma}[definition]{Lemma}
  \newtheorem{cor}[definition]{Corollary}
  \newtheorem*{theorem*}{Theorem}
  \newtheorem*{lemma*}{Lemma}

  \theoremstyle{remark}
  \newtheorem*{claim}{Claim}
  \newtheorem*{notation}{Notation}

\DeclareMathOperator{\id}{id}
  
  \DeclareMathOperator{\Ker}{Ker}

  \DeclareMathOperator{\image}{Im}

  \DeclareMathOperator{\colim}{colim}

        \DeclareMathOperator{\V}{\mathcal{V}}
                \DeclareMathOperator{\U}{\mathcal{U}}

\begin{document}

\begin{center}

\vspace*{2.5cm}

\huge

\textbf{Coarse Homotopy Theory and Shape Theory}

\vspace{2cm}

\normalsize

I n a u g u r a l d i s s e r t a t i o n 

\vspace{.5cm}

zur 

\vspace{.5cm}

Erlangung des akademischen Grades

\vspace{.5cm}

Dr. rer. nat. (Doctor rerum naturalium) 

\vspace{.5cm}

der

\vspace{.5cm}

Mathematisch-Naturwissenschaftlichen Fakultät

\vspace{.5cm}

der

\vspace{.5cm}

Universität Greifswald

\end{center}

\vspace{1.5cm}

\begin{flushright}

vorgelegt von \qquad \qquad \qquad \qquad \qquad \qquad 

\vspace{.4cm}

Felix Lange \, \,\,\qquad \qquad \qquad \qquad \qquad \qquad 

\end{flushright}

\vspace{1.5cm}

\begin{center}
Greifswald, 21.04.2026
\end{center}

\newpage

\thispagestyle{empty}

\null 

\vfill

The present document is a slightly modified version of the my dissertation, which was submitted to and accepted by the University of Greifswald, and will be referred to as \textit{thesis} throughout. \\

Universität Greifswald, Walther-Rathenau-Straße 47, 17489 Greifswald, Germany. I acknowledge financial support by the Deutsche Forschungsgemeinschaft DFG through Priority Programme SPP 2026 "Geometry at Infinity" (EN1163/5-1, project number 441426261, Macroscopic invariants of manifolds).\\

\vspace{1cm}

\newpage

\tableofcontents

\newpage

\begin{table}
\caption{Table of Notation}
\begin{tabularx}{\textwidth}{p{0.22\textwidth}X}
\toprule
  Parts I,II \\
 $cX$ & Euclidean cone of a compact metric space\\
 $Z_{h}$ & Maximal $\frac{1}{2^{h+1}}$-separated subset of $X$\\
 $\mathcal{U}_{h}$ & Open cover $\cup_{z \in Z_{h}}B(z,\frac{1}{2^h})$ of $X$\\
$|\mathcal{U}|$ & Nerve of the open cover $\mathcal{U}$\\
$\phi_{h}$ & Simplicial maps $\phi_{h}: |\mathcal{U}_{h}| \rightarrow |\mathcal{U}_{h-1}|$ induced by refinement of open covers\\
 $\phi_{l,s}$ &  $\phi_{s+1}\circ \dots \circ \phi_{l}: |\mathcal{U}_{l}| \rightarrow |\mathcal{U}_{s}|$ for $l>s$\\
 $\langle gst(u) \rangle $ &  Open geometric star of a vertex $u\in M_{X}$\\
$\mathscr{C}_{\Delta} $ &  Geometric core of the geometric simplex $\Delta$\\
$\mathcal{N}_{\Delta}$ &  Geometric neighbourhood of the geometric simplex $\Delta$\\
 $M_{X}$ & Inverse mapping telescope of $\{|\mathcal{U}_{h}|, \phi_{h+1}; \mathbb{N}_{0}\}$\\
$p$ & Projection to the height variable $p: cX\rightarrow [1,\infty)$, $p: M_{X}\rightarrow [1,\infty)$  \\
$p'$ & Projection to the $x$-variable $p': cX\rightarrow X$, $p': M_{X}\rightarrow \coprod_{h\in \mathbb{N}_{0}} |\mathcal{U}_{h}| $ \\ 
$b^H_{\omega_{0},\omega}$ & Change of base ray homomorphism $\pi_{n}^c(Z,\omega_{0})\rightarrow \pi_{n}^c(Z,\omega)$ via the homotopy $H: c[0,1]\rightarrow Z$  \\
  $i_{h}$ & $i_{h}: X \times (2^h,2^{h+1}] \rightarrow |\mathcal{U}_{h}| \times (2^h,2^{h+1}]$ induced by a partition of unity subordinate to $\mathcal{U}_{h}$\\
  $i$ & $i: cX\rightarrow M_{X}$ defined piecewise by $\{i_{h}\}_{h\in \mathbb{N}_{0}}$\\
 $R_{h}$ & Coarse map $R_{h}: |\mathcal{U}_{h}| \times (2^h,2^{h+1}] \rightarrow X \times (2^h,2^{h+1}]$ sending vertices to $Z_h$  \\
  $R$ & $R: M_{X}\rightarrow cX$ defined piecewise by $\{R_{h}\}_{h\in \mathbb{N}_{0}}$\\
$s_{r}$ & Shrinking map $s_{r}: cX\rightarrow cX$ dividing heights by $2^r$\\
$\mathfrak{q}_{a}$ & Deformation retract of $M_{X}$ to height $a$\\
$\mathbf{H}^r$ & Homotopy $I_{p}(cX)\rightarrow cX$ between $\id$ and $s_{r}$\\
$\mathbf{H}^{r,r'}$ & Homotopy $I_{p}(cX)\rightarrow cX$ between $s_{r}$ and $s_{r'}$\\
$D_{r}$ & Deformation retract $D_{r}: M_{X}\rightarrow M_{X}$ dividing heights by $2^r$\\
$\mathbf{D}^r$ & Homotopy $I_{p}(M_{X})\rightarrow M_{X}$ between $\id$ and $D_{r}$\\
$\mathbf{D}^{r,r'}$ & Homotopy $I_{p}(M_{X})\rightarrow M_{X}$ between $D_{r}$ and $D_{r'}$\\
$c_{a}X$ & Points in $cX$ with height $\geq a$\\
$c_{[a,b]}X$ & Points in $cX$ with height in $[a,b]$\\
$M_{h}$ & Points in $M_{X}$ with height $\geq 2^h+1$\\
$M_{[a,b]}$ & Points in $M_{X}$ with height in $[a,b]$\\
 $Y \oset{c}{-} L_i = L^c_{i}$ & CW-complement of $L_i$ in $Y$\\
  Part III \\
$X$ & Proper metric space\\
 $Z_{h}$ & Maximal $\frac{2^h}{2}$-separated subset of $X$\\
 $\mathcal{U}_{h}$ & Open cover $\cup_{z \in Z_{h}}B(z,2^h)$ of $X$\\
  $i_{h}$  &  $i_{h}: X\rightarrow |\mathcal{U}_{h}|$ induced by partition of unity subordinate to $\mathcal{U}_{h}$\\
 $\phi_{h,h+1}$ & Simplicial maps $\phi_{h,h+1}: |\mathcal{U}_{h}| \rightarrow |\mathcal{U}_{h+1}|$ induced by refinement of open covers\\
 $\phi_{s,l}$ &  $\phi_{l-1,l}\circ \dots \circ \phi_{s,s+1}: |\mathcal{U}_{s}| \rightarrow |\mathcal{U}_{l}|$ for $l>s$\\
 $i_{h}$  &  $i_{h}: X\rightarrow |\mathcal{U}_{h}|$ induced by partition of unity subordinate to $\mathcal{U}_{h}$\\
 $R_{h}$ & Coarse map $R_{h}: |\mathcal{U}_{h}|\rightarrow X$ sending vertices to $Z_h$\\ 
\bottomrule
\end{tabularx}
\end{table}

\newpage

\textit{Throughout this thesis, we write "$\cup_{U\in \mathcal{U}} U$ is an open cover of $X$" rather than "$\{U\in \mathcal{U}\}$ is an open cover of $X$". We also use the same symbol $f$ for a map $f: X\rightarrow Y$ and its induced map $f_{\ast}$ on homotopy groups, when the context is clear enough to do so.}

\newpage

\section{Introduction}


%


\subsection{Part I}

Euclidean cones are studied in coarse geometry because they provide a canonical way to turn compact metric spaces $X$ into coarse spaces $cX$ that are amenable to large-scale analysis. Coarse geometry studies spaces "from far away": the geometry of $cX$ at arbitrarily large scales reflects fine geometric and topological properties of $X$. Loosely speaking, the cone stretches the metric so that features of $X$ occuring at smaller and smaller scales appear further out along the cone, relating coarse invariants such as asymptotic dimension and coarse cohomology to classical invariants of $X$ \cite{engel2021coronasproperlycombablespaces}. Moreover, cone constructions have been adapted to produce important counterexamples in coarse geometry; for instance, Roe’s warped cone construction over a compact space with a group action yields spaces which cannot be coarsely embedded into a Hilbert space \cite{RoeWarped05}, and provide counterexamples to the coarse Baum–Connes conjecture \cite{kitsios2025coarsebaumconneswarpedcones}.\\

The Euclidean cone enjoys several nice structural features. The construction is natural with respect to continuous maps: maps between compact metric spaces induce coarse maps between their cones. Moreover, different reasonable choices of metrics on $X$ or $cX$ yield coarsely homotopy equivalent spaces. This robustness allows the cone to serve as a canonical large-scale model for compacta arising as boundaries or coronas, such as the visual boundaries of hyperbolic or CAT(0) spaces, where the cone of the boundary is coarsely homotopy equivalent to the original space and provides a concrete metric realisation of its geometry at infinity. \\

Let $X$ be a compact metric space which is isometrically embedded into a normed space $\iota: X\rightarrow L$ (eg. via the Kuratowski embedding). We define the cone $cX \subset L \times \mathbb{R}$ as follows:
\begin{align*}
cX:= \{(h\iota(x),h)\in L \times [1,\infty) \,|\, x\in X\}
\end{align*}
We equip $cX$ with the metric induced by the ambient norm on $L\times \mathbb{R}$ defined by $\|(h\iota(x),h)\| = \|h\iota(x)\|_{L}+ |h|$.  Up to coarse homotopy equivalence, $cX$ agrees with other standard constructions in the literature. \\

Let $x_0\in X$ and let $\omega:[1,\infty)\rightarrow cX, h\mapsto  (hx_0,h)$ be the standard parametrisation of the base ray $c\{x_0\}$. The $n$-th coarse homotopy group $\pi_{n}^{c}(cX,\omega)$ is defined as the coarse homotopy classes of coarse maps 
\begin{align*}
f: (c[0,1]^n,c\partial [0,1]^n) \rightarrow (cX,\omega)
\end{align*}
$\pi_{0}^{c}(cX,\omega)$ is defined as the pointed set of coarse homotopy classes of coarse base rays.\\

The motivation for studying coarse homotopy groups comes from the question: what is the coarse homotopy type of a cone of a compact metric space, and are there algebraic invariants to detect and classify this? The analogous question in usual homotopy theory has led historically to the development of a expansive theory with a rich tapestry of interesting results: cellular approximation, Whitehead theorem, and later on to the development of surgery theory.\\

The definition of coarse homotopy groups was introduced by Norouzizadeh in his PhD thesis \cite{norouzizadeh2010some}, in which he showed that $\pi_{n}^{c}(cS^n,c\{x_0\}) \cong \pi_{n}(S^n,x_0)$ for all $n\geq 0$. Mitchener, Norouzizadeh and Schick extended this result in \cite{mitchener2020coarse} to finite simplicial complexes $(K,x_0)$. They prove the isomorphism 
\begin{align*}
\pi_{n}^{c}(cK,c\{x_0\})\cong \pi_{n}(K,x_0)
\end{align*}
for all $n\geq 0$ by approximating a coarse map $f: (c[0,1]^n,c\partial [0,1]^n)\rightarrow (cK,c\{x_0\})$ by a simplicial map in the same coarse homotopy class, and showing that it is coarse-Lipschitz homotopic to the cone of a Lipschitz map $g: ([0,1]^n,\partial [0,1]^n)\rightarrow (K,x_0)$. \\

Restricting now to the coarse fundamental group, Weighhill proves in \cite{weighill2025liftingcoarsehomotopies} that $\pi^{c}_{1}(cM,c\{x_0\}) \cong \pi_{1}(M,x_0)$ for $M$ a compact manifold with non-positive sectional curvature. He develops a coarse version of covering space theory; using the fact that $ \pi_{1}(M,x_0)$ acts geometrically on the universal cover $\tilde{M}$, he obtains a group action of $\pi_{1}(M,x_0)$ on $c\tilde{M}$ with quotient $cM$. This induces a short exact sequence
\begin{align*}
0\rightarrow \pi^{c}_{1}(c\tilde{M}, c \{\tilde{x}_{0}\}) \rightarrow  \pi_{1}^{c}(cM,c\{x_0\}) \rightarrow  \pi_{1}(M,x_0) \rightarrow 0
\end{align*}
where $\tilde{x}_{0}$ is any lift of $x_0$. It is then proven that $\pi^{c}_{1}(c\tilde{M}, c \{\tilde{x}_{0}\}) = 0$. This is done by showing that every class is trivial by retracting along geodesics to the base ray. The assumption that $M$ has non-positive sectional curvature means that for $\gamma_{1},\gamma_{2}$ geodesic rays in $\tilde{M}$, $t\in [0,1]$ and $q,q'>0$, the inequality 
\begin{align*}
d(\gamma_{1}(tq),\gamma_{2}(tq'))\leq \max \{d(\gamma_{1}(0),\gamma_{2}(0)),d(\gamma_{1}(q),\gamma_{2}(q')) \}
\end{align*}
holds. This inequality is used multiple times to show that certain homotopies are coarse.  \\

The results in \cite{mitchener2020coarse} and \cite{weighill2025liftingcoarsehomotopies} both use the fine geometry of $cX$, and rely heavily on the fact that $X$ is locally well-behaved enough so that $cX$ is a path metric space. For an arbitrary compact metric space $X$ which is not a manifold or a simplicial complex, $cX$ may not be a path metric space. Therefore, current proofs available for computing $\pi_{n}^{c}(cX,c\{x_0\})$ cannot be generalisable to these cases. \\ 

In fact, after some preliminary investigation, I was sceptical that the coarse homotopy groups of $cX$ for an arbitrary compact metric space $X$ had anything to do with the specific choice of metric at all, because of the following: suppose we have a compact topological space $X$ which is metrisable. For any two choices of metrics $(X,d),(X,d')$ which induce the same topology, the identity maps in both directions are uniformly continuous homeomorphisms. By the proof of Lemma \ref{unicoarse}, there exist shrinking maps $S: c(X,d)\rightarrow c(X,d)$ and $S': c(X,d')\rightarrow c(X,d')$ such that $\id S$ and $\id S'$ are coarse. The compositions $(\id S)\circ(\id S')$ and $(\id S')\circ(\id S)$ are coarsely homotopic to the identity. By Lemma \ref{arrow1}, a coarse homotopy equivalence induces an isomorphism on coarse homotopy groups: $\pi_{n}^{c}(c(X,d),c\{x_0\}) \cong \pi_{n}^{c}(c(X,d'),c\{x_0\}) $ for all $n\geq 0$. Therefore the coarse homotopy groups should only depend on the \textit{topological} structure of $X$. \\

I was also unconvinced that the coarse homotopy groups of $cX$ should always coincide with the ordinary homotopy groups of $X$. The reasoning for this is that coarse maps $f: (c[0,1]^n,c\partial [0,1]^n)\rightarrow (cX,c\{x_0\})$ are allowed to be discontinuous, as long as the discontinuity is controlled by a uniform distance bound, and one can modify the image of $f$ by an arbitrary but uniformly bounded distance and remain in the same coarse homotopy class. Therefore for locally "bad" compact metric spaces $X$, like the Warsaw circle, I suspected that even after taking the cone, the coarse homotopy groups would not be able to detect the local topological structure but would rather measure the homotopy type of an approximation of $X$.\\

The idea of how to compute the coarse homotopy groups of the cone of a general compact metric space came from the following observation. Consider the Hawaiian earring $\mathbb{H}$ as an infinite wedge of shrinking circles (with centre $(\frac{1}{2^{j}},0)$, radius $\frac{1}{2^{j}}$) embedded in $\mathbb{R}^2$ with wedge point $x_0=(0,0)$ and the cone $c\mathbb{H}$ as a subset of $\mathbb{R}^3$:
\begin{align*}
\mathbb{H} = \bigvee_{j \in \mathbb{N}} S^1_{j} \quad \quad \quad c\mathbb{H} = \{(hx,h)\in \mathbb{R}^3\,|\, x\in \mathbb{H}, h\in [1,\infty)\}
\end{align*}  
I noticed that the set 
\begin{align*}
\bigcup_{j\in \mathbb{N}} S^1_{j} \times [1,2^j] \subset c\mathbb{H}
\end{align*}
is bounded distance away from the base ray $c\{x_0\}$. This means that for the purposes of coarse homotopy theory, the $j$-th circle is invisible until height $2^j$. Therefore, $c\mathbb{H}$ restricted to height $2^h$ looks like a wedge of $h$ circles. This led to the conjecture that $\pi_{1}^{c}(c\mathbb{H},c\{x_0\})$ should be isomorphic to $\varprojlim_{h} \pi_{1}(\vee_{j=1}^{h} S^1_{j}, x_0) = \varprojlim_{h} F_{h}$. \\

Following the strategy employed by Norouzizadeh in his thesis, I managed to prove the existence of a surjective map $\pi_{1}^{c}(c\mathbb{H},c\{x_0\})\rightarrow \varprojlim_{h} F_{h}$, but the kernel was difficult to compute. The problem is that working with arbitrary coarse maps is extremely cumbersome. Therefore I sought an approach that was less computational and generalisable to all compact metric spaces. \\

This is done by introducing the object $M_{X}$, which is the inverse mapping telescope of the nerves of a fixed, cofinal sequence of open covers $\{\mathcal{U}_{h}\}_{h\in \mathbb{N}}$ by balls of radius $\frac{1}{2^h}$ with simplicial refinement maps $\phi_{h+1,h}: |\mathcal{U}_{h+1}|\rightarrow |\mathcal{U}_{h}|$ as gluing maps. The idea for this was inspired by Wright, who used a similar construction to prove the coarse Baum-Connes conjecture for spaces with finite asymptotic dimension \cite{wright2005coarse}. Although his mapping telescope comes from an anti-\v{C}ech sequence, the underlying geometric idea is the same: instead of considering a coarse space and coarse maps, find a nicer space where computations can be done using continuous maps, and use these results to make statements about the original problem. The inverse mapping telescope $M_{X}$, although only a CW complex, is sufficiently well-behaved to define a version of (geometric) simplicial approximation, since all the gluing maps are simplicial and are separated by discrete intervals. It is also a path metric space, so the coarse structure is easy to work with. There are canonical maps $i: cX\rightarrow M_{X}$, defined via partitions of unity subordinate to $\mathcal{U}_{h}$, and $R: M_{X}\rightarrow cX$, which sends vertices to centres of balls in $\mathcal{U}_{h}$. After precomposing with a shrinking map $s: cX\rightarrow cX$, $i\circ s$ is a coarse homotopy equivalence with inverse $R$. \\

Using (geometric) simplicial approximation, I proved the isomorphism:
\begin{align*}
\pi_{n}^{c}(cX,c\{x_0\}) \cong \pi_{n}^{e}(M_{X},c\{x_0\})
\end{align*}
where $\pi_{n}^{e}(M_{X},c\{x_0\})$ denotes the end homotopy groups (also known as the strong or Steenrod homotopy groups) of $M_{X}$. These homotopy groups have been studied previously, and fit inside a $\varprojlim^1$ sequence. We therefore obtain the main theorems of Part I, Theorem \ref{maintheorem} and Theorem \ref{maintheorem2}: 

\begin{theorem*} Let $X$ be a compact metric space, $x_0\in X$, $\omega$ the standard parametrisation of $c\{x_0\}$, $n\geq 1$. There is a short exact sequence:
\begin{align*}
0\rightarrow {\varprojlim_h}^1  \pi_{n+1}(|\mathcal{U}_{h}|,x_0) \rightarrow  \pi_{n}^c(cX,\omega) \rightarrow   \varprojlim_{h} \pi_n (|\mathcal{U}_{h}|,x_0)  \rightarrow 0
\end{align*}
\end{theorem*}

\begin{theorem*}  Let $X$ be a compact metric space, $x_0\in X$, $\omega$ the standard parametrisation of $c\{x_0\}$. There is a short exact sequence of pointed sets:
\begin{align*}
 {\varprojlim_h}^1  \pi_{1}(|\mathcal{U}_{h}|,x_0) \hookrightarrow  \pi^{c}_{0}(cX,\omega) \twoheadrightarrow \varprojlim_{h} \pi_0 (|\mathcal{U}_{h}|,x_0) 
\end{align*}
\end{theorem*}

The term  $\varprojlim_{h} \pi_{n}(|\mathcal{U}_{h}|,x_0)$ is the $n$-th \v{C}ech homotopy group (set for $n=0$) $\check{\pi}_n(X,x_0)$, which is a shape invariant of $X$. The connection between coarse homotopy groups to something previously studied was a pleasant surprise. I was curious as to the exact relationship between the coarse homotopy type of $cX$ and the shape of $X$, with the goal of proving a coarse version of the Whitehead theorem for cones. This was a sensible idea since there are several versions known of the Whitehead theorem for pointed shape maps $F: S(X,x_0)\rightarrow S(Y,y_0)$ which induce an isomorphism pro-$\pi_{n}(F):$pro-$\pi_{n}(X,x_0)\rightarrow$ pro-$\pi_{n}(Y,y_0)$ in pro-Gp. This line of inquiry led to the results of Part II:

\subsection{Part II}

Let $(X,x_0), (Y,y_0)$ be compact metric spaces with base points. Let $\omega,\tau$ be the standard parametrisations of the base rays $c\{x_0\},c\{y_0\}$ respectively. We have the following series of implications; blue arrows are proven in this thesis, red arrows are known theorems from the literature, and black arrows come from definitions.

\begin{enumerate}
\item A coarse homotopy equivalence $\mathfrak{a}:(cX,\omega)\rightarrow (cY,\mathfrak{a}\omega)$ induces an isomorphism in coarse homotopy groups. This is a corollary of Lemma \ref{arrow1}. 
\item A coarse map $\mathfrak{a}:(cX,\omega)\rightarrow (cY,\mathfrak{a}\omega)$, defines a proper, cellular, locally Lipschitz map $\psi:=\widetilde{is_{l}\mathfrak{a}R}: (M_{X},\omega)\rightarrow (M_{Y}, \widetilde{is_{l}\mathfrak{a}\omega})$ which is the (geometric) cellular approximation of $i_{Y}s_{l}\mathfrak{a}R_{X}$. This is shown in Proposition \ref{psiexistence}. If $\mathfrak{a}$ induces an isomorphism on coarse homotopy groups, $\psi$ induces an isomorphism on end homotopy groups.  
\item A proper cellular map $\psi: (M_{X}, \omega)\rightarrow (M_{Y},\psi\omega)$ which induces an isomorphism on end homotopy groups has a proper homotopy inverse, assuming that $X,Y$ are connected and have finite shape dimension. The statement is found in Lemma \ref{weakimpliesphe} and is a result of Corollary \ref{inverse2}.
\item If a coarse map $\mathfrak{a}: cX \rightarrow cY$ has a coarse homotopy inverse $\mathfrak{b}: cY\rightarrow cX$, then $\psi:=\widetilde{i_{Y}s_{l}\mathfrak{a}R_{X}}$ is a proper homotopy equivalence with inverse $\xi:= \widetilde{i_{X}s_{r}\mathfrak{b}R_{Y}}$. This is proven in Theorem \ref{SSBornCoarse2}. 
\item A proper homotopy equivalence $\psi: M_{X}\rightarrow M_{Y}$ defines a coarse homotopy equivalence $R_{Y}(\psi' S_{\psi'})i_{X}s_{X}: cX\rightarrow cY$, where $\psi'$ is a locally Lipschitz map in the proper homotopy class of $\psi$ and $S_{\psi'}: M_{X}\rightarrow M_{X}$ is a shrinking map. This is proven in Theorem \ref{SSBornCoarse}. 
\item A proper homotopy equivalence $\psi: (M_{X},\omega)\rightarrow (M_{Y},\psi \omega)$ induces a shape equivalence. This is Lemma \ref{phemeansshape}. 
\item A shape morphism $F: \mathcal{X}\rightarrow \mathcal{Y}$ induces a continuous, proper, cellular map $\psi:M_{\mathcal{X}}\rightarrow M_{\mathcal{Y}}$, where $M_{\mathcal{X}}, M_{\mathcal{Y}}$ are inverse mapping telescopes associated with the shape expansions $\mathcal{X},\mathcal{Y}$ respectively. Pre- and postcomposing with the proper homotopy equivalences $M_{X}\simeq M_{\mathcal{X}}, M_{\mathcal{Y}}\simeq M_{Y}$ respectively, we obtain a proper, cellular map $\psi': M_{X}\rightarrow M_{Y}$. If $F$ is a pointed shape morphism, then $\psi': (M_{X},\omega)\rightarrow (M_{Y},\tau)$ is a map of pairs. In this case, if $F$ induces an isomorphism on pro-$\pi_{n}$ for all $n\geq 0$, then the map $\psi'$ induces an isomorphism on end homotopy groups. This is Theorem \ref{wseimplieswphe1} and Corollary \ref{wseimplieswphe2}. 
\item The proper homotopy inverse from point $3.$ can be chosen to be a map of pairs. 
\end{enumerate}

\begin{figure}[H]
\centering
   \includegraphics[width=0.7 \linewidth]{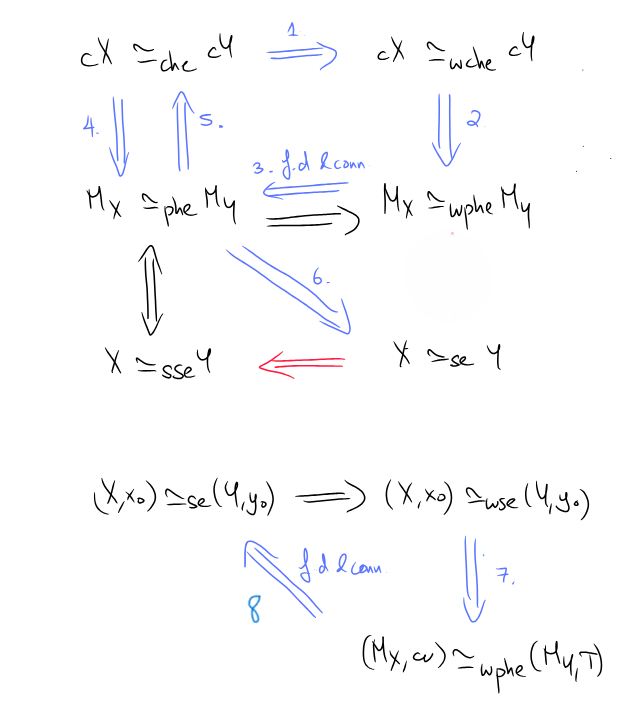}
  \label{fig:arrows}
\end{figure}

Following the top square clockwise gives us the coarse Whitehead theorem (Theorem \ref{unbasedwhitehead}):

\begin{theorem*} Let $(X,x_0)$ and $(Y,y_0)$ be connected, compact metric spaces with finite shape dimension. Let $\omega$ be the standard parametrisation of the base ray $c\{x_0\}$. Let $\mathfrak{a}:(cX,\omega)\rightarrow (cY,\mathfrak{a}\omega)$ be a coarse map, which induces an isomorphism on coarse homotopy groups $\mathfrak{a}: \pi_{n}^{c}(cX,\omega)\rightarrow \pi_{n}^{c}(cY,\mathfrak{a}\omega)$ for all $n\geq 0$. Then $\mathfrak{a}$ is a coarse homotopy equivalence, ie. there exists $\mathfrak{b}: cY\rightarrow cX$ such that $\mathfrak{b}\circ \mathfrak{a} \simeq id_{cX}$ and $\mathfrak{a}\circ \mathfrak{b} \simeq id_{cY}$.
\end{theorem*} 

Following the bottom triangle clockwise recovers a classical theorem in shape theory: the shape Whitehead theorem for compact, connected metric spaces (Corollary \ref{wseimplieswphe2}). The relationship between the two diagrams is known to be a difficult problem. For example, pointed shape map $F:S(X,x_0)\rightarrow S(Y,y_0)$ which has an unpointed shape inverse does not have to have a pointed shape inverse \cite{dydak1980pointed}. \\

Assembling arrows $4$ and $5$ into a concise statement, we obtain Theorems \ref{SSBornCoarse} and \ref{SSBornCoarse2}. 

\begin{theorem*} There is a full functor 
\begin{align*}
\mathcal{F}: \mathbf{StrongShape} \rightarrow \mathbf{hCoarse}
\end{align*}
where $\mathbf{StrongShape}$ has $(W,X)$ as objects ($W$ is a compact $AR$ and $X$ is a (closed) $Z$-set in $W$) and morphisms proper homotopy classes of continuous, proper maps $\psi: W_X \setminus X\rightarrow W_{Y}\setminus Y$. $\mathbf{hCoarse}$ denotes the homotopy category of coarse spaces with bornology induced from the coarse structure. Morphisms in $\mathbf{hCoarse}$ are coarse homotopy classes of coarse bornological maps.  
\end{theorem*}

\begin{theorem*} Two objects $(W_{X},X),(W_{Y},Y)$ are equivalent in $\mathbf{StrongShape}$ if and only if their images $\mathcal{F}(W_{X},X),\mathcal{F}(W_{Y},Y)$ are equivalent in $\mathbf{hCoarse}$. 
\end{theorem*}

The functor $\mathcal{F}$ assigns $(W_{X},X)$ to the Euclidean cone $cX$ of the compact metric space $X$, and morphisms $[\psi]: M_{X}\rightarrow M_{Y}$ (after composing with equivalences $M_{X}\simeq W_{X}$ and $M_{Y}\simeq W_{Y}$) to the coarse homotopy class of $[R_{Y}\psi' S_{\psi'} i_{X} s_{X}]$. For details, see Definition \ref{DefF}.\\

To summarise, we establish a correspondence between coarse homotopy equivalences of Euclidean cones of compact metric spaces and proper homotopy equivalences of their mapping telescopes. This provides a concrete realisation of the principle that the coarse homotopy type of Euclidean cones encodes exactly the strong shape type of the underlying compact metric spaces, yielding explicit maps and homotopies realising the equivalence. To my knowledge, this has not appeared previously in the literature.\\

\subsection{Part III}

The contents of Part III were created after a literature search to see if the methods developed in Parts I and II could be useful in solving other problems in current research. \\

The following theorems are proven in a recent paper \cite{ashley2025interactions}:

\begin{theorem*} (Theorem $7.4$ in \cite{ashley2025interactions}) Let $X$ be a proper geodesic metric space. There exists a natural surjection $\pi^c_{0}(X)\rightarrow \mathcal{E}nds(X)$. 
\end{theorem*}

\begin{theorem*} (Theorem $9.9$ in \cite{ashley2025interactions}) Let $X$ be a locally finite geometric tree. There exists a natural isomorphism $\pi^c_{0}(X)\cong \mathcal{E}nds(X)$.
\end{theorem*}

I was able to recover these results, as well as generalise them to statements about $\pi_{n}^{c}(X,\omega)$ ($\omega$ is a coarse base ray) for all $n\geq 0$ and $X$ an arbitrary connected, proper metric space, using a similar proof method to the $\varprojlim^1$ sequence for cones. Rather than a \v{C}ech sequence, we take an anti-\v{C}ech sequence $\{\mathcal{U}_{h}\}_{h\in \mathbb{N}_{0}}$ of open covers of radius $2^h$ centred at points $z \in Z_{h}$ in a maximal $\frac{2^h}{2}$-separated subset of $X$. There are canonical maps $i_{h}: X\rightarrow |\mathcal{U}_{h}|$ and $R_{h}: |\mathcal{U}_{h}|\rightarrow X$ for each $h\in \mathbb{N}_{0}$, defined analogoulsy to $i,R$ from Part I.  From the simplicial bonding maps $\phi_{h,h+1}: |\mathcal{U}_{h}|\rightarrow |\mathcal{U}_{h+1}|$ we obtain a directed system of end homotopy groups $\{\pi_{n}^{e}(|\mathcal{U}_{h}|,\phi_{h_0,h} \widetilde{i_{h_0}\omega}), (\phi_{h,h+1})_{\ast}; \mathbb{N}_{\geq h_0}\}$.

\begin{theorem*} (Theorem \ref{main1}) There is a natural isomorphism 
\begin{align*}
\Theta: \pi_{n}^{c}(X,\omega) &\rightarrow \varinjlim_{h\in \mathbb{N}_{\geq h_0}}\pi^{e}_{n}(|\mathcal{U}_{h}|, \phi_{h_0,h} \widetilde{i_{h_0}\omega})\\
[f] &\mapsto b_{\widetilde{i_{h}\omega},\phi_{h_0,h} \widetilde{i_{h_0}\omega}}[\widetilde{i_{h} f}]
\end{align*}
\end{theorem*}

\begin{theorem*} (Theorem \ref{main2}) There is a natural isomorphism 
\begin{align*}
\theta: \pi_{0}^{c}(X,\omega) &\rightarrow \varinjlim_{h\in \mathbb{N}_{\geq h_0}}\pi^{e}_{0}(|\mathcal{U}_{h}|, \widetilde{i_{h}\omega})\\
[\tau] &\mapsto [\widetilde{i_{h} \tau}]
\end{align*}
\end{theorem*}

The main technical tool is once again simplicial approximation of a homotopy class $[f]$ after mapping it into an approximating space $|\mathcal{U}_{h}|$ of $X$. There are significantly fewer technical details to these proofs compared to those in Parts I and II, since $|\mathcal{U}_{h}|$, unlike the mapping telescope, is an actual simplicial complex; the reader is encouraged to browse through Part III first if they wish to become acquainted with the ideas of this thesis without too much frustration. \\

The $\varprojlim^1$-sequence of end homotopy groups yields us Theorem \ref{mainthm}:

\begin{theorem*} (Theorem \ref{mainthm}) Let $X$ be a connected proper metric space and $[\omega]\in \pi_{0}^{c}(X)$ a fixed base ray such that $\widetilde{i_{h}\omega}$ exists for all $h\geq h_0$.  Let $\mathcal{L}_{h}:=\{L_{j,h}\}_{j\in \mathbb{N}_{0}}$ be a finite filtration of $|\mathcal{U}_{h}|$ such that $\widetilde{i_{h}\omega}$ is well-parametrised with respect to $\mathcal{L}_{h}$. Let $n\geq 1$. There is a short exact sequence of groups 
 \begin{align*}
0\rightarrow  \varinjlim_{h} {\varprojlim_{j}}^{1} \pi_{n+1}(L_{j,h}^c,\widetilde{i_{h}\omega}(j))  \rightarrow  \pi_{n}^c(X,\omega) \rightarrow  \varinjlim_{h} \varprojlim_{j} \pi_{n}(L_{j,h}^c,\widetilde{i_{h}\omega}(j)) \rightarrow 0
\end{align*}
For $n=0$, there is a short exact sequence of pointed sets
\begin{align*}
\varinjlim_{h} {\varprojlim_{j}}^{1} \pi_{1}(L_{j,h}^c,\widetilde{i_{h}\omega}(j))  \hookrightarrow \pi^{e}_{0}(Y,\omega) \twoheadrightarrow  \varinjlim_{h} \varprojlim_{j} \pi_{0}(L_{j,h}^c,\widetilde{i_{h}\omega}(j))
\end{align*}
\end{theorem*}

The results of \cite{ashley2025interactions} can be obtained as special cases of Theorem \ref{mainthm}. If $X$ is geodesic, then the bonding maps $\phi_{h,h+1}$ are quasi-isometries. Since the ends of a space is a quasi-isometry invariant, the direct limit on the right collapses to $h=h_0$. $i_{h_0}$ is also a quasi-isometry, so 
\begin{align*}
\varinjlim_{h} \varprojlim_{j} \pi_{0}(L_{j,h}^c,\widetilde{i_{h}\omega}(j))\cong \mathcal{E}nds(X,\omega)
\end{align*} 
which gives us the surjection in \cite{ashley2025interactions}. If $X$ is a locally finite geometric tree, then for all $j\in \mathbb{N}$, $L_{j,h}^c$ is homotopy equivalent to a set of points which in the limit corresponds to the ends of $|\mathcal{U}_{h}|$. This boils down to the fact that metric balls in a tree are convex, so the intersection $\cap_{i=0}^{m} B(z_i,2^h)$ is a non-empty subtree of $X$. We obtain 

\begin{lemma*} (Lemma \ref{geomtree}) Let $X$ be a locally finite geometric tree, $\omega: [1,\infty)\rightarrow X$ a coarse base ray. The coarse homotopy groups of $X$ are

\[ \pi_{n}^{c}(X,\omega) \cong \begin{cases} 
        0 & n\geq 1 \\ 
        \mathcal{E}nds(X,\omega) & n=0\\
       \end{cases}
    \]
\end{lemma*}

There is also an example constructed in \cite{ashley2025interactions} where  $\pi^c_{0}(X)\rightarrow \mathcal{E}nds(X)$ is not injective, a $1$-dimensional ladder $X$ with rungs of length $n^2$ at heights $n\in \mathbb{N}$. Using the $\varprojlim^1$-sequence, I found an algebraic description of the kernel and computed the coarse homotopy groups for $X$. I  generalised this to all ladders with rung lengths which grow unboundedly; $\pi_{0}^{c}$ is the $\varprojlim^1$ of an inverse system of free groups on countably many generators, and all higher homotopy groups vanish:

\begin{lemma*} (Lemma \ref{laddersize}) Let $X$ be a ladder: a $1$-dimensional simplicial complex defined by two rays $\alpha,\alpha'$ joining at a distinguished vertex at height $0$, and rungs of length $l_{m}\in \mathbb{N}$  joining the rays at heights $h_{m} \in \mathbb{N}$. Let $X$ be given the path metric. Assume that $\{h_{m}\}_{m\in \mathbb{N}}$ is a strictly increasing sequence and that $\lim_{m\rightarrow \infty} l_{m}=\infty$. Then the coarse homotopy groups of $X$ are:

\[ \pi_{n}^{c}(X,\alpha) \cong \begin{cases} 
       0 &n\geq 1  \\
        {\varprojlim}^1_{m}  F_{(\mathbb{N}_{0})_{\geq m}} & n=0 \\ 
       \end{cases}
    \]
where the bonding maps $\psi_{m+1,m}: F_{(\mathbb{N}_{0})_{\geq m+1}} \rightarrow F_{(\mathbb{N}_{0})_{\geq m}}$ are the canonical inclusion. 
\end{lemma*}

These ladders also provide us with an important counterexample. Let $Y$ be a ladder with rungs of length $2^m$ at heights $2^m$ for $m\geq 10$. The Lipschitz map $Y\rightarrow X$ which sends the rung at height  $2^m$ to the rung at height $m$ does not have a coarse homotopy inverse.

\begin{theorem*} (Theorem \ref{fail}) There exist $1$-dimensional simplicial complexes $Y,X$ and a proper Lipschitz map $\mathfrak{a}: Y\rightarrow X$ which induces an isomorphism on all coarse homotopy groups, but is not a coarse homotopy equivalence.
\end{theorem*}

This counterexample shows that coarse homotopy groups are distinguishing invariants only for cone-like spaces: the proof of the coarse Whitehead theorem for cones relies heavily on the fact that one can "shrink" arbitrarily far down the cone to make any continuous map coarse. Obviously this argument fails if there is no way to shrink without destroying the topology. The goal to discover better invariants for the coarse homotopy type of a space is an interesting subject for possible future research.

\subsection{Applications to geometric group theory}

The non-believer may wonder why we even bother generalising to $cX$, for $X$ not a CW complex.  Do these pathological spaces ever appear in "nature"? It turns out that they do appear, with overwhelming probability (in a precise sense), as boundaries of groups.\\

Consider a finitely generated group $G$ equipped with the coarse structure induced from the word metric. If $G$ is hyperbolic, CAT$(0)$, systolic, or hierarchically hyperbolic (and additionally proper and finite-dimensional) (Examples 3.27 in \cite{engel2021coronasproperlycombablespaces}) it is shown by Fukaya-Oguni \cite{Fukaya_2018} that $G$ is coarsely homotopy equivalent to the cone of its Gromov boundary $\partial G$ (or a suitable replacement). More precisely, 

\begin{theorem} (Theorem $1.1$ in \cite{Fukaya_2018}) Let $X$ be a proper coarsely convex space. Then $X$ is coarsely homotopy equivalent to $O(\partial X)$, the open cone over the ideal boundary of $X$. 
\end{theorem}

If $G$ acts geometrically on $X$, then $G$ is quasi-isometric to $X$. The open cone $O\partial X$ in Fukaya's construction is coarsely equivalent to $c(\partial X)$. The results of this thesis imply that two such groups $G_{1},G_{2}$ with boundaries $\partial G_{1}, \partial G_{2}$ (which are always compact and metrisable) are coarsely homotopy equivalent if and only if their boundaries $\partial G_1, \partial G_2$ are shape equivalent. This quite a weak condition for the following reasons:

\begin{enumerate}
\item A shape equivalence $S(\partial G_{1})\simeq S(\partial G_{2})$ does not have to be induced from a continuous map of spaces $\partial G_{1}\rightarrow \partial G_{2}$. \cite{zdravkovska1981example}
\item The shape classification of compact metric spaces is strictly weaker than the homotopy classification, coinciding on the subcategory of spaces homotopy equivalent to finite simplicial complexes. 
\end{enumerate}

It was already clear before the writing of this thesis that a coarse homotopy equivalence is much weaker than a quasi-isometry. For example, from the work of Bunke and Engel \cite{bunke2020homotopy}, a homotopy equivalence between compact metric spaces induces a coarse homotopy equivalence on their cones. However, we know that if $G_1,G_2$ are hyberbolic, a quasi-isometry $\mathfrak{a}:G_1\rightarrow G_2$ extends to a \textit{homeomorphism} on the Gromov boundaries $\partial \mathfrak{a}: \partial G_1\rightarrow \partial G_2$. The results of this thesis imply that the coarse homotopy classification of hyperbolic groups is even weaker than what previous research suggested: a shape equivalence $\partial G_{1} \simeq_{se} \partial G_{2}$ is sufficient to induce a coarse homotopy equivalence $G_1 \simeq_{che} G_{2}$. Whereas for quasi-isometries, the existence of a homeomorphism $\partial G_{1} \rightarrow \partial G_{2}$ does not imply that groups $G_1,G_2$ are quasi-isometric (the homoeomorphism must additionally respect the quasi-conformal structures on the boundaries \cite{paulin1996groupe}). Interestingly, in the CAT$(0)$ case, a quasi-isometry between groups does not even necessarily induce a homeomorphism on their boundaries, only a shape equivalence (if the groups are torsion-free \cite{bestvina1996local}). The precise relationship between these concepts is complicated and the subject of active research. 

\begin{prop} (\cite{alonso1991notes}) Let $G$ be a group generated by a finite set $S$ and suppose that the Cayley graph $\Gamma(G,S),d_{S}$ is $\delta$-hyperbolic. Let $d\geq 4\delta +2$ be an integer. Then $P_{d}(G,S)$ is contractible.
\end{prop}

$P_{d}(G,S)$ is a finite-dimensional simplicial complex with a natural left action of $G$. If $G$ is torsion-free, the action is free and transitive, with compact quotient $BG = P_{d}(G,S)/G$. It is therefore a model for the classifying space $EG$. Since $P_{d}(G,S)$ is coarsely equivalent (even quasi-isometric) to $G$ and $\partial G$ is a $Z$-set in $P_{d}(G)\cup \partial G$ \cite{bestvina1991boundary}, our results show that the coarse homotopy type of such classifying spaces is uniquely determined by its proper homotopy type, which is determined completely by the shape of $\partial G$.
\begin{align*}
G\simeq_{qi} P_{d}(G)\simeq_{phe} M_{\partial G} \simeq_{che} c(\partial G) \simeq_{che} G
\end{align*} 
where $qi$ stands for quasi-isometry, $phe$ stands for proper homotopy equivalence, and $che$ for coarse homotopy equivalence. \\

We now discuss some technicalities present in the coarse Whitehead theorem (Theorem \ref{unbasedwhitehead}) with an outlook of their impact on applicability of results. Let $\partial G_1,\partial G_2$ be the boundaries of suitable (in the sense of \cite{Fukaya_2018}) groups $G_1,G_2$.  Recall that Theorem \ref{unbasedwhitehead} allows us to detect when a weak coarse homotopy equivalence $\mathfrak{a}:G_1\rightarrow G_2$ is a coarse homotopy equivalence. It requires the assumptions of finite (shape) dimensionality and connectedness of the boundaries.\\

\textbf{Finite shape dimension}: For a compact metric space $\partial G$, the shape dimension is equivalent to the topological dimension. Corollary $8.10$ of \cite{Fukaya_2018} states that if $G$ admits a finite model for its classifying space $BG$, then 
\begin{align*}
cd(G) = sd (\partial G)+1
\end{align*}
where $cd(G)$ denotes the cohomological dimension of $G$. Even if $cd(G) = \infty$, there are still some things known about the dimension of $\partial G$, which we summarise in the table below: the boundary type is the suitable definition so that $G\simeq_{che} c(\partial G)$. A metric space $(X,d)$ has asymptotic dimension at most $n$ ($asdim(X)\leq n$) if for every $R>0$ there exists a uniformly bounded cover $\mathcal U$ of $X$ such that every ball $B(x,R)$ intersects at most $n+1$ elements of $\mathcal U$.

\begin{center}
\begin{tabular}{|p{3.5cm}|p{3.5cm}|p{6.5cm}|}

\hline

\textbf{Group class} 
& \textbf{Boundary type} 
& \textbf{$asdim(G)$ and $sd(\partial G)$} \\
\hline

Hyperbolic groups 
& Gromov boundary  
& $asdim(G) = sd(\partial G) + 1<\infty$ \cite{BuyaloLebedeva2007} \\
\hline

CAT(0) groups 
& Visual boundary  
& Relation between $asdim(G)$ and $sd(\partial G)$ is open \\
\hline

Systolic groups 
& Systolic boundary 
& $sd(\partial G)< \infty$ \cite{Osajda2009} \\
\hline

Relatively hyperbolic groups (+ conditions)
& Bowditch boundary 
& Open problem; all known examples $sd(\partial G)< \infty$ \\
\hline
\end{tabular}
\end{center}

\textbf{Connectedness} The proof for coarse Whitehead theorem only works for connected spaces. This is equivalent to the group $G$ having one end. By Stalling's theorem, a group either has one end, two ends, or infinitely many ends. If it has two ends, it is virtually $\mathbb{Z}$. It has infinitely many ends if it can be written as an amalgamated product of free groups or as a HNN extension. Otherwise it has one end. These theorems can be found in the fun book "Trees" by Serre \cite{serre2002trees}. Therefore the restriction on $\partial G$ being connected is not too severe.\\

We now answer the potential question of the non-believer. A random group $G_{ran}$ (for a definition, see Example 8 in Subsection \ref{goodexamples}) is infinite, hyperbolic, $1$-ended, and has boundary homeomorphic to the Menger sponge $\mathbb{M}$, a topologically $1$-dimensional fractal embedded in $\mathbb{R}^3$ with infinitely many holes, zero volume and infinite surface area. We can compute the coarse homotopy groups of $G_{ran}$ quite easily. $\mathbb{M}$ is the inverse limit of a sequence of $4$-valent graphs $\{X_h\}_{h\in \mathbb{N}}$, where each $X_h$ has $2\cdot 4^h$ vertices. We have that $\pi_1(X_{h+1}) = \pi_{1}(X_h) \ast F_{3}\ast \dots \ast F_{3}$  with one factor of $F_{3}$ for each vertex of $X_{h}$. The bonding morphisms $\pi_{1}(X_{h+1})\rightarrow \pi_{1}(X_h)$ trivialise these additional factors. Therefore,  

\[ \pi_{n}^{c}(G_{ran},\omega) \cong \begin{cases} 
       0 &n\geq 2  \\
        {\varprojlim}_{h}  F_{c_{h}} & n=1 \\ 
        \{\ast\} & n=0\\
       \end{cases}
    \]

where $c_{h} =2\cdot 4^{h}+1$, and the bonding maps $\psi_{h+1,h}: F_{c_{h+1}} \rightarrow F_{c_{h}}$ are the canonical projection.

\newpage

\part{Coarse Homotopy Groups}

\section{Definitions and setup} We begin with some preliminaries about coarse structures and develop the definition of $\mathbf{BornCoarse}$ with metric spaces as the motivating example. We introduce a canonical coarse structure on the Euclidean cone of a compact metric space, and a simplicial structure on the cone of $[0,1]^n$. This leads to the geometric definition of coarse homotopy groups. We cite some elementary properties, such as functoriality, the existence of a long exact sequence, and that a coarse homotopy equivalence induces an isomorphism on coarse homotopy groups. We then define the change of base ray homomorphisms necessary for the rest of the thesis. 

\subsection{Coarse structures}

We begin with a sequence of definitions collected from \cite{bunke2020homotopy} and   \cite{mitchener2020coarse}.\\

A \textit{coarse structure} on a set $X$ is a distinguished collection, $\mathcal{C}$, of subsets of the product $X\times X$  (called \textit{entourages}) such that:
\begin{itemize}
\item Any finite union of entourages is an entourage. Any subset of an entourage is an entourage.
\item The union of all entourages is the entire space $X \times X$.
\item The inverse of an entourage $V$:
\begin{align*}
V^{-1} = \{(y,x)\in X\times X \,|\, (x,y)\in V\}
\end{align*}
is an entourage. 
\item The composition of entourages $V_1$ and $V_2$:
\begin{align*}
V_1\circ V_2 = \{(x,z)\in X\times X\,|\,(x,y)\in V_1, (y,z)\in V_{2} \text{ for some } y\in X\}
\end{align*}
is an entourage. 
\item The diagonal $\Delta = \{(x,x) \,|\, x\in X\}$ is an entourage. 
\end{itemize}

A space equipped with a coarse structure is called a \textit{coarse space}. \\

Let $X$ be a set. A \textit{bornology} on $X$ is a subset $\mathcal{B}$ of the power set $\mathcal{P}(X)$ which is closed under forming finite unions, taking subsets, and which contains all one-point sets. The elements of $\mathcal{B}$ are called the \textit{bounded sets} of $X$. A map between sets equipped with bornologies is said to be \textit{proper} if pre-images of bounded sets are bounded.\\

For a coarse entourage $V \in \mathcal{C}$ and bounded set $B\in \mathcal{B}$, we define the $V$\textit{-thickening of} $B$ as
\begin{align*}
V[B] := \{ x\in X\,|\, \exists b\in B, (x,b)\in V\}
\end{align*}

 A bornology $\mathcal{B}$ and a coarse structure $\mathcal{C}$ on the same set are said to be compatible if $\mathcal{B}$ is stable under forming thickenings. A \textit{bornological coarse space} is a triple $(X,\mathcal{C},\mathcal{B})$ of a set with a coarse structure $\mathcal{C}$ and bornology $\mathcal{B}$ such that the coarse structure and the bornology are compatible.\\

Let $(X,d)$ be a proper metric space (ie. the closures of sets of finite diameter are compact). This will be the main example of a bornological coarse space considered in this thesis. The \textit{bounded coarse structure} $\mathcal{C}_{d}$ on $X$ is the coarse structure formed by defining the entourages to be subsets of $R$-neighbourhoods of the diagonal:
\begin{align*}
D_{R}= \{(x,y)\in X\times X \,|\, d(x,y)<R\}
\end{align*}
for $R\in \mathbb{R}_{> 0}$. The \textit{metric bornology} $\mathcal{B}_{d}$ on $X$ consists of sets which are bounded with respect to the metric. This bornology is compatible with the coarse structure since for a bounded set $B$, 
\begin{align*}
D_{R}[B] = \{x\in X\,|\, \exists y\in B, d(x,y) <R\}
\end{align*}
is a bounded set of diameter $< R+diam(B)$. We say that $\mathcal{C}_{d},\mathcal{B}_{d}$ are induced by the metric $d$ and call $(X,\mathcal{C}_{d}, \mathcal{B}_{d})$ the bornological coarse space associated to $(X,d)$. \\

We now consider a map $f:(X, \mathcal{C},\mathcal{B})\rightarrow (Y, \mathcal{C}',\mathcal{B}')$ between bornological coarse spaces.  We say that $f$ is

\begin{itemize}
\item  \textit{controlled} if for every $V\in \mathcal{C}$ we have $(f\times f)(V)\in \mathcal{C}'$. 
\item \textit{proper} if for every $B'\in \mathcal{B}'$ we have $f^{-1}(B') \in \mathcal{B}$. \item \textit{coarse} if it is controlled and proper. 
\item  \textit{bornologous} if for every $B\in \mathcal{B}$ we have $f(B)\in \mathcal{B}'$. 
\end{itemize}

\begin{definition} ($\mathbf{BornCoarse}$) A set equipped with compatible bornological and coarse structures is an object of the category $\mathbf{BornCoarse}$ of bornological coarse spaces. A morphism between bornological coarse spaces is a coarse map. 
\end{definition}

Consider a map between proper metric spaces with the induced bornological coarse structures $f: (X,\mathcal{C}_{d},\mathcal{B}_{d})\rightarrow (Y,\mathcal{C}_{d'}, \mathcal{B}_{d'})$. $f$ is proper iff it is metrically proper. $f$ is controlled iff for every $R>0$ there exists an $S>0$ such that 
\begin{align*}
d_{X}(x,y)< R \implies d_{Y}(f(x),f(y))<S
\end{align*}
A controlled map is in this case automatically bornologous. \\

If $X$ is additionally a path metric space, then by a standard Arzelà–Ascoli argument, it is geodesic. Therefore $f$ being controlled is equivalent to its being \textit{large-scale Lipschitz}: there exists a $L>0$ such that 
\begin{align*}
d_{Y}(f(x),f(y))< L d_{X}(x,y)+L 
\end{align*}
for all $x,y\in X$. Note that this is true even if we remove the assumption of $X$ being proper, since a path metric space is always quasi-geodesic (see Definition $1.4.10$ in \cite{nowak2023large}).  \\

$\mathbf{BornCoarse}$ enjoys a number of nice categorical properties (see Chapter $2$ of \cite{bunke2020homotopy}), most of which we will ignore. The two constructions we do need are (some) colimits and products:\\

Let $(X,\mathcal{C},\mathcal{B})$ be a bornological coarse space and $V\in \mathcal{C}$ an entourage. We consider the bornological coarse space 
\begin{align*}
X_V:= (X,\mathcal{C}\langle \{V\} \rangle, \mathcal{B})
\end{align*}
where $\mathcal{C}\langle \{V\} \rangle$ denotes the coarse structure generated by $\{V\}$, ie. the minimal coarse structure containing $\{V\}$. In $\mathbf{BornCoarse}$ we have the isomorphism 
\begin{align*}
X\cong \colim_{V\in \mathcal{C}} X_{V}
\end{align*}

Let $(X,\mathcal{C},\mathcal{B})$ and $(Y,\mathcal{C}',\mathcal{B}')$ be bornological coarse spaces. We define 
\begin{align*}
(X,\mathcal{C},\mathcal{B}) \otimes (Y,\mathcal{C}',\mathcal{B}'):= (X\times Y,\mathcal{C}\times \mathcal{C}', \mathcal{B}\times \mathcal{B}')
\end{align*}
where $\mathcal{C}\times \mathcal{C}' := \mathcal{C}\langle \{V\times V' \,|\, V\in \mathcal{C} \text{ and } V'\in \mathcal{C}'\}\rangle $ and analogously for $\mathcal{B}\times \mathcal{B}'$. Note that if $X,Y$ are proper metric spaces, then the bornology and coarse structure of $X\otimes Y$ are induced by the product metric on $X\times Y$. \\

Let $f_{0},f_{1}: X\rightarrow Y$ be two coarse maps between coarse bornological spaces. We say that $f_0$ and $f_1$ are \textit{close} if $f_0\times f_1 (\Delta)$ is an entourage of $Y$. We say that $f: X\rightarrow Y$ is a \textit{coarse equivalence} if there exists a morphism $g: Y\rightarrow X$ in $\mathbf{BornCoarse}$ such that $f\circ g$ and $g\circ f$ are close to the respective identities. $X$ and $Y$ are called \textit{coarsely equivalent} if there is a coarse equivalence $f: X\rightarrow Y$.  \\

We wish to define when two coarse maps $f_0:X\rightarrow Y$ and $f_{1}: X\rightarrow Y$ are coarsely homotopic. The idea is that coarse homotopies have to end eventually, but the end will be allowed to depend on the given point in the bornological coarse space (and may go to infinity as one goes to infinity). This is measured by a controlled, bornologous map $p: X \rightarrow [1,\infty)$. 

\begin{definition} Let $X$ be a bornological coarse space and let $p: X\rightarrow [1,\infty)$ be a controlled, bornologous map. We define the cylinder
\begin{align*}
I_p(X) := \{(x,t)\in X \times [0,\infty)\,|\, t\leq p(x)\}
\end{align*}
with inclusions $i_0: X\rightarrow I_p(X)$ and $i_1: X\rightarrow I_p(X)$ defined by $i_0(x) = (x,0)$ and $i_1(x) = (x,p(x))$. We give $I_{p}(X)$ the coarse structure and bornology induced from $X \otimes [0,\infty)$.
\end{definition}

Let $X$ and $Y$ be bornological coarse spaces. 
\begin{itemize}
\item A \textit{coarse homotopy} is a coarse map $H: I_p(X)\rightarrow Y$ for some controlled, bornologous map $p: X\rightarrow [1,\infty)$.
\item We call coarse maps $f_0:X\rightarrow Y$ and $f_1: X\rightarrow Y$ \textit{coarsely homotopic} if there is a coarse homotopy $H: I_p(X)\rightarrow Y$ such that $f_0 = H\circ i_0$ and $f_1 = H\circ i_1$. 
\item We say that $f: X\rightarrow Y$ is a \textit{coarse homotopy equivalence} if there exists a morphism $g: Y\rightarrow X$ in $\mathbf{BornCoarse}$ such that $f\circ g$ and $g\circ f$ are coarsely homotopic to the respective identities. $X$ and $Y$ are called \textit{coarsely homotopy equivalent} if there is a coarse homotopy equivalence $f: X\rightarrow Y$. 
\end{itemize}

\begin{lemma} (Theorem $2.4$ in \cite{mitchener2020coarse}) The notion of two coarse maps being coarsely homotopic is an equivalence relation.
\end{lemma}

Note that we use a slightly different definition of coarse homotopy to \cite{mitchener2020coarse}. They let $q: X\rightarrow [0,\infty)$ be a proper, controlled map and define the coarse cylinder as
\begin{align*}
I'_{q}(X): = \{(x,t)\in X\times [0,\infty)\,|\,t\leq q(x) +1\}
\end{align*} 
with inclusions $i_{0}: X\rightarrow I'_{q}(X)$ and $i_1: X\rightarrow I'_{q}(X)$ defined by $i_{0}(x) = (x,0)$ and $i_{1}(x,q(x)+1)$. It is clear by letting $q = p-1$ the two definitions are equivalent. The reason for our definition is so that $I_{p}(c[0,1]^n)$ can be later identified with $c[0,1]^{n+1}$. The definition in \cite{mitchener2020coarse} also requires that $p$ is proper, but this condition is never used in the proof that coarse homotopy is an equivalence relation, nor in any subsequent discussion. \\

We remark that it \textit{is} necessary to assume that the map $p$ defining a coarse cylinder is controlled, as well as bornologous, in order to work with general bornological coarse spaces. This is because \cite{mitchener2020coarse} use the definition that a set $B$ is bounded iff the inclusion $B\rightarrow X$ is close to a constant map, ie. there exists a point $x \in X$ such that $B \times \{x\}$ is an entourage. With this bornology, controlled maps are automatically bornologous. To see this, let $f: X\rightarrow Y$ be a controlled map into an arbitrary bornological coarse space, $B$ a bounded set in $X$. We have that $(f\times f)(B \times \{x\}) = f(B) \times \{f(x)\} = V \in \mathcal{C}_{Y}$ is an entourage. By definition $f(B) = V[\{f(x)\}]$. This is a bounded set since a bornology contains all one-point sets, and the bornology on $Y$ is compatible with $\mathcal{C}_{Y}$. \\

There is also another definition of coarse homotopy introduced by Bunke and Engel in \cite{bunke2020homotopy}, but this is also equivalent to ours above (see Remark $2.8$ in  \cite{weighill2025liftingcoarsehomotopies} for details). \\

The notion of maps being coarsely homotopic is closed under pre- and post-composition with coarse, bornologous maps:

\begin{prop} \label{properties} Let $X,Y,Z$ be bornological coarse spaces with controlled and bornologous projections $p_{X},p_{Y}$. let $f: X\rightarrow Y$ and $g: X\rightarrow Z$ be coarse. Assume additionally that $f$ is bornologous.
\begin{enumerate}
\item If $f\simeq f'$ then $g\circ f\simeq g\circ f'$
\item If $g\simeq g'$ then $g\circ f\simeq g'\circ f$
\end{enumerate}
\end{prop}

\begin{proof}
\begin{enumerate}
\item There exists a coarse homotopy $H: I_{p_{X}}(X)\rightarrow Y$ from $f$ to $f'$. The homotopy $g\circ H: I_{p_{X}}(X)\rightarrow Z$ is a coarse homotopy between $g\circ f$ and $g\circ f'$.
\item There exists a coarse homotopy $G: I_{p_{Y}}(Y)\rightarrow Z$ between $g$ and $g'$. We define
\begin{align*}
H: I_{p_{Y}f}(X)&\longrightarrow Z\\
(x,p_{Y}f(x)t) &\longmapsto G(f(x),p_{Y}f(x)t)
\end{align*}
$H$ is a homotopy between $G(f(x),0) = g\circ f$ and $G(f(x),p_{Y}f(x)) = g'\circ f$. Since $f$ and $p_{Y}$ are controlled and bornologous, the object $I_{p_{Y}f}(X)$ is well-defined and $H$ is coarse.
\end{enumerate}
\end{proof}

\subsection{Cones of compact metric spaces} 

In this subsection we define the cone of a compact metric space and compare this definition to various other constructions in the coarse geometry literature. Our construction relies on the fact that one can always isometrically embed a metric space $X$ into a normed space $C_{b}(X)$. 

\begin{definition} \label{Kuratowski} Let $(X,d)$ be a metric space, $x_0\in X$. Consider the Banach space $C_{b}(X)$ of all bounded functions $X\rightarrow \mathbb{R}$ equipped with the supremum norm. The Kuratowski embedding is defined as
\begin{align*}
\Psi: X&\rightarrow C_{b}(X)\\
\Psi(x)(y) &= d(x,y) - d(x_0,y)
\end{align*}
for $x,y\in X$. It is an isometry onto its image.
\end{definition}

\begin{definition} Let $X$ be a compact metric space which is isometrically embedded into a normed space $\iota: X\rightarrow L$. We define the \textit{Euclidean cone} of $X$, $cX \subset L \times \mathbb{R}$, as follows:
\begin{align*}
cX:= \{(h\iota(x),h)\in L \times [1,\infty) \,|\, x\in X\}
\end{align*}We equip $cX$ with the metric induced by the ambient norm on $L\times \mathbb{R}$ defined by $\|(h\iota(x),h)\| = \|h\iota (x)\|_{L}+ |h|$. 
\end{definition}

$cX$ is a proper metric space. The definition of $cX$ is independent up to bi-Lipschitz equivalence of the choice of normed space $L$, the particular representative in a bi-Lipschitz class of metrics on $X$, and the choice of embedding $\iota$. 

\begin{lemma}Let $\tilde{d}$ be another choice of metric on $X$, which is bi-Lipschitz equivalent to $d$. Let $\tilde{\iota}: X\rightarrow \tilde{L}$ be an isometric embedding of $(X,\tilde{d})$ into another normed space $\tilde{L}$. The cones $cX$ and $\tilde{c}X$, constructed from the data $(X,d,L,\iota)$ and $(X,\tilde{d},\tilde{L},\tilde{\iota})$ respectively, are bi-Lipschitz equivalent.  
\end{lemma}

\begin{proof}
Let $C_1,C_2>0$ be constants such that $C_1 d(x,y)\leq \tilde{d}(x,y)\leq C_{2} d(x,y)$ for all $x,y\in X$. Consider the map 
\begin{align*}
c(\tilde{\iota}\circ \iota^{-1}): cX&\rightarrow \tilde{c}X\\
(h\iota(x),h)&\rightarrow (h\tilde{\iota}(x),h)
\end{align*}
Let $(hx,h)$ and $(sy,s)$ be two points in $cX$. Assume that $h\geq s$. 
\begin{align*}
\|h\tilde{\iota}(x)-s\tilde{\iota}(y)\|_{\tilde{L}} + |h-s| \leq h\|\tilde{\iota}(x)-\tilde{\iota}(y)\|_{\tilde{L}} + |h-s| \|\tilde{\iota}(y)\|_{\tilde{L}} + |h-s|\\
\leq h\tilde{d}(x,y)+(\tilde{D}+1) |h-s| \leq hC_{2}d(x,y)+ (\tilde{D}+1)|h-s| \\
= hC_{2}\|\iota(x)-\iota(y)\|_{L} + (\tilde{D}+1)|h-s| \\
\leq C_2 (\|h\iota(x)-s\iota(y)\|_{L} + |h-s| \|\iota(y)\|_{L}) + (\tilde{D}+1) |h-s|\\
\leq (C_{2} + C_{2}D+\tilde{D} +1)(\|h\iota(x)-s\iota(y)\|_{L} + |h-s|)
\end{align*}
where $D = \sup_{x\in X} \|\iota(x)\|_{L}$ and $\tilde{D} = \sup_{x\in X} \|\tilde{\iota}(x)\|_{\tilde{L}}$. Therefore $c(\tilde{\iota}\circ \iota^{-1})$ is $(C_{2} + C_{2}D+\tilde{D} +1)$-Lipschitz. The Lipschitz constant of the inverse $c(\iota\circ \tilde{\iota}^{-1})$ is $\frac{1}{C_{1}} + D+ \frac{1}{C_{1}} \tilde{D}+1$ by an analogous calculation. 
\end{proof}

Therefore we can always assume that $L = C_{b}(X)$ and $\iota = \Psi$ is the Kuratowski embedding with $\Psi(x_0) = 0$. Additionally, by scaling the metric on $X$ such that $d(x_0,y)<\frac{1}{4}$ for all $y\in X$, we can assume that $\|\Psi(x)\| < \frac{1}{4}$ for all $x\in X$.\\

\begin{notation}
Let $(hx,h)\in cX$. Throughout this thesis, we refer to $h$ as the height variable and $x$ as the $X$-variable. We denote by $p: cX\rightarrow [1,\infty)$ the height projection and $p': cX\rightarrow X$ the projection onto $X$. We will use the notation:
\begin{align*}
c_{[a,b]} X &:=\{(hx,h)\in cX \,|\, h\in [a,b]\}\\
c_{a} X &:=  \{(hx,h)\in cX \,|\, h\in [a,\infty)\}
\end{align*} 
\end{notation}

Observe that our cones are \textit{not} topological cones, as they are truncated at height $1$. The original definitions induced by Mitchener,  Norouzizadeh and Schick in  \cite{mitchener2020coarse} define the metric cone as follows:
\begin{align*}
c'X:= \{(hx,h)\in L \times [0,\infty) \,|\, x\in X\}
\end{align*} 

The reason for defining $cX$ as we do is purely for the purpose of making calculations easier, and occasionally we add the tip back in when it suits the particular computation. Since $cX$ and $c'X$ are coarsely equivalent, this has no effect on any outcomes. (See Subsection \ref{diffdef} in the Appendix for a proof of equivalence of definitions.) \\

It is proven in Lemma $3.3$ in \cite{mitchener2020coarse} that in the restricted case where $X$ is a finite simplicial complex piecewise-linearly embedded via the map $i$ into a subset of the unit sphere of a real Hilbert space $L \cong \mathbb{R}^n$, $c'(i(X))$ is bi-Lipschitz (and therefore coarsely) equivalent to another standard definition in the literature, the metric cone with spherical base
\begin{align*}
C'(i(X)) = \{tx\,|\, t\geq 0, x\in i(X)\subset L\}
\end{align*}

Lemma $3.3$ in \cite{mitchener2020coarse} also states: if we equip $c(i(X))$ with either the subspace metric obtained as a restriction of the metric on $\mathbb{R}^n \times \mathbb{R}$, or with the induced path metric, then the identity map $\id_{c(i(X))}$ applied when changing metrics is a bi-Lipschitz homeomorphism. Additionally, the bi-Lipschitz class of $c(i(X))$ does not depend on the chosen piecewise linear embedding $i$.\\

It is important to note that these statements are \textit{not} true if $X$ is not a finite simplicial complex, as the following example shows.

\begin{example} (The Warsaw circle) Consider the Warsaw circle as a subset of $\mathbb{R}^2$:
\begin{align*}
S_{W} &:=  \{(x,sin(\frac{\pi}{x}))\,|\, 0<x<1\} \cup (\{0\} \times [-1,1]) \cup \\
& \{(x,y)\,|\, (x-\frac{1}{2})^2 +(y+\frac{15}{16})^2= (\frac{17}{16})^2, y\leq 0\}
\end{align*}

\begin{figure}
\centering
  \centering
  \includegraphics[width=0.4\linewidth]{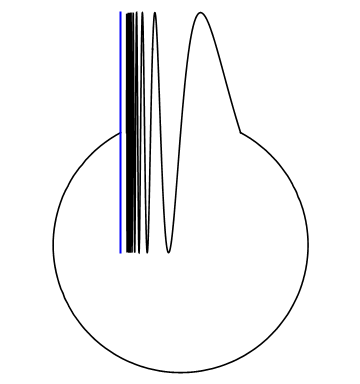}
  \caption{\cite{brazasWarsaw} The Warsaw circle.}
  \label{fig:warsaw}
\end{figure}
The cone $cS_{W}$ is a subset of $\mathbb{R}^2 \times \mathbb{R} = \mathbb{R}^3$, with metric induced from the ambient $L^1$ norm. By rescaling, we can assume that the diameter of $S_{W}\subset \mathbb{R}^2$ is less than $1$. The metric space $(cS_{W},\|\|)$ is not bi-Lipschitz equivalent to $cS_{W}$ equipped with the induced path metric $d_{path}$. \\

To see this, let $(hx,h),(hy,h) \in cS_{W}$ and $\gamma(t)= (h(t)x(t),h(t))$ a path between them. Call $\tau(t) = (x(t),1)$ the projection of $\gamma$ to height $1$. This is a path between $(x,1)$ and $(y,1)$ which has length at most the length of $\gamma$:
 \begin{align*}
\|\gamma'(t)\| = \|h'(t)x(t)+h(t)x'(t)\| + |h'(t)| \geq h(t)\|x'(t)\|+ (1-\|x(t)\|)|h'(t)| \geq \|x'(t)\|\\ = \|\tau'(t)\|\\
\implies \int \|\gamma'(t)\| dt \geq  \int \|\tau'(t)\| dt
 \end{align*}

 The distance between $(hx,h),(hy,h)$ with respect to the path metric on $cS_{W}$ therefore satisfies
 \begin{align*}
d_{path}((hx,h), (hy,h)) \geq d^{path}_{S_{W}} (x,y)
 \end{align*}
Consider the identity map $\id: (cS_{W}, \|\|) \rightarrow (cS_{W},d_{path})$. We will show that $\id$ is not coarse. Consider the metric entourage $D_{2}$ and the following sequence of pairs 
\begin{align*}
(x_{h}, y_{h})=((0,h), (h(\frac{1}{h}, 0),h)) \in D_{2} \subset cS_{W} \times cS_{W}
\end{align*}
for $h\in \mathbb{N}$. Let $\eta^{h}$ denote the shortest path in $S_{W}$ from $(0,0)$ to $(\frac{1}{h},0)$ with length $l(\eta^{h})$. Since we have to travel "around" the circle, and the length of the curve $(x,sin(\frac{\pi}{x}))$ approaches infinity as $x\rightarrow 0$, we get that
\begin{align*}
d_{path}(x_{h},y_{h}) \geq l(\eta^{h}) \xrightarrow{h\rightarrow \infty} \infty
\end{align*}
This shows that the identity map is not even controlled, much less a bi-Lipschitz equivalence.  
\end{example}

This example shows that as soon as we leave the world of finite simplicial complexes, strange things can occur. It turns out that the Warsaw circle is not so badly behaved after all: later on we will show that $cS_{W}$ is actually coarsely homotopy equivalent to $cS^1 = \mathbb{R}^2$. We will also show that for any compact metric space $X$, $cX$ is coarsely homotopy equivalent to a path metric space $M_{X}$, but as seen above, in many cases this is not the path metric space induced from the norm on $C_{b}(X) \times \mathbb{R}$. \\

There is yet another, more abstract construction of cones over uniform bornological coarse spaces by Bunke and Engel, which agrees (up to coarse equivalence) with our definitions when we restrict to metric spaces. \\

Let $X$ be a set equipped with a uniform structure $\mathcal{U}$, and a coarse structure $\mathcal{C}$. We say that $\mathcal{U}$ is compatible with $\mathcal{C}$ if $\mathcal{U}\cap \mathcal{C} \neq \emptyset$. 

\begin{definition} ($\mathbf{UBC}$) A bornological coarse space with an additional compatible uniform structure is called a \textit{uniform bornological coarse space}. We let $\mathbf{UBC}$ denote the category of uniform bornological coarse spaces and proper, controlled, and uniformly continuous maps. 
\end{definition}

As an example, let be $X$ be a compact metric space. It has unique coarse, uniform, and bornological structures induced by the metric which are mutually compatible. Continuous maps between compact metric spaces are always proper, controlled, and uniformly continuous. \\

Let $X$ be a uniform bornological coarse space and let $\mathcal{T}$ denote the uniform structure of $X$. We give $\mathcal{P}(X\times X)$ and $\mathcal{T}$ a partial order which is the opposite of the inclusion relation. An order preserving map $\phi: [0,\infty)\rightarrow \mathcal{P}(X\times X)$ is called $\mathcal{T}$-admissible if for every uniform entourage $U\in \mathcal{T}$ there exists a $t>0$ such that $\phi(s)\subset U$ for all $s\geq t$. 

\begin{definition} (Definition $8.5$ in \cite{bunke2020coarse}) Let $X$ be a uniform bornological coarse space and let $\phi:[0,\infty)\rightarrow \mathcal{P}(X\times X)$ be a $\mathcal{T}$-admissible function. We let $(\tilde{\mathcal{O}}_{\phi}(X), \mathcal{C}_{\phi},\mathcal{B})$ be the bornological coarse space defined as follows:
\begin{enumerate}
\item The underlying set of $\tilde{\mathcal{O}}_{\phi}(X)$ is $[0,\infty) \times X$.
\item The bornology $\mathcal{B}$ of $\tilde{\mathcal{O}}_{\phi}(X)$ is generated by the subsets $[0,n] \times B$ for all $n$ in  $\mathbb{N}$ and bounded subsets $B$ of $X$.
\item The coarse structure  $\mathcal{C}_{\phi}$ of $\tilde{\mathcal{O}}_{\phi}(X)$ is generated by entourages of the form $V\cap U_{\phi}$, where $V$ is a coarse encourage of $[0,\infty)\otimes X$ 
\begin{align*}
U_{\phi}:= \{((h,x),(s,y))\in ([0,\infty) \times X) \times ([0,\infty) \times X) \,|\,(x,y)\in \phi(\max\{h,s\})\}
\end{align*}
 \end{enumerate}
 \end{definition}

\begin{lemma} Let $X$ be a compact metric space, $R>0$. Let $\phi_{R}$ denote the $\mathcal{T}$-admissible function defined as $\phi_{R}(t) = D_{R/t}$ for $t>0$ and $\phi_{R}(0)=X\times X$. $(c'X, \mathcal{C}_{d}, \mathcal{B}_{d})$ is coarsely equivalent to $(\colim_{\phi_{R}}\tilde{\mathcal{O}}_{\phi_{R}}(X), \colim_{R} \mathcal{C}_{\phi_{R}},\mathcal{B})$.
\end{lemma}

\begin{proof} 
$\colim_{\phi_{R}}\tilde{\mathcal{O}}_{\phi_{R}}(X)$ has the same underlying set and bornology as each $\tilde{\mathcal{O}}_{\phi_{R}}(X)$. The coarse structure of $\colim_{\phi_{R}}\tilde{\mathcal{O}}_{\phi_{R}}(X)$ is generated by all the entourages $\{V\cap U_{\phi_{R}}\}_{V,R}$ where $V$ is a coarse entourage of $[0,\infty)\otimes X$ and $R> 0$. \\

We identify points in $\colim_{\phi_{R}}\tilde{\mathcal{O}}_{\phi_{R}}(X)$ with their images in $c'X$ under the map
\begin{align*}
\id: \colim_{\phi_{R}}\tilde{\mathcal{O}}_{\phi_{R}}(X) &\rightarrow c'X\\
(h,x) &\mapsto (hx,h)
\end{align*}
It suffices to show that for all $R>0$ there exists an $S>0$ such that $D_{R}\in \mathcal{C}_{\phi_S}$, and that for all pairs $\{V,R\}$ with $V$ a coarse entourage of $[0,\infty)\otimes X$ and $R>0$, there exists an $S>0$ such that $V\cap U_{\phi_{R}}\subset D_{S}$.\\

Let $(hx,h),(sy,s) \in D_{R}$ be two points in $c'X$. We have that $|h-s|+\|hx-sy\| <R$. Assume that $h\geq s$. We have
\begin{align*}
\|x-y\| \leq \frac{1}{h} (\|hx-sy\| + |h-s|\|y\|) < \frac{R}{h}
\end{align*}
since we can assume $\|y\|<\frac{1}{2}$ by compactness. This shows that $D_{R}\subset U_{\phi_{R}}$. We can let $V = D_{R} \times D_{R}$. Therefore $D_{R}\in \colim_{R} \mathcal{C}_{\phi_{R}} $.\\

Conversely, let $V = D_{R'} \times D_{R'}$ be a base coarse entourage of $[0,\infty) \otimes X$ for some $R'>0$. Let $R>0$. 
\begin{align*}
V\cap U_{\phi_{R}} = \{((h,x),(&s,y))\in ([0,\infty)\times X)^2 \,|\, \|x-y\| < \min\{\frac{R}{\max\{h,s\}}, R'\}, |h-s|< R'\} \\
d_{c'X} ((hx,h),(sy,s)) &= |h-s|+\|hx-sy\| \leq h\|x-y\| + 2|h-s| \|y\| < R+ 2R'
\end{align*}
Therefore $V \cap U_{\phi_R} \subset D_{R+2R'}$ and so is a coarse entourage of $c'X$. \\

It is clear that $\mathcal{B}_{d}\subset \mathcal{B}$. Suppose that $(hx,h),(sy,s)\in [0,n]\times B \in \mathcal{B}$. Again, we can assume that $B$ has diameter less than $1$.  We have
\begin{align*}
|h-s|+\|hx-sy\| \leq |h-s|+h\|x-y\|+|h-s|\|y\| < 3n
\end{align*}
hence, $[0,n] \times B$ is bounded with respect to the metric on $c'X$. 
\end{proof}

\begin{remark} If the metric on $c'X$ obtained from the ambient norm is coarsely equivalent to $c'X$ with the induced path metric (eg. if $X$ a finite simplicial complex), then the coarse structure $\mathcal{C}_{d}$ is generated by one entourage $D_{1}$. In this case, the proof of the lemma shows that $c'X$ is coarsely equivalent to $\tilde{\mathcal{O}}_{\phi_1}(X)$. 
\end{remark}

\begin{lemma} The cone $\colim_{\phi_{R}}\tilde{\mathcal{O}}_{\phi_{R}}(X)$  is coarsely homotopy equivalent to $\tilde{\mathcal{O}}_{\phi_1}(X)$.  
\end{lemma}

\begin{proof}We first prove that the coarse structure $\mathcal{C}_{\phi_1}$ of $\tilde{\mathcal{O}}_{\phi_1}(X)$ contains all sets of the form $B \times B$, where $B$ is a bounded set. By definition $B\subset [0,n] \times X$ for some $n\in [0,\infty)$. Consider $(h,x),(s,y)\in B$. Let $V_{n} = ([0,n] \times X)^2$. We have 
\begin{align*}
((h,x),(0,x)), ((0,x), (0,y)), ((0,y),(s,y)) \in V_{n} \cap U_{\phi_1}
\end{align*}  
which gives us that $((h,x),(s,y))\in (V_{n} \cap U_{\phi_1})^{\circ 3}$. We show that the logarithm defined as 
\begin{align*}
Log:  \colim_{\phi_{R}}\tilde{\mathcal{O}}_{\phi_{R}}(X) &\rightarrow  \tilde{\mathcal{O}}_{\phi_1}(X)\\
(t,x) &\mapsto (\log (t), x)
\end{align*}
for $t\geq 1$ and $(t,x)\mapsto (0,x)$ for $t\leq 1$, is a coarse homotopy equivalence. We use the same symbol $\log$ to denote the extension to $t\leq 1$. $Log$ is proper since $Log^{-1}([0,n]\times X) = [0,\exp(n)] \times X$. To show it is controlled, let $V\cap U_{\phi_{R}}$ be an entourage in $ \colim_{\phi_{R}}\tilde{\mathcal{O}}_{\phi_{R}}(X)$, $V=D_{R'}\times D_{R'}$. There exists a $t_{R}>1$ such that $\frac{R}{t}<\frac{1}{\log t}$ for all $t> t_R$. Consider the bounded set $B_{t_R}: = \{(t,x) \in [0,\infty) \times X\,|\, t\leq t_R\}$. Let $((h,x),(s,y)) \in V\cap U_{\phi_{R}}$ be in the complement of $B_{t_R} \times B_{t_{R}}$. We have that
\begin{align*}
\|x-y\| &< \min\{\frac{R}{\max\{h,s\}},R'\} < \min\{\frac{1}{\log{\max\{h,s\}}},R'\}\\
&=\min\{\frac{1}{\max\{\log(h),\log(s)\}},R'\}\\
|\log(h)-\log(s)|&< |h-s|<R'
\end{align*}
Therefore $(Log(h,x),Log(s,y))\in V\cap U_{\phi_1}$. For points in $B_{t_{R}}$we have $(Log \times Log) (B_{t_{R}}\times B_{t_{R}})\subset  (B_{\log(t_{R})}\times B_{\log(t_{R})})$. This shows that $(Log \times Log)(V\cap U_{\phi_{R}})\subset (V\cap U_{\phi_1})\cup (B_{\log(t_{R})}\times B_{\log(t_{R})})\in \mathcal{C}_{\phi_{1}}$. Therefore $Log$ is coarse. \\

We show that $Log$ is a coarse homotopy inverse to the set-wise identity map $\id: \tilde{\mathcal{O}}_{\phi_1}(X) \rightarrow \colim_{\phi_{R}}\tilde{\mathcal{O}}_{\phi_{R}}(X)$. Let $p$ be the projection to the height variable (this is a controlled, bornologous map). Define the homotopy
\begin{align*}
H: I_{p}(\colim_{\phi_{R}}\tilde{\mathcal{O}}_{\phi_{R}}(X))&\rightarrow \colim_{\phi_{R}}\tilde{\mathcal{O}}_{\phi_{R}}(X) \\
(h,x,ht) &\mapsto (th + (1-t)\log h, x)
\end{align*}
$H$ is coarse: let $W =((D_{R}\times D_{R}) \cap U_{\phi_{R}}) \times D_{R}$ be a base entourage in $ I_{p}(\colim_{\phi_{R}}\tilde{\mathcal{O}}_{\phi_{R}}(X))$. Let $((h,x,ht),(s,y,sl))\in W$. We have that
\begin{align*}
&|(1-t) \log h +ht - (1-l)\log s - ls| \leq |(1-t)\log h - (1-l) \log s| + |ht-sl|\\
&\leq |(1-t) \log h - (1-t) \log s| + |l-t| \log s  + |ht-sl|\\
&\leq 2|h-s| + 2|ht-sl| \\
&<4R\\
&\|x-y\| < \min\{\frac{R}{\max\{h,s\}},R\} \leq \min\{\frac{R}{\max\{(1-t) \log h +ht,(1-l)\log s - ls\}},4R\}
\end{align*}
This shows that $(H(h,x,ht), H(s,y,sl))\in (D_{4R} \times D_{4R})\cap U_{\phi_{R}}$. Therefore $H$ is controlled.\\

Suppose that $[0,n] \times X$ is a bounded set in $\colim_{\phi_{R}}\tilde{\mathcal{O}}_{\phi_{R}}(X)$. We have that 
\begin{align*}
H^{-1}([0,n] \times X) = \{(h,x,ht)\in  I_{p}(\colim_{\phi_{R}}\tilde{\mathcal{O}}_{\phi_{R}}(X))\,|\, ((1-t) \log h + th, x)\in [0,n] \times X\}\\
\subset [0,\exp(n)] \times X \times [0,\exp(n)]
\end{align*}
An analagous calculation shows that $Log \circ \id: \tilde{\mathcal{O}}_{\phi_1}(X)\rightarrow \tilde{\mathcal{O}}_{\phi_1}(X)$ is also coarsely homotopic to the identity. 
\end{proof}

To help us actually compute some coarse homotopy groups, we need a specific, geometric simplicial structure on $c[0,1]^n$, where $[0,1]^n$ is viewed as a subset of $\mathbb{R}^n$. Therefore, we introduce some notions from PL topology and geometry, and define the simplicial cone of a finite simplicial complex embedded into Euclidean space. Doing this means we have some useful theorems from the PL world at our disposal.  \\

\begin{definition} \begin{itemize}
\item A \textit{realisation} of a simplex $\sigma:= [v_{0},\dots, v_{m}]$ is a set in $\mathbb{R}^N$ which is the convex hull of $(m+1)$-affinely independent vertices $v_{0},\dots, v_{m} \in \mathbb{R}^N$. 
\item Let $K$ be an abstract simplicial complex and let $\iota_K: vert K \rightarrow \mathbb{R}^N$ be a map defined on the vertices of $K$ such that for every simplex $\sigma = [v_{0},\dots, v_{m}]$ in $K$, the set $\{\iota_K(v_{0}),\dots,\iota_K(v_{m})\}$ is affinely independent in $\mathbb{R}^N$. In this case we can extend $\iota_K$ linearly in barycentric coordinates to the geometric realisation of $\sigma$ by sending its image to the convex hull of $\{\iota_K(v_{0}),\dots, \iota_K(v_{m})\}$.  We obtain a map  $\iota_{|K|}: |K|\rightarrow \mathbb{R}^N$, which we call a \textit{simplicial embedding} if it is injective. We identify $K$ with its geometric realisation $|K|$ and its image under $\iota_{|K|}$ and we say that $K= |K|= \iota_{|K|}(|K|)$ is \textit{simplicially embedded}.
\item A map between simplicially embedded simplicial complexes, $K,K'$, with embeddings $\iota_{|K|},\iota_{|K'|}: |K|,|K'|\rightarrow \mathbb{R}^N$, is a map $f: \iota_{|K|}(|K|) \rightarrow \iota_{|K'|}(|K'|)$. $f$ is called \textit{simplicial} if the composition 
\begin{align*}
|K|\xrightarrow{\iota_{|K|}} \iota_{|K|}(|K|) \xrightarrow{f} \iota_{|K'|}(|K'|) \xrightarrow{(\iota_{|K'|})^{-1}} |K'|
\end{align*}
is a simplicial map. 
\end{itemize}
\end{definition}

Given a simplicially embedded simplicial complex $K$, we can give $|K|$ the pullback metric $\iota^{*}_{|K|} d_{std}$ from the norm on $\mathbb{R}^N$, or the induced path metric $d_{path}$. From now on, if $K$ is simplicially embedded, we freely switch notation between $K,|K|$ and $\iota_{|K|}(|K|)$, eg. we simply write $f:K\rightarrow K'$ to denote a map between simplicially embedded simplicial complexes. Similarly, for any abstract simplicial complex $K$ we identify $K$ with its geometric realisation $|K|$.  

\begin{definition} A simplicial complex $K$ is called \textit{locally ordered} if there is a partial order on its vertices which restricts to a total order (or \textit{orientation}) on the vertices of each simplex of $K$. 
\end{definition}

For a locally ordered simplicial complex $K$ recall that there is a canonical triangulation of $K\times [0,1]$ with the simplices in $K\times \{0\}$ and $K\times \{1\}$ coming from the triangulation of $K$. In addition, for any oriented simplex $\langle v_{0},\dots v_{m}\rangle$ of $K$ and $0\leq j\leq m$ we have a new simplex spanned by $\{v_{0}, \dots v_{j},v_{j}',\dots v_{m}'\}:=\{(v_0,0),\dots,(v_j,0),(v_j,1),\dots,(v_m,1)\}$ where the primes indicate copies of vertices at $K \times \{1\}$. $K\times [0,1]$ has a local order inherited from $K$: we let $(v_i,k)\leq(v_j,l)$ if and only if $k\leq l$ and $i\leq j$. This restricts to a total order on each simplex of $K\times [0,1]$. Observe that if $\langle (v_0,h_0),\dots, (v_m, h_m) \rangle$ is an oriented simplex of $K \times [0,1]$, then the image of its projection to $K$ spans an oriented simplex, after possibly deleting a vertex that has two pre-images. 
\begin{notation} We denote by $[v_0,\dots, v_m]$ the realisation of a simplex (which is a set in $\mathbb{R}^N$) and use $\langle v_0,\dots, v_m \rangle$ for an oriented simplex with orientation obtained from the order on the vertices. 
\end{notation}

The following defintions about subdivision are taken from \cite{mitchener2020coarse}.

\begin{figure}[H]
\centering
  \centering
  \includegraphics[width=0.5\linewidth]{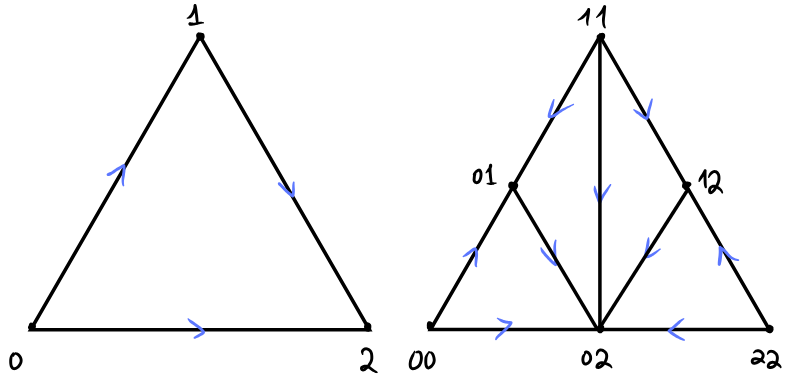}
  \caption{The $1$st standard subdivision $S(\sigma)$ (right) of $\sigma = \langle 0,1,2 \rangle$ (left) with order $0<1<2$.}
  \label{fig:ssubd0}
\end{figure}

\begin{definition} Let $\sigma:= \langle v_{0},\dots, v_{n}\rangle \subset \mathbb{R}^N$ be an oriented realised simplex. We define its \textit{standard subdivision} $S(\sigma)$ as the simplicial complex with vertices $v_{i,j}:= \frac{v_{i}+v_{j}}{2}$ for $0\leq i\leq j\leq n$. On this set of vertices we define a partial order setting $(i,j)\leq (k,l)$ if and only if $k\leq i\leq j\leq l$. The simplices of the standard subdivision are spanned by increasing sequences of vertices, making $S(\sigma)$ locally ordered. See Figure   \ref{fig:ssubd0}. \\

Given a locally ordered simplicial complex $K$, define a standard subdivision $S(K)$ by applying the standard decomposition to each simplex to obtain a simplicial decomposition of the whole simplicial complex. 
\end{definition}

\begin{figure}
\centering
  \centering
  \includegraphics[width=0.4\linewidth]{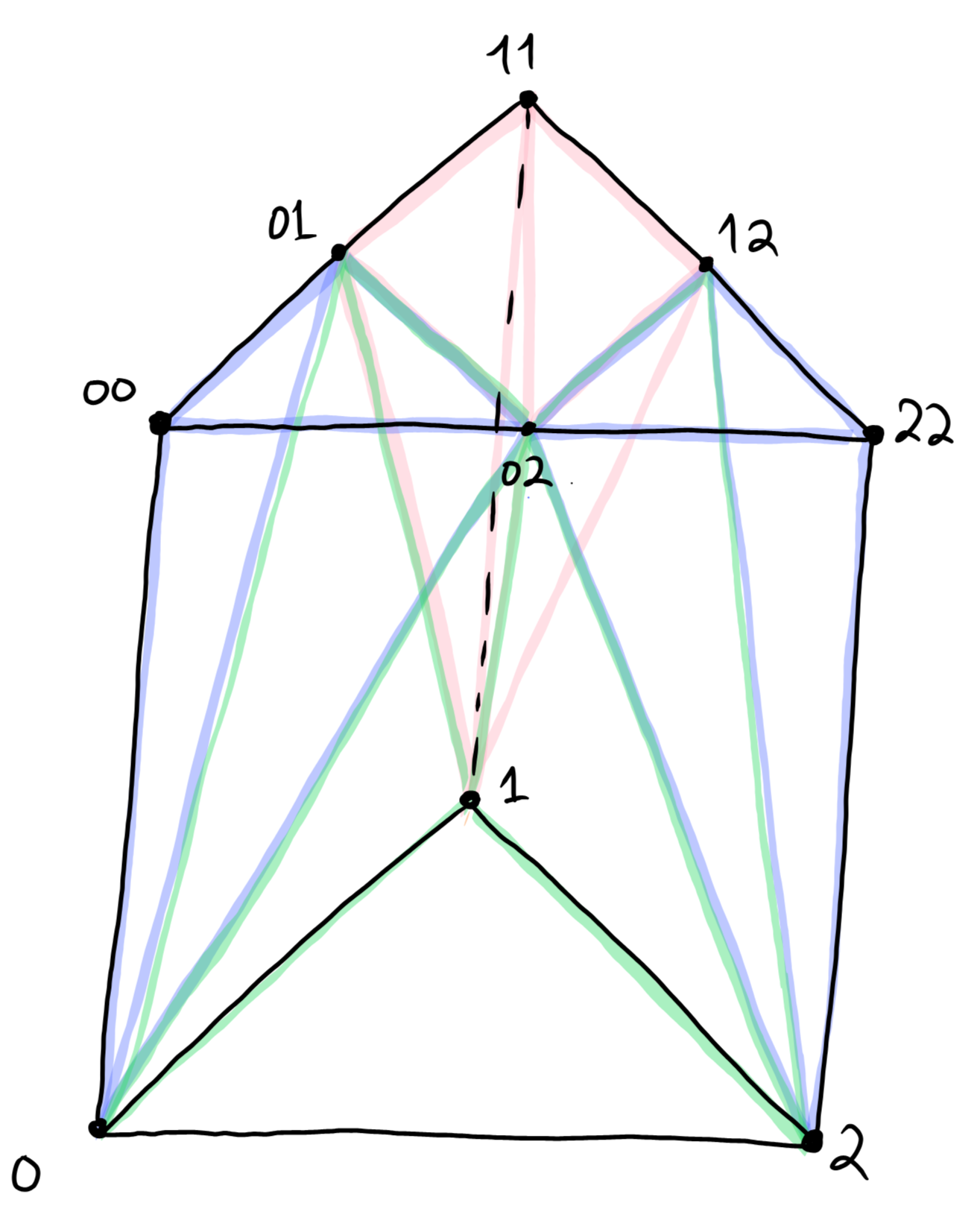}
  \caption{The standard product subdvision of $\langle 0,1,2\rangle \times [0,1]$. Red: simplicies that have $1$ as a vertex. Blue: simplices that have $0$ or $2$ as a vertex. Green: simplices that have $\langle 0,1\rangle$ or $\langle1,2 \rangle$ as an edge. There is an additional uncoloured simplex $\langle 0,1,2,02\rangle $.}
  \label{fig:egS}
\end{figure}

We also require a triangulation of $S(K) \times [0,1]$: 

\begin{definition} Let $K$ be a locally ordered simplicial complex. We define the \textit{standard product subdivision} of $K\times [0,1]$ as the triangulation of $K \times [0,1]$ which restricts to the given triangulation on $K\times \{0\}$ and to the standard subdivision $S(K) \times \{1\}$ on the other end. The additional simplices are as follows: whenever $u_{l}\leq \dots \leq u_{0} \leq v_{0} <\dots< v_{k} \leq w_{0}\leq \dots \leq w_{l}$ are vertices of a simplex of $K$ such that $(u_0,w_0)< \dots <(u_{l},w_{l})$ in the standard subdivision, we get a simplex in the standard product subdivision of $K \times [0,1]$ spanned by 
\begin{align*}
(v_0,0),\dots,(v_{k},0), ((u_0,w_0),1), \dots, ((u_l,w_l),1)
\end{align*}
See Figure \ref{fig:egS}.
\end{definition}

\begin{remark} The convex hull of $[(u_0,w_0), \dots, (u_{l},w_{l})] \times \{1\}$ and $[v_0,\dots,v_{k}] \times \{0\}$ is precisely the union of the convex hulls of $[(u_0,w_0), \dots (u_l,w_{l})] \times \{1\}$ and the simplices in the standard subdivision of $[v_0,\dots,v_{k}] \times \{0\}$. 
\end{remark}

\begin{lemma} Let $K$ be a locally ordered simplicial complex and $I = [0,1]$ the unit interval with one edge and two vertices, with order $0<1$. Equip $K\times I$ with the product simplicial structure. Then 
\begin{align*}
S^k(K\times I) = S^k(K) \times S^k(I)
\end{align*}
for all $k\in \mathbb{N}_{0}$. 
\end{lemma}

\begin{proof} Assume the statement is true for all $j\leq k$. Note that $S^{k+1}(I)$ can be identified with $S^{k}(I)\cup S^{k}(I^{-})$ where $I^{-1}$ is $I=[0,1]$ given the opposite orientation. Here, $\{1\}$ in the first copy is glued to $\{0\}$ in the second. We have 
\begin{align*}
S^{k+1}(K\times I) = S^{k} (S(K) \times S(I)) = S^{k}(S(K) \times I \cup S(K) \times I^{-1})\\
= S^{k+1}(K) \times S^{k}(I) \cup S^{k+1}(K) \times S^k(I^{-1})  = S^{k+1}(K) \times S^{k+1}(I)
\end{align*}
It remains to prove $S(K\times I) = S(K) \times S(I)$. It suffices to restrict  $\sigma \times I$, where $\sigma = \langle 0,1,\dots, m-1, m\rangle$ is a simplex in $K$, since subdivision commutes with inclusion of subcomplexes. Denote by $\langle 0',1',\dots, (m-1)',m'\rangle$ the copy of $\sigma$ at level $\{1\}$. See Figure \ref{fig:ssubd}.\\

\begin{claim} There is a copy of $S(\sigma)$ on the slices $\{0\},\{\frac{1}{2}\},\{1\}$. 
\end{claim}

The claim for $\{0\},\{1\}$ is obvious. Suppose that $(u,v)$ is a vertex in $S(\sigma)$. The vertex $(u,v')$ is a vertex in $S(\sigma \times I)$ on the slice $\{\frac{1}{2}\}$. We let $f: S(\sigma)\rightarrow S(\sigma \times I)$ be the extension of this vertex map to simplices: if $\langle (u_{0},v_{0}),\dots, (u_{l}, v_{l})\rangle$ is a simplex in $S(\sigma)$ we have:
\begin{align*}
u_{l}\leq u_{l-1}\leq \dots \leq u_{0} \leq v_{0} \leq \dots \leq v_{l-1} \leq v_{l}
\end{align*}
By the definition of the order on $\sigma \times I$ we have $u_{l}\leq u_{l-1}\leq \dots \leq u_{0} \leq v'_{0} \leq \dots \leq v'_{l-1} \leq v'_{l}$, so $\langle (u_{0},v'_{0}),\dots, (u_{l}, v'_{l})\rangle$ is a simplex in $S(\sigma \times I)$ on the slice $\{\frac{1}{2}\}$. \\

We now construct simplicial map $g: S(\sigma) \times I \rightarrow S(\sigma \times I)$ which is a bijection onto the subcomplex $S(\sigma) \times [\frac{1}{2}, 1]$. Let $\tau = \langle ((u_{0},v_{0}),0),\dots, ((u_{j}, v_{j}),0), ((u_{j},v_{j}),1), \dots, ((u_{l}, v_{l}),1)\rangle$ be a simplex in $S(\sigma) \times I$. We define
\begin{align*}
g(\tau)= \langle (u'_{0},v'_{0}),\dots, (u'_{j}, v'_{j}), (u_{j},v'_{j}), \dots, (u_{l}, v'_{l}) \rangle
\end{align*}
This is a well-defined simplicial map because 
\begin{align*}
u_{l}\leq \dots \leq u_{j}\leq u'_{j} \leq \dots \leq u'_{0} \leq v'_{0} \leq \dots \leq v'_{l}
\end{align*}
It is obviously a bijection onto the image. The construction for the subcomplex $S(\sigma) \times [0,\frac{1}{2}]$ is analogous. 
\begin{align*}
\tilde{g}(\tau)= \langle (u_{0},v_{0}),\dots, (u_{j}, v_{j}), (u_{j},v'_{j}), \dots, (u_{l}, v'_{l}) \rangle
\end{align*}
Therefore $S(K\times I) = S(K) \times S(I)$.

\begin{remark} $S(I)$ has a canonical order $0<1<\frac{1}{2}$ and the product simplicial structure of $S(K) \times S(I)$ is with respect to this order. That is why $g$ sends $S(\sigma)\times \{0\}$ to $S(\sigma \times \{1\})$ and $S(\sigma)\times \{1\}$ to $S(\sigma \times \{\frac{1}{2}\})$. 
\end{remark}

\begin{figure}
\centering
  \centering
  \includegraphics[width=0.5\linewidth]{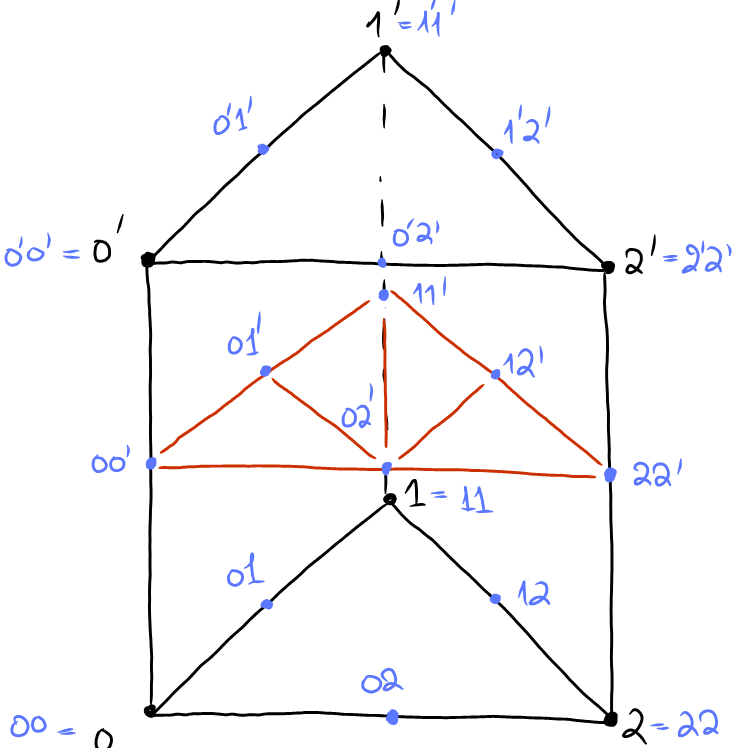}
  \caption{The $1$st standard subdivision of $\sigma \times I = \langle 0,1,2 \rangle \times I$ (with orientations and some simplices omitted). The original simplicial complex $\sigma \times I$ is denoted in black. Vertices in $S(\sigma \times I)$ are blue. The simplex $S(\sigma) \times \{\frac{1}{2}\}$ is denoted in red.}
  \label{fig:ssubd}
\end{figure}

\end{proof}

We now define the simplicial structure of $cK$, for $K$ a finite simplicial complex simplicially embedded into $\mathbb{R}^N$.

\begin{definition} \label{triang} Let $K \subset \mathbb{R}^N$ be a finite simplicial complex simplicially embedded. Write $c(K) \subset \mathbb{R}^N \times [0,\infty)$ as the union of the infinitely many copies of $K\times [0,1]$ given as $ZK(n):= \{(hx,h)\,|\, x\in K, h\in [n,n+1]\}$ for $n\in \mathbb{N}$. We define a simplicial structure on $c(K)$ as follows: use the $k$-th standard subdivision of $K$ on $nK \times \{n\}$ for $n\in \mathbb{N}$ with $2^k\leq n< 2^{k+1}$ and the product simplicial structure on $ZK(n)$ compatible with the given simplicial structure on the top and bottom. We call this simplicial structure of $c(K)$ the cone of the simplicial structure on $K$. 
\end{definition}

\begin{example} Image \ref{fig:examplecone} shows the simplicial structure of $c[0,1]$, where $[0,1]$ has one edge and orientation $0<1$.
\begin{figure}
\centering
  \centering
  \includegraphics[width=0.5\linewidth]{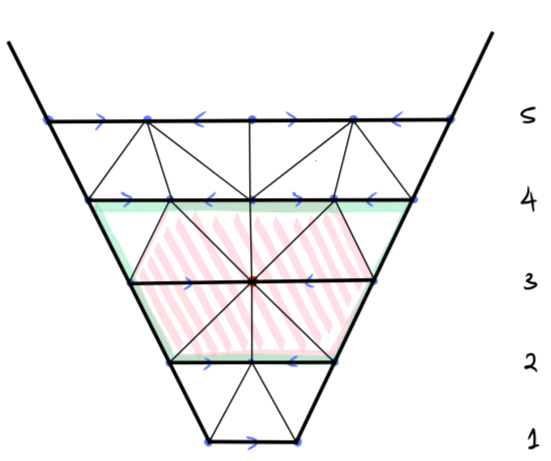}
  \caption{The cone of $[0,1]$ for $h\in [1,5]$.}
  \label{fig:examplecone}
\end{figure}
\end{example}

For $K$ a finite simplicial complex simplicially embedded into $\mathbb{R}^N$, there is an embedding of $K\times [0,1]$ into $\mathbb{R}^{N+1}$ in the obvious way. Recall that the open star of a vertex $v \in K$ consists of all the interiors of simplices which have $v$ as a vertex. We denote this by $\langle st(v) \rangle$ and its closure as $|st(v)|$. The next technical lemma will be used later on when considering homotopies of maps from $cK$.

\begin{lemma} \label{technicalcone} Let $K$ be a finite simplicial complex and $(hx,ht,h)$ be a vertex in $c(K \times [0,1])$. The closed star of $(hx,ht,h)$ is contained within $|st(hx,h)| \times I$ where $|st(hx,h)|$ denotes the star in $cK$ with the coned structure of $K$, and $I$ is an interval of length $4$ centred at $ht$. 
\end{lemma}

\begin{proof} Let $h\in \mathbb{N}$ be fixed. We have that $x \in vert S^{k}(K)$ and $t \in vert S^{k}([0,1])$ for some $k$. Since at most one of the heights $h-1,h+1$ has a different simplicial structure to $S^{k}(K\times [0,1])$ it suffices to consider $|st(x,t,0)|$ and $|st(x,t,1)|$ in the standard product subdivision of $(S^{k-1}(K) \times S^{k-1}(I)) \times [0,1]$. \\

There are $3$ cases:\\

\textbf{Case 1:} $[h-1,h+1] \subset (2^{k-1},2^{k})$ for some $k\in \mathbb{N}$.\\

The simplicial structure restricted to each integer height is just $S^{k-1}(K \times [0,1]) = S^{k-1}(K) \times S^{k-1}([0,1])$. Suppose we have a simplex 
 \begin{align*}
\sigma = (v_0,0),\dots,(v_{j},0), (v_{j},1), \dots, (v_{l},1)
\end{align*}
in $S^{k-1}(K) \times S^{k-1}([0,1]) \times [h-1,h+1]$ with $(x,t,h)$ as a vertex. By construction $\langle v_{0},\dots,v_{j},\dots, v_{l}\rangle$ is a simplex in $S^{k-1}(K) \times S^{k-1}([0,1])$ with $(x,t)$ as a vertex. Consider now the projection 
\begin{align*}
q\sigma = sp\{(qv_0,0),\dots,(qv_{j},0), (qv_{j},1), \dots, (qv_{l},1)\}
\end{align*}
to the level $S^{k-1}(K) \times \{t\}$, where $sp$ denotes the span of a set of vertices. Since the order is preserved under projection, we have that the vertices $\{qv_{0},\dots,qv_{j},\dots, qv_{l}\}$ span an oriented simplex in $S^{k-1}(K)$ (after possibly deleting a duplicate vertex). From this it follows that $q\sigma$ is a simplex in $S^{k-1}(K) \times \{t\} \times [h-1,h+1]$ with $(x,t,h)$ as a vertex. Therefore $|st(x,t,h)| \subset |st(x,h)| \times I$, where $I$ is an interval of length $\frac{2}{2^{k-1}}$. To go back to the metric on $c(K \times [0,1])$ we multiply the length of $I$ by a constant which is at most $2^{k}-1<2^{k}$. Therefore $|st(hx,ht,h)| \subset |st(hx,h)| \times I$ where $I$ is an interval of length at most $4$. \\

\textbf{Case 2:} $h = 2^{k}$ for some $k\in \mathbb{N}$. We identify the interval $[2^{k}-1,2^k]$ with $[0,1]$. \\

Let 
 \begin{align*}
\sigma = \langle (v_0,0),\dots,(v_{j},0), ((u_0,w_0),1)), \dots, ((u_l,w_l),1) \rangle
\end{align*}
be a simplex with $(x,t,1)$ as a vertex. By definition we have $u_{l}\leq \dots \leq u_{0} \leq v_{0} <\dots< v_{j} \leq w_{0}\leq \dots \leq w_{l}$. Recall that  $\langle v_0,\dots,v_{j}\rangle$ is a simplex in the product simplicial structure of $S^{k-1}(K) \times I$ for an interval $I$ of length $\frac{1}{2^{k-1}}$, and that $ \langle (u_0,w_0)), \dots, (u_l,w_l)\rangle$ is a simplex in $S^{k}(K) \times I$ for an interval of length $\frac{1}{2^k}$. Let $qv_{0},\dots,qv_{j}$ denote the image of $v_0,\dots,v_{j}$ under the projection to the level $t$. Analogously, denote by  $(qu_0,qw_0), \dots, (qu_l,qw_l)$ the projection of $(u_0,w_0), \dots, (u_l,w_l)$ to $S^{k}(X) \times \{t\}$. As before, we may have to remove a copy if two vertices project to the same point. By the canonical order on the product simplicial structure, we have that 
 \begin{align*}
qu_{l}\leq \dots \leq qu_{0} \leq qv_{0} <\dots< qv_{j} \leq qw_{0}\leq \dots \leq qw_{l}
\end{align*}

ie. the order is preserved under projection. Therefore
\begin{align*}
q\sigma = sp\{(qv_0,0),\dots,(qv_{j},0), ((qu_0,qw_0),1)), \dots, ((qu_l,qw_l),1)\}
\end{align*} 
 is a simplex in $c(K)\times \{t\}$ with $(x,t,1)$ as a vertex. If $v_{0},\dots,v_{j}$ are  contained within a single slice of $S^{k-1}(K)$, then $(u_0,w_0), \dots, (u_l,w_l)$ may span an interval of length $\frac{1}{2^{k}}$. If  $v_{0},\dots,v_{j}$ spans an interval of length $\frac{1}{2^{k-1}}$, then $t$ must be a vertex in $S^{k}(I)$ but not in $S^{k-1}(I)$, and $(u_0,w_0), \dots, (u_l,w_l)$ must be contained within $S^{k}(K) \times \{t\}$. In either case, $\sigma$ is contained in the set $q\sigma \times I\subset |st(x,1)| \times I$, where $I$ is an interval of length less than $\frac{2}{2^{k}}$ centred at $t$. The metric on $c(K\times [0,1])$ means that this interval is scaled by at most $2^{k}$: therefore $I$ has length less than $\frac{2^{k}}{2^{k-1}}=2$. \\

For simplices $\sigma$ that lie in the height interval $[2^{k},2^{k}+1]$, we reduce to the situation in case $1$. So we have $\sigma \subset |st(hx,h)| \times I$ where $I$ has length less than $4$. Putting it altogether, we obtain $|st(hx,ht,h)|\subset |st(hx,h)| \times I$, where $I$ is an interval centred at $ht$ of length less than $4$.  \\
 
\textbf{Case 3:}  $h = 2^{k}-1$ for some $k\in \mathbb{N}$.\\

This is the same as in case $2$: let $\sigma$ be a simplex with $(x,t,0)$ as a vertex. Here, $\sigma$ is contained in the set $q\sigma \times I \subset |st(x,0)| \times I$, where $I$ is an interval of length less than $\frac{2}{2^{k-1}}$ centred at $t$. The metric on $c(K \times [0,1])$ means that the interval is scaled at most $2^{k}$. Therefore $I$ has length less than $\frac{2^{k}}{2^{k-2}} = 4$. \\

For simplices that lie in the height interval $[2^{k}-2,2^{k}-1]$, we reduce to the situation in case $1$. This concludes the proof.\\

For an illustration of the star of a vertex in case $3$, see Figure \ref{fig:examplecone}. Here, $K = \{\ast\}$. The star of the red vertex $(3\ast, 3(01),3)$ is shaded in red and the set $|st(hx,h)| \times I = |st(3\ast,3)| \times I$ is outlined in green.\\
\end{proof}

\begin{remark}If the reader is unable to visualise the argument in the previous proof, we spell out now an example in excruciating detail. Suppose we are in the first subcase of case $2$. First, let all the $v$'s in $S^{k-1}(K) \times S^{k-1}(I)$ lie on one level (in the $I$-coordinate). That means that in the subdivision $S^{k}(K) \times S^{k}(I)$, we could have multiple levels. For example, we could have a simplex of the form 
 \begin{align*}
\sigma = \langle (v_0,0),\dots,(v_{j},0), ((u_0,w_0),1)), \dots,  ((u_m,w_m),1),  ((u_m,w'_m),1), \dots, ((u_l,w'_l),1)\rangle
\end{align*}
where the sequence of vertices in the original structure has the order  $u_{l}\leq \dots \leq u_{0} \leq v_{0} <\dots< v_{j} \leq w_{0}\leq \dots  \leq w_m< w'_{m}\leq \dots \leq w'_{l}$ for some $m$. The primes indicate a copy at a different level. Then when we project we get two copies of $w_{m}$, even though $w_{m}<w'_{m}$. So we delete the extra copy. Recall that by the definition of the product simplicial structure, removing the primes on the $w$'s does not change the order. This gives us the canonical orientation on the projection. Observe that the convex hull of the vertices defining $\sigma$ are subset of the convex hull of the projection $p\sigma$ plus the copy of $p\sigma$ at the primed level. This gives us the statement that we want. \\

Now if the $v$'s lie on two separate levels, then the changeover from the un-primed to the primed variables happens in the middle of the $v$'s, ie. we have a simplex of the form 
 \begin{align*}
\sigma = \langle (v_0,0),\dots,(v_{m},0), (v'_{m},0),\dots, (v'_{j},0), ((u_0,w'_0),1)), \dots, ((u_l,w'_l),1) \rangle
\end{align*}
for some $m$. We have the order $u_{l}\leq \dots \leq u_{0} \leq v_{0} <\dots<v_{m}<v'_m <\dots< v'_{j} \leq w'_{0} \leq \dots \leq w'_{l}$.  \\

Looking at this, we see that all the $u$'s are un-primed and all the $w$'s are primed. That means that the simplex $((u_0,w'_0),1)), \dots, ((u_l,w'_l),1)$ lies on one level - namely the middle between the unprimed and the primed levels. We delete the copy of $v_{m}$ after projection. 
\end{remark}

We now cite some results obtained by Michener, Norouzizadeh and Schick about cones of finite simplicial complexes embedded into Euclidean space. These will be used throughout the thesis. For proofs, refer to \cite{mitchener2020coarse}.

\begin{lemma} \label{Lbound}(Lemma $1.11.$ in \cite{mitchener2020coarse}) Let $\sigma:=[v_0,\dots,v_{n}]\subset \mathbb{R}^N$ be a realisation of an $n$-simplex, spanned by the vectors $v_0,\dots,v_n$ in general position. Let $w_0,\dots,w_n$ be vertices of a realised $k$-simplex $\tau\subset \mathbb{R}^M$. There is a unique affine linear map $f:\sigma\rightarrow \tau$ sending $v_i$ to $w_i$. The Lipschitz constant of $f$ is bounded above by $c(n,k,w)\max\{|w_i-w_j|\}$ where $c(n,k,w)$ depends on the dimensions $n$ and $k$ of the simplices and in addition on a lower bound $w$ on the width of $\sigma$ defined to be the shortest distance from any vertex of $\sigma$ to the opposite face.
\end{lemma}

\begin{definition} Two realised simplices are strongly similar if one can be obtained from the other by translation and multiplication by a positive constant. 
\end{definition}

\begin{definition} A simplicial complex $K$ is called \textit{uniformly bounded} if
\begin{enumerate} 
\item  The lengths of the edges are contained in a compact interval $[a,b]$ with $0<a<b<\infty$.
\item There are only finitely many strong similarity types of simplices in $K$.
\end{enumerate}
\end{definition}

\begin{lemma} \label{bound} (Lemma $1.16$ in \cite{mitchener2020coarse}) Let $X$ be a finite simplicial complex simplicially embedded into $\mathbb{R}^N$. There are only finitely many strong similarity types in the simplicial structure of $c(X)$. Moreover, the lengths of the edges are contained in a compact interval $[a,b]$ with $0<a<b<\infty$. In particular, there is a positive lower bound on the width of the simplices and an upper bound on the diameter. 
\end{lemma}

\subsection{Coarse homotopy groups}

In order to define coarse homotopy groups, we need a coarse analogue of a basepoint. Let $[1,\infty)$ be given the standard metric $d(t,s) = |t-s|$. 

\begin{definition} Let $X$ be a bornological coarse space. A \textit{base ray} of $X$ is a coarse map $\omega: [1,\infty)\rightarrow X$. If $Y$ is another bornological coarse space with base ray $\tau$, then a coarse map $f: X\rightarrow Y$ is called \textit{base ray preserving} if $\tau=f\circ \omega$. 
\end{definition}

\begin{definition} Let $X$ be a bornological coarse space. We define the \textit{$0$-th coarse homotopy set} $\pi_0^c(X)$ to be the set of coarse homotopy classes of coarse maps from $[1,\infty)$ to $X$.
\end{definition}

The idea of the $0$-th coarse homotopy set is that it counts the "components at infinity". 

\begin{example} $\pi_0^c(\mathbb{R}^n)$ has two elements if $n=1$ and one element if $n\geq 2$. This is not trivial and followed from the main theorem of \cite{mitchener2020coarse}.
\end{example}

If $X$ has at least one base ray $\omega$, ie. $\pi_0^c(X) \neq \emptyset$, then we can consider the pointed set $\pi^c_{0}(X,\omega)$ of coarse homotopy classes of coarse maps from $[1,\infty)$ to $X$ with base point $[\omega]$. \\

A \textit{bornological coarse pair} is a pair of bornological coarse spaces $(X,A)$ along with a bornologous, coarse map $k_A: A\rightarrow X$. 
\begin{definition} Let $(X,A)$ and $(Y,B)$ be bornological coarse pairs. A \textit{coarse map of pairs} $f:(X,A)\rightarrow (Y,B)$ is a commutative diagram 

\begin{figure}[H]
\center
\begin{tikzcd}
&A \arrow {r}{f_{|A}} \arrow{d}{k_A}
&B \arrow{d}{k_B}
  \\
&X \arrow {r}{f} 
&Y

\end{tikzcd}
\end{figure}

\end{definition}

As an example, let $X$ and $Y$ be bornological coarse spaces with base rays $\omega,\tau$. $f: (X,\omega)\rightarrow (Y,\tau)$ is base ray preserving exactly when the diagram commutes.

\begin{figure}[H]
\center
\begin{tikzcd}
&{[1,\infty)} \arrow {r}{\id} \arrow{d}{\omega}
&{[1,\infty)} \arrow{d}{\tau}
  \\
&X \arrow {r}{f} 
&Y
\end{tikzcd}
\end{figure}

\begin{definition} Let $f,g:(X,A)\rightarrow (Y,B)$ be coarse maps. A \textit{coarse homotopy of pairs} between $f$ and $g$ is a coarse map of pairs $H: (I_p(X), I_{pk_{A}}(A))\rightarrow (Y,B)$ so that $H \circ i_{0} = f$ and $H\circ i_{1} = g$. A \textit{relative coarse homotopy} between $f$ and $g$ is a coarse homotopy of pairs between $f$ and $g$ such that $H(a,t) = f(a)$ for all $a\in A$ and $t\leq pk_{A}(a)$.
\end{definition}

\begin{definition} Let $X$ be a bornological coarse space with base ray $\omega:[1,\infty)\rightarrow X$. For $n\geq 1$ define the \textit{$n$-th coarse homotopy group} $\pi_n^c(X,\omega)$ to be the set of relative coarse homotopy classes of maps 
\begin{align*}
f: (c([0,1]^n), c(\partial [0,1]^n))\rightarrow (X,\omega[1,\infty))
\end{align*}
such that $f|_{c(\partial[0,1]^n)} = \omega \circ p$, where $p:c([0,1]^n)\rightarrow [1,\infty), (hx,h)\mapsto h$ denotes the height variable of the cone. \\

For a bornological coarse pair $k_A: A\rightarrow X$ with base ray $\omega:[1,\infty)\rightarrow A$ we define the relative $n$-th coarse homotopy "group" $\pi_{n}^{c}(X,A, \omega)$ to be the set of relative coarse homotopy classes of maps 
\begin{align*}
f: (c([0,1]^n), c(\partial [0,1]^n), c(\partial_{+}[0,1]^n))\rightarrow (X,A,\omega [1,\infty))
\end{align*}
such that $f|_{c(\partial_{+}[0,1]^n)}= \omega \circ p$ where $\partial_{+}[0,1]^n := \{(x_1,\dots,x_n)\in \partial [0,1]^n\,|\,x_n>0\}$. 
\end{definition}

The notion of relative here means relative to $c(\partial_{+}[0,1]^n)$. The homotopy restricted to $c(\partial [0,1]^n \times [0,1])$ has image in $A$ and is fixed on $c(\partial_{+}[0,1]^n \times [0,1])$ but does not need to be fixed on all of $c(\partial [0,1]^n \times [0,1])$. We frequently refer to such homotopies as "relative and base ray preserving" and call two maps $f,g$ "base ray preserving relatively coarsely homotopic" if there is a relative and base ray presreving coarse homotopy between them. \\

The reason why can fix $p$ to be the standard height projection is because of the following:

\begin{lemma} \label{p_0} Let $X$ be a compact metric space, $cX$ the Euclidean cone of $X$, thought of as a bornological coarse space. Let $p_0:cX\rightarrow [1,\infty)$ be the height projection. Let $q: cX\rightarrow [1,\infty)$ be any controlled map. Then any coarse homotopy $H:I_{q}(cX)\rightarrow Y$ between $f:cX\rightarrow Y$ and $g: cX\rightarrow Y$ gives rise to a coarse homotopy $\overline{H}: I_{p_0}(cX)\rightarrow Y$ between $f$ and $g$. 
\end{lemma}
\begin{proof} Let $(h\mathbf{x},h)\in cX$. We have that $d_{cX}((h\mathbf{x},h),((h+1)\mathbf{x},h+1)) = \|\mathbf{x}\|+1\leq \frac{5}{4}$. Let $L>0$ be such that $|q(x)-q(y)|<L$ for all $x,y\in cX$ with $d_{cX}(x,y)<\frac{5}{4}$. Then
\begin{align*}
|q(h\mathbf{x},h)-q(\mathbf{x},1)| \leq |q(h\mathbf{x},h)-q(\left \lfloor{h}\right\rfloor\mathbf{x}, \left \lfloor{h}\right\rfloor )| + \sum_{j=2}^{\left \lfloor{h}\right \rfloor } |q(j\mathbf{x},j)-q((j-1)\mathbf{x},j-1)|<Lh 
\end{align*}
Let $C = \sup_{x\in X} q(\mathbf{x},1)$. 
\begin{align*}
|q(h\mathbf{x},h)| \leq |q(h\mathbf{x},h)-q(\mathbf{x},1)| + |q(\mathbf{x},1)| < Lh + C
\end{align*}
Let $q' = L p_0 + C$. We can extend the homotopy $H$ to $H': I_{q'}(cX)\rightarrow Y$ by extending "constantly" for the additional time, ie. $H'((h\mathbf{x},h),t) = g(h\mathbf{x},h)$ if $((h\mathbf{x},h),t) \in I_{q'}(X)\setminus I_{q}(X)$. Finally, there is a coarse equivalence $\Psi: I_{p_0}(cX)\rightarrow I_{q'}(cX)$ with 

\[ \Psi(x,t) = \begin{cases} 
          (x,(C+L)t) & t\in [0,1] \\
          (x,C+Lt)& t\in [1,p_{0}(x)] \\ 
       \end{cases}
    \]

For coarseness, there is nothing to check except in the time coordinate, but $\Psi$ is piecewise linear there, so it is controlled. Properness and the existence of an inverse are also easy to check.
\end{proof} 

\begin{remark} An analogous argument can be made for the mapping telescope $M_{X}$ defined later on. \end{remark}

The advantage of  Lemma \ref{p_0} is we can always assume that the map $p$ defining a coarse cylinder is the height projection. We can canonically identify
\begin{align*}
I_{p}(c[0,1]^n) &= c[0,1]^{n+1}\\
(hx_1,\dots, hx_n,h,ht)&\mapsto (hx_1,\dots, hx_n,ht,h)
\end{align*} 
They are isometric if we choose the standard $L^1$-norm on $\mathbb{R}^{N}$ and the product metric on $I_{p}(c[0,1]^n)$, and bi-Lipschitz equivalent for any other sensible choice. Analogously, $I_{p}(cX) = c(X \times [0,1])$. Note that because of this identification, we will occasionally swap the last two variables throughout the rest of the thesis, writing parameters in the order most convenient for exposition, eg. putting the time variable last when we want to emphasise that a map is a homotopy. \\

The product of two elements is defined in the standard way. Let $n\geq 1$ and $f,g$ be representatives of $[f],[g]\in \pi_n^{c}(X,\omega)$. The element $[f]\cdot [g]$ is represented by 
\[ f \cdot g(hx_1,hx_2,\dots, hx_n, h) = \begin{cases} 
          f(2hx_1,hx_2,\dots, hx_n, h) &  x_1\in [0, \frac{1}{2}]\\
         g(h(2x_1-1),hx_2,\dots, hx_n, h) &  x_1 \in [\frac{1}{2},1]\\ 
       \end{cases}
    \]

This operation turns $\pi_n^{c}(X,\omega)$ into a group, which is abelian if $n\geq 2$. For $n\geq 2$, the same formula defines a group structure for the relative homotopy groups, which are abelian if $n\geq 3$. The unit is represented by the map $\omega \circ p: c[0,1]^n \rightarrow X$. \\

Coarse homotopy groups are functorial under coarse maps: 

\begin{prop} \label{functorial} Let $(X,A)$ and $(Y,B)$ be bornological coarse pairs with base rays $\omega,\tau$ and let $\mathfrak{a}: (X,A)\rightarrow (Y,B)$ be a base ray preserving coarse map. Then there is an induced homomorphism
\begin{align*}
\mathfrak{a}_{*}: \pi_n^c(X,A,\omega)\rightarrow \pi_{n}^c(Y,B,\tau)
\end{align*}
defined by the formula $\mathfrak{a}_*[f]=[\mathfrak{a} \circ f]$. Furthermore, if $\mathfrak{a}, \mathfrak{a}': (X,A,\omega)\rightarrow (Y,B,\tau)$ are base ray preserving relatively coarsely homotopic, then the homomorphisms $\mathfrak{a}_*$ and $\mathfrak{a}'_*$ are equal.
\end{prop}

\begin{proof} The map $\mathfrak{a} f$ is coarse. If $f$ is coarsely homotopic to $g$, then $\mathfrak{a}f$ is coarsely homotopic to $\mathfrak{a}g$. By Proposition \ref{properties}, if $\mathfrak{a}$ is coarsely homotopic to $\mathfrak{a'}$ then $\mathfrak{a}f$ is coarsely homotopic to $\mathfrak{a}'f$, by a relative and base ray preserving homotopy. We have used here that the bornology on $c[0,1]^n$ is induced by the metric, so the controlled map $f$ is automatically bornologous. 
\end{proof}

An analogous statement is true for pointed sets. If $\mathfrak{a}:(X,\omega)\rightarrow (Y,\tau)$ is a base ray preserving coarse map between bornological coarse spaces with base rays, there is an induced map $\mathfrak{a}_{*}:\pi^c_0(X,\omega)\rightarrow \pi^c_0(Y,\tau)$ by sending $[\beta]\in \pi^c_0(X,\omega)$ to $[\mathfrak{a}\beta]$. Base ray preserving homotopic maps induce the same morphism. 

\begin{prop} \label{longexact} If $(X,A,\omega)$ is a bornological coarse space with map $k: A\rightarrow X$, and base ray $\omega:[1,\infty)\rightarrow A$, the analogue of the usual construction in topology defines a long exact sequence of coarse homotopy groups or pointed sets
\begin{align*}
\dots \rightarrow \pi_2^c(A,\omega)\xrightarrow{k_*} \pi_2^c(X,\omega)\xrightarrow{r} \pi_2^c(X,A,\omega) \xrightarrow{\partial} \pi_1^c(A,\omega)\xrightarrow{k_*} \pi_1^c(X,\omega)\xrightarrow{r} \pi_1^c(X,A,\omega)\\
\xrightarrow{\partial}  \pi_0^c(A,\omega)\xrightarrow{k_*} \pi_0^c(X,\omega)
\end{align*}
\end{prop}
We only prove the proposition for $A\subset X$ with the restricted bornological coarse structure. (The general statement is left as an exercise.) Consider $[0,1]^n$ as a metric space equipped with a (any) norm-induced metric. Equip $[0,1]^n \times I$ with the product metric. 

\begin{lemma} Let $G: [0,1]^n \times I \rightarrow [0,1]^n$ be a Lipschitz homotopy with Lipschitz constant $C$. Define the cone $cG$ as 
\begin{align*}
cG: I_{p} (c[0,1]^n)= c[0,1]^{n+1} &\rightarrow c[0,1]^n
\\ (hx,h,ht) &\mapsto (h G(x,t),h)
\end{align*}
This is a proper Lipschitz homotopy between $cG_0$ and $cG_1$. 
\end{lemma}

\begin{proof} A bounded set $B$ in $c[0,1]^n$ has a maximum height $h'$, and by construction a point at height $h$ stays at height $h$, so the preimage of $B$ is contained in the set $c_{[1,h']} [0,1]^{n+1}$, which is bounded. To prove that $cG$ is controlled, suppose we have two points $(hx,h,ht)$ and $(sy,s,sl)$ with distance less than $R$. Obviously then $\|hx-sy\|<R$, $|h-s|<R$, $|ht-sl|<R$. We compute:
\begin{align*}
d(cG(hx,h,ht), cG(sy,s,sl))= |h-s| + \|hG(x,t)-s G(y,l)\| \\
<R + \|h G(x,t)-h G(y,l)\|+ \|hG(y,l)-sG(y,l)\| \\
<R +Ch(\|x-y\|+ |t-l|) + RK_n \\
\leq R+ RK_{n}+ C(\|hx-sy\|+ \|sy-hy\|+ |ht-sl|+ |sl-hl|) \\
< R+ RK_n + C(R+RK_n+R+R)
< \tilde{C} R
\end{align*}
where $K_n$ is a constant depending on the choice of metric on $[0,1]^n$ and $\tilde{C}$ is some other constant after absorbing everything, dependent only on $K_n$ and $C$. 
\end{proof}

Note that this proof works just as well for arbitrary compact metric spaces $X,Y$: a Lipschitz homotopy $G: X\times [0,1]\rightarrow Y$ defines a proper Lipschitz homotopy $cG: c(X \times [0,1])\rightarrow cY$.

\begin{lemma} (Compression criterion) \label{compression} Suppose that
\begin{align*}
f: (c([0,1]^n), c(\partial [0,1]^n), c(\partial_{+}[0,1]^n))\rightarrow (X,A,\omega[1,\infty))
\end{align*}
is a coarse map. Then $[f] = 0 \in \pi_n^c(X,A,\omega)$ if and only if it is coarsely homotopic relative to $c(\partial [0,1]^n)$ to a map with image contained in $A$.
\end{lemma}

\begin{proof} First, suppose that $[f]$ is homotopic to a map $[g]$ as in the lemma. Then $[f]=[g]\in \pi_n^c(X,A,\omega)$. Now, choose a Lipschitz homotopy $H: [0,1]^n \times [0,1]\rightarrow [0,1]^n$ from the identity to a constant map with value $s = (x_1,\dots,x_{n-1},1)$ with the property that $H(\partial_{+}[0,1]^n \times [0,1])\subset \partial_{+}[0,1]^n$ (for example, deformation retract everything to level $1$ and then deformation retract to $s$). Then $g \circ cH$ is a coarse homotopy with image in $A$ between $g$ and the constant map $\omega \circ p$ to the base ray . By construction $g \circ cH$ is $\omega \circ p$ when restricted to $c(\partial_{+}[0,1]^n \times [0,1])$. Therefore, it is a homotopy of triples $(c([0,1]^n), c(\partial [0,1]^n), c(\partial_{+}[0,1]^n))\rightarrow (X,A,\omega[1,\infty))$. Thus $[g]=0$. \\

Conversely, suppose that $[f] = 0 \in  \pi_n^c(X,A,\omega)$. This means that there is a homotopy
\begin{align*}
F: I_p(c[0,1]^n) \rightarrow X
\end{align*} 
of triples between $f$ and the constant map to the base ray. Observe that the "end" of the cylinder along with the cylinder restricted to $c(\partial[0,1]^n \times [0,1])$
\begin{align*}
\{(hx,h,h)\,|\, x\in [0,1]^n\} \cup \{ (hx,h) \times [0,h] \,|\,x\in \partial[0,1]^n\}
\end{align*}
is itself a copy of $c([0,1]^n)$ which is Lipschitz homotopic to $c([0,1]^n \times \{0\})$ in $I_p(c[0,1]^n)$, fixing $c(\partial[0,1]^n \times \{0\})$. This shows that $f$ is homotopic rel $c(\partial[0,1]^n)$ to a map with image in $F(\{(hx,h,h)\,|\, x\in [0,1]^n\} \cup \{ (hx,h) \times [0,h] \,|\,x\in \partial[0,1]^n\}) \subset A$. Concretely, we have a Lipschitz deformation retract $D: [0,1]^n \times [0,1] \times [0,1] \rightarrow [0,1]^n \times [0,1]$ relative to $\partial [0,1]^n \times \{0\}$ such that $D_0= \id$ and $D_{1}([0,1]^n \times [0,1])\subset ([0,1]^n \times \{1\})\cup(\partial [0,1]^n \times [0,1])$,  define $q: I_p(c[0,1]^n)\rightarrow [1,\infty)$ as $q(hx,h,s) = h$ and we have a homotopy of the homotopy cylinder 
\begin{align*}
H: c([0,1]^n \times [0,1] \times [0,1])=I_{q}(I_{p}(c[0,1]^n)) \rightarrow X 
\end{align*} 
which is defined as $F\circ cD$. The restriction of $H$ to $c([0,1]^n \times \{0\} \times [0,1])$ is the desired coarse homotopy from $f$ to a map with image in $A$.  
\end{proof}

Now we prove the exactness of the long exact sequence in coarse homotopy groups:

\begin{proof}  (of Proposition \ref{longexact}) The first three points apply to $n \geq 1$. 
\begin{enumerate}
\item $\image(\partial)=\Ker(k_*)$: Suppose we have $[f]\in \pi_n^c(X,A,\omega)$, represented by 
\begin{align*}
f:(c([0,1]^n), c(\partial [0,1]^n), c(\partial_{+}[0,1]^n))\rightarrow (X,A,\omega[1,\infty))
\end{align*}
and denote by $g$ the restriction of $f$ to $ c([0,1]^{n-1} \times \{0\})$, ie $\partial[f]=[g]$. Choose a Lipschitz deformation retract $H$ of $[0,1]^n$ to $s \in [0,1]^{n-1} \times \{1\}$. Then $f\circ cH$ restricted to $c([0,1]^{n-1} \times \{0\} \times [0,1])$ is a base ray preserving homotopy between $g$ and the map $\omega \circ p$. Hence $[g] = 0$ in $\pi_{n-1}^c(X,\omega)$. \\

Conversely, suppose that $k_*([g])=0$, which means that there is a coarse homotopy $H: I_p(c[0,1]^{n-1})\rightarrow X$ rel $c(\partial[0,1]^{n-1})$ of $g$ to the map $\omega \circ p$. Put $f(h(x_1,\dots,x_{n-1},x_n),h) = H(h(x_1,\dots,x_{n-1}),h,hx_n)$. This is a map $f$ with $[f]\in \pi_n^c(X,A,\omega)$ which restricts to $g$ on $c([0,1]^{n-1})$, showing that $[g]=\partial [f]$.
\item $\image(k_*) = \Ker{(r)}$: This is exactly the contents of Lemma \ref{compression}.
\item $\image{(r)} = \Ker(\partial)$: If $f$ represents a class in $\pi_n^c(X,\omega)$ then thinking about it as an element of $\pi_n^c(X,A,\omega)$ means immediately that its restriction to $c([0,1]^{n-1} \times \{0\})$ is the constant map $\omega \circ p$ to the base ray, showing that $\partial \circ r = 0$.\\

Conversely, suppose the restriction of $f$  to $c([0,1]^{n-1} \times \{0\})$ is homotopic rel $c(\partial [0,1]^{n-1} \times \{0\})$ with image in $A$ to the map $\omega \circ p$. Then we can "stack" this homotopy onto the bottom of $f$ to obtain a homotopy of triples between $f$ and map $g$, where $g$ restricted to $c(\partial[0,1]^n)$ is the constant map $\omega \circ p$ to the base ray. Therefore $[f]=[g] \in \image (r)$. 
\item Exactness at $\pi^{c}_0(A,\omega)$: Let $[\beta]\in \pi^c_{0}(A,\omega)$. If there is a coarse homotopy $H: I_{p}([1,\infty))=c(\{\ast\} \times [0,1])\rightarrow X$ between $\beta$ and $\omega$, then $[H]\in \pi^c_1(X,A,\omega)$ with $\partial[H]=[\beta]$. Conversely, for $[f]\in  \pi^c_1(X,A,\omega)$, we have that $f$ restricted to $c(\{0\})$ is a base ray in $A$ which is coarsely homotopic in $X$ to $\omega$. 
\end{enumerate}
\end{proof}

We show that coarse homotopy groups commute with taking some colimits. 

\begin{definition} A \textit{based bornological coarse space} is a quadruple  $(X,\mathcal{C},\mathcal{B},\omega)$ consisting of a bornological coarse space $(X,\mathcal{C},\mathcal{B})$ and a coarse map $\omega: [1,\infty)\rightarrow X$. 
$\textbf{BornCoarse}_{*}$ is the category of based bornological coarse spaces. A morphism $f: (X,\mathcal{C},\mathcal{B},\omega) \rightarrow  (Y,\mathcal{C}',\mathcal{B}',\tau)$ is a base ray preserving coarse map. 
\end{definition}

We have shown in Proposition \ref{functorial} (take $A=\emptyset$) that coarse homotopy groups (resp. sets) are functors: $\textbf{BornCoarse}_{*}\rightarrow \mathbf{C}$, where $\mathbf{C}$ is the category of pointed sets for $n=0$, groups for $n=1$, or abelian groups for $n\geq 2$. We now introduce the idea of $u$-continuity, which is used in \cite{bunke2020coarse} to define the axioms of a local homology theory.

\begin{definition}(Definition $3.11$ in \cite{bunke2020coarse}) Let $E:\mathbf{UBC}\rightarrow \mathbf{C}$ be a functor. We say that $E$ is $u$-continuous if for every uniform bornological space $X$ the canonical morphism
\begin{align*} 
\colim_{V\in \mathcal{C}} E(X_{V})\rightarrow E(X)
\end{align*} is an equivalence. 
\end{definition}

We use a variant of this definition.\\

Let $(X,\mathcal{C},\mathcal{B},\omega)$ be a bornological coarse space with base ray. Since $\omega: [1,\infty) \rightarrow X$ is coarse, we have that $(\omega \times \omega)D_{1}$ is an entourage. Let $\mathcal{C}_{\omega}\subset \mathcal{C}$ be the subset of $\mathcal{C}$ consisting of entourages that contain $(\omega \times \omega)D_{1}$. Denote by $(X_{V}, \mathcal{C}\langle V\rangle, \mathcal{B},\omega)$ the based coarse bornological space with coarse structure generated by the entourage $V\in \mathcal{C}_{\omega}$. If $V\subset W$ then the set-wise identity $i_{V,W}:X_{V}\rightarrow X_{W}$ is coarse and base ray preserving. We have the equivalence 
\begin{align*}
\colim_{V \in \mathcal{C}_{\omega}} (X_{V},\omega) \simeq (X,\omega)
\end{align*}
in $\textbf{BornCoarse}_{*}$.

\begin{definition} Let $E: \textbf{BornCoarse}_{*}\rightarrow \mathbf{C}$ be a functor. We say that $E$ is \textit{$u$-continous} if for every based bornological coarse space $(X,\omega)$, the canonical morphism 
\begin{align*}
\colim_{V \in \mathcal{C}_{\omega}} i_{V}: \colim_{V \in \mathcal{C}_{\omega}} E(X_{V},\omega)\rightarrow E(X,\omega)
\end{align*}
induced by the set-wise identity maps $i_{V}:  (X_{V},\omega) \rightarrow (X,\omega)$ is an equivalence.
\end{definition}

\begin{prop} The functor $\pi^{c}_{n}:  \textbf{BornCoarse}_{*} \rightarrow \mathbf{C}$ where $\mathbf{C}$ is the category of pointed sets for $n=0$, groups for $n=1$, or abelian groups for $n\geq 2$, is $u$-continuous. 
\end{prop}

\begin{proof} We define the inverse to $\colim_{V\in \mathcal{C}_{\omega}}i_{V}$. Let $[f]\in \pi_{n}^{c}(X,\omega)$. For $n\geq 1$, $[f]$ is represented by a coarse map of pairs $f:(c[0,1]^n,c\partial[0,1]^n)\rightarrow (X,\omega)$. We have that $(f \times f)D_{1}\subset U$ for some $U\in \mathcal{C}$. Since $(\omega \times \omega) D_{1}\subset (f\times f)D_{1}$ (take points in $c\{(0,0)\}$) we have that $U\in \mathcal{C}_{\omega}$. Denote by $[f_{U}]\in \pi_{n}^{c}(X_{U},\omega)$ the homotopy class represented by $f$. We send $[f]$ to the image of $[f_{U}]$ in the colimit. 
\begin{align*}
j: \pi_{n}^{c}(X,\omega)&\rightarrow \colim_{V\in \mathcal{C}_{\omega}} \pi_{n}^{c} (X_{V},\omega)\\
[f]&\mapsto [f_{U}]
\end{align*}
This is well defined: if $U\subset W$ then $i_{U,W}[f_{U}] = [f_{W}]$ because all maps are induced by the identity. If $[f]=[g] \in \pi_{n}^{c}(X,\omega)$ then there exists a coarse homotopy $H: (c[0,1]^{n+1}, c(\partial[0,1]^n \times [0,1])) \rightarrow (X,\omega)$ such that $H_{|c([0,1]^{n}\times \{0\})} =f, H_{|c([0,1]^{n}\times \{1\})} = g$.  There exists a $W\in \mathcal{C}_{\omega}$ such that $(H \times H)D_{1}\subset W$. This gives us 
\begin{align*}
j[f] = [f_{W}] = [g_{W}] = j[g]
\end{align*}
Both compositions $j\circ \colim_{V\in \mathcal{C}_{\omega}} i_{V}, \colim_{V\in \mathcal{C}_{\omega}} i_{V} \circ j$ are the identity on representatives. \\

An analogous statement holds for $n=0$.
\end{proof}

\subsection{Base ray considerations}

We introduce in this section the change of base ray homomorphism $b$, and show that a coarse homotopy equivalence induces an isomorphism on coarse homotopy groups. Note that throughout this section, we will use $\|\cdot\|$ to denote the $L^1$-norm $\|\cdot \|_{1}$ on $[0,1]^n\subset \mathbb{R}^{n}$. 

 \begin{definition} \label{DefBaserayChange} Identify the unit cube $[0,1]^n$ with the set $[-1,1]^n\subset \mathbb{R}^n$ in the obvious way and take the Euclidean cone $c[-1,1]^n\subset \mathbb{R}^{n+1}$. Let $f: (c[-1,1]^n, c\partial[-1,1]^n)\rightarrow (Z,\omega_0)$ be a coarse map representing a class $[f]\in \pi_n^c(Z,\omega_0)$. Let $\omega: [1,\infty)\rightarrow Z$ be another base ray in $Z$ such that $\omega_0$ is coarsely homotopic to $\omega$: ie. there exists a coarse map $H: I_p[1,\infty)\rightarrow Z$ restricting to $\omega_0$ and $\omega$ on $c\{0\}$ and $c\{1\}$ respectively. We define 
 \begin{align*}
 b_{\omega_0,\omega}^H: \pi_n^c(Z, \omega_0) \rightarrow \pi_n^c(Z,\omega)
 \end{align*} as follows:\\

  Let $U = cW := \{(hx_1,\dots,hx_n,h) \subset c[-1,1]^n \,|\, -\frac{1}{2}\leq x_i \leq \frac{1}{2} \, \forall 1\leq i\leq n\}$
be a smaller cone inside $c[-1,1]^n$. Let $\mathfrak{r}: U\rightarrow c[-1,1]^n$ be the height-wise radial multiplication by $2$, $(hx_1,\dots,hx_n,h)\mapsto (2hx_1,\dots,2hx_n,h)$. Define $b^H_{\omega_0,\omega}[f]|_{U}$ as the composition $f \circ \mathfrak{r}$. On $c[-1,1]^n\setminus U$ we have 
\begin{align*}
b^H_{\omega_0,\omega}[f]|_{c[-1,1]^n\setminus U}: c[-1,1]^n\setminus U &\rightarrow Z\\
(h\mathbf{x},h) &\mapsto H(h, h(2\|\mathbf{x}\|_{\infty} - 1))
\end{align*}
That is, on height $h$  we travel radially outwards along the path between $\omega_0(h)$ and $\omega(h)$ as defined by the coarse homotopy $H$. See Figure \ref{fig:defb}. \\

\begin{figure}[H]
\centering
  \centering
  \includegraphics[width=0.4\linewidth]{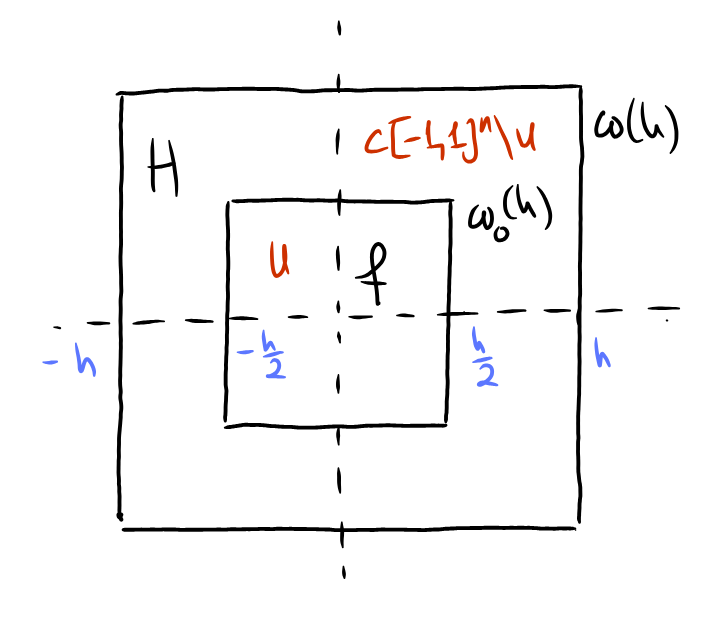}
  \caption{A representative of $b^H_{\omega_0,\omega}[f]$ restricted to height $h$.}
  \label{fig:defb}
\end{figure} 
 
We call $b^H_{\omega_0,\omega}$ the \textit{change of base ray homomorphism} between $\omega_{0}$ and $\omega$. We omit the indices $\omega_0,\omega,H$ from the notation when it is obvious from context. 

 \end{definition}
  
 \begin{definition} Let $H$ be a coarse homotopy between $\omega_0$ and $\omega_1$, $H'$ a coarse homotopy between $\omega_1$ and $\omega_2$. We define the concatenation of homotopies $H\ast H' :c[0,1] \rightarrow Z$ as 

\[ H\ast H'(h,ht) = \begin{cases} 
          H(h, 2ht) & t \in [0,\frac{1}{2}] \\
          H'(h, h(2t-1)) & t\in [\frac{1}{2},1] \\ 
       \end{cases}
    \]
 \end{definition}
This is just a rescaled version of running first $H$, then $H'$. It is a coarse homotopy between the rays $\omega_0$ and $\omega_2$. An analogous formula holds for homotopies $H,H': c([0,1]^n \times [0,1])\rightarrow Z$, $n\geq 1$.

\begin{prop}\label{baseddefined} $b^H_{\omega_0,\omega}$ is well defined: 
\begin{enumerate}
\item $b^H_{\omega_0,\omega}[f]\in \pi_{n}^{c}(Z,\omega)$
\item $b^H_{\omega_0,\omega}[f]$ is independent of choice of representative $f$
\item $b^H_{\omega_0,\omega}$ is a group homomorphism 
\item $b^H_{\omega_0,\omega}$ is independent of homotopy class of $H$
\item $b_{\omega_1,\omega_2}^{H'}b_{\omega_0,\omega_1}^{H}= b_{\omega_0,\omega_2}^{H\ast H'}$
\item $b_{\omega_0,\omega_0}^{\id} = \id$
\item Let $\mathfrak{a}: X\rightarrow Y$ be a coarse map between bornological coarse spaces. Then the diagram commutes

\begin{figure}[H]
\center
\begin{tikzcd}
&\pi_n^{c}(X,\omega_0)\arrow {r}{\mathfrak{a}_{*}} \arrow{d}{b_{\omega_0,\omega}^H}
&\pi_n^{c}(Y,\mathfrak{a}\omega_0) \arrow{d}{b_{\mathfrak{a}\omega_0,\mathfrak{a}\omega}^{\mathfrak{a}H}}
  \\
&\pi_{n}^c(X,\omega) \arrow {r}{\mathfrak{a}_{*}} 
&\pi_{n}^c(Y,\mathfrak{a}\omega)

\end{tikzcd}
\end{figure}

\end{enumerate}
\end{prop}

\begin{remark} Piecing statements $4,5,6$ together we get that the change of base ray homomorphism is an isomorphism: $b_{\omega,\omega_0}^{H^{-1}}b_{\omega_0,\omega}^{H}= b_{\omega_0,\omega_0}^{H\ast H^{-1}}=b_{\omega_0,\omega_0}^{\id}= \id$ and $b_{\omega_0,\omega}^{H}b_{\omega,\omega_0}^{H^{-1}} = b_{\omega,\omega}^{H^{-1}\ast H} = b_{\omega,\omega}^{\id}=  \id$. 
\end{remark}

\begin{proof} 
\begin{enumerate}
\item $\mathfrak{r}$ is clearly proper and has Lipschitz constant $2$. $b[f]|_{U}$ is coarse, as a composition of proper Lipschitz and coarse maps. To show that $b[f]|_{c[-1,1]^n\setminus U}$ is coarse, we just compute: 
 \begin{align*}
 d_{c[-1,1]^n}(h,h(2\|\mathbf{x}\|_{\infty}-1))(s, s(2\|\mathbf{y}\|_{\infty}-1)) \\
 \leq 2|\|h\mathbf{x}\|_{\infty}-\|s\mathbf{y}\|_{\infty}|+2|h-s|\\
 \leq 2\|h\mathbf{x}-s\mathbf{y}\|_{\infty}+2|h-s| \\
 \leq 2 d_{c[-1,1]^n}((h\mathbf{x}, h,ht),(s\mathbf{y} ,s,sl))
 \end{align*} 
 $b[f]|_{c[-1,1]^n\setminus U}$ is a composition of proper Lipschitz and coarse maps, therefore coarse. \\
 
 Clearly $b[f]$ sends $c\partial[-1,1]^n$ to the base ray $\omega$. Properness follows because the pre-image of a bounded set $K$ is $b[f]^{-1}(K)=  b[f]|_{U}^{-1}(K) \cup b[f]_{c[-1,1]^n\setminus U}^{-1}{K} = \mathfrak{s}(f^{-1}{K})\bigcup H^{-1}(K)\times \partial[-1,1]^n$, where $\mathfrak{s}: c[-1,1]^n \rightarrow U$ is the height-wise radial multiplication by $\frac{1}{2}$ and $H^{-1}(K)\times \partial[-1,1]^n$ refers to the points $(hx_1,\dots,hx_n,h)$ in $c[-1,1]^n\setminus U \cong c([0,1] \times \partial [-1,1]^n)$ which have appropriate norm to be in the preimage of $K$. Since taking products and forming unions of bounded sets is bounded, $b[f]^{-1}{(K)}$ is bounded. Since the maps agree on their intersection and $c[-1,1]^n$ is a path metric space, $b[f]$ is globally coarse and therefore $b[f]\in \pi_{n}^{c}(Z,\omega)$.\\
\item Let $g$ be coarsely homotopic to $f$ via $G':I_{p}(c[-1,1]^n) = c([-1,1]^n \times [0,1])\rightarrow Z$ which preserves the base ray. We can run $G'$ on the cylinder $c(W \times [0,1])$ of $U$, and "fill in" with the constant homotopy which changes base ray on the remainder. Concretely, on $(U\times [0,\infty)) \cap I_{p}(c[-1,1]^n)$ we define $G$ as the composition $G' \circ (\mathfrak{r} \times \id)$ and on $(c[-1,1]^n \setminus U) \times [0,\infty)$ we project onto $c[-1,1]^n \setminus U$ and apply $b[f]|_{c[-1,1]^n\setminus U} =b[g]|_{c[-1,1]^n\setminus U}$. See Figure \ref{fig:bpic7}. 

\begin{figure}[H]
\centering
  \centering
  \includegraphics[width=0.5\linewidth]{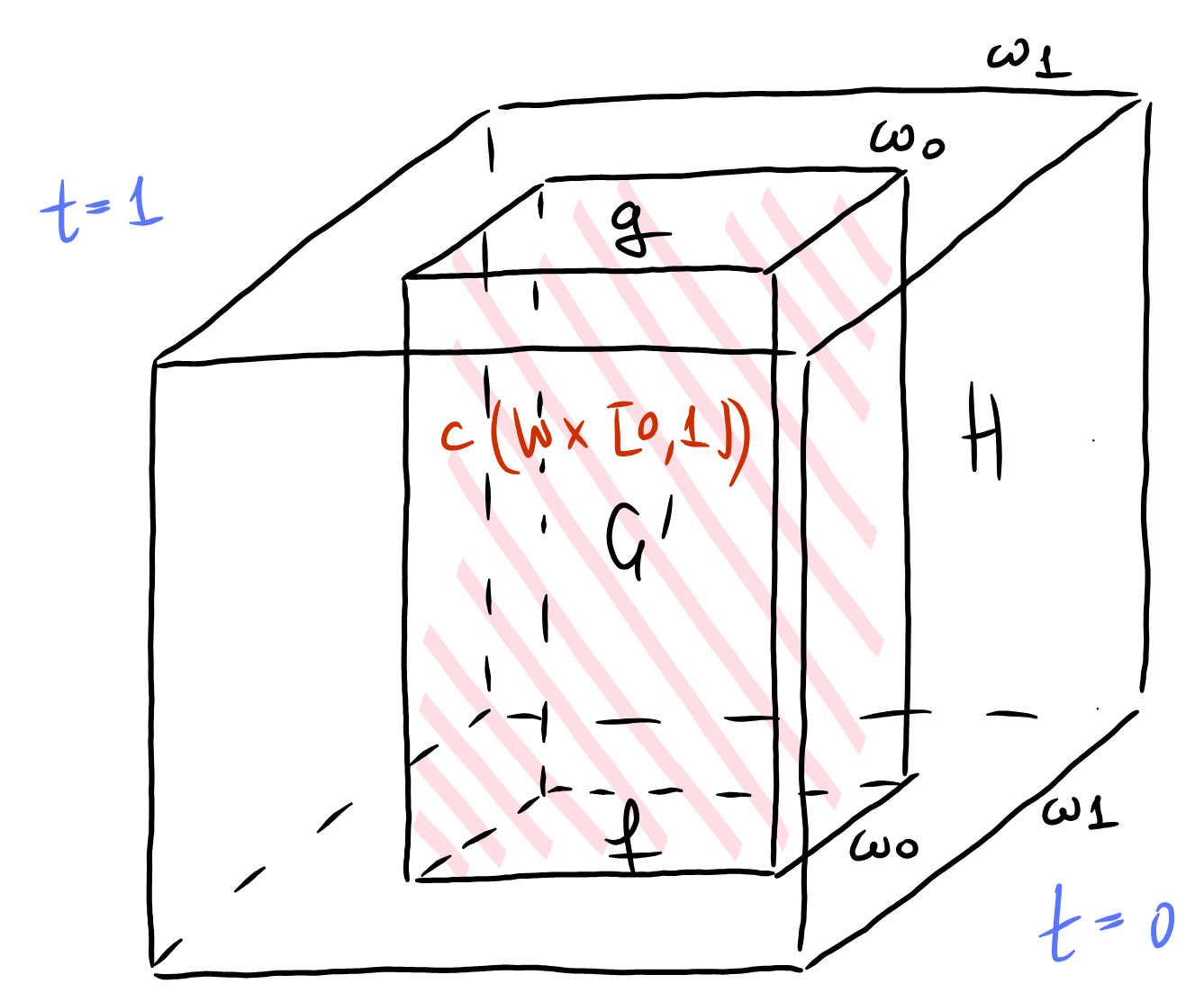}
  \caption{The coarse homotopy $G$ between $b[f]$ and $b[g]$ on height $h$.}
  \label{fig:bpic7}
\end{figure}

 Again, the restriction to each component is clearly controlled and the maps agree on the intersection and so is globally controlled. To show properness, observe that for a bounded set $K$ the preimage is $G^{-1}(K) = (\mathfrak{s}\times \id)(G'^{-1}(K)) \cup ((H^{-1}(K)\times \partial[-1,1]^n)\times [0,\infty) \cap I_p(c[-1,1]^n))$. This set is bounded since $((H^{-1}(K)\times \partial[-1,1]^n)\times [0,\infty) \cap I_p(c[-1,1]^n))$ is contained within $(H^{-1}(K)\times \partial[-1,1]^n)\times [0,h]$ for some $h>0$. Hence $G$ is proper.

\item The proof is basically taking the cone of the proof for ordinary homotopy groups. First we show that for $[f]\in \pi_n^{c}(Z,\omega_0)$, replacing half the cone with the constant map to the base ray gives us $f'$ which is coarsely homotopic to $f$. Let $V$ be $c([-1,0]\times [-1,1]^{n-1})\subset c[-1,1]^n$ and define:
\begin{align*}
f'_{|V}(hx_1,\dots,hx_n,h) &= f(2h (x_1+\frac{1}{2}),hx_2,\dots,hx_n,h)\\
f'_{|V^c}(h\mathbf{x},h) &= \omega_0(h)
\end{align*} 
 Let $Y:=\{(h\mathbf{x},h,ht)\in I_{p}(c[-1,1]^n)\,|\,x_1\leq t\}$. For the homotopy $J$ between $f'$ and $f$, we define
\begin{align*}
J_{|Y}(h\mathbf{x},h,ht) &= f(\frac{2h}{1+t}(x_1-\frac{t-1}{2}),hx_2,\dots,hx_n,h)\\
J_{|Y^{c}}(h\mathbf{x},h,ht)&= \omega_0(h)
\end{align*}
See Figure \ref{fig:bpic2}.

\begin{figure}[H]
\centering
  \centering
  \includegraphics[width=0.5\linewidth]{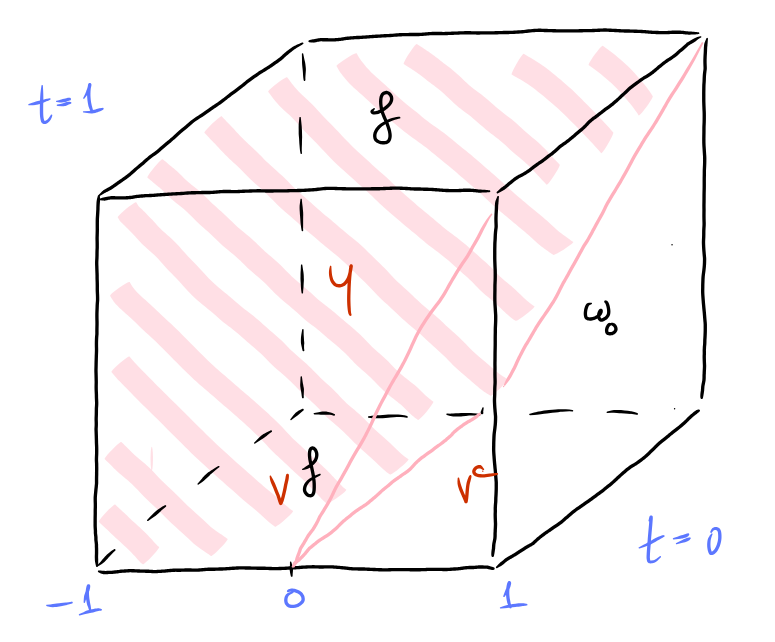}
  \caption{The coarse homotopy $J$ between $f'$ and $f$, restricted to height $h$.}
  \label{fig:bpic2}
\end{figure}

To compute the Lipschitz constants, let $(h\mathbf{x},h,ht)$ and $(s\mathbf{y},s,sl)$ be two points in $Y$. We have the following estimates:
\begin{align*}
&|\frac{2h}{1+t}(x_1-\frac{t-1}{2}) - \frac{2s}{1+l}(y_1-\frac{l-1}{2})|\\
&\leq |\frac{2hx_1}{1+t} - \frac{2sy_1}{1+l}| + |\frac{h (t-1)}{1+t} - \frac{s(l-1)}{1+l}|\\
&\leq |\frac{2hx_1}{1+t} - \frac{2sy_1}{1+t}| + |\frac{2sy_1}{1+t} - \frac{2sy_1}{1+l}|
+ |\frac{h(t-1)}{1+t}- \frac{h (l-1)}{1+l}| + |\frac{h(l-1)}{1+l} - \frac{s(l-1)}{1+l}|\\
&\leq 2|hx_1-sy_1| + 2s|t-l| + 2(|th-ls|+|h-s|) + |\frac{l-1}{1+l}||h-s|\\
&\leq 2|hx_1-sy_1| + 4(|th-ls|+|h-s|) + |h-s| 
\end{align*}
Therefore
\begin{align*}
&|\frac{2h}{1+t}(x_1-\frac{t-1}{2}) - \frac{2s}{1+l}(y_1-\frac{l-1}{2})| + \sum_{i=2}^n|hx_i-sy_i|+|h-s|\\
&\leq 2|hx_1-sy_1| + 4(|th-ls|+|h-s|) + 2|h-s| + \sum_{i=2}^n|hx_i-sy_i|\\
&\leq 2 \|h\mathbf{x}-s\mathbf{y}\| + 4|th-ls|+6|h-s|
\end{align*}
This shows that $H$ is a composition of controlled and Lipschitz maps, therefore controlled. Properness is clear.\\

Let $[f],[g]\in \pi_n^c(Z,\omega_0)$. We denote by $f+0$ the map $f'$ as defined above and by $0+g$ the analogous construction with $V'= c([0,1]\times [-1,1]^{n-1})\subset c[-1,1]^n$. We now construct a coarse homotopy $G$ between $b[f+0]+b[0+g]$ and $b[f+g]$. See Figure \ref{fig:bpic3}.

\[ G(h\mathbf{x},h,ht) = \begin{cases} 
          b(f+0) (\frac{2h}{1+t}(x_1-\frac{t-1}{2}),hx_2,\dots,hx_n,h) & x_1\in [-1,0] \\
          b(0+g) (\frac{2h}{1+t}(x_1+\frac{t-1}{2}),hx_2,\dots,hx_n,h) & x_1\in [0,1] \\ 
       \end{cases}
    \]

\begin{figure}
\centering
  \centering
  \includegraphics[width=0.8\linewidth]{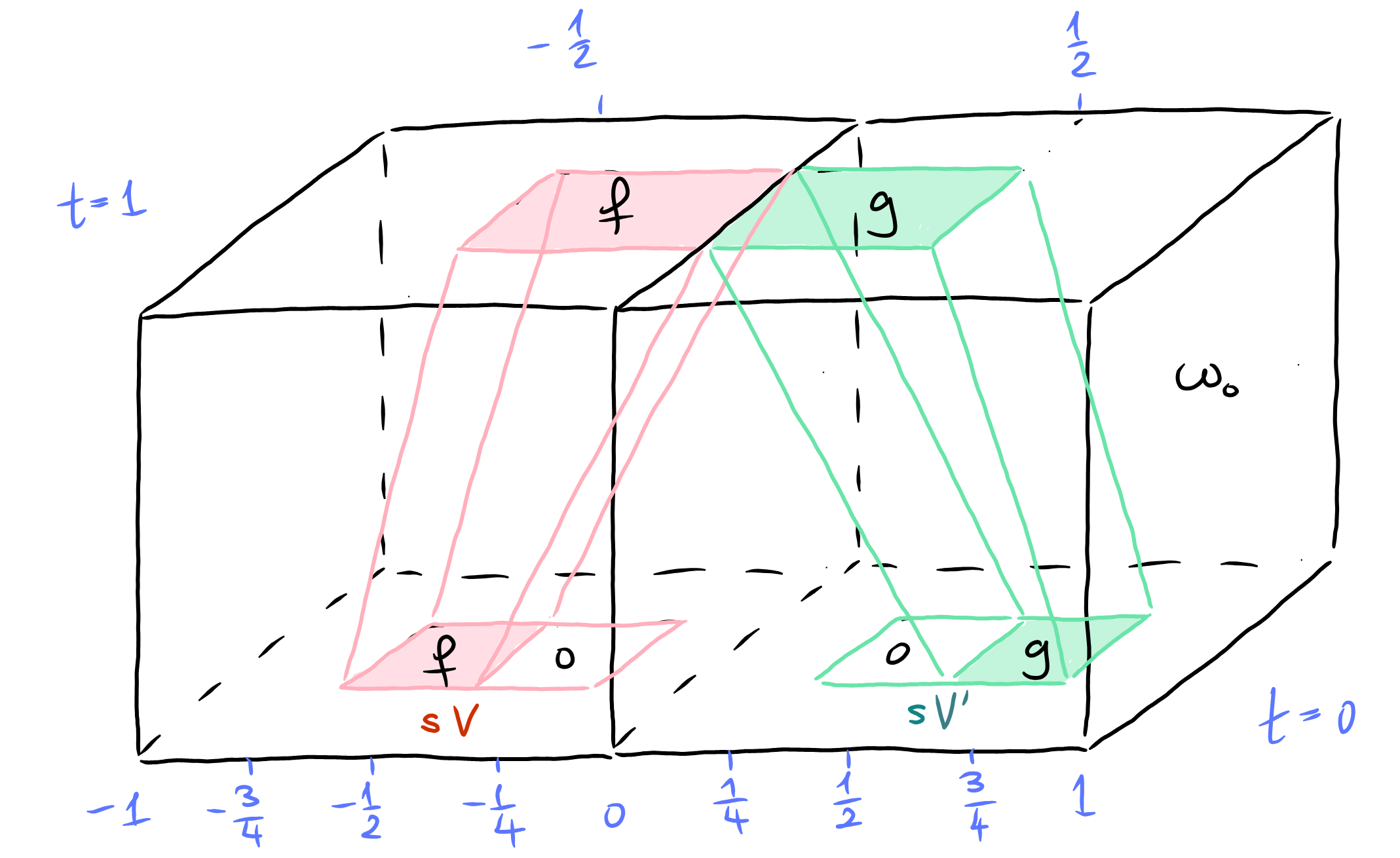}
  \caption{The coarse homotopy $G$ between $b[f+0]+b[0+g]$ and $b[f+g]$, restricted to height $h$.}
  \label{fig:bpic3}
\end{figure}

The calculation of the Lipschitz constants is exactly the same as above, and properness is also clear. The composition $(b[(J_{f})^{-1}]+b[(J_{g})^{-1}])\ast G$ (abusing notation) is a coarse homotopy between $b[f]+b[g]$ and $b[f+g]$.
\item Let $H'$ be homotopic to $H$ via a coarse homotopy $G': I_p(c[0,1]) = c([0,1]^2) \rightarrow Z$ restricting to $H$ on $c([0,1] \times \{0\})$ and $H'$ on $c([0,1] \times \{1\})$ respectively. On the cylinder of $(U\times [0,\infty)) \cap I_{p}(c[-1,1]^n)$ we define $G$ by projecting onto $U$ and applying the composition $f\circ \mathfrak{r}$. This is controlled, since projection has Lipschitz constant $1$. On its complement we define
\begin{align*}
G: ((c[-1,1]^n \setminus U) \times [0,\infty)) \cap I_{p}(c[-1,1]^n) &\rightarrow Z
\\(h\mathbf{x}, h,ht)&\mapsto G'(h, h(2\|\mathbf{x}\|_{\infty} - 1),ht)
\end{align*}
See Figure  \ref{fig:bpic1}. Since $G'$ is controlled, so is $G$ restricted to $((c[-1,1]^n \setminus U) \times [0,\infty) )\cap I_{p}(c[-1,1]^n)$. For properness, for a bounded set $K$ the preimage is $ (\mathfrak{s}\times \id)(f^{-1}(K) \times [0,\infty) \cap I_{p}(c[-1,1]^n)) \cup (G'^{-1}(K) \times \partial[-1,1]^n)$ which is bounded. 

\begin{figure}[H]
\centering
  \centering
  \includegraphics[width=0.5\linewidth]{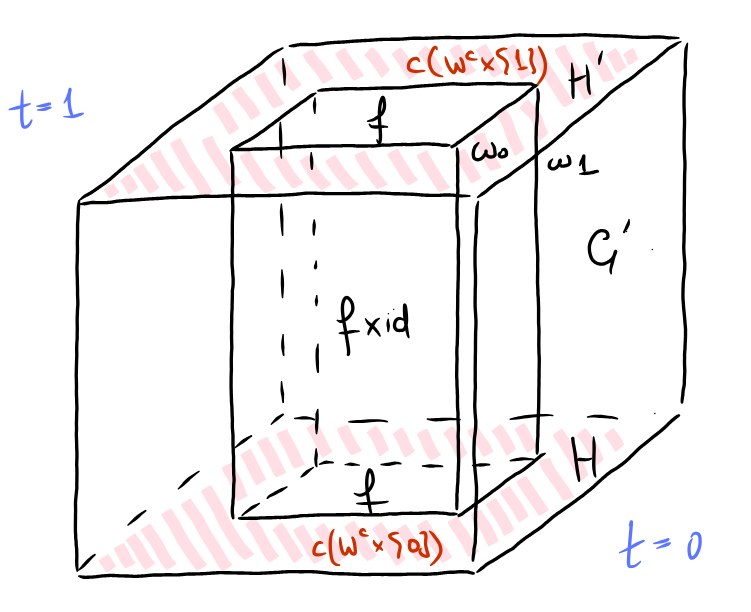}
  \caption{The coarse homotopy $G$ between $b^{H}_{\omega_{0},\omega}[f]$ and $b^{H'}_{\omega_{0},\omega}[g]$, restricted to height $h$.}
  \label{fig:bpic1}
\end{figure}

\item Let $[f]\in \pi_n^{c}(Z,\omega_0)$. For a fixed height $h$ on $c[-1,1]^n$ the description of  $b_{\omega_1,\omega_2}^{H'}b_{\omega_0,\omega_1}^{H}([f])$ is $f$ on points with supremum norm less than $\frac{h}{4}$, the homotopy $H$ between $\frac{h}{4}$ and $\frac{h}{2}$, and the homotopy $H'$ between $\frac{h}{2}$ and $h$. The description of $b_{\omega_0,\omega_2}^{H\ast H'}$ is $f$ on points with norm less than $\frac{h}{2}$, $H$ between $\frac{h}{2}$ and $\frac{3h}{4}$ and $H'$ between $\frac{3h}{4}$ and $h$. Let $\varphi_1: [-1,1]\rightarrow[-1,1]$ be the function which takes values
\begin{align*}
(0,0),(\pm \frac{1}{4},\pm \frac{1}{2}),(\pm \frac{1}{2},\pm \frac{3}{4}),(\pm1,\pm1)
\end{align*} 
with linear interpolation in between. Let $\varphi_t(x) = (1-t)x+t\varphi_1(x)$ be a homotopy between the identity and $\varphi_1$. We now write down a homotopy between $b_{\omega_1,\omega_2}^{H'}b_{\omega_0,\omega_1}^{H}[f]$ and $b_{\omega_0,\omega_2}^{H\ast H'}[f]$ . Define
\begin{align*}
G': I_{p}(c[-1,1]^n)=c([-1,1]^n \times [0,1]) &\rightarrow c[-1,1]^n \\
(hx_1,\dots,hx_n,h,ht)&\mapsto (h\varphi_t(x_1),\dots,h\varphi_t(x_n),h)
\end{align*}
for $t\in [0,1]$. Let $G$ be the composition $b_{\omega_1,\omega_2}^{H'}b_{\omega_0,\omega_1}^{H}[f]\circ G'$, which clearly restricts to $b_{\omega_1,\omega_2}^{H'}b_{\omega_0,\omega_1}^{H}[f]$ on $c([-1,1]^n\times \{0\})$ and $b_{\omega_0,\omega_2}^{H\ast H'}[f]$ on $c([-1,1]^n\times \{1\})$. See Figure 
\ref{fig:bpic4}.\\

\begin{figure}
\centering
  \centering
  \includegraphics[width=0.5\linewidth]{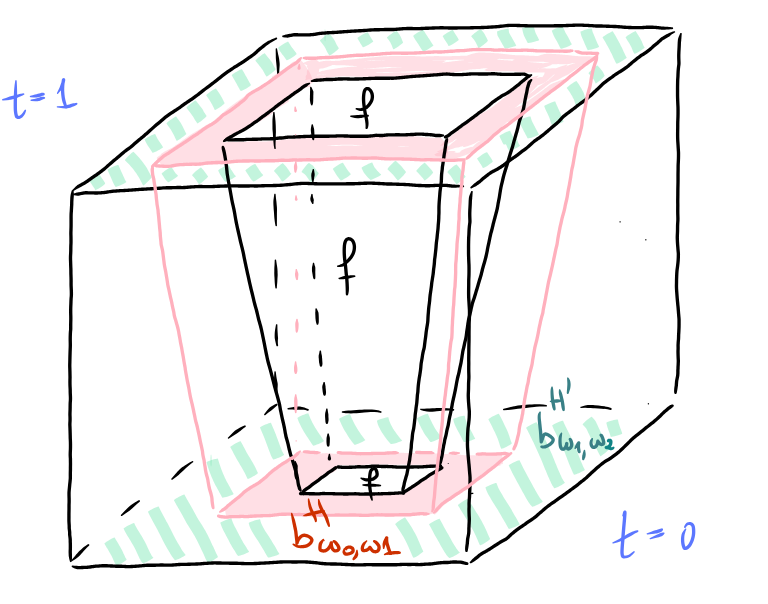}
  \caption{The coarse homotopy $G$ between $b_{\omega_1,\omega_2}^{H'}b_{\omega_0,\omega_1}^{H}[f]$ and $b_{\omega_0,\omega_2}^{H\ast H'}[f]$, restricted to height $h$.}
  \label{fig:bpic4}
\end{figure}

To show that this is controlled involves a straightforward calculation: 
\begin{align*}
|\varphi_{t}(x)-\varphi_{l}(y)| &\leq |\varphi_t(x)-\varphi_t(y)|+|\varphi_t(y)-\varphi_l(y)| \\
|\varphi_t(x)-\varphi_t(y)| &= |(1-t)x+t\varphi_1(x)-(1-t)y-t\varphi_1(y)|\\
&\leq |x-y| + |\varphi_1(x)-\varphi_1(y)|\leq 3|x-y|\\
|\varphi_t(y)-\varphi_l(y)| &= |(1-t)y+t\varphi_1(y)-(1-l)y-l\varphi_1(y)| \leq 2|t-l|
\end{align*}
for $x,y\in [-1,1]$. For $(hx_1,\dots,hx_n,h,ht)$ and $(sy_1,\dots,sy_n,s,sl)$ we have the estimate
\begin{align*}
 \sum_{i=1}^{n} |h\varphi_{t}(x_i)-s\varphi_{l}(y_i)| &\leq \sum_{i=1}^{n} |h\varphi_{t}(x_i)-h\varphi_{l}(y_i)|+|h\varphi_{l}(y_i)-s\varphi_{l}(y_i)|\\
&\leq \sum_{i=1}^{n} 3h(|t-l|+|x_i-y_i|) + |h-s| \\
&\leq \sum_{i=1}^{n}3(|ht-sl|+|hx_i-sy_i|)+7|h-s|\\
&=3n|ht-sl|+\|h\mathbf{x}-s\mathbf{y}\|+7n|h-s|
\end{align*}
Hence $G$ is controlled. For a bounded set $K$ the preimage is contained within $(b_{\omega_1,\omega_2}^{H'}b_{\omega_0,\omega_1}^{H}[f]^{-1}(K) \times [0,\infty)) \cap I_{p}(c[-1,1]^n)$ which is bounded. 
\item The homotopy is defined as follows: let $V\subset I_{p}(c[-1,1]^n)$ be the set
\begin{align*}
V:=\{(h\mathbf{x} ,h,ht)\in I_{p}(c[-1,1]^n)\,|\, 2\|\mathbf{x}\|_{\infty}-1\leq t\}
\end{align*}
Then we define:
\begin{align*}
G_{|V}: V&\rightarrow Z\\
(h\mathbf{x},h,ht)&\mapsto f(\frac{2h\mathbf{x}}{t+1},h)
\end{align*}
and on the complement we define $G$ by taking each height to $\omega_0(h)$. See Figure \ref{fig:bpic5}.

\begin{figure}[H]
\centering
  \centering
  \includegraphics[width=0.5\linewidth]{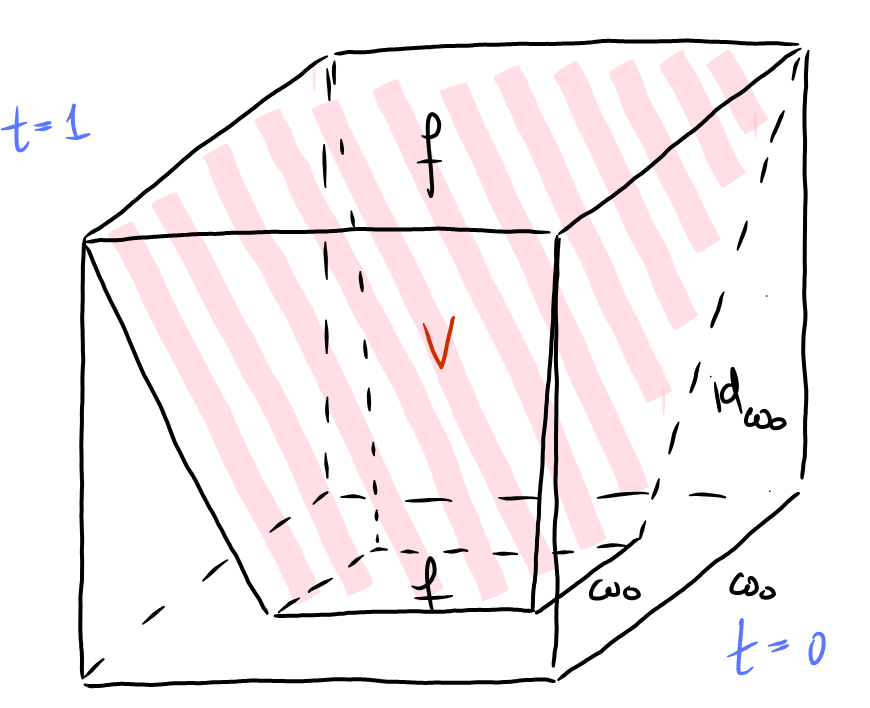}
  \caption{The coarse homotopy $G$ between $b^{id_{\omega_0}}_{\omega_0,\omega_0}[f]$ and $[f]$, restricted to height $h$.}
  \label{fig:bpic5}
\end{figure}

Let us compute the constants once again, for two points $(h\mathbf{x},h,ht)$ and $(s\mathbf{y},s,sl)$ in $V$.
\begin{align*}
\|\frac{2h\mathbf{x}}{t+1}-\frac{2s\mathbf{y}}{l+1}\| &= 2\|\frac{h\mathbf{x}}{t+1}-\frac{s\mathbf{x}}{t+1}+\frac{s\mathbf{x}}{t+1}-\frac{s \mathbf{y}}{t+1}+\frac{s\mathbf{y}}{t+1} - \frac{s \mathbf{y}}{l+1}\|\\
&\leq 2(|h-s|\|\mathbf{x}\| + s \|\mathbf{x}-\mathbf{y}\|+s\|\mathbf{y}\| |\frac{1}{t+1}-\frac{1}{l+1}|)\\
&\leq 2(n|h-s| + n|h-s| + \|h\mathbf{x}-s \mathbf{y}\|+ ns |t-l|)\\
&\leq 6n|h-s|+2\|h\mathbf{x}-s \mathbf{y}\|+2n |ht-sl|
\end{align*}
Hence $G$ is controlled. Properness is shown in a similar way to the previous points. 

\item The diagram already commutes on the level of maps: take any representative for $[f]\in \pi_{n}^{c}(X,\omega_0)$ and the result is that on $U$ we have $\mathfrak{a}\circ f \circ \mathfrak{r}$ and on the complement of $U$ we run along from $\mathfrak{a}\omega_0$ to $\mathfrak{a}\omega$ with the homotopy $\mathfrak{a}\circ H$. 

\end{enumerate}

\end{proof}

  \begin{lemma} \label{arrow1} Let $X,Y$ be bornological coarse spaces and $\omega: [1,\infty)\rightarrow X$ a coarse base ray. Let $\mathfrak{a}:X\rightarrow Y$ be a coarse homotopy equivalence. The induced map $\mathfrak{a}_{*}: \pi_{n}^c(X,\omega)\rightarrow \pi_n^c(Y,\mathfrak{a}\circ \omega)$ is an isomorphism for all $n\geq 0$. 
 \end{lemma}

\begin{proof} Let $\mathfrak{b}: Y\rightarrow X$ be coarse homotopy inverse, that is, a coarse map such that $\mathfrak{b}\circ \mathfrak{a}: X\rightarrow X$ and $\mathfrak{a}\circ \mathfrak{b}:Y\rightarrow Y$ are coarsely homotopic to $\id_X$ and $\id_Y$ respectively. Call these homotopies $H: I_{p_{X}}(X)\rightarrow X$ and $\tilde{H}: I_{p_{Y}}(Y)\rightarrow Y$. We show that the composition 
\begin{align*}
 \pi_{n}^c(X,\omega)\xrightarrow{\mathfrak{a}_{*}} \pi_n^c(Y,\mathfrak{a} \omega) \xrightarrow{\mathfrak{b}_{*}}\pi_n^c(X,\mathfrak{b}\mathfrak{a}\omega)  \xrightarrow{b_{\mathfrak{b}\mathfrak{a}\omega,\omega}^{H}} \pi_{n}^c(X,\omega)
\end{align*}
is the identity, which gives us that $\mathfrak{a}_{*}$ is injective. Let $n\geq 1$. The image of $[f]\in \pi_n^c(X,\omega)$ under the composition is represented by the following map: on $U$ we have $\mathfrak{b}\circ \mathfrak{a} \circ f \circ \mathfrak{r}$ and on the complement of $U$ we travel radially outwards along the path between $\mathfrak{b}\mathfrak{a}\omega$ and $\omega$ as given by the homotopy $H$. Let us define a homotopy between this map and the identity. The idea is the same as the proof of point $6$ in Proposition \ref{baseddefined}. Let $V\subset I_{p}(c[-1,1]^n)$ be the set 
\begin{align*}
V:=\{(h\mathbf{x} ,h,ht)\in I_{p}(c[-1,1]^n)\,|\, 2\|\mathbf{x}\|_{\infty}-1\leq t\}
\end{align*}
Then we define:
\begin{align*}
G_{|V}: V&\rightarrow X\\
(h\mathbf{x},h,ht)&\mapsto H(f(\frac{2h\mathbf{x}}{t+1},h), p_{X}(f(\frac{2h\mathbf{x}}{t+1},h))t)
\end{align*}
That is, on each time slice $ht$ we scale, apply $f$, then take the point along the homotopy $H$ with appropriate rescaling with $p_{X}: X\rightarrow [1,\infty)$.  On the complement of $V$ we just "finish" the path between $\mathfrak{b}\mathfrak{a}\omega$ and $\omega$ as given by the homotopy $H$. See Figure \ref{fig:bpic6}. The concrete formula for $G|_{V^c}$ is:
\begin{align*}
G|_{V^c}: I_{p}(c[-1,1]^n)\setminus V &\rightarrow X
\\(h\mathbf{x}, h,ht)&\mapsto H(\omega(h),p_{X}(\omega(h))(2\|\mathbf{x}\|_{\infty}-1))
\end{align*}

\begin{figure}
\centering
  \centering
  \includegraphics[width=0.5\linewidth]{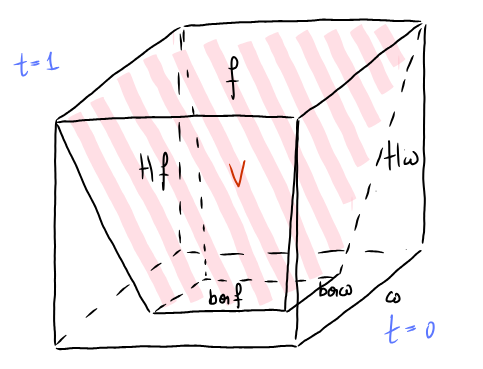}
  \caption{The coarse homotopy $G$ between $b_{\mathfrak{b}\mathfrak{a}\omega,\omega}^{H}\mathfrak{b}_{*}\mathfrak{a}_{*}[f]$ and $[f]$, restricted to height $h$.}
  \label{fig:bpic6}
\end{figure}

Observe that the images agree on the intersection of $V$ and the closure of its complement: on $V$ for a fixed time $ht$ and points satisfying $2\|\mathbf{x}\|_{\infty} = t+1$ we have that $\|\frac{2h\mathbf{x}}{t+1}\|_{\infty} = h$, and so these points get mapped to $H(\omega(h), p_{X}(\omega(h))t)$. The formula for $G|_{V^c}$ also gives for points where $2\|\mathbf{x}\|_{\infty}-1= t$ the result $G|_{V^c}(h\mathbf{x}, h,ht) = H(\omega(h), p_{X}(\omega(h))t)$, as claimed.\\

To show that this map is controlled, we again consider points first in $V$. Since $H$ is coarse, it takes entourages to entourages. It therefore is sufficient to compute coarse estimates on $I_{p_{X}}(X)$.  Let $(h\mathbf{x},h,ht)$ and $(s\mathbf{y}, s,sl)$ be in $V$. We already know that:
\begin{align*}
\|\frac{2h\mathbf{x}}{t+1}-\frac{2s\mathbf{y}}{l+1}\| + |h-s| \leq (6n+1)|h-s|+2\|h\mathbf{x}-s \mathbf{y}\|+2n |ht-sl|
\end{align*}
from point $6$ of Proposition \ref{baseddefined}. Since the composition of controlled maps is controlled and $p_{X}$ is controlled, so is $p_{X}\circ f$. Let $D$ be the large-scale Lipschitz constant for $p_{X}\circ f$.
\begin{align*}
|pf(&\frac{2h\mathbf{x}}{t+1},h)t-pf(\frac{2s\mathbf{y}}{l+1},s)l| \\
&=|pf(\frac{2h\mathbf{x}}{t+1},h)t-pf(\frac{2h\mathbf{x}}{t+1},h)l+pf(\frac{2h\mathbf{x}}{t+1},h)l-pf(\frac{2s\mathbf{y}}{l+1},s)l| \\
&\leq |t-l||pf(\frac{2h\mathbf{x}}{t+1},h)| +|pf(\frac{2h \mathbf{x}}{t+1},h)-pf(\frac{2s\mathbf{y}}{l+1},s)|\\
&\leq  |t-l||pf(\frac{2h\mathbf{x}}{t+1},h)| +  D(\|\frac{2h\mathbf{x}}{t+1}-\frac{2s\mathbf{y}}{l+1}\| + |h-s|) +D
\end{align*}
To bound the value of $|pf(\frac{2h\mathbf{x}}{t+1},h)|$ with relation to $h$, we consider the point $(\frac{2\mathbf{x}}{t+1},1)$, the projection of $(\frac{2h\mathbf{x}}{t+1},h)$ to the base of the cone. We have
\begin{align*}
\|(\frac{2h\mathbf{x}}{t+1},h)-(\frac{2\mathbf{x}}{t+1},1)\| = \frac{2(h-1)}{t+1} \|\mathbf{x}\| + (h-1)\\
|pf(\frac{2h\mathbf{x}}{t+1},h)| \leq |pf(\frac{2h\mathbf{x}}{t+1},h)-pf(\frac{2\mathbf{x}}{t+1},1)| + |pf(\frac{2\mathbf{x}}{t+1},1)|
\end{align*}
Since $f(x_0,1) = \omega(1)$ for a point $x_0\in \partial [0,1]^n$ and $(\frac{2\mathbf{x}}{t+1},1)$ is distance at most $2n$ from $(x_0,1)$, there exists a constant $E$ such that $|pf(\frac{2\mathbf{x}}{t+1},1)|\leq E$ independent of $\mathbf{x}$ and $t$.
\begin{align*}
|pf(\frac{2h\mathbf{x}}{t+1},h)| \leq D(\frac{2n}{t+1}(h-1)+(h-1)) + D + E \leq (2n+1)Dh + D+E
\end{align*}
Putting the computation together:
\begin{align*}
|pf(&\frac{2h\mathbf{x}}{t+1},h)t-pf(\frac{2s\mathbf{y}}{l+1},s)l| \\
&\leq  |t-l|((2n+1)Dh + D+E) +  D((6n+1)|h-s|+2\|h\mathbf{x}-s \mathbf{y}\|+2n |ht-sl|) +D\\
&\leq C(|h-s|+\| h\mathbf{x}-s\mathbf{y}\| + |ht-sl|) +C
\end{align*}
For some $C>0$. These estimates are sufficient to show that $G$ is controlled on $V$. \\

Let $B$ be the large-scale Lipschitz constant for $\omega$, $B'$ for $p\circ \omega$. On the complement of $V$ we that $|\omega(h)-\omega(s)| \leq B |h-s| + B$ and 
\begin{align*}
|p\omega(h)(2&\|\mathbf{x}\|_{\infty}-1) - p\omega(s)(2\|\mathbf{y}\|_{\infty}-1)| \\
&\leq 2|p\omega(h)\|\mathbf{x}\|_{\infty}-p\omega(s)\|\mathbf{y}\|_{\infty}|+|p\omega(h)-p\omega(s)| \\
&\leq (1+2\|\mathbf{y}\|_{\infty})|p\omega(h)-p\omega(s)| +  2p\omega(h)|\|\mathbf{x}\|_{\infty}-\|\mathbf{y}\|_{\infty}| \\
&\leq 3(B'|h-s|+B') +2p\omega(h) \|\mathbf{x}-\mathbf{y}\|_{\infty}\\
&\leq 3(B'|h-s|+B') +2(B'h+B'+p\omega(1)) \|\mathbf{x}-\mathbf{y}\|_{\infty}\\
&\leq 3(B'|h-s|+B') +2B'(\|h\mathbf{x}-s\mathbf{y}\|_{\infty}+|h-s|)+2(B'+p\omega(1))\\
&\leq C(\|h\mathbf{x}-s\mathbf{y}\|_{\infty}+|h-s|)+C
\end{align*}
For some $C>0$. These estimates are sufficient to show that $G$ is controlled on $V^{c}$.\\

To show that $G$ is proper, let $K\subset X$ be a bounded subset. Consider first $V^c$. For a fixed height $h$ in $V^c$, its image under $G$ has nonempty intersection with $K$ if $H((\omega(h)\times [0,\infty)) \cap I_{p}(X)) \cap K \neq \emptyset$. Since $H$ and $\omega$ are proper, the $h's$ for which this is true is contained within a bounded set $[1,l]$. Hence $V^c \cap G^{-1}(K)$ is bounded. For $V$, a similar argument holds. $G|_V$ is the composition of a rescaled version of the map $(h\mathbf{x},h,ht)\mapsto (f(h\mathbf{x},h),p_{X}f(h\mathbf{x},h)t)$ and the homotopy $H$. Since $H$ is proper, $H^{-1}(K)$ is contained within $Z\times [0,\infty)\cap I_{p}(X)$ for $Z$ bounded. Now since $f$ is proper, $G^{-1}(K)\cap V$ is contained within a bounded set of heights $[1,l']$. Taking the maximum of $l,l'$ we see that $G^{-1}(K)$ is bounded. \\ 

We now show that the composition 
\begin{align*}
\pi_n^c(Y,\mathfrak{a} \omega) \xrightarrow{\mathfrak{b}_{*}}\pi_n^c(X,\mathfrak{b}\mathfrak{a}\omega) \xrightarrow{\mathfrak{a}_{*}} \pi_{n}^{c}(Y,\mathfrak{a}\mathfrak{b}\mathfrak{a}\omega) \xrightarrow{b_{\mathfrak{a}\mathfrak{b}\mathfrak{a}\omega,\mathfrak{a}\omega}^{\tilde{H}}} \pi_{n}^{c}(Y,\mathfrak{a}\omega)
\end{align*}
is the identity, which gives us that $\mathfrak{a}_{*}$ is surjective. Let $[f]\in \pi_n^c(Y, \mathfrak{a}\omega)$. The image of $[f]$ under the composition is represented by the following map: on $U$ we have $\mathfrak{a}\circ \mathfrak{b}\circ f \circ \mathfrak{r}$ and on the complement of $U$ we travel radially outwards along the path between $\mathfrak{a}\mathfrak{b}\mathfrak{a}\omega$ and $\mathfrak{a}\omega$ given by image of the base ray $\mathfrak{a}\omega$ under the homotopy $\tilde{H}$. The homotopy between this and the original $f$ is exactly the same idea as the construction of $G$, and the proof of its coarseness is also identical. By point $6$ of Proposition \ref{baseddefined}, the diagram

\begin{figure}[H]
\center
\begin{tikzcd}
&\pi_n^{c}(X,\mathfrak{b}\mathfrak{a}\omega)\arrow {r}{\mathfrak{a}_{*}} \arrow{d}{b_{\mathfrak{b}\mathfrak{a}\omega,\omega}^H}
&\pi_n^{c}(Y,\mathfrak{a}\mathfrak{b}\mathfrak{a}\omega) \arrow{d}{b_{\mathfrak{a}\mathfrak{b}\mathfrak{a}\omega,\mathfrak{a}\omega}^{\mathfrak{a}H}}
  \\
&\pi_{n}^c(X,\omega) \arrow {r}{\mathfrak{a}_{*}} 
&\pi_{n}^c(Y,\mathfrak{a}\omega)

\end{tikzcd}
\end{figure}

commutes. Since the vertical arrows are isomorphisms, $\mathfrak{a}_{*}: \pi_{n}^{c}(X,\omega)\rightarrow \pi_n^{c}(Y,\mathfrak{a}\omega)$ is also surjective. Combining the statements, $\mathfrak{a}_{*}$ is an isomorphism. \\

For $n=0$, both compositions
\begin{align*}
 \pi_{0}^c(X,\omega)\xrightarrow{\mathfrak{a}_{*}} \pi_0^c(Y,\mathfrak{a} \omega) \xrightarrow{\mathfrak{b}_{*}}\pi_0^c(X,\mathfrak{b}\mathfrak{a}\omega)  \\
\pi_0^c(Y,\mathfrak{a} \omega) \xrightarrow{\mathfrak{b}_{*}}\pi_0^c(X,\mathfrak{b}\mathfrak{a}\omega) \xrightarrow{\mathfrak{a}_{*}} \pi_{0}^{c}(Y,\mathfrak{a}\mathfrak{b}\mathfrak{a}\omega) 
\end{align*}

are the identity. This is because for a base ray $[\tau] \in \pi^c_{0}(X,\omega)$, the map $\mathfrak{b}\mathfrak{a}\tau$ is coarsely homotopic to $\tau$. The other composition is analogous.

 \end{proof}

This shows that coarse homotopy groups are invariant under coarse homotopy equivalence. This, combined with u-continuity, means we can make statements about the universal cone $\tilde{\mathcal{O}}(X)$ over a compact metric space $X$. This is the bornological coarse bornological space defined as 
\begin{align*}
\tilde{\mathcal{O}}(X):= \colim_{\phi} \tilde{\mathcal{O}}_{\phi}(X)
\end{align*}
for all $\mathcal{T}$-admissible functions $\phi: [0,\infty) \rightarrow \mathcal{P}(X\times X)$.\\

\begin{cor} The Euclidean cone $cX$ has the same weak coarse homotopy type as $\tilde{\mathcal{O}}(X)$. 
\end{cor}

\begin{proof} We know that $cX,c'X$ and $\tilde{\mathcal{O}}_{\phi_1}(X)$ are coarsely homotopy equivalent. It is shown in the proof of Lemma $8.7$ in \cite{bunke2020coarse}  that if a $\mathcal{T}$-admissible function $\phi$ takes values in $\mathcal{T}$ (and not just $\mathcal{P}(X\times X)$), is monotone, and satisfies $\phi(0) = X \times X$, then for every $\mathcal{T}$-admissible function $\phi'$ such that $\phi(t) \subset \phi'(t)$ for all $t\in [0,\infty)$ the canonical inclusion induced by the set-wise identity 
\begin{align*}
k^{\phi'}_{\phi} : \tilde{\mathcal{O}}_{\phi}(X) \rightarrow \tilde{\mathcal{O}}_{\phi'}(X)
\end{align*}
is a coarse homotopy equivalence. Furthermore, we have that 
\begin{align*}
\tilde{\mathcal{O}}(X) \simeq \colim_{\phi'\geq \phi} \tilde{\mathcal{O}}_{\phi}(X)
\end{align*}
is an equivalence in $\mathbf{BornCoarse}$. The function $\phi=\phi_1$ satisfies these conditions.\\

Let $\omega$ be a coarse base ray in $cX$. Consider the based bornological coarse space $(cX,\omega)$. By invariance under coarse homotopy equivalence and u-continuity, we obtain 
\begin{align*}
\pi^c_{n}(cX,\omega) \cong \pi^c_{n}(c'X,\omega) \cong \pi^c_{n}(\tilde{\mathcal{O}}_{\phi_1}(X), Log\omega) \cong \pi^c_{n}(\colim_{\phi'\geq \phi_1} \tilde{\mathcal{O}}_{\phi_1}(X),Log \omega) \\
=  \pi^c_{n}(\tilde{\mathcal{O}}(X),Log\omega)
\end{align*}
for all $n\geq 0$. 
\end{proof}

\newpage
\section{Proof of $\varprojlim^1$ sequence}

The main theorem of this chapter is:

\begin{theorem} \label{maintheorem} Let $X$ be a compact metric space, $x_0\in X$, $\omega$ the standard parametrisation of $c\{x_0\}$, $n\geq 1$. There is a short exact sequence 
\begin{align*}
0\rightarrow {\varprojlim_h}^1  \pi_{n+1}(|\mathcal{U}_{h}|,x_0) \rightarrow  \pi_{n}^c(cX,\omega) \rightarrow  \check{\pi}_n(X,x_0) \rightarrow 0
\end{align*}
corresponding to a sequence of open covers $\mathcal{U}_{h}:=\cup_{z\in Z_h} B(z,\frac{1}{2^h})$ of metric balls centered on a maximal $\frac{1}{2^{h+1}}$-separated subset $Z_{h}$ of $X$. 
\end{theorem}

In the following special case we have the stronger result: 

\begin{cor} \label{good} Let $X$ be a compact metric space such that the open covers $\mathcal{U}_{h}$ are good for sufficiently large $h$. Let $x_0\in X$, $\omega$ the standard parametrisation of $c\{x_0\}$, $n\geq 1$. There is an isomorphism
\begin{align*}
\pi_{n}^{c}(cX,\omega)\cong \pi_{n}(X,x_0)
\end{align*}
\end{cor}

\subsection{Geometric setup}

Before we prove the theorem, let us go through some preliminaries. Recall that $X$ can be embedded isometrically into $C_{b}(X)$ via the Kuratowski embedding $\Psi$ (Definition \ref{Kuratowski}) such that $\|\Psi(x)\| < \frac{1}{4}$ for all $x\in X$ and $\Psi(x_0) = 0$. We identify $X$ with its image in $C_{b}(X)$. 

\begin{definition} Let $(X,d)$ be a metric space. A subset $Z\subset X$ is $C$-separated if $d(x,y)> C$ for all $x,y\in Z$. 
\end{definition}

Let $Z_0 = \{x_0\}$. For $h\in \mathbb{N}$ choose maximal $\frac{1}{2^{h+1}}$-separated subsets $Z_h$ such that $Z_{h}\subset Z_{h+1}$ for all $h$. Consider the open cover $\cup_{z\in Z_h} B(z,\frac{1}{2^h})=: \mathcal{U}_{h}$. The nerves of these open covers $|\mathcal{U}_h|$ correspond to the Rips complexes $R_{<}(Z_h,\frac{1}{2^{h-1}})$. There are simplicial maps $\phi_{h+1}: |\mathcal{U}_{h+1}| \rightarrow |\mathcal{U}_{h}|$ defined on vertices by the identity on $Z_h\subset Z_{h+1}$ and by sending new points in $Z_{h+1}$ to some point $z'$ in $Z_h$ with distance $\leq \frac{1}{2^{h+1}}$ away. These are well-defined since $\mathcal{U}_{h+1}$ is a refinement of $\mathcal{U}_h$: a ball of radius $\frac{1}{2^{h+1}}$ around $z\in Z_{h+1}$ is completely contained in $B(z',\frac{1}{2^h})$ since for $d(z,y)<\frac{1}{2^{h+1}}$ we have $d(y,z')\leq d(y,z)+ d(z,z') < \frac{1}{2^{h+1}} +\frac{1}{2^{h+1}}= \frac{1}{2^{h}}$. These induce morphisms
\begin{align*}
\phi_{h+1}: \pi_{n} (|\mathcal{U}_{h+1}|,x_0)\rightarrow \pi_{n}(|\mathcal{U}_{h}|,x_0)
\end{align*}
and we define $\phi_{l,s}: \pi_{n} (|\mathcal{U}_{l}|,x_0)\rightarrow \pi_{n}(|\mathcal{U}_{s}|,x_0)$ for $l>s$ as $\phi_{s+1} \circ \dots \circ \phi_{l-1} \circ \phi_{l}$. The inverse limit $\varprojlim_{h\in \mathbb{N}_{0}}\pi_{n}(|\mathcal{U}_h|,x_0)$ can be identified with the \v{C}ech (or shape) homotopy group $\check{\pi}_n(X,x_0)=: \varprojlim_{\mathcal{U}}\pi_{n}(|\mathcal{U}|,x_0)$ (see Definition \ref{defcech}). This is because $\{\mathcal{U}_h\}_{h\in \mathbb{N}_0}$ is cofinal with respect to all open covers of $X$. $X$ is compact, so for any open cover $\mathcal{U}$ there exists a finite subcover with some Lebesgue number $\varepsilon>0$. Choose $\frac{1}{2^h}< \varepsilon$ and $\mathcal{U}_h$ is a refinement of $\mathcal{U}$. We return to a more detailed discussion of the shape of $X$ later. \\

Now we define the inverse mapping telescope, which is the main geometric object used in the proof. Consider the mapping cylinder of $\phi_{h}: |\mathcal{U}_{h}| \rightarrow |\mathcal{U}_{h-1}|$. 
\begin{align*}
(|\mathcal{U}_{h}| \times [2^{h},2^{h+1}]) \cup_{|\mathcal{U}_{h}| \times \{2^h\}} |\mathcal{U}_{h-1}| 
\end{align*}

We glue these mapping cylinders together. 

\begin{definition} (The inverse mapping telescope $M$) 
\begin{align*}
M:=|\mathcal{U}_0| \times [1,2] \cup_{\phi_1} |\mathcal{U}_1| \times [2,4] \cup_{\phi_2} \dots \cup_{\phi_h} |\mathcal{U}_h| \times [2^h,2^{h+1}] \cup_{\phi_{h+1}} \dots
\end{align*}
Each simplex $\sigma$ in $|\mathcal{U}_{h}|$ is given the spherical metric (scaled such that the length of edges is $1$) and $M$ is equipped with the largest path metric bounded above by the product metric on each $\sigma \times [2^h,2^{h+1}]$: see Definition \ref{metricdef}. 
\end{definition}

As with the case of $cX$, we have canonical projections $p: M_{X}\rightarrow [1,\infty)$ and $p': M_{X}\rightarrow \coprod_{h}|\mathcal{U}_{h}|$. \\


\begin{definition} (Geometric simplices) Let $\{u_0,\dots,u_m\}$ be vertices in $M$. We say that $[u_0,\dots, u_m]$ is a \textit{geometric simplex} if $\{u_{0},\dots,u_m\}$ is a subset of the vertices in $\sigma \times I$ for some simplex $\sigma$ in $|\mathcal{U}_{h}|$. For $\tau,\sigma$ geometric simplices, we say that $\tau$ is a \textit{geometric face} of $\sigma$ if the vertices defining $\tau$ are a subset of the vertices defining $\sigma$, ie. $vert(\tau)\subset vert (\sigma)$. 
\end{definition}

\begin{definition} (Geometric cylinders) A \textit{geometric cylinder} is a geometric simplex of the form $\sigma \times \{t\}$ or $\sigma \times I$ for some simplex $\sigma$ in $|\mathcal{U}_{h}|$. 
\end{definition}

\begin{figure}
\centering
\begin{minipage}{.5\textwidth}
  \centering
  \includegraphics[width=.9\linewidth]{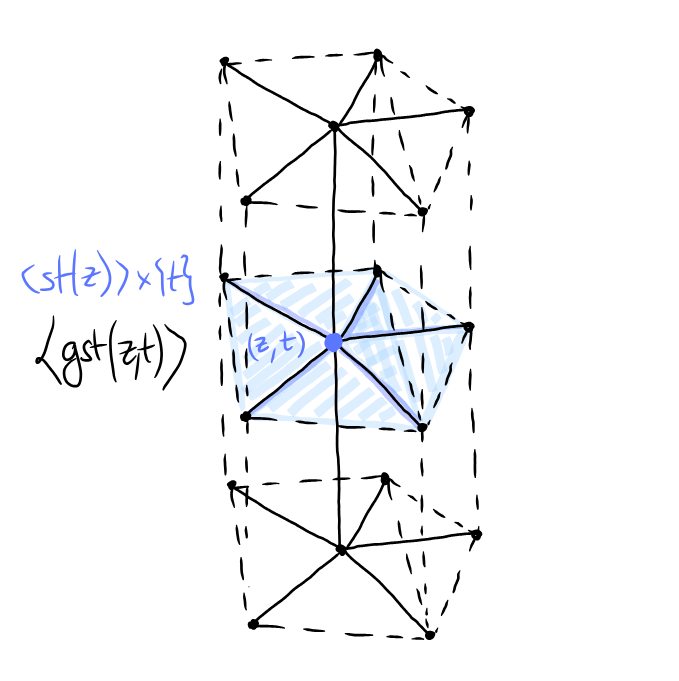}
  \caption{An open geometric star $\langle gst(z,t)\rangle$ for $t=2^h+k$ for $1\leq k<2^h$.}
  \label{fig:starnoglue}
\end{minipage}%
\begin{minipage}{.5\textwidth}
  \centering
  \includegraphics[width=.76\linewidth]{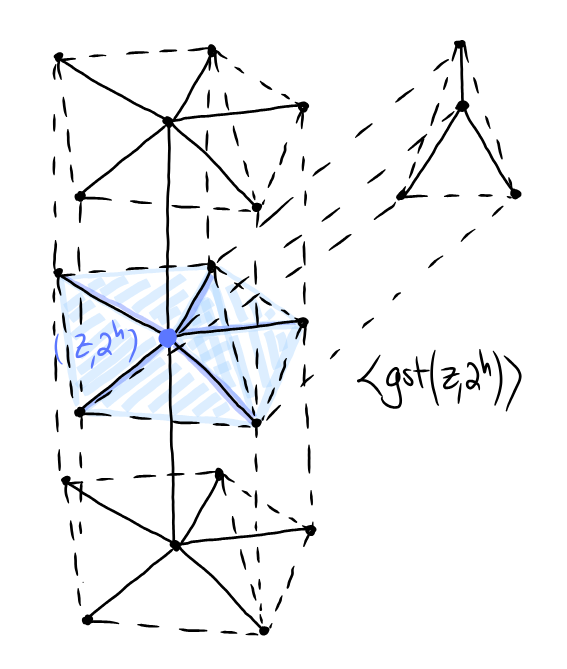}
  \caption{An open geometric star $\langle gst(z,2^h) \rangle$.}
  \label{fig:starglue}
\end{minipage}
\end{figure}

\begin{definition} (Geometric stars) Let $u = (z, 2^h + k)$ be a vertex in $M$. We define the \textit{open geometric star} of $u$ in $M$, denoted by $\langle gst(u)\rangle$, as follows. If $1\leq k< 2^{h}$ (ie. $u$ is not in a gluing component) then we let
\begin{align*}
\langle gst(u)\rangle = \langle st(z)\rangle \times (2^h+k-1, 2^h + k+1)
\end{align*}
where $\langle st(z)\rangle$ denotes the open star of $z$ in $|\mathcal{U}_{h}|$ (see Figure \ref{fig:starnoglue}). If $u=(z',2^h)$ then we define
\begin{align*}
\langle gst(u)\rangle = (\cup_{z\in \phi_{h}^{-1}(z')} \langle st(z)\rangle \times [2^h, 2^h+1)) \bigcup (\langle st(z')\rangle \times (2^{h}-1,2^h])
\end{align*} 
where $\langle st(z)\rangle$ denotes the open star in $|\mathcal{U}_{h}|$ and $\langle st(z')\rangle$ the open star in $|\mathcal{U}_{h-1}|$. (see Figure \ref{fig:starglue}).  The \textit{closed geometric star} $|gst(u)|$ is defined as the closure of $\langle gst(u)\rangle$.
\end{definition}

The idea of this definition is that the geometric star of $u$ is the union of its stars over all possible choices of product CW-structure given the fixed vertex set and the geometric information of which vertices form a simplex in each $|\mathcal{U}_{h}|$.

\begin{figure}
    \centering
  \includegraphics[width=0.4\linewidth]{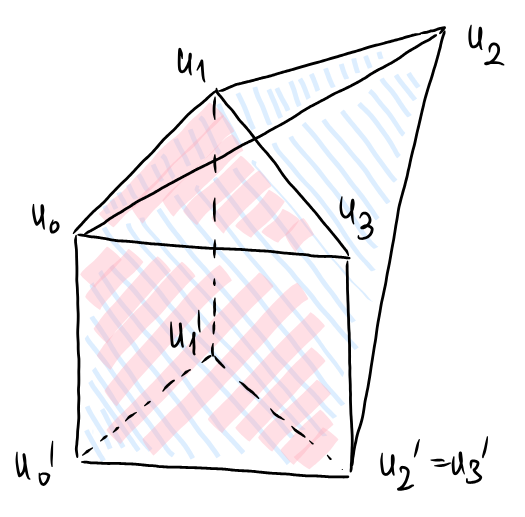}
  \caption{$K:= ([u_0,u_1,u_2] \times I)\cup ([u_0,u_1,u_3] \times I)\cup [u'_0, u'_1, u'_{2}]$. Geometric simplices of $K$ are all subsets of $\{u_0,u_1,u_2,u'_0, u'_1, u'_{2}\}$ and all subsets of $\{u_0,u_1,u_3,u'_0, u'_1, u'_{2}\}$.}
  \label{fig:geomsimpl}
\end{figure}

\begin{figure}
    \centering
  \includegraphics[width=0.4\linewidth]{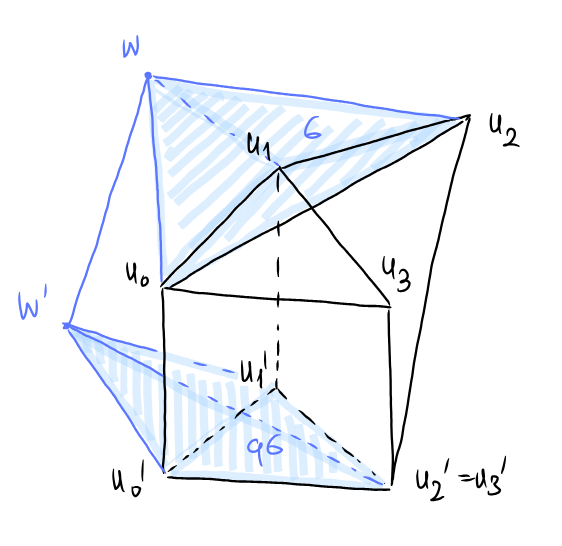}
  \caption{$\sigma = [w,u_0,u_1,u_2], q\sigma = [w',u'_0, u'_{1},u'_{2}]$. Let $K' = K\cup (\sigma \times I \cup q\sigma)$. $w\in \mathcal{N}_{\Delta}$}
  \label{fig:geomneigh}
\end{figure}

\begin{remark} (Description of geometric simplices and stars in $M$) Let us restrict for a moment to heights in the interval $[2^h, 2^{h+1}]$. In $M$, this corresponds to the subset
\begin{align*}
(|\mathcal{U}_{h}| \times [2^h,2^{h+1}]) \cup_{\phi_{h}} |\mathcal{U}_{h-1}|
\end{align*}
A geometric simplex can be of the form $\sigma \times \{2^{h}+k \}$ where $\sigma$ is a simplex in $|\mathcal{U}_{h}|$ and $1\leq k\leq 2^h$, or a simplex in $|\mathcal{U}_{h-1}|$ if $k=0$. Alternatively it is a subset of the vertices $\{u_0,\dots, u_j, u'_0, \dots, {u}'_{j}\}$ where the $u_i = (z_i,2^{h} + k)$ are vertices in $|\mathcal{U}_{h}| \times \{2^{h} + k\}$ for $1\leq k\leq 2^{h}$ and $u'_i = (z_i,2^{h} + k-1)$ refers to the same vertex in $|\mathcal{U}_{h}| \times \{2^{h} + k-1\}$ (which are identified with their images $(\phi_{h}(z_i),2^h) \in|\mathcal{U}_{h-1}| \times \{2^{h}\}$ if $k=1$). Therefore, a set of vertices $\{u_0,\dots,u_m\}$ forms a geometric simplex in $(|\mathcal{U}_{h}| \times [2^h,2^{h+1}]) \cup_{\phi_{h}} |\mathcal{U}_{h-1}|\subset M$ exactly when their height variables differ by at most $1$ and the following condition is satisfied: if $\{u_0,\dots,u_m\}\subset |\mathcal{U}_{h}| \times [2^h+1,2^{h+1}]$ or $\{u_0,\dots,u_m\}\subset |\mathcal{U}_{h-1}| \times \{2^h\}$ then their projections to $ |\mathcal{U}_{h}|$ (resp. $|\mathcal{U}_{h-1}|$) must span a simplex. Otherwise the vertices are of the form $\{u_0,\dots,u_k,w_{k+1},\dots,w_m\}\subset  (|\mathcal{U}_{h}| \times [2^h,2^{h}+1]) \cup_{\phi_{h}} |\mathcal{U}_{h-1}|$ where $u_i = (z_i,2^{h}+1)$ for $0\leq i\leq k$ and $w_{i} = (z'_i,2^h)$ for $k+1\leq i\leq m$. They define a geometric simplex exactly when there exist vertices $z_{k+1},\dots, z_m$ such that $z_i\in \phi^{-1}_{h}(z'_i)$ for $k+1\leq i\leq m$ and $\{z_0,\dots,z_k,z_{k+1},\dots,z_m\}$ spans a simplex $\sigma\subset |\mathcal{U}_{h}|$. In this case, $\{u_0,\dots,u_k,w_{k+1},\dots,w_m\}$ is a subset of $\sigma \times I = sp\{z_0,\dots,z_k,z_{k+1},\dots,z_m\} \times [2^h,2^h+1]$, where $sp(\cdot)$ denotes "the simplex spanned by". See Figure \ref{fig:geomsimpl}.\\

\begin{notation} From now on, we denote by $M_{[a,b]}$ the mapping telescope with heights in $[a,b]$. For example, we refer to the set  $(|\mathcal{U}_{h}| \times [2^h,2^{h+1}]) \cup_{\phi_{h}} |\mathcal{U}_{h-1}| \subset M$ as $M_{[2^h,2^{h+1}]}$. 
\end{notation}

\begin{lemma} \label{geomsexist} For a set of vertices $\{u_0,\dots u_m\}$, a geometric simplex $[u_0,\dots, u_m]$ exists exactly when $\cap_{i=0}^{m} \langle gst(u_i)\rangle \neq \emptyset$
\end{lemma}  

\begin{proof} Let $u_i = (z_i,t_i)$ for $0\leq i\leq m$. Observe that $\langle gst(z_i,t_i) \rangle \cap \langle gst(z_j, t_j)\rangle \neq \emptyset$ only if $|t_i-t_j| \in \{0,1\}$ for all $0\leq i,j\leq m$. If $t_i = 2^h+k$ for some $1\leq k < 2^{h}$ and all $i$ then this condition reduces to $\cap_{i=0}^{m}\langle st(z_i) \rangle \neq \emptyset$, which is equivalent to the existence of a simplex $sp\{z_0,\dots,z_m\}$ in $|\mathcal{U}_{h}|$. \\

 If $t_i=2^h$ for all $i$ and $\{z_0,\dots,z_{m}\}$ spans a simplex then clearly their geometric stars intersect. Conversely, if $\cap_{i=0}^{m} \langle gst(z_i,t_i) \rangle \neq \emptyset$ that means that either $\cap_{i=0}^{m}\langle st(z_i) \rangle \times (2^h-1, 2^h] \neq \emptyset$ or $\cap_{i=0}^{m}(\cup_{z\in \phi_{h}^{-1}(z_i)} \langle st(z) \rangle \times [2^h, 2^h +1)) \neq \emptyset$. In either case, we conclude that $\cap_{i=0}^{m}\langle st(z_i) \rangle \neq \emptyset$ and therefore $\{z_0,\dots,z_{m}\}$ spans a simplex in $|\mathcal{U}_{h-1}|$. \\

If $t_i\in \{2^h+k, 2^h+k+1\}$ for some $1\leq k <2^{h}$ and all $i$ we are reduced again to the condition that $\cap_{i=0}^{m} \langle st(z_i) \rangle \neq \emptyset$. Finally, if $t_i\in \{2^h,2^h+1\}$ and $\cap_{i=0}^{m} \langle gst(z_i,t_i)\rangle \neq \emptyset$ it means that 
\begin{align*}
\cap_{t_i=2^h+1} \langle st(z_i) \rangle \times (2^h, 2^h+1] \bigcap \cap_{t_j=2^h} \cup_{y\in \phi_{h}^{-1}(z_j)} \langle st(y) \rangle \times [2^h, 2^h+1) \neq \emptyset 
\end{align*}
Staring at this for a while, one concludes that the definition of a geometric simplex is satisfied. Namely, there exist vertices $y_j\in \phi_{h}^{-1}(z_j)$ for $t_{j} = 2^h$ such that the vertices $\cup_{t_{i}=2^h+1} z_i \bigcup \cup_{t_j=2^h} y_j$ span a simplex in $|\mathcal{U}_{h}|$. The converse follows from the definition of geometric simplices. See Figure \ref{fig:geomneigh}. 
 \end{proof}
\end{remark}

\begin{figure}
    \centering
  \includegraphics[width=0.4\linewidth]{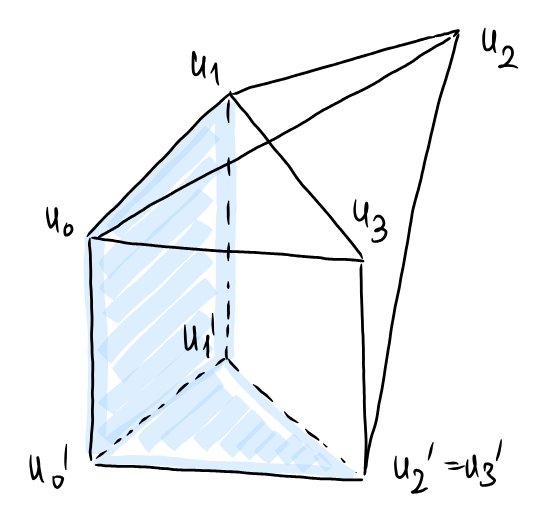}
  \caption{The geometric core $\mathscr{C}_{\Delta} = ([u_0,u_1] \times I)\cup [u'_0,u'_1,u'_{2}]$ of $\Delta:=[u_0,u_1,u'_{2}]$ in $K$. Geometric cylinders $\eta_{\Delta}$ are $([u_0,u_1,u_{2}] \times I) \cup [u'_0,u'_1,u'_{2}]$ and $([u_0,u_1,u_{3}] \times I)\cup [u'_0,u'_1,u'_{2}]$.}
  \label{fig:simplapprox1}
\end{figure}

\begin{definition} (Cores, spanning cylinders and neighbourhoods) Let $\Delta = [u_0,\dots u_m]$ be a geometric simplex.
\begin{itemize}
\item The \textit{geometric core} of $\Delta$, denoted by $\mathscr{C}_{\Delta}$, is defined as follows: if $\Delta$ lies on one height, let $\mathscr{C}_{\Delta} = \Delta$. Otherwise, $\Delta = [x_0,\dots,x_k,w_{k+1},\dots,w_m]$, and we define 
\begin{align*}
\mathscr{C}_{\Delta}= ([x_0,\dots,x_k] \times I )\cup_{[x_0,\dots,x_k] \times \{0\}} sp\{q(x_0),\dots, q(x_k), w_{k+1},\dots,w_m\} 
\end{align*}
where $q$ denotes the projection to the bottom of $\Delta$, ie. $q: |\mathcal{U}_{h}| \times \{2^h+k\} \rightarrow |\mathcal{U}_{h}| \times \{2^h+k-1\}$ for some $h\in \mathbb{N}_{0}$ and $1\leq k\leq 2^h$. Note that $q$ can be identified with the identity for $k>1$ and the gluing map $\phi_{h}$ for $k=1$. 
\item A \textit{geometric cylinder spanned by} $\Delta$, denoted by $\eta_{\Delta}$, is defined as follows:  if $\Delta$ lies on one height, let $\eta_{\Delta} = \Delta$. For $\Delta = [x_0,\dots,x_k,w_{k+1},\dots,w_m]$ we let 
\begin{align*}
\eta_{\Delta}= ([x_0,\dots,x_k, y_{k+1},\dots, y_{m}] \times I )\cup sp\{q(x_0),\dots, q(x_k), w_{k+1},\dots,w_m\} 
\end{align*}
for a choice of $y_i\in q^{-1}(w_i)$ for $k+1\leq i\leq m$. See Figure \ref{fig:simplapprox1}. \\

Let $\mathcal{G}_{\Delta}$ denote the set of geometric cylinders spanned by $\Delta$.
\item The \textit{geometric neighbhourhood} of $\Delta$, denoted by $\mathcal{N}_{\Delta}$, is defined as follows: if $\Delta = [(z_0,2^{h}+k),\dots, (z_m,2^{h}+k)]$ lies on one height, we let 
 \[ \mathcal{N}_{\Delta} := \begin{cases} 
          cl(\cap_{i=0}^{m} \langle st(z_i) \rangle) \times (2^{h}+k-1, 2^{h}+k+1) & 1\leq k< 2^h \\
          \big( cl(\cap_{i=0}^{m} \langle st(z_i) \rangle) \times (2^h-1, 2^h]\big)\bigcup \\
\big( cl(\cap_{i=0}^{m}\cup_{y\in \phi_{h}^{-1}(z_{i})} \langle st(y) \rangle) \times [2^h, 2^h+1)\big)& k=0 \\ 
       \end{cases}  \] 
Note that in this case
\begin{align*}
\bigcap_{i=0}^m \langle gst(u_i) \rangle \subset \mathcal{N}_{\Delta} \subset cl(\bigcap_{i=0}^m \langle gst(u_i) \rangle )
\end{align*}
 If $\Delta$ does not lie on one height, we define $\mathcal{N}_{\Delta}$ as the closure of the intersection of geometric stars
\begin{align*}
cl(\bigcap_{i=0}^{m} \langle gst(u_i) \rangle)
\end{align*}
See Figure \ref{fig:simplapprox1}.
\end{itemize}
\end{definition}

\begin{notation} A geometric simplex can be written in the form $[\sigma,\tau]$, where $\sigma$ is a simplex in some $|\mathcal{U}_{h}| \times \{2^{h}+k\}$ for $1\leq k\leq 2^{h}$ with $\tau$ is a simplex at height $2^{h}+k-1$. If the geometric simplex lies on one height, then either $\sigma$ or $\tau$ is empty. If $\tau\neq \emptyset$, we can write $\mathscr{C}_{[\sigma,\tau]}= (\sigma \times I)\cup sp\{q\sigma,\tau\}$. Likewise $\eta_{[\sigma,\tau]} = (sp\{\sigma,\delta\} \times I) \cup sp\{q\sigma,\tau\}$ for a choice of simplex $\delta \in q^{-1}\tau$. This is an abuse of notation, since we identify $|\mathcal{U}_{h}| \times \{2^{h}+k\} \times [0,1]$ with $|\mathcal{U}_{h}| \times [2^{h}+k-1, 2^h+k]$. \\

We say that a geometric simplex has type $1$ if all its vertices lie on one height, and type $2$ otherwise.  \\

For a simplex $\Delta = [\sigma,\tau]$ of type $2$, we denote by $\mathcal{I}_{0}\mathcal{N}_{\Delta}$ the subset of $\mathcal{N}_{\Delta}$ without the simplices at the height of $\tau$. Likewise, we denote by $\mathcal{I}_{1}\mathcal{N}_{\Delta}$ the subset of $\mathcal{N}_{\Delta}$ without the simplices at the height of $\sigma$.
\end{notation}

The following are easy observations:

\begin{itemize}
\item$ \bigcup_{\eta_{\Delta}\in \mathcal{G}_{\Delta}} \eta_{\Delta} \subset \mathcal{N}_{\Delta}$.
\item$ \mathscr{C}_{\eta_{\Delta}} = \eta_{\Delta}$.
\item $\mathscr{C}_{\Delta} \subset \bigcap_{\eta_{\Delta}\in \mathcal{G}_{\Delta}} \eta_{\Delta}$.
\item If $\Delta'$ is a face of $\Delta$, we have $\mathscr{C}_{\Delta'}\subset \mathscr{C}_{\Delta}$ .
\item If $\Delta'$ is a face of $\Delta$, we have $\mathcal{N}_{\Delta}\subset cl(\mathcal{N}_{\Delta'})$. If $\Delta',\Delta$ are both type $1$ or both type $2$ we have $\mathcal{N}_{\Delta}\subset \mathcal{N}_{\Delta'}$.
\item For a simplex $\Delta = [\sigma,\tau]$ of type 2, a point $(x,t) \in (0,1]$ lies in $\mathcal{N}_{\Delta}$ if and only if $x\in \mu$, $\mu$ a simplex which contains $sp\{\sigma, \delta\}$ as a face, for some simplex $\delta \in q^{-1}(\tau)$. A point $(x,0)$ lies in $\mathcal{N}_{\Delta}$ if and only if $x\in \mu$, $\mu$ a simplex which contains $sp\{q\sigma, \tau\}$ as a face. 
\item For a simplex $\Delta = [\sigma,\tau]$ of type 2 with $\Delta' = [\sigma',\emptyset]$ (resp. $[\emptyset, \tau']$) a face of $\Delta$, we have $\mathcal{I}_{0} \mathcal{N}_{\Delta} \subset \mathcal{N}_{\Delta'}$ (resp. $\mathcal{I}_{1}\mathcal{N}_{\Delta} \subset \mathcal{N}_{\Delta'}$).
\end{itemize}
 We now prove some useful properties about neighbourhoods and cores. 

\begin{lemma} \label{usefulprops} Let $\Delta=[\sigma,\tau]$ be a geometric simplex and let $\Delta'=[\sigma',\tau']$ be a face of $\Delta$, ie. $vert(\tau') \subset vert(\tau), vert(\sigma')\subset vert \sigma$. 
\begin{enumerate}
\item Let $\Delta$ be equipped with the path metric induced from the standard ($L^1$) metric on $\sigma, \tau$, and the product metric on $\sigma \times I$. There is a $2$-Lipschitz deformation retract $H^{\mathscr{C}_{\Delta}}: \mathscr{C}_{\Delta} \times [0,1] \rightarrow \mathscr{C}_{\Delta}$ of $\mathscr{C}_{\Delta}$ to a point. 
\item There is a deformation retract $H^{\mathcal{N}_{\Delta}}: \mathcal{N}_{\Delta} \times [0,1]\rightarrow \mathcal{N}_{\Delta}$ of $\mathcal{N}_{\Delta}$ to a point. 
\item Suppose that $\Delta$ is a geometric simplex with $\Delta'$ as a face. The deformation retract of $\mathcal{N}_{\Delta'}$ restricted to $\mathcal{N}_{\Delta} \cap \mathcal{N}_{\Delta'}$ has image in $\mathcal{N}_{\Delta}$. 
\end{enumerate}
\end{lemma}

\begin{figure}
\centering
  \centering
  \includegraphics[width=0.5\linewidth]{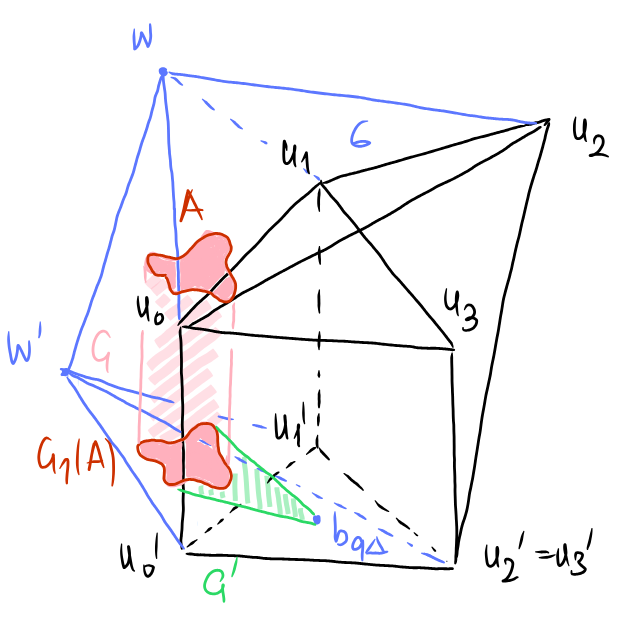}
  \caption{The image of the set $A \times [0,1]$ under $H^{\mathcal{N}_{\Delta}}$. $\sigma = [w,u_0,u_1,u_2]$,\\$q\sigma = [w',u'_0, u'_{1},u'_{2}]$, $\Delta=[u_{0},u_1, u'_{2}] $}
  \label{fig:homotopygeosimpl}
\end{figure}

\begin{proof}
\begin{enumerate}
\item For a simplex of type $1$, we simply take the standard deformation retract of a simplex to its barycentre. Otherwise, we define the homotopies
\begin{align*}
H: \mathscr{C}_{\Delta} \times [0,1]&\rightarrow \mathscr{C}_{\Delta}\\
([x,s],t)&\mapsto [x,(1-t)s]\\
([y,0],t)&\mapsto [y,0]\\
H': sp \{q\sigma,\tau\} \times [0,1]&\rightarrow sp \{q\sigma,\tau\}\\
(y,t)&\mapsto (1-t)y + t b_{sp \{q\sigma,\tau\}}
\end{align*}
where $b_{sp \{q\sigma,\tau\}}$ denotes the barycentre of $sp \{q\sigma,\tau\}$. The composition $H\ast (H'\circ (H_{1} \times \id))$ is a deformation retract of $\mathscr{C}_{\Delta}$ to the barycentre of  $sp \{q\sigma,\tau\}$. It is $2$-Lipschitz in the time variable and $1$-Lipschitz in the $\mathscr{C}_{\Delta}$-variable (since the bonding maps $\phi_{h}$ are distance-decreasing). 

\item If $\Delta = [u_0, \dots, u_m]$ is type $1$, with $u_i=(z_i,t)$ for all $0\leq i\leq m$, then there is first an obvious deformation retract of $\mathcal{N}_{\Delta}$ to  $cl(\cap_{i=0}^{m} \langle st(z_{i}) \rangle) \times \{t\}$, which we denote by $G$. Then we can deformation retract to the barycentre of $[z_{0},\dots, z_m] \times \{t\}$ by a linear homotopy $G'$. \\

If $\Delta$ is type $2$, then it has the form $[x_0,\dots,x_{j},w_{j+1},\dots, w_{m}]$, where the $w$'s lie on a lower height. In this case $\mathcal{N}_{\Delta}$ deformation retracts to the bottom slice (which we denote by $G$), and then to the barycentre of the projection 
\begin{align*}
[q(x_0),\dots, q(x_j),w_{j+1},\dots, w_{m}]
\end{align*}
via a homotopy $G'$, where $q$ is the identity or $\phi$, depending on the height (see Figure \ref{fig:homotopygeosimpl}). To see that this works, let $(x,t)$ be a point in $\mathcal{N}_{\Delta}$ for $t\in (0,1]$. $x$ lies in a simplex $\tau$ contained within the closure of 
\begin{align*}
\cap_{i=0}^{j} \langle st(x_{i}) \rangle \bigcap \cap_{i=j+1}^{m} \langle st(y_i)\rangle
\end{align*} 
for some $y_i \in q^{-1}(w_i)$ for $j+1\leq i\leq m$. This is equivalent to $\tau$ having $[x_0,\dots,x_j, y_{j+1}, \dots, y_m]$ as a face. Since the bonding maps $\phi$ are simplicial, we have that image of $\tau$ under the deformation retract to the bottom slice $t=0$ always lies in the closure of 
\begin{align*}
\cap_{i=0}^{j} \langle st(q(x_{i})) \rangle \bigcap \cap_{i=j+1}^{m} \langle st(w_i)\rangle
\end{align*}
 So, the barycentre of $[q(x_0),\dots, q(x_j),w_{j+1},\dots, w_{m}]$ can be joined to $q(x)$ by a straight line within a simplex. \\ 

If $\mathcal{N}_{\Delta}$ is equipped with the path metric induced from the standard, product metric on each simplex $\sigma \times I$, then $H^{\mathcal{N}_{\Delta}}: = G\ast (G'\circ G_{1})$ is $2$-Lipschitz.

\item If $\Delta = [u_0,\dots, u_m]$ is type $1$, then $u_i = (z_i,t)$ for all $0\leq i\leq m$ and $vert \Delta' \subset \{u_0,\dots, u_m\}$. Let $H^{\mathcal{N}_{\Delta'}} = G\ast (G'\circ G_{1})$ be the concatenation of two homotopies described in point $2.$ above.  It is clear that $G_{|\mathcal{N}_{\Delta} \times [0,1]}(\mathcal{N}_{\Delta} \times [0,1]) \subset \mathcal{N}_{\Delta}$. For $G'$, the straight line between $ x \in cl(\cap_{i=0}^{m} \langle st(z_{i}) \rangle) \times \{t\}$ and the barycentre of $\Delta'$ is contained within $cl(\cap_{i=0}^{m} \langle st(z_{i})\rangle) \times \{t\}$ since $x$ lies within a simplex $\tau$ which contains $\Delta'$ as a face, and simplices are convex. The argument for $\Delta'$ being type $2$ is analogous. \\

Suppose that $\Delta = [\sigma, \tau]$ is type $2$ and $\Delta' = [\sigma',\emptyset]$ is type $1$. Then $\mathcal{N}_{\Delta} \cap \mathcal{N}_{\Delta'} = \mathcal{I}_{0}\mathcal{N}_{\Delta} \cap \mathcal{N}_{\Delta'} = \mathcal{I}_{0}\mathcal{N}_{\Delta}$. It is easy to see that $G_{| \mathcal{I}_{0}\mathcal{N}_{\Delta} \times [0,1]}$ has image in $ \mathcal{I}_{0}\mathcal{N}_{\Delta}$. For a point $(x,t)\in \mathcal{I}_{0}\mathcal{N}_{\Delta}$, with $x\in \mu$, $\mu$ a simplex containing $sp\{\sigma, \delta\}$ for some simplex $\delta\in q^{-1}(\tau)$, $G_{1}(x,t)\in sp\{\sigma, \delta\}$. As before, the straight line between $G_{1}(x,t)$ and $b_{\sigma'}$ lies in $\mu \subset \mathcal{I}_{0}\mathcal{N}_{\Delta}$. The argument for  $\Delta' = [\emptyset,\tau']$ is analogous.

\end{enumerate}
\end{proof}

\begin{remark} The reason why we had to define geometric neighbourhoods $\mathcal{N}_{\Delta}$ differently depending on the type of $\Delta$ is because if $\Delta$ is type $1$ at height $2^h+1$, the closure $cl(\mathcal{N}_{\Delta})$ includes a gluing slice at height $2^h$, so does not deformation retract to height $2^{h}+1$. 
\end{remark}

We now digress into a short discussion about shrinking maps, which will be useful in many parts of the rest of the arguments:

\begin{definition} (Shrinking maps) Let $\gamma: [1,\infty)\rightarrow [1,\infty)$ be a coarse map such that $\gamma(h)\leq h$. We define $s_{\gamma}: cX\rightarrow cX$ as
\begin{align*}
s_{\gamma}: cX&\rightarrow cX\\
(hx,h)&\mapsto (\gamma(h)x, \gamma(h))
\end{align*}
We call $s_{\gamma}$ the \textit{shrinking map associated to $\gamma$}. 
\end{definition}

\begin{lemma} \label{shrinkindep} Let $s_{\gamma},s_{\gamma'}$ be two shrinking maps. There is a linear coarse homotopy $G: I_{p}(cX)\rightarrow cX$ between $s_{\gamma}$ and $s_{\gamma'}$.
\end{lemma}

\begin{proof}
 Since $[1,\infty)$ is a path metric space, $\gamma,\gamma'$ are large-scale Lipschitz. So there exist constants $C_{\gamma}, C_{\gamma'}>1$ such that $|\gamma(t)-\gamma(s)|<C_{\gamma}|t-s|+C_{\gamma}$ for all $t,s\in [1,\infty)$ and similarly for $\gamma'$. For simplicity's sake, define $\eta(h,t) = t \gamma'(h) + (1-t) \gamma(h)$. Consider the homotopy
\begin{align*}
G: I_{p}(cX)&\longrightarrow cX
\\ (hx,ht,h) &\longmapsto (\eta(h,t)x, \eta(h,t))
\end{align*}
This is a homotopy between $s_{\gamma}$ (when $t=0$) and $s_{\gamma'}$ (when $t=1$). We now show that $G$ is a coarse map. 
\begin{align*}
d_{cX}(G(hx,ht,h), G(sy,sl,s)) &\leq \|\eta(h,t)x-\eta(h,t)y\|+ |\eta(h,t)-\eta(s,l)| (\|y\| +1)\\
&< \eta(h,t)\|x-y\| + 2|\eta(h,t)-\eta(s,l)|\\
 \eta(h,t)\|x-y\|  &\leq  h\|x-y\| \leq \|hx-sy\|+|h-s|\|y\| \\
 &< \|hx-sy\|+|h-s|\\
 |\eta(h,t)-\eta(s,l)| &\leq |\eta(h,t)-\eta(s,t)| + |\eta(s,t)- \eta(s,l)|\\
|\eta(h,t)-\eta(s,t)| &\leq t|\gamma'(h) - \gamma'(s)| + (1-t) |\gamma(h) - \gamma(s)|\\
&<(C_{\gamma'}+ C_{\gamma})|h-s| + (C_{\gamma'}+ C_{\gamma})\\
|\eta(s,t)- \eta(s,l)| &\leq |t-l|(\gamma'(s) + \gamma(s))\\
&<|t-l|(C_{\gamma'}s + C_{\gamma'} + C_{\gamma}s +C_{\gamma})\\
&\leq C_{\gamma'}|st-sl| + C_{\gamma}|st-sl|+(C_{\gamma'} + C_{\gamma})\\
& \leq (C_{\gamma'}+ C_{\gamma})(|ht-sl| +|h-s|) + (C_{\gamma'} + C_{\gamma})
\end{align*}
These estimates are sufficient to show that $G$ is controlled. $G$ is proper because $\eta: [1,\infty) \times [0,1]\rightarrow [1,\infty)$ is proper: for a bounded interval $[1,h']$ there exists a $l$ such that $\gamma(h)>h'$ and $\gamma'(h)>h'$ for all $h>l$. Therefore $\eta^{-1}[1,h'] \subset [1,l] \times [0,1]$ and $G^{-1}(c_{[1,h']}X) \subset  c_{[1,l]}(X \times [0,1])$. Therefore $G$ is coarse.  
\end{proof}

In particular, it means that any shrinking map $s_{\gamma}$ is coarsely homotopic to the identity. Observe that the proof also applies to any $\gamma,\gamma'$ coarse and not necessarily shrinking, since there exist constants $D_{\gamma}>C_{\gamma}$, $D_{\gamma'}>C_{\gamma'}$ such that $\gamma(h)< D_{\gamma}h + D_{\gamma}$ and $\gamma'(h)< D_{\gamma'} h+ D_{\gamma'}$ for all $h\in [1,\infty)$. So we have 
\begin{align*}
\eta(h,t) < t (D_{\gamma'}h + D_{\gamma'}) + (1-t)(D_{\gamma}h+D_{\gamma}) < Dh + D
\end{align*} 
where $D$ is the maximum of the two constants. The coarseness of $G$ follows. This concludes the digression. \\

We now define the maps $i: cX\rightarrow M_{X}$ and $R:M_{X}\rightarrow cX$. Let $\rho_{h}= \{\rho_{U}\}_{U\in \mathcal{U}_{h}}$ be a partition of unity subordinate to the cover $\mathcal{U}_{h}$. Define $\varphi_{h}: X\rightarrow |\mathcal{U}_{h}|$ by sending a point $x\in X$ to its image $\sum_{U\in \mathcal{U}_{h}}\rho_{U}(x) [U]$ in barycentric co-ordinates. 

\begin{definition} Let $i: cX\rightarrow M$ be a map defined piecewise as follows:
\begin{align*}
i_h:X \times (2^{h},2^{h+1}] &\rightarrow |\mathcal{U}_{h}| \times (2^{h},2^{h+1}]\\
(tx,t) &\mapsto (\varphi_{h}(x),t)
\end{align*}
Denote by $i$ as the gluing of these maps. 
\end{definition}

Observe that $i$ is not continuous. Nonetheless, it is useful for the following reason:

{\begin{lemma} \label{weirdness} Let $U\subset cX$ be a set of diameter $<Q$. There exists a shrinking map $s$ such that $i\circ s(U)$ is contained within the geometric star of a vertex in $M$. 
\end{lemma}
\begin{proof} Let $r\in \mathbb{N}$ be large enough such that $\frac{Q}{2^r}<\frac{1}{4}$. Let $s_r: cX\rightarrow cX$ be the shrinking map that divides everything by $2^r$:

\[ s_r(tx,t) = \begin{cases} 
          (\frac{t}{2^r}x,\frac{t}{2^r}) & t \geq 2^r \\
          (x,1)& t\in [1,2^r] \\ 
       \end{cases}
    \]

This is obviously $\frac{1}{2^r}$-Lipschitz, so $diam s_r(U)<\frac{Q}{2^r}$. From this we know that $s_r(U)$ intersects at most one gluing component. Assume first that $s_r(U)\subset X \times [2^h,2^{h+1}]$ for some $h\in \mathbb{N}_{0}$. Let $(tx,t),(sy,s)\in s_r(U)$. By a simple calculation
\begin{align*}
t\|x-y\| = \|tx-sy+sy-ty\| \leq \|tx-sy\| + |t-s| \|y\|\\
\|x-y\| \leq \frac{1}{t}(\|tx-sy\| + |t-s| \|y\|) < \frac{Q}{2^{h+r}}<\frac{1}{2^{h+2}}
\end{align*}
Hence, the points projected onto $X$ are distance less than the Lebesgue number of $\mathcal{U}_{h}$ (which is $\frac{1}{2^{h+1}}$) and therefore $i\circ s_r(U)$ is contained within the geometric star of a vertex in $|\mathcal{U}_{h}| \times [2^h,2^{h+1}]\subset M_{[2^h,2^{h+1}]}$. \\

Suppose now that $s_r(U) \subset X \times (2^{h}-1,2^{h}+1)$ with $(tx,t)\in s_r(U)\cap c_{[2^{h}, 2^{h+1})}X$. By the same calculation as before, we have that $\|x-y\| < \frac{1}{2^{h+2}}$ for all $(sy,s)\in s_{r}(U)$. Therefore the projection of $s_r(U)$ onto $X$ has diameter $<\frac{1}{2^{h+1}}$.  So there exists a $z\in Z_{h}$ such that $y\in B(z,\frac{1}{2^{h}})$ for all $(sy,s)\in s_{r}(U)$. We have that $i(sy,s) = i_h(sy,s)=(\varphi_{h}(y),s)$ or $i(sy,s) = i_{h-1}(sy,s) = (\varphi_{h-1}(y),s)$ depending on the height $s$. Let $z'= \phi_{h}(z)$. It is now clear that $i(sy,s)\in \langle gst([B(z', \frac{1}{2^{h-1}})], 2^h) \rangle$ where $[B(z', \frac{1}{2^{h-1}})]$ denotes the vertex in $|\mathcal{U}_{h-1}|$.
\end{proof}

\begin{notation} From now on we identify points $z\in Z_{h}$ with the vertices $[B(z,2^h)]$ in $|\mathcal{U}_{h}|$.
\end{notation}
\begin{definition}\label{R}Let $R: M\rightarrow cX$ be defined piecewise 
\begin{align*}
R_{h}:|\mathcal{U}_{h}| \times (2^{h},2^{h+1}] &\rightarrow X \times (2^{h},2^{h+1}]
\end{align*}
by sending points $(z,t), z\in vert |\mathcal{U}_{h}|$ to $(tz,t)$ for $z\in Z_{h} \subset X$. For $x\in |\mathcal{U}_{h}| \times \{t\}$, $x$ in the interior of some simplex $\sigma\subset |\mathcal{U}_{h}|$, we send $(x,t)$ to a point $(y,t)$ where $y\in \cap_{z\in vert(\sigma)}B(z,2^h)$. 
\end{definition}

The resulting map $R$  is again not continuous, but it is coarse: for this, it suffices to prove the distance estimate for points lying in a common geometric cylinder, since $M$ is a path metric space induced from the spherical metric on each $|\mathcal{U}_{h}| \times \{t\}$ (for a discussion on why the spherical metric is necessary, see Subsection \ref{sphs} in the Appendix). Let $\sigma \times I \subset |\mathcal{U}_{h}| \times [2^{h}+1,2^{h+1}]$ be a geometric cylinder. Vertices are sent to points whose projection to $X$ are distance $<\frac{1}{2^{h-1}}$ from each other. Interior points are sent to points whose projection are distance $< \frac{1}{2^{h}}$ from the image of any vertex. Hence, the set $p'R(\sigma \times I)$ has diameter $<\frac{1}{2^{h-1}}$. Suppose that $(tx,t),(sy,s)$ are two points in $R(\sigma \times I)$. We have
\begin{align*}
\|tx-sy\| +|t-s| \leq \|tx-sx\| +\|sx-sy\| +|t-s| =|t-s|\|x\| +s \|x-y\| +|t-s|\\
 < \|x\| + 2^{h+1} \frac{1}{2^{h-1}} + 1 < \frac{21}{4} 
\end{align*} 
Now let $\sigma \times I = [z_0,\dots,z_m, \phi_{h}(z_0), \dots, \phi_{h}(z_m)] \subset (|\mathcal{U}_{h}| \times [2^{h},2^{h}+1]) \cup_{\phi_h} |\mathcal{U}_{h-1}|$ be a geometric cylinder.  All points in $R(\sigma \times I)$, after projecting to $X$, are distance  $\leq \frac{1}{2^{h-1}}$ from some $\phi_{h}(z_i)$. The vertices  $\phi_{h}(z_i)$ are sent to distance $< \frac{1}{2^{h-2}}$ from one another. Hence, the set $p'R(\sigma \times I)$ has diameter $<\frac{1}{2^{h-2}} + \frac{2}{2^{h-1}} = \frac{1}{2^{h-3}}$. Suppose that $(tx,t),(sy,s)$ are two points in $R(\sigma)$. We have
 \begin{align*}
\|tx-sy\| +|t-s| <  \|x\| + 2^{h+1} \frac{1}{2^{h-3}} + 1 < \frac{69}{4} 
\end{align*} 
Hence, $R(\sigma \times I)$ has diameter $< \frac{69}{4}$. This shows that $R$ is controlled. It is height preserving, and therefore proper.\\

Let $(tx,t) \in X \times (2^h,2^{h+1}]$. The composition $R \circ i$ sends $(tx,t)$ to a point $(tz,t)$ which is at most distance $\frac{1}{2^{h-1}} \cdot 2^{h+1} = 4$ away.  It is therefore close to the identity.\\

The composition $i\circ R: M\rightarrow M$ sends vertices $(z,t)\in M$ to something in the star $\langle st(z)\rangle \times \{t\}$, and it sends interiors of simplices $(x,t)$, $x\in \tau$ into the intersection of open stars of its vertices $\cap_{z \in vert \tau} \langle st(z) \rangle \times \{t\}$. Therefore, $i\circ R$ distance less the maximum diameter of a simplex at a single height in $M$, which is $1$, away from the identity.  \\

\subsection{Isomorphism to Lipschitz end homotopy groups}

Let us now define the Lipschitz end homotopy groups.

\begin{definition} (Coarse-Lipschitz homotopy) Let $X$ be a proper metric space with Lipschitz projection map $p: X\rightarrow [1,\infty)$ and a proper Lipschitz base ray $\omega: [1,\infty)\rightarrow X$. We define the \textit{coarse-Lipschitz cylinder}
\begin{align*}
I_{p}(X) := \{(x,t)\in X\times [0,\infty) \,|\, t\leq p(x)\}
\end{align*}
with inclusions $i_0: X\rightarrow I_{p}(X)$ and $i_{1}: X\rightarrow I_{p}(X)$ defined by $i_0(x) = (x,0)$ and $i_1(x) = (x,p(x))$.\\

Let $X$ and $Y$ be proper metric spaces with projections. A \textit{coarse-Lipschitz homotopy} between two proper Lipschitz maps $f_0, f_1: X\rightarrow Y$ is a proper Lipschitz map $H: I_{p}(X)\rightarrow Y$ such that $f_0 = H\circ i_0$ and $f_1 = H\circ i_1$.
\end{definition}

As usual, for the cone $cX$ of a compact metric space $X$, we can assume that $p$ is the height projection. 

\begin{definition} (Lipschitz end homotopy groups) Let $X$ be a proper metric space with proper Lipschitz base ray $\omega: [1,\infty)\rightarrow X$. For $n\geq 1$ we define the \textit{$n$-th Lipschitz end homotopy group} $\pi_n^{L,e}(X,\omega)$ to be the set of relative coarse-Lipschitz homotopy classes of proper Lipschitz maps 
\begin{align*}
f: (c[0,1]^n,c\partial [0,1]^n)\rightarrow (X,\omega[1,\infty))
\end{align*}
such that $f|_{c\partial[0,1]^n} = \omega\circ p$ where $p: c[0,1]^n \rightarrow [1,\infty); (hx,h)\mapsto h$ denotes the height variable of the cone. \\

 We define the \textit{$0$-th coarse-Lipschitz homotopy set} $\pi_0^{L,e}(X)$ to be the set of coarse-Lipschitz homotopy classes of proper Lipschitz maps from $[1,\infty)$ to $X$.
\end{definition}

Many of the properties of coarse homotopy groups are directly applicable to Lipschitz end homotopy groups: the group multiplication is the same formula, they are functorial under proper Lipschitz maps with coarse-Lipschitz homotopic maps inducing the same morphism. There is a long exact sequence for pairs $(X,A)$ with proper Lipschitz map $k: A \rightarrow X$ and proper Lipschitz base ray $\omega:[1,\infty) \rightarrow A$. The results of Proposition \ref{baseddefined} apply verbatim if one exchanges "coarse map" with "proper Lipschitz map" and "coarse homotopy" with "coarse-Lipschitz homotopy", since all the other homotopies used in the proofs are coarse-Lipschitz. The details are left as an exercise to the reader. We will not be needing many of these properties as the coarse-Lipschitz homotopy groups themselves are only an intermediary construction in our proofs. \\

Consider now the base ray $\omega: [1,\infty)\rightarrow cX$ defined as $\omega(t)=(tx_0,t) = (0,t)$. We call this the \textit{standard parametrisation of the base ray $c\{x_0\}$}. By construction, $x_0\in Z_{h}$ for all $h\in \mathbb{N}_{0}$ so there is a copy of $\{x_0\} \times [1,\infty) $ in $M$. We also denote this base ray by $\omega$. In both cases, $\omega: [1,\infty)\rightarrow cX,M$ is an isometric embedding.  

\begin{theorem} \label{isoomega}There is an isomorphism  
\begin{align*}
\lambda:  \pi_{n}^{L,e}(M,\omega)\rightarrow \pi_{n}^c(cX,\omega) 
\end{align*} 
with inverse $\theta:  \pi_{n}^c(cX,\omega)\rightarrow  \pi_{n}^{L,e}(M,\omega)$. 
\end{theorem}

The main goal of this subsection is to prove this theorem. Let us define all the relevant maps. There is a homomorphism $\lambda: \pi_{n}^{L,e}(M,\omega) \rightarrow \pi_{n}^c(cX,\omega)$ defined in the following way: a class $[f]\in  \pi_{n}^{L,e}(M,\omega)$ is represented by a proper Lipschitz map, which is in particular coarse, and a coarse-Lipschitz homotopy between maps is also a coarse homotopy. So $[f]\in \pi_{n}^{c}(M,\omega)$. The map $R: M\rightarrow cX$ is coarse, and satisfies $R(\omega)=\omega$, so it induces a morphism $R: \pi_n^c(M,\omega)\rightarrow \pi_{n}^c(cX,\omega)$. We define $\lambda$ as the composition 

\begin{equation*}
	\begin{tikzcd}[row sep=0.5em, column sep=0.9em]
		\lambda: \pi^{L,e}_{n}(M,\omega)  \arrow[r] & \pi_{n}^c(M,\omega)  \arrow[r] & \pi_{n}^c(cX,\omega) \\
		{[g]} \arrow[r,|->] & {[g]} \arrow[r,|->]& {[R\circ g]} 
	\end{tikzcd}
\end{equation*}

Let us now define $\theta: \pi_n ^{c}(cX,\omega)\rightarrow  \pi^{L,e}_{n}(M,\omega)$. Let $[f]\in \pi_n^c(cX,\omega)$. Fix a simplicial structure on $[0,1]^n$. Recall that $c[0,1]^n$ is uniformly bounded, ie. there are only finitely many strong similarity types of simplices; there is a positive lower bound on their widths, and their diameters are bounded between $a,b$ for $0<a< b<\infty$ (with respect to the metric induced by any norm on $\mathbb{R}^{n+1}$). Since $f$ is controlled, there exists a $S>0$ such that $diam f (|st(v)|) < S$ for all closed stars $|st(v)|$ of vertices $v \in c[0,1]^n$. By Lemma \ref{weirdness} there is a shrinking map $s_{r}$ such that for all vertices $v$, $is_{r}  f (|st(v)|)\subset \langle gst(u_v) \rangle$ for some $u_v \in M$. We say that $is_{r}f$ satisfies the geometric star condition. \\

We can define a Lipschitz map $\widetilde{is_{r}f}: c[0,1]^n \rightarrow M$ by letting $\widetilde{is_{r}f}(v)=u_v$ on vertices and extending skeleta by skeleta as follows. Assume that $\widetilde{is_{r}f}$ is defined on the $(k-1)$-skeleton. Let $\Delta^{k}=[v_0,\dots,v_k]$ be a $k$-simplex. By the geometric star condition $is_{r}f(\Delta^{k})\subset \cap_{i=0}^{k} \langle gst(u_{v_{i}}) \rangle$. From Lemma \ref{geomsexist}, we have that $[u_{v_{0}}, \dots, u_{v_{k}}]$ is a geometric simplex. By abuse of notation, denote the core of this geometric simplex by $\mathscr{C}_{\Delta^{k}}:= \mathscr{C}_{[u_{v_{0}}, \dots, u_{v_{k}}]}$. Assume that $\widetilde{is_{r}f}|_{\partial \Delta^{k}}(\partial \Delta^k) \subset \mathscr{C}_{\Delta^{k}}$. By Lemma \ref{usefulprops}, there is a deformation retract $H^{\mathscr{C}_{\Delta^{k}}}$ of $\mathscr{C}_{\Delta^{k}}$ to a point. Therefore we can extend $\widetilde{is_{r}f}$ to $\Delta^k$. Let $\Delta^{k+1} = [w_0,\dots,w_{k+1}]$ be a $(k+1)$-simplex with $\Delta^k \subset \Delta^{k+1}$ as a face. Since $[u_{v_{0}}, \dots, u_{v_{k}}]\subset [u_{w_0},\dots,u_{w_{k+1}}]$ we have that $\widetilde{is_{r}f}(\Delta^k)\subset \mathscr{C}_{\Delta^k} \subset \mathscr{C}_{\Delta^{k+1}}$. Therefore the induction step is met, and there is a globally defined map $\widetilde{is_{r}f}: c[0,1]^n \rightarrow M$, which we call the geometric simplicial approximation of $is_{r}f$. \\

To show that $\widetilde{is_{r}f}$ is Lipschitz, we give an explicit description of the extension to the $k$-skeleton. See Figure \ref{fig:lipschitzfilling}. Assume that for every $k-1$-simplex $\Delta^{k-1}$, $\widetilde{is_{r}f}_{|\Delta^{k-1}}: \Delta^{k-1} \rightarrow \mathscr{C}_{\Delta^{k-1}}$ is Lipschitz with respect to the norm on $\Delta^{k-1}\subset \mathbb{R}^{n+1}$ and the path metric on $\mathscr{C}_{\Delta^{k-1}}$, with a uniform Lipschitz constant $L_{k-1}$ depending only on the dimension $k-1$ and constants related to the uniform boundedness of $c[0,1]^n$. Let $\Delta^{k} = [v_0,\dots,v_k]$ be a $k$-simplex. If $\mathscr{C}_{\Delta^{k}}$ is type $1$, geometric and ordinary simplicial aproximation coincide. $\mathscr{C}_{\Delta^{k}}$ is a simplex of dimension at most $k$, equipped with the spherical metric. The composition of simplicial maps
\begin{align*}
\widetilde{is_{r}f}: \Delta^{k} \rightarrow \Delta^k_{std} \rightarrow \mathscr{C}_{\Delta^{k}}
\end{align*}
where $\Delta^k_{std}$ denotes the standard $k$-simplex, has Lipschitz constant depending only on $k$ and the lower bound $w$ on the width of all simplices in $c[0,1]^n$. \\

 If $\mathscr{C}_{\Delta^{k}} = (\sigma \times I) \cup sp\{q\sigma, \tau\}$ is type $2$, recall that $H^{\mathscr{C}_{\Delta^{k}}}$ is the composition $H\ast (H'\circ H_1 \times \id)$ where $H$ deformation retracts to $sp\{q\sigma, \tau\}$, and $H'$ to the barycentre of $sp\{q\sigma, \tau\}=sp\{qu_{v_{0}}, \dots, qu_{v_{k}}\}$. \\

\begin{figure}
    \centering
  \includegraphics[width=0.8\linewidth]{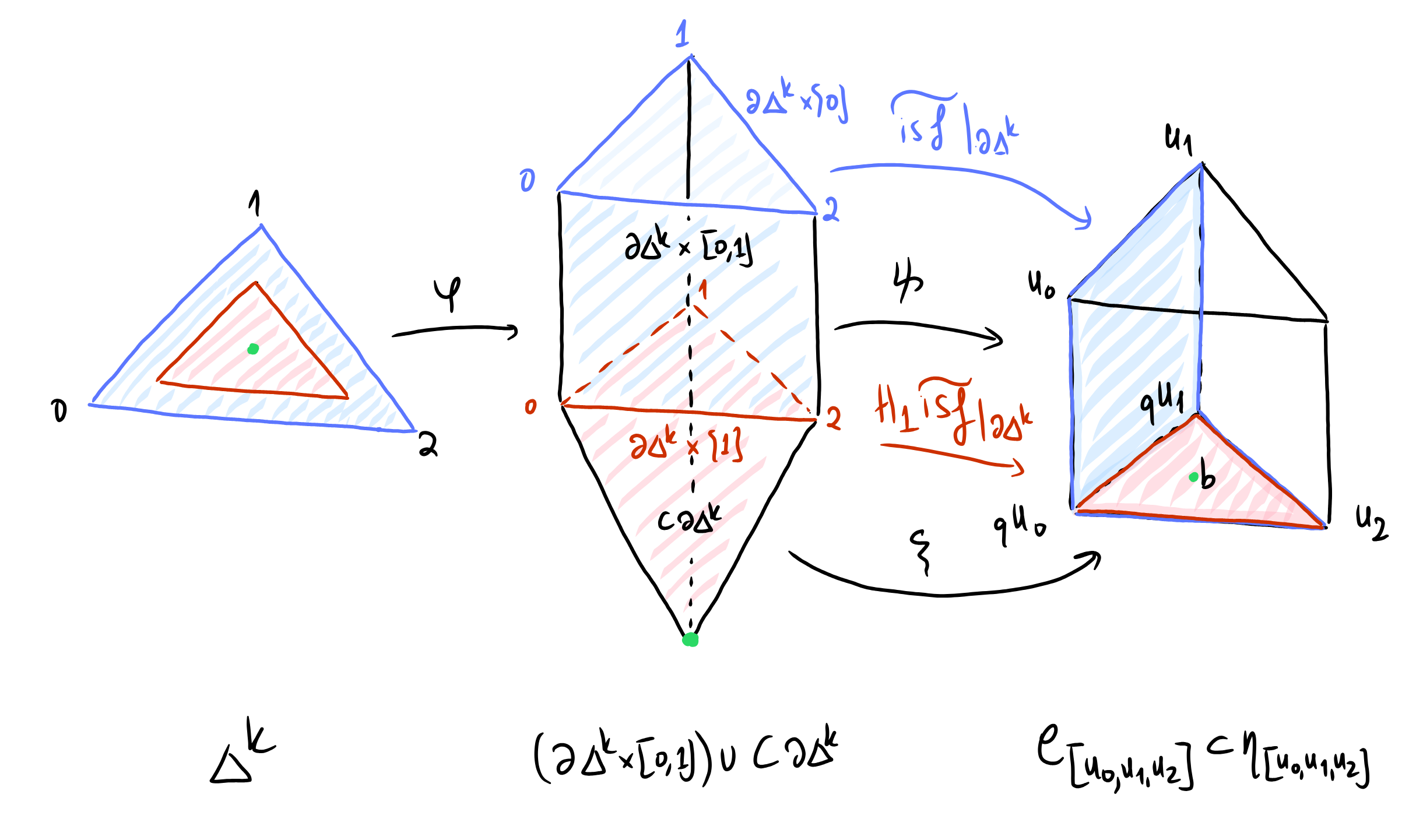}
  \caption{The Lipschitz extension of $\widetilde{is_{r}f}$ to $\Delta^{k}$ for $\Delta^{k}=[0,1,2]$ given the vertex map $\widetilde{is_{r}f}(j)=u_j$.}
  \label{fig:lipschitzfilling}
\end{figure}

Every point in $\Delta^{k}$ can be written in co-ordinates as $(x,t):=(1-t)b_{\Delta^{k}}+t x$ where $b_{\Delta^{k}}$ is the barycentre of $\Delta^{k}$ and $x\in \partial \Delta^{k}$.  By a linear translation, we can assume that $b_{\Delta^{k}} = 0$ so that $d((x,t),(y,s)) = \|tx-sy\|$.  We can write $\Delta^k = (\Delta^k\setminus  (\Delta^{k})^{'}) \cup (\Delta^{k})^{'}$ where $(\Delta^{k})':=\{(x,t)\in \Delta^{k}\,|\, t\in [0,\frac{1}{2}]\}$ is a smaller copy of $\Delta^k$ in the interior. Consider the metric space  $(\partial{\Delta^k} \times [0,1]) \cup C\partial\Delta^{k}$, where $\partial{\Delta^k} \times [0,1]$ is equipped with the product metric and $C\partial\Delta^{k}:= (\partial \Delta^{k} \times [0,1]) / (\partial \Delta^{k} \times \{0\})$ denotes the topological cone with metric 
\begin{align*}
d_{C\partial\Delta^{k}}((x,t),(y,s)) = \|tx-sy\| + |t-s|
\end{align*}
The subspaces $\partial \Delta^{k} \times \{1\}\subset \partial{\Delta^k} \times [0,1]$ and $\partial \Delta^{k} \times \{1\} \subset C\partial\Delta^{k}$ are identified.  There is a PL homeomorphism $\varphi: \Delta^k=(\Delta^k\setminus  (\Delta^{k})^{'}) \cup (\Delta^{k})^{'}\rightarrow (\partial{\Delta^k} \times [0,1]) \cup C\partial\Delta^{k}$. \\

$\varphi$ is Lipschitz: for points $(x,t),(y,s)\in (\Delta^k\setminus  (\Delta^{k})^{'})$ we have 
\begin{align*}
d_{\partial{\Delta^k} \times [0,1]}((x,t),(y,s)) = \|x-y\| + |t-s| \leq 2(\|tx-sy\|+|t-s|\|y\|) + |t-s|\\
 < \delta_1(\|tx-sy\|+|t-s|)
\end{align*}
where $\delta_1>0$ is a constant depending only on the diameter of $\Delta^{k}$. For points in $(\Delta^{k})^{'}$, since $\varphi_{|(\Delta^{k})^{'}}$ is PL, it is $\delta_{2}$-Lipschitz for some $\delta_{2}>0$ . Because $\Delta^{k}$ is geodesic, the Lipschitz constant for $\varphi$ is $\max \{\delta_1,\delta_2\}$. Since $c[0,1]^n$ is uniformly bounded, we can choose a $\delta>0$ so that $\varphi$ is $\delta$-Lipschitz for all $k$-simplices (maximise over the diameters for $\delta_{1}$, and maximise over the smallest-diameter elements in the finitely many strong similarity types for $\delta_{2}$).\\

The deformation retract $H:\mathscr{C}_{\Delta^{k}} \times [0,1]\rightarrow \mathscr{C}_{\Delta^{k}}$ is $1$-Lipschitz and $H_{1}(\mathscr{C}_{\Delta^{k}})\subset sp\{q\sigma,\tau\}$. We define
\begin{align*}
\psi &:  \partial{\Delta^k} \times [0,1] \rightarrow \mathscr{C}_{\Delta^{k}} \\
\psi &= H\circ (\widetilde{is_{r}f}_{|\partial \Delta^{k}} \times \id_{[0,1]})
\end{align*}
The path metric on $\partial \Delta^{k}$ is bi-Lipschitz equivalent to the norm-induced metric when considering $\partial \Delta^{k}$ as a subspace of $\Delta^{k}$. Since there are only finitely many strong similarity types of $\Delta^{k}$ and the inclusion $\cup_{\Delta^{k-1}\subset \partial \Delta^{k}}\mathscr{C}_{\Delta^{k-1}}\subset \mathscr{C}_{\Delta^{k}}$ is $1$-Lipschitz, the Lipschitz constant of $\widetilde{is_{r}f}_{|\partial \Delta^{k}}$, and therefore $\psi$, depends only on $L_{k-1}$. We define 
\begin{align*}
\xi &: C\partial \Delta^{k} \rightarrow sp\{q\sigma,\tau\}\\
\xi &(x,t) = (1-t)b_{sp\{q\sigma,\tau\}} + tH_{1}(\widetilde{is_{r}f}_{|\partial \Delta^{k}}(x)) 
\end{align*}
where $b_{sp\{q\sigma,\tau\}}$ is the barycentre of $sp\{q\sigma,\tau\}$. We let the extension of $\widetilde{is_{r}f}$ to $\Delta^{k}$ be defined as $(\psi\cup \xi)\circ \varphi$. It is obviously Lipschitz when restricted to $\Delta^k\setminus  (\Delta^{k})^{'}$. To show that $\xi$ is Lipschitz, let $c,c'$ be the Lipschitz constants obtained from changing between the standard ($L^1$) and spherical metrics on $sp\{q\sigma,\tau\}$ ie. $\frac{1}{c'} \|w-z\| \leq d^s(w,z)\leq c \|w-z\| $ where $d^{s}$ denotes the spherical metric. Since $sp\{q\sigma,\tau\}$ is at most a $k$-dimensional simplex, $c,c'$ depend only on $k$. Let $L'_{k}$ be the Lipschitz constant of $\widetilde{is_{r}f}_{|\partial \Delta^{k}}$.
\begin{align*}
d^{s}(\xi(x,t),\xi (y,s))&=d^{s}((1-t)b_{sp\{q\sigma,\tau\}} + tH_{1}(\widetilde{is_{r}f}(x)), (1-s)b_{sp\{q\sigma,\tau\}} + sH_{1}(\widetilde{is_{r}f}(y)))\\
& \leq  c \| t(H_{1}(\widetilde{is_{r}f}(x))-b_{sp\{q\sigma,\tau\}}) -s(H_{1}(\widetilde{is_{r}f}(y))-b_{sp\{q\sigma,\tau\}})\| \\
& \leq c |t-s|\|H_{1}(\widetilde{is_{r}f}(x)) -b_{sp\{q\sigma,\tau\}}\| + cs\|H_{1}(\widetilde{is_{r}f}(x)) -H_{1}(\widetilde{is_{r}f}(y))\| \\
&\leq c |t-s| +cc' s d^{s}(H_{1}\widetilde{is_{r}f}(x),H_{1}\widetilde{is_{r}f}(y)) \\
& \leq c |t-s| +cc' sd_{\mathscr{C}_{\Delta^{k}}}(\widetilde{is_{r}f}(x),\widetilde{is_{r}f}(y)) \\
&\leq c |t-s| +cc's L'_k\|x-y\|\\
& \leq c |t-s| +cc'L'_k(\|tx-sy\| + |t-s|)\\
& \leq 2\max\{c,cc'L'_k\} d_{C\partial\Delta^{k}}((x,t),(y,s))
\end{align*}
Therefore $\xi$ is Lipschitz, with Lipschitz constant depending only on $L_{k-1}$ and the dimension $k$.\\

This shows that we can choose a constant $L_{k}$ such that $\widetilde{is_{r}f}_{|\Delta^{k}}: \Delta^{k} \rightarrow \mathscr{C}_{\Delta^k}$ is $L_{k}$-Lipschitz for all $k$-simplices. Since $c[0,1]^n$ is finite-dimensional, by induction, there exists a $L>0$ such that $f_{|\Delta}: \Delta \rightarrow \mathscr{C}_{\Delta}$ is $L$-Lipschitz for every simplex $\Delta\subset c[0,1]^n$. Moreover, the inclusion $\mathscr{C}_{\Delta} \subset M$ is $1$-Lipschitz. Therefore $\widetilde{is_{r}f}$ is uniformly locally Lipschitz (see Defintion \ref{uniformll}) with respect to the set of simplices in $c[0,1]^n$.\\

Consider now $c[0,1]^n$ as a path metric space. This is equivalent to $c[0,1]^n$ with the metric induced from the Euclidean norm, since $c[0,1]^n$ is convex. As a path metric space, it is geodesic, and each geodesic segment contained entirely in one simplex has length equal to the norm-induced distance. These facts are sufficient to show that $\widetilde{is_{r}f}: c[0,1]^n \rightarrow M$ is globally Lipschitz. It is proper because for $a\in \mathbb{N}$,
\begin{align*}
(\widetilde{is_{r}f})^{-1}M_{[1,a]} \subset  (is_{r}f)^{-1} M_{[1,a+1]}\subset f^{-1}(c_{[1,2^r(a+1)]}X)
\end{align*}
which is a bounded set. \\

Observe that since $f(j\partial[0,1]^n \times \{j\}) = \{(jx_0,j)\}$ for all $j\in \mathbb{N}$ and $s_{r}(jx_0,j) = (\frac{j}{2^r}x_0,\frac{j}{2^r})$ for an appropriately large $r$, we have $i\circ s_{r} \circ f(|st(v)|)\subset \langle gst(x_0,\omega_{r}(j))\rangle $ for any vertex $v\in j\partial[0,1]^n \times \{j\}$, where $\omega_r(j)$ is an integer with distance $\leq \frac{1}{2}$ to $\frac{j}{2^r}$. We can extend $\omega_{r}$ linearly to a map $\omega_{r}: [1,\infty) \rightarrow \{x_0\} \times [1,\infty)$. We can choose the simplicial approximation such that $\widetilde{is_{r}f}: (c[0,1]^n, c\partial[0,1]^n)\rightarrow (M,\omega_r)$, where $\widetilde{is_{r}f}_{|c\partial[0,1]^n} = \omega_r\circ p$. Note that the linear homotopy $H^{r,0}:  c(\{\ast\}\times [0,1]) \rightarrow c\{x_0\}$ between $\omega_r$ and $\omega$ is proper and Lipschitz. Let $b^{H^{r,0}}_{\omega_r,\omega}:  \pi^{L,e}_{n}(M,\omega_r)\rightarrow \pi^{L,e}_{n}(M,\omega)$ be the change of base ray homomorphism. We define 

\begin{align*}
\theta: \pi_n ^{c}(cX,\omega)&\rightarrow \pi^{L,e}_{n}(M,\omega)\\
[f]&\mapsto b^{H^{r,0}}_{\omega_r,\omega}[\widetilde{is_{r}f}]
\end{align*}

There are quite a number of things to check to see this actually works, which are contained in the following three technical propositions and lemmas (\ref{baserayindep}, \ref{geosa} and \ref{indepwelldefined}).

\begin{lemma} \label{baserayindep} (Reparametrisations of $\omega$) Let $\gamma: [1,\infty)\rightarrow \{x_0\}\times [1,\infty)$ and $\gamma': [1,\infty)\rightarrow \{x_0\} \times [1,\infty)$ be coarse (resp. proper Lipschitz) maps which reparametrise $\omega$. 
\begin{enumerate}
\item There is a linear coarse (resp. coarse-Lipschitz) homotopy $G^{\gamma,\gamma'}: c([0,1])\rightarrow \{x_0\} \times [1,\infty)$ between $\gamma$ and $\gamma'$. 
\item Let $G:c([0,1])\rightarrow  \{x_0\} \times [1,\infty)$ be any coarse (resp. coarse-Lipschitz) homotopy between $\gamma$ and $\gamma'$. Let $\eta$ and $\eta'$ denote the linear coarse (resp. coarse-Lipschitz) homotopies between $\gamma$ and $\omega$, $\gamma'$ and $\omega$ respectively. There is a coarse (resp. coarse-Lipschitz) homotopy $F: c([0,1]^2)\rightarrow \{x_0\} \times [1,\infty)$ between $\eta$ and $\eta'$, which restricts to $G$ and $\id_{\omega}$ on  $c(\{0\} \times [0,1])$ and $c(\{1\} \times [0,1])$ respectively. 
\end{enumerate} 
\end{lemma}

\begin{proof}
In this proof, the height variable $h$ of $c[0,1]^n$ is always ordered last. We can identify $\{x_0\} \times [1,\infty)$ with its height, since $\omega$ is an isometric embedding. 
\begin{enumerate}
\item This is the same computation as Lemma \ref{shrinkindep} by letting $G(hs,h)=\eta(h,s) = s\gamma'(h) + (1-s)\gamma(h)$. The same distance estimates work for the Lipschitz case, just without the additive constant.  
\item We define the homotopy
\begin{align*}
F: c([0,1]^2)&\longrightarrow \{x_0\} \times [1,\infty)\\
(ht,hs,h)&\longmapsto (1-t) G(hs,h) + th
\end{align*}
For $s=0$ we have the linear homotopy between $G(0,h) = \gamma$ and $\omega$; for $s=1$ we have a homotopy between $G(h,h) = \gamma'$ and $\omega$; for $t=0$ we have $G(hs,h)$; for $t=1$ we have the constant base ray $\omega$. \\

To show that this is coarse, let $C$ be the large-scale Lipschitz constant for $G$. For two points $(ht,hs,h), (h't',h's',h')$ we have
\begin{align*}
|F(ht,hs,h)-&F(h't',h's',h')| = |(1-t)G(hs,h) -(1-t')G(h's',h')+ th - t'h'| \\
&\leq (1-t)|G(hs,h) -G(h's',h')| + |t-t'||G(h's',h')| + |ht-h't'|\\
&\leq C(|hs-h's'| + |h-h'|)+C + |t-t'||G(h's',h')| + |ht-h't'|
\end{align*}
Consider the point $(0,1)$, which is distance $sh' + (h'-1) < 2h'$ away from $(h's',h')$. We have
\begin{align*}
|G(h's',h')| &\leq |G(h's',h')-G(0,1)|+ |G(0,1)| < 2Ch' + C + \gamma(1) \\
&= (3C+ \gamma(1))h'\\
|t-t'| |G(h's',h')| &< |th'-t'h'|(3C+ \gamma(1)) \leq (3C+ \gamma(1))( |ht-h't'| + |h-h'|)
\end{align*}
These estimates show that $F$ is controlled. To see that $F$ is proper, observe that since $G$ is proper, for $h'\in [1,\infty)$ there exists a $l>h'$ such that $G(hs,h) > h'$ for all $h>l$. Therefore $F^{-1}([0,h'])\subset c_{[0,l]}([0,1]^2)$. We conclude that $G$ is coarse.\\

The same distance estimates work for the Lipschitz case, just without the additive constant.  
\end{enumerate}
\end{proof}

The next proposition (\ref{geosa}) proves some basic properties about geometric simplicial approximation.  Many of these statements and proofs are analogous to the case of ordinary simplicial approximation. We first observe that the construction of $\widetilde{is_{r}f}$ works for any simplicial complex $K$ and any map $f: K\rightarrow M$ which satisfies the geometric star condition:

\begin{definition} Let $K$ be a simplicial complex, $M$ the inverse mapping telescope associated to the compact metric space $X$. We call a map $g: K\rightarrow M$ \textit{geometrically simplicial} if for every simplex $\sigma \subset K$, $\cup_{v \in vert \sigma} g(v)$ spans a geometric simplex, and $g(\sigma) \subset \mathscr{C}_{ [\cup_{v \in vert \sigma} g(v)]}$. \\

 Let $f: K\rightarrow M$ be a (possibly non-continuous) map which satisfies the \textit{geometric star condition}: for all vertices $v\in K$, there exists a vertex $u_{v}\in M$ such that $f|st(v)|\subset \langle gst(u_{v})\rangle $. The \textit{geometric simplicial approximation} of $f$ is a geometrically simplicial map $\tilde{f}$ which is defined on vertices as $\tilde{f}(v)=u_{v}$ and extends to the interior of simplices by Lipschitz deformation retracts of geometric cores.\\
\end{definition}

\begin{definition} \label{uniformll} Let $Z$ be a proper metric space with a collection of subsets $\mathcal{L}=\{L_{i}\}_{i\in \mathbb{N}_{0}} \subset Z$ such that $Z= \cup_{i\in \mathbb{N}_{0}} L_{i}$. Let $Y$ be a proper metric space. A map $f: Z\rightarrow Y$ is called \textit{locally Lipschitz} with respect to $\mathcal{L}$ if $f_{|L_{i}}$ is Lipschitz for all $i\in \mathbb{N}_{0}$. We say that $f$ is \textit{uniformly locally Lipschitz} if we can choose the same Lipschitz constant for all $i\in \mathbb{N}_{0}$.
\end{definition}

Note that we have not specified a metric on $L_{i}$. Usually this is taken as the subspace metric from $Z$.  \\

If $K$ is a simplicial complex with path metric induced from an embedding into a normed space, $\mathcal{L}$ the collection of simplices in $K$ with the subspace metric, then we have shown that the geometric simplicial approximation $\tilde{f}: K\rightarrow M$ is locally Lipschitz with respect to $\mathcal{L}$. If $K$ is additionally uniformly bounded and finite-dimensional, then $\tilde{f}$ is globally Lipschitz. \\

To prove the next proposition, we require the definition of the standard metric on a simplicial complex $K$. Recall that an abstract simplicial complex $K$ consists of a vertex set $V$ and a collection $\mathcal{S}$ of non-empty, finite subsets of $V$, which is closed under taking subsets.

\begin{definition} \label{metstd}(Standard metric on geometric realisation) Let $K= (V,\mathcal{S})$ be an abstract simplicial complex. On the product space $[0,1]^V$ of functions $\alpha: V\rightarrow [0,1]$ we consider the subspace
\begin{align*}
|K| = \{\alpha\,|\, \{v\,|\,\alpha(v)>0\} \in \mathcal{S}, \sum_{v\in V} \alpha(v) = 1\}
\end{align*}
The standard metric on $|K|$, denoted by $d_{std}$, is defined as 
\begin{align*}
d_{std}(\sum_{v\in V}\alpha(v)\delta_v,\sum_{v\in V}\beta(v) \delta_v) = \sum_{v \in V} |\alpha(v)-\beta(v)|
\end{align*}
\end{definition}

\begin{prop} (Properties of geometric simplicial approximation) \label{geosa} Let $K$ be a locally finite, finite-dimensional simplicial complex, $L$ a subcomplex of $K$. Assume each simplex of $K$ is equipped with a norm-induced metric such that $K$ is uniformly bounded and give $K$ the induced path metric. Let $f:K\rightarrow M$ be a (possibly non-continuous) map which satisfies the geometric star condition and denote by $\tilde{f}:K\rightarrow M$ a geometric simplicial approximation of $f$. 
\begin{enumerate}
\item Let $\tilde{f}':K\rightarrow M$ be another choice of geometric simplicial approximation. $\tilde{f}'$ is homotopic to $\tilde{f}$ by a Lipschitz homotopy. If $\tilde{f}'_{|L} = \tilde{f}_{|L}$ then this homotopy fixes the subcomplex $L$. 
\item Let $bs(K)$ be the barycentric subdivision of $K$. The geometric simplicial approximation of $f$ with respect to $bs(K)$, denoted $bs\tilde{f}$, exists and is Lipschitz homotopic to $\tilde{f}$. 
\item If $f$ is a continuous (resp. Lipschitz) map, then $f$ is homotopic to $\tilde{f}$ by a continuous (resp. locally Lipschitz) homotopy. If $f_{|L} = \tilde{f}_{|L}$ then this homotopy fixes the subcomplex $L$. 
\item Let $\overline{K}$ be another locally finite, finite-dimensional simplicial complex. Assume each simplex of $\overline{K}$ is equipped with a norm-induced metric such that $\overline{K}$ is uniformly bounded. Give $\overline{K}$ the path metric. Assume that $L$ is also a subcomplex of $\overline{K}$ and there is a bi-Lipschitz homeomorphism $\varphi: (K,L)\rightarrow (\overline{K},L)$ such that $\varphi_{|L}=\id_{L}$. Suppose that that $f\circ \varphi^{-1}$ satisfies the geometric star condition with respect to $\overline{K}$. Let $\tilde{f}': \overline{K}\rightarrow M$ denote its geometric simplicial approximation. $\tilde{f}$ is Lipschitz homotopic to $\tilde{f}'\circ \varphi$. If $\tilde{f}_{|L} = (\tilde{f}'\circ \varphi)_{|L}$ then this homotopy fixes the subcomplex $L$.  

\end{enumerate}
Additionally, if $f$ is proper, then all maps in the above points as well as any homotopy between them, are proper. 
\end{prop}

\begin{remark} Let $K= c(i(X))$ be the cone of the finite simplicial complex $(X,d_{std})$, piecewise linearly embedded into $\mathbb{R}^n$ via $i: X\rightarrow \mathbb{R}^n$. Consider another choice of bi-Lipschitz triangulation for $X$, ie. a bi-Lipschitz homeomorphism $\psi: (X,d_{std})\rightarrow (\overline{X},d_{std})$, and a PL embedding $j: \overline{X}\rightarrow \mathbb{R}^m$. For $\overline{K}= c(j(\overline{X}))$, we have that $c(j\circ \psi \circ i^{-1}): K \rightarrow \overline{K}$ is a bi-Lipschitz homeomorphism. Additionally, the path metrics on $K,\overline{K}$ are bi-Lipschitz equivalent to the metric induced from the ambient norms. 
\end{remark}

\begin{proof}
\begin{enumerate}

\item Suppose we have $f:K\rightarrow M$ and two simplicial approximations $\tilde{f},\tilde{f}'$ of $f$. Give $K \times [0,1]$ the product simplicial structure and metric. Taking products with $[0,1]$ only adds finitely many additional strong similarity classes: $K\times [0,1]$ is uniformly bounded. We identify $K$ with $K \times \{0\} \subset K\times [0,1]$. A simplex in $K\times [0,1]$ has the form $\sigma = [v_0, \dots, v_{j-1}, v_j, v'_j, v'_{j+1} \dots, v'_{m}]$ (or some subset thereof) for some $0\leq j\leq m$, where $\Delta^m= [v_0,\dots, v_m]$ is a simplex in $K$ and $v'_i$ denotes a copy of $v_i$ on $K \times \{1\}$.\\

Let $pr: K\times [0,1] \rightarrow K$ denote the projection to $K$. Consider the map $f\circ pr: K \times [0,1] \rightarrow M$. We show that a simplicial approximation of $f\circ pr$ is a homotopy between $\tilde{f}$ and $\tilde{f}'$. \\

Let $v$ be a vertex in $K \times \{0\}$. We have that 
\begin{align*}
f \circ pr(|st(v,K\times [0,1])|) \subset f(|st(v,K)|) \subset \langle gst(\tilde{f}(v)) \rangle 
\end{align*} 
For a vertex $v' \in K \times \{1\}$, we have
\begin{align*}
f \circ pr(|st(v',K\times [0,1])|) \subset f(|st(v,K)|) \subset \langle gst(\tilde{f}'(v)) \rangle 
\end{align*} 
This shows that $f\circ pr$ has a simplicial approximation, defined on vertices as $\widetilde{f\circ pr}(v)= \tilde{f}(v)$ and $\widetilde{f\circ pr}(v')= \tilde{f}'(v)$. since $K \times [0,1]$ is uniformly bounded, $\widetilde{f\circ pr}$ is locally Lipschitz with respect to the set of simplices, each equipped with the norm-induced metric. The path metric on $K \times [0,1]$ is the same as the product metric: therefore $\widetilde{f\circ pr}$ is globally Lipschitz. The geometric simplicial approximation of a proper map is proper. It is clear that $\widetilde{f\circ pr}$ restricts to $\tilde{f}$ on $K \times \{0\}$ and $\tilde{f}'$ on $K \times \{1\}$. Therefore $\widetilde{f\circ pr}$ is a proper, Lipschitz homotopy between $\tilde{f}$ and $\tilde{f}'.$\\

If $\tilde{f}_{|L} = \tilde{f}'_{|L}$ we can choose $H$ to fix the subcomplex $L$ by inducting only over simplices not contained in $L \times [0,1]$.

\item Let $bs(K)$ be the barycentric subdivision of $K$. Let $\Delta=[v_0,\dots,v_m]$ be a simplex in $K$, $b_{\Delta}$ the barycentre of $\Delta$. Let $u_{v_0},\dots ,u_{v_m}$ be the vertices in $M$ such that $f (|st(v_i)|) \subset \langle gst(u_{v_i})\rangle$ for all $i$, so that $\tilde{f}$ is the Lipschitz extension of $\tilde{f}(v_i) = u_{v_i}$ to the interior of simplices. We can define $bs(\tilde{f})$ by letting $bs(\tilde{f})(v_i) = u_{v_i}$ for all $0\leq i \leq m$ and $bs(\tilde{f})(b_{\Delta}) = u_{v_j}$ for any $j$. This is because $|st(v_i,bs(K))|\subset|st(v_i,K)|$ and $|st(b_{\Delta}, bs(K))| \subset \cap_{i=0}^{m} |st(v_i)|$, so the vertex map $bs(\tilde{f})$ extends to the interior of simplices and gives us a simplicial approximation of $f$ with respect to $bs(K)$. Since $bs(K)$ is uniformly bounded, $bs(\tilde{f})$ is Lipschitz.  \\

Recall that in $bs(K)$, $\Delta$ is divided into $(m+1)!$ many $m$-simplices, corresponding to an ordering $v_{i_{0}}, \dots, v_{i_{m}}$ of the vertices in $\Delta$. Geometrically, this corresponds to the simplex $\delta=[v_{i_{0}}, b_{[v_{i_{0}},v_{i_{1}}]},\dots,b_{\Delta}]$, ie. the vertices are barycentres of an increasing flag of subcomplexes of $\Delta$. The image of $\delta$ under $\tilde{f}$ is contained within $\tilde{f}(\Delta) \subset \mathscr{C}_{\tilde{f}\Delta}$, ie. the geometric core defined by the vertices $\{u_{v_0},\dots, u_{v_{m}}\}$. By construction, the image of $\delta$ under $bs(\tilde{f})$ is also contained within $\mathscr{C}_{\tilde{f}\Delta}$. Let $\mathcal{S}$ be the set of simplices in $K$, $bs(\mathcal{S})$ the set of simplices in $bs(K)$. Denote by $\iota: bs(\mathcal{S})\rightarrow \mathcal{S}$ the map which sends a simplex $\delta$ in $bs(K)$ to the unique simplex $\Delta$ in $K$ which contains it.\\

Give $bs(K) \times [0,1]$ the product simpicial structure. We want to extend the map $H$ which we fix to be $\tilde{f}$ on $bs(K) \times \{0\}$ and $bs(\tilde{f})$ on $bs(K) \times \{1\}$.  Assume that $H$ is defined on the $(k-1)$-skeleton. Let $\sigma$ be a $k$-simplex in $bs(K) \times [0,1]$, not contained in $bs(K) \times \{0,1\}$. Assume that for every face $\tau \subset \partial \sigma$, $H_{|\tau}$ has image in the core of the geometric simplex spanned by the vertices $\cup_{v \in vert \iota(pr \tau)} \tilde{f}(v)$, which we denote (by abusing notation) as $\mathscr{C}_{\tau}$. Note that this automatically holds for faces in $bs(K) \times \{0,1\}$. Since $\iota (pr \tau) \subset \iota (pr \sigma)$, we have that $\cup_{\tau \subset \partial \sigma} \mathscr{C}_{\tau} \subset \mathscr{C}_{\sigma}$.\\ 

As in the construction of the geometric simplicial approximation, $\mathscr{C}_{\sigma}$ is contractible, so we can extend $H$ to the interior of $\sigma$.  There is only one small modification: since $H_{|\partial \sigma}$ is not necessarily simplicial, for $\mathscr{C}_{\sigma}$ of type $1$, we consider the composition 
\begin{align*}
H_{|\sigma}: \sigma \rightarrow C\partial \sigma \xrightarrow{\xi} \mathscr{C}_{\sigma}
\end{align*} 
where $\xi(x,t) = (1-t)b_{\mathscr{C}_{\sigma}} + t H_{|\partial \sigma}(x)$. 
To show that this extension is Lipschitz, we assume that for every simplex $\tau$ (equipped with the norm-induced metric) of dimension $\leq k-1$, the restriction $H_{|\tau}: \tau \rightarrow \mathscr{C}_{\tau}$ is uniformly Lipschitz with Lipschitz constant $L_{k-1}$. Since $\mathscr{C}_{\sigma}$ is at most $k$-dimensional and $bs(K) \times [0,1]$ is uniformly bounded, by previous arguments, there exists a uniform Lipschitz constant $L_{k}$ such that $H_{|\sigma}: \sigma \rightarrow \mathscr{C}_{\sigma}$ is $L_{k}$-Lipschitz for all $k$-simplices $\sigma$. \\

By induction, there is a globally defined homotopy $H: bs(K) \times [0,1] \rightarrow M$ which is uniformly locally Lipschitz with respect to the set of simplices in $bs(K) \times [0,1]$. By previous arguments, $H$ is globally Lipschitz. It is proper because $\tilde{f}$ is proper. Therefore, it is a proper, Lipschitz homotopy between $\tilde{f}$ and $bs(\tilde{f})$. \\

\begin{remark}
More generally, this proof works for the barycentric subivision of $K$ relative to a subcomplex $L$ (induct only over simplices not in $L \times [0,1]$), or the standard subdivision $S(K)$ of $K$.
\end{remark}

\item Consider a map $f: K\rightarrow M$ which satisfies the geometric star condition. If $f$ is already Lipschitz, then a more careful version of the argument in point $2.$ of this lemma shows that $\tilde{f}$ is locally Lipschitz homotopic to $f$. One has to be more careful because geometric neighbourhoods are generally more complicated than geometric cores. For example, for a simplex $\Delta^{m}$ in $K$, the image $f(\Delta^{m})$ may intersect a height slice, in which case $\{u_{v_{0}}, \dots, u_{v_{m}}\}$, where $\tilde{f}(v_i) = u_{v_i} = (z_{v_i}, t)$, lie on a single height $t$, and define a geometric simplex of type $1$, whose neighbourhood contains $f(\Delta^{m})$. Alternatively $f(\Delta^{m})$  is fully contained within one open height interval, and $\{u_{v_{0}}, \dots, u_{v_{m}}\}$ may define a geometric simplex of type $1$ or $2$. \\

 In the product simplicial structure on $K \times [0,1]$ we want to extend the map $H$ which we fix to be $f$ on $K\times \{0\}$ and $\tilde{f}$ on $K\times \{1\}$. Denote by $pr$ the projection $pr: K\times [0,1] \rightarrow K$.  Let  $\sigma = [v_0, \dots, v_{j-1}, v_j, v'_j, v'_{j+1}\dots, v'_{m}]$ (or some subset of vertices thereof) be a simplex in $K\times [0,1]$ which does not lie entirely in $K \times \{0,1\}$, where $v'_{i}$ denotes a copy of $v_i \in K \times \{0\} = K$ on $K \times \{1\}$. We have that $H(v_i) = f(v_i)$ for $0\leq i\leq j$ and $H(v'_i) = u_{v_i}$ for $j\leq i\leq m$.  \\

Consider now the simplicial neighbourhood of $pr(\sigma)$ in $K$, denoted by $N(pr \sigma)$, consisting of all the closed simplices which contain $pr(\sigma)$ as a face. By assumption we have
\begin{align*}
f (pr \sigma) \subset f \bigcap_{v \in vert N(pr \sigma)} |st(v)| \subset  \bigcap_{v \in vert N(pr \sigma)} \langle gst(\tilde{f}(v)) \rangle 
\end{align*}
Therefore, the geometric simplex spanned by the vertices $\cup_{v \in vert N(pr \sigma)}\tilde{f}(v)$ exists. We denote the geometric neighbourhood of $[\cup_{v \in vert N(pr \sigma)}\tilde{f}(v)]$  by $\mathcal{M}_{\sigma}$. Assume that for every face $\tau \subset \partial \sigma$,  $H_{|\tau}$ lies within the intersection of geometric neighbourhoods $\mathcal{M}_{\sigma} \cap \mathcal{N}_{\tilde{f}pr(\sigma)} $, where $\mathcal{N}_{\tilde{f}pr(\sigma)}$ denotes the geometric neighbhourhood of the geometric simplex spanned by $\cup_{v\in vert pr(\sigma)} \tilde{f}(v)$. Note that this is automatically true for faces in $K \times \{0,1\}$: for $\tau\subset K \times \{0\}$ (which we identify with $pr(\tau) \subset K$) we have that 
\begin{align*}
f(\tau) \subset f(\bigcap_{v \in vert N(\tau)} |st(v)|) \subset \bigcap_{v\in vert N(\tau)}\langle gst (\tilde{f}(v))\rangle \subset \bigcap_{v\in vert N(pr \sigma)}\langle gst (\tilde{f}(v))\rangle \\
\subset \mathcal{M}_{\sigma} \cap \mathcal{N}_{\tilde{f}pr(\sigma)}
\end{align*} 
Likewise, for $\tau\subset K \times \{1\}$ we have that 
\begin{align*}
\tilde{f}(\tau) \subset \mathscr{C}_{\tilde{f}\tau} \subset  \mathscr{C}_{\tilde{f} N(pr \sigma)}\cap \mathscr{C}_{\tilde{f} pr \sigma} \subset  \mathcal{M}_{\sigma} \cap \mathcal{N}_{\tilde{f}pr(\sigma)}
\end{align*} 
Recall that $\mathcal{N}_{\tilde{f}pr(\sigma)}$ deformation retracts to a point, such that the homotopy restricted to $\mathcal{M}_{\sigma}\cap \mathcal{N}_{\tilde{f}pr(\sigma)}$ has image in $\mathcal{M}_{\sigma}\cap \mathcal{N}_{\tilde{f}pr(\sigma)}$. Therefore, we can extend $H$ to the interior of $\sigma$. Now, for the induction step, suppose $\sigma \subset \partial\sigma'$ is a face on the boundary of $\sigma'$. We have the following sequence of inclusions:
\begin{align*}
pr \sigma \subset pr \sigma' \subset N(pr \sigma') \subset N(pr \sigma)
\end{align*}
This gives us the inclusions
\begin{align*}
\mathcal{M}_{\sigma} \subset cl(\mathcal{M}_{\sigma'})\\
 \mathcal{M}_{\sigma'}\subset cl(\mathcal{N}_{\tilde{f}pr \sigma'})\\
\mathcal{N}_{\tilde{f}pr \sigma'}\subset cl(\mathcal{N}_{\tilde{f}pr \sigma})
\end{align*}
where we can remove the closures if the two geometric neighbouroods are of the same type. \\

There is at most one change of type in this chain, but even if there were, we know that a point in $\mathcal{N}_{\tilde{f} pr \sigma} \cap \mathcal{M}_{\sigma}$ is contained in $\mathcal{I}_{\ast}\mathcal{M}_{\sigma} $, where $\ast \in \{0,1\}$ corresponds to the height of $\tilde{f} pr \sigma$: ie. we cannot lie in the boundary component in a neighbourhood of a simplex of type $2$ not present in the neighbourhood of a simplex of type $1$. Therefore, we have the inclusions
\begin{align*}
\mathcal{N}_{\tilde{f} pr \sigma} \cap \mathcal{M}_{\sigma} \subset \mathcal{N}_{\tilde{f} pr \sigma}\cap \mathcal{M}_{\sigma'} \subset \mathcal{N}_{\tilde{f} pr \sigma'} \cap \mathcal{M}_{\sigma'}
\end{align*}
as required.\\

Therefore, by induction on dimension, there is a homotopy $H: K\times [0,1] \rightarrow M$ between $f$ and $\tilde{f}$. $H$ is locally Lipschitz with respect to the set of simplices in $K \times [0,1]$. The reason that $H$ may not be globally Lipschitz is because a simplex in a geometric neighbourhood may have arbitrarily large dimension. If $f_{|L} = \tilde{f}_{|L}$ we can choose $H$ to fix the subcomplex $L$ by inducting only over simplices not contained in $L \times [0,1]$. If $f$ is proper, then $H$ is proper.\\

\item  Let $K,\overline{K},L, \varphi, \tilde{f}, \tilde{f'}$ be as in the proposition. We show that the relative simplicial approximation theorem (\cite{zeeman1964relative}) is applicable to the the map $\varphi: K\rightarrow \overline{K}$ since $\varphi$ is Lipschitz and $K,\overline{K}$ are uniformly bounded: the cover by open stars in $\overline{K}$ has strictly positive Lebesgue number and the diameter of simplices in $K$ approaches $0$ with iterated barycentric subdivision. To understand this proof, it is required to first read Subsection \ref{relgeosimplsect} in the Appendix for the notation used. \\

First we show that the cover by open stars in $\overline{K}$ has strictly positive Lebesgue number. To that end, we consider a different metric on $\overline{K}$. Let $d'$ be the path metric induced from the standard metric $d_{std}$ (see Definition \ref{metstd}) on $\overline{K}$. We denote this metric space by $(\overline{K},d')$. \\

\begin{claim}
For any $n$-dimensional abstract simplicial complex $K$, the simplicial complex $(|K|,d')$ has Lebesgue number $\frac{1}{n+1}$. 
\end{claim}
To see this, let $U$ be a set of diameter less than $\frac{1}{n+1}$ and $x\in U$. $x = \sum_{v \in vert \sigma} x_{v}v$ is contained in the interior of a unique simplex $\sigma \in K$. We choose the closest vertex $y$ in $\sigma$ with respect to $d_{std}$. By the definition of the metric, we have that $x_{y}\geq \frac{1}{n+1}$. Let $z = \sum_{v\in vert K}z_{v}v$ be a point in $U$. $z$ has $y$-coordinate strictly greater than $0$, and therefore lies in the open star of $y$. \\ 

Now we consider $(\overline{K},d)$ with our given metric. Since $d$ is uniformly bounded, the identity map $\id: (\overline{K},d)\rightarrow (\overline{K},d')$ is Lipschitz, with Lipschitz constant $C$. Therefore the Lebesgue number of the covering by open stars in $(\overline{K},d)$ is $\frac{1}{C(n+1)}$.\\

Let $D$ be the Lipschitz constant for $\varphi$. Choose $\varepsilon>0$ small enough such that $D\varepsilon<1$. By Corollary \ref{starconv} there exists a $q$ large enough such that $|st(v,K_{q})|\subset N_{\varepsilon} |st(v,L)|$ for all vertices $v$ in L. Fix this $q$. As in the proof of the relative simplicial approximation theorem, there exists a simplicial map $\mathfrak{h}:K_{q}'\rightarrow K_{q}$ which is a simplicial approximation of the identity map on $K$. Since there are only finitely many strong similarity types of simplices in $K'_{q}$ and $K_{q}$ (each type of simplex only leads to a finite number of new types of simplices  in the subdivision) by Lemma \label{Lbound} the map $\mathfrak{h}$ is Lipschitz. Therefore, the composition $\varphi\circ \mathfrak{h}: K_{r}\rightarrow \overline{K}$ is a Lipschitz map for all $r\geq q$. By the assumptions on $\overline{K}$, the Lebesgue number of the cover of open stars of vertices in $\overline{K}$ is strictly positive.  Using the fact that $\varphi\circ \mathfrak{h}$ is Lipschitz, we obtain that the open cover $\beta:=(\varphi\circ \mathfrak{h})^{-1}(\mathcal{W})$, where $\mathcal{W}$ denotes the star covering of $\overline{K}$, has strictly positive Lebesgue number. Let $\alpha_{r}$ denote the star covering of $K_{r}$. The subset 
\begin{align*}
\alpha'_r:= \{|st(v,K_r)| \,|\, v\in V_r\}
\end{align*}
satisfies the property that as $r \rightarrow \infty$, mesh $\alpha'_r\rightarrow 0$ (again, since $K$ is uniformly bounded). Choose $s$ such that the mesh of $\alpha'_s$ is less than the Lebesgue number of $\beta$. Therefore $\alpha'_s$ refines $\beta$. By the construction of $\mathfrak{h}$, $\alpha_{s}$ also refines $\beta$. \\

Therefore $\varphi\circ \mathfrak{h}: K_{s}\rightarrow \overline{K}$ satisfies the star condition and has a simplicial approximation $g$. Additionally, $g_{|L} = \varphi_{|L} = \id_{L}$  and $g$ is homotopic to $\varphi$ keeping $L$ fixed. From the proof of the relative simplicial approximation theorem, we also have that $\mathfrak{h}|st(y,K)|\subset |st(y,K)|$ for any vertex $y\in L$. \\

Let $\tilde{f}_{s}: K_{s}\rightarrow M$ denote the simplicial approximation of $f$ with respect to $K_{s}$. We show that $\tilde{f_{s}}$ and $\tilde{f}' \circ g$ are both simplicial approximations of $f\circ \mathfrak{h}$.  The statement for $\tilde{f}' \circ g$ is obvious:
\begin{align*}
f\varphi^{-1}\varphi\mathfrak{h}|st(v)| \subset f\varphi^{-1} |st(g(v))| \subset  \langle gst(\tilde{f}'g(v))\rangle 
\end{align*}
 For $\tilde{f_{s}}$, let $x$ be a vertex in $K_{s}$. If $x$ lies in the set $cl(K_{q}-N(L,K_{q}))$ then $\mathfrak{h}$ is the identity on $|st(x,K_{s})|$ for any $s>q+1$, so $f\mathfrak{h}|st(x,K_{s})| = f|st(x,K_{s})| \subset \langle gst( \tilde{f_{s}} (x))\rangle$. Similarly if $x\in L$ then $\mathfrak{h}|st(x,K_{s})|\subset |st(x,L)|$, so $f\mathfrak{h}|st(x,K_{s})| \subset f|st(x,L)| \subset \langle gst(\tilde{f}(x))\rangle=\langle gst(\tilde{f_{s}}(x))\rangle$. Finally, if $x$ lies in the interior of  $ N(L,K_{q})- L$ then $x$ lies in the interior of some simplex $\sigma$ in $K_{q}-L$  with $y$ a vertex of $\sigma \cap L$. By construction we can always choose $\tilde{f_{s}}(x) = \tilde{f}(y)$. We also have
\begin{align*}
\mathfrak{h}|st(x,K_{s})| \subset \mathfrak{h}|st(y,K)| \subset |st(y,K)|
\end{align*}
and therefore $f\mathfrak{h}|st(x,K_{s})|\subset f(|st(y,K)|) \subset \langle gst(\tilde{f}(y))\rangle = \langle gst(\tilde{f_{s}}(x))\rangle$. Therefore $\tilde{f_{s}}$ is a simplicial approximation of $f\circ \mathfrak{h}$. Since $\tilde{f_{s}}$ and $\tilde{f}' \circ g$ restricted to $L$ are equal, the homotopy between them fixes $L$. \\

We have the following sequence of homotopies 
\begin{align*}
\tilde{f}\simeq_{L} \tilde{f_{s}} \simeq_{L} \tilde{f}' \circ g \simeq_{L} \tilde{f}'\circ \varphi
\end{align*}
where the symbol $\simeq_{L}$ stands for a homotopy keeping $L$ fixed. All the homotopies are Lipschitz. Since the homotopy $ g \simeq_{L} \varphi$ is bounded, all homotopies are proper if $f$ is proper.\\

The absolute case is a special case of the relative statement for the subcomplex $L=\emptyset$. \\
\end{enumerate} 
\end{proof}

We now show that $\theta$ is well-defined and independent of choices. Most of these statments are obtained as corollaries of Proposition \ref{geosa}. We identify the set  $c[0,1]^n$ with $c[-1,1]^n$ occasionally to make change of base ray homomorphisms easier to write.

\begin{lemma} (Well-definedness of $\theta$) \label{indepwelldefined}
\begin{enumerate}
\item Independence of choice of simplicial approximation.
\item Independence of choice of shrinking map $s_{r}$ for sufficiently large $r$.
\item Barycentric subdivision of $c[-1,1]^n$.
\item Independence of coarse homotopy class of $f$.
\item $\theta$ is a group homomorphism. 
\end{enumerate}
\end{lemma}

\begin{proof}

\begin{enumerate}
\item Independence of choice of simplicial approximation:

Suppose we have $f: (c[-1,1]^n,c\partial [-1,1]^n)\rightarrow (cX,\omega)$ and two simplicial approximations $\widetilde{is_{r}f},\widetilde{is_{r}f}'$ of $i s_{r} f$. By Proposition \ref{geosa} they are homotopic via a proper Lipschitz homotopy, which in particular can be used to define a coarse-Lipschitz homotopy by extending constantly. Denote this extended homotopy by $H: c([-1,1]^n \times [0,1])\rightarrow M$.\\

Let $\omega_r, \omega'_r$ be the base rays associated to $\widetilde{is_{r}f}$ and $\widetilde{is_{r}f}'$ respectively. Observe that $H_{|c(\partial[-1,1]^n \times [0,1])}$ is a coarse-Lipschitz homotopy between $\omega_r$ and $\omega'_r$, which we denote by $H^{\omega_{r},\omega'_r}$. Let $H^{\omega_r,\omega},H^{\omega'_{r},\omega}$ be linear homotopies between $\omega_r,\omega$ and between $\omega'_{r},\omega$ respectively. By Lemma \ref{baserayindep} there is a proper Lipschitz homotopy $F: c([0,1]^2)\rightarrow \{x_0\} \times [1,\infty)$ between $H^{\omega_r,\omega}$ and $H^{\omega'_{r},\omega}$, which restricts to $H^{\omega_{r},\omega'_r}$ and $\id$ on $c(\{0\} \times [0,1])$ and $c (\{1\} \times [0,1])$ respectively. \\ 
 
Recall that the set $U = cW$ is defined as 
 \begin{align*}
 U:= \{(hx_1,\dots,hx_{n},h)\subset c[-1,1]^n\,|\, -\frac{1}{2}\leq x_i\leq \frac{1}{2} \,|\,\forall 1\leq i\leq n\}
 \end{align*}
 We construct a map $G: c([-1,1]^n \times [0,1])\rightarrow M$ as follows: 
 \begin{align*}
 G_{|c(W\times [0,1])}: c(W\times [0,1]) &\longrightarrow M\\
 (h\mathbf{x},ht, h)&\longmapsto H(2h\mathbf{x},ht,h)\\
 G_{|c(([-1,1]^n \setminus W)\times [0,1])}: c(([-1,1]^n \setminus W)\times [0,1]) &\longrightarrow M\\
 (h\mathbf{x},ht,h) &\longmapsto F(h(2\|\mathbf{x}\|_{\infty}-1),ht,h)
 \end{align*}
 Basically on $c(W \times [0,1])$ we have the given homotopy $H$ and on the complement we have the homotopy between the two changes of base rays. See Figure \ref{fig:isolip1}. $G_{|c(W\times [0,1])}$ is a composition of a $2$-Lipschitz map and a Lipschitz map $H$, so it is Lipschitz. We know from previous computation that the map $(h\mathbf{x},ht,h)\mapsto (h(2\|\mathbf{x}\|_{\infty}-1),ht,h)$ is $2$-Lipschitz.  $G_{|c(([-1,1]^n \setminus W)\times [0,1])}$ is a composition of a $2$- Lipschitz and a Lipschitz map, thus it is Lipschitz. Since the two maps agree on their common intersection, $G$ is globally Lipschitz. Properness is clear. Therefore, $b_{\omega_r,\omega}[\widetilde{is_{r}f}]=b_{\omega_r',\omega}[\widetilde{is_{r}f}']\in \pi_{n}^{L,e}(M,\omega)$.

\begin{figure}
\centering
  \centering
  \includegraphics[width=0.5\linewidth]{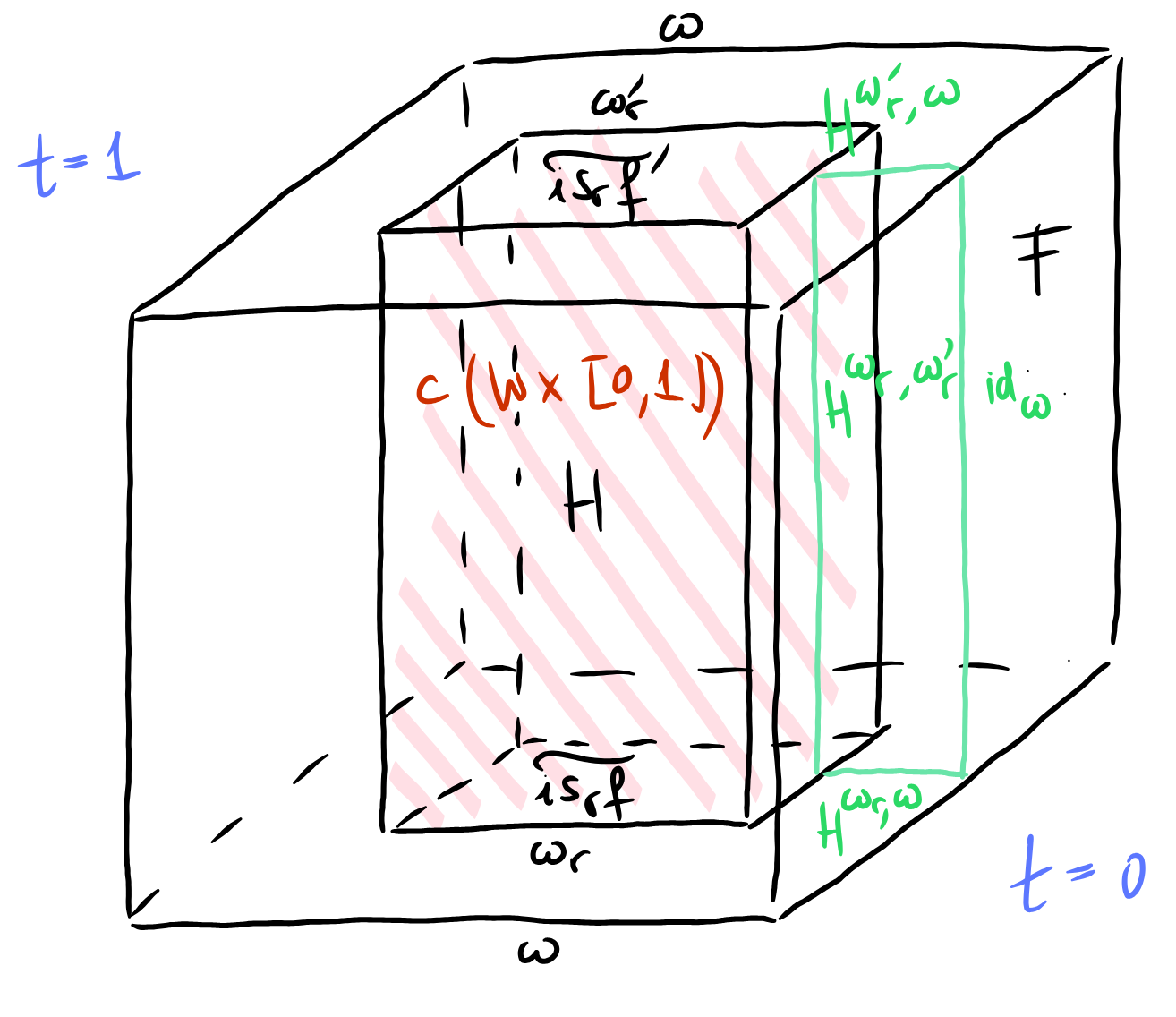}
  \caption{The coarse homotopy between $b_{\omega_{r},\omega}[\widetilde{is_{r}f}]$ and $b_{\omega_{r}',\omega}[\widetilde{is_{r}f}']$}
  \label{fig:isolip1}
\end{figure}

\item Independence of choice of shrinking map $s_{r}$ for sufficiently large $r$:

Let $s_{1}: cX\rightarrow cX$ be the shrinking map that divides by $2$. By Lemma \ref{shrinkindep} there exists a coarse homotopy $\mathbf{H}^{1}: I_{p}(cX)\rightarrow cX$ between $\id$ and $s_{1}$.  Choose a simplicial structure on $c([0,1]^{n+1})$ such that it induces our fixed simplicial structure of $c[0,1]^n$ on $c([0,1]^n \times \partial [0,1])$ (this is achieved by choosing a simplicial structure on $[0,1]^{n+1}$ so that induces our fixed simplicial structure of $[0,1]^n$ on $[0,1]^{n} \times \partial [0,1]$ and running the coning construction of finite simplicial complexes). 

Then we define:
\begin{align*}
G: c([0,1]^n\times [0,1])&\rightarrow cX\\
(h\mathbf{x},h,ht)&\mapsto \mathbf{H}^{1}(f(h\mathbf{x},h), pf(h\mathbf{x},h)t)
\end{align*}
That is, on each time slice $ht$ we apply $f$, then take the point along the homotopy $\mathbf{H}^{1}$ with appropriate rescaling with $p: cX\rightarrow [1,\infty)$. $G$ is a homotopy from $f$ to $s_{1} f$. To show that $G$ is coarse, let $(h\mathbf{x},h,ht), (s\mathbf{y},s,sl)$ be two points with distance less than $Q$. Since $f$ is coarse, $f(h\mathbf{x},h), f(s\mathbf{y},s)$ lie in an entourage. $p\circ f: c[0,1]^n\rightarrow [1,\infty)$ is large-scale Lipschitz: let $A$ be the large-scale Lipschitz constant for it. 
\begin{align*}
|pf(h\mathbf{x},h)t-pf(s\mathbf{y},s)l|&\leq |pf(h\mathbf{x},h)t-pf(s\mathbf{y},s)t|+|pf(s\mathbf{y},s)t-pf(s\mathbf{y},s)l|\\
&< AQ+A + |t-l| (|pf(s\mathbf{y},s) - pf (s \partial [0,1]^n,s)|+s)\\
&< AQ+A+ |t-l|(sA+A+s)\\
&\leq AQ+2A+(A+1)(|ht-ls|+|h-s|)\\
&< CQ +C
\end{align*} 
for some constant $C>0$. Therefore $(f(h\mathbf{x},h), pf(h\mathbf{x},h)t), (f(s\mathbf{y},s), pf(s\mathbf{y},s))$ lie in an entourage. Since $\mathbf{H}^{1}$ is controlled, their images lie in an entourage. This shows that $G$ is controlled. For properness, let $K$ be a bounded set. Observe that since $\mathbf{H}^{1}$ is proper, $(\mathbf{H}^{1})^{-1}(K)$ is contained within $Z\times [0,\infty)\cap I_{p}(cX)$ for $Z$ bounded. Now since $f$ is proper, $G^{-1}(K)$ is contained within a bounded set of heights $[1,l]$. Therefore $G$ is proper.  

\begin{remark} This computation also yields the result that $s_{\gamma}\circ f$ is homotopic to $f$ for any choice of shrinking map $s_{\gamma}$. This can be concluded abstractly by Proposition \ref{properties} and Lemma \ref{p_0} but it is  useful in later arguments to have an explicit formula for a homotopy.
\end{remark}

As before, there exists a $S>0$ such that $G(|st(v)|)$ has diameter $<S$ for all vertices $v\in c[0,1]^{n+1}$, and a $r$ large enough so that $i s_{r} G$ has a simplicial approximation $\widetilde{is_{r}G}$. Observe that $\widetilde{is_{r}G}_{|c([0,1]^n \times \{0\})}$ and $\widetilde{is_{r}G}_{|c([0,1]^n \times \{1\})}$ are simplicial approximations of $i s_{r}  f$ and $i  s_{r} s_{1} f = i s_{r+1} f$ respectively (since for a vertex, the star considered in the subcomplex $c([0,1]^n \times \partial[0,1])$ is contained within the star considered in $c[0,1]^{n+1}$. After changing base rays with the same argument as point $1$ of this lemma, $b_{\omega_{r},\omega}[\widetilde{is_{r}f}] = b_{\omega_{r+1},\omega}[\widetilde{is_{r+1}f}] \in \pi_{n}^{L,e}(M,\omega)$. By replacing $r$ by $r+1$, we can inductively show that  $b_{\omega_{r'},\omega}[\widetilde{is_{r'}f}] = b_{\omega_{r},\omega} [\widetilde{is_{r}f}]$ for all $r'\geq r$. 

\item Barycentric subdivision of $c[-1,1]^n$:

Suppose we take the $k$-th barycentric subdivision of the fixed simplicial structure on $c[-1,1]^n$, denoted by $bs^k(c[-1,1]^n)$. Assume that $r$ is sufficiently large as in point $2.$ of this lemma. By Proposition \ref{geosa} the simplicial approximation of $is_{r}f$ with respect to $bs^k(c[-1,1]^n)$ exists (denoted by $(\widetilde{is_{r}f})_{k}$) and is proper Lipschitz homotopic to $\widetilde{is_{r}f}$ via a homotopy $H$. Let $J^{k}: c(\partial[-1,1]^n) \times [0,1]\rightarrow M$ denote $H$ restricted to the boundary. Let $B^{J^{k}}(\widetilde{is_{r}f})_{k}$ be the map which is $(\widetilde{is_{r}f})_{k}$ when restricted to the set $U$, and the homotopy $J^{k}$ (rescaled) on the set $c[-1,1]^n \setminus U$, which can be identified topologically with $c(\partial[-1,1]^n \times [0,1])$. 

We have that
\begin{align*}
b_{\omega_r,\omega}[\widetilde{is_{r}f}] = b_{\omega_r,\omega}[B^{J^{k}}(\widetilde{is_{r}f})_{k}] 
\end{align*}

The reason we use a different symbol $B$ here is because $(\widetilde{is_{r}f})_{k}$ restricted to $c \partial [-1,1]^n$ does not need to take the form $\omega'\circ p$ for a base ray $\omega': [1,\infty)\rightarrow M$. 

\item Independence of coarse homotopy class of $f$:

Suppose that $[g]=[f]$ in $\pi_n^c(cX,\omega)$. This means that there is a coarse homotopy $H: (c[0,1]^{n+1}, c(\partial [0,1]^{n}\times [0,1]))\rightarrow (cX,\omega)$ which restricts to $f$ and $g$ on $c([0,1]^n \times \{0\})$ and $c([0,1]^n \times \{1\})$ respectively. There exists a $S>0$ such that $H(|st(v)|)$ has diameter $<S$ for all vertices $v\in c[0,1]^{n+1}$, and a $r$ large enough so that $i s_{r}  H$ has a simplicial approximation $\widetilde{is_{r}H}$.  We can assume that $r$ is larger than the constants for $f$ and $g$ from point $2.$ of this lemma. $\widetilde{is_rH}_{|c([0,1]^n \times \{0\})}$ and $\widetilde{is_{r}H}_{|c([0,1]^n \times \{1\})}$ are simplicial approximations of $i  s_{r}  f$ and $i  s_{r}  g$ respectively.  This means that $[\widetilde{i  s_{r}  f}] = [\widetilde{i  s_{r}  g}] \in \pi_{n}^{L,e}(M,\omega_r)$ and after changing base rays we obtain $\theta[f] = \theta[g]$. 

\item This comes from the fact that $\theta$ is inverse to a group homomorphism. 
\end{enumerate}
\end{proof}

Now we resume the proof of Theorem \ref{isoomega}. We show that $\theta$ and $\lambda$ are inverse to each other. First we compute the composition $\lambda \circ \theta: \pi_{n}^{c}(cX,\omega)\rightarrow \pi_{n}^{c}(cX,\omega)$. Let $f: (c[-1,1]^n,c\partial[-1,1]^{n})\rightarrow (cX,\omega)$ be a coarse map representing a class in $\pi_{n}^{c}(cX,\omega)$. We show that $\lambda\circ \theta[f]$ is represented by a map which is coarsely homotopic to $f$. The composition is represented by $R \circ b^{H^{r,0}}_{\omega_r,\omega}[\widetilde{i s_{r} f}]$, which, as a reminder, just $R \circ \widetilde{i s_{r} f}\circ \mathfrak{r}$  when restricted to the set $U = \{(h\mathbf{x},h)\,|\, \mathbf{x}\in [-\frac{1}{2}, \frac{1}{2}]^n\}$ and the linear homotopy from $R\omega_{r}=\omega_r$ to $R \omega=\omega$ on its complement via $H^{r,0}$, ie.

\begin{align*}
R \circ b^{H^{r,0}}_{\omega_r,\omega}[\widetilde{i s_{r} f}]_{|c[-1,1]^n \setminus U}: c[-1,1]^n\setminus U &\longrightarrow cX\\
(h\mathbf{x},h) &\longmapsto H^{r,0} (h(2\|\mathbf{x}\|_{\infty}-1), h)
\end{align*}
See Figure \ref{fig:isolip2}.\\

\begin{figure}[H]
\centering
  \centering
  \includegraphics[width=0.6\linewidth]{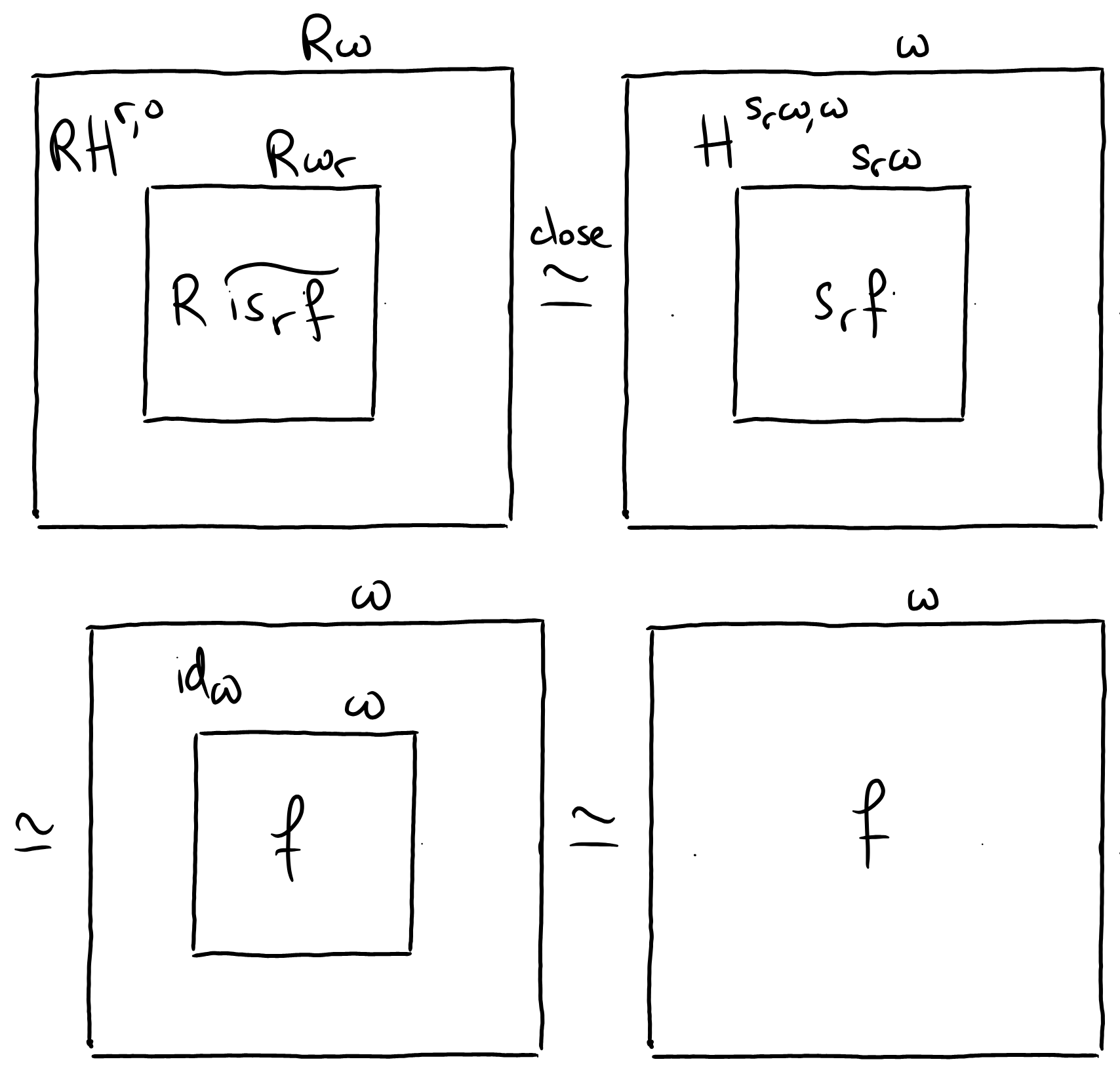}
  \caption{A representative of $\lambda\theta[f]$.}
  \label{fig:isolip2}
\end{figure}

$R\widetilde{is_{r}f}$ is close to $Ris_{r}f$ since $is_{r}f$ is close to its geometric simplicial approximation, and $R$ is a coarse map. We know that $Ri$ is close to the identity, and therefore on $U$, we can replace the map $R \circ \widetilde{i s_{r} f}\circ \mathfrak{r} \cdot$ by $s_{r} f\circ \mathfrak{r}$.  On the complement of $U$, observe that 
\begin{align*}
H^{r, 0}(h(2\|\mathbf{x}\|_{\infty}-1), h) &= (2-2\|\mathbf{x}\|_{\infty})\omega_{r}(h)+h(2\|\mathbf{x}\|_{\infty}-1) \\
H^{s_{r}\omega, \omega}(h(2\|\mathbf{x}\|_{\infty}-1), h) &=  (2-2\|\mathbf{x}\|_{\infty})\frac{h}{2^r}+h(2\|\mathbf{x}\|_{\infty}-1)
\end{align*}
where we have identified $\omega(h) = (hx_0,h)$ with its height variable. Since $|\omega_{r}(h) - \frac{h}{2^r}|<1$, we can replace $H^{r,0}$ on the complement of $U$ with $H^{s_r\omega,\omega}$, which is close. Putting it together, $R \circ b_{\omega_r,\omega}[\widetilde{i s_{r} f}] = b_{s_{r}\omega,\omega}[s_{r}f]$. Since $s_{r}f$ is coarsely homotopic to $f$, we have that 
\begin{align*}
b_{s_{r}\omega,\omega} [s_{r}f]=b_{\omega,\omega} [f] = [f]
\end{align*}
Therefore $\lambda\circ \theta = \id$.\\

Now we compute $\theta\circ \lambda:\pi_{n}^{L,e}(M,\omega)\rightarrow \pi_{n}^{L,e}(M,\omega)$. Let $f: (c[-1,1]^n,c\partial[-1,1]^n)\rightarrow (M,\omega)$ be a proper Lipschitz map representing a class in $\pi_{n}^{L,e}(M, \omega)$. There exists a $k$ large enough such that the $k$-th barycentric subdivision of $c[-1,1]^n$ satisfies the property that for every vertex $v\in bs^k(c[-1,1]^n)$, $f(|st(v)|)\subset \langle gst(z_v,t_v) \rangle $ where $(z_v,t_v)$ is the closest vertex to the image of $v$. We can do this because the covering of open stars in $M$ has Lebesgue number  $\frac{1}{4}>0$ and $f$ is Lipschitz. For the class $[Rf]$, choose $r$ sufficiently large as in point $2$. of Lemma \ref{indepwelldefined}. $i s_{r} R f$ has a simplicial approximation $(i s_{r} R f)^{\sim}$. We show that a simplicial approximation $(i s_{r} R f)^{\sim}_{k}$ of $(i s_{r} R f)$ with respect to $bs^{k}(c[-1,1]^n)$ is homotopic to $f$. See Figure \ref{fig:isolip3}. \\

\begin{figure}[H]
\centering
  \centering
  \includegraphics[width=0.6\linewidth]{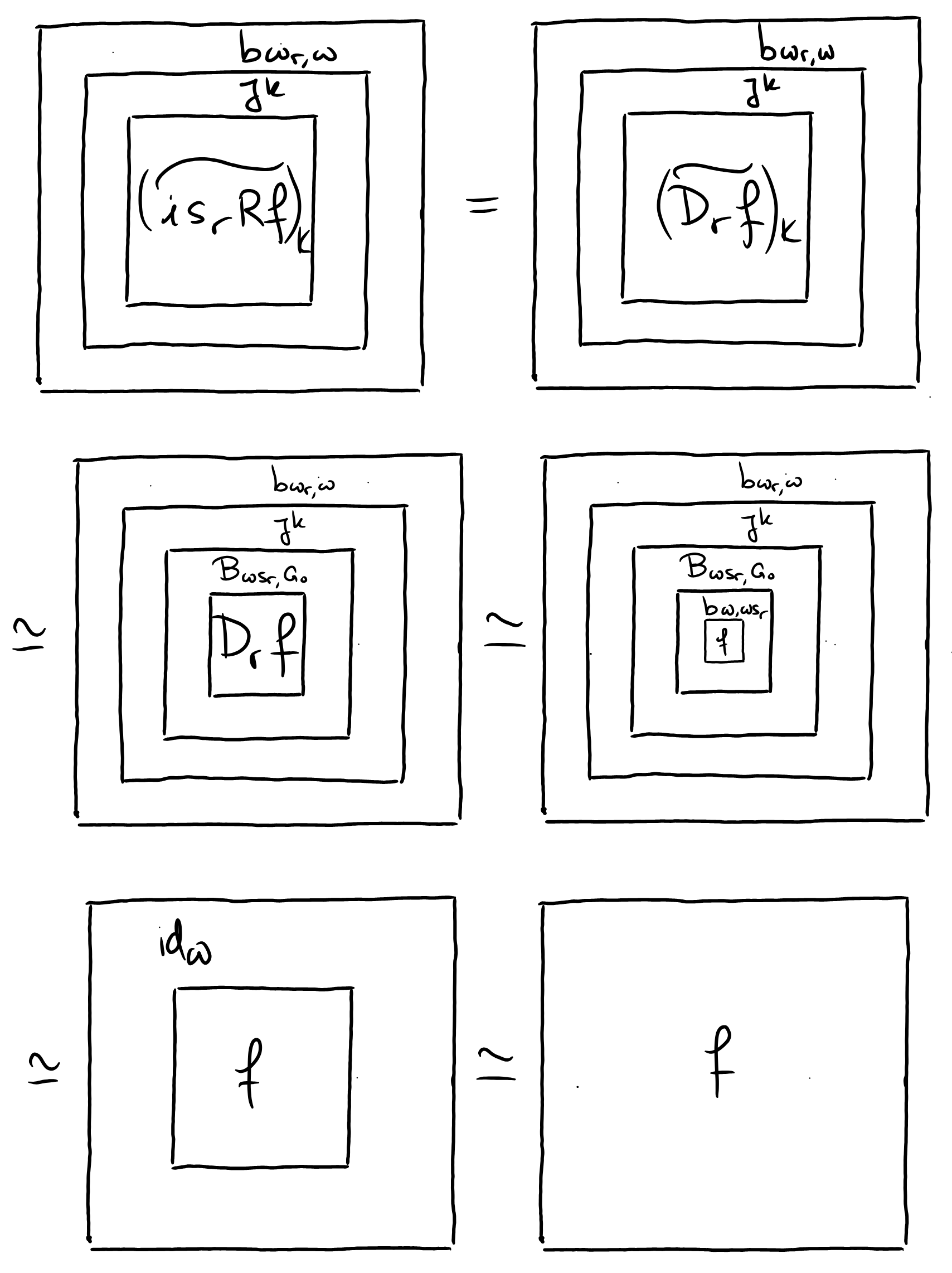}
  \caption{A representative of $\theta\lambda[f]$}
  \label{fig:isolip3}
\end{figure}

Let $(z_{v},t_{v})$ be a vertex in $|\mathcal{U}_{h}| \times (2^{h},2^{h+1})$. The image of $R \langle gst(z_{v}, t_{v})\rangle$ consists of points whose projections in $X$ live inside the open ball around $z_{v}$ and whose heights in $cX$ live in the interval $(t_{v}-1,t_{v}+1)$. It is clear now that these projections also live in the open ball around $\phi_{h,h-r}(z_{v})$ and their heights divided by $2^{r}$ are in the interval $(p_{r}(t_{v})-1,p_{r}(t_{v})+1)$ where $p_{r}(t_{v})$ is the closest integer to $\frac{t_{v}}{2^r}$.

Therefore we have
\begin{align*}
is_{r}R \langle gst(z_v,t_v) \rangle \subset \langle gst(\phi_{h,h-r}(z_v),p_{r}(t_v)) \rangle 
\end{align*}  
where $\phi_{h,h-r}(z_v)$ is identified with its image $\phi_{h,h-r-1}(z_v)$ if $p_{r}(t_v) = 2^{h-r}$. \\

On the other hand, let $D_{r}$ denote the map on the telescope $M$ which divides heights by $2^{r}$. It is defined as follows: for $(x,t)\in |\mathcal{U}_{h}| \times (2^h,2^{h+1}]$ with $h\geq r$ we have $D_{r}(x,t) =  (\phi_{h,h-r}(x),\frac{t}{2^{r}})$. Otherwise we have $D_{r}(x,t) = (x_0,1)$. There is a deformation retract $\mathbf{D}^{r}: M\times [0,1]\rightarrow M$ from the identity to $D_{r}$ and we have 
\begin{align*}
D_{r}\langle gst(z_{v},t_{v})\rangle \subset \langle gst(\phi_{h,h-r}(z_{v}),p_{r}(t_{v})) \rangle  
\end{align*}

If $t_{v}=2^h$, recall that from the definition of geometric stars we have
\begin{align*}
&is_{r}R (\langle gst(z_{v},t_{v}\rangle) \\
&= is_{r}(\cup_{z\in \phi_{h}^{-1}(z_{v})} R \langle st(z) \rangle \times [2^h, 2^h+1)) \bigcup is_{r} (R\langle st(z_{v})\rangle \times (2^{h}-1,2^h])\\
&=i(\cup_{z\in \phi_{h}^{-1}(z_{v})} R \langle st(z) \rangle \times [2^{h-r}, 2^{h-r}+\frac{1}{2^r})) \bigcup i (R\langle st(z_{v})\rangle \times (2^{h-r}-\frac{1}{2^r},2^{h-r}])\\
&\subseteq i(\cup_{z\in \phi_{h}^{-1}(z_{v})} R \langle st(z) \rangle \times [2^{h-r}, 2^{h-r}+1)) \bigcup i (R\langle st(z_{v})\rangle \times (2^{h-r}-1,2^{h-r}])\\
&\subseteq (\cup_{z\in \phi_{h}^{-1}(z_{v})}  \langle st(\phi_{h,h-r}(z)) \rangle \times [2^{h-r}, 2^{h-r}+1)) \bigcup \langle st(\phi_{h-1,h-r-1} (z_{v}))\rangle \times (2^{h-r}-1,2^{h-r}])\\
&\subseteq (\cup_{z\in \phi_{h-r}^{-1}(\phi_{h-1,h-r-1}(z_{v}))}  \langle st(z) \rangle \times [2^{h-r}, 2^{h-r}+1)) \bigcup \langle st(\phi_{h-1,h-r-1} (z_{v}))\rangle \times (2^{h-r}-1,2^{h-r}])\\
&= \langle gst(\phi_{h-1,h-r-1}(z_{v}), p_{r}(t_{v}))\rangle 
\end{align*} 
This comes from the fact that if $z\in \phi_{h}^{-1}(z_v)$ then  
\begin{align*}
\phi_{h,h-r}(z)=\phi_{h-1,h-r}(z_{v})\in \phi_{h-r}^{-1}(\phi_{h-1,h-r-1}(z_v))
\end{align*}
In this case we have also that 
\begin{align*}
D_{r}\langle gst(z_{v},t_{v})\rangle \subset \langle gst(\phi_{h-1,h-r-1}(z_{v}),p_{r}(t_{v})) \rangle  
\end{align*}
This shows that $is_{r}Rf$ and $D_{r}f$ have the same geometric simplicial approximation with respect to $bs^{k}(c[-1,1]^n)$ and that $\theta\circ \lambda[f]$ is represented by $b_{\omega_{r},\omega}[B^{J^{k}}(\widetilde{is_{r}Rf})_{k}] = b_{\omega_r,\omega} [B^{J^{k}}(\widetilde{D_{r}f})_{k}]$. We know that $(\widetilde{D_{r}f})_{k}$ is locally Lipschitz homotopic to $D_{r}f$. The homotopy restricted to the boundary, which we denote by $G: c\partial [-1,1]^n \times [0,1]$, has image contained in the base ray $\omega$. We denote $G_{0}= {(\widetilde{D_{r}f})_{k}}_{|c\partial[-1,1]^n}$, $G_1 = (D_{r}f)_{|c\partial [-1,1]^n}= \omega s_{r}$. Since $\omega$ is a $1$-dimensional subcomplex, $G$ is Lipschitz. From the proof of Lemma \ref{clip} in the next subsection, by precomposing with a shrinking map we have a coarse-Lipschitz homotopy between $(\widetilde{D_{r}f})_{k}$  and $D_{r}f$, via a homotopy that restricts to $G$ on the boundary. This gives us
\begin{align*}
(\widetilde{D_{r}f})_{k} \simeq B_{G_{0}',G_{0}} (\widetilde{D_{r}f})'_{k}\simeq B_{\omega s_{r},G_{0}}b_{(\omega s_{r})',\omega s_{r}} (D_{r}f)'  \simeq B_{\omega s_{r},G_{0}} D_{r}f
\end{align*} 
rel $c\partial[-1,1]^n$. Putting it together, we obtain
\begin{align*}
\theta\circ\lambda[f] &=b_{\omega_{r},\omega}[B^{J^{k}}(\widetilde{is_{r}Rf})_{k}]= b_{\omega_{r},\omega}[B^{J^{k}}(\widetilde{D_{r}f})_{k}]=b_{\omega_{r},\omega}[B^{J^{k}}B_{\omega s_{r},G_{0}} D_{r}f] \\
&= b_{\omega_{r},\omega}[B^{J^{k}}B_{\omega s_{r},G_{0}} b_{\omega,\omega s_{r}}f] 
\end{align*}

Observe that the map $b_{\omega_{r},\omega}B^{J^{k}}B_{\omega s_{r},G_{0}} b_{\omega,\omega s_{r}}: c(\partial[-1,1]^n  \times [0,1])\rightarrow M$ is a homotopy with image in $\omega$. One can show that 
\begin{align*}
b_{\omega_{r},\omega}B^{J^{k}}B_{\omega s_{r},G_{0}} b_{\omega,\omega s_{r}} \simeq b^{\id}_{\omega,\omega}
\end{align*}
by a linear homotopy. Therefore, $\theta\circ\lambda[f] =[f]$. 

\subsection{End homotopy groups}

The goal of this section is to prove the following proposition:
\begin{prop} There is an isomorphism \label{natiso}
\begin{align*}
\pi_{n}^{L,e}(M,\omega) \cong \pi_{n}^{e}(M,\omega)
\end{align*} 
\end{prop}

We begin with an introduction to end homotopy groups $\pi^e_{n}$. 

\begin{definition}
Let $X$ be a CW complex and $S\subset X$. The \textit{carrier} of $S$ is the intersection of all subcomplexes of $X$ which contain $S$. It is a subcomplex of $X$, denoted $C(S)$. A CW complex $X$ is called \textit{strongly locally finite} if 
\begin{align*}
\{C(e)\,|\, e \text{ is a cell of } X\}
\end{align*}
 is a locally finite cover of $X$. 
\end{definition}

\begin{example} \label{CW} $M$ can be given the structure of a CW complex, since the gluing maps $\phi_{h}$ are simplicial (and therefore cellular). We take the product CW-structure on each $|\mathcal{U}_{h}| \times [2^{h},2^{h+1}]$ with vertices at integer heights, delete the cells in $|\mathcal{U}_{h}| \times \{2^h\}$, and redefine the characteristic maps to be compatible with $\phi_{h}$. The cells  $e$ in $M$ have the form $\sigma \times \{t\}$ or $\sigma \times I$ where $\sigma$ is a closed simplex in some $|\mathcal{U}_{h}|$. We have that $C(e) = e$ (topologically) since each closed cell is a subcomplex. The cover by closed cells is locally finite, because each $|\mathcal{U}_{h}| \times [0,1]$ is finite. Therefore $M$ is strongly locally finite. 
\end{example}

\begin{definition} Let $Y$ be a strongly locally finite, path-connected CW complex with continuous proper base ray $\omega: [1,\infty)\rightarrow Y$. For $n\geq 1$ we define the \textit{$n$-th end homotopy group} $\pi_n^e(Y,\omega)$ to be the set of relative, continuous proper homotopy classes of continous proper maps 
\begin{align*}
f: ([0,1]^n \times [1,\infty), \partial [0,1]^n \times [1,\infty))\rightarrow (Y,\omega[1,\infty))
\end{align*}
such that $f_{|\partial [0,1]^n \times [1,\infty)} = \omega\circ p$, where $p: [0,1]^n \times [1,\infty)\rightarrow [1,\infty); (x,h)\mapsto h$. \\

We define the \textit{$0$-th end homotopy set} $\pi^e_0(Y)$ to be the set of continous proper homotopy classes of continous proper maps from $[1,\infty)$ to $Y$. 
\end{definition}

End homotopy groups satisfy all the expected properties: the existence of a group multiplication, functoriality under continuous proper maps with continous, properly homotopic maps inducing the same morphism; relative end homotopy groups are defined in the obvious way, and there is a long exact sequence for pairs $(Y,A)$ with continous proper map $k: A\rightarrow X$ and continous proper base ray $\omega: [1,\infty)\rightarrow A$. There is a change of base ray homomorphism, and the results of Proposition \ref{baseddefined} apply if one exchanges "coarse map" with "continous, proper map" and "coarse homotopy" with "continous, proper homotopy". \\

In the literature (Page $418$ of \cite{geoghegan2007topological}), $\pi_n^e(Y,\omega)$ is called the $n$-th strong homotopy group, or the $n$-th Steenrod homotopy group. The name "strong" refers to strong shape theory, which is discussed in Subsection \ref{strongshapetheory}. I choose to use the name "end homotopy groups" because they fit inside a $\varprojlim^1$ sequence, and $\varprojlim^1$ and $\varprojlim$ of a inverse system of groups does not depend on the first $j$  elements for any $j$ finite. \\

The pair $(c[0,1]^n, c\partial [0,1]^n)$ can be identified topologically with $([0,1]^n \times [1,\infty), \partial [0,1]^n \times [1,\infty))$.  If $[f]=[g]\in \pi_{n}^{L,e}(Y,\omega)$, a coarse-Lipschitz homotopy between them
\begin{align*}
H:(c[0,1]^{n+1},c(\partial[0,1]^n \times [0,1]))\rightarrow (Y,\omega)\\
H_{|c([0,1]^n\times \{0\})}=f \quad H_{|c([0,1]^n\times \{1\})}= g
\end{align*}
defines a proper homotopy just by scaling, ie. $G: (c[0,1]^n \times [0,1], c\partial[0,1]^n \times [0,1]) \rightarrow (Y,\omega)$; $(hx,h,t)\mapsto H(hx,h,ht)$. Therefore, every coarse-Lipschitz class $[f]\in \pi_n^{L,e}(Y,\omega)$ determines a class $[f]\in\pi^e_{n}(Y,\omega)$ by forgetting the Lipschitz property. We denote the forgetful map by $\mathfrak{d}:  \pi_n^{L,e}(Y,\omega)\rightarrow \pi^e_{n}(Y,\omega)$. We show that $\mathfrak{d}$ is an isomorphism for $Y=M$ the inverse mapping telescope.\\

Surjectivity: let $[f]\in \pi_{n}^{e}(M,\omega)$. The outline of the proof is to homotope $f$ to be locally Lipschitz via (geometric) simplicial approximation, and then find a globally Lipschitz representative by rescaling on the cone $c[0,1]^n$. For this we need some auxillary definitions and theorems: 

\begin{definition} \label{derivedcomplex} Let $K$ be a finite simplicial complex, $L$ a subcomplex of $K$. Let $(K \,rel \,L)'$ denote the barycentric derived complex of $K$ relative to $L$, which is obtained from $K$ by subdividing barycentrically all simplices of $K-L$. In particular $(K \,rel \,L)'$ contains $L$ as a subcomplex. Therefore we can define inductively 
\begin{align*}
K_0 &= K\\
K_r &= (K_{r-1} \,rel \,L)'
\end{align*}
\end{definition}

$K-L$ here refers to all simplices in $K$ not contained in $L$. For a more combinatoric description of $K_{r}$, see Subsection \ref{relgeosimplsect}. Adapting the proof in \cite{zeeman1964relative} to geometric simplicial approximation we obtain the following theorem:

\begin{theorem}\label{relgeosimpl} Let $K$ be a finite simplicial complex, $L$ a subcomplex of $K$. Let $M$ be the inverse mapping telescope of a compact metric space, and $f: |K|\rightarrow M$ a continuous map such that the restriction $f|_{L}$ is a geometrically simplicial map from $L$ to $M$. Then there exists an integer $s$ and a geometrically simplicial map $g: K_{s}\rightarrow M$ such that $g_{|L} = f_{|L}$ and $g$ is homotopic to $f$ keeping $L$ fixed.\\

Additionally, given an arbitrary neighbourhood $U$ of $|L|$ in $|K|$, we can choose $g$ and the homotopy such that
\begin{enumerate}
\item If $x\notin U$, the homotopy of $x$ is a path from $f(x)$ to $g(x)$ with image contained within the neighbourhood $ \mathcal{N}$ of a geometric simplex.
\item If $x\in |L|$, the homotopy leaves $x$ fixed at $f(x)$.
\item If $x\in U - |L|$, the homotopy of $x$ is contained in $gst^2(\mathscr{C},M)\cup \mathcal{N}$, where $\mathscr{C}$ is the core of a geometric simplex, $gst^2(\mathscr{C},M) := \cup_{v\in vert \mathscr{C}} \langle gst(v) \rangle$, and $ \mathcal{N}$ is the neighbourhood of a geometric simplex.
\end{enumerate}
\end{theorem}

An analogous statement applies if $K$ is uniformly bounded and $f$ is Lipschitz (see the proof of point $4$ of Proposition \ref{geosa}). The proof of Theorem \ref{relgeosimpl} is contained in Subsection \ref{relgeosimplsect}. \\

Now we go back to the proof of surjectivity. Consider the set $c_{[1,2]}[0,1]^n$. Applying the simplicial approximation theorem relative to $c\partial [0,1]^n$, there exists a homotopy $\overline{H}_1$ on $c_{[1,2]}[0,1]^n$ from $f$ to a simplicial map $\overline{f_1}$ such that $\overline{f_1}_{|c\partial [0,1]^n}= f_{|c\partial [0,1]^n}$. Consider the set $c_{[2, 3]}[0,1]^n$ and the subcomplex 
\begin{align*}
(2[0,1]^n \times \{2\})\cup (3[0,1]^n \times \{3\}) \cup (\partial [0,1]^n \times [2,3])
\end{align*}
Since subcomplexes of CW complexes satisfy the homotopy extension property, there is an extension of $\overline{H}_1$ to $c_{[2,3]}[0,1]^n$ which restricts to $\overline{H}_1$ on $(2[0,1]^n \times \{2\})$ and is constant on $ (3[0,1]^n \times \{3\}) \cup (\partial [0,1]^n \times [2,3])$. Extending this via the constant homotopy on the rest of $c[0,1]^n$, we obtain a homotopy $H_1: c[0,1]^n\times [0,1]\rightarrow M$ to a map $f_1$, which is simplicial on $c_{[1,2]}[0,1]^n \cup c\partial [0,1]^n$.\\

Suppose that $H_1\ast \dots \ast H_{k}$ has already been defined and the endpoint of the homotopy $f_k: (c[0,1]^n, c\partial [0,1]^n)\rightarrow (M,\omega)$ is simplicial when restricted to $c_{[1,2^{k}]}[0,1]^n\cup c\partial[0,1]^n$. By relative simplicial approximation, we have a homotopy $\overline{H}_{k+1}: c_{[2^{k},2^{k+1}]}[0,1]^n \times [0,1]\rightarrow M$ relative to $2^{k} [0,1]^n \times \{2^{k}\}\cup c\partial [0,1]^n$ and a simplicial $\overline{f}_{k+1}$. As before, $c_{[2^{k+1},2^{k+1}+1]}[0,1]^n$ has $(2^{k+1}[0,1]^n \times \{2^{k+1}\})\cup ((2^{k+1}+1)[0,1]^n \times \{2^{k+1}+1\}) \cup (\partial [0,1]^n \times [2^{k+1}, 2^{k+1}+1])$ as a subcomplex, so we extend $\overline{H}_{k+1}$ to $c_{[2^{k+1},2^{k+1}+1]}[0,1]^n$ and then to all of $c[0,1]^n$ via the constant homotopy.\\

The infinite concatenation of homotopies $H:= H_{1}\ast \dots \ast H_{k}  \ast H_{k+1}\ast \dots$ is well-defined since any point in $c[0,1]^n$ is only affected by at most $2$ of the $H_k$. Call $\tilde{f}$ the endpoint of the homotopy $H$. It is the same as $f$ on $c \partial[0,1]^n$ and is Lipschitz when restricted to any $c_{[2^{k}, 2^{k+1}]}[0,1]^n$ with Lipschitz constant $a_{k}$. By the properness of $f$, for every $h\in \mathbb{N}$ there exists a $l_{h}\in \mathbb{N}$ such that $f(c_{2^{l_h}}[0,1]^n)\subset M_{[h,\infty)}$. It is an easy exercise to check that $H(c_{2^{l_h}}([0,1]^n \times [0,1])) \subset M_{[h,\infty)}$ for all $h\in \mathbb{N}$. Therefore the homotopy is proper. \\

It should be noted that the resulting $\tilde{f}$ is simplicial with respect to a simplicial structure that gets progressively subdivided the further we go up the cone. \\

We now find a globally Lipschitz representative for $\tilde{f}$. Choose constants $m_h$ such that $2^{m_h} > a_h$ and $m_h \leq m_{h+1}$ for all $h\in \mathbb{N}_{0}$. The projection
\begin{align*}
q_{h}: [0,1]^n \times [2^{h+m_{h}}, 2^{h+1+m_{h}}] &\rightarrow [0,1]^n \times [2^{h}, 2^{h+1}]\\
 (tx,t)&\mapsto (\frac{t}{2^{m_{h}}}x, \frac{t}{2^{m_{h}}})
\end{align*}
has Lipschitz constant $\frac{1}{2^{m_{h}}}$, so $\tilde{f}\circ q_h$ is $1$-Lipschitz. We define $f'$ to be $\tilde{f}\circ q_h$ wherever it is defined, and since $q_{h+1}([0,1]^n \times \{2^{h+1+m_{h+1}}\}) = q_{h}([0,1]^n \times \{2^{h+1+m_h}\})$ we can extend this in the obvious way to a globally $1$-Lipschitz map on $c[0,1]^n$. Unfortunately $f'$ does not preserve the base ray: let $\omega'$ be the linear extension of the function $2^{h+m_{h}}\mapsto 2^{h}, 2^{h+1+m_{h}}\mapsto 2^{h+1}$ to $[1, \infty)\rightarrow [1,\infty)$. By construction $[f']\in \pi_{n}^{L,e}(M,\omega')$. Since the distance between $\omega'(h)$ and $\omega(h)$ is at most $h$, the linear homotopy 
\begin{align*}
G:I_{p}([1,\infty))\rightarrow [1,\infty)
\end{align*}
between $\omega'$ and $\omega$ is $2$-Lipschitz. Hence, $b_{\omega',\omega}[f'] \in \pi_{n}^{L,e}(M,\omega)$ and $b_{\omega',\omega}(f')$ is properly homotopic to $\tilde{f}$. Therefore, $\mathfrak{d}$ is surjective. \\

We also have the following lemma, to be useful in the proof of injectivity:
\begin{lemma} \label{clip} Suppose that $\tilde{f}$ is already Lipschitz. Then $f'$ is coarse-Lipschitz homotopic to $\tilde{f}$. 
\end{lemma}
\begin{proof} We define the homotopy as follows:
\begin{align*}
H:I_{p}(c[0,1]^n)&\rightarrow M\\
(hx,h,ht)&\mapsto \tilde{f}([th+(1-t)\omega'(h)]x, th+(1-t)\omega'(h))
\end{align*}
Let $L$ be the Lipschitz constant for $\tilde{f}$ and $(hx,h,ht), (sy,s,sl)\in I_{p}(c[0,1]^n)$. We compute
\begin{align*}
d(H(hx,h,ht),H(sy,s,sl)) \leq L(\|[ht+(1-t)\omega'(h)]x-[sl+(1-l)\omega'(s)]y\| \\
+ |ht+(1-t)\omega'(h)-sl-(1-l)\omega'(s)|)\\
\|[ht+(1-t)\omega'(h)]x-[sl+(1-l)\omega'(s)]y\| \leq [ht-(1-t)\omega'(h)]\|x-y\|\\
+\|ht+(1-t)\omega'(h)y-sl-(1-l)\omega'(s)y\|\\
\leq h\|x-y\| + |ht+(1-t)\omega'(h)-sl-(1-l)\omega'(s)|\|y\|
\end{align*}
\begin{align*}
|ht+(1-t)\omega'(h)-sl-(1-l)\omega'(s)|\leq |ht-sl| + |(1-t)\omega'(h)-(1-t)\omega'(s)| \\
+ |(1-t)\omega'(s)- (1-l)\omega'(s)|\\
\leq |ht-sl| + |\omega'(h)-\omega'(s)|+|t-l||\omega'(s)| \leq  |ht-sl| + |h-s|+s|t-l|\\
\leq  2|ht-sl| + 2|h-s|\\
d(H(hx,h,ht),H(sy,s,sl)) \leq L(\|hx-sy\|+5|h-s| + 4|ht-sl|)
\end{align*}
\end{proof}

\begin{remark} This lemma shows that if $f : (c[0,1]^n,\omega) \rightarrow (M,\omega_{f})$ and $g: (c[0,1]^n,\omega) \rightarrow (M,\omega_{g})$ are Lipschitz maps that are locally Lipschitz homotopic, such that the restriction of the homotopy to $ c \partial[0,1]^n$ is a Lipschitz homotopy $b_{\omega_{f},\omega_{g}}$, then $f$ and $g$ are coarse-Lipschitz homotopic, via a homotopy that restricts to $b_{\omega_{f},\omega_{g}}$ on the boundary. This is because
\begin{align*}
f\simeq b_{\omega_{f}',\omega_{f}} f' \simeq b_{\omega_{g}',\omega_{g}} g' \simeq g
\end{align*}
where $\simeq$ denotes coarse-Lipschitz homotopic.\\
\end{remark}

To show injectivity, suppose that $f,g$ are Lipschitz maps which are properly homotopic. Since the covering of open geometric stars in $M$ has Lebesgue number $\frac{1}{4}>0$ and $c[0,1]^n$ is uniformly bounded, the requirements for relative geometric simplicial approximation is met. Therefore, there exist integers $q\geq 1$ (which can be chosen freely), $s>q+1$, and simplicial maps $\mathfrak{h}: (c[0,1]^n)'_{q}\rightarrow (c[0,1]^n)_{q}$ and $\tilde{f},\tilde{g}: (c[0,1]^n)_{s}\rightarrow M$ such that: $\tilde{f}_{|c \partial[0,1]^n} = \tilde{g}_{|c \partial[0,1]^n} = \omega$, $\tilde{f},\tilde{g}$ are locally Lipschitz homotopic to $f \mathfrak{h},g \mathfrak{h}$ rel $c\partial [0,1]^n$, and $\mathfrak{h}$ is Lipschitz homotopic to $\id$ rel  $c\partial [0,1]^n$.\\

We first show that $[f] = [\tilde{f}]$ by constructing a coarse-Lipschitz homotopy between them. To that end, there exists a locally Lipschitz homotopy $G: c[0,1]^{n+1}\rightarrow M$ between $f,\tilde{f}$ which is the base ray $\omega$ on $G_{|c (\partial [0,1]^n \times [0,1])}$:
\begin{align*}
\tilde{f} \simeq f\mathfrak{h} \simeq f
\end{align*}
where all homotopies are rel. $c \partial[0,1]^n$. By the remark above, we can scale to a globally Lipschitz homotopy $G'$ between $f'$ and $\tilde{f'}$, where $G'_{|c (\partial [0,1]^n \times [0,1])} = \omega'$. Therefore
\begin{align*}
[f] = b_{\omega',\omega} [f'] = b_{\omega',\omega}[\tilde{f}'] = [\tilde{f}]
\end{align*}
where we have used that both $f$ and $\tilde{f}$ are Lipschitz. As before, there is a proper homotopy between $\tilde{f}$ and $\tilde{g}$ fixing the boundary, ie.
\begin{align*}
H: c[0,1]^{n+1} & \rightarrow M\\
H_{|c([0,1]^n \times \{0\})} = \tilde{f}, &\quad H_{|c([0,1]^n \times \{1\})} = \tilde{g}
\end{align*}
The homotopy restricted to $c(\partial[0,1]^n \times [0,1])$ is fixed to be $\omega$. We can choose a simplicial structure on $c[0,1]^{n+1}$ so that $H_{|c([0,1]^n \times\partial [0,1])}$ is simplicial (take the $s$-th barycentric subdivision of $c[0,1]^{n+1}$ rel $c(\partial [0,1]^n \times [0,1])$). As in the proof of surjectivity, we do simplicial approximation relative to $c(\partial [0,1]^{n+1})=c([0,1]^n \times \partial [0,1]) \cup c(\partial[0,1]^n \times [0,1])$ to obtain a locally Lipschitz homotopy $\tilde{H}$ then scale to a globally Lipschitz homotopy $H''$ between $\tilde{f}''$ and $\tilde{g}''$. This means that $b_{\omega'',\omega}[\tilde{f}''] = b_{\omega'',\omega} [\tilde{g}'']$. Putting it all together, we obtain
\begin{align*}
[f] = [\tilde{f}] = b_{\omega'',\omega}[\tilde{f}''] = b_{\omega'',\omega} [\tilde{g}''] = [\tilde{g}] = [g]
\end{align*}
which gives us injectivity of $\mathfrak{d}$. 
\subsection{$\varprojlim^1$ sequence}

\begin{definition} Let $Y$ be a strongly locally finite CW complex which is path-connected. A \textit{finite filtration} $\mathcal{L} = \{L_{i}\}_{i\in \mathbb{N}}$ of $Y$ consists of a sequence of finite, full subcomplexes $L_{1}\subset L_{2}\subset \dots$ such that $Y = \cup_{i\in \mathbb{N}} L_{i}$. A base ray $\omega: [1,\infty)\rightarrow Y$ is \textit{well-parametrised} with respect to $\mathcal{L}$ if $\omega(i) \subset Y \oset{c}{-} L_i$ for all $i$, where $Y \oset{c}{-} L_i$ refers to the CW-complement of $L_i$, ie. the largest subcomplex of $Y$ which has $0$-skeleton the vertices of $Y\setminus L_i$. 
\end{definition}

\begin{theorem*}(Theorem \ref{lim1absolutegeq1}) Let $Y$ be a strongly locally finite CW complex which is path-connected. If $\mathcal{L} = \{L_i\}_{i\in \mathbb{N}}$ is a finite filtration of $Y$ and $\omega$ is well-parametrised wrt $\mathcal{L}$, then there is a natural short exact sequence of groups
\begin{align*}
0\rightarrow {\varprojlim}^1  \pi_{n+1}(Y \oset{c}{-} L_i,\omega(i))  \xrightarrow{\overline{a}} \pi^{e}_{n}(Y,\omega) \xrightarrow{\overline{b}}  \varprojlim \pi_{n}(Y \oset{c}{-} L_i,\omega(i)) \rightarrow 0
\end{align*} 
for $n\geq 1$. 
\end{theorem*}

This theorem will be proven in Part II.\\

Consider the finite filtration $\{M_{[1,2^{h}]}\}_{h\in \mathbb{N}}$ of $M$. The CW complement of $M_{[1,2^{h}]}$ is $M_{[2^{h}+1,\infty)}=: M_{h}$. The base ray $\omega'(t) = \omega(2^{t}+1)$ is well-parametrised with respect to $\{M_{[1,2^{h}]}\}_{h\in \mathbb{N}}$ and properly homotopic to $\omega$. We have that $\pi_n^{e}(M,\omega)\cong \pi^e_{n}(M,\omega')$. This gives us the short exact sequence

\begin{figure}[H]
\center
\begin{tikzcd} [sep = .5 cm]
0 \arrow [r] & {\varprojlim}_{h}^1  \pi_{n+1}(M_h,\omega(2^h+1))  \arrow [r, "\overline{a}"]  & \pi^{e}_{n}(M,\omega) \arrow [r, "\overline{b}"]& \varprojlim_h \pi_{n}(M_h,\omega(2^h+1))\arrow [r] &0 
\end{tikzcd}
\end{figure}

\begin{lemma} There are isomorphisms 
\begin{align*}
{\varprojlim}_{h}^1  \pi_{n+1}(M_h, \omega(2^h+1)) &\cong {\varprojlim}_{h}^1   \pi_{n+1}(|\mathcal{U}_{h}|,x_0)\\
{\varprojlim}_{h}  \pi_{n}(M_h, \omega(2^h+1)) &\cong {\varprojlim}_{h}   \pi_{n}(|\mathcal{U}_{h}|,x_0)
\end{align*}
for all $n\geq 1$. 
\end{lemma}

\begin{proof} We define isomorphisms $\alpha^{h}: \pi_{n}(M_h,\omega(2^h+1))\rightarrow \pi_{n}(|\mathcal{U}_{h}|,x_0)$ as follows. The projection $\mathfrak{q}_{2^{h}+1}: M_{h}\rightarrow |\mathcal{U}_{h}| \times \{2^h+1\}$ is a deformation retract, hence a homotopy equivalence. Let $\alpha^h[f]$ be the image of $[f]$ under $\mathfrak{q}_{2^{h}+1}$ in $\pi_{n}(|\mathcal{U}_{h}| \times \{2^h+1\}, \{x_0\} \times \{2^h+1\}) \cong \pi_{n}(|\mathcal{U}_{h}|, \{x_0\})$. \\

The diagram

\begin{figure}[H]
\center
\begin{tikzcd}
&\pi_{n}(M_{h+1},\omega(2^{h+1}+1)) \arrow{r}{\alpha^{h+1}}  \arrow{d}{b_{\eta_{h+1,h}}}
&\pi_{n}(|\mathcal{U}_{h+1}|,x_0) \arrow{d}{\phi_{h+1}} 
  \\
&\pi_{n}(M_{h},\omega(2^{h}+1))\arrow {r}{\alpha^h}
&\pi_{n}(|\mathcal{U}_{h}|,x_0)
\end{tikzcd}
\end{figure}

commutes for all $h\in \mathbb{N}$ where $b_{\eta_{h+1,h}}$ denotes the change of base point homomorphism via the linear path $\eta_{h+1,h}(t)= (x_0,t(2^h+1)+(1-t)(2^{h+1}+1))$. To see the commutativity, observe that a deformation retract of $|\mathcal{U}_{h+1}| \times \{2^{h+1}+1\} \subset M_{h+1} \subset M_{h}$ to $|\mathcal{U}_{h}| \times \{2^h+1\}$ exactly corresponds to the map $\phi_{h+1}: |\mathcal{U}_{h+1}| \rightarrow |\mathcal{U}_{h}|$ with the obvious identifications.\\

So the two directed systems $\{\pi_n(M_{h},\omega(2^h+1)), b_{\eta_{h+1,h}}; \mathbb{N} \}$ and $\{\pi_{n}(|\mathcal{U}_{h}|,x_0), \phi_{h+1}; \mathbb{N}\}$ are pro-isomorphic, hence their $\varprojlim_{h}$ and $\varprojlim^1_{h}$ are isomorphic. Denote for future use the isomorphisms as 
\begin{align*}
c:{\varprojlim_h}^1  \pi_{n+1}(M_h, \omega(2^h+1)) &\rightarrow {\varprojlim_h}^1   \pi_{n+1}(|\mathcal{U}_{h}|,x_0) \\
d: \varprojlim_{h} \pi_{n}(M_{h},\omega(2^h+1)) &\rightarrow \varprojlim_{h} \pi_{n}(|\mathcal{U}_{h}|,x_0)
\end{align*}
\end{proof}

\begin{proof} Finally we actually get to the proof of the short exact sequence. As a reminder, we want the following
\begin{align*}
0\rightarrow {\varprojlim_h}^1  \pi_{n+1}(|\mathcal{U}_{h}|,x_0) \rightarrow \pi_{n}^{c}(cX,\omega) \rightarrow  \check{\pi}_n(X,x_0) \rightarrow 0
\end{align*}

The previous lemma and theorem give us a short exact sequence
\begin{align*}
0\rightarrow {\varprojlim_h}^1  \pi_{n+1}(|\mathcal{U}_{h}|,x_0) \xrightarrow{\overline{a} \circ c^{-1}}\pi_{n}^{e}(M,\omega) \xrightarrow{d\circ \overline{b}}\varprojlim_{h} \pi_{n}(|\mathcal{U}_{h}|,x_0) \rightarrow 0
\end{align*}
From previous results we have $\pi^{e}_{n}(M,\omega) \cong \pi^{L,e}_{n}(M,\omega) \cong \pi_{n}^{c}(cX,\omega)$ and $\varprojlim_{h} \pi_{n}(|\mathcal{U}_{h}|,x_0) \cong \check{\pi}_n(X,x_0)$, which concludes the proof. 
\end{proof}

\subsection{$\pi_0$ case}

This subsection deals with the $\pi_0$ case, which is analogous to $n\geq 1$ but easier because there are no change of base ray homomorphisms. Recall that $\pi_0^{c}(cX,\omega)$ is the pointed set of coarse homotopy classes of maps from $[1,\infty)$ to $cX$ with base point $[\omega] = [c\{x_0\}]$ and that $\pi_0^{L,e}(M,\omega)$ (resp. $\pi_{0}^{e}(M,\omega)$) is the set of proper coarse-Lipschitz homotopy classes of proper Lipschitz maps (resp. continuous proper homotopy classes of continuous proper maps) from $[1,\infty)$ to $M$ with base point $[\omega]$.

\begin{theorem} Let $X$ be a compact metric space with base point $x_0\in X$. Let $\omega$ be the standard parametrisation of the base ray $c\{x_0\}$. There is an isomorphism of pointed sets
\begin{align*}
 \pi_{0}^c(cX,\omega) \cong \pi^{L,e}_{0}(M,\omega) 
\end{align*}
\end{theorem}

\begin{proof} The map $\lambda: \pi_0^{L,e}(M)\rightarrow \pi_0^c(cX)$ is the same as before, ie $\lambda[\tau] = [R\tau]$. For $\theta$, choose the simplicial structure on $[1,\infty)$ with vertices at integer heights. As before, $[\tau]\in \pi_0^c(cX)$ has a simplicial approximation $\widetilde{i s_{r} \tau}\in \pi_n^{L,e}(M)$ for sufficiently large $r$. The independence of choice of simplicial approximation, independence of subdivision, independence of $r$, and independence of representative in the coarse homotopy class of $[\tau]$ follow by previous discussion. Recall that $is_{r}\tau$ is close to its geometric simplicial approximation and that $iR$ is close to the identity. So 
\begin{align*}
R\widetilde{is_{r}\tau} \simeq Ris_{r}\tau \simeq s_{r}\tau \simeq \tau
\end{align*}Hence $\lambda\circ \theta= \id$. For $\theta \circ \lambda$, choose a small enough barycentric subdivision of $[1,\infty)$, denoted by $bs^{k}([1,\infty))$, such that $\tau(|st(v)|)\subset \langle gst(z_v,t_v) \rangle$ for all vertices $v\in bs^{k}([1,\infty))$. We obtain 
\begin{align*}
(\widetilde{is_{r}R \tau})_{k}= (\widetilde{D_{r} \tau})_{k}\simeq D_{r}\tau \simeq \tau
\end{align*}
Hence $\theta\circ \lambda = \id$. We have that $\theta[\omega] = [\widetilde{i\omega}]=[\omega]$ and $[R\omega]=[\omega]$ so we obtain an isomorphism of pointed sets. 
\end{proof}

\begin{lemma} There is an isomorphism 
\begin{align*}
\pi_0^{L,e}(M,\omega)\cong \pi_0^{e}(M,\omega)
\end{align*}
\end{lemma}
\begin{proof} The map $\mathfrak{d}$ is as before. $[\tau]\in \pi_0^{e}(M)$ can be represented by a map $\tau: [1,\infty)\rightarrow M$. Choose the simplicial structure on $[1,\infty)$ to have vertices at integers. There exists a subdivision of $[1,\infty)$, a homotopy $\overline{H_1}$ on $[1,2]$ from $\tau$ to a simplicial map $\overline{\tau_1}$, and an extension $H_1$ of $\overline{H_1}$ to $[1,\infty)$ which restricts to $\overline{H_1}$ on $[1,2]$ and is the identity on $[3,\infty)$. Suppose that $H_1\ast \dots \ast H_{k}$ has already been defined and the endpoint of the homotopy $\tau_k: [1, \infty)\rightarrow M$ is simplicial when restricted to $[1,2^{k}]$. By relative simplicial approximation we have a homotopy $\overline{H}_{k+1}$ on $[2^{k},2^{k+1}]$ relative to $\{2^{k}\}$ and a simplicial $\overline{\tau}_{k+1}$. We extend $\overline{H}_{k+1}$ to a homotopy $H_{k+1}$ on $[1,\infty)$, keeping it the identity outside of $[2^{k},2^{k+1}+1]$. The homotopy $H:=H_1\ast \dots \ast H_{k}\ast H_{k+1}\ast \dots$ is well-defined. Call $\tilde{\tau}$ the endpoint of $H$. It is Lipschitz when restricted to any $[2^{k}, 2^{k+1}]$ with Lipschitz constant $a_{k}$. $H$ is proper since for every $h\in \mathbb{N}$ there exists a $l_h\in \mathbb{N}$ such that $H([2^{l_h},\infty)\times [0,1]) \subset M_{[h,\infty)}$.\\

As before we can find a globally Lipschitz representative for $\tilde{\tau}$ by choosing constants $m_h$ such that $2^{m_h}>a_h$ so that $\tilde{\tau}\circ q_h$ is $1$-Lipschitz, where $q_h: [2^{h+m_h}, 2^{h+1+m_h}]\rightarrow [2^{h},2^{h+1}]$ are the shrinking maps defined previously. By the same gluing procedure we extend $\tilde{\tau}\circ q_h$ to a globally $1$-Lipschitz map $\tau': [1,\infty)\rightarrow M$. This shows that $\mathfrak{d}$ is surjective.\\

For injectivity, suppose that $\tau$ and $\eta$ are Lipschitz maps which are properly homotopic. There exist simplicial approximations $\tilde{\tau},\tilde{\eta}$ with respect to a common barycentric subdivision. $\tilde{\tau},\tilde{\eta}$ are locally Lipschitz homotopic (and therefore coarse-Lipschitz homotopic, by rescaling) to $\tau$ and $\eta$ respectively. We have a proper homotopy
\begin{align*}
H: c[0,1]\rightarrow M\\
H|_{c\{0\}} = \tilde{\tau}, \quad H|_{c\{1\}} = \tilde{\eta}
\end{align*}
By relative simplicial approximation and rescaling, we obtain a Lipschitz homotopy $H'$ between $\tilde{\tau}'$ and $ \tilde{\eta}'$ ($\tilde{\tau}',\tilde{\eta}'$ are just $\tilde{\tau}$ and $\tilde{\eta}$ rescaled) and coarse-Lipschitz homotopies
\begin{align*}
\tau \simeq \tilde{\tau}\simeq \tilde{\tau}'\simeq \tilde{\eta}' \simeq \tilde{\eta} \simeq \eta
\end{align*}
\end{proof}

\begin{theorem} \label{lim10} (Prop $16.1.4$  in \cite{geoghegan2007topological}). Let $Y$ be a strongly locally finite path-connected CW complex, $[\omega] \in \pi_{0}^{e}(Y)$. If $\mathcal{L} = \{L_i\}$ is a finite filtration of $Y$ and $\omega$ is well-parametrised wrt $\mathcal{L}$, there is a natural short exact sequence of pointed sets
\begin{align*}
 {\varprojlim}^1  \pi_{1}(Y \oset{c}{-} L_i,\omega(i))  \xhookrightarrow{\overline{a}} \pi^{e}_{0}(Y,\omega) \xtwoheadrightarrow{\overline{b}}  \varprojlim \pi_{0}(Y \oset{c}{-} L_i,\omega(i)) 
\end{align*}
where $\overline{a}$ is injective and $\overline{b}$ is surjective. 
\end{theorem}

The proof is found in Chapter $16$ of \cite{geoghegan2007topological}. Applying this to $M$, we obtain the short exact sequence
\begin{align*}
{\varprojlim}_{h}^1  \pi_{1}(M_{h},\omega(2^h+1))  \xhookrightarrow{\overline{a}} \pi^{e}_{0}(M,\omega) \xtwoheadrightarrow{\overline{b}}   \pi_{0}(M_{h},\omega(2^h+1)) 
\end{align*}
The rest of the argument follows exactly as in the $n\geq 1$ case; all the maps and homotopies are the same, substituting $n=0$ and removing changes of base point. We obtain the $\varprojlim^1$-sequence in the $\pi_0$ case:

\begin{theorem} \label{maintheorem2} Let $X$ be a compact metric space, $x_0\in X$, $\omega$ the standard parametrisation of $c\{x_0\}$. There is a short exact sequence of pointed sets:
\begin{align*}
 {\varprojlim_h}^1  \pi_{1}(|\mathcal{U}_{h}|,x_0) \hookrightarrow  \pi^{c}_{0}(cX,\omega) \twoheadrightarrow \varprojlim_{h} \pi_0 (|\mathcal{U}_{h}|,x_0) 
\end{align*}
\end{theorem}
 
\newpage
\section{Applications}

In this section, we investigate corollaries of Theorem \ref{maintheorem} and compute some examples.\\

Let $X$ be a paracompact Hausdorff space. Let $O(X,x_{0})$ be the directed set of open covers of $X$ with a distinguished element containing $x_{0}$. It consists of pairs $(\mathcal{U},U_{0})$ where $\mathcal{U}$ is an open cover of $X$, $U_{0}\in \mathcal{U}$ and $x\in U_{0}$. We say that $(\mathcal{V}, V_{0})$ refines $(\mathcal{U},U_{0})$, denoted $(\mathcal{V}, V_{0}) \succeq (\mathcal{U},U_{0})$, if $\mathcal{V}$ refines $\mathcal{U}$ and $V_{0} \subset U_{0}$. \\

The nerve of a cover $(\mathcal{U}, U_{0})\in O(X,x_{0})$ is the abstract simplicial complex $N(\mathcal{U})$ whose vertex set is $\mathcal{U}$ and vertices $A_{0}, \dots, A_{k}\in \mathcal{U}$ span a $k$-simplex if $\cap_{i=0}^{k} A_{i}\neq \emptyset$. The vertex $U_{0}$, which we often denote as $x_{0}$, is taken to be the basepoint of the geometric realisation $|N(\mathcal{U})|$. Whenever $(\mathcal{V}, V_{0})$ refines $(\mathcal{U},U_{0})$, there is a simplicial map $\phi_{\mathcal{V},\mathcal{U}}: |N(\mathcal{V})|\rightarrow |N(\mathcal{U})|$ defined by sending a vertex $V\in \mathcal{V}$ to a vertex $U\in \mathcal{U}$ such that $V\subset U$. We make the convention that $\phi_{\mathcal{V},\mathcal{U}}(V_{0}) = U_{0}$. $\phi_{\mathcal{V},\mathcal{U}}$ is unique up to based homotopy. Therefore, the induced homomorphism on homotopy groups (resp. pointed sets)
\begin{align*}
\phi_{\mathcal{V},\mathcal{U}}: \pi_{n}(|N(\mathcal{V})|,V_{0}) \rightarrow \pi_{n}(|N(\mathcal{U})|,U_{0})
\end{align*}
is (up to coherent isomorphism) independent of the choice of simplicial map.

\begin{definition} \label{defcech} The \textit{$n$-th shape homotopy group} (resp. pointed set) is the inverse limit 
\begin{align*}
\check{\pi}_{n}(X,x_{0}) = \varprojlim_{\mathcal{U}\in O(X,x_{0})} \{\pi_{n}(|N(\mathcal{U})|, U_{0}), \phi_{\mathcal{U}\mathcal{V}}\}
\end{align*}
\end{definition}

If $X$ is a compact metric space with basepoint $x_{0}\in X$, then the sequence of metric open covers $\{(\mathcal{U}_{h},x_{0}), \phi_{h+1}; \mathbb{N}\}$ we constructed is cofinal in $O(X,x_{0})$. Therefore, $\check{\pi}_{n}(X,x_{0}) \cong \varprojlim_{h} \{\pi_{n}(|\mathcal{U}_{h}|, x_0), \phi_{h+1}\}$, as previously discussed.\\

More generally, consider an inverse sequence $\{(X_{k},x_{k}),\phi_{k+1};\mathbb{N}\}$ of spaces having the homotopy type of a simplicial complex and pointed maps $\phi_{k+1}: (X_{k+1},x_{k+1})\rightarrow (X_{k},x_{k})$, with the property that $(X,x_0) = \varprojlim_{k} (X_{k},x_{k})$ in pointed topological spaces. We will show later (see Subsection \ref{inversesc}) that 
\begin{align*}
\check{\pi}_{n}(X,x_{0}) &\cong \varprojlim_{k} \pi_{n}(X_{k}, x_k)\\
{\varprojlim_h}^1 \pi_{n+1}(|\mathcal{U}_{h}|,x_0) &\cong {\varprojlim_k}^1 \pi_{n+1}(X_{k}, x_k)
\end{align*}
for all $n\in \mathbb{N}_{0}$.

\subsection{Examples} \label{goodexamples}

We begin with a proof of Corollary \ref{good}. Let $X$ be a compact metric space with $x_{0}\in X$ and assume there exists a $N\in \mathbb{N}$ such that $\mathcal{U}_{h}$ is good for all $h>N$. We show that $X\simeq |\mathcal{U}_{h}|$ and that the bonding maps $\phi_{h+1}: |\mathcal{U}_{h+1}|\rightarrow |\mathcal{U}_{h}|$ are homotopy equivalences for $h>N$. With this, it is easy to see that $\varprojlim_{h} \pi_{n}(|\mathcal{U}_{h}|,x_0) \cong \pi_{n}(X,x_0)$ and that $\varprojlim^1_{h}\pi_{n+1}(|\mathcal{U}_{h}|,x_0) = 0$ for all $n\in \mathbb{N}_{0}$.

\begin{proof} (of Corollary \ref{good}) For the full definitions used in this proof, refer to Chapter $4.$G of \cite{hatcher2002algebraic}. We have the following diagram

\begin{figure}[H]
\center
\begin{tikzcd}
&X  \arrow{r}{s_{h+1}} \arrow{dr}{s_{h}}
&\Delta X_{\mathcal{U}_{h+1}} \arrow{d}{\varphi_{h+1}} \arrow{r}{a_{h+1}}
& \Gamma_{\mathcal{U}_{h+1}} = bs|\mathcal{U}_{h+1}| \arrow{d}{\phi_{h+1}}
  \\
&
&\Delta X_{\mathcal{U}_{h}} \arrow{r}{a_{h}}
& \Gamma_{\mathcal{U}_{h}} = bs|\mathcal{U}_{h}|
\end{tikzcd}
\end{figure}

The space $\Delta X_{\mathcal{U}_{h}}$ is the realisation of the complex of spaces associated to the open cover $\mathcal{U}_{h}$. It is the quotient space obtained from the disjoint union of all the products $U_{i_0}\cap \dots \cap U_{i_n} \times \Delta^n$, with subscripts ranging over sets of $n+1$ distinct indices with $n\geq 0$. There are identifications over the faces of $\Delta^n$ using the inclusions $U_{i_0}\cap \dots \cap U_{i_n} \hookrightarrow U_{i_0}\cap \dots \cap \widehat {U_{i_{j}}} \cap \dots \cap U_{i_n}$. \\

Let $\{\rho^{h}_{U}\}_{U\in \mathcal{U}}$ be a partition of unity subordinate to $\mathcal{U}_{h}$. The map $s_{h}: X\rightarrow \Delta X_{\mathcal{U}_{h}}$ is defined as
\begin{align*}
s_{h}(x) = (x,\sum_{U\in \mathcal{U}_{h}} \rho^{h}_{U}(x) [U])
\end{align*}
where $x$ is considered as a point in $\cap_{x\in U} U$. $\varphi_{h+1}$ is the product of the inclusion $\cap_{i=0}^{n} U_{i} \subset \cap_{i=0}^{n} \phi_{h+1}(U_{i})$ and the simplicial map $\phi_{h+1}:[U_{0},\dots,U_{n}] \rightarrow[\phi_{h+1}(U_{0}),\dots,\phi_{h+1}(U_{n})]$. $a_{h}$ is induced by sending each intersection to a point while keeping the barycentric coordinates. $bs|\mathcal{U}_{h}|$ is the barycentric subdivision of $|\mathcal{U}_{h}|$ (Hatcher emphasises the subdivided structure to construct a homotopy inverse to $a_{h}$). Observe that $a_{h}\circ s_{h}: X\rightarrow bs|\mathcal{U}_{h}|=|\mathcal{U}_{h}|$ is just the standard map into the nerve induced by a partition of unity subordinate to $\mathcal{U}_{h}$. \\

We now show that this diagram commutes up to homotopy. The left triangle commutes up to homotopy because 
\begin{align*}
\varphi_{h+1} \circ s_{h+1}(x) &= (x,\sum_{V\in \mathcal{U}_{h+1}} \rho^{h+1}_{V}(x) [\phi_{h+1} V])
\end{align*}
Since $\varphi_{h+1} \circ s_{h+1}(x)$ and $s_{h}(x)$ lie in a common simplex, they can be joined by a linear homotopy. \\

The rectangle commutes by following the definitions. The maps $s_{h+1}$ and $s_{h}$ are homotopy equivalences by Proposition 4G.2 and the maps $a_{h+1}$ and $a_{h}$ are homotopy equivalences by Corollary 4G.3. in \cite{hatcher2002algebraic}. Therefore $\varphi_{h+1}$ and $\phi_{h+1}$ are also homotopy equivalences for $h>N$, and we have
\begin{align*}
\pi_{n}^{c}(cX,\omega)\cong \pi_{n}(X,x_0)
\end{align*}
for all $n\in \mathbb{N}_{0}$, as claimed. 
\end{proof}

\begin{remark} The conditions of Corollary \ref{good} are satisfied for $X$ a compact manifold or a finite simplicial complex. This recovers the results of \cite{mitchener2020coarse}. 
\end{remark}

Now we look at some interesting non-trivial examples:

\begin{enumerate}
\item \textbf{Cantor set}: The Cantor set $\mathcal{C}$ is obtained by iteratively deleting the open middle third from a set of line segments, ie. let
\begin{align*}
C_{0} &:= [0,1]\\
C_{n} &:= \frac{C_{n-1}}{3} \cup (\frac{2}{3} +\frac{C_{n-1}}{3}) \quad n\geq 1
\end{align*}
and define $\mathcal{C}:= \cap_{n=0}^{\infty} C_{n}$. This is a compact metric space with metric induced from $[0,1]$. Equivalently one can use the $p$-adic metric on $2^\mathbb{N}$: given two sequences $(x_n), (y_n)\in 2^\mathbb{N}$, the distance between them is $d((x_n),(y_n)) = 2^{-k}$, where $k$ is the smallest index such that $x_{k}\neq y_{k}$. \\

Let $x_0$ denote the sequence corresponding to $0$. $|\mathcal{U}_{h}|$ is a disjoint union of $2^{h}$ points, and is homotopy equivalent to $C_{h}$ ($2^h$ disjoint intervals). Under this identification, the bonding maps $\phi_{h+1}: C_{h+1} \rightarrow C_{h}$, up to homotopy, map each interval in $C_{h+1}$ to the interval in $C_{h}$ that contains it. Now it is easy to see that $\varprojlim_{h}\pi_{0}(|\mathcal{U}_{h}|,x_0) \cong \mathcal{C}$ and $\pi_{0}^{c}(c\mathcal{C},\omega) = \mathcal{C}$. All higher coarse homotopy groups vanish.
\item \textbf{Warsaw circle}: Recall that the Warsaw circle $S_{W}\subset \mathbb{R}^2$ is defined as:
\begin{align*}
S_{W} &:=  \{(x,sin(\frac{\pi}{x}))\,|\, 0<x<1\} \cup (\{0\} \times [-1,1]) \cup \\
& \{(x,y)\,|\, (x-\frac{1}{2})^2 +(y+\frac{15}{16})^2= (\frac{17}{16})^2, y\leq 0\}
\end{align*}
with metric induced from the ambient norm. $S_W$ is a connected, compact metric space weakly homotopy equivalent to a point. \\

Let $x_0\in S_{W}$. The Warsaw circle $(S_{W},x_0)$ is pointed shape equivalent to $(S^1, s)$. To see this, consider the inverse system $\{N_{k},\iota_{k+1,k};\mathbb{N}\}$ of $\frac{1}{k}$-neighbourhoods of $S_{W}$ in $\mathbb{R}^2$. By the discussion in Subsection \ref{inversesc}, this system is a pointed shape expansion of $(S_{W},x_0)$. It is easy to see that for sufficiently large $k$, the $\frac{1}{k}$-neighbourhood of $S_W$ is homotopy equivalent to a circle (see Figure \ref{fig:warsaw}), and the inclusion maps $\iota_{k+1,k}$ induce the identity on $\pi_{1}(S^1,s) \cong \mathbb{Z}, \pi_{0}(S^1,s) = \{[s]\}$. Therefore $\pi_{1}^{c}(cS_{W},\omega)\cong \mathbb{Z}$, $\pi_{0}^{c}(cS_{W},\omega)\cong \{[\omega]\}$ and all higher coarse homotopy groups vanish. 

\item \textbf{Dyadic solenoid}: The dyadic solenoid $\Lambda$ is a compact metrisable space which is the limit of the inverse sequence $\{X_{k}, p_{k+1};\mathbb{N}\}$ where each $X_{k}$ is a copy of $S^1$ and each bonding map $p_{k+1}: X_{k+1}\rightarrow X_{k}$ is given by the squaring map $e^{i\theta}\mapsto e^{2i\theta}$ (viewing $S^1$ as the unit circle in the complex plane). Consider the point $1\in \mathbb{C}$. We have that $(\Lambda, x_0) = \varprojlim_{k} (X_{k},1)$ for some $x_0\in \Lambda$. Therefore 
\begin{align*}
\check{\pi}_{n}(\Lambda,x_{0}) &\cong \varprojlim_{k} \pi_{n}(X_{k}, 1)\\
{\varprojlim_h}^1 \pi_{n+1}(|\mathcal{U}_{h}|,x_0) &\cong {\varprojlim_k}^1 \pi_{n+1}(X_{k}, 1)
\end{align*}
for all $n\in \mathbb{N}_{0}$. $\varprojlim_{k} \pi_{1}(X_{k}, 1)$ is the inverse limit of the sequence 
\begin{align*}
\mathbb{Z}\xleftarrow{\cdot 2} \mathbb{Z} \xleftarrow{\cdot 2} \mathbb{Z} \xleftarrow{\cdot 2} \mathbb{Z} \xleftarrow{\cdot 2} \dots
\end{align*}
where $\cdot 2$ denotes multiplication by $2$. Since there is no integer divisible by every power of $2$, we obtain $\varprojlim_{k} \pi_{1}(X_{k}, 1)=0$.\\

 Since $S^1$ is connected, we have $\varprojlim_{k} \pi_{0}(X_{k}, 1) = \{\ast\}$. However, $\pi_{0}^{c}(c\Lambda,\omega) \cong \varprojlim^1_{k}\pi_{1}(X_{k},1) \neq 0$. In fact, $\varprojlim^1_{k}\pi_{1}(X_{k},1)$ has cardinality $\mathbb{R}$ (Example $17.7.2$ in \cite{geoghegan2007topological}). Therefore

\[ \pi_{n}^{c}(c\Lambda,\omega) \cong \begin{cases} 
        0 & n\geq1 \\ 
        \mathbb{R} & n=0\\
       \end{cases}
    \]

\item \textbf{Shrinking wedge of spaces}: If $X = \vee_{j\in \mathbb{N}} X_{j}$ is a shrinking wedge of path-connected, compact, based spaces with wedge point $x_0$, we have that $(X,x_0) = \varprojlim_{k} \{(X_{\leq k},x_0), \phi_{k+1,k};\mathbb{N}\}$ in pointed topological spaces, where  $X_{\leq k} = \vee_{j=1}^{k}X_k$ and the bonding maps $\phi_{k+1,k}: X_{\leq k+1}\rightarrow X_{\leq k}$ collapse of the last space $X_{k+1}$ to $x_0$. The shape homotopy groups are
\begin{align*}
\check{\pi}_{n}(X,x_0) \cong \varprojlim_{k} \{\pi_{n}(X_{\leq k},x_0), \phi_{k+1,k}\}
\end{align*}
The inclusion map $\iota_{k,k+1}: X_{\leq k}\rightarrow X_{\leq k+1}$ induces a canonical section of $\phi_{k+1,k}: \pi_{n}(X_{\leq k+1})\rightarrow \pi_{n}(X_{\leq k})$. Therefore $\varprojlim^1 \{\pi_{n}(X_{\leq k},x_0), \phi_{k+1,k}\}$ is trivial for $n\geq 1$. This means that there is an isomorphism 
\begin{align*}
\pi_{n}^{c}(cX,\omega)\cong \check{\pi}_{n}(X,x_0)
\end{align*}
for $n\geq 0$.

\item \textbf{Hawaiian earring}: The Hawaiian earring $\mathbb{H}$ is defined as the infinite wedge of shrinking circles embedded in $\mathbb{R}^2$ with wedge point $x_0 = (0,0)$. 
\begin{align*}
\mathbb{H}:= \bigvee _{j\in \mathbb{N}} S^1_{j}
\end{align*}
By the previous example, we have that $\pi_{n}^{c}(c\mathbb{H},\omega)\cong \check{\pi}_{n}(\mathbb{H},x_0)$ for $n\in \mathbb{N}_{0}$. \\

$\varprojlim_{k} \{\pi_{n}(S^{1}_{\leq k},x_0), \phi_{k+1,k}\} = 0$ for $n\geq 2$. For $n=1$ it is the inverse limit of $\{F_{k},\phi_{k+1,k}\}$, with bonding maps 
\begin{align*}
\phi_{k+1,k}: F_{k+1}&\rightarrow F_{k}\\
w(a_1,\dots,a_k, a_{k+1}) &\mapsto w(a_1,\dots,a_k,1)
\end{align*} 
$F_{k}$ denotes the free group on the generating set $G_{k} = \{a_{1},\dots,a_{k}\}$ and $\phi_{k+1,k}$ takes a word in $G_{k+1}\cup G_{k+1}^{-1}$ and maps it to the same word with $a_{k+1}^{\pm 1}$ replaced by the identity. We obtain

\[ \pi_{n}^{c}(c\mathbb{H},\omega) \cong \begin{cases} 
        0 & n\geq 2 \\ 
 \varprojlim_k F_{k} & n= 1 \\
        \{\ast\} & n=0\\
       \end{cases}
    \]

\item  \textbf{Shrinking wedge of spheres}: The infinite earring space  $\mathbb{E}_m$ is a shrinking wedge of $m$-spheres ($m\geq 2$) with wedge point $x_0$. The shape homotopy groups are
\begin{align*}
\check{\pi}_n(\mathbb{E}_{m},x_0) \cong \bigoplus_{1\leq j\leq \frac{n-1}{m-1}} (\pi_{n}(S^{mj-j+1}))^{\mathbb{N}}
\end{align*} 
for $n \geq 2$ \cite{brazas2025v}.\\

Since $\mathbb{E}_{m} = \vee_{j\in \mathbb{N}} S^{m}_{j}$ satisfies the conditions of Example $4$ we have that the $\varprojlim^1$ term vanishes and $\pi_{n}^{c}(c\mathbb{E}_{m},\omega)\cong \check{\pi}_n(\mathbb{E}_{m},x_{0})$ for $n\geq 0$. Since $\pi_{1}(\vee_{j=1}^{k}S^m_{j})=0$ and  $\pi_{0}(\vee_{j=1}^{k}S^m_{j})=\{\ast\}$ for all $k\in \mathbb{N}$ we obtain

\[ \pi_{n}^{c}(c\mathbb{E}_{m},\omega) \cong \begin{cases} 
        \bigoplus_{i\leq j\leq \frac{n-1}{m-1}} (\pi_{n}(S^{mj-j+1}))^{\mathbb{N}} &n\geq 2  \\
        0 & n=1 \\ 
        \{\ast\} & n=0\\
       \end{cases}
    \]
\item \textbf{Strange wedge}: \textit{This example is a warning that the coarse homotopy groups of $cX$ can be different depending on the choice of $x_0\in X$.} Consider $X = \Lambda \vee S^1$, the wedge of the dyadic solenoid and a circle with wedge point $x_0$. Let $P = (S^1,\ast)\vee (S^1,\ast)$ and consider the map $f$ which takes the paths $\omega_1,\omega_2,\omega_3$ to the paths $\omega_1,\omega_2\omega_1\omega_2,\omega_3$ respectively. Similarly, let $g$ be the map taking the paths $\omega_1,\omega_2,\omega_3$ to $\omega_1\omega_2\omega_1,\omega_2,\omega_3$. See Figure \ref{fig:strangewedge}.

\begin{figure}[H]
\centering
  \centering
  \includegraphics[width=0.4\linewidth]{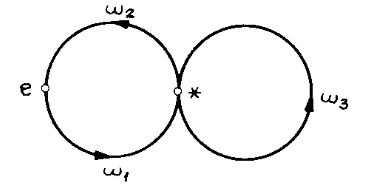}
  \caption{\cite{mardesic1982shape} The space $P = (S^1,\ast)\vee (S^1,\ast)$.}
  \label{fig:strangewedge}
\end{figure}

$(X,x_{0})$ is the inverse limit of the sequence 
\begin{align*}
(P,\ast)\xleftarrow{f} (P,\ast)\xleftarrow{g} (P,\ast) \xleftarrow{f} (P,\ast) \xleftarrow{g} \dots
\end{align*}
Analogously, the inverse limit of 
\begin{align*}
(P,e)\xleftarrow{f} (P,e)\xleftarrow{g} (P,e) \xleftarrow{f} (P,e) \xleftarrow{g} \dots
\end{align*}
is homeomorphic to $(X,x_1)$ for a suitable $x_1\in X$. We have that $\pi_{n}^c(cX,c\{x_0\})=\pi_{n}^c(cX,c\{x_1\})=0$ for $n\geq 2$. By definition $\pi_0^c(cX,c\{x_0\})\cong \pi_{n}^c(cX,c\{x_1\})$ since it is only a pointed set. However, it can be shown that
\begin{align*}
\pi_{1}^c(cX,c\{x_0\}) \cong \check{\pi}_{1}(X,x_0) \cong \mathbb{Z}\\
\pi_{1}^c(cX,c\{x_1\}) \cong \check{\pi}_{1}(X,x_1) =0
\end{align*}
See Example 4 of Ch. II §3.3 in \cite{mardesic1982shape} for the details of this computation.

\item \textbf{Random groups}: Let $F_n$ be the free group on $n\geq 2$ generators $s_{1},\dots, s_{n}$. For any integer $L$ let $R_{L}\subset F_{n}$ be the set of reduced words of length $L$ in these generators. 

\begin{definition} Let $0<d<1$. A random set of relators at density $d$, at length $L$ is a $\left \lfloor{(2n-1)^{dL}}\right\rfloor$-tuple of elements of $R_{L}$, randomly picked among all elements of $R_{L}$. A random group at density $d$, at length $L$ is the group $G$ presented by $\langle S|R \rangle $, where $S = \{s_1,\dots, s_n\}$ and $R$ is a random set of relators at density $d$, at length $L$. 
\end{definition} 

We say that a property occurs with overwhelming probability at density $d$ if for fixed $d$, the probability of occurence tends to $1$ as $L\rightarrow \infty$. 

\begin{theorem} (Section 9.B in \cite{gromov1992asymptotic}) A random group is with overwhelming probability 
\begin{enumerate}
\item Trivial or $\mathbb{Z}/2\mathbb{Z}$ at density greater than $\frac{1}{2}$.
\item Word-hyperbolic, with aspherical presentation complex, at density less than $\frac{1}{2}$.
\end{enumerate}
\end{theorem}

\begin{theorem} (Theorem $1.4$ in \cite{dahmani2011random}) Let $0<d<\frac{1}{2}$. Then with overwhelming probability the boundary of a random group at density $d$ is the Menger sponge. 
\end{theorem}

Combining this with Theorem $1.1$ of \cite{Fukaya_2018}, we obtain that the generic (finitely presented, infinite) group is coarsely homotopy equivalent to the cone of the Menger sponge $c\mathbb{M}$. The Menger sponge $\mathbb{M}$ is a topologically $1$-dimensional fractal constructed from the unit cube by subdividing it into $27$ smaller cubes, removing the central cube and the six cubes centred on each face; and then inductively repeating the same procedure with the $20$ remaining smaller cubes (see Figure  \ref{fig:Menger}).

\begin{figure}[H]
\centering
  \centering
  \includegraphics[width=0.5\linewidth]{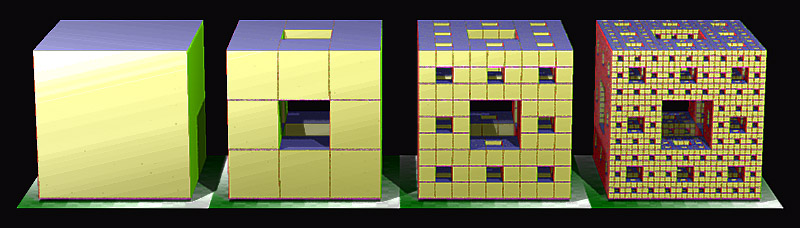}
  \caption{\cite{WikiCommonsSponge} The first four iterations of the construction of the Menger sponge.}
  \label{fig:Menger}
\end{figure}

We now follow the paper \cite{fischer2013word}. $(\mathbb{M},x_0)$ can be alternatively constructed as the inverse limit of a sequence of $4$-valent graphs $\{X_h,x_h\}_{h\in \mathbb{N}}$, where each $X_h$ has $2\cdot 4^h$ vertices and $x_h$ is the top-most vertex. See Figure \ref{fig:GraphsMenger}.

\begin{figure}[H]
\centering
  \centering
  \includegraphics[width=0.7\linewidth]{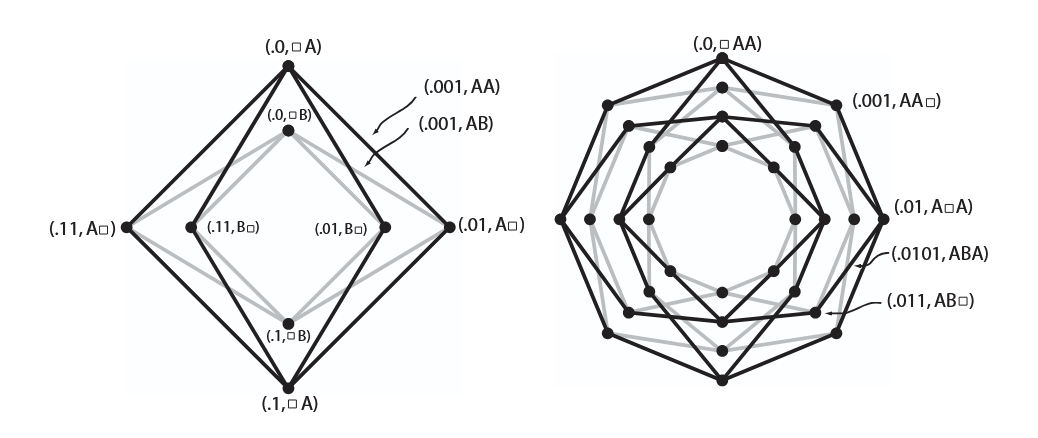}
  \caption{\cite{fischer2013word} The graphs $X_1$ and $X_2$.}
  \label{fig:GraphsMenger}
\end{figure}

 We have that $\pi_1(X_{h+1},x_{h+1}) = \pi_{1}(X_h,x_h) \ast F_{3}\ast \dots \ast F_{3}$  with one factor of $F_{3}$ for each vertex of $X_{h}$. The bonding morphisms $\pi_{1}(X_{h+1},x_{h+1})\rightarrow \pi_{1}(X_h,x_h)$ trivialise these additional factors. The surjectivity of the bonding morphisms imply the vanishing of $\varprojlim^1 \pi_{1}$, and all higher homotopy groups vanish. Therefore,  

\[ \pi_{n}^{c}(G_{ran},\omega) \cong \begin{cases} 
       0 &n\geq 2  \\
        {\varprojlim}_{h}  F_{c_{h}} & n=1 \\ 
        \{\ast\} & n=0\\
       \end{cases}
    \]
where $c_{h} =2\cdot 4^{h}+1$, and the bonding maps $\psi_{h+1,h}: F_{c_{h+1}} \rightarrow F_{c_{h}}$ are the canonical projection.

\end{enumerate}

\subsection{Relation to Spanier groups}

We relate the coarse homotopy groups $\pi_{n}^c(cX,\omega)$ the kernel of the group homomorphism $\Psi_{n}:\pi_{n}(X,x_0) \rightarrow \check{\pi}_{n}(X,x_0)$, which consists of elements of $\pi_{n}(X,x_0)$ which cannot be detected by polyhedral approximations to $X$. We will show that there exists a coning map $\mathfrak{c}$ such that the composition 
\begin{align*}
\pi_{n}(X,x_0)\xrightarrow{\mathfrak{c}}\pi_{n}^{c}(cX,\omega)\rightarrow \check{\pi}_{n}(X,x_0)
\end{align*}
where the second arrow comes from the $\varprojlim^1$ sequence, coincides with $\Psi_{n}$.\\

\begin{prop} \label{unicoarse} Let $X$ be a compact metric space, $x_0\in X$, $n\geq 1$. There is a canonical group homomorphism $\mathfrak{c}: \pi_{n}(X,x_0)\rightarrow \pi_{n}^{c}(cX,\omega)$ obtained by coning a representative of $[f]\in \pi_{n}(X,x_0)$. 
\end{prop}
\begin{proof} This is subtler than may first appear: it is not true in general that any continuous map $f: ([0,1]^n,\partial[0,1]^n)\rightarrow (X,x_0)$, where $X$ is a compact metric space, can be continuously homotoped to a Lipschitz map, so it may not even be possible to find a representative $f$ of $[f]\in \pi_{n}(X,x_0)$ such that the cone $cf: (c[0,1]^n,c\partial[0,1]^n)\rightarrow (cX,\omega)$ is coarse. \\

However, the fact that any continuous map between compact metric spaces is uniformly continuous rescues us.  That means there is a function $\beta: \mathbb{R}_{\geq 0} \rightarrow \mathbb{R}_{\geq 0}$ with the following properties:
\begin{itemize}
\item $\beta(0) = 0$. 
\item $\beta$ is monotonically increasing. 
\item If $x,y\in [0,1]^n$ with $d_{[0,1]^n}(x,y)\leq r$ then $d_{X}(f(x),f(y))\leq \beta(r)$. 
\item $\beta$ is continuous. 
\end{itemize}

We let
\begin{align*}
\beta(r) := \max_{d_{[0,1]^n}(x,y)\leq r} d_{X}(f(x),f(y))
\end{align*}
The first three points are obvious. Assume that $\beta$ is not left continuous. Then there exists an $r_0$ and $\varepsilon>0$ such that $\beta(r) < \beta(r_0) -\varepsilon$ for all $r<r_0$. The function
\begin{align*}
d_{X}(f\times f):[0,1]^n \times [0,1]^n &\rightarrow [0,\infty)\\
(x,y) &\mapsto d_{X}(f(x),f(y))
\end{align*}
is continuous and $\{(x,y)\in [0,1]^n \times [0,1]^n\,|\,d_{[0,1]^n}(x,y)\leq r_0\}$ is a compact set, so it assumes its maximum. This maximum must be for points $x_{0},y_{0}$ with distance equal to $r_0$, or it would contradict the assumption. Since $[0,1]^n$ is a geodesic path metric space, there exists a sequence $\{y_n\}_{n\in \mathbb{N}}$ with $d_{[0,1]^n}(x_0,y_n)<d_{[0,1]^n}(x_0,y_0)$ for all $n$, such that $\lim_{n\rightarrow \infty} d_{[0,1]^n}(x_0,y_n)= r_0$ and $\lim_{n\rightarrow \infty} d_{X}(f(x_0),f(y_n)) =  d_{X}(f(x_0), f(y_0))=\beta(r_0)$, contradiction. \\

Right continuity: the sets $K_{r}:= \{(x,y)\in [0,1]^n \times [0,1]^n \,|\, d(x,y)\leq r\}$ for $r\geq r_0$ are compact and $\cap_{r\geq r_0} K_{r} = K_{r_0}$. The function $d_{X}(f\times f)$ attains its maximum on each $K_{r}$: there exist $(x_r,y_r) \in K_{r}$ with $\beta(r) = d_{X}(f(x_r),f(y_r))$. By passing to a subsequence, we can assume that $\lim_{r\rightarrow r_{0}} (x_{r},y_{r}) = (x_0,y_0) \in \cap_{r\geq r_0} K_{r}= K_{r_0}$. We obtain that $\lim_{r\rightarrow r_{0}} \beta(r) = \beta(r_0)$. \\ 

For ease of notation, define
\begin{align*}
\gamma(t):= [\beta(\frac{1}{\log(t)+1})]^{-1} +1
\end{align*}

Observe that $\gamma: [1,\infty)\rightarrow [1,\infty)$ is a cotinuous, monotonically increasing function and that $\lim_{t\rightarrow \infty} \gamma(t) = \infty$. $\gamma$ might not be coarse. 

\begin{lemma} There exists a shrinking map $sh_{\gamma}: [1,\infty)\rightarrow [1,\infty)$ such that $\gamma\circ sh_{\gamma}$ is still continuous, monotonically increasing, and is coarse. 
\end{lemma}

\begin{proof} Let $h\in \mathbb{N}_{0}$. Consider the restriction $\gamma_{|[2^{h},2^{h+1}]}: [2^h,2^{h+1}]\rightarrow [1,\infty)$, which we denote by $\gamma_{h}$. Since $\gamma_{h}$ is uniformly continuous, there exists a constant $d_{h}$ such that  $|\gamma_{h}(t)-\gamma_{h}(s)|<1$ for $t,s\in [2^h,2^{h+1}]$ and $|t-s|<d_{h}$. Let $N_{h}$ denote the open $d_{h}$-neighbourhood of the diagonal. The function 
\begin{align*}
([2^{h}, 2^{h+1}] \times [2^{h}, 2^{h+1}])\setminus N_{h} &\longrightarrow [0,\infty)\\
(t,s)&\longmapsto \frac{|\gamma(t)-\gamma(s)|}{|t-s|}
\end{align*}
is continuous on a compact set, so is bounded by some constant $m_{h}$. Choose a constant $l_h$ such that $2^{l_{h}}> m_{h}$. Without loss of generality we can assume that $\{l_{h}\}_{h\in \mathbb{N}_{0}}$ is an increasing sequence and $l_{0} \geq 1$. \\

We take $sh_{\gamma}$ to be the linear extension of a function which scales by $\frac{1}{2^{l_{h}}}$ on the set $[2^{h+l_{h}}, 2^{h+l_{h}+1}]$. This is a continuous shrinking function, coarse-Lipschitz homotopic to the identity. For any two points $2^{l_{h}}t, 2^{l_{h}} s \in [2^{h+l_{h}}, 2^{h+l_{h}+1}]$, the composition $\gamma\circ sh_{\gamma}$ takes them to $\gamma(t),\gamma(s)$. If $|t-s|<d_{h}$ then $|\gamma sh_{\gamma}(2^{l_{h}}t)-\gamma sh_{\gamma} (2^{l_{h}}s)|<1$. Otherwise, we have
\begin{align*}
|\gamma sh_{\gamma}(2^{l_{h}}t)-\gamma sh_{\gamma} (2^{l_{h}}s)| = |\gamma(t)-\gamma(s)| < 2^{l_{h}}|t-s|
\end{align*}
For two points in the interval $[2^{{h}+l_{h-1}}, 2^{h+l_{h}}]$, the composition $\gamma\circ sh_{\gamma}$ sends them to the same point, $\gamma(2^{h})$. Now suppose we have any two $t,s \in [1,\infty)$ with $t\geq s$. We have
\begin{align*}
|\gamma sh_{\gamma}(t)-\gamma sh_{\gamma}(s)|&\leq |\gamma sh_{\gamma}(t)-\gamma sh_{\gamma}(\lfloor{t}\rfloor)|+|\gamma sh_{\gamma}(\lfloor{t}\rfloor)-\gamma sh_{\gamma}(\lceil s\rceil)|+|\gamma sh_{\gamma}(\lceil s\rceil)-\gamma sh_{\gamma}(s)|\\
&< 1+\max\{|t-s|,1\}+1 \\
&\leq|t-s|+3
\end{align*}
This shows that $\gamma\circ sh_{\gamma}$ is controlled. It is clear that $\gamma\circ sh_{\gamma}$ is continuous, montonically increasing, and proper. Additionally we have that
\begin{align*}
\gamma sh_{\gamma}(t) = \gamma sh_{\gamma}(t) - \gamma sh_{\gamma}(1)+ \beta(1)^{-1} +1
 < t + C
  \end{align*} 
for all $t\in [1,\infty)$, where $C = 1 + \beta(1)^{-1}$. 

\end{proof}

Now we have a series of technical claims:

\begin{claim} The map 
\begin{align*}
cf\circ a_{f}: c[0,1]^n &\longrightarrow cX\\
(hx,h)&\longmapsto (\gamma sh_{\gamma}(h)f(x), \gamma sh_{\gamma}(h)) 
\end{align*}
is coarse. 
\end{claim}

We just compute. Let $(tx,t), (sy,s) \in c[0,1]^n$ with distance less than $Q>0$ and $t\geq s$. This means $|t-s|<Q$ and $\|x-y\|< \frac{AQ}{t}$ for $A = 1+diam ([0,1]^n)$. 
\begin{align*}
d_{cX}(cfa_{f}(tx,t), cfa_{f}(sy,s)) &\leq \|\gamma sh_{\gamma}(t)f(x)-\gamma sh_{\gamma}(t)f(y)\| + |\gamma sh_{\gamma}(t)-\gamma sh_{\gamma}(s)| (\|f(y)\| +1)\\
&< \gamma(t) \beta(\frac{AQ}{t})+ (Q+3)(diam (X) +1)\\
&= ([\beta(\frac{1}{\log(t)+1})]^{-1} +1) \beta (\frac{AQ}{t})+ (Q+3)(diam (X) +1)
\end{align*}
Since $\frac{(\log(t)+1)AQ}{t} \rightarrow 0$ as $t\rightarrow \infty$, there is a $t_0$ such that for all $t\geq t_{0}$, $\frac{1}{\log(t)+1} > \frac{AQ}{t}$. Additionally choose $t_0$ large enough so that $\beta(\frac{AQ}{t})<1$ for all $t\geq t_{0}$. Since $\beta$ is monotonically increasing, we obtain 
\begin{align*}
([\beta(\frac{1}{\log(t)+1})]^{-1} +1) \beta (\frac{AQ}{t}) < 1+1 = 2
\end{align*}
Hence for two points $(tx,t), (sy,s)$ with $t\geq s, t_0$ we have a coarse estimate. If both $t,s \leq  t_0$, then they lie in the set $B^{c[0,1]^n}_{t_{0}}:=\{(tx,t) \in c[0,1]^n\,|\, t\leq t_0\}$, which is bounded. We have
\begin{align*}
(cf a_{f} \times cf a_{f})(B^{c[0,1]^n}_{t_0} \times B^{c[0,1]^n}_{t_0}) \subset B^{cX}_{\gamma sh_{\gamma}(t_0)} \times B^{cX}_{\gamma sh_{\gamma}(t_0)}
\end{align*}
We simply add $diam (B^{cX}_{\gamma sh_{\gamma}(t_0)})$ to the distance estimate for these points. So $cf \circ a_{f}$ is controlled. It is obviously proper since $\gamma sh_{\gamma}$ is proper. Therefore, $cf \circ a_{f}$ is coarse. \end{proof}

\begin{claim} The coarse homotopy class of $cf \circ a_f$ is independent of the choices of $\beta$ and $sh_{\gamma}$. 
\end{claim}

\begin{proof} Suppose we now have a $\beta', \gamma'$. Suppose further that  $sh_{\gamma'}: [1,\infty)\rightarrow [1,\infty)$ is a shrinking map so that $\gamma' \circ sh_{\gamma'} $ is still continuous, monotonically increasing and is coarse, ie. for every $Q>0$ there exists a constant $C_{\gamma',Q}$ such that 
\begin{align*}
|t-s|<Q \implies |\gamma' sh_{\gamma'}(t)- \gamma' sh_{\gamma'}(s)|<C_{\gamma', Q}
\end{align*}
Additionally, because $[1,\infty)$ is a path metric space, $\gamma' sh_{\gamma'}$ is large-scale Lipschitz. So there exists a constant $D_{\gamma'}$ such that $|\gamma' sh_{\gamma'}(t)| < D_{\gamma'} t + D_{\gamma'}$ for all $t\in [1,\infty)$. \\

Let $\eta(h,t) = t \gamma sh_{\gamma}(h) + (1-t) \gamma'sh_{\gamma'}(h)$. Consider the homotopy
\begin{align*}
G: c([0,1]^n \times [0,1]) &\longrightarrow cX
\\ (hx,ht,h) &\longmapsto (\eta(h,t)f(x), \eta(h,t))
\end{align*}
This is a homotopy between $cf\circ a_{f}$ (when $t=1$) and $cf\circ a'_{f}$ (when $t=0$). We now show that $G$ is a coarse map. Suppose we have two points $(hx,ht,h), (ly,ls,l)$ with distance $<Q$. Assume that $h\geq l$. As before, this means that $|h-l|<Q$ and $|t-s|, \|x-y\|< \frac{AQ}{h}$ for some constant $A\geq 1$. 

\begin{align*}
d_{cX}(G(hx,ht,h), G(ly,ls,l)) &\leq \|\eta(h,t)f(x)-\eta(h,t)f(y)\|+ |\eta(h,t)-\eta(l,s)| (\|f(y)\| +1)\\
&\leq \eta(h,t)\|f(x)-f(y)\| + |\eta(h,t)-\eta(l,s)|(diam(X)+1)\\
\eta(h,t)\|f(x)-f(y)\| &=  t\gamma sh_{\gamma}(h)\|f(x)-f(y)\|  + (1-t) \gamma'sh_{\gamma'}(h) \|f(x)-f(y)\|\\
 &< t\gamma sh_{\gamma}(h)\beta(\frac{AQ}{h})  + (1-t) \gamma'sh_{\gamma'}(h) \beta'(\frac{AQ}{h})\\
 & <2t + 2(1-t) = 2
\end{align*}
for all $h$ greater than some $h_0$. 
\begin{align*}
|\eta(h,t)-\eta(l,s)| &\leq |\eta(h,t)-\eta(l,t)| + |\eta(l,t)- \eta(l,s)|\\
|\eta(h,t)-\eta(l,t)| &\leq t|\gamma sh_{\gamma}(h) - \gamma sh_{\gamma}(l)| + (1-t) |\gamma'sh_{\gamma'}(h) - \gamma'sh_{\gamma'}(l)|\\
&<(Q+3) + C_{\gamma',Q}\\
|\eta(l,t)- \eta(l,s)| &\leq |t-s|(\gamma sh_{\gamma}(l) + \gamma'sh_{\gamma'}(l))\\
&< \frac{AQ}{h}(l + C + D_{\gamma'} l + D_{\gamma'})\\
&< AQ(1+ C+ 2D_{\gamma'})
\end{align*}
 If both $h,l \leq  h_0$, we have
\begin{align*}
(G \times G)(B^{c[0,1]^{n+1}}_{h_0} \times B^{c[0,1]^{n+1}}_{h_0}) \subset B^{cX}_{\max \{\gamma sh_{\gamma}(h_0), \gamma' sh_{\gamma'}(h_0) \}} \times B^{cX}_{\max \{\gamma sh_{\gamma}(h_0), \gamma' sh_{\gamma'}(h_0) \}}
\end{align*}
This is sufficient to show that $G$ is controlled. $G$ is obviously proper.  
\end{proof}

\begin{claim} The coarse homotopy class of $cf\circ a_{f}$ is independent of the choice of representative $f$. 
\end{claim}

\begin{proof} Let $[f] =[g]$ in $\pi_{n}(X,x_0)$. This means there is a homotopy
\begin{align*}
H: ([0,1]^n \times [0,1], \partial [0,1]^n \times [0,1]) \rightarrow (X,x_0)
\end{align*}
which restricts to $f$ on $[0,1]^n \times \{0\}$ and to $g$ on $[0,1]^n \times \{1\}$.  Let $\beta_{H}$ and $\gamma_{H}$ be the functions associated with $H$ and consider the composition 
\begin{align*}
cH\circ a_{H}:c([0,1]^n \times [0,1]) &\rightarrow cX\\
(h(x,t), h) &\mapsto (\gamma_{H} sh_{\gamma_{H}}(h) H(x,t), \gamma_{H} sh_{\gamma_{H}} (h))
\end{align*}
By the same computation as above, the map $cH\circ a_{H}$ is coarse, and $cH \circ a_{H}$ provides a coarse homotopy between $cf\circ a_{f}$ and $cg\circ a_{g}$, since $\beta_{H}$ is also an appropriate function for $f$ and $g$ (by the previous claim, $cf \circ a_{f}$ and $cg\circ a_{g}$ are independent of the choice of $\beta$). Additionally, the homotopy restricted to $c(\partial[0,1]^n \times [0,1])$ has image in $\omega[1,\infty)$. So we are done. 
\end{proof}

\begin{definition} The \textit{coning homomorphism} $\mathfrak{c}: \pi_{n}(X,x_0)\rightarrow \pi_{n}^{c}(cX,\omega)$ is defined by taking a representative $f$ of $[f]\in \pi_{n}(X,x_0)$, and sending it to $b_{\gamma sh_{\gamma} \omega,\omega}[cf\circ a_{f}]$. 
\end{definition}

This is a group homomorphism. Let $[f],[g]\in \pi_{n}(X,x_0)$. We can choose the same $\beta = \max\{\beta_{f},\beta_{g}\}$ and shrinking map $a$ for both $f$ and $g$. It follows that 
\begin{align*}
\mathfrak{c}[f]\cdot \mathfrak{c}[g]= b_{\gamma sh_{\gamma} \omega,\omega}[cf\circ a]\cdot b_{\gamma sh_{\gamma} \omega,\omega}[cg\circ a]=b_{\gamma sh_{\gamma} \omega,\omega}([cf\circ a]\cdot[cg\circ a])\\
=b_{\gamma sh_{\gamma} \omega,\omega}[c(f\cdot g)\circ a]
=\mathfrak{c}([f]\cdot [g]) 
\end{align*}

\begin{remark} There are easier ways to prove the existence of a shrinking map $a_{f}$ such that $cf\circ a_{f}$ is coarse, by using the fact that $cf$ is uniformly continous on compact subsets and that the coarse structure of $c[0,1]^n$ is generated by one entourage (eg. see the proof of Theorem \ref{endcoarse}). However, the reason we have included the proof above is because it works for all compact metric spaces $K$ and continuous maps $f:K\rightarrow X$ even if $cK$ is not a path metric space. This is because we can replace $\beta(r) := \max_{d_{K}(x,y)\leq r} d_{X}(f(x),f(y))$ with
\begin{align*}
\beta_1(r):=\frac{1}{r}\int^{2r}_{r} \beta(s)ds 
\end{align*}
which is a continous, monotonically increasing function with $\beta_1(0)=0$.\\

The rest of the proof of Proposition \ref{unicoarse} follows through with $\beta_1$. For some intuition, going up the cone expands linearly along the modulus of continuity $\beta_1$, but the shrinking map is logarithmic, so for any distance scale, the shrunken map is eventually coarse. This uses the property that $cf$ is the cone of a single map, rather than some arbitrary continuous map between cones. We obtain the following corollary:
\end{remark}

\begin{cor} Let $f:K\rightarrow X$ be a continous map between compact metric spaces. There exists a shrinking map $a_{f}: cK\rightarrow cK$ such that $cf\circ a_{f}: cK\rightarrow cX$ is coarse.
\end{cor}

To define the homomorphism $\Psi_{n}: \pi_{n}(X,x_0)\rightarrow \check{\pi}_{n}(X,x_0)$, we require a small modification to the definition of $\check{\pi}_{n}(X,x_0)$. \\

Let $X$ be a paracompact Hausdorff space. Let $\Lambda$ be the subset of $O(X,x_{0})$ consisting of pairs $(\mathcal{U}, U_{0})$ where $\mathcal{U}$ is a normal open cover of $X$ and such that there is a partition of unity $\{\rho_{U}\}_{U\in \mathcal{U}}$ subordinate to $\mathcal{U}$ with $\rho_{U_{0}}(x_0) = 1$. It is known that every open cover of a paracompact Hausdorff space $X$ is normal. Additionally, if $(\mathcal{U}, U_{0})\in O(X,x_{0})$, we can refine $(\mathcal{U}, U_{0})$ to a cover $(\mathcal{V}, V_{0})$ such that $V_{0}$ is the only element of $\mathcal{V}$ containing $x_{0}$. Therefore, $\Lambda$ is cofinal in $O(X, x_{0})$ and we can write
\begin{align*}
\check{\pi}_{n}(X,x_{0}) \cong \varprojlim_{\mathcal{U}\in \Lambda} \{\pi_{n}(|N(\mathcal{U})|, U_{0}), \phi_{\mathcal{V}\mathcal{U}}\}
\end{align*}

To explain what this means, consider our fixed sequence of open covers $(\mathcal{U}_{h}, x_0)$. For each $h\in \mathbb{N}$, we can modify $\mathcal{U}_{h}$ as follows: if $B(z,\frac{1}{2^h})$ contains $x_{0}$, and $z\neq x_{0}$, we replace it by $B'(z,\frac{1}{2^h}):=B(z,\frac{1}{2^h})\setminus \{x_0\}$. Otherwise, $B'(z,\frac{1}{2^h}) = B(z,\frac{1}{2^h})$. The open cover
\begin{align*}
\mathcal{U}'_{h} := \bigcup_{z\in Z_{h}}B'(z,\frac{1}{2^h})
\end{align*}
has the same nerve $|\mathcal{U}'_{h}|$ as $|\mathcal{U}_{h}|$. This is because the intersection of finitely many open sets is open. Assuming that the singleton $\{x_0\}$ is not open (otherwise all the homotopy groups $\pi_{n}(|\mathcal{U}_{h}|,x_0)$ would be trivial for large $h$), removing the point $\{x_0\}$ does not affect whether intersections are empty or not. The bonding maps $\phi_{h+1,h}: |\mathcal{U}_{h+1}|\rightarrow |\mathcal{U}_{h}|$ are also unaffected. \\

The map $i: cX\rightarrow M_{X}$ can be modified to $\hat{i}: cX\rightarrow M_{X}$ by taking partitions of unity $\{\rho_{h}'\}_{h\in \mathbb{N}}$ subordinate to $\{\mathcal{U}'_{h}\}_{h\in \mathbb{N}}$ instead, and this also has no effect on the homomorphism $\pi_{n}^{c}(cX,\omega)\rightarrow \pi_{n}^{L,e}(M_{X},\omega)$ since we can always choose the geometric simplicial approximation $\widetilde{is_{r}f}$ to be the same as $\widetilde{\hat{i}s_{r}f}$ for any $[f]\in \pi_{n}^{c}(cX,\omega)$ and $r\in \mathbb{N}$. This is because the only case to consider is if, for a vertex $v$ in $c[0,1]^n$, $p' s_{r}f |st(v)|$ has non-trivial intersection with $\{x_0\}$. In that case, we can always define $u_{v}$ to be $(hx_{0},h)$ where $h$ is the closest integer to $p s_{r}f(v)$. Therefore, $is_{r}f$ and $\hat{i}s_{r}f$ satisfy the same geometric star condition. \\

The purpose of this is purely for the convenience of avoiding change of base ray homomorphisms. Observe that the image of $c\partial[0,1]^n$ under $\hat{i}s_{r}f$ is $s_{r} \omega$, since $\rho'_{h}(x_0) = x_{0}$ for all $h\in \mathbb{N}_{0}$. Therefore $\widetilde{\hat{i}s_{r}f}$ is homotopic to $\hat{i}s_{r}f$ by a homotopy whose image when restricted to $c\partial[0,1]^n$ is contained within the image of the base ray $\omega$.\\

For the rest of this subsection, we assume that $X$ is path-connected and locally path-connected. 

\begin{definition} \label{psispanier} Let $X$ be a paracompact Hausdorff space and $n\geq 1$. Let $\{\rho_{U}\}_{U\in \mathcal{U}}$ be a locally finite partition of unity subordinate to $\mathcal{U}$ such that $\rho_{U_{0}}(x_0) = 1$. Let 
\begin{align*}
\varphi_{\mathcal{U}}: (X,x_0) &\rightarrow (|N(\mathcal{U})|,U_{0}) \\
x &\mapsto \sum_{U\in \mathcal{U}} \rho_{U}(x) [U]
\end{align*}
be the map which sends a point to its barycentric co-ordinates in $|N(\mathcal{U})|$. We define $\Psi_{n}$ as the group homomorphism 
\begin{align*}
\Psi_{n}: \pi_{n}(X,x_0) &\rightarrow \check{\pi}_{n}(X,x_0) \cong \varprojlim_{\mathcal{U}\in \Lambda} \{\pi_{n}(|N(\mathcal{U})|, U_{0}), \phi_{\mathcal{V}\mathcal{U}}\}\\
[f] &\mapsto \{[\varphi_{\mathcal{U}}\circ f]\}_{\mathcal{U}\in \Lambda}
\end{align*}
 \end{definition}
 
 This is well-defined: the sum in $\varphi_{\mathcal{U}}(x)$ is finite since the covers are locally finite; since $\mathcal{U}\in \Lambda$ we have that $\varphi_{\mathcal{U}}(x_0) = U_{0}$. Additionally, if $(\mathcal{V},V_{0})\succeq (\mathcal{U}, U_{0})$ then $\phi_{\mathcal{V}\mathcal{U}} \circ \varphi_{\mathcal{V}}$ is homotopic to $\varphi_{\mathcal{U}}$ by a base point preserving homotopy, so $\Psi_{n}$ is a map into the direct limit. \\

Roughly speaking, the kernel of $\Psi_{n}$ is a measure of how much homotopical information is lost by approximating $X$ by polyhedra: $[f]\in \pi_{n}(X,x_0)\setminus \ker{\Psi_{n}}$ if and only if there exist a polyhedron $K$ and a map $q: (X,x_0)\rightarrow (K,k_0)$ such that $q_{*}[f]\neq 0$ in $\pi_{n}(K,k_0)$. 

\begin{prop} \label{samepsi} The composition
\begin{align*}
\pi_{n}(X,x_0)\xrightarrow {\mathfrak{c}}\pi_{n}^{c}(cX,\omega)\xrightarrow{\cong} \pi_{n}^{L,e}(M_{X},\omega)\xrightarrow{\cong} \pi_{n}^{e}(M_{X},\omega) \rightarrow \varprojlim_{h} \pi_{n} (|\mathcal{U}_{h}|,x_0)
\end{align*}
coincides with $\Psi_{n}$. 
\end{prop}
\begin{proof} Let $[f] \in \pi_{n}(X,x_0)$, where $f:([0,1]^n,\partial [0,1]^n)\rightarrow (X,x_0)$ is a representative of $[f]$. Ignoring change of base rays for now, the map $\pi_{n}(X,x_0)\rightarrow \pi_{n}^c(cX,\omega) \xrightarrow{\cong} \pi_{n}^{L,e}(M_{X},\omega)$ sends $[f]$ to the cone $[cf \circ a_{f}]$ of $f$, and then to $[\widetilde{\hat{i} s_{r} cf a_{f}}]$, for sufficiently large $r$. Observe that $(\hat{i} \circ s_{r} \circ cf \circ a_{f})_{|[0,1]^n \times (\gamma sh_{\gamma})^{-1}[2^{h}+2^r,2^{h+1}]}$ is just, by following the definitions,  
\begin{align*}
(\varphi'_{h-r} \times \id _{[2^{h-r} +1, 2^{h-r+1}]})\circ (\id_{X} \times \frac{1}{2^r} \cdot)\circ (f\times \id_{[2^h+2^r,2^{h+1}]}) \circ (\id_{[0,1]^n} \times \gamma sh_{\gamma})\\
= (\varphi'_{h-r}\circ f) \times (\frac{1}{2^r} \cdot \circ \gamma sh_{\gamma})
\end{align*} where $\varphi'_{h-r}: X\rightarrow |\mathcal{U}_{h-r}| $ is the canonical map induced by a partition of unity subordinate to $\mathcal{U}'_{h-r}$. By the same argument as in the proof of independence of choice of (geometric) simplicial approximation, $(\hat{i}  s_{r}  cf  a_{f})_{|[0,1]^n \times (\gamma sh_{\gamma})^{-1}\{2^{h}+2^r\}}$ and its simplicial approximation are continuously homotopic, with image above height $2^{h-r}+1$ in $M_{X}$. Additionally, this homotopy has image in $\omega$ when restricted to $\partial [0,1]^n \times (\gamma sh_{\gamma})^{-1}\{2^{h}+2^r\}$. Therefore the last map $\pi_{n}^{e}(M_{X},\omega) \rightarrow \varprojlim_{h} \pi_{n} (|\mathcal{U}_{h}|,x_0)$ takes us to the element 
\begin{align*}
\{[ \widetilde{\hat{i} s_{r} cf a_{f}}_{|[0,1]^n \times (\gamma sh_{\gamma})^{-1}\{2^h+2^r\}}]\}_{h> r} = \{[ \varphi'_{h-r}\circ f ] \}_{h> r} \in \{ \pi_{n}(|\mathcal{U}_{h-r}|,x_0)\}_{h> r}
\end{align*}
in $ \varprojlim_{h}\pi_{n}(|\mathcal{U}_{h}|,x_0)$, as claimed.  The changes of base ray homomorphisms we have omitted for the purpose of clarity do not affect the result.
\end{proof}

With this identification, we can say things about when the coning map $\mathfrak{c}:\pi_{n}(X,x_0)\rightarrow \pi_{n}^{c}(cX,\omega)$  is an isomorphism, and develop some intuition about what the coarse homotopy groups actually measure. Recall that when $\varprojlim_h^1\pi_{n+1}(|\mathcal{U}_{h}|,x_0)$ vanishes, we have $\pi_{n}^{c}(cX,\omega) \cong \check{\pi}_{n}(X,x_0)$. In this case, the kernel of $\mathfrak{c}$ coincides exactly with $\ker(\Psi_n)$. 

\begin{definition} Let $n\geq 1$ and $\alpha: ([0,1],0)\rightarrow (X,x_0)$ be a path and $U$ be an open neighbourhood of $\alpha(1)$ in $X$. Define
\begin{align*}
[\alpha] \ast \pi_{n}(U) = \{b^{\alpha^{-1}}_{\alpha(1),x_0}[f]\in \pi_{n}(X,x_0) \,|\, f([0,1]^n)\subset U, f(\partial [0,1]^n) = \alpha(1)\}
\end{align*}
where $b^{\alpha^{-1}}_{\alpha(1),x_0}: \pi_{n}(X,\alpha(1))\rightarrow \pi_{n}(X,x_0)$ is the change of base point homomorphism via the path $\alpha^{-1}$. 
\end{definition}
The set $[\alpha] \ast \pi_{n}(U)$ is a subgroup of $\pi_{n}(X,x_0)$. 

\begin{definition} (Spanier groups) Let $n\geq 1$, $\mathcal{U}$ be an open cover of $X$, and $x_0\in X$. The \textit{n-th Spanier group of $(X,x_0)$ with respect to} $\mathcal{U}$ is the subgroup $\pi_{n}^{Sp}(\mathcal{U},x_0)$ of $\pi_{n}(X,x_0)$ generated by the subgroups $[\alpha] \ast \pi_{n}(U)$ for all pairs $(\alpha,U)$ with $\alpha(1) \in U$ and $U\in \mathcal{U}$:
\begin{align*}
\pi_{n}^{Sp}(\mathcal{U},x_0) = \langle [\alpha] \ast \pi_{n}(U)\,|\, U\in \mathcal{U}, \alpha(1)\in U \rangle 
\end{align*}
The \textit{n-th Spanier group of $(X,x_0)$} is the intersection
\begin{align*}
\pi_{n}^{Sp}(X,x_0) = \bigcap_{\mathcal{U}\in O(X)} \pi_{n}^{Sp}(\mathcal{U},x_0)
\end{align*}
where $O(X)$ is the directed system of open covers of $X$. 
\end{definition}

\begin{definition} \begin{itemize}
\item Let $n\in \mathbb{N}_{0}$. A space $X$ is $LC^n$ if for every neighbourhood $U$ of any point $x\in X$, there exists a neighbourhood $V$ of $x$ in $U$ such that every map $f: S^{k}\rightarrow V, 0\leq k \leq n$ is nullhomotopic in $U$. 
\item A space $X$ is \textit{semi-locally $\pi_{n}$-trivial} if for every $x\in X$ there exists an open neighbourhood $U$ of $x$ such that every map $S^n\rightarrow U$ is null-homotopic in $X$. 
\end{itemize}
\end{definition}

Observe that $LC^n\implies LC^{n-1}$ and semilocally $\pi_n$-trivial.

\begin{theorem} \label{kernelspanier} (Thm $1.1$ in \cite{aceti2023elements}) Let $n\geq 1$ and $x_0\in X$. If $X$ is paracompact, Hausdorff, and $LC^{n-1}$ then $\pi_{n}^{Sp}(X,x_0)=\ker{\Psi_n}$. 
\end{theorem}

\begin{theorem} \label{coningiso} (Thm $1.2$ in \cite{aceti2023elements}) Let $n\geq 1$ and $x_0\in X$. If $X$ is paracompact, Hausdorff, $LC^{n-1}$ and semi-locally $\pi_n$-trivial, then $\Psi_n$ is an isomorphism.
\end{theorem}

\begin{theorem} (Corollary $3.11$ in  \cite{aceti2023elements}) Let $n\geq 1$. For any based space $(X,x_0)$ we have $\pi_{n}^{Sp}(X,x_0) \subset \ker (\Psi_{n})$. 
\end{theorem}

In summary: for compact metric spaces $X$, the kernel of $\Psi_{n}$ exactly corresponds to those elements $[f]$ of $\pi_{n}(X,x_0)$ such that 
\begin{align*}
\mathfrak{c} [f] \in  {\varprojlim_h}^1 \pi_{n+1}(|\mathcal{U}_{h}|,x_0) \leq \pi_{n}^{c}(cX,\omega)
\end{align*}

If $X$ is path-connected, locally path-connected and $LC^{n-1}$, then this coincides with the Spanier group. If $X$ is additionally semi-locally $\pi_{n}$-trivial, then the coning homomorphism $\mathfrak{c}$ is injective, and  $ \pi_{n}(X,x_0)$ can be identified with a subgroup of $\pi_{n}^{c}(cX,\omega)$. \\

If we assume that $\varprojlim_h^1 \pi_{n+1}(|\mathcal{U}_{h}|,x_0)$ vanishes, the coning map induces a group homomorphism $\overline{\mathfrak{c}}: \pi_{n}(X,x_0)/\pi_{n}^{Sp}(X,x_0) \rightarrow \pi_{n}^{c}(cX,\omega)\cong \check{\pi}_{n}(X,x_0)$. If $X$ is additionally $LC^{n-1}$, then $\overline{\mathfrak{c}}$ is injective. If $X$ is additionally semi-locally $\pi_{n}$-trivial, then $\mathfrak{c}:\pi_{n}(X,x_0) \rightarrow \pi_{n}^{c}(cX,\omega)$ is an isomorphism. \\

Note that if there exists a $N\in \mathbb{N}$ such that $\mathcal{U}_{h}$ is good for all $h>N$, then $X$ is $LC^{n}$ for all $n\in \mathbb{N}_{0}$. To see this, take any point $x\in X$ and $U$ a neighbourhood of $x$. There exists $\varepsilon>0$ such that $B(x,\varepsilon)\subset U$. Choose $h>N$ large enough so that $\frac{1}{2^h}+\frac{1}{2^{h+1}}<\varepsilon$. Then there is an $z\in Z_{h}$ such that $ B(z,\frac{1}{2^h})$ is contractible and $x\in B(z,\frac{1}{2^h})\subset U$. \\

Now for some interesting examples:

\begin{enumerate}
\item \textbf{Hawaiian earring}: Recall that $\pi_{n}^{c}(c\mathbb{H},\omega) \cong \check{\pi}_{n}(\mathbb{H},x_0)$ for all $n\geq0$. The conditions of Theorem \ref{kernelspanier} apply to the Hawaiian earring, but the conditions of Theorem \ref{coningiso} do not, since it is famously not semi-locally simply connected. Nonetheless, $\Psi_{1}$ is injective (which is not easy to prove), and so $\pi_{1}^{Sp}(\mathbb{H},x_0) = 0$. To see why $\Psi_{1}$ is not surjective, observe than an element of $\check{\pi}_{1}(\mathbb{H},x_0) = \varprojlim_{h}F_h$ can be represented by a sequence $(w_h)_{h\in \mathbb{N}}$ where
\begin{itemize}
\item $w_{h}\in F_{h}$ is a reduced word in the free group on the generators $a_1,\dots, a_h$.
\item Removing the letter $a_h^{\pm}$ from $w_h$ gives the word $w_{h-1}$.
\end{itemize}

Suppose we had a word $(w_h)_{h\in \mathbb{N}}$ such that for some $k$, the number of appearances of $a_k$ in $w_{h}$ grows without bound as $h\rightarrow \infty$. This cannot be represented by a continuous loop $\gamma: [0,1]\rightarrow \mathbb{H}$ since for any such $\gamma$, the length of each appearance of $a_k$, by uniform continuity, has length bounded below by a strictly positive constant; contradiction. This suggests that the fundamental group of the Hawaiian earring is exactly those sequences with the additional property that for each $k\geq 1$, the number of times the letter $a_k$ appears in $w_{h}$ is eventually constant as $h\rightarrow \infty$. Denote this subgroup of $\varprojlim_{h} F_{h}$ as $\#_{\mathbb{N}}\mathbb{Z}$. 

\begin{theorem} (Theorem 4.1 in \cite{morgan1986van}) $\Psi_{1}: \pi_{1}(\mathbb{H},x_{0})\rightarrow \varprojlim_{h} F_{h}$ embeds $\pi_{1}(\mathbb{H}, x_{0})$ isomorphically into $\#_{\mathbb{N}}\mathbb{Z}$. 
\end{theorem}

The intuitive reason why every element of $\varprojlim_{h} F_{h}$ has a representative as a coarse homotopy class of $\pi_{1}^{c}(c\mathbb{H},\omega)$ is because the additional height variable allows for "infinitely" long paths in $\mathbb{H}$. 
\item \textbf{Griffiths twin cone}: The Griffiths twin cone $\mathbb{G}$ is obtained by taking two copies of the cone of the Hawaiian earring and gluing them together at the "bad" point on the base. It is a $2$-dimensional, path-connected, locally path-connected, compact metric space. See Figure \ref{fig:griffithscone}. 

\begin{figure}
\centering
  \centering
  \includegraphics[width=0.4\linewidth]{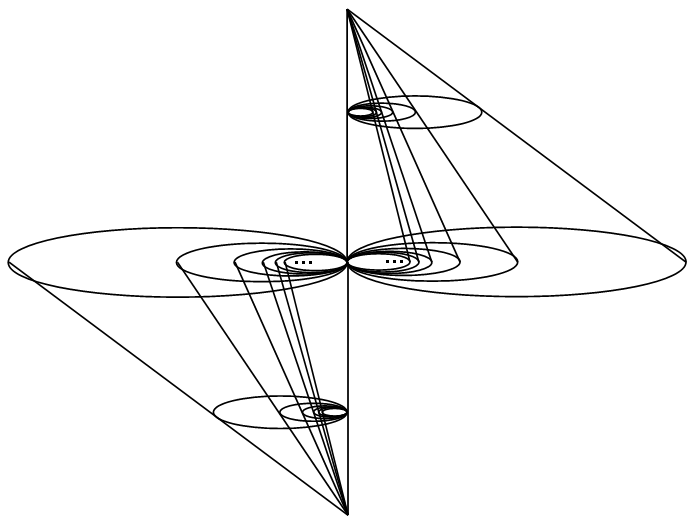}
  \caption{\cite{brazasGriffithsTwinCone}  The Griffiths twin cone $\mathbb{G}$.}
  \label{fig:griffithscone}
\end{figure}

It can be embedded into $\mathbb{R}^3$ via the alternative construction: for each integer $n\geq 1$, let $S^1_n$ be the circle of radius $\frac{1}{n}$ centred at $(\frac{1}{n}, 0)$ with basepoint $x_0 = (0,0)$. Let $\mathbb{H}_{o} =  \vee_{j\in \mathbb{N}} S^1_{2j-1}$ be the odd circles and $\mathbb{H}_{e} = \vee_{j\in \mathbb{N}} S^1_{2j}$ be the even circles so that $\mathbb{H} = \mathbb{H}_{o} \vee \mathbb{H}_{e}$ and $\mathbb{H}_{o} \cap \mathbb{H}_{e} = \{x_0\}$. Embed $\mathbb{H}\subset \mathbb{R}^3$ in the $xy$-plane. Let $C\mathbb{H}_{o} \subset \mathbb{R}^3$ be the cone over the odd circles with vertex $(0,0,1)$ and $C \mathbb{H}_{e} \subset \mathbb{R}^3$ be the cone over the even circles with vertex $(0,0,-1)$. Then $\mathbb{G} \cong C\mathbb{H}_{o} \vee_{\{x_{0}\}} C\mathbb{H}_{e}$. See Figure \ref{fig:griffiths2}. \\

\begin{figure}
\centering
  \centering
  \includegraphics[width=0.35\linewidth]{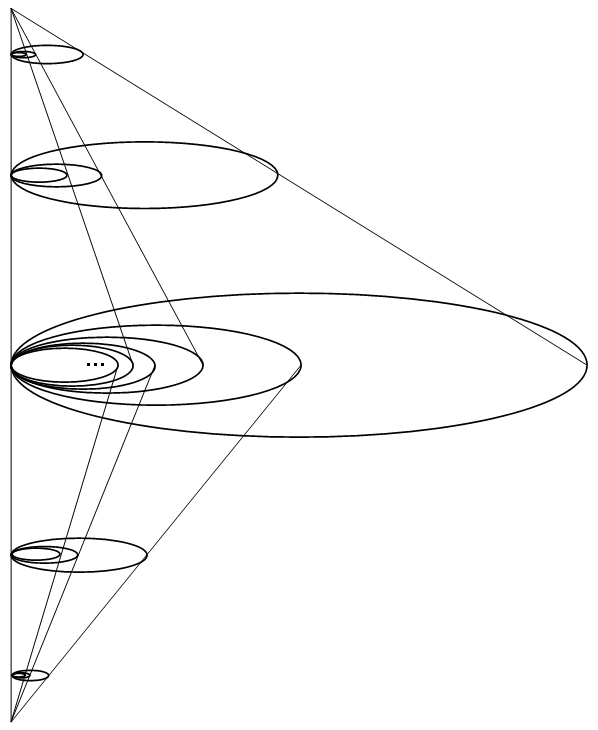}
  \caption{\cite{brazasGriffithsTwinCone}  An alternate embedding of $\mathbb{G}$ into $\mathbb{R}^3$.}
  \label{fig:griffiths2}
\end{figure}

The inclusions $\iota_o: \mathbb{H}_{o}\rightarrow \mathbb{H}$ and $\iota_{e}: \mathbb{H}_{e}\rightarrow \mathbb{H}$ induce isomorphisms on $\pi_{1}$, since $\mathbb{H}_{o}$ and $\mathbb{H}_{e}$ are homeomorphic copies of $\mathbb{H}$. Let $H_{o}$ and $H_{e}$ be these images in $\pi_{1}(\mathbb{H}, x_0)$ respectively. By the van Kampen Theorem
\begin{align*}
\pi_{1} (\mathbb{G},x_0) \cong \pi_{1}(\mathbb{H}, x_0)/N
\end{align*}
where $N$ is the normal closure of $H_{o}\cup H_{e}$. The higher homotopy groups of $\mathbb{G}$ vanish. \\

To compute the shape homotopy groups of $\mathbb{G}$, observe that $(\mathbb{G},x_0)$ is the inverse limit of the sequence 
\begin{align*}
\{C(\bigvee_{j=1}^k S^1_{2j-1})\vee_{\{x_0\}} C(\bigvee_{j=1}^k S^1_{2j}),\phi_{k+1,k};\mathbb{N}\}
\end{align*}
where the maps $\phi_{k+1,k}:C(\vee_{j=1}^{k+1} S^1_{2j-1})\vee_{\{x_0\}} C(\vee_{j=1}^{k+1} S^1_{2j})\rightarrow C(\vee_{j=1}^k S^1_{2j-1})\vee_{\{x_0\}} C(\vee_{j=1}^k S^1_{2j}) $ are induced by the collapse of the last wedge of cones. Since all these spaces are contractible, we obtain
\begin{align*}
\pi_{n}^{c}(c\mathbb{G},\omega) \cong \check{\pi}_{n}(\mathbb{G},x_0) = 0
\end{align*}
for all $n$. \\

In this example, $\Psi_{n}:\pi_{n}(\mathbb{G},x_0) \rightarrow \check{\pi}_{n}(\mathbb{G},x_0)$ makes every homotopy class trivial, ie. \textit{no} class can be detectable by polyhedra. We consider now the coarse homotopy groups of $c\mathbb{G}$. Intuitively speaking, the reason that there are any non-trivial homotopy classes in $\pi_{1}(\mathbb{G},x_0)$ at all, despite its being the wedge sum of two contractible spaces, lies in the fact that there are two cone points, and so if the path is complicated enough, the obvious null-homotopy takes too long to fit into an interval. This, however, is not a problem when looking at its image in the coarse homotopy groups under $\mathfrak{c}$, because the infinite height direction gives us the space we need.  \\

\item  \textbf{Shrinking wedge of spheres}: Recall that $ \pi_{n}^{c}(c\mathbb{E}_{m},\omega)\cong \check{\pi}_n(\mathbb{E}_{m},x_0)$ for all $n\geq 0$, and that 
\begin{align*}
\check{\pi}_n(\mathbb{E}_{m},x_0) \cong \bigoplus_{1\leq j\leq \frac{n-1}{m-1}} (\pi_{n}(S^{mj-j+1}))^{\mathbb{N}}
\end{align*} 
for $n \geq 2$. \\

$\mathbb{E}_{m}$ is $LC^{m-1}$ but is not semi-locally $\pi_{m}$-trivial. However, $\Psi_{m}: \pi_{m}(\mathbb{E}_{m},x_0)\rightarrow \check{\pi}_{m}(\mathbb{E}_{m},x_0)$ is an isomorphism where both groups are canonically isomorphic to $\mathbb{Z}^{\mathbb{N}}$. It is shown in  \cite{brazas2025v} that $\Psi_{n}$ is split surjective for $m\leq 2n-1$. 
\end{enumerate}

\newpage
\part{Relation to Shape Theory}

Shape theory is designed to extend homotopy-theoretic methods to compact metric spaces which are not homotopy equivalent to CW complexes. It replaces homotopy equivalence with shape equivalence, a weaker relation that captures topological information invariant under approximation by polyhedra. \\

Every compact metric space $X$ admits a shape expansion: that is, there exists an inverse system $\mathcal{X}:=\{X_{\alpha},\phi_{\alpha',\alpha}; \mathbb{N}\}$ of finite polyhedra associated to $X$ which is unique up to pro-isomorphism in the category pro-hCW. This allows us to define the shape of $X$ as the pro-homotopy type of a system. Shape pro-homotopy groups arise naturally from this construction, and coincide with the classical homotopy groups when $X$ is a CW complex. Shape-equivalent spaces have isomorphic pro-homotopy groups.\\

In this thesis, we also adopt the framework of strong shape theory developed via $Z$-sets and proper homotopy. In this setting, the end homotopy groups (also known as strong or Steenrod homotopy groups) provide invariants of the strong shape of a space.

\section{Proper homotopy theory}

This section introduces some basic definitions from proper homotopy theory and proves a version of the Whitehead theorem for continuous, proper maps inducing isomorphisms on end homotopy groups. 

\subsection{Pro-completion of a category}  \label{progp}

 We begin with a small digression into pro-completions. Loosely speaking, the pro-completion of a category $\mathcal{C}$, denoted by pro-$\mathcal{C}$, is a way of formally adding inverse limits to $\mathcal{C}$.  The objects of pro-$\mathcal{C}$ are cofiltered inverse systems of objects in $\mathcal{C}$ indexed over a small cofiltered category $\mathcal{A}$, and there is a canonical, fully faithful embedding 
\begin{align*}
\iota: \mathcal{C}\hookrightarrow \text{pro-}\mathcal{C}
\end{align*}
which sends an object to the constant inverse system (over any small cofiltered category $\mathcal{A}$). pro-$\mathcal{C}$ admits all small cofiltered limits, and satisfies the following universal property: for any category $\mathcal{D}$ with small cofiltered limits, the composition with $\iota$ induces an equivalence of categories
\begin{align*}
\mathrm{Fun}_{\mathrm{cofilt}}(\text{pro-}\mathcal{C},\mathcal{D}) \simeq \mathrm{Fun}(\mathcal{C},\mathcal{D})
\end{align*} 
where  $\mathrm{Fun}_{\mathrm{cofilt}}$ denotes the category of functors that preserve small cofiltered limits. \\

We will mostly set aside the categorical properties of pro-completions in favour of a hands-on approach. Let $\mathcal{C}$ be the category of groups, pointed sets, or hTop (where the objects are topological spaces and morphisms are homotopy classes of continuous maps).

\begin{definition} A \textit{directed set} $(\mathcal{A},\leq)$ is a nonempty set together with a reflexive and transitive binary relation, with the additional property that every pair of elements has an upper bound.
\end{definition}

$(\mathcal{A},\leq)$ forms a small cofiltered category by letting objects be the elements of $\mathcal{A}$ and letting there be a unique morphism $a\rightarrow b$ between two objects if $a\geq b$. Every small cofiltered category admits a cofinal functor from a directed set: one can replace a cofiltered indexing category by a directed set without changing the pro-object (p.$71$ of \cite{edwards2006cech}). This is why we will define pro-objects using directed sets instead of arbitrary cofiltered categories.

\begin{definition} (inv-$\mathcal{C}$) Let $\mathcal{A}$ be a directed set. The objects of inv-$\mathcal{C}$ are inverse systems $\{X_{\alpha}, \phi_{\alpha',\alpha}; \mathcal{A}\}$ in $\mathcal{C}$, where $\phi_{\alpha',\alpha}: X_{\alpha'}\rightarrow X_{\alpha}$ are morphisms whenever $\alpha\leq \alpha'$. A morphism of inv-$\mathcal{C}$ from $\mathcal{X}=\{X_{\alpha},\phi_{\alpha',\alpha}; \mathcal{A}\}$ to $\mathcal{Y} = \{Y_{\beta}, \psi_{\beta',\beta}; \mathcal{B}\}$ consists of a function $ \rho:\mathcal{B}\rightarrow \mathcal{A}$ and for each $\beta\in \mathcal{B}$, a morphism of $\mathcal{C}$, $q_{\beta}: X_{ \rho(\beta)}\rightarrow Y_{\beta}$, such that whenever $\beta\leq \beta'\in \mathcal{B}$ there exists $\alpha\in \mathcal{A}$ with $ \rho(\beta)\leq \alpha$ and $ \rho(\beta')\leq \alpha$ making the following diagram commute in $\mathcal{C}$.

\begin{figure}[H]
\center
\begin{tikzcd}[row sep=scriptsize, column sep=scriptsize] 
& X_{\alpha}\arrow[dl, "\phi_{\alpha, \rho(\beta)}"'] \arrow[dr, "\phi_{\alpha,  \rho(\beta')}"]   & & & \\ 
X_{ \rho(\beta)}  \arrow[dd, "q_{\beta}"'] & & X_{ \rho(\beta')}\arrow[dd, "q_{\beta'}"] \\ \\
 Y_{\beta}  & & Y_{\beta'} \arrow[ll, "\psi_{\beta',\beta}"]
 \end{tikzcd}
 \end{figure}
 
 There is an obvious definition of composition in inv-$\mathcal{C}$, and the identity morphism of $\{X_{\alpha},\phi_{\alpha',\alpha}; \mathcal{A}\}$ consists of $\id_{\mathcal{A}}$, with each $q_{\alpha} = \id_{X_{\alpha}}$. 
 
\end{definition}

\begin{definition} (pro-$\mathcal{C}$) The category pro-$\mathcal{C}$ is the quotient category of inv-$\mathcal{C}$ which has the same objects. If $\mathcal{X}$ and $\mathcal{Y}$ are objects in inv-$\mathcal{C}$, define an equivalence relation on the set of morphisms of inv-$\mathcal{C}$ from $\mathcal{X}$ to $\mathcal{Y}$ by: $( \rho, \{q_{\beta}\}) \sim ( \rho', \{q'_{\beta}\})$ if for each $\beta\in \mathcal{B}$ there exists a $\alpha\in \mathcal{A}$ with $\alpha \geq  \rho(\beta)$ and $\alpha\geq  \rho'(\beta)$ such that the following diagram commutes in $\mathcal{C}$:

\begin{figure}[H]
\center
\begin{tikzcd}[row sep=scriptsize, column sep=scriptsize] 
& X_{\alpha} \arrow[dl, "\phi_{\alpha, \rho(\beta)}"'] \arrow[dr, "\phi_{\alpha,  \rho'(\beta)}"]   & & & \\ 
X_{ \rho(\beta)}  \arrow[dr, "q_{\beta}"'] & & X_{ \rho'(\beta)}\arrow[dl, "q'_{\beta}"] \\ 
 &Y_{\beta}  & &&  
 \end{tikzcd}
 \end{figure}
 
 The equivalence classes so defined are the morphisms of pro-$\mathcal{C}$ from $\mathcal{X}$ to $\mathcal{Y}$. Taking equivalence classes is compatible with composition of morphisms so this is a well-defined category. Two objects are called \textit{pro-isomorphic} if there is an isomorphism in pro-$\mathcal{C}$ from one to the other.  

\end{definition}

The definition of pro-$\mathcal{C}$ gives the flexibility we need to change the index set $\mathcal{A}$. For example, if $\mathcal{A'}$ is cofinal in $\mathcal{A}$, $\mathcal{X}= \{X_{\alpha},\phi_{\beta,\alpha}; \mathcal{A}\}$ gives rise to a cofinal subsystem $\mathcal{X}'=\{X_{\alpha'}, \phi_{\beta',\alpha'}; \mathcal{A}'\}$ in which one retains only the objects and bonds indexed by elements of $\mathcal{A}'$. There is a restriction morphism $\mathcal{X}\rightarrow \mathcal{X}'$ in inv-$\mathcal{C}$ defined by $\mathcal{A'}\hookrightarrow \mathcal{A}$ and $\{\id_{X_{\alpha'}}\,|\,\alpha'\in \mathcal{A}'\}$. 

\begin{lemma} The restriction morphism induces an isomorphism in pro-$\mathcal{C}$.
\end{lemma}

\begin{proof} Let us construct a candidate for an inverse $\mathcal{X}'\rightarrow \mathcal{X}$. We define $\rho:\mathcal{A}\rightarrow \mathcal{A}'$ via

\[
\rho(\alpha)=
\begin{dcases*}
\alpha
   & if  $\alpha\in \mathcal{A}'$\, \\[1ex]
\tilde{\alpha}> \alpha
   & else\,
\end{dcases*}
\]

for some $\tilde{\alpha}\in \mathcal{A}'$. This exists since $\mathcal{A}'\subset \mathcal{A}$ is cofinal. As morphisms we choose the bonding morphisms of $\mathcal{X}$, ie. $\phi_{\rho(\alpha),\alpha}: X_{\rho(\alpha)}\rightarrow X_{\alpha}$. To check that this defines a morphism in inv-$\mathcal{C}$ we have to check the commutativity of the pentagon diagram. Assume that $\alpha\leq \alpha'$. Since $\mathcal{A}'$ is a directed set, there exists a $\beta\in \mathcal{A}'$ such that $\beta\geq \rho(\alpha), \rho(\alpha')$. The diagram becomes

\begin{figure}[H]
\center
\begin{tikzcd}[row sep=scriptsize, column sep=scriptsize] 
& X_{\beta}\arrow[dl] \arrow[dr]   & & & \\ 
X_{\rho(\alpha)}  \arrow[dd] & & X_{\rho(\alpha')}\arrow[dd] \\ \\
 X_{\alpha}  & & X_{\alpha'} \arrow[ll]
 \end{tikzcd}
 \end{figure}
 
 All the maps appearing are the bonding maps of $\mathcal{X}$ and thus the diagram commutes. \\
 
Now we check that this is an inverse in pro-$\mathcal{C}$. The composition $\mathcal{X}'\rightarrow\mathcal{X}\rightarrow \mathcal{X}'$ in inv-$\mathcal{C}$ is already the identity, so it is the identity in pro-$\mathcal{C}$. We first calculate the composition $\mathcal{X}\rightarrow\mathcal{X}'\rightarrow \mathcal{X}$ in inv-$\mathcal{C}$. The map of directed sets $\mathcal{A}\rightarrow \mathcal{A}$ is the composition of $\rho$ and the inclusion $\mathcal{A}'\subset \mathcal{A}$. The maps $X_{\rho(\alpha)}\rightarrow X_{\alpha}$ are the bonding maps of $\mathcal{X}$. We check that morphism is equivalent to the identity. Let $\alpha\in \mathcal{A}$, $\beta = \rho(\alpha)$. Since $\rho(\alpha)\geq \alpha$ by construction we have $\beta\geq \rho(\alpha), \id(\alpha)$. The diagram

\begin{figure}[H]
\center
\begin{tikzcd}[row sep=scriptsize, column sep=scriptsize] 
& X_{\beta} \arrow[dl, "\id_{X_{\beta}}" '] \arrow[dr, "\phi_{\beta,\alpha}"]   & & & \\ 
X_{\rho(\alpha)}  \arrow[dr, "\phi_{\beta,\alpha}" '] & & X_{\id(\alpha)}\arrow[dl, "\id_{X_{\alpha}}"] \\ 
 &X_{\alpha}  & &&  
 \end{tikzcd}
 \end{figure}

commutes. This shows that the composition is equivalent to the identity. 
 
\end{proof}

Let $\mathcal{C} = $ Gp, the category of groups. We use the symbol $0$ to denote both the trivial group and the trivial group homomorphism. 

\begin{lemma} Let $\mathcal{Z}$ be a directed set. The constant system $0:=\{0_{\zeta},0_{\zeta',\zeta}; \mathcal{Z}\}$ is a $0$-object in pro-Gp. 
\end{lemma}

\begin{proof} We check that $0$ is initial and terminal. Let $\mathcal{X}:=\{X_{\alpha}, \phi_{\alpha',\alpha}; \mathcal{A}\}$ and $\mathcal{Y}:=\{Y_{\beta}, \varphi_{\beta',\beta};\mathcal{B}\}$ be two objects. We will show that there is a unique morphism $0\rightarrow \mathcal{Y}$ (showing that $0$ is initial) and a unique morphism $\mathcal{X}\rightarrow 0$ (showing that $0$ is terminal). 
\begin{itemize}
\item $0$ is initial:\\

For $0\rightarrow \mathcal{Y}$ choose some map $\rho:\mathcal{B}\rightarrow \mathcal{Z}$. Then the maps $q_{\beta}: 0_{\rho(\beta)}\rightarrow Y_{\beta}$ have to be the $0$-morphism. Given $\beta\leq \beta'$ in $\mathcal{B}$ choose a $\zeta\in \mathcal{Z}$ such that $\zeta\geq \rho(\beta), \rho(\beta')$ which exists since $\mathcal{Z}$ is a directed set. Then the diagram

\begin{figure}[H]
\center
\begin{tikzcd}[row sep=scriptsize, column sep=scriptsize] 
& 0_{\zeta}\arrow[dl] \arrow[dr]   & & & \\ 
0_{\rho(\beta)}  \arrow[dd] & & 0_{\rho(\beta')}\arrow[dd] \\ \\
 Y_{\beta}  & & Y_{\beta'} \arrow[ll]
 \end{tikzcd}
 \end{figure}
 
 obviously commutes since both paths are $0$. This shows that there is a morphism $0\rightarrow \mathcal{Y}$. To show uniqueness, let $\rho': \mathcal{B}\rightarrow \mathcal{Z}$ be another map. Again the group homomorphisms $q'_{\beta}:0_{\rho'(\beta)}\rightarrow Y_{\beta}$ have to be the $0$-morphism. We now check that $(\rho,\{q_{\beta}\})\sim (\rho',\{q'_{\beta}\})$. Let $\beta\in \mathcal{B}$ and choose a $\zeta\in \mathcal{Z}$ with $\zeta\geq \rho(\beta), \rho'(\beta)$. Then the diagram

\begin{figure}[H]
\center
\begin{tikzcd}[row sep=scriptsize, column sep=scriptsize] 
& 0_{\zeta} \arrow[dl] \arrow[dr]   &&& \\ 
0_{\rho(\beta)}  \arrow[dr] & & 0_{\rho'(\beta)}\arrow[dl] \\ 
 &Y_{\beta}  & &&  
 \end{tikzcd}
 \end{figure}

commutes because both paths are $0$.

\item $0$ is terminal:\\

This is basically the same. To define $(\rho,\{q_{\zeta}\}): \mathcal{X}\rightarrow 0$ choose any map $\rho: \mathcal{Z}\rightarrow \mathcal{A}$. All the group homomorphisms $q_{\zeta}$ must be $0$ and the pentagon commutes since $0$ is the target (ie. in the bottom left corner), so both paths need to be the $0$-morphism. Given another $(\rho',\{q'_{\zeta}\})$, the diamond commutes since we have $0$ again as the target (ie. at the bottom).
\end{itemize}
\end{proof}

For a concrete example, we consider objects in pro-Gp indexed over the natural numbers. What does it mean for a pro-group $\mathcal{X}:=\{X_{i},\phi_{j,i};\mathbb{N}\}$ to be isomorphic to $0:=\{0_{i}, 0_{j,i};\mathbb{N}\}$? It means that there are morphisms $\mathcal{X}\xrightarrow{\rho} 0 \xrightarrow{\sigma} \mathcal{X}$ so that $\id_{\mathcal{X}} \sim \sigma\rho$. $\sigma \rho$ can be represented by the morphism $(\id_{\mathbb{N}},0_{X_{i}})$. The identity morphism $\id_{\mathcal{X}}$ can be represented by $(\id_{\mathbb{N}},\id_{X_{i}})$. The condition that $(\id_{\mathbb{N}},\id_{X_{i}}) \sim (\id_{\mathbb{N}},0_{X_{i}})$ is: for every $i\in \mathbb{N}$ there exists a $j\geq i$ such that the diagram

\begin{figure}[H]
\center
\begin{tikzcd}[row sep=scriptsize, column sep=scriptsize] 
& X_{j} \arrow[dl, "\phi_{j,i}" '] \arrow[dr, "\phi_{j,i}"]   & & & \\ 
X_{i}  \arrow[dr, "\id" '] & & X_{i}\arrow[dl, "0"] \\ 
 &X_i  & &&  
 \end{tikzcd}
 \end{figure}
commutes $\iff$ for every $i\in \mathbb{N}$ there exists a $j\geq i$ such that the bonding map $\phi_{j,i}=0$. 

\begin{example} \label{relexample1} Consider the object $\{\pi_{n}(L_i^{c}, (L_{i}\cap A)^c, \omega(i)), \psi_{i+1}; \mathbb{N}\}$ (see Theorem \ref{rellimone}) in pro-Gp, where  the bonding maps $\psi_{i+1}: \pi_{n}(L_{i+1}^{c}, (L_{i+1}\cap A)^c, \omega(i+1))\rightarrow \pi_{n}(L_i^{c}, (L_{i}\cap A)^c, \omega(i))$ are the composition of the inclusion and a change of base point $b_{\omega(i+1),\omega(i)}(\iota_{i+1,i})_{\ast}$ . We see that $\{\pi_{n}(L_i^{c}, (L_{i}\cap A)^c, \omega(i)), \psi_{i+1}; \mathbb{N}\}=0$ $\iff$ for every $i\in \mathbb{N}$ there exists a $j\geq i$ such that the morphism
\begin{align*}
\psi_{j,i}:\pi_{n}(L_{j}^c, (L_j\cap A)^c, \omega(j))\xrightarrow{(\iota_{j,i})_{\ast}} \pi_{n}(L_{i}^c, (L_i\cap A)^c, \omega(j)) \xrightarrow{b_{\omega(j),\omega(i)}}  \pi_{n}(L_{i}^c, (L_i\cap A)^c, \omega(i))
\end{align*}
where $(\iota_{j,i})_{\ast}$ denotes the morphism induced by inclusion, and $b_{\omega(j),\omega(i)}$ denotes the change of base point homomorphism using the path in $\omega$ between $\omega(j)$ and $\omega(i)$, is $0$. Since the change of base point homomorphism is an isomorphism, this means that $(\iota_{j,i})_{\ast}=0$.  
\end{example}

We now relate the condition of $\mathcal{X}$ being pro-isomorphic to $0$ to the vanishing of its $\varprojlim$ and $\varprojlim^1$.

\begin{prop}\label{protrivial} Let $\mathcal{X}=\{X_{i},\phi_{j,i};\mathbb{N}\}$ be an object in pro-Gp. Assume that $X_{i}$ is countable for all $i\in \mathbb{N}$. Then the following are equivalent: 
\begin{enumerate}
\item $\mathcal{X}=\{X_{i},\phi_{j,i};\mathbb{N}\}$ is pro-isomorphic to $0$. 
\item $\varprojlim_{i}^1 \{X_{i}\}$ and $\varprojlim_{i} \{X_{i}\}$  are trivial. 
\end{enumerate}
\end{prop}

The implication $1\implies 2$ comes from the fact that $\varprojlim$ and $\varprojlim^1$ are functors pro-Gp$\rightarrow$ Gp. $2\implies 1$ is a combination of the following two lemmas:

\begin{lemma} (Theorem 11.3.2 in \cite{geoghegan2007topological}) If $\varprojlim_{i}^1 \{X_{i}\}$ is trivial and each group $X_i$ is countable then $\{X_i\}$ is Mittag-Leffler, ie. for each $i\in \mathbb{N}$ there exists a $j\geq i$ such that for all $k\geq j$, $\image \phi_{j,i} = \image \phi_{k,i}$. 
\end{lemma} 

\begin{lemma} \label{alglimone} (Theorem 3.4 in \cite{melikhov2009steenrod}) Let $\mathcal{X}=\{X_{i},\phi_{j,i};\mathbb{N}\}$ be an inverse sequence of pointed sets. If  $\{X_i\}$ is Mittag-Leffler and $\varprojlim_{i} \{X_{i}\}$ is trivial, then for each $i$ there exists a $j>i$ such that $\phi_{j,i}: X_{j}\rightarrow X_{i}$ is trivial. 
\end{lemma}

\begin{remark} The groups in Example \ref{relexample1} are countable: by Proposition $11.4.3$ of \cite{geoghegan2007topological}, a path-connected CW complex having locally finite type is countable. By Theorem C. in \cite{wall1965finiteness} the homotopy groups of a countable CW complex are countable. Hence the absolute homotopy groups of $L_i^c$ and $(L_i\cap A)^c$ are countable. By the long exact sequence, so are the relative homotopy groups $\pi_{n}(L_{i}^c, (L_i\cap A)^c, \omega(i))$. 
\end{remark}
\subsection{Proof of $\varprojlim^1$ sequence for $\pi^e_{n}$}

\begin{theorem}\label{lim1absolutegeq1} Let $Y$ be a strongly locally finite CW complex which is path-connected. If $\mathcal{L} = \{L_i\}_{i\in \mathbb{N}}$ is a finite filtration of $Y$ and $\omega$ is well-parametrised wrt $\mathcal{L}$, then there is a natural short exact sequence of groups
\begin{align*}
0\rightarrow {\varprojlim}^1  \pi_{n+1}(Y \oset{c}{-} L_i,\omega(i))  \xrightarrow{\overline{a}} \pi^{e}_{n}(Y,\omega) \xrightarrow{\overline{b}}  \varprojlim \pi_{n}(Y \oset{c}{-} L_i,\omega(i)) \rightarrow 0
\end{align*} 
for $n\geq 1$. 
\end{theorem}

We let $L_i^c:=Y\oset{c}{-}L_i$. 

\begin{proof} Recall that $\{\pi_{n}(L^c_{i},\omega(i)), \psi_{i+1};\mathbb{N}\}$ is a system in pro-Gp, where the bonding maps $\psi_{i+1}: \pi_{n}(L^c_{i+1},\omega(i+1))\rightarrow \pi_{n}(L^c_{i},\omega(i))$ are the composition of the inclusion and a change of base point $b_{\omega(i+1),\omega(i)}(\iota_{i+1,i})_{\ast}$.\\

Let $[f]\in \pi_{n}^e(Y,\omega)$. A representative $f$ is a continuous, proper map
\begin{align*}
f: ([0,1]^n \times [1,\infty), \partial [0,1]^n \times [1,\infty)) \rightarrow (Y,\omega)
\end{align*}
For every $i\in \mathbb{N}$ there exists a $j_{i}\geq i$ such that $f_{|[0,1]^n \times [j_{i},\infty)}$  has image in $L^c_{i}$. We can assume that $\{j_{i}\}_{i\in \mathbb{N}}$ is a strictly increasing sequence. Let  $f_{j_{i}}:=f_{|[0,1]^n \times \{j_{i}\}}$ and $b_{j_i,i}:= b_{\omega(j_i),\omega(i)}$. We define $\overline{b}: \pi^{e}_{n}(Y,\omega) \rightarrow  \varprojlim \pi_{n}(Y \oset{c}{-} L_i,\omega(i))$ as follows: 
\begin{align*}
\overline{b}: \pi^{e}_{n}(Y,\omega) &\rightarrow  \varprojlim \pi_{n}(L^c_{i},\omega(i))\\
[f] &\mapsto \{b_{j_i,i}[f_{j_i}]\}_{i\in \mathbb{N}}
\end{align*}
The restriction of $f$ to $[0,1]^n \times [j_{i},j_{i+1}]$ gives us a homotopy between $f_{j_{i+1}}$ and $f_{j_i}$ which is $b_{j_{i+1},j_i}$ on the boundary. We obtain
\begin{align*}
b_{i+1,i}b_{j_{i+1},i+1}[f_{j_{i+1}}] = b_{i+1,i}b_{j_i,i+1}b_{j_{i+1},j_i} [f_{j_{i+1}}] =b_{i+1,i}b_{j_i,i+1} [f_{j_i}] = b_{j_i,i} [f_{j_i}] 
\end{align*} 
where the inclusion maps $\iota$ have been omitted. Therefore $\{b_{j_i,i}[f_{j_i}]\}_{i\in \mathbb{N}} \in  \varprojlim \pi_{n}(L^c_{i},\omega(i))$. \\
 
Suppose we have another choice of $\{j'_{i}\}_{i\in \mathbb{N}}$. For a fixed $i$, assume that $j'_{i}\geq j_i$. We have that 
\begin{align*}
b_{j'_{i},i} [f_{j'_{i}}] =b_{j_i,i} b_{j'_{i},j_{i}} [f_{j'_{i}}]= b_{j_i,i} [f_{j_i}]
\end{align*}

Now if $[f]=[g]$ there exists a continous, proper homotopy  
\begin{align*}
H:[0,1]^n \times [1,\infty) \times [0,1] \rightarrow Y
\end{align*}
such that $H_{| [0,1]^n \times [1,\infty) \times \{0\}}=f$, $H_{| [0,1]^n \times [1,\infty) \times \{1\}} = g$ and $H_{|\partial [0,1]^n \times [1,\infty) \times [0,1]} = \omega$. There exists a strictly increasing sequence $\{j_i\}_{i\in \mathbb{N}}$, such that $j_i\geq i$ and $H_{|[0,1]^n \times [j_i,\infty) \times [0,1]}$ has image in $L^c_{i}$ for all $i\in \mathbb{N}$. We choose the sequence $\{j_i\}_{i\in \mathbb{N}}$ for both $f$ and $g$. Then $b_{j_i,i}[f_{j_i}] = b_{j_i,i} [g_{j_i}]$ for all $i\in \mathbb{N}$ which gives us $\overline{b}[f] = \overline{b}[g]$.  This shows that $\overline{b}$ is well-defined. $\overline{b}$ is obviously a group homomorphism. \\

$\overline{b}$ is surjective: let $\{[f_{i}]\}_{i\in \mathbb{N}} \in \varprojlim \pi_{n}(L^c_{i},\omega(i))$. There exist homotopies $H^{i,i+1}: [0,1]^n \times [i,i+1] \rightarrow L^c_{i}$ such that 
\begin{align*}
H^{i,i+1}_{|[0,1]^n \times \{i\}} = f_i\\
H^{i,i+1}_{|[0,1]^n \times \{i+1\}} = f_{i+1}\\
H^{i,i+1}_{|\partial[0,1]^n \times [i,i+1]} = b_{i,i+1}:=b_{\omega(i),\omega(i+1)}=(b_{i+1,i})^{-1}
\end{align*}
We glue these homotopies together to a map $H: ([0,1]^n \times [1,\infty), \partial [0,1]^n \times [1,\infty))\rightarrow (Y,\omega)$. We have that $\overline{b}[H] = \{[f_i]\}_{i\in \mathbb{N}}$. \\

We now define $\overline{a}: {\varprojlim}^1  \pi_{n+1}(Y \oset{c}{-} L_i,\omega(i)) \rightarrow \pi^e_{n}(Y,\omega)$. Let $([f_i])_{i\in \mathbb{N}} \in \prod_{i\in \mathbb{N}}\pi_{n+1}(L^c_{i},\omega(i))$. We can identify $[0,1]^{n+1}$ with $[0,1]^n \times [0,\frac{1}{2}] = [0,1]^n \times [i-\frac{1}{4}, i+\frac{1}{4}]$ for all $i\in \mathbb{N}$. We define the map $G$ piecewise as follows:
\begin{align*}
G: [0,1]^n \times [i-\frac{1}{4}, i+\frac{1}{4}] &\rightarrow Y\\
 (x,t) &\mapsto f_{i}(x,t) \quad i\geq 2\\
[0,1]^n \times [i+\frac{1}{4}, i+\frac{3}{4}] &\rightarrow Y\\
 (x,t) &\mapsto \omega(2t-i-\frac{1}{2}) \quad i\geq 1\\
[0,1]^n \times  [1,1+ \frac{1}{4}] &\rightarrow Y\\
(x,t) &\mapsto f_{1}(x,t)
\end{align*}
This glues together to a continous, proper map. See Figure   \ref{fig:limproof1}. The base ray $\omega' = G_{|\partial [0,1]^n \times [1,\infty)}$ is a reparametrisation of (and therefore properly homotopic to) $\omega$. We define $a([f_i])_{i\in \mathbb{N}} = b_{\omega',\omega} [G]$. It is clear that if $([f_i])_{i\in \mathbb{N}}=([g_i])_{i\in \mathbb{N}}$ then $a([f_i])_{i\in \mathbb{N}}=a([g_i])_{i\in \mathbb{N}}$ and that $a: \prod_{i\in \mathbb{N}}\pi_{n+1}(L^c_{i},\omega(i)) \rightarrow \pi^e_{n}(Y,\omega)$ is a group homomorphism.

\begin{figure}[H]
\centering
  \centering
  \includegraphics[width=0.6\linewidth]{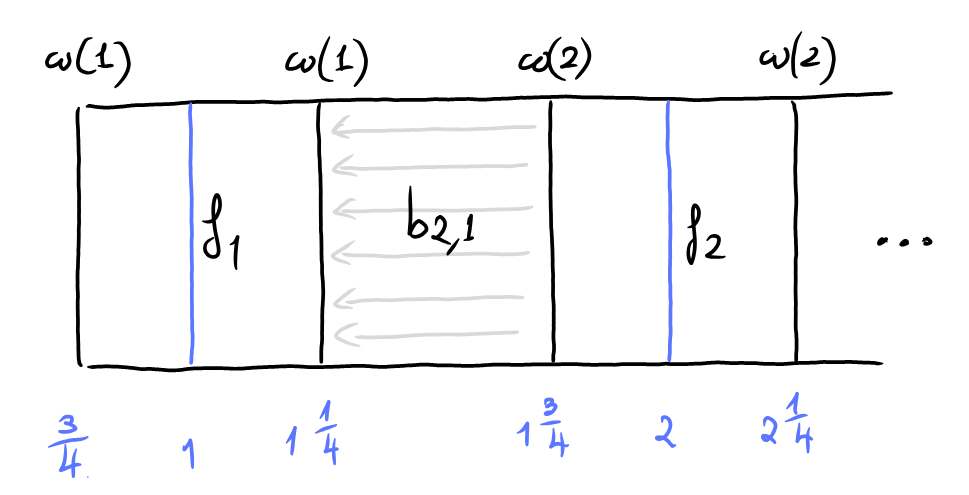}
  \caption{The map $G$}
  \label{fig:limproof1}
\end{figure}

We show that the kernel of $a = \image \varphi$, where $\varphi$ is the group homomorphism 
\begin{align*}
\varphi: \prod_{i\in \mathbb{N}}\pi_{n+1}(L^c_{i},\omega(i))&\rightarrow \prod_{i\in \mathbb{N}}\pi_{n+1}(L^c_{i},\omega(i))\\
([f_i])_{i\in \mathbb{N}}&\mapsto ([f_i]- \psi_{i+1}[f_{i+1}])_{i\in \mathbb{N}}
\end{align*}
Suppose that $([g_i])_{i\in \mathbb{N}} \in \image \varphi$, ie. there exists $([f_i])_{i\in \mathbb{N}}$ such that $[g_i] = [f_i]-\psi_{i+1}[f_{i+1}]$ for all $i\in \mathbb{N}$. We prove that $[G]$ is properly homotopic to the identity by picture. Consider the set $[0,1]^n \times [i,i+\frac{3}{4}]$ for $i\in \mathbb{N}$, labelled red in Figure \ref{fig:limproof2}.
\begin{figure}[H]
\centering
  \centering
  \includegraphics[width=0.7\linewidth]{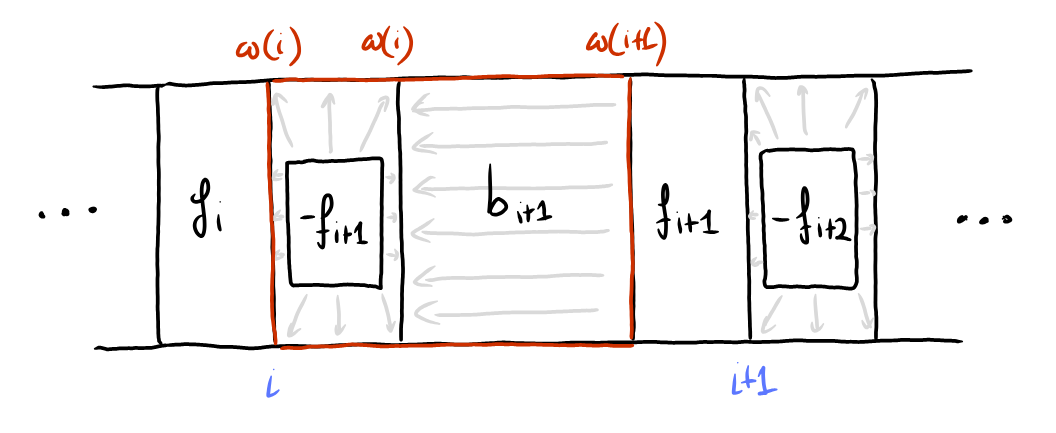}
  \caption{A representative of $[G]$, where $a([g_i])_{i\in \mathbb{N}} = b_{\omega',\omega}[G]$.}
  \label{fig:limproof2}
\end{figure}

There is a homotopy $Q: [0,1]^n \times [i,i+\frac{3}{4}] \times [0,1]$ defined by Figure $\ref{fig:homotopyQ}$. Observe that $Q_1 = G_{|[0,1]^n \times [i,i+\frac{3}{4}]}$ and that $Q$ is constant on the boundary $[0,1]^n \times \partial [i,i+\frac{3}{4}]$, and therefore can be constantly extended to all of $[0,1]^n \times [1,\infty)$, which we denote also by $Q$.

\begin{figure}[H]
\centering
  \centering
  \includegraphics[width=0.85\linewidth]{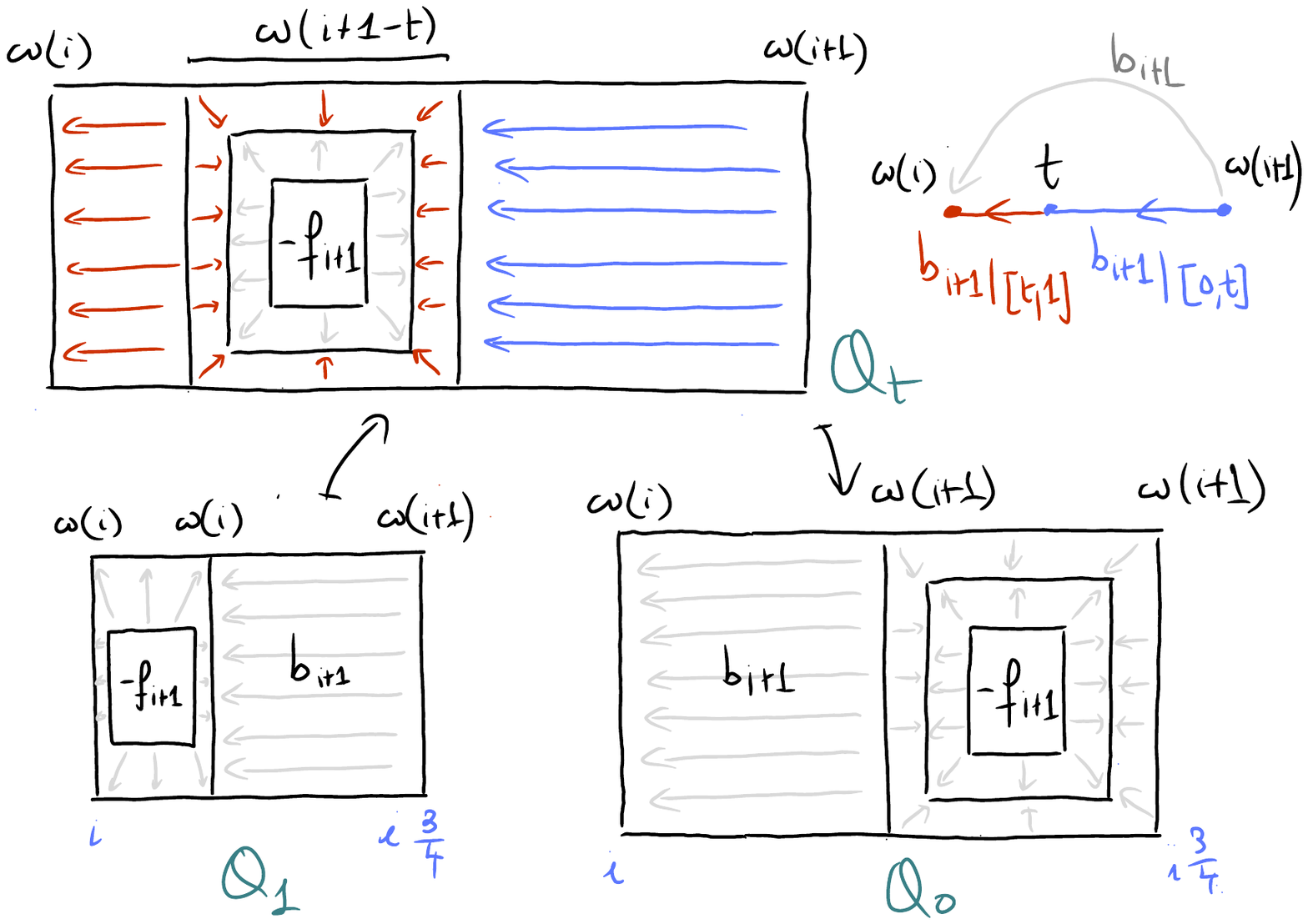}
  \caption{Diagram of the homotopy $Q$.}
  \label{fig:homotopyQ}
\end{figure}

It is clear that $Q_{0}$ restricted to the set $[0,1]^n \times [i,i+1]$ is homotopic to the change of base point $b_{i,i+1}$, via a homotopy $Q'$ that can be extended constantly to $[0,1]^n \times [1,\infty)$. We denote the concatenation of $Q^{-1}\ast Q'$ as $J_i$. Let $J: [0,1]^n \times [1,\infty) \times [0,1]$ be the infinite concatenation of homotopies $J_1\ast J_2\ast \dots$, which is well-defined since each compact set is affected by only finitely many of the $J_i$. $J$ is a proper homotopy between $G$ and the constant map to the base ray $\omega$. The homotopy is not base ray preserving, but $J_{|\partial [0,1]^n \times [1,\infty) \times [0,1]}$ always has image in $\omega$, and any two changes of base ray obtained from reparametrisation are homotopic. Therefore $b_{\omega',\omega} G$ is homotopic to the constant map to $\omega$. This shows that $\image \varphi \subset \ker a$.\\

For the other inclusion, we prove an auxilliary lemma about $\varprojlim^1$. Let $\{G_{i},\psi_{i+1};\mathbb{N}\}$ be an inverse system of (not necessarily abelian) groups. $\varprojlim^1 G_{i}$ is defined a quotient of the product $\prod_{i\in \mathbb{N}} G_{i}$ under the equivalence relation $\sim$.
\begin{align*}
(f_i)_{i\in \mathbb{N}} \sim (f'_{i})_{i\in \mathbb{N}} \iff \exists (g_i)_{i\in \mathbb{N}}: f_i = g_{i} f'_{i} \psi_{i+1}(g_{i+1})^{-1} \quad \forall i\in \mathbb{N}
\end{align*}
If the groups $G_i$ are all abelian, this agrees with the definition of $\varprojlim^1 G_{i}$  as the cokernel of the group homomorphism $\varphi$. 
\begin{lemma} \label{equlim1} Let $\{n_i\}_{i\in \mathbb{N}}$ be a strictly increasing sequence such that $n_i>i$ for all $i$. Let $(f_i)_{i\in \mathbb{N}}\in \prod_{i\in \mathbb{N}} G_{i}$. The element $(f'_i)_{i\in \mathbb{N}}$, defined as 
\begin{align*}
f'_{i} = \prod_{j\in [n_i,n_{i+1}-1]} \psi_{j,i} f_{j}
\end{align*}
is equivalent to $(f_i)_{i\in \mathbb{N}}$.
\end{lemma}

\begin{proof}
This is just a computation. Let 
\begin{align*}
g_{i} = \prod_{j\in [i,n_{i}-1]} \psi_{j,i} f_{j}
\end{align*}
Then we have 
\begin{align*}
g_{i}f'_{i} \psi_{i+1,i} (g_{i+1})^{-1} &=  \prod_{j\in [i,n_{i}-1]} \psi_{j,i} f_{j} \prod_{j\in [n_i,n_{i+1}-1]} \psi_{j,i} f_{j}( \psi_{i+1,i} (\prod_{j\in [i+1,n_{i+1}-1]} \psi_{j,i+1} f_{j})^{-1})\\
&=\prod_{j\in [i,n_{i+1}-1]} \psi_{j,i} f_{j} \prod_{j\in [n_{i+1}-1,i+1]} \psi_{j,i} (f_{j})^{-1}\\
&= f_{i}
\end{align*}
\end{proof}

We now show that $\ker a \subset \image \varphi$. Suppose we have $([f_{i}])_{i\in \mathbb{N}} \in \ker  a$, ie. the map $G$ is properly homotopic to the constant map to $\omega'$, via a base ray preserving homotopy $H: [0,1]^n \times [1,\infty) \times [0,1] \rightarrow Y$. For every $i$, choose a $n_{i}>i$ such that $H_{|[0,1]^n \times [n_{i}-\frac{1}{4}, \infty) \times [0,1]}$ has image in $L^c_{i}$. We can assume $\{n_{i}\}_{i\in \mathbb{N}}$ is strictly increasing. By Lemma \ref{equlim1}, the element 
\begin{align*}
([g_{i}])_{i\in \mathbb{N}}=(\sum_{j\in [n_{i},n_{i+1}-1]} \psi_{j,i} [f_{j}])_{i\in \mathbb{N}}
\end{align*}
is equivalent to $([f_{i}])_{i\in \mathbb{N}}$.

\begin{figure}[H]
\centering
  \centering
  \includegraphics[width=0.65\linewidth]{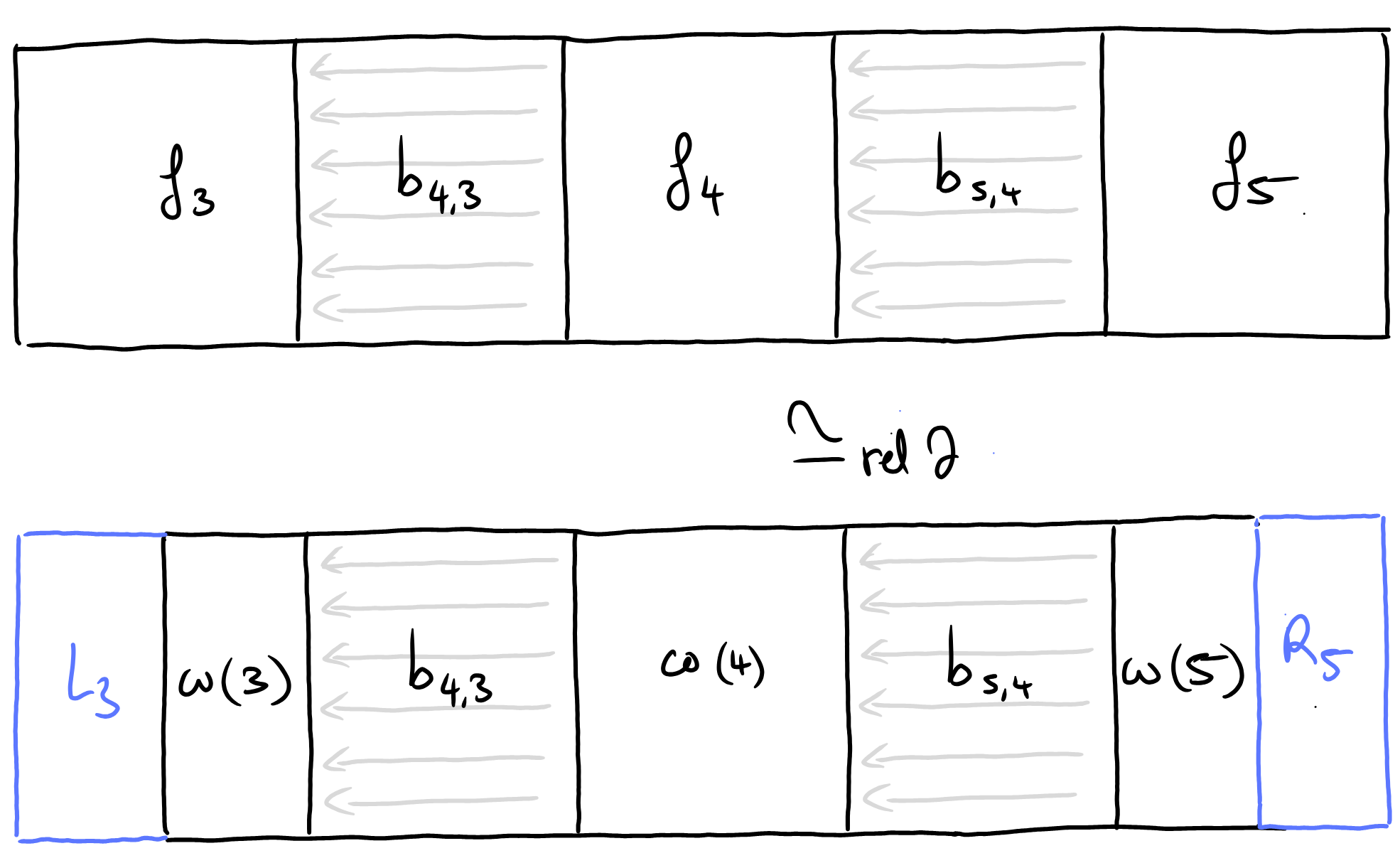}
  \caption{The two images are homotopic relative to the boundary, via the homotopy $H$ restricted to the set $[0,1]^n \times [3-\frac{1}{4},5+\frac{1}{4}] \times [0,1]$.}
  \label{fig:limproof3}
\end{figure}

\begin{figure}[H]
\centering
  \centering
  \includegraphics[width=0.75\linewidth]{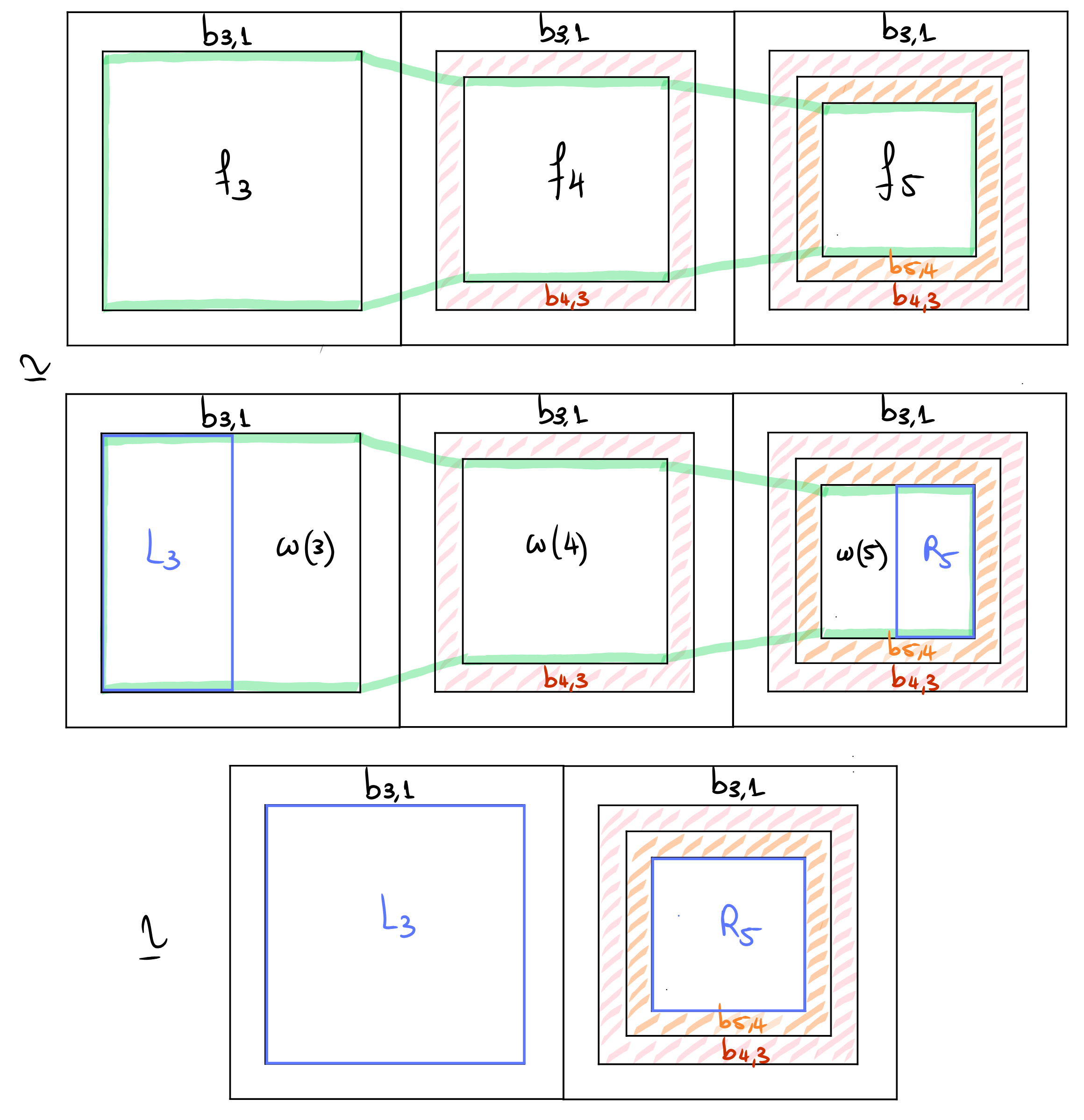}
  \caption{The homotopy between $\psi_{3,1}[f_{3}]+\psi_{4,1}[f_{4}] + \psi_{5,1}[f_5]$ and $\psi_{3,1}[L_{3}]+ \psi_{5,1} [R_{5}]$ The region outlined in green represents the homotopy in Figure \ref{fig:limproof3}.}
  \label{fig:limproof4}
\end{figure}

\begin{figure}[H]
\centering
  \centering
  \includegraphics[width=0.3\linewidth]{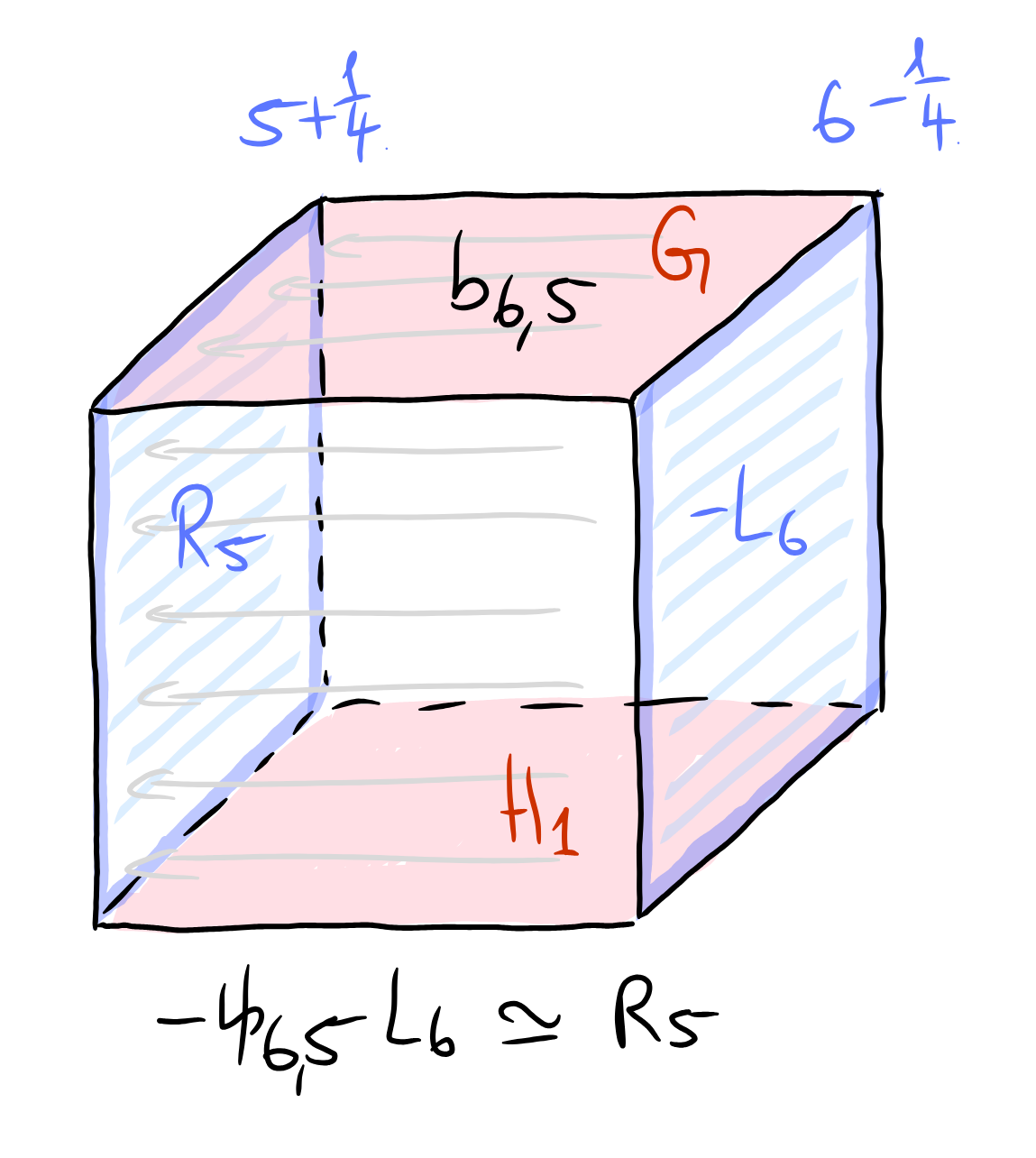}
  \caption{Diagram showing the homotopy  $[R_{n_{i+1}-1}] =- \psi_{n_{i+1},n_{i+1}-1}[L_{n_{i+1}}]$ for $i=1,n_2 = 6$.}
  \label{fig:limproof6}
\end{figure}

For each $i\in \mathbb{N}$, we have that $[g_i] = \psi_{n_i,i}[L_{n_i}] + \psi_{n_{i+1}-1,i} [R_{n_{i+1}-1}]$, where $[L_{n_{i}}] \in \pi_{n+1}(L^c_{i},\omega(n_i))$ is the homotopy class of  $H_{|[0,1]^n \times \{n_{i}-\frac{1}{4}\} \times [0,1]}$.  Similarly, $ [R_{n_{i+1}-1}]\in \pi_{n+1}(L^c_{i},\omega(n_{i+1}-1))$ is the homotopy class of $H_{|[0,1]^n \times \{n_{i+1}-1+\frac{1}{4}\} \times [0,1]}$. See Figure \ref{fig:limproof3} and Figure \ref{fig:limproof4} for a diagram of the proof for the example $[g_1]$, where $n_1=3,n_2=6$. By Figure \ref{fig:limproof6} we have that $[R_{n_{i+1}-1}] =- \psi_{n_{i+1},n_{i+1}-1}[L_{n_{i+1}}]$. Define $[h_i] = \psi_{n_{i},i}[L_{n_i}]$. We have that  
\begin{align*}
[h_i]-\psi_{i+1,i}[h_{i+1}]& = \psi_{n_{i},i}[L_{n_i}] - \psi_{i+1,i}\psi_{n_{i+1},i+1}[L_{n_{i+1}}]\\
&= \psi_{n_{i},i}[L_{n_i}] - \psi_{i+1,i}\psi_{n_{i+1}-1,i+1}\psi_{n_{i+1},n_{i+1}-1}[L_{n_{i+1}}]\\
&= \psi_{n_{i},i}[L_{n_i}] +\psi_{i+1,i}\psi_{n_{i+1}-1,i+1}[R_{n_{i+1}-1}]\\
& = \psi_{n_{i},i}[L_{n_i}] +\psi_{n_{i+1}-1,i}[R_{n_{i+1}-1}]\\
&= [g_i]
\end{align*}
This shows that $([g_i])_{i\in \mathbb{N}} \in \image \varphi$. Therefore, $\ker a = \image \varphi$, so $a$ descends to an injective group homomorphism $\overline{a}: {\varprojlim}^1  \pi_{n+1}(Y \oset{c}{-} L_i,\omega(i)) \rightarrow \pi^e_{n}(Y,\omega)$.\\

To see the sequence is exact, observe that $\image \overline{a} \subset \ker \overline{b}$ comes from Figure   \ref{fig:limproof1}: the blue lines indicate the image of $([f_i])_{i\in \mathbb{N}}$ under $\overline{b}\circ \overline{a}$. Each ${f_{i}}_{| [0,1]^n \times \{i\}}$ is nullhomotopic in $L^c_{i}$ via ${f_{i}}_{| [0,1]^n \times [i,i+\frac{1}{4}]}$. \\

Suppose now that $[G'] \in \pi^e_{n}(Y,\omega)$ lies in $\ker \overline{b}$. There exists a strictly increasing sequence $\{j_i\}_{i\in \mathbb{N}}$ such that $G'_{| [0,1]^n \times \{j_i\}}$ is nullhomotopic in $L^c_{i}$. Let $\gamma$ be the linear extension of the function $i\mapsto j_i$. By pre- and post-composing with $\id_{[0,1]^n}\times \gamma$ and a change of base ray, we obtain $[G']=[b_{\gamma,\omega} G' (\id_{[0,1]^n} \times \gamma)] =:[G]\in \pi^e_{n}(Y,\omega)$ with the property that $[G_{| [0,1]^n \times \{i\}}] \in \pi_{n}(L^c_{i},\omega(i))$ is nullhomotopic via a homotopy $H^{i}:  [0,1]^n \times [0,1] \rightarrow L^c_{i}$. Consider the restriction $G_{|[0,1]^n \times [i,i+1]}$, which we denote by $G_{i}$. We can add some change of base point "wings" to the boundary $\partial [0,1]^n \times [i,i+1]$ to turn this into an element which we denote by $g_i$ (see Figure \ref{fig:limproof7}) with the property that ${g_{i}}_{|\partial [0,1]^n \times [i,i+1]}=\omega(i)$. We can stack the sequence $-H^{i}+g_i+\psi_{i+1,i}H^{i+1} + B_{i+1,i}$ together, where $B_{i+1,i}$ is the change of base point $b_{i+1,i}$ with wings attached, and $\psi_{i+1,i}$ is written suggestively to denote the wings attached to $H^{i+1}$. We have that
\begin{align*}
[f_{i}] := [-H^{i}+g_i+\psi_{i+1,i}H^{i+1} + B_{i+1,i}] \in \pi_{n}(L^c_{i},\omega(i))
\end{align*}
for $i\in \mathbb{N}$. The element $\overline{a}([f_{i}])_{i\in \mathbb{N}}$ is represented by Figure \ref{fig:limproof8}. This element is homotopic to $G$ by iteratively performing the homotopy in Figure \ref{fig:limproof8}. (For $i=1$, we can rescale $[f_1]$ so that $-[H_1]$ takes up exactly half the interval). Note that since the boundary $\partial[0,1]^n \times [1,\infty)$ always has image in $\omega$ during the homotopies, we can freely re-parametrise the base ray via $b,b'$. Therefore $\ker \overline{b}\subset \image \overline{a}$. This completes the proof of the short exact sequence. We omit the proof for naturality, since we only need the easier statement for the inclusion of a subcomplex. 

\begin{figure}[H]
\centering
  \centering
  \includegraphics[width=0.5\linewidth]{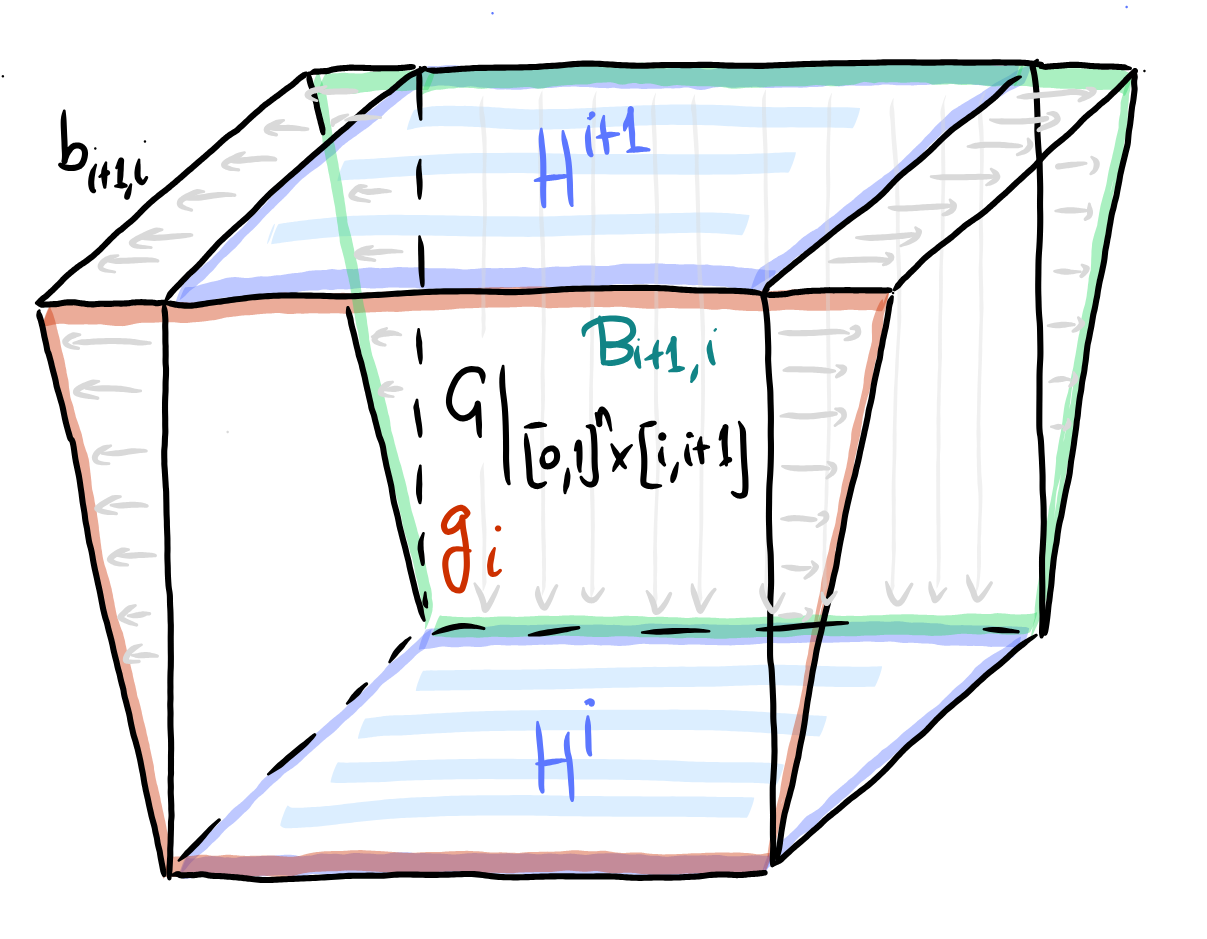}
  \caption{The segments $H^{i},H^{i+1}$ (blue), $g_i$ (red), and $B_{i+1,i}$ (green).}
  \label{fig:limproof7}
\end{figure}

\begin{figure}
\centering
  \centering
  \includegraphics[width=\linewidth]{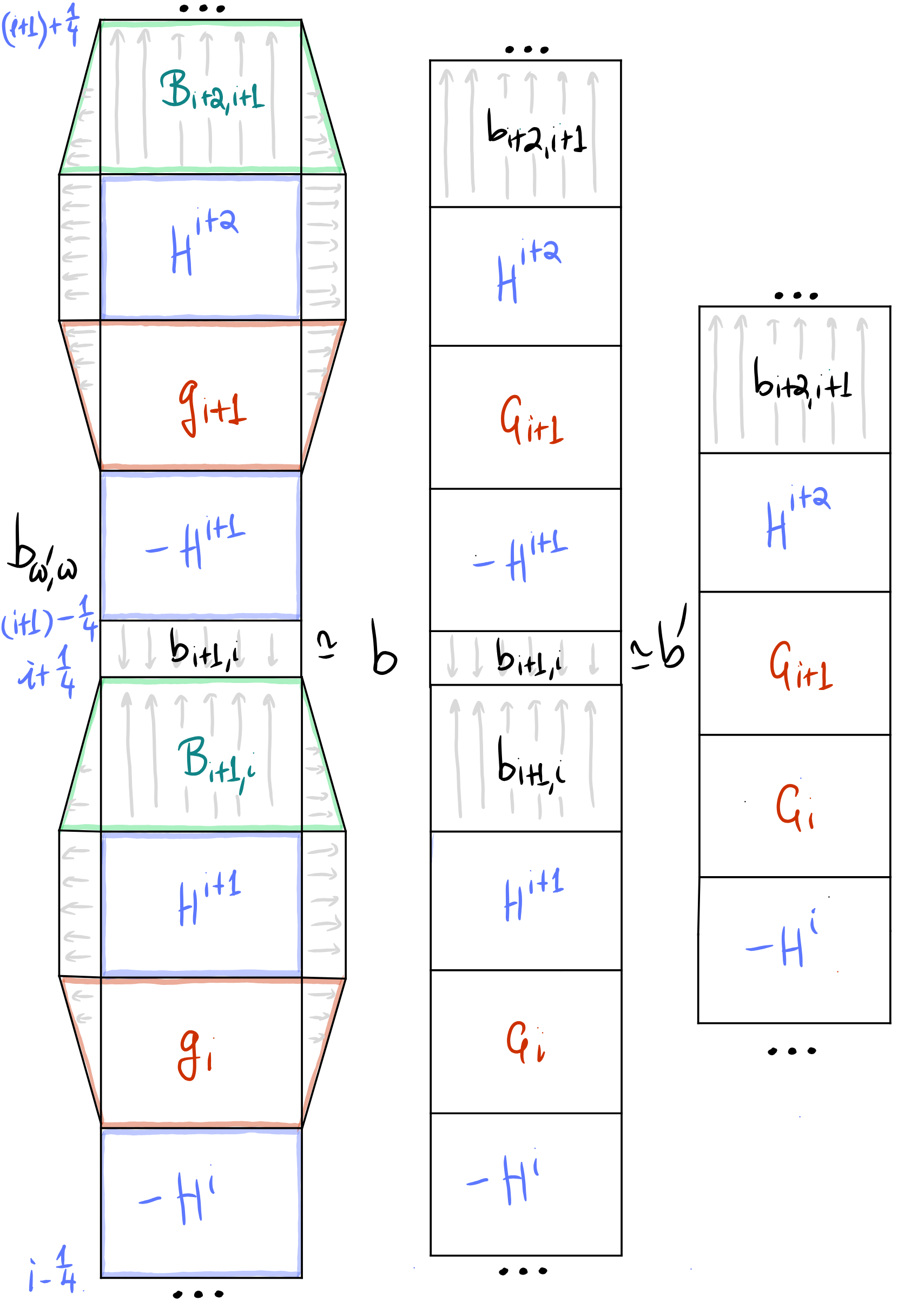}
  \caption{The homotopy between $\overline{a}([f_{i}])_{i\in \mathbb{N}}$ and $[G]$ when restricted to heights in $[i-\frac{1}{4}, (i+1) + \frac{1}{4}]$. It is  relative to $[0,1]^n \times \partial[i-\frac{1}{4}, (i+1) + \frac{1}{4}] $. $b,b'$ denote changes of base ray with image in $\omega$. }
  \label{fig:limproof8}
\end{figure}
\end{proof} 

\begin{theorem} \label{rellimone} Let $Y$ be a strongly locally finite CW complex which is path-connected, $A$ a path-connected subcomplex of $Y$. Let $\omega: [1,\infty)\rightarrow A$ be a base ray. If $\mathcal{L} = \{L_i\}$ is a finite filtration of $Y$ and $\omega$ is well-parametrised wrt $\mathcal{L}$, then there is a natural short exact sequence of groups 
\begin{align*}
0\rightarrow {\varprojlim}^1  \pi_{n+1}(L_i^c, (L_i\cap A)^c,\omega(i)) \xrightarrow{\overline{a}} \pi_{n}^e(Y,A,\omega) \xrightarrow{\overline{b}}  \varprojlim \pi_{n}(L_i^c, (L_i\cap A)^c,\omega(i))  \rightarrow 0
\end{align*} 
for $n\geq 2$, where $L_i^c:=Y\oset{c}{-}L_i$ and $(L_i\cap A)^c := A\oset{c}{-}(L_i\cap A) = A\cap(Y\oset{c}{-}L_i)$. For $n=1$, there is a natural short exact sequence of pointed sets 
\begin{align*}
 {\varprojlim}^1  \pi_{2}(L_i^c, (L_i\cap A)^c,\omega(i)) \xhookrightarrow{\overline{a}} \pi^{e}_{1}(Y,A,\omega) \xtwoheadrightarrow{\overline{b}}  \varprojlim \pi_{1}(L_i^c, (L_i\cap A)^c,\omega(i)) 
\end{align*}
where $\overline{a}$ is injective and $\overline{b}$ is surjective. 
\end{theorem}

\begin{proof} The situation for $n\geq 2$ is analagous to Theorem \ref{lim1absolutegeq1}: we can define a change of base ray homomorphism $b_{\omega',\omega}: \pi_n^{e}(Y,A,\omega')\rightarrow \pi_n^{e}(Y,A,\omega)$ by coning the construction for ordinary relative homotopy groups, and $\overline{a},\overline{b}$ are defined in the obvious way. All arguments from the proof of Theorem \ref{lim1absolutegeq1} work with slightly different diagrams.\\

For $n=1$, ${\varprojlim}^1  \pi_{n+1}(L_i^c, (L_i\cap A)^c,\omega(i)) $ is not a group, only a pointed set. This is because relative $\pi_2$ is not necessarily abelian, so the map $\varphi$ is not a group homomorphism. However, observe that the proof of Theorem \ref{lim1absolutegeq1} never actually uses the fact that higher homotopy groups are abelian: for example, Lemma \ref{equlim1} applies to arbitrary inverse systems of groups, and all the homotopies constructed apply to sequential elements. Therefore, it suffices to replace all the sums by products, and minus signs by inverses. The details are left as an exercise to the reader. 
\end{proof}

\subsection{Proper Whitehead theorem}

\begin{definition}  Let $Y$ be a finite-dimensional, strongly locally finite, path-connected CW complex, $A$ a subcomplex, $\mathcal{L}$ a finite filtration of $Y$. Let $n\geq 0$. We say $(Y,A)$ is \textit{$n$-connected at infinity} if for every $0\leq k\leq n$ and for every $i$ there exists $j\geq i$ such that every map $(D^{k},S^{k-1})\rightarrow (Y\oset{c}{-}L_j, A\oset{c}{-}(L_{j}\cap A))$ is homotopic rel $S^{k-1}$ in $Y\oset{c}{-}L_i$ to a map whose image lies in $A$.  $(Y,A)$ is called \textit{properly $n$-connected} if $(Y,A)$ is $n$-connected and $n$-connected at infinity. 
\end{definition}

For example, $(Y,A)$ is $0$-connected at infinity if the inclusion morphism on the set of ends $\iota_{*}: \varprojlim \pi_{0}((L_i\cap A)^c)\rightarrow \varprojlim \pi_{0}(L_i^c)$ is surjective.

\begin{prop} (Prop $17.1.1$ in \cite{geoghegan2007topological}) \label{propeq} The following are equivalent:
\begin{enumerate}
\item $A$ is a proper strong deformation retract of $Y$.
\item $\iota: A\rightarrow Y$ is a proper homotopy equivalence. 
\item $(Y,A)$ is properly $n$-connected for all $n$ such that $Y - A$ contains an $n$-cell. 
\end{enumerate}
\end{prop}

\begin{proof} For our purposes, we only need the implication $3\implies 1$. For this, we consider the identity map $\id_{Y}$. The pair $(Y,A)$ is properly $0$-connected, so $\id_{|Y^{0}\cup A}$ can be homotoped to have image in $A$. Call this homotopy $H^0: (Y^{0}\cup A) \times I\rightarrow Y$. \\

Assume that the homotopy $H^{k-1}$ has already defined on $Y^{k-1}\cup A$ with $(H^{k-1})_{0} = \id_{Y^{k-1}\cup A}$ and $(H^{k-1})_{1}(Y^{k-1}\cup A)\subset A$. Let $e^k$ be a $k$-cell in $Y - A$ with characteristic map $\varphi: (D^k,S^{k-1})\rightarrow (Y^k\cup A,Y^{k-1}\cup A)$. Let $\psi$ be the composite
\begin{align*}
\psi:(D^{k} \times \{0\}\cup S^{k-1} \times I, S^{k-1} \times \{1\})\xrightarrow{\varphi \cup (\varphi_{|S^{k-1}} \times \id)} \\
((Y^k \cup A)\times \{0\} \cup (Y^{k-1} \cup A) \times I, (Y^{k-1} \cup A) \times \{1\}) \xrightarrow{\id\cup H_{k-1}} (Y,A)
\end{align*} 
Observe that there is a homotopy equivalence $(D^{k} \times \{0\} \cup S^{k-1}\times I, S^{k-1} \times \{1\}) \simeq (D^{k},S^{k-1})$. Since $(Y,A)$ is $0$-connected we can use the homotopy rel $S^{k-1}$ to "fill in" the map $\psi$, ie. there is an extension of $\psi$ to
\begin{align*}
\overline{\psi}: (D^{k} \times I, S^{k-1} \times I) \rightarrow (Y,A)
\end{align*}
such that the image of $D^{k} \times \{1\}$ lies in $A$. Doing this for all $k$-cells of $Y - A$ simultaneously we obtain an extension of $H^{k-1}$ to a homotopy $H^k$
\begin{align*}
H^{k}: ((Y^{k} \cup A) \times I, (Y^{k} \cup A) \times \{1\})\rightarrow (Y, A)
\end{align*}
Repeating this process we obtain the desired homotopy 
\begin{align*}
H: Y \times I &\rightarrow Y\\
H_{|Y \times \{0\}} = \id_Y \quad &H(Y \times \{1\}) \subset A
\end{align*}
To show that we can choose $H$ to be proper, we need that $Y$ is finite-dimensional.  In defining $H$ we always choose the extension of $\psi$ to lie in $ Y\oset{c}{-}L_{j}$ for the largest possible $j$. Let $l$ be the maximal dimension of cells in $Y - A$ and let $i\in \mathbb{N}$ be fixed. Let $i_l$ to be the number such that all maps $(D^{l},S^{l-1})\rightarrow (Y\oset{c}{-}L_{i_l}, A\oset{c}{-}(L_{i_l}\cap A))$ are homotopic rel $S^{l-1}$ in $Y\oset{c}{-}L_i$ to a map whose image lies in $A$. Let $i_{l-1}$ be the number such that all maps $(D^{l-1},S^{l-2})\rightarrow (Y\oset{c}{-}L_{i_{l-1}}, A\oset{c}{-}(L_{i_{l-1}}\cap A))$ are homotopic rel $S^{l-1}$ in $Y\oset{c}{-}L_{i_l}$ to a map whose image lies in $A$. Repeating this process we obtain $L_{i}\subset L_{i_{l}}\subset L_{i_{l-1}} \subset \dots \subset L_{i_1}\subset L_{i_0}$.\\

Claim: any cell $e\in Y\oset{c}{-}L_{i_{0}}$ has the property that $H(e\times I)\subset Y\oset{c}{-}L_{i}$. To show this, observe that the $0$-skeleton's homotopy lies in $Y\oset{c}{-}L_{i_{1}}$, the $1$-skeleton's homotopy in $Y\oset{c}{-}L_{i_{2}}$, and so on, until we get that the $l$-skeleton's homotopy  has image in $Y\oset{c}{-}L_{i}$. Therefore $H^{-1}(L_i)\subset L_{i_0} \times I$ which shows that $H$ can be chosen to be proper. 
\end{proof}

We now relate proper $n$-connectedness of $(Y,A)$ to some group-theoretic properties. \\

\begin{theorem} \label{deforetract} Let $Y$ be a finite-dimensional, strongly locally finite, path-connected CW complex, $A$ a subcomplex with base ray $\tau_{1}$, $\{L_i\}_{i\in \mathbb{N}}$ a finite filtration of of $Y$. Assume that $\varprojlim \pi_{0}((L_i\cap A)^c,\tau_{1}(i))$ is a finite set with well-parametrised representatives $\{\tau_{1},\dots,\tau_{t}\}$ in $\pi_{0}^{e}(A,\tau_{1})$ and that $\iota_{*}: \varprojlim \pi_{0}((L_i\cap A)^c,\tau_{1}(i))\rightarrow \varprojlim \pi_{0}(L_i^c,\tau_{1}(i))$ is bijective. Suppose that $\iota:A\rightarrow Y$ induces an isomorphism $\iota_{*}: \pi_{n}^{e}(A,\tau) \xrightarrow{\simeq} \pi_{n}^{e}(Y,\tau)$ for all $n\geq 0$ and all $\tau \in \{\tau_{1},\dots,\tau_{t}\}$.  Then there exists a $m\in \mathbb{N}$ and a proper homotopy 
\begin{align*}
H: Y\oset{c}{-}L_m \times [0,1] \rightarrow Y
\end{align*}
from the identity to a map whose image lies in $A$. 
\end{theorem}

\begin{proof}

For $n\geq 1$ we have the following commutative diagram
\begin{figure}[H]
 \advance\leftskip-2cm
\begin{tikzcd}
&
&
& \pi_{n+1}^e(Y,A,\tau) \arrow{d}{\partial}
&
&
\\
& 0 \arrow{r}
&\varprojlim^1 \pi_{n+1}((L_i\cap A)^c,\tau(i)) \arrow{r} \arrow{d}
&\pi_{n}^e(A,\tau) \arrow{d}{\simeq} \arrow{r}
&\varprojlim \pi_{n}((L_i\cap A)^c,\tau(i)) \arrow{d} \arrow{r}
& 0
  \\
& 0 \arrow{r}
&\varprojlim^1 \pi_{n+1}(L_i^c,\tau(i)) \arrow{r}  \arrow{d}
&\pi_{n}^e(Y,\tau) \arrow{d} \arrow{r}
&\varprojlim \pi_{n}(L_i^c,\tau(i)) \arrow{d} \arrow{r}
& 0
  \\
& 0 \arrow{r}
&\varprojlim^1 \pi_{n+1}(L_i^c, (L_i\cap A)^c,\tau(i)) \arrow {r}
&\pi_{n}^e(Y,A,\tau)  \arrow{r} \arrow{d}{\partial}
&\varprojlim \pi_{n}(L_i^c, (L_i\cap A)^c,\tau(i))  \arrow{r}
& 0
\\
&
&
& \pi_{n-1}^e(A,\tau) 
&
&
\end{tikzcd}
\end{figure}

where the horizontal arrows come from the $\varprojlim^1$ sequence of end homotopy groups (Theorems \ref{lim1absolutegeq1} and \ref{rellimone}), and the middle vertical arrows come from the long exact sequence of relative end homotopy groups. All the vertical arrows are induced by inclusion maps and are the literal identity on representatives; hence the diagram commutes strictly. (Note that for $n=1$ we have to replace the bottom horizontal sequence with a short exact sequence of pointed sets.) Since $\pi_{n}^e(A,\tau)\simeq \pi_{n}^e(Y,\tau)$ for all $n\geq 0$ and $\tau\in \{\tau_{1},\dots, \tau_{t}\}$ we obtain 
\begin{align*}
\pi_{n}^e(Y,A,\tau) = 0 \\
\varprojlim \pi_{n}(L_i^c, (L_i\cap A)^c,\tau(i)) = 0\\
{\varprojlim}^{1} \pi_{n+1}(L_i^c, (L_i\cap A)^c,\tau(i)) = 0
\end{align*}
for all $n\geq 1$ and $\tau \in \{\tau_{1},\dots,\tau_{t}\}$. By Proposition \ref{protrivial} we have that  $\{\pi_{n}(L_{i}^c,(L_i\cap A)^c,\tau(i)),\psi_{i+1,i};\mathbb{N}\}$ is pro-isomorphic to $0$ for $n\geq 2$. We also have the commutative diagram (Theorem \ref{lim10}) of pointed sets

\begin{figure}[H]
\begin{tikzcd}
\centering
&\varprojlim^1 \pi_{1}((L_i\cap A)^c,\tau(i)) \arrow[r, hook] \arrow{d}
&\pi_{0}^e(A,\tau) \arrow{d}{\simeq}  \arrow[r, twoheadrightarrow]
&\varprojlim \pi_{0}((L_i\cap A)^c,\tau(i)) \arrow{d}{\simeq} 
  \\
&\varprojlim^1 \pi_{1}(L_i^c,\tau(i))\arrow[r, hook] 
&\pi_{0}^e(Y,\tau)\arrow[r, twoheadrightarrow]
&\varprojlim \pi_{0}(L_i^c,\tau(i)) 
\end{tikzcd}
\end{figure}

where the horizontal sequences are exact. Hence, we obtain an isomorphism
\begin{align*}
{\varprojlim}^{1} \pi_{1}((L_i\cap A)^c,\tau(i))\simeq {\varprojlim}^{1} \pi_{1}(L_i^c,\tau(i))
\end{align*} 
The next point is surprisingly subtle. We cannot directly apply Proposition \ref{protrivial} to  $\{\pi_{1}(L_{i}^c,(L_i\cap A)^c,\tau(i))\}$ since relative fundamental groups are only pointed sets, not groups. However, we have the following lemma:

\begin{lemma} (Lemma $3.7$ in \cite{melikhov2009steenrod}) Let $G_{i}$ be an inverse sequence of groups and $H_i$ an inverse sequence of their subgroups. 
\begin{enumerate}
\item If the inverse sequence of the pointed sets $G_{i}/H_{i}$ of right cosets satisfies the Mittag-Leffler condition, then ${\varprojlim}^{1} H_{i}\rightarrow {\varprojlim}^{1} G_{i}$ is surjective.
\item If all the $G_i$ are countable, the converse holds.
\end{enumerate} 
\end{lemma}

We can consider $H_{i}:=\pi_{1}((L_i\cap A)^c,\tau(i))$ as a subgroup of $G_{i}:=\pi_{1}(L_i^c,\tau(i))$ by identifying it with its image under the inclusion morphism. Since ${\varprojlim}^{1} H_{i}\rightarrow {\varprojlim}^{1} G_{i}$ is surjective, the inverse sequence of pointed sets $\{G_{i}/H_{i}, \psi_{i+1}; \mathbb{N}\}$ satisfies the Mittag-Leffler condition. The short exact sequence of groups $0\rightarrow H_{i}\rightarrow G_{i}\rightarrow G_{i}/H_{i}\rightarrow 0$ commutes with the bonding maps. By Theorem 3.1 $d')$ of \cite{melikhov2009steenrod}, there is a $5$-term exact sequence of pointed sets 
\begin{align*}
\{\ast\}\rightarrow \varprojlim H_{i}\rightarrow \varprojlim G_{i}\rightarrow \varprojlim (G_{i}/H_{i})\rightarrow {\varprojlim}^{1}H_{i}\rightarrow {\varprojlim}^{1}G_{i}
\end{align*}
This gives us $\varprojlim (G_{i}/H_{i}) = \{\ast\}$. By Lemma \ref{alglimone} we obtain that for every $i$ there exists a $j>i$ such that $\psi_{j,i}: G_{j}/H_{j}\rightarrow G_{i}/H_{i}$ is the trivial morphism. \\

We now show that $(Y,A)$ is $n$-connected at infinity for all $n\geq 0$. For this we need the assumption that $Y$ has finitely many ends. 

\begin{prop} (Prop. 13.4.8. and Addendum 13.4.9. in \cite{geoghegan2007topological}) \label{pathcomp} Let $Y$ have $t$ ends. There exists a $q\in \mathbb{N}$ such that for all $j\geq q$ there exists a path-connected finite subcomplex $K_j$ of $Y$ containing $L_j$ such that $Y\oset{c}{-}K_{j}$ has exactly $t$ path components. 
\end{prop}

Since $K_j$ is compact, there is a $k\geq j$ such that $K_{j}\subset L_k$. Therefore for any point $y$ in $Y\oset{c}{-}L_{k}$ there is a path from $y$ to $\tau(k)$ in $Y\oset{c}{-}L_{j}$ for some $\tau \in \{\tau_1,\dots,\tau_{t}\}$. We can concatenate this path with the segment $\tau_{[j,k]}^{-1}$ to obtain a path in $Y\oset{c}{-}L_{j}$ from $y$ to $\tau(j)$. 

\begin{remark} We require that $Y$ has finitely many ends, because if it has infinitely many, even though a $K_j$ can still be chosen so that $Y\oset{c}{-}K_{j}$ has only unbounded path components, $K_{j}$ is not necessarily compact, so the properness of the later constructed homotopy $H$ fails. 
\end{remark}

Now let $n\geq 1$ be fixed and $i\in \mathbb{N}$, $i\geq q$. For every $\tau\in \{\tau_1,\dots,\tau_{t}\}$ there exists a $j_{\tau}\geq i$ such that $\iota_{*}:\pi_{n}(L_{j_{\tau}}^c,(L_{j_{\tau}}\cap A)^c, \tau(j_{\tau}))\rightarrow \pi_{n}(L_i^c,(L_i\cap A)^c,\tau(j_{\tau}))$ is the $0$ morphism. We let $j:=\max_{\tau} j_{\tau}$. The diagram

 \begin{figure}[H]
\center
\begin{tikzcd}[row sep=scriptsize, column sep=scriptsize] 
& \pi_{n}(L_{j}^c,(L_j\cap A)^c, \tau(j))\arrow[dl, "\iota"'] \arrow[dr, "\iota"]   & & & \\ 
\pi_{n}(L_{j_{\tau}}^c,(L_{j_{\tau}}\cap A)^c, \tau(j))  \arrow[dd, "b"] & & \pi_{n}(L_{i}^c,(L_i\cap A)^c, \tau(j)) \arrow[dd, "b"] \\ \\
 \pi_{n}(L_{j_{\tau}}^c,(L_{j_{\tau}}\cap A)^c, \tau(j_{\tau})) \arrow[rr, "\iota=0"]  & & \pi_{n}(L_{i}^c,(L_i\cap A)^c, \tau(j_{\tau}))
 \end{tikzcd}
 \end{figure}
 
commutes because all the maps are inclusions $\iota$ or changes of base point $b$ using the path segment $\tau_{[j_{\tau},j]}^{-1}$. Since the $b$'s are isomorphisms, we obtain that 
 \begin{align*}
\iota_{*}:\pi_{n}(L_j^c,(L_j\cap A)^c, \tau(j))\rightarrow \pi_{n}(L_i^c,(L_i\cap A)^c,\tau(j))
\end{align*}
is the $0$ morphism for all $\tau\in \{\tau_1,\dots,\tau_{t}\}$.\\

Now let $k$, $k\geq j$ be the number which satisfies the condition obtained from Proposition \ref{pathcomp}, ie. every point in  $L_{k}^{c}$ is connected to some $\tau(j)$ by a path in $L_{j}^c$. Suppose we now have a map:
\begin{align*}
f: (D^n,S^{n-1},s)\rightarrow (L_{k}^c,(L_k\cap A)^c, f(s))
\end{align*}
This defines a homotopy class $[f]\in \pi_{n}(L_{k}^c,(L_k\cap A)^c, f(s))$. Now we have the diagram
 
 \begin{figure}[H]
\center
\begin{tikzcd}[row sep=scriptsize, column sep=scriptsize] 
& \pi_{n}(L_{k}^c,(L_k\cap A)^c, f(s))\arrow[dl, "\iota"'] \arrow[dr, "\iota"]   & & & \\ 
\pi_{n}(L_{j}^c,(L_j\cap A)^c, f(s))  \arrow[dd, "b"] & & \pi_{n}(L_{i}^c,(L_i\cap A)^c, f(s)) \arrow[dd, "b"] \\ \\
 \pi_{n}(L_{j}^c,(L_j\cap A)^c, \tau(j)) \arrow[rr, "\iota=0"]  & & \pi_{n}(L_{i}^c,(L_i\cap A)^c, \tau(j))
 \end{tikzcd}
 \end{figure}
 
 which commutes since all the maps are inclusions $\iota$ or changes of base point $b$ using the path from $f(s)$ to $\tau(j)$ constructed above. Therefore, $\iota_{*}: \pi_{n}(L_{k}^c,(L_k\cap A)^c, f(s))\rightarrow \pi_{n}(L_{i}^c,(L_i\cap A)^c, f(s)) $ is the $0$-map. This shows that every map $(D^n,S^{n-1})\rightarrow (L_{k}^c,(L_{k}\cap A)^c)$ can be homotoped rel $S^{n-1}$ in $L_{i}^c$ to a map whose image lies in $A$. \\

This shows that $(Y,A)$ is $n$-connected at infinity for all $n\geq 1$. $(Y,A)$ is $0$-connected at infinity because $\iota: A\rightarrow Y$ is surjective on the set of ends. If additionally $(Y,A)$ is $n$-connected for all $n\geq 0$, (eg. if both $A$ and $Y$ are contractible), then by Proposition \ref{propeq} we obtain a proper strong deformation retract of $Y$ to $A$. Without this assumption, there still exists a proper deformation retract outside a compact subcomplex, defined in the following way:\\

Let $l$ be the maximal dimension of cells in $Y - A$. Let $q_l$ be the number such that all maps $(D^{l},S^{l-1})\rightarrow (Y\oset{c}{-}L_{q_l}, A\oset{c}{-}(L_{q_l}\cap A))$ are homotopic rel $S^{l-1}$ to a map whose image lies in $A$. Let $q_{l-1}$ be the number such that all maps $(D^{l-1},S^{l-2})\rightarrow (Y\oset{c}{-}L_{q_{l-1}}, A\oset{c}{-}(L_{q_{l-1}}\cap A))$ are homotopic rel $S^{l-2}$ in $Y\oset{c}{-}L_{q_{l}}$ to a map whose image lies in $A$. Repeating this process we obtain $ L_{q_{l}}\subset L_{q_{l-1}}\subset \dots \subset L_{q_1}\subset L_{q_0}$. Let $m:=q_{0}$. \\

There is a proper homotopy $H: Y\oset{c}{-}L_m \times [0,1] \rightarrow Y$ from the identity map to a map whose image lies in $A$.  To see this, we simply follow the proof of Proposition \ref{propeq}, this time only considering cells in $(Y-A)\cap (Y\oset{c}{-}L_m)$. Our choice of $m$ means that $H$ can be extended skeleta by skeleta to be defined on all of $Y\oset{c}{-}L_m$. Furthermore, as in the proof of Proposition \ref{propeq}, we can choose $H$ to be proper by picking for every $i\in \mathbb{N}$, $i\geq m$ a sequence of natural numbers $i\leq i_{l}\leq i_{l-1}\leq \dots\leq i_{0}$ such that any cell $e\in Y\oset{c}{-}L_{i_0}$ has the property that $H(e\times [0,1]) \subset Y\oset{c}{-}L_i$.\\

This concludes the proof of Theorem \ref{deforetract}.
 \end{proof}

Now we extend to the case where $\psi: X\rightarrow Y$ is a proper cellular map which induces an isomorphism on end homotopy groups and the set of ends by a standard trick in topology involving the mapping cylinder $M_{\psi}$. 

\begin{cor} \label{inverse} Let $X$ and $Y$ be finite-dimensional, strongly locally finite, path-connected CW complexes with finitely many ends. Let $\omega_{1}, \dots, \omega_{t} \in \pi_{0}^{e}(X,\omega_1)$ be proper rays which represent the ends of $X$. Let $\psi:X\rightarrow Y$ be a proper cellular map which induces an isomorphism $\psi: \pi_{n}^{e}(X,\omega)\rightarrow \pi_{n}^{e}(Y,\psi \omega)$ for all $n\geq 0$ and all $\omega \in \{\omega_1,\dots,\omega_{t}\}$, and a bijection on the set of ends. Then there exist finite subcomplexes $Q\subset Y$ and $L\subset X$ and a map
\begin{align*}
\xi: (Y\oset{c}{-}Q)\rightarrow X
\end{align*}
such that 
\begin{align*}
\xi \circ \psi_{|X\oset{c}{-}L} \simeq \id_{|X\oset{c}{-}L} \\
\psi \circ \xi \simeq \id_{|Y\oset{c}{-}Q}
\end{align*}
If $\psi$ additionally induces an isomorphism on ordinary homotopy groups, $Q=\emptyset$, $L=\emptyset$, and $\xi$ is a proper homotopy inverse to $\psi$. 
\end{cor}

\begin{proof} We consider the mapping cylinder
\begin{align*}
M_{\psi}:= X\times [0,1] \cup_{\psi}Y
\end{align*}
There are canonical maps 
\begin{align*}
i: X\rightarrow M_{\psi}\\
j: Y\rightarrow M_{\psi}\\
j^{-1}: M_{\psi}\rightarrow Y
\end{align*}
where $i$ denotes inclusion into $X\times \{0\}$, $j$ is also the inclusion, and $j^{-1}$ is the deformation retract of $M_{\psi}$ to $Y$. $i$ and $j$ are obviously proper, and $j^{-1}$ is proper since the preimage of a compact set $K\subset Y$ is $\psi^{-1}(K) \times [0,1] \cup_{\psi} K$. By standard arguments, $j^{-1}$ is a proper homotopy inverse to $j$ and therefore induces an isomorphism on all end homotopy groups. Let $\omega$ be a base ray in $X$. Observe that the composition 
\begin{align*}
(X,\omega) \xrightarrow{i} (M_{\psi}, \omega \times \{0\}) \xrightarrow{j^{-1}} (Y,\psi \omega)
\end{align*}
 is exactly the map $\psi$, which by assumption induces an isomorphism on end homotopy groups. Therefore 
 \begin{align*}
 i_{*}: \pi_{n}^{e}(X,\omega) \rightarrow \pi_{n}^{e}(M_{\psi}, \omega \times \{0\})
 \end{align*}
 is an isomorphism for all $n\geq 0$.\\ 

We can choose a finite filtration on the mapping cylinder as follows. Choose finite filtrations $\{L_k\}_{k\in \mathbb{N}}, \{Q_{l}\}_{l\in \mathbb{N}} $ on $X$ and $Y$ respectively so that $\omega$ is well-parametrised with respect to $\{L_k\}_{k\in \mathbb{N}}$. Then for every $k$ there exists a $l_k$ such that $\psi(L_k)\subset Q_{l_{k}}$. Take $\{l_{k}\}_{k\in \mathbb{N}}$ to be increasing. We define
\begin{align*}
M_{k} := L_k \times I \cup Q_{l_{k}}
\end{align*}
$\mathcal{M}:=\{M_{k}\}_{k\in \mathbb{N}}$ is a finite filtration for $M_{\psi}$, and $\omega \times \{0\}$ is well-parametrised with respect to $\mathcal{M}$. \\

By Theorem \ref{deforetract} there is a $m\in \mathbb{N}$ and a homotopy $H:(M_{\psi}\oset{c}{-} M_{m}) \times [0,1] \rightarrow M_{\psi}$ such that $H_{1}(M_{\psi} \oset{c}{-}M_{m}) \subset X\times \{0\}$. Let $r = H_{1}$. Since $\psi$ is proper, there exists a $k_{m}\in \mathbb{N}$, $k_{m}>m$ such that $\psi(X \oset{c}{-} L_{k_{m}}) \subset Y \oset{c}{-} Q_{l_{m}}$. We have that the homotopy $G: M_{\psi}\times [0,1]\rightarrow M_{\psi}$ from $j\circ j^{-1}$ to $\id_{M_{\psi}}$, when restricted to $(M_{\psi}\oset{c}{-}M_{k_{m}}) \times [0,1]$, has image in $M_{\psi}\oset{c}{-} M_{m}$. \\
 
 \begin{figure}[H]
\center
 \begin{tikzcd}
Y\rar{j}  & M_{\psi} \rar{r}\arrow[bend right]{l}[swap]{j^{-1}}   & X\arrow[bend right]{l}[swap]{i} \arrow[bend left]{ll}[black]{\psi}
\end{tikzcd}
 \end{figure}
 
Let $K=M_{m}, K'= M_{k_{m}}$. Let $Q,L$ be finite subcomplexes chosen such that $j(Y\oset{c}{-}Q) \subset  (M_{\psi}\oset{c}{-}K')$, $i(X\oset{c}{-}L)\subset  (M_{\psi}\oset{c}{-}K')$, and  $\psi(X\oset{c}{-}L)\subset (Y\oset{c}{-}Q)$.  We define a homotopy inverse to $\psi$ as $\xi:=r\circ j_{|Q^c}:(Y\oset{c}{-}Q)\rightarrow X$. To check that this is actually what we want, we compute
 \begin{align*}
 \psi \circ r\circ j_{|Q^c} = j^{-1} \circ i \circ r\circ j_{|Q^{c}} \simeq j^{-1} \circ \id_{|K'^{c}} \circ j_{|Q^c} \simeq j^{-1} \circ j_{|Q^c} = \id_{|Q^c}\\
 r\circ j_{|Q^c} \circ \psi_{|L^c} = r\circ j \circ j^{-1} \circ i_{|L^c} \simeq r\circ \id_{K'^c} \circ i_{|L^c} = r \circ i_{|L^c} 
 = \id_{|L^c}
 \end{align*}
 All the homotopies above are proper, as they are a composition of proper maps and proper homotopies. This completes the proof. Note that if additionally $\psi$ induces an isomorphism on ordinary homotopy groups, we can let $K=\emptyset$, $Q = \emptyset$, $L = \emptyset$ and obtain a proper homotopy inverse $\xi: Y\rightarrow X$. 
\end{proof}

\begin{remark} It should be noted that the homotopies 
\begin{align*}
\xi \circ \psi_{|X\oset{c}{-}L} \simeq \id_{|X\oset{c}{-}L} \\
\psi \circ \xi \simeq \id_{|Y\oset{c}{-}Q}
\end{align*}
do not have to remain in $X\oset{c}{-}L$ and $Y\oset{c}{-}Q$ respectively.
\end{remark}

\begin{cor} \label{inverse2} Let $X$ and $Y$ be finite-dimensional, strongly locally finite, path-connected CW complexes with finitely many ends. Let $\omega_{1}, \dots, \omega_{t} \in \pi_{0}^{e}(X,\omega_1)$ be proper rays which represent the ends of $X$. Let $\psi:X\rightarrow Y$ be a proper cellular map which induces an isomorphism $\psi: \pi_{n}^{e}(X,\omega)\rightarrow \pi_{n}^{e}(Y,\psi \omega)$ for all $n\geq 0$ and all $\omega \in \{\omega_1,\dots,\omega_{t}\}$, and a bijection on the set of ends. Assume further that $\psi$ also induces an isomorphism on ordinary homotopy groups. Then there is a proper homotopy inverse $\xi$ to $\psi$. 
\end{cor}

\begin{proof} This follows from Proposition \ref{propeq}, Theorem \ref{deforetract} and Corollary \ref{inverse}. The combination of  $\psi: \pi_{n}^{e}(X,\omega)\rightarrow \pi_{n}^{e}(Y,\psi \omega)$ being an isomorphism  for all $n\geq 0$,  being an isomorphism on $\pi_n$ for all $n\geq 0$ and inducing a bijection on the finite set of ends, imply that $i: X\rightarrow M_{\psi}$ is properly $n$-connected for all $n\geq 0$. 
\end{proof}

\begin{remark} It may seem strange that we only require isomorphisms $\psi: \pi_{n}^{e}(X,\omega)\rightarrow \pi_{n}^{e}(Y,\psi \omega)$ for a fixed set of representatives $\omega \in \{\omega_1,\dots,\omega_{t}\}$ for the ends of $X$, rather than all proper homotopy classes of proper rays $\pi_{0}^{e}(X)$, as this set can be much larger. This has to do with the relationship between shape theory and strong shape theory (Subsection \ref{strongshapetheory}).
\end{remark}
\newpage

\section{Shape Theory} 

This chapter introduces some basic notions about shape theory and proves some further results about coarse homotopy groups. There are two main constructions of shape theory. The first, originally introduced by Borsuk, uses systems of neighbourhoods and fundamental sequences, and only applies to compact metrisable spaces. The second, attributed to Marde\u{s}i\'{c}, Segal and Morita,  uses inverse systems, and is applicable to arbitrary topological spaces. We will work with the inverse system contruction. 

\subsection{Inverse system construction}\label{inversesc}

 Let hTop be the homotopy category of topological spaces, where the objects are topological spaces and the morphisms between them are homotopy classes of continous functions. The homotopy category hCW of spaces homotopy equivalent to CW complexes, is a subcategory of hTop. We can construct the categories inv-hTop and pro-hTop as discussed in Subsection \ref{progp}. Then inv-hCW (resp. pro-hCW) is a subcategory of inv-hTop (resp. pro-hTop). 

\begin{definition} \begin{itemize}
\item A metrisable topological space $Y$ is an \textit{absolute neighbourhood retract} (ANR) if for any embedding $Y\subset Z$ as a closed subspace in a metrisable topological space $Z$, $Y$ is a neighbourhood retract of $Z$.
\item A metrisable topological space $Y$ is an \textit{absolute retract} (AR) if for any embedding $Y\subset Z$ as a closed subspace in a metrisable topological space $Z$, $Y$ is a retract of $Z$. 
\end{itemize}
\end{definition}

\begin{prop} (Lemma 1.1 in \cite{morita1975shapes}) For a topological space $X$, the following conditions are equivalent:
\begin{enumerate}
\item $X$ has the homotopy type of a CW complex. 
\item $X$ has the homotopy type of a simplicial complex with the metric topology.
\item $X$ has the homotopy type of an ANR. 
\end{enumerate}
\end{prop}

We call such an $X$ an ANR-space. Taking this into account, we can say that hCW is the homotopy category of spaces having the homotopy type of an ANR. In shape theory, only the properties of ANRs are used, so we refer to inverse systems in inv-hCW as ANR-systems. 

\begin{prop} (Proposition II.7.2 in \cite{sze1965theory}) A metrisable topological space is an absolute retract exactly when it is a contractible topological space and an absolute neighbourhood retract. 
\end{prop}

\begin{definition} \label{systemassocX} An inverse ANR-system $\mathcal{X} = \{X_{\alpha}, [\phi_{\alpha',\alpha}]; \mathcal{A}\} \in$ inv-hCW is said to be \textit{associated with a topological space} $X$ if there are homotopy classes of continuous maps $[p_{\alpha}]: X\rightarrow X_{\alpha}$, $\alpha\in \mathcal{A}$, satisfying the following conditions:
\begin{enumerate}
\item $[p_{\alpha}] = [\phi_{\alpha',\alpha}] [p_{\alpha'}]$ for all $\alpha \leq \alpha'$.
\item For any ANR-space $P$ and any homotopy class $[f]: X\rightarrow P$, there exist an $\alpha\in  \mathcal{A}$ and a homotopy class $[h_{\alpha}]: X_{\alpha}\rightarrow P$ such that $[h_{\alpha}][p_{\alpha}] = [f]$. 
\item If $[h] [p_{\alpha}] = [h'] [p_{\alpha}]$ for two homotopy classes $[h],[h']: X_{\alpha}\rightarrow P$, then there exists an index $\alpha'\geq \alpha$ such that $[h][\phi_{\alpha',\alpha}] = [h'] [\phi_{\alpha',\alpha}]$. 
\end{enumerate} 
We all such a system $\mathcal{X}$ a \textit{shape expansion} of $X$.
\end{definition}

Analagously, we say that an inverse ANR-system $\mathcal{X} = \{(X_{\alpha},x_{\alpha}), [\phi_{\alpha,\alpha'}]; \mathcal{A}\} \in$ inv-hCW$_{\ast}$ is associated to the pointed topological space $(X,x_0)$, if all the maps and homotopies in the above definition are pointed. We call $\mathcal{X}$ a \textit{pointed shape expansion} of $X$. 

\begin{prop} \label{prosystem} For any topological space $X$, there exists an associated inverse ANR-system. Moreover, any two ANR-systems associated to $X$ are isomorphic in pro-hCW by a canonical isomorphism. 
\end{prop}

\begin{proof} Let $\{\mathcal{U}_{\alpha} \,|\, \alpha\in \mathcal{A}\}$ be the set of normal, locally finite open covers of $X$. If a cover $\mathcal{U}_{\alpha'}$ is a refinement of $\mathcal{U}_{\alpha}$, then we write $\alpha'\geq\alpha$. Let $X_{\alpha}$ be the geometric realisation of the nerve of $\mathcal{U}_{\alpha}$; this is a simplicial complex with the metric topology. For each $\alpha\in \mathcal{A}$, there exists a continous map $p_{\alpha}: X\rightarrow X_{\alpha}$ (which can be defined via a partition of unity subordinate to $\mathcal{U}_{\alpha}$) satisfying the condition
\begin{align*}
p_{\alpha}^{-1}\langle st(u,X_{\alpha}) \rangle  \subset U
\end{align*}
where $u$ is the vertex of $X_{\alpha}$ corresponding to the element $U$ of the cover $\mathcal{U}_{\alpha}$. A map satisfying this condition is said to be a  \textit{canonical projection}. All canonical projections are homotopic to one another. If $\alpha'\geq \alpha$ then there exists a simplicial map $\phi_{\alpha',\alpha}: X_{\alpha'}\rightarrow X_{\alpha}$ (coming from the refinement of open covers) such that $\phi_{\alpha',\alpha}(u) = v$ implies $U\subset V$, where $u$ and $v$ are the vertices of the simplicial complex $X_{\alpha'}$ and $X_{\alpha}$ corresponding to the elements $U$ and $V$ of the covers $\mathcal{U}_{\alpha'}$ and $\mathcal{U}_{\alpha}$ respectively. The map $\phi_{\alpha',\alpha}$ is called a \textit{canonical bonding map}. Any two canonical bonding maps are homotopic to one another. Moreover
\begin{align*}
[\phi_{\alpha',\alpha}] [p_{\alpha'}] = [p_{\alpha}]
\end{align*}
for all $\alpha'\geq\alpha$. By Theorem 4.3 of \cite{morita1975vcech}, $\{X_{\alpha},[\phi_{\alpha',\alpha}];\mathcal{A}\}$ is associated to X.\\

To show uniqueness, suppose that $\mathcal{Y}:=\{Y_{\beta},\varphi_{\beta',\beta};\mathcal{B}\}$ is another inverse ANR-system associated to $X$, with projections $[q_{\beta}]: X\rightarrow Y_{\beta}$. For $\beta\in \mathcal{B}$ there exist $\rho(\beta)\in \mathcal{A}$ and $h_{\beta}: X_{\rho(\beta)}\rightarrow Y_{\beta}$ such that $[h_{\beta}][p_{\rho(\beta)}] = [q_{\beta}]$. Analogously, for $\alpha\in \mathcal{A}$ there exist $\eta(\alpha)\in \mathcal{B}$ and $k_{\alpha}: Y_{\eta(\alpha)}\rightarrow X_{\alpha}$ such that $[k_{\alpha}][q_{\eta(\alpha)}] = [p_{\alpha}]$. The pair $(\rho,\{h_{\beta}\}): \mathcal{X}\rightarrow \mathcal{Y}$ is a morphism in inv-hCW. This is because for any $\beta'\geq \beta$ there exists an $\alpha'\geq \rho(\beta),\rho(\beta')$ such that
the diagram
 \begin{figure}[H]
\center
\begin{tikzcd}[row sep=scriptsize, column sep=scriptsize] 
& X\arrow[dl]  \arrow[dddl] \arrow[dddddl] \arrow[dddr] \arrow[dddddr] \arrow[dr]   & & & \\ 
X_{\alpha'}  \arrow[dd, "\phi_{\alpha',\rho(\beta)}"'] & &X_{\alpha'} \arrow[dd, "\phi_{\alpha',\rho(\beta')}"] \\ \\
 X_{\rho(\beta)}\arrow[dd, "h_{\beta}"']  & & X_{\rho(\beta')} \arrow[dd, "h_{\beta'}"]\\ \\
Y_{\beta}  & &Y_{\beta'} \arrow[ll, "\varphi_{\beta',\beta}"] 
 \end{tikzcd}
 \end{figure}
commutes up to homotopy, where the unlabelled arrows are canonical projections. We obtain 
\begin{align*}
[\varphi_{\beta',\beta}][h_{\beta'}][\phi_{\alpha',\rho(\beta')}][p_{\alpha'}]=[h_{\beta}][\phi_{\alpha',\rho(\beta)}][p_{\alpha'}]
\end{align*} Since $\mathcal{Y}$ is associated to $X$, there exists $\alpha''\geq \alpha$ such that 
\begin{align*}
[\varphi_{\beta',\beta}][h_{\beta'}][\phi_{\alpha',\rho(\beta')}][\phi_{\alpha'',\alpha'}]=[h_{\beta}][\phi_{\alpha',\rho(\beta)}][\phi_{\alpha'',\alpha'}]
\end{align*}
This shows that $(\rho,\{h_{\beta}\})$ is a morphism in inv-hCW. Analogously, $(\eta, \{k_{\alpha}\}): \mathcal{Y}\rightarrow \mathcal{X}$ is a morphism in inv-hCW. \\

We show that the composition $(\eta,\{k_{\alpha}\})\circ (\rho,\{h_{\beta}\}) = (\rho\circ\eta, \{k_{\alpha}\circ h_{\eta(\alpha)}\})$ is equivalent to the identity. The diagram

 \begin{figure}[H]
\center
\begin{tikzcd}[row sep=scriptsize, column sep=scriptsize] 
& X\arrow[dl, "p_{\rho\eta(\alpha)}"']  \arrow[dddl, "q_{\eta(\alpha)}"] \arrow[dddr, "p_{\alpha}"']\arrow[dr, "p_{\rho\eta(\alpha)}"]   & & & \\ 
X_{\rho\eta(\alpha)}  \arrow[dd, "h_{\eta(\alpha)}"'] & &X_{\rho\eta(\alpha)} \arrow[dd, "\phi_{\rho\eta(\alpha),\alpha}"] \\ \\
 Y_{\eta(\alpha)} \arrow[rr, "k_{\alpha}"]  & & X_{\alpha}
 \end{tikzcd}
 \end{figure}
commutes up to homotopy. We have that $[k_{\alpha}][h_{\eta(\alpha)}][p_{\rho\eta(\alpha)}]=[\phi_{\rho\eta(\alpha),\alpha}][p_{\rho\eta(\alpha)}]$. Therefore, there exists $\alpha'\geq \rho\eta(\alpha)$ such that 
\begin{align*}
[k_{\alpha}][h_{\eta(\alpha)}] [\phi_{\alpha',\rho\eta(\alpha)}] = [\phi_{\rho\eta(\alpha),\alpha}][\phi_{\alpha',\rho\eta(\alpha)}] = [\phi_{\alpha',\alpha}]
\end{align*}
This shows that $(\rho\circ\eta, \{k_{\alpha}\circ h_{\eta(\alpha)}\})$ is equivalent to the identity. The composition $ (\rho,\{h_{\beta}\})\circ (\eta,\{k_{\alpha}\})$ is also equivalent to the identity, by the same argument. Therefore, the systems $\mathcal{X},\mathcal{Y}$ are pro-isomorphic. 
\end{proof}

Analagously, for every pointed topological space $(X,x_0)$, there exists an associated pointed inverse ANR-system, by taking normal, locally finite open covers of $X$ with distinguished element $U_{0}\in \mathcal{U}_{\alpha}$ where $x_0\in U_{0}$. The maps $p_{\alpha}$ are defined via partitions of unity subordinate to $\mathcal{U}'_{\alpha}$, where $\mathcal{U}'_{\alpha}$ denotes the modified open cover as discussed in the proof of Lemma \ref{samepsi}. Any two pointed ANR-systems associated to $(X,x_0)$ are isomorphic in pro-hCW$_{\ast}$. 

\begin{lemma} \label{cofinalsys} Let $\mathcal{X} = \{X_{\alpha}, [\phi_{\beta,\alpha}]; \mathcal{A}\}$ be an inverse ANR-system associated to $X$ and suppose $\mathcal{X}' = \{X_{\alpha'}, [\phi_{\beta',\alpha'}]; \mathcal{A}'\}$ is cofinal in $\mathcal{X}$. Then $\mathcal{X}'$ is also associated to $X$. 
\end{lemma}

\begin{proof} For any homotopy class $[f]: X\rightarrow P$ there exists $\alpha \in \mathcal{A}$ and $[h_{\alpha}]: X_{\alpha}\rightarrow P$ such that $[h_{\alpha}][p_{\alpha}] = [f]$. By cofinality, there is an $\alpha'\geq \alpha$ in $\mathcal{A}'$. We let $h'_{\alpha'}:=h_{\alpha}\phi_{\alpha',\alpha}$. Then $[h'_{\alpha'}][p_{\alpha'}] =[h_{\alpha}][\phi_{\alpha',\alpha}][p_{\alpha'}]= [h_{\alpha}][p_{\alpha}] = [f]$, as required. \\

Suppose we have two homotopy classes $[h],[k]: X_{\alpha'} \rightarrow P$ with $[h][p_{\alpha'}]= [k][p_{\alpha'}]$. There exists an index $\beta \in \mathcal{A}$ such that $[h][\phi_{\beta,\alpha'}] = [k] [\phi_{\beta,\alpha'}]$. By cofinality there is a $\beta'\in \mathcal{A}'$ with $\beta'\geq \beta$. We obtain 
\begin{align*}
[h][\phi_{\beta',\alpha'}] = [h][\phi_{\beta,\alpha'}][\phi_{\beta',\beta}] = [k] [\phi_{\beta,\alpha'}][\phi_{\beta',\beta}] = [k][\phi_{\beta',\alpha'}]
\end{align*}
\end{proof}

The analagous statement is true for pointed inverse ANR-systems.\\

Let us restrict to the case where $X$ is a compact metric space. The ANR-system $\{|\mathcal{U}_{h}|,[\phi_{h+1}];\mathbb{N}\}$ with projections $\varphi_{h}: X\rightarrow |\mathcal{U}_{h}|$ is cofinal in $\{X_{\alpha},[\phi_{\alpha',\alpha}];\mathcal{A}\}$. Therefore, it is an inverse system associated to $X$. Additionally, it is pointed (with modified projections $\varphi'_{h}$), and cofinal in the pointed inverse ANR-system of normal, locally finite open covers of $X$ with distinguished element. So, it is a pointed shape expansion of $(X,x_0)$.

\begin{notation} From now on we omit the brackets to indicate homotopy class, eg. it is understood that $\phi_{h+1}: |\mathcal{U}_{h+1}|\rightarrow |\mathcal{U}_{h}|$ denotes the homotopy class of $\phi_{h+1}$. 
\end{notation}

There are two other helpful ways to obtain an ANR-system associated to $X$. 

\begin{enumerate}
\item We can embed the compact metric space $X$ as a closed subset of an AR $P$ (eg. via the Kuratowski embedding into the Banach space $C_{b}(X)$). Consider the set $\{N\}$ of all open neighbourhoods $N$ of $X$ in $P$. The inverse system $\{N,\iota_{N',N}\}$ of ANRs with inclusions $\iota_{N',N}$ as bonding maps, denoted by $U(X,P)$, is called the  \textit{complete neighbourhood system} of $X$ in $P$. $\{N,\iota_{N',N}\}$ is an inverse system associated to $X$ with inclusion maps $\iota_{N}: X\rightarrow N$ as the canonical projections. If we equip $P$ with a metric, then the inverse sequence  $\{N_{k}, \iota_{k+1,k};\mathbb{N}\}$ of $\frac{1}{k}$-neighbourhoods of $X$ is cofinal in $\{N,\iota_{N',N}\}$. 
\item Let $X=\varprojlim_{k} X_{k}$ be the limit (in Top) of an inverse sequence $\mathcal{X}=\{X_{k},\phi_{k+1,k};\mathbb{N}\}$, where $X_{k}$ have the homotopy type of ANRs. The canonical projections $p_{k}: X\rightarrow X_{k}$, in this case, come from the definition of $X$ as a limit. 
\end{enumerate}

By Theorems $1.4$ and $1.9$ of \cite{morita1975shapes}, these two inverse ANR-systems are associated to $X$. If $x_0\in X$, then $\{(N_{k},x_0), \iota_{k+1,k};\mathbb{N}\}$ is a pointed shape expansion of $(X,x_0)$. Analagously, if $(X,x_0)=\varprojlim_{k} (X_{k},x_{k})$ is the limit (in Top$_{\ast}$) of an inverse sequence $\mathcal{X}=\{(X_{k},x_k),\phi_{k+1,k};\mathbb{N}\}$ with canonical projections $p_{k}: (X,x_0)\rightarrow (X_{k},x_{k})$, then $\mathcal{X}$ is a pointed shape expansion of $(X,x_0)$.\\

Consider now a continous map $f:X\rightarrow Y$ between topological spaces. Let $\mathcal{X}=\{X_{\alpha},\phi_{\alpha'\alpha};\mathcal{A}\},\mathcal{Y}=\{Y_{\beta},\varphi_{\beta',\beta},\mathcal{B}\}$ be shape expansions of $X$ and $Y$ respectively, with canonical projections $p_{\alpha},q_{\beta}$. We construct a morphism $S(f): \mathcal{X}\rightarrow \mathcal{Y}$ as follows. For every $\beta\in \mathcal{B}$, there exists a $\rho(\beta)\in \mathcal{A}$ and $h_{\beta}: X_{\rho(\beta)}\rightarrow Y_{\beta}$ such that $q_{\beta}\circ f= h_{\beta}\circ p_{\rho(\beta)}$ up to homotopy. The pair $(\rho,\{h_{\beta}\}):  \mathcal{X}\rightarrow \mathcal{Y}$ is a morphism in inv-hCW. It can be checked that any other choice of $(\rho',\{h'_{\beta}\})$ is equivalent to $(\rho,\{h_{\beta}\})$ in pro-hCW. Moreover, homotopic maps  $f,g$ induce (by construction) equivalent morphisms in pro-hCW. The morphism $(\rho,\{h_{\beta}\}):\mathcal{X}\rightarrow \mathcal{Y}$ constructed in Proposition \ref{prosystem} between two shape expansions of $X$ satisfies $S(\id_{X})=(\rho,\{h_{\beta}\})$.  \\

\begin{definition} (The category Shape) Let Shape be the category whose objects are pairs $(\mathcal{X},X)$, where $\mathcal{X}$ is a shape expansion of $X$. A morphism $(\mathcal{X},X)\rightarrow (\mathcal{Y},Y)$ in Shape is a morphism $f: \mathcal{X} \rightarrow \mathcal{Y}$ in pro-hCW. 
\end{definition}

\begin{definition} (The shape functor) Let $X$ be an object in hTop. The \textit{shape functor} $S$ sends $X$ to $(\mathcal{X},X)$, where $\mathcal{X}$ is a shape expansion of $X$. To define $S$ on morphisms, let $X,Y\in $ hTop, $f: X\rightarrow Y$ a homotopy class of continuous maps and let $\mathcal{X},\mathcal{Y}$ be ANR-systems associated to $X,Y$. We define $S(f):= (\rho,\{h_{\beta}\}):  \mathcal{X}\rightarrow \mathcal{Y}$. We obtain a functor
\begin{align*}
S: \text{hTop} \rightarrow \text{Shape}
\end{align*}
We say that two spaces $X,Y$ are \textit{shape equivalent} if $S(X),S(Y)$ are equivalent objects in Shape. 
\end{definition}

Analagously, there is a pointed shape category Shape$_{\ast}$ and a pointed shape functor $S:$ hTop$_{\ast}$ $\rightarrow$ Shape$_{\ast}$. 

\begin{remark} $S$ is dependent on a choice of shape expansion $\mathcal{X}$ for each $X\in$ hTop. A different set of choices results in a functor $S'$ which is naturally isomorphic to $S$. 
\end{remark}

\begin{prop} (Page 49 of  \cite{mardesic1982shape}) The shape functor $S:$ hTop $\rightarrow$ Shape restricted to the subcategory hCW of hTop is full and faithful. 
\end{prop}

This means that while the shape classification of spaces is in general weaker than the homotopy classification, on the class hCW, both classifications coincide.

\begin{definition} Let $(X,x_0)$ be a pointed topological space and $\mathcal{X} = \{(X_{\alpha}, x_\alpha), \phi_{\alpha',\alpha}; \mathcal{A}\}$ a pointed shape expansion of $(X,x_0)$. Applying the functors $\pi_{n}$ for $n\in \mathbb{N}_{0}$ to this system, we obtain the homotopy pro-groups (or pointed sets) of the pointed space $(X,x_0)$
\begin{align*}
\text{pro-}\pi_{n}(\mathcal{X},x_0) = \{\pi_{n}(X_{\alpha},x_{\alpha}), \phi_{\alpha',\alpha}; \mathcal{A}\}
\end{align*}
Let $F=(\rho,\{q_{\beta}\}): \mathcal{X} \rightarrow \mathcal{Y}$ be a pointed shape morphism between pointed shape expansions $\{(X_{\alpha}, x_\alpha), \phi_{\alpha',\alpha}; \mathcal{A}\},\{(X_{\alpha}, x_\alpha), \psi_{\beta',\beta}; \mathcal{B}\}$ of $(X,x_0),(Y,y_0)$. By applying the functors $\pi_{n}$ to $\{q_{\beta}\}$, we obtain the induced morphism of pro-groups (or pointed sets)
\begin{align*}
\text{pro-}\pi_{n}(F):\text{pro-}\pi_{n}(\mathcal{X},x_0)\rightarrow \text{pro-}\pi_{n}(\mathcal{Y},y_0) 
\end{align*}
This gives us a functor
\begin{align*}
\text{pro-}\pi_{n}: \text{Shape}_{\ast} \rightarrow \text{pro-}\mathcal{C}
\end{align*}
where $\mathcal{C}$ is the category of groups for $n\geq 1$ and pointed sets for $n=0$. 
\end{definition}

We write the pro-homotopy groups of $X$ as $\text{pro-}\pi_{n}(X,x_0):=\text{pro-}\pi_{n}(S(X),x_0)$. This is well-defined up to canonical natural isomorphism. Post-composing with the functor $\varprojlim:$ pro-$\mathcal{C}\rightarrow \mathcal{C}$ yields the shape homotopy groups of $(X,x_0)$. (Definition \ref{defcech}).\\

\subsection{Proper homotopy theory and Shape}

An important theorem in shape theory is the shape Whitehead theorem:

\begin{theorem} \label{shapewhitehead} (Theorem C in \cite{morita1974hurewicz}) A pointed shape morphism $F:S(X,x_0)\rightarrow S(Y,y_0)$ of finite shape-dimensional, connected topological spaces is a pointed shape equivalence if and only if it induces an isomorphism pro-$\pi_{n}(F)$:pro-$\pi_{n}(X,x_0)\rightarrow$ pro-$\pi_{n}(Y,y_0)$ of homotopy pro-groups for all $n\geq 0$. 
\end{theorem}

In this subsection, we will present a proof of this theorem (Corollary \ref{wseimplieswphe2}) when $X,Y$ are compact metric spaces. In doing so, we connect the proper homotopy theory of the inverse mapping telescopes $M_{\mathcal{X}}, M_{\mathcal{Y}}$ associated to pointed shape expansions of $(X,x_0)$ and $(Y,y_{0})$ to the shapes $S(X,x_0),S(Y,y_0)$. To do this, we have to assume that each space $X_{h}$ in a pointed shape expansion $\mathcal{X}=\{(X_{h},x_{h}),\phi_{h+1,h};\mathbb{N}\}$ is an actual simplicial complex, $x_0$ is a vertex, and each $\phi_{h+1,h}$ is simplicial. \\

\begin{lemma} Let $\mathcal{X}=\{(X_{h},x_{h}),\phi_{h+1,h};\mathbb{N}\}$ be a pointed ANR-expansion of $(X,x_{0})$. Let $(X'_{h},x'_{h})$ be a finite simplicial complex and $f_{h}: (X'_{h},x'_{h})\rightarrow (X_{h},x_{h})$ a pointed homotopy equivalence with inverse $f_{h}^{-1}$. Then the inverse system
\begin{align*}
\mathcal{X}=\{(X'_{h},x'_{h}), f_{h}^{-1}\phi_{h+1,h}f_{h+1};\mathbb{N}\}
\end{align*}
is associated to $X$. The analogous statement for the unpointed case holds. 
\end{lemma}

The proof of the lemma is left as an exercise in following the definitions. We can always assume that $x'_{h}$ is a vertex of $X'_{h}$, by stellar subdivision of the unique simplex whose interior contains $x'_{h}$. By simplicial approximation, we can assume that $f_{h}^{-1}\phi_{h+1,h}f_{h+1}$ is base point preserving and simplicial. \\

Recall that we have a fixed pointed shape expansion $\{(|\mathcal{U}_{h}|,x_0),\phi_{h+1};\mathbb{N}\}$ of $(X,x_0)$, with inverse mapping telescope $M_{X}$. Since we are now only interested in proper homotopy, the scaling on the height variable of $M_{X}$ does not matter. We can therefore identify it topologically with 
\begin{align*}
M_{X}:= |\mathcal{U}_{0}| \cup_{{\phi}_{1,0}} |\mathcal{U}_{1}| \times [1,2] \cup_{{\phi}_{2,1}} |\mathcal{U}_{2}| \times [2,3] \cup_{{\phi}_{3,2}} \dots
\end{align*}

Indeed, for any pointed shape expansion $\mathcal{X}=\{(X_{h},x_{h}),\phi_{h+1,h};\mathbb{N}\}$ of $(X,x_0)$, we can define 
\begin{align*}
M_{\mathcal{X}}:= \{\ast\} \cup_{\phi_{1,0}} X_{1} \times [1,2] \cup_{\phi_{2,1}} X_{2} \times [2,3] \cup_{\phi_{3,2}} \dots 
\end{align*}
where $\phi_{1,0}: X_{1} \times \{1\} \rightarrow \{\ast\}=:X_{0}$ is the unique map to a point. There is a base ray $\omega: [1,\infty)$ associated to $x_0$ which takes the value $\omega(t) = (x_{t},t)$. Assuming that all $X_{h}$ are finite simplicial complexes, $M_{\mathcal{X}}$ is a strongly locally finite, path-connected CW complex.\\

Let $(X,x_0),(Y,y_0)$ be pointed compact metric spaces with pointed shape expansions $\mathcal{X}=\{(X_{h},x_{h}),\phi_{h+1,h};\mathbb{N}\}, \mathcal{Y}=\{(Y_{h},y_{h}),\varphi_{h+1,h};\mathbb{N}\}$. Let $\omega$ and $\tau$ be the base rays associated to $x_0$ and $y_0$ respectively. Suppose that $\psi: M_{\mathcal{X}}\rightarrow M_{\mathcal{Y}}$ is a continuous, proper map. 

\begin{lemma} \label{functorP} $\psi: M_{\mathcal{X}}\rightarrow M_{\mathcal{Y}}$ induces a morphism 
\begin{align*}
\mathcal{P}(\psi):=(\rho,\{q_{h}\}): \{X_{h},\phi_{h+1,h};\mathbb{N}\} \rightarrow \{Y_{h},\varphi_{h+1,h};\mathbb{N}\}
\end{align*}
in pro-hCW. If $\psi: (M_{\mathcal{X}},\omega) \rightarrow (M_{\mathcal{Y}},\tau)$ is a map of pairs (ie. $\psi(\omega) \subset \tau$), then this morphism is pointed. 
\end{lemma}

\begin{proof} Define a finite filtration $\{L_h\}_{h\in \mathbb{N}}$ of $M_{\mathcal{X}}$ by letting $L_h := (M_{\mathcal{X}})_{[1,h]}$ be everything with height in the interval $[1,h]$. Analogously, let $Q_h:=(M_{\mathcal{Y}})_{[1,h]}$. As before, the CW complement of $L_{h}$, denoted by $L_h^{c}:=M_{\mathcal{X}} \oset{c}{-} L_{h}$, consists of all points with height $\geq h+1$, and deformation retracts onto the slice $X_{h} \times \{h+1\}$. Analogously, $Q^c_h$ deformation retracts to the slice $Y_{h}\times \{h+1\}$. \\

Let $i\in \mathbb{N}$ be given. Since $\psi$ is proper, there exists a $\rho(i)$ such that $\psi(L_{\rho(i)}^c)\subset Q_{i}^c$. We define $q_{i}$ as the homotopy class of 
\begin{align*}
q_{i}: X_{\rho(i)}&\longrightarrow Y_{i} \quad q_{i}:= \mathfrak{q}^{Y}_{i} \circ \psi_{|X_{{\rho(i)}} \times \{\rho(i)+1\}}
\end{align*}
where $\mathfrak{q}^Y_{i}$ refers to the deformation retract of $Q^c_{i}\subset M_{\mathcal{Y}}$ to height $i+1$.  We can always choose $\rho$ to be an increasing sequence. \\

We check that $(\rho,\{q_{h}\})$ actually defines a morphism in pro-hCW. Let $i\leq j$. We choose $\alpha = \rho(j)$ and obtain the diagram:

\begin{figure}[H]
\center
\begin{tikzcd}[row sep=scriptsize, column sep=scriptsize] 
& X_{\rho(j)} \arrow[dl, "\phi_{\rho(j),\rho(i)}" '] \arrow[dr, "\id"]   & & & \\ 
X_{\rho(i)} \arrow[dd, "q_{i}"'] & & X_{\rho(j)}\arrow[dd, "q_{j}"] \\ \\
 Y_{i}  & & Y_{j} \arrow[ll, "\varphi_{j,i}"]
 \end{tikzcd}
 \end{figure}

To show commutativity, observe that $\varphi_{ji} q_{j}$ is obtained from the map $\mathfrak{q}^Y_{i} \mathfrak{q}^Y_{j}  \psi_{|X_{\rho(j)} \times \{\rho(j)+1\}}$. $q_i \phi_{\rho(j),\rho(i)}$ is obtained by the map $\mathfrak{q}^Y_{i} \psi_{|X_{\rho(i)}\times \{\rho(i)+1\}} \mathfrak{q}^X_{\rho(i)}$, where $\mathfrak{q}^X_{\rho(i)}$ refers to the deformation retract of $L^c_{\rho(j)} \subset M_{\mathcal{X}}$ to the height slice $\{\rho(i)+1\}$. Since  $\psi \mathfrak{q}^X_{\rho(i)}$ is homotopic to $\psi$ we have that 
\begin{align*}
\mathfrak{q}^Y_{i} \psi \mathfrak{q}^X_{\rho(i)} \simeq \mathfrak{q}^Y_{i}  \psi = \mathfrak{q}^Y_{i} \mathfrak{q}^Y_{j} \psi
\end{align*}
as required.  By a similar calculation, one can show that a different choice of $\rho$ results in an equivalent morphism in pro-hCW. \\

If $\psi(\omega)\subset \tau$, then we have $q_{i} (x_{\rho(i)},\rho(i)+1)=\mathfrak{q}^{Y}_{i} \circ \psi(x_{\rho(i)},\rho(i)+1) = (y_{i},i+1)$. Therefore $(\rho,\{q_{h}\})$ is a pointed morphism. 
\end{proof}

%

%

Assume that $\psi$ has a partially defined proper homotopy inverse $\xi: (M_{\mathcal{Y}}\oset{c}{-} Q) \rightarrow M_{\mathcal{X}}$ in the sense of Corollary \ref{inverse}. Then we have the following:

\begin{lemma} \label{phemeansshape} The map $\mathcal{P}(\xi): \mathcal{Y}\rightarrow \mathcal{X}$ induced by $\xi: (M_{\mathcal{Y}}\oset{c}{-} Q) \rightarrow M_{\mathcal{X}}$ has the property that $\mathcal{P}(\xi) \circ \mathcal{P}(\psi)= \id_{\mathcal{X}}$ and $\mathcal{P}(\psi)\circ \mathcal{P}(\xi)= \id_{\mathcal{Y}}$ in pro-hCW. \\

If $\psi: (M_{\mathcal{X}},\omega)\rightarrow (M_{\mathcal{Y}},\tau)$ and $\xi: (M_{\mathcal{Y}}\oset{c}{-} Q, (M_{\mathcal{Y}}\oset{c}{-} Q) \cap \tau)\rightarrow (M_{\mathcal{X}}, \omega)$ are maps of pairs, such that $\xi \circ \psi \simeq \id_{|L^c}$ and $\psi\circ \xi \simeq \id_{|Q^c}$ as maps of pairs, then $\mathcal{P}(\xi)\circ \mathcal{P}(\psi)= \id_{\mathcal{X}}$ and $\mathcal{P}(\psi)\circ \mathcal{P}(\xi) = \id_{\mathcal{Y}}$ in pro-hCW$_{\ast}$.   
\end{lemma} 

From the lemma we conclude that $X$ and $Y$ are shape equivalent (resp. pointed shape equivalent). We write $X\simeq_{se} Y$ (resp. $(X,x_0)\simeq_{se} (Y,y_0)$). 

\begin{proof}

$\xi$ induces a morphism in pro-hCW $\mathcal{P}(\xi): \mathcal{Y}\rightarrow \mathcal{X}$ in the same way as $\psi$. We show that both compositions are equivalent to the identity. 
\begin{itemize}
\item $\mathcal{X}\xrightarrow{\mathcal{P}(\xi)\mathcal{P}(\psi)} \mathcal{X}$:\\
The composition $\mathcal{P}(\xi)\mathcal{P}(\psi)$ is defined as follows. Let $i\in \mathbb{N}$ be given. Since $\xi$ is proper there exists a $j$ such that $\xi(Q_{j}^c)\subset L_{i}^c$. We can assume that $Q_{j}^c\subset Q^c$. Since $\psi$ is proper there exists a $k\geq i$ such that $\psi(L_{k}^c)\subset Q_{j}^c$. So we have 
\begin{align*}
\xi\psi(L_{k}^{c})\subset \psi(Q_{j}^c) \subset L_{i}^c
\end{align*}
Choose additionally $k$ large enough so that $L^c_{k}\subset L^{c}$ and the homotopy $H$ from $\xi\psi_{|L_{k}^c} $ to $\id_{|L_{k}^c}$ has $H(L_{k}^c \times [0,1])\subset L_{i}^c$. We let $\rho(i) = k$ and define $q_{i}$ by taking the bottom slice of $L_{k}^c$, applying $\xi\psi$, then retracting to the bottom of $L_{i}^c$, ie. 
\begin{align*}
q_{i}: X_{k}&\longrightarrow {X}_{i} \quad q_{i}:= \mathfrak{q}^{X}_{i} \xi \psi_{|X_{k} \times \{k+1\}}
\end{align*}
The composition $\mathcal{P}(\xi)\mathcal{P}(\psi): \mathcal{X}\rightarrow \mathcal{X}$ is represented by the morphism $(\rho,\{q_{h}\})$. \\

Consider $\mathfrak{q}^{X}_{i} H$ restricted to ${X}_{k} \times \{k+1\}$. We have that $\mathfrak{q}^{X}_{i}\xi \psi_{|{X}_{k} \times \{k+1\}}$ is homotopic to $\mathfrak{q}^{X}_{i}$. Thus the diagram (under appropriate identifications)
\begin{figure}[H]
\center
\begin{tikzcd}[row sep=scriptsize, column sep=scriptsize] 
& X_{k} \arrow[dl, "\id" '] \arrow[dr, "\phi_{k,i}"]   & & & \\ 
X_{k}  \arrow[dr, "q_i \simeq \mathfrak{q}^{X}_{i} = \phi_{k,i}" '] & & {X}_{i}\arrow[dl, "\id"] \\ 
 &X_{i} & &&  
 \end{tikzcd}
 \end{figure}
 commutes. This shows that $\mathcal{P}(\xi)\mathcal{P}(\psi)$ is equivalent to $\id_{\mathcal{X}}$ in pro-hCW. 

\item $\mathcal{Y}\xrightarrow{\mathcal{P}(\psi)\mathcal{P}(\xi)} \mathcal{Y}$:\\
The argument for $\mathcal{P}(\psi)\mathcal{P}(\xi) = \id_{\mathcal{Y}}$ follows analogously. 
\end{itemize}

If $\psi,\xi$ are maps of pairs, then $\mathcal{P}(\psi),\mathcal{P}(\xi)$ are morphisms in pro-hCW$_{\ast}$ by the same formulas. If the homotopies $H,G$ realising $\xi \circ \psi \simeq \id_{|L^c}$ and $\psi\circ \xi \simeq \id_{|Q^c}$ satisfy $H(\omega \times [0,1])\subset \omega$ and $G(\tau \times [0,1]) \subset \tau$, then all morphisms and homotopies in the proof are pointed. 
\end{proof}

\begin{remark} Clearly if $\xi$ is an actual proper homotopy inverse then Lemma \ref{phemeansshape} also holds. However, a partially-defined inverse is sufficient, since we have the freedom of going up the telescope in our choice of indices. 
\end{remark}

In summary, we have shown that a proper map $M_{\mathcal{X}}\rightarrow M_{\mathcal{Y}}$ induces a shape morphism $\mathcal{X}\rightarrow \mathcal{Y}$, with a proper homotopy equivalence inducing a shape equivalence. We now construct the correspondence in the other direction.

\begin{lemma} \label{modify} Let $\mathcal{X}=\{X_{h}, \phi_{h+1,h}; \mathbb{N}\}$ and $\mathcal{Y}=\{Y_{h}, \varphi_{h+1,h}; \mathbb{N}\}$ be objects in pro-hCW. A morphism $(\rho,\{q_{h}\}): \mathcal{X}\rightarrow \mathcal{Y}$ is equal in pro-hCW to a morphism $(\rho',\{q'_{h}\})$ with the properties
\begin{itemize}
\item $\rho': \mathbb{N}\rightarrow \mathbb{N}$ is a strictly increasing function. 
\item The diagram

\begin{figure}[H]
\center
\begin{tikzcd}
&X_{\rho'(h)}   \arrow{d}{q'_{h}}
&&X_{\rho'(h+1)} \arrow{d}{q'_{h+1}} \arrow[ll, "\phi_{\rho'(h+1)\rho'(h)}"']
  \\
& Y_{h}
&&Y_{h+1} \arrow{ll}{\varphi_{h+1,h}}
\end{tikzcd}
\end{figure}
commutes up to homotopy for all $h\in \mathbb{N}$. 
\end{itemize}
The analogous statement is true for objects and morphisms in pro-hCW$_{\ast}$. 
\end{lemma}

\begin{proof}By the definitions, there exists a $j\geq \rho(1),\rho(2)$ such that
\begin{align*}
\varphi_{2,1}q_2\phi_{j,\rho(2)} \simeq q_1 \phi_{j,\rho(1)}
\end{align*}
Define $\rho'(1) = \rho(1)$, $q'_{1}=q_{1}$, and $\rho'(2) = j$ and $q'_{2}: X_{\rho'(2)}\rightarrow Y_{2}$ as 
\begin{align*}
q'_{2}:= q_{2} \phi_{j,\rho(2)}.  
\end{align*}
Then the diagram

\begin{figure}[H]
\center
\begin{tikzcd}
&X_{\rho(1)} = X_{\rho'(1)}  \arrow{d}{q'_{1}}
&&X_{\rho'(2)} \arrow{d}{q'_{2}=q_{2} \phi_{\rho'(2)\rho(2)}} \arrow[ll, "\phi_{\rho'(2)\rho'(1)}"']
  \\
& Y_{1}
&&Y_{2} \arrow{ll}{\varphi_{2,1}}
\end{tikzcd}
\end{figure}
commutes up to homotopy by construction. Suppose for the induction step that $\rho'(h)$ and $q'_h: X_{\rho'(h)}\rightarrow Y_{h}$ has been defined and satisfes the properties in the lemma. Consider $\rho(h+1)$ and $q_{h+1}: X_{\rho(h+1)}\rightarrow Y_{h+1}$. From the definitions there exists $j\geq \rho(h), \rho(h+1)$ such that 
\begin{align*}
\varphi_{h+1,h}q_{h+1}\phi_{j,\rho(h+1)} \simeq q_h \phi_{j,\rho(h)}
\end{align*}
If we choose any $j'\geq j$ the diagram

\begin{figure}[H]
\center
\begin{tikzcd}[row sep=scriptsize, column sep=scriptsize] 
& X_{j'} \arrow[dl, "\phi_{j\rho(h)}\phi_{j'j}\simeq\phi_{j'\rho(h)}" '] \arrow[dr, "\phi_{j'\rho(h+1)}\simeq\phi_{j\rho(h+1)}\phi_{j'j}"]   & & & \\ 
X_{\rho(h)} \arrow[dd, "q_h"'] & &X_{\rho(h+1)}\arrow[dd, "q_{h+1}"] \\ \\
 Y_{h}  & & Y_{h+1}\arrow[ll, "\varphi_{h+1,h}"]
 \end{tikzcd}
 \end{figure}

also commutes up to homotopy since the bonding maps $\phi$ are compatible. Therefore we can assume that $j' > \rho'(h)$ and we let $\rho'(h+1):=j'$ and $q'_{h+1}= q_{h+1} \phi_{j' \rho(h+1)}$. The diagram

\begin{figure}[H]
\center
\begin{tikzcd}
&X_{\rho'(h)}   \arrow{d}{q'_{h}}
&&X_{\rho'(h+1)} \arrow{d}{q'_{h+1}} \arrow[ll, "\phi_{\rho'(h+1)\rho'(h)}"']
  \\
& Y_{h}
&&Y_{h+1} \arrow{ll}{\varphi_{h+1,h}}
\end{tikzcd}
\end{figure}

commutes up to homotopy since  
\begin{align*}
q'_{h}\phi_{\rho'(h+1)\rho'(h)}\simeq q_{h}\phi_{\rho'(h)\rho(h)}\phi_{\rho'(h+1)\rho'(h)}\simeq q_{h} \phi_{\rho'(h+1)\rho(h)}\simeq \varphi_{h+1,h} q_{h+1} \phi_{j' \rho(h+1)}\\
\simeq \varphi_{h+1,h}q'_{h+1} 
\end{align*}

Therefore $(\rho',\{q'_{h}\})$ is defined and satisfies the desired properties. To check that it is actually a morphism in pro-hCW requires that for all $i\leq j$ there exists a $k\geq \rho'(i), \rho'(j)$ such that the diagram

\begin{figure}[H]
\center
\begin{tikzcd}[row sep=scriptsize, column sep=scriptsize] 
& X_{k} \arrow[dl, "\phi_{k \rho'(i)}" '] \arrow[dr, "\phi_{k\rho'(j)}"]   & & & \\ 
X_{\rho'(i)} \arrow[dd, "q'_{i}"'] & &X_{\rho'(j)}\arrow[dd, "q'_{j}"] \\ \\
 Y_{i}  & & Y_{j}\arrow[ll, "\varphi_{ji}"]
 \end{tikzcd}
 \end{figure}
 
 is homotopy commutative. To this end, choose $k=\rho'(j)$. Since $\rho'(j) \geq \rho'(l), \rho(l)$ for all $i\leq l\leq j$ we have the following homotopy commutative diagram:

\begin{figure}[H]
\center
\begin{tikzcd}[row sep=scriptsize, column sep=scriptsize] 
&& & X_{\rho'(j)} \arrow[d] \arrow[dr] \arrow[dll]  & & & \\ 
& X_{\rho(i)} \arrow[dd, "q_{i}"']  &\dots & X_{\rho(j-1)} \arrow[dd, "q_{j-1}"']  &X_{\rho(j)}\arrow[dd, "q_{j}"] \\ \\
  & Y_{i}   & \dots \arrow[l] & Y_{j-1}\arrow[l]   & Y_{j}\arrow[l]
 \end{tikzcd}
 \end{figure}

 where all the unlabelled arrows are bonding morphisms. This, combined with the definition of $q'_{h}$ as $q_{h}\phi_{\rho'(h)\rho(h)}$ for all $h\in \mathbb{N}$ gives us the desired result. \\
 
It remains to show that $(\rho,\{q_{h}\}) =(\rho',\{q'_{h}\})$ in pro-hCW. To do this, we need that for all $i\in \mathbb{N}$ there exists a $j\geq \rho(i), \rho'(i)$ such that the diagram

\begin{figure}[H]
\center
\begin{tikzcd}[row sep=scriptsize, column sep=scriptsize] 
& X_{j} \arrow[dl, "\phi_{j\rho(i)}" '] \arrow[dr, "\phi_{j\rho'(i)}"]   & & & \\ 
X_{\rho(i)}  \arrow[dr, "q_{i}" '] & & X_{\rho'(i)}\arrow[dl, "q'_{i}"] \\ 
 &Y_{i} & &&  
 \end{tikzcd}
 \end{figure}
 
 commutes up to homotopy. We simply pick $j=\rho'(i)$ and the diagram commutes by construction. \\

If $\mathcal{X},\mathcal{Y}\in$ pro-hCW$_{\ast}$ the same proof works by letting all maps and homotopies be pointed. 
\end{proof}

We say that an ANR-system $\mathcal{X} = \{(X_{h},x_0), \phi_{h+1,h}; \mathbb{N}\}$, where all $X_{h}$ are simplicial complexes, has dimension $\leq n $ if all $X_{h}$ have dimension $\leq n$.

\begin{theorem} \label{wseimplieswphe1} Let $\mathcal{X}:=\{(X_{h},x_0), \phi_{h+1,h}; \mathbb{N}\},\mathcal{Y}:=\{(Y_{h}, y_0), \varphi_{h+1,h}; \mathbb{N}\}$ be inverse systems associated with the compact, based metric spaces $(X,x_{0})$ and $(Y,y_{0})$. A pointed shape morphism $(\rho,\{q_{h}\}): \mathcal{X}\rightarrow \mathcal{Y}$ induces a continuous, proper, cellular map of pairs $\psi: (M_{\mathcal{X}}, \omega)\rightarrow (M_{\mathcal{Y}}, \tau)$ such that $\mathcal{P}(\psi)= (\rho,\{q_{h}\})$.\\

Assume that $X,Y$ are connected and that $\mathcal{X},\mathcal{Y}$ are finite-dimensional expansions. $\psi$ has a proper homotopy inverse of pairs $\xi$ if and only if the induced map on end homotopy groups $\psi: \pi_{n}^{e}(M_{\mathcal{X}},\omega)\rightarrow \pi_{n}^{e}(M_{\mathcal{Y}}, \psi \omega)$ is an isomorphism for all $n\geq 0$. This occurs exactly when the induced map pro-$\pi_{n}(\rho,\{q_{h}\}):$ pro-$\pi_{n}(X,x_0)\rightarrow$ pro-$\pi_{n}(Y,y_0)$ on pro-homotopy groups is an isomorphism for all $n\geq 0$. 
\end{theorem}

\begin{proof} We can assume that $(\rho,\{q_{h}\})$ has the form as in Lemma \ref{modify}.  Since $\rho$ is a strictly increasing function the object $\mathcal{X}':=\{X_{\rho(h)}, \phi_{\rho(h+1)\rho(h)}; \mathbb{N}\}$ is pro-isomorphic to $\mathcal{X}=\{X_{h}, \phi_{h+1, h}; \mathbb{N}\}$ and the morphism $(\rho,\{q_{h}\})$ restricts to $\mathcal{X}'\rightarrow \mathcal{Y}$. We define the inverse mapping telescopes $M_{\mathcal{X}}, M_{\mathcal{X}'}, M_{\mathcal{Y}}$ as follows:
\begin{align*}
M_{\mathcal{X}}&:= \{\ast\} \cup X_{1} \times [0,1] \cup_{\phi_{2,1}} X_{2} \times [0,1] \cup_{\phi_{3,2}} \dots\\
M_{\mathcal{X}'}&:= \{\ast\} \cup X_{\rho(1)} \times [0,1] \cup_{\phi_{\rho(2)\rho(1)}} X_{\rho(2)} \times [0,1] \cup_{\phi_{\rho(3)\rho(2)}} \dots \\
M_{\mathcal{Y}}&:= \{\ast\} \cup Y_{1} \times [0,1] \cup_{\varphi_{2,1}} Y_{2} \times [0,1] \cup_{\varphi_{3,2}} \dots
\end{align*}
There is a proper homotopy equivalence $r: M_{\mathcal{X}}\rightarrow M_{\mathcal{X}'}$, obtained by collapsing all intermediate cylinders, with inverse $\iota: M_{\mathcal{X}'}\rightarrow M_{\mathcal{X}}$ . Therefore it suffices to construct a map $\psi': M_{\mathcal{X}'}\rightarrow M_{\mathcal{Y}}$ with the desired properties. $\psi'$ is defined piecewise for $h\in \mathbb{N}$ by 
\begin{align*}
\psi'_{h}: X_{\rho(h)} \times [0,1] \rightarrow Y_{h} \times [0,1] \cup_{\phi_{h,h-1}} Y_{h-1} \times \{1\}
\end{align*}
 \[ \psi'_{h}(x,t) := \begin{cases} 
          [(q_{h}(x), t)]& t\in [\frac{1}{2},1] \\
          [(q_{h}(x), 2t-\frac{1}{2})]& t\in [\frac{1}{4},\frac{1}{2}]\\ 
	 [(H^h(x,4t),1)]& t\in [0,\frac{1}{4}]\\ 
       \end{cases}  \] 

where $H^h: X_{\rho(h)} \times [0,1]\rightarrow Y_{h-1}$ is a homotopy between $H^h_{0} = q_{h-1} \phi_{\rho(h)\rho(h-1)}$ and $H^h_1 = \varphi_{h,h-1} q_{h}$. This glues to a continuous, proper map $\psi'$. Let $\psi:= \psi' \circ r$. Since all maps and homotopies are pointed, $\psi: (M_{\mathcal{X}},\omega)\rightarrow (M_{\mathcal{Y}},\tau)$ is a map of pairs. By construction we have $\mathcal{P}(\psi) = (\rho,\{q_{h}\})$. Additionally, by the relative simplicial approximation theorem, we can assume that all $q_{h}, H^h$ are simplicial. Then $\psi'$ is cellular map.\\

For the second half of the theorem, we require the following claim: 

\begin{claim} Let $X$ be a compact metric space. $X$ is connected exactly when the inverse mapping telescope $M_{X}$ has one end. 
\end{claim}

Suppose $X$ is connected. Then every $|\mathcal{U}_{h}|$ is path-connected. Therefore the inverse mapping telescope $M_{X}$ has the property that $M_{X} \oset{c}{-} L_i$ is path-connected for all $i\in \mathbb{N}$, which implies that $M_{X}$ has one end. Conversely, if $X$ not connected there exist $U,V$ open and disjoint such that $X = U\cup V$. Since the sequence of open covers $\mathcal{U}_{h}$ is cofinal, there exists a $h'$ large enough such that $\mathcal{U}_{h}$ is refinement of $U\cup V$ for all $h\geq h'$. This shows that $M_{X}$ has at least $2$ ends. \\

The number of ends of $M_{X}$ is a shape invariant of $X$. Therefore for any shape expansion $\mathcal{X}$ of $X$, the inverse mapping telescope $M_{\mathcal{X}}$ has one end. \\

Let $\omega' = r\omega$. Consider the mapping cylinder of $\psi'$, which we denote by 
\begin{align*}
M_{\psi'}:= M_{\mathcal{X}'} \times [0,1] \cup_{\psi'} M_{\mathcal{Y}}
\end{align*}
This is a finite-dimensional, strongly locally finite (this uses the properness of $\psi'$), path-connected CW complex, with $M_{\mathcal{X}'}\times \{0\}$ as a subcomplex. The commutative diagrams in the proof of \ref{deforetract} with $Y = M_{\psi'}$, $A = M_{\mathcal{X}'} \times \{0\}$ and $\{L_{i}\}_{i\in \mathbb{N}}$ a finite filtration of   $M_{\psi'}$, show that $\psi':(M_{\mathcal{X}'},\omega)\rightarrow (M_{\mathcal{Y}},\tau)$ induces an isomorphism on end homotopy groups for all $n\geq 0$ exactly when the induced morphisms
\begin{align*}
{\varprojlim}^1 \pi_{n+1}((L_i\cap A)^c,\omega'(i))& \rightarrow {\varprojlim}^1 \pi_{n+1}(L^c_i,\omega'(i)) \\
\varprojlim \pi_{n}((L_i\cap A)^c,\omega'(i))  &\rightarrow \varprojlim \pi_{n}(L^c_i,\omega'(i)) 
\end{align*}
are isomorphisms for all $n\geq 0$. \\

By post-composing with the deformation retract $M_{\psi'}\rightarrow M_{\mathcal{Y}}$, this is equivalent (after suitable identifications) to the condition that 
\begin{align*}
{\varprojlim}^1 \text{pro-}\pi_{n+1}(\rho,\{q_{h}\})&: {\varprojlim}^1 \text{pro-}\pi_{n+1}(X,x_0)\rightarrow {\varprojlim}^1 \text{pro-}\pi_{n+1}(Y,y_0) \\
 \varprojlim \text{pro-}\pi_{n}(\rho,\{q_{h}\})&: \varprojlim \text{pro-}\pi_{n}(X,x_0)\rightarrow \varprojlim \text{pro-}\pi_{n}(Y,y_0)
\end{align*}
are isomorphisms for all $n\geq 0$.  This occurs exactly when pro-$\pi_{n}(\rho,\{q_{h}\}):$ pro-$\pi_{n}(X,x_0)\rightarrow$ pro-$\pi_{n}(Y,y_0)$ is an isomorphism for all $n\geq 0$, by applying the following lemma:

\begin{lemma} (Lemma 3.8 in \cite{melikhov2009steenrod}) A collection of homomorphisms $f_{i}:G_{i}\rightarrow H_{i}$ between inverse sequences of countable groups commuting with the bonding maps is a pro-isomorphism if and only if $\varprojlim f_{i}: \varprojlim G_{i}\rightarrow \varprojlim H_{i}$ and ${\varprojlim}^1 f_{i}: {\varprojlim}^1 G_{i}\rightarrow {\varprojlim}^1 H_{i}$ are bijections.  
\end{lemma}

and using the fact that connectedness of $X,Y$ implies pro-$\pi_{0}$ is the trivial system. \\

By Corollary \ref{inverse} and Corollary \ref{inverse2}, there is a proper homotopy $G: M_{\psi'} \times [0,1]\rightarrow M_{\psi'}$ such that $G_{0} = \id$ and $G_{1}(M_{\psi'})\subset  M_{\mathcal{X}'} \times \{0\}$. Let $K: = ((\omega' \times [0,1])\cup \tau)$.  We can assume $G$ satisfies $G(K \times [0,1]) \subset K$ by the following argument: since $\psi'_{|\omega'}: \omega' \rightarrow \tau$ is a proper homotopy equivalence, there exists a homotopy $G^{\omega'}: K \times [0,1] \rightarrow K$ such that $G^{\omega'}_{0} = \id$ and $G^{\omega'}_{1}(K) \subset \omega' \times \{0\}$. We can extend $G^{\omega'}$ to $G$ now by inducting over cells not in $K\times [0,1]$ in the usual fashion. \\

 This gives us a continuous proper map of pairs $\xi': (M_{\mathcal{Y}},\tau) \rightarrow (M_{\mathcal{X'}},\omega')$ which is a proper homotopy inverse of pairs to $\psi'$. The composition $\xi:=\iota \circ \xi' : (M_{\mathcal{Y}},\tau) \rightarrow (M_{\mathcal{X}},\omega)$ is a proper homotopy inverse of pairs to $\psi$. 
\end{proof}

Before we prove the shape Whitehead theorem, we need to discuss dimension.

\begin{definition} (Lebesgue covering dimension) (Definition $2.1$ in \cite{DRANISHNIKOV2018429}) The \textit{multiplicity} of an open cover $\mathcal{V}$ is the maximal number of elements of $\mathcal{V}$ having a common point. The \textit{Lebesgue covering dimension} of a topological space $X$ does not exceed $n$ ($\dim X\leq n$) if for every open cover $\mathcal{U}$ of $X$ there is an open refinement $\mathcal{V}$ of $\mathcal{U}$ with multiplicity $\leq n+1$. Define $\dim X$ as the minimal number such that $\dim X\leq n$ and $\infty$ otherwise.
\end{definition}

\begin{definition} (Shape dimension) (Definition $11$ in \cite{gevorgyan2021shape}) A topological space $X$ has \textit{shape dimension} $sd X\leq n$ if there exists an associated ANR-system $\mathcal{X}=\{X_{\alpha}, \phi_{\alpha',\alpha}; \mathcal{A}\}$ such that each $X_{\alpha}$, $\alpha \in \mathcal{A}$ is \textit{homotopy dominated} by an ANR of dimension $\leq n$, ie. there exists an ANR $Y_{\alpha}$ of dimension $\dim Y_{\alpha}\leq n$ and maps $f_{\alpha}: X_{\alpha}\rightarrow Y_{\alpha}$ and $g_{\alpha}: Y_{\alpha} \rightarrow X_{\alpha}$ such that
\begin{align*}
g_{\alpha}f_{\alpha} \simeq \id_{X_{\alpha}}
\end{align*} 
Define $sd X$ as the minimal number such that $sd X\leq n$ and $\infty$ otherwise.
\end{definition}

\begin{theorem} (Theorem $2.2$ of \cite{DRANISHNIKOV2018429}) For a compact metric space $X$, its shape dimension and Lebesgue covering dimension coincide.
\end{theorem} 

From this we obtain the following proposition:

\begin{prop} \label{finitedim} Let $(X,x_0)$ be a compact metric space with pointed shape expansion $\mathcal{X}:=\{(|\mathcal{U}_{h}|,x_0),\phi_{h+1,h};\mathbb{N}\}$.  Let $X$ have finite shape dimension $sd X = n$. Then there exists a pointed shape expansion $\mathcal{X'}:=\{(\mathcal{V}_{h}, V_{0h}),\tilde{\phi}_{h+1,h}; \mathbb{N}\}$ of dimension $\leq n$ obtained from a cofinal sequence of open covers with distinguished element containing $x_0$. There is a proper homotopy equivalence of pairs $\zeta_{\mathcal{X'}\mathcal{X}}: (M_{\mathcal{X}'}, \tau) \rightarrow (M_{X},\omega)$. 
\end{prop}

\begin{proof}  
Let $\{\mathcal{U}'_{h}\}_{h\in \mathbb{N}}$ be our fixed sequence of (modified) open covers of $X$. For ease of notation in this proof, we drop the primes on $\mathcal{U}_{h}$ entirely. It is understood that $\mathcal{U}_{h}$ is the modified sequence of open covers (ie. $B(x_0,\frac{1}{2^h})$ is the only set in $\mathcal{U}_{h}$ containing $x_0$). We also stop writing $|\mathcal{U}_{h}|$ and just use $\mathcal{U}_{h}$ even when we mean the geometric realisation of the nerve.\\

Let $\V_{0} = \{\ast\}$. Let $\psi_{0}: \U_{h_1}:=\U_1\rightarrow \V_0$ be the obvious map to a point. Choose $\V_{1}$ to be a refinement of $\U_1$ with multiplicity $\leq n+1$.  By a further modification if necessary,  we can assume that  $\V_{1}$ has a distinguished element $V_{01}$ as the unique $V\in \V_{1}$  containing $x_0$.  There is a simplicial refinement map  $\varphi_1: \V_1\rightarrow \U_1$ which sends the vertex $V_{01}$ to $x_0$. We define $\tilde{\phi}_{1,0}: \V_1\rightarrow \V_0$ to be $\psi_{0}\circ \varphi_1$. Since $\{\U_{h}\}_{h\in \mathbb{N}}$ is cofinal, there exists $h_2 > 1$ such that $\U_{h_2}$ refines $\V_1$ with refinement map $\psi_1: \U_{h_2}\rightarrow \V_1$. We have that $\psi_1(x_0) = V_{01}$. We can choose $\V_{2}$ with multiplicity $\leq n+1$ and distinguished element $V_{02}$ so that it refines $\U_{h_2}$ with pointed refinement map $\varphi_{2}: \V_{2}\rightarrow \U_{h_2}$. We define $\tilde{\phi}_{2,1}: \V_2\rightarrow \V_1$ as $\psi_1\circ \varphi_2$. Continuing iteratively we have
\begin{align*}
\psi_{j-1}\circ \varphi_j=:\tilde{\phi}_{j,j-1}: \V_{j} \rightarrow \V_{j-1}
\end{align*}
for all $j\geq 1$ where 
\begin{align*}
\varphi_{j}&: \V_{j}\rightarrow \U_{h_{j}}\\
\psi_{j-1}&: \U_{h_{j}} \rightarrow \V_{j-1} 
\end{align*}
are pointed, simplicial refinement maps.\\

Since $\{(|\mathcal{U}_{h}|,x_0),\phi_{h+1,h};\mathbb{N}\}$ is cofinal in the system of locally finite, normal covers of $X$ with distinguished element, so is $\mathcal{X}'=\{(\mathcal{V}_{h}, V_{0h}),\tilde{\phi}_{h+1,h}; \mathbb{N}\}$. Therefore, by Lemma \ref{cofinalsys}, $\mathcal{X}'$ is a pointed shape expansion of $(X,x_0)$. \\

We construct the inverse mapping telescope  $M_{\mathcal{X}'}$ as
\begin{align*}
M_{\mathcal{X}'}:= \V_{0} \cup_{\tilde{\phi}_{1,0}} \V_{1} \times [0,1] \cup_{\tilde{\phi}_{2,1}} \V_{2} \times [0,1] \cup_{\tilde{\phi}_{3,2}} \dots 
\end{align*}
Let $\tau$ be the base ray of $M_{\mathcal{X}'}$ corresponding to the distinguished points $V_{0h}$ in the simplicial complex $\mathcal{V}_{h}$. \\

Likewise, we can identify $M_{\mathcal{X}}=M_{X}$ topologically with 
\begin{align*}
M_{X}:= \U_{0} \cup_{{\phi}_{1,0}} \U_{1} \times [0,1] \cup_{{\phi}_{2,1}} \U_{2} \times [0,1] \cup_{{\phi}_{3,2}} \dots
\end{align*}
with base ray $\omega$ corresponding to $x_0$. Observe that $\varphi_j\circ \psi_{j}: \U_{h_{j+1}}\rightarrow \U_{h_{j}}$ is pointed homotopic to $\phi_{h_{j+1}, h_{j}}$ for all $j\geq 1$. \\

Our goal now is to extend the algebraic data $\varphi_{j}$ and $\psi_{j}$ to a proper homotopy equivalence $\zeta_{\mathcal{X}'\mathcal{X}}: (M_{\mathcal{X}'},\tau) \rightarrow (M_X,\omega)$ of pairs, defined on the whole telescope by "filling in" the gaps. To do this, we define $\zeta_{\mathcal{X}'\mathcal{X}}$ piecewise:
\begin{align*}
\zeta^{j}_{\mathcal{X}'\mathcal{X}}: \V_{j+1} \times [0,\frac{1}{2}]\cup_{\tilde{\phi}} \V_{j} \times [\frac{1}{2},1] \longrightarrow \U_{h_{j+1}} \times [0,\frac{1}{2}] \cup_{\phi} (\bigcup_{h_{j}<h<h_{j+1}} \U_{h} \times [0,1]) \cup_{\phi} \U_{h_j} \times [\frac{1}{2},1] 
\end{align*}

The definition is best explained with a picture (Figure \ref{fig:prophom1}), but it does the following on specific slices: 
\begin{align*}
\zeta^{j}_{\mathcal{X}'\mathcal{X}}:& \V_{j+1} \times \{\frac{1}{2}\} \rightarrow \U_{h_{j+1}} \times \{\frac{1}{2}\}  \quad   &:=\varphi_{j+1}\\
&\V_{j+1} \times \{\frac{3}{8}\} \rightarrow \U_{h_{j+1}-1} \times \{1\}  &:= \phi_{h_{j+1}, h_{j+1}-1}\varphi_{j+1}\\
&\V_{j+1} \times \{\frac{1}{4}\} \rightarrow \U_{h_{j}} \times \{1\} &:= \phi_{h_{j+1}, h_{j}}\varphi_{j+1}\\
&\V_{j+1} \times \{0\} \rightarrow \U_{h_j} \times \{\frac{1}{2}\} &:= \varphi_{j} \psi_{j}\varphi_{j+1} \\
&\V_{j} \times [\frac{1}{2},1] \rightarrow \U_{h_j} \times \{\frac{1}{2}\} &:= \varphi_j
\end{align*}
On the segments $\V_{j+1} \times [\frac{3}{8}, \frac{1}{2}]$ and $\V_{j+1} \times [\frac{1}{4}, \frac{3}{8}]$ we fill in with the obvious deformation retract on $M_{X}$. On $\V_{j+1} \times [0,\frac{1}{4}]$ we fill in with a fixed choice of linear homotopy $\phi_{h_{j+1}, h_{j}}\simeq \varphi_{j} \psi_{j}$. One can check that all this together is a well-defined map compatible with gluing. The $\zeta^{j}_{\mathcal{X}'\mathcal{X}}$ are also compatible with each other on boundary components. On the segment $\V_{1} \times [0,\frac{1}{2}]$ we take the cone of the map $\varphi_1$. We therefore obtain a globally defined map $\zeta_{\mathcal{X}'\mathcal{X}}$, which is continuous and proper. \\

\begin{figure}[H]
\centering
  \centering
  \includegraphics[width=0.85\linewidth]{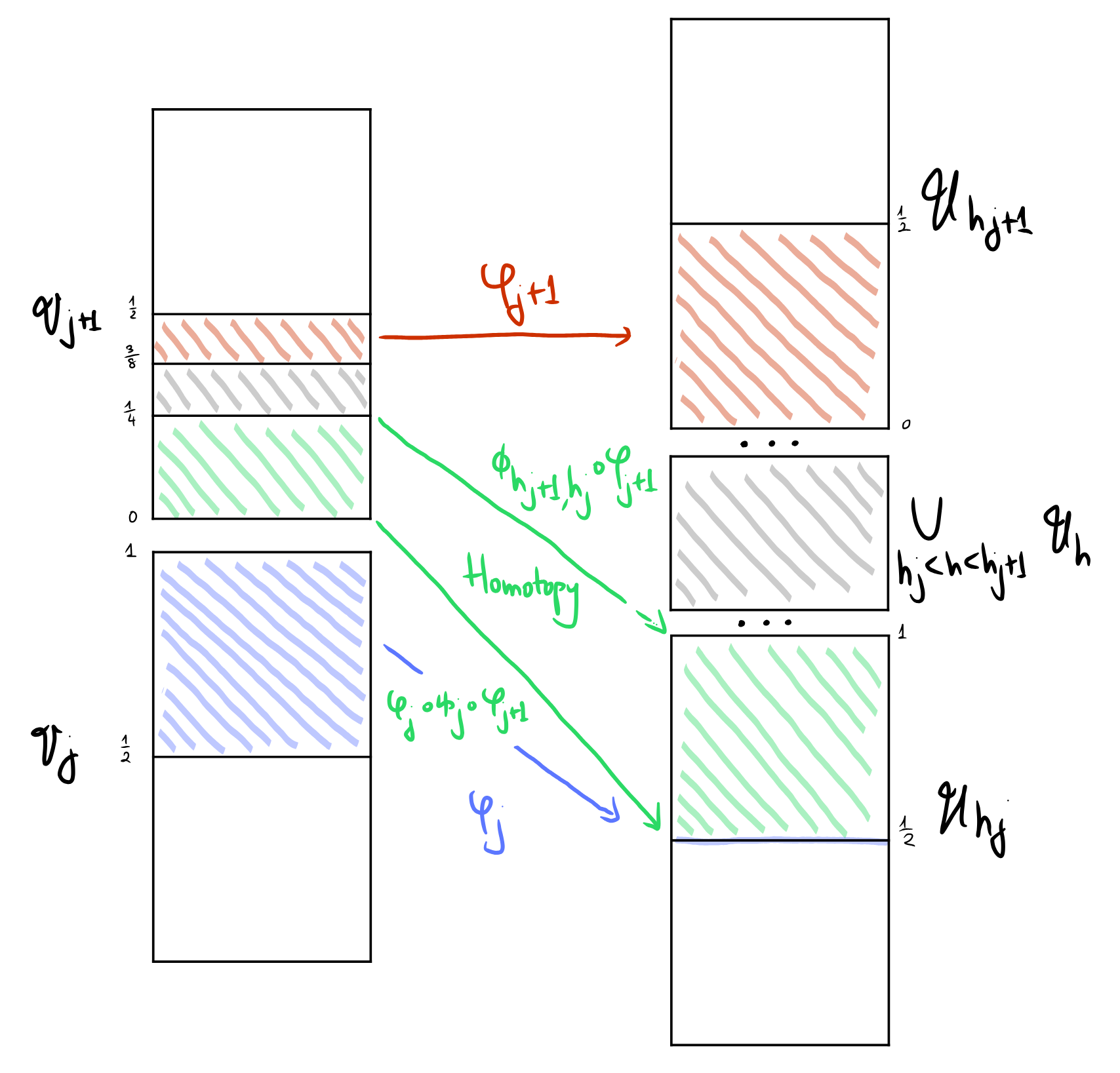}
  \caption{The definition of $\zeta^j_{\mathcal{X}'\mathcal{X}}$.  }
  \label{fig:prophom1}
\end{figure}

We now define $\zeta_{\mathcal{X} \mathcal{X}'}: M_{X} \rightarrow M_{\mathcal{X}'}$ similarly. See Figure \ref{fig:prophom2}. 

\begin{align*}
\zeta^{j}_{\mathcal{X} \mathcal{X}'}:  \U_{h_{j+1}} \times [0,\frac{1}{2}] \cup_{\phi} (\bigcup_{h_{j}<h<h_{j+1}} \U_{h} \times [0,1]) \cup_{\phi} \U_{h_j} \times [\frac{1}{2},1] \longrightarrow \V_{j} \times [0,\frac{1}{2}]\cup_{\tilde{\phi}} \V_{j-1} \times [\frac{1}{2},1] 
\end{align*}

\begin{align*}
\zeta^{j}_{\mathcal{X} \mathcal{X}'}: & \U_{h_{j+1}} \times \{\frac{1}{2}\} \rightarrow \V_{j} \times \{\frac{1}{2}\}  \quad   &:=\psi_{j}\\
&\U_{h_{j+1}} \times \{\frac{1}{4}\} \rightarrow  \V_{j-1} \times \{1\}  &:= \tilde{\phi}_{j,j-1} \psi_{j}\\
&\U_{h_{j+1}} \times \{0\} \rightarrow \V_{j-1} \times \{\frac{1}{2}\} &:= \psi_{j-1} \phi_{h_{j+1},h_j} \\
&\U_{h_{j}} \times [\frac{1}{2}, 1] \rightarrow \V_{j-1} \times \{\frac{1}{2}\} &:= \psi_{j-1}\\
\end{align*}

\begin{figure}
\centering
  \centering
  \includegraphics[width=0.85\linewidth]{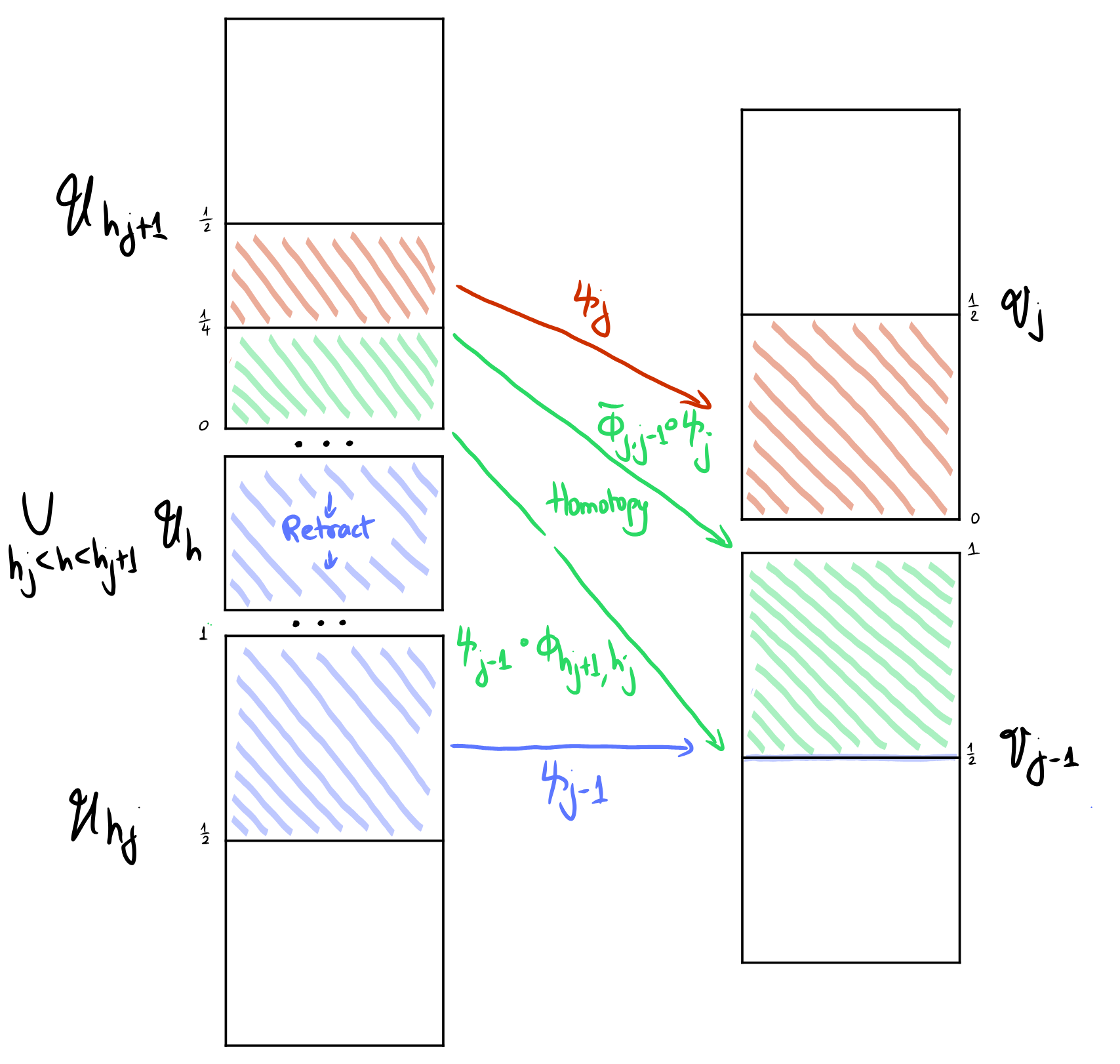}
  \caption{The definition of $\zeta^j_{\mathcal{X} \mathcal{X}'}$.}
  \label{fig:prophom2}
\end{figure}

On the segment $\U_{h_{j+1}} \times [\frac{1}{4}, \frac{1}{2}]$ we fill in with the obvious deformation retract on $M_{\mathcal{X}'}$. On $\U_{h_{j+1}} \times [0,\frac{1}{4}]$ we fill in with a linear homotopy $\tilde{\phi}_{j,j-1}\psi_{j}=\psi_{j-1}\varphi_{j} \psi_{j} \simeq \psi_{j-1} \phi_{h_{j+1},h_j}$ coming from a fixed choice of linear homotopy $\phi_{h_{j+1}, h_{j}}\simeq \varphi_{j} \psi_{j}$. On the segment $\bigcup_{h_{j}<h<h_{j+1}} \U_{h} \times [0,1]$ we deformation retract to the slice $\U_{h_{j}} \times \{1\}$, then apply $\psi_{j-1}$. Again, this is well-defined and compatible with gluing. The maps $\zeta^{j}_{\mathcal{X} \mathcal{X}'}$ are also compatible with each other on boundary components. The map extends via the constant map to $\V_{0}$ on everything below the slice $\U_{h_{2}} \times \{\frac{1}{4}\}$. We therefore obtain a globally defined map $\zeta_{\mathcal{X} \mathcal{X}'}$, which is continuous and proper. \\

Now we compute the composition $\zeta_{\mathcal{X} \mathcal{X}'}\zeta_{\mathcal{X}'\mathcal{X}}: M_{\mathcal{X}'}\rightarrow M_{\mathcal{X}'}$ (see Figure \ref{fig:prophom3}). This is given, piecewise, by the composition $\xi^j_{\mathcal{X}, \mathcal{X}'}\xi^j_{\mathcal{X}',\mathcal{X}}$:

\begin{align*}
\xi^j_{\mathcal{X}, \mathcal{X}'}\xi^j_{\mathcal{X}',\mathcal{X}}: \V_{j+1} \times [0,\frac{1}{2}]\cup_{\tilde{\phi}} \V_{j} \times [\frac{1}{2},1] \longrightarrow \V_{j} \times [0,\frac{1}{2}]\cup_{\tilde{\phi}} \V_{j-1} \times [\frac{1}{2},1] 
\end{align*}

\begin{align*}
\xi^j_{\mathcal{X}, \mathcal{X}'}\xi^j_{\mathcal{X}',\mathcal{X}}: & \V_{j+1} \times \{\frac{1}{2}\} \rightarrow \V_{j} \times \{\frac{1}{2}\}\quad   &= \psi_{j}\varphi_{j+1}\\
 &\V_{j+1} \times \{\frac{7}{16}\} \rightarrow \V_{j-1} \times \{1\}  &= \psi_{j-1} \varphi_{j}\psi_{j}\varphi_{j+1}\\
 &\V_{j+1} \times \{\frac{3}{8}\} \rightarrow \V_{j-1} \times \{\frac{1}{2}\}   &=\psi_{j-1}\phi_{h_{j+1}, h_{j}} \varphi_{j+1}\\
 &\V_{j+1} \times \{\frac{1}{4}\} \rightarrow \V_{j-1} \times \{\frac{1}{2}\}&= \psi_{j-1}\phi_{h_{j+1}, h_{j}} \varphi_{j+1}\\
 &\V_{j+1} \times \{0\} \rightarrow \V_{j-1} \times \{\frac{1}{2}\} &= \psi_{j-1} \varphi_{j} \psi_{j} \varphi_{j+1}\\
 &\V_{j} \times [\frac{1}{2},1] \rightarrow \V_{j-1} \times \{\frac{1}{2}\} &= \psi_{j-1}\varphi_j
\end{align*}

\begin{figure}
\centering
  \centering
  \includegraphics[width=0.8\linewidth]{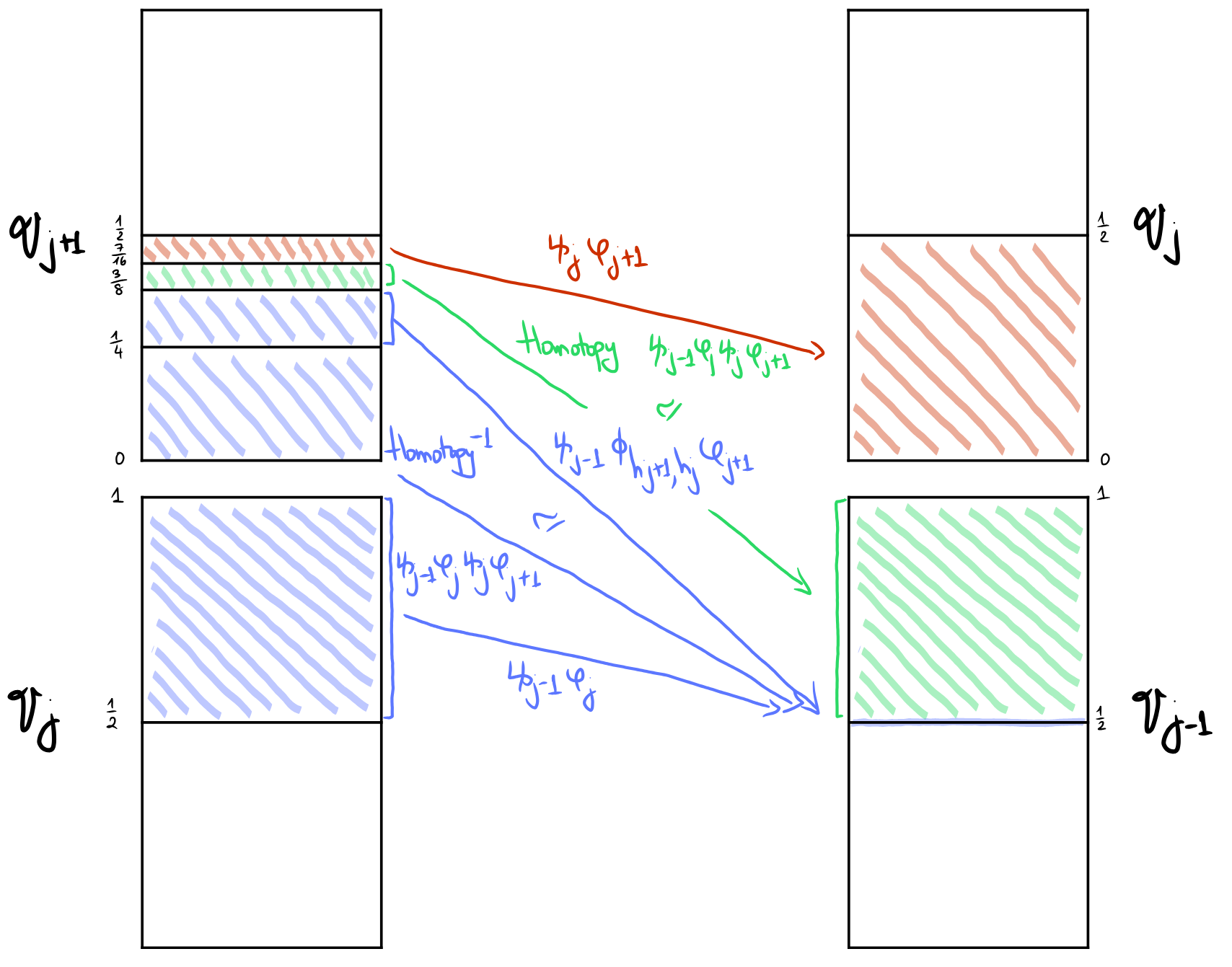}
  \caption{The composition $\zeta^j_{\mathcal{X} \mathcal{X}'}\zeta^j_{\mathcal{X}'\mathcal{X}}$ for $j\geq 2$. }
  \label{fig:prophom3}
\end{figure}

Looking at this carefully, we observe that the map between the slices $\V_{j+1} \times \{\frac{7}{16}\}$ and $\V_{j+1} \times \{0\}$ consist precisely of a homotopy between $\psi_{j-1} \varphi_{j}\psi_{j}\varphi_{j+1}$ and $\psi_{j-1}\phi_{h_{j+1}, h_{j}} \varphi_{j+1}$, a constant segment, then the inverse homotopy. Therefore, this is homotopic to a map which is defined to be $\psi_{j-1} \varphi_{j}\psi_{j}\varphi_{j+1}$ on the entire segment. Remembering that $\psi_{j-1}\varphi_{j} = \tilde{\phi}_{j,j-1}$ we obtain a map which is just a scaled version of the deformation retract of $M_{\mathcal{X}'}\rightarrow M_{\mathcal{X}'}$ obtained by going down a level via the bonding map. This is properly homotopic to the identity map. Therefore, $\zeta_{\mathcal{X} \mathcal{X}'}\zeta_{\mathcal{X}'\mathcal{X}}\simeq \id_{M_{\mathcal{X}'}}$. \\

Now we compute the composition $\zeta_{\mathcal{X}'\mathcal{X}}\zeta_{\mathcal{X} \mathcal{X}'}: M_{X}\rightarrow M_{X}$  (see Figure \ref{fig:prophom4}). This is given, piecewise, by the composition $\xi^{j-1}_{\mathcal{X}',\mathcal{X}} \xi^j_{\mathcal{X}, \mathcal{X}'}$:

\begin{align*}
\xi^{j-1}_{\mathcal{X}',\mathcal{X}} \xi^j_{\mathcal{X}, \mathcal{X}'}:  \U_{h_{j+1}} \times [0,\frac{1}{2}] \cup_{\phi} (\bigcup_{h_{j}<h<h_{j+1}} \U_{h} \times [0,1]) \cup_{\phi} \U_{h_j} \times [\frac{1}{2},1] \longrightarrow \\
 \U_{h_{j}} \times [0,\frac{1}{2}] \cup_{\phi} (\bigcup_{h_{j-1}<h<h_{j}} \U_{h} \times [0,1]) \cup_{\phi} \U_{h_{j-1}} \times [\frac{1}{2},1] 
\end{align*}

\begin{align*}
\xi^{j-1}_{\mathcal{X}',\mathcal{X}} \xi^j_{\mathcal{X}, \mathcal{X}'}:& \U_{h_{j+1}} \times \{\frac{1}{2}\} \rightarrow \U_{h_{j}} \times \{\frac{1}{2}\}\quad   &= \varphi_{j}\psi_{j}\\
&\U_{h_{j+1}} \times \{\frac{7}{16}\} \rightarrow \U_{h_{j}-1} \times \{1\}  &= \phi_{h_{j}, h_{j}-1} \varphi_{j} \psi_{j}\\
&\U_{h_{j+1}} \times \{\frac{3}{8}\} \rightarrow \U_{h_{j-1}} \times \{1\}   &=\phi_{h_{j}, h_{j-1}} \varphi_{j} \psi_{j}\\
&\U_{h_{j+1}} \times \{\frac{1}{4}\} \rightarrow \U_{h_{j-1}} \times \{\frac{1}{2}\}&= \varphi_{j-1}\psi_{j-1}\varphi_{j}\psi_{j}\\
&\U_{h_{j+1}} \times \{0\} \rightarrow \U_{h_{j-1}} \times \{\frac{1}{2}\} &= \varphi_{j-1}\psi_{j-1}\phi_{h_{j+1},h_{j}}\\
&\U_{h_{j}} \times [\frac{1}{2},1] \rightarrow \U_{h_{j-1}} \times \{\frac{1}{2}\} &= \varphi_{j-1} \psi_{j-1}\\
\end{align*}

We want to homotope this map to the just a scaled version of the deformation retract of $M_{X}\rightarrow M_{X}$ obtained by going down a level via the bonding map. The obstacle is that $\varphi_{j} \psi_{j}$ is not the bonding map $\phi_{h_{j+1}, h_{j}}$, but is homotopic to it by a linear homotopy. We claim that we can glue together a homotopy $M_{X}\times [0,1] \rightarrow M_{X}$ piecewise, which results in what we want. \\

It is clear what to do on the segments outside of $\U_{h_{j+1}} \times [0, \frac{3}{8}]$: we take the linear homotopy between $\varphi_{j} \psi_{j}$ and $\phi_{h_{j+1}, h_{j}}$, and post-compose with the relevant bonding map. On the segment $\U_{h_{j+1}} \times [0, \frac{3}{8}]$ we wish to homotope what we currently have to a map which is $\phi_{h_{j}, h_{j-1}} \phi_{h_{j+1}, h_{j}}$ on the entire segment. \\

\begin{figure}[H]
\centering
  \centering
  \includegraphics[width=0.8\linewidth]{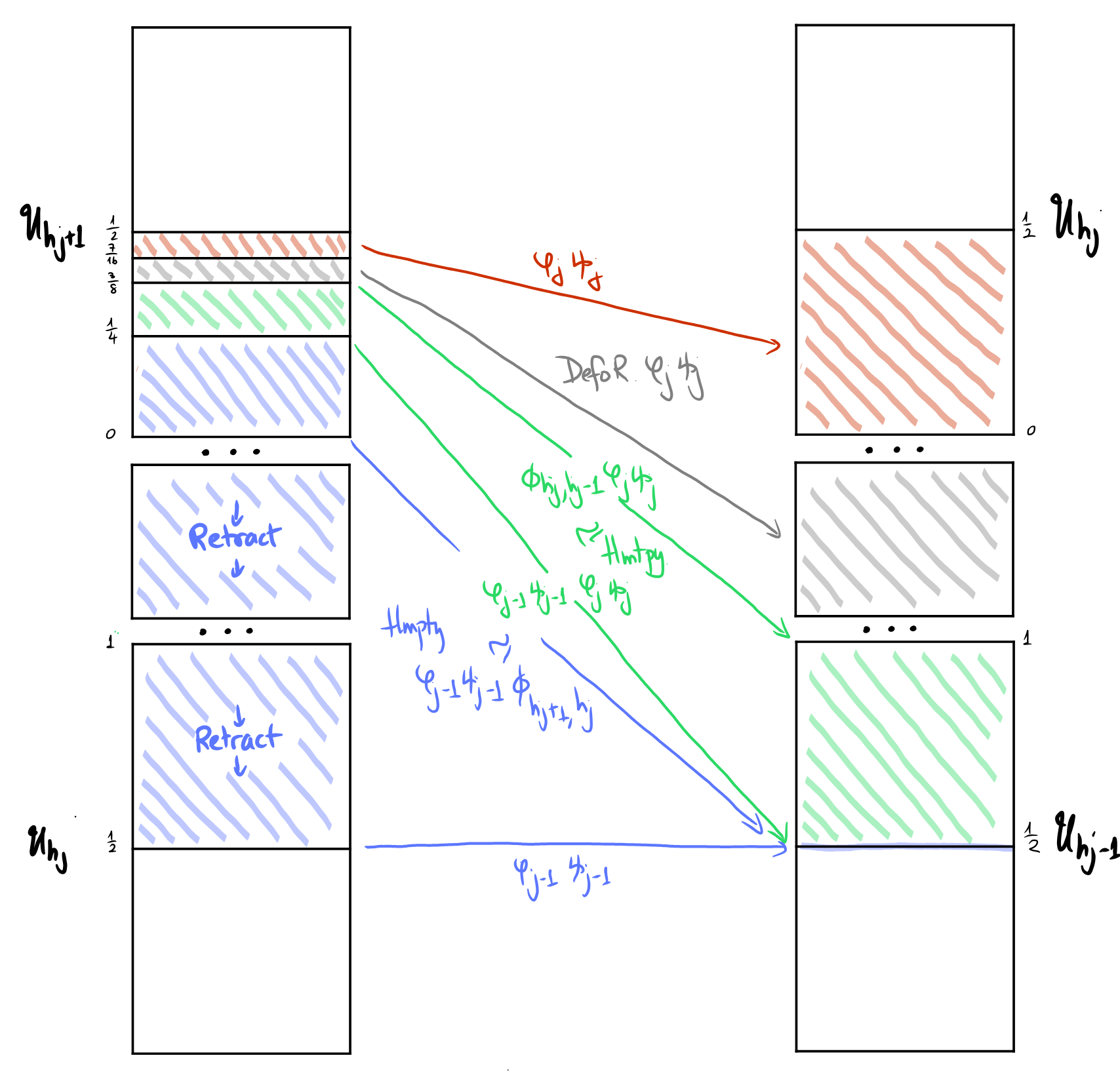}
  \caption{The composition $\xi^{j-1}_{\mathcal{X}',\mathcal{X}} \xi^j_{\mathcal{X}, \mathcal{X}'}$ for $j\geq 2$.}
  \label{fig:prophom4}
\end{figure}

Let us just consider the interval $[\frac{1}{4}, \frac{3}{8}]$ as an example of what to do. Let $f_t: \U_{h_{j+1}} \times [0,1] \rightarrow \U_{h_{j-1}}$ be the linear homotopy between $f_0=\phi_{h_j, h_{j-1}}\varphi_j\psi_j$ and $f_1=\phi_{h_j, h_{j-1}} \phi_{h_{j+1}, h_{j}}$. Let $g_t: \U_{h_{j+1}} \times [0,1]\rightarrow \U_{h_{j-1}}$ be the linear homotopy between $g_0=\varphi_{j-1}\psi_{j-1}\varphi_{j}\psi_j$ and $g_1=\phi_{h_j, h_{j-1}} \phi_{h_{j+1}, h_{j}}$. We can define a homotopy 
\begin{align*}
H:\U_{h_{j+1}} \times [0,1] \times [0,1] &\longrightarrow \U_{h_{j-1}}\\
H(s,t) &:= (1-s) f_t + sg_t
\end{align*}
in barycentric coordinates, where $s$ is the re-parametrised height variable of the interval $[\frac{1}{4},\frac{3}{8}]$. Observe that this is well-defined since the following condition is met: for any simplex $\sigma\subset \U_{h_{j+1}}$:
\begin{align*}
\bigcap_{v\in vert \sigma}\bigcap_{i=0,1} \langle st(f_i(v))\rangle  \cap \langle st(g_i(v))\rangle \neq \emptyset
\end{align*}
This comes from the fact that all the maps involved are simplicial refinement maps. So everything lies within a common simplex and we can form linear combinations in barycentric co-ordinates. Furthermore, $H$ restricted to $s=0,1$ are $f_t$ and $g_t$ respectively, and $H$ restricted to $t=0,1$ are the original map we had on $\U_{h_{j+1}} \times [\frac{1}{4}, \frac{3}{8}]$ and the map $\phi_{h_j, h_{j-1}} \phi_{h_{j+1}, h_{j}}$ on the entire segment respectively. \\

We do the same thing on the segment $\U_{h_{j+1}} \times [0, \frac{1}{4}]$ and observe that since all the homotopies are linear, they glue together to a well-defined homotopy on the entirety of $M_X \times [0,1]$. Therefore, the composition $\zeta_{\mathcal{X}'\mathcal{X}}\zeta_{\mathcal{X} \mathcal{X}'}: M_{X}\rightarrow M_{X}$ is homotopic to a deformation retract obtained by going down a level. As before, this is homotopic to the identity. \\

The computations above work for $j\geq 2$. For the segment $\mathcal{V}_{2} \times [0,\frac{1}{2}]\cup \mathcal{V}_{1} \times [0,1]$, the composition $\zeta_{\mathcal{X} \mathcal{X}'}\zeta_{\mathcal{X}'\mathcal{X}}$ on slices is
\begin{align*}
\zeta_{\mathcal{X} \mathcal{X}'}\zeta_{\mathcal{X}'\mathcal{X}} &: \V_{2} \times \{\frac{1}{2}\} \rightarrow \V_{1} \times \{\frac{1}{2}\}\quad   &= \psi_{1}\varphi_{2} =\tilde{\phi}_{2,1}\\
&\V_{2} \times \{\frac{7}{16}\} \rightarrow \V_{0}  &= \text{constant map to point}
\end{align*}
Everything below $\V_{2} \times \{\frac{7}{16}\}$ also gets sent to the point $\V_{0}$. There is a homotopy from $\zeta_{\mathcal{X} \mathcal{X}'}\zeta_{\mathcal{X}'\mathcal{X}}$ to the identity which glues compatibly with the homotopy on the rest of the telescope. For the segments below the slice $\mathcal{U}_{h_{2}} \times \{\frac{1}{2}\}$ we get
\begin{align*}
\zeta_{\mathcal{X}',\mathcal{X}} \zeta_{\mathcal{X}, \mathcal{X}'}: \U_{h_{2}} \times \{\frac{1}{2}\} &\rightarrow \U_{1} \times \{\frac{1}{2}\}\quad   &= \varphi_{1}\psi_{1}\\
\U_{h_{2}} \times \{\frac{1}{4}\} &\rightarrow \U_{0}  &=  \text{constant map to point}
\end{align*}
Everything below $\U_{h_{2}} \times \{\frac{1}{4}\}$ is also sent to $\U_{0}$.
The cone of the linear homotopy between $\varphi_{1}\psi_{1}$ and $\phi_{h_2,1}$ glues compatibly with the homotopy constructed on the rest of $M_{X}$, as well as the homotopy from the identity to the deformation retract obtained by going downwards.\\

In conclusion, we have constructed a proper homotopy equivalence between mapping telescopes $\zeta_{\mathcal{X}',\mathcal{X}}: M_{\mathcal{X}'}\rightarrow M_{X}$ with proper homotopy inverse $\zeta_{\mathcal{X},\mathcal{X}'}$. By construction, all maps and homotopies are of pairs. 
\end{proof}

\begin{remark} Observe that the restriction to the base ray ${\zeta_{\mathcal{X}',\mathcal{X}}}_{|\tau}: \tau \rightarrow \omega$ is already cellular. By the proper cellular approximation theorem (Theorem 10.1.14 in \cite{geoghegan2007topological}), $\zeta_{\mathcal{X}',\mathcal{X}}$ is properly homotopic, rel $\tau$, to a proper cellular map. Similarly, $\zeta_{\mathcal{X},\mathcal{X}'}$ is properly homotopic, rel $\omega$ to a proper, cellular map. \\

The analogous statement of Lemma \ref{finitedim} is true for the unpointed case, by forgetting all base points and base rays.  
\end{remark}

\begin{cor} \label{wseimplieswphe2} (The shape Whitehead theorem) A shape morphism $F:S(X,x_{0})\rightarrow S(Y, y_0)$ of finite shape-dimensional, connected, compact metric spaces is a shape equivalence if and only if pro-$\pi_{n}(F):$pro-$\pi_{n}(X,x_0)\rightarrow$ pro-$\pi_{n}(Y,y_0)$ is an isomorphism for all $n\geq 0$. 
\end{cor}

\begin{proof} By Proposition \ref{finitedim} we can assume that $\mathcal{X}=\{(X_{h},x_{h}),\phi_{h+1,h};\mathbb{N}\}$ and $\mathcal{Y}=\{(Y_{h},y_{h}),\varphi_{h+1,h};\mathbb{N}\}$ are finite-dimensional, pointed shape expansions of $(X,x_0)$ and $(Y,y_0)$ respectively, where all $X_{h},Y_{h}$ are simplicial complexes and all bonding maps are pointed and simplicial. We can assume that $F= (\rho,\{q_{h}\}): \mathcal{X}\rightarrow \mathcal{Y}$. Let $(M_{\mathcal{X}},\omega)$ and $(M_{\mathcal{Y}},\tau)$ be the inverse mapping telescopes of $\mathcal{X}$ and $\mathcal{Y}$, with base rays $\omega,\tau$ associated to $x_0$ and $y_0$ respectively. By Theorem \ref{wseimplieswphe1} there is a continuous, proper, cellular map of pairs $\psi: (M_{\mathcal{X}},\omega) \rightarrow (M_{\mathcal{Y}},\tau)$ with $\mathcal{P}(\psi) = (\rho,\{q_{h}\})$. By the assumptions of the corollary, $\psi$ is a proper homotopy equivalence of pairs, with inverse $\xi$. By Lemma \ref{phemeansshape} we have that $\mathcal{P}(\xi)\circ F = \id_{\mathcal{X}}$ and $F\circ \mathcal{P}(\mathcal{\xi)} = \id_{\mathcal{Y}}$ in pro-hCW$_{\ast}$. Therefore $F$ is a pointed shape equivalence. \\

The other direction comes from the fact that pro-$\pi_{n}$ is a functor.  
\end{proof}

\begin{remark}
The proof of Corollary \ref{wseimplieswphe2} requires that $X,Y$ are compact metric spaces, because we have used the proper homotopy theory of a specific geometric object (the inverse mapping telescope). The general shape Whitehead theorem (Theorem \ref{shapewhitehead}), applicable to arbitrary topological spaces, was proven by Morita (\cite{morita1974hurewicz}) using similar ideas: construct an inverse piece by piece by using the mapping cylinders of each $q_{\beta}: X_{\rho(\beta)}\rightarrow Y_{\beta}$ associated to the shape morphism $F=(\rho,\{q_{\beta}\}): \{X_{\alpha}, \phi_{\alpha',\alpha}; \mathcal{A}\} \rightarrow \{Y_{\beta}, \varphi_{\beta',\beta}; \mathcal{B}\}$. 
\end{remark}

The assumption of finite-dimensionality of $X,Y$ cannot be omitted. 

\begin{example} (The Kahn compact space) Consider the space $X_0$ defined as the CW complex obtained by attaching a $(2p+1)$-cell to the sphere $S^{2p}$ by a map of prime degree $p$. For each $n\geq 0$, define $X_{n+1}$ inductively as the $(2p-2)$-fold suspension of $X_{n}$:
\begin{align*}
X_{n+1} = \Sigma^{2p-2}(X_{n})
\end{align*}
The bonding maps $f_{n}: X_{n}\rightarrow X_{n-1}$ for $n>1$ are defined inductively 
\begin{align*}
f_{n} = \Sigma^{2p-2}(f_{n-1})
\end{align*}
beginning with the map $f_{1}: X_{1}\rightarrow X_{0}$ defined in a special way. The Kahn compact space $K$ is the limit of the inverse sequence 
\begin{align*}
X_{0}\xleftarrow{f_{1}} X_{1} \xleftarrow{f_{2}} X_{2} \leftarrow \dots
\end{align*}
This is a compact metric space of infinite shape dimension. Since $\pi_{k}(\Sigma^{n(2p-2)},\ast)=0$ for $n(2p-2)\geq k$, it follows that pro-$\pi_{k}(K,\ast)=0$ for all $k$. However, $K$ does not have the shape of a point because the composition 
\begin{align*}
f_{1}\circ f_{2}\circ \dots \circ f_{n}: X_{n}\rightarrow X_{0}
\end{align*}
is not homotopic to a constant for any $n$. See \cite{draper1976example} for more details. 
\end{example}

\subsection{Strong shape theory} \label{strongshapetheory}

\begin{definition} ($Z$-sets)
Let $Y$ be a locally compact metrisable space which is path-connected. 
\begin{itemize}
\item A \textit{compactification} of $Y$ is a compact metrisable space $W$ which contains $Y$ as a dense open subset. The nowhere dense compact set $X:=W\setminus Y$ is called the \textit{compactifying space}.
\item A closed subset $D\subset W$ is called a \textit{$Z$-set} in $W$ if for every open set $U\subset W$, the inclusion map $\iota: U\setminus D \rightarrow U$ is a homotopy equivalence. If the compactifying space $X$ is a $Z$-set in $W$, we say that $W$ is a \textit{$Z$-set compactification} of $Y$. 
\end{itemize}
\end{definition}

Let $X$ be a compact metric space with shape expansion $\mathcal{X}=\{X_{h},\phi_{h+1,h};\mathbb{N}\}$, where $X_{h}$ are finite simplicial complexes and $\phi_{h+1,h}$ are simplicial maps. The inverse mapping telescope $M_{\mathcal{X}}$ is a strongly locally finite CW complex, and $X \cup M_{\mathcal{X}}$ (with the appropriate topology) is a $Z$-set compactification of $M_{\mathcal{X}}$ with compactifying space $X$. See Example 17.5.6 in \cite{geoghegan2007topological} for details. 

\begin{definition} (The category $\mathbf{StrongShape}$)
Let $\mathbf{StrongShape}$ be the category whose objects are pairs $(W,X)$ where $W$ is a compact $AR$ and $X$ is a (closed) $Z$-set in $W$; a morphism $(W_1,X_1)\rightarrow (W_2,X_2)$ in $\mathbf{StrongShape}$ is a proper homotopy class of proper maps $Y_1\rightarrow Y_2$ where $Y_{i} = W_{i}\setminus X_{i}$. 
\end{definition}

There is a functor from $\mathbf{StrongShape}$ to Shape, defined as follows. Let $(W,X)$ be an object of $\mathbf{StrongShape}$. Let $\{U_{\alpha},\iota_{\alpha',\alpha};\mathcal{A}\}$ be the directed set of all compact ANR neighbourhoods of $X$ in $W$ with inclusion morphisms $\iota_{\alpha',\alpha}: U_{\alpha'}\rightarrow U_{\alpha}$ for $U_{\alpha'} \subset U_{\alpha}$. This is a shape expansion of $X$. In fact, there exists a cofinal subsystem $\{U_{h},\iota_{h+1,h};\mathbb{N}\}$ with the property that $X=\cap_{h\in \mathbb{N}} U_{h}$.  \\

A morphism $(W_1,X_1)\rightarrow (W_{1},X_{2})$ is represented by a continuous, proper map $\psi: W\setminus X_{1}\rightarrow W_{2}\setminus X_{2}$. Let $U(W_1,X_{1})=\{U_{h},\iota^1_{h+1,h};\mathbb{N}\}$ and $U(W_{2},X_{2})=\{V_{h},\iota^2_{h+1,h};\mathbb{N}\}$ be cofinal ANR-systems associated with $X_1,X_{2}$ respectively. Since $X_{1},X_{2}$ are $Z$-sets, we can replace $U_{h}$ with $U_{h}\setminus X_{1}$ (resp. $V_{h}$ with $V_{h}\setminus X_{2}$). For every $i\in \mathbb{N}$ there exists a $\rho(i)$ such that $\psi (U_{\rho(i)}\setminus X_{1}) \subset V_{i}\setminus X_{2}$. Let $q_{i}:= \psi_{|U_{\rho(i)}\setminus X_{1}}$ The pair $\mathcal{P}(\psi):=(\rho,\{q_{h}\}): U(W_{1},X_{1})\rightarrow U(W_{2},X_{2})$ is a morphism in pro-hCW. 

\begin{definition} (The functor $\mathcal{P}$). The functor $\mathcal{P}: \mathbf{StrongShape} \rightarrow$ Shape is defined on objects as $\mathcal{P}(W,X) := (U(W,X),X)$ and on morphisms $\psi: W\setminus X_{1}\rightarrow W_{2}\setminus X_{2}$ as $\mathcal{P}(\psi) := (\rho,\{q_{h}\})$.  
\end{definition}

\begin{remark} If $W_1= X_1\cup M_{\mathcal{X}_1}$ and $W_{2} = X_{2} \cup M_{\mathcal{X}_{2}}$ then 
\begin{align*}
U(W_1,X_{1})&=\{L^c_{h},\iota^1_{h+1,h};\mathbb{N}\}\\
 U(W_2,X_{2})&=\{Q^c_{h},\iota^2_{h+1,h};\mathbb{N}\}
\end{align*}
where $L_{h}:=(M_{\mathcal{X}_1})_{[1,h]}$ and $Q_{h}:=(M_{\mathcal{X}_2})_{[1,h]}$. In this case, the functor $\mathcal{P}$ agrees (after suitable identifications) with the $\mathcal{P}$ from Lemma \ref{functorP}.
\end{remark}

The functor $\mathcal{P}$ in general is full but not faithful (see Example 17.7.2 in \cite{geoghegan2007topological}). 

\begin{definition} We say that two compact metric spaces $X,Y$ are \textit{strong shape equivalent} $X\simeq_{sse} Y$ if there exists a proper homotopy equivalence $\psi: M_{\mathcal{X}}\rightarrow M_{\mathcal{Y}}$ for shape expansions $\mathcal{X},\mathcal{Y}$ of $X$ and $Y$.
\end{definition}

We have the following theorem:

\begin{theorem} \label{SSandS} (Theorem $5.2.9$ in \cite{edwards2006cech}) Let $\mathcal{X}$ and $\mathcal{Y}$ be inverse ANR-systems. Let  $\psi: M_{\mathcal{X}}\rightarrow M_{\mathcal{Y}}$ be a continuous proper map such that the induced morphism $\mathcal{P}(\psi): \mathcal{X}\rightarrow \mathcal{Y}$ is invertible in pro-hCW. There is a proper homotopy equivalence $\psi': M_{\mathcal{X}}\rightarrow M_{\mathcal{Y}}$ such that $\mathcal{P}(\psi')=\mathcal{P}(\psi)$ in pro-hCW.
\end{theorem}

By Theorem \ref{wseimplieswphe1} (which also applies to the unpointed case) a shape morphism $(\rho, \{q_{h}\}): \mathcal{X}\rightarrow \mathcal{Y}$ induces a proper, cellular map $\psi: M_{\mathcal{X}}\rightarrow M_{\mathcal{Y}}$ such that $\mathcal{P}(\psi) = (\rho, \{q_{h}\})$. Therefore, we conclude that two compact metric spaces are shape equivalent $X\simeq_{se}Y$ if and only if they are strong shape equivalent $X\simeq_{sse}Y$. 

\begin{cor} Let $X$ be a compact metric space. Any two inverse mapping telescopes of different shape expansions of $X$ are properly homotopy equivalent. 
\end{cor}

The properties of the functor $\mathcal{P}$, as well as the exact relationship between Shape, Shape$_{\ast}$, and $\mathbf{StrongShape}$ are complicated have been omitted for the sake of brevity. 
\newpage
\section{Coarse Whitehead theorem}

We begin with the following setup: let $(X,x_0)$ and $(Y,y_0)$ be connected, compact metric spaces with finite shape dimension. Let $\omega$ be the standard parametrisation of the base ray $c\{x_0\}$ and $\mathfrak{a}:(cX,\omega)\rightarrow (cY,\mathfrak{a}\omega)$ a coarse map. Let $M_X$ and $M_Y$ be the inverse mapping telescopes for the spaces $X,Y$. The goal of this chapter is to show that $\mathfrak{a}$ is a coarse homotopy equivalence under the assumption that $\mathfrak{a}$ induces isomorphisms on all coarse homotopy groups (Theorem \ref{unbasedwhitehead}).

\subsection{Connections to coarse homotopy groups}

In the following two technical propositions, we prove the isomorphism between coarse homotopy groups $\pi^c_n(cX,\tau)$ and end homotopy groups $\pi^e_{n}(M,\tau)$ for arbitrary base rays $\tau$ (Proposition \ref{isogeneralray}) and construct from the coarse map $\mathfrak{a}: cX\rightarrow cY$ a proper, locally Lipschitz map $\psi: M_{X}\rightarrow M_{Y}$, which induces morphisms on the end homotopy groups of $M_{X}$ and $M_{Y}$ (Proposition \ref{psiexistence}).

\begin{prop} \label{isogeneralray}  Let $\tau$ be a representative of $[\tau]\in \pi_{0}^{c}(cX,\omega)$. There exists a $l\in \mathbb{N}$, sufficiently large, such that the composition 
\begin{align*}
\theta \circ s_{l}: \pi_{n}^c(cX,\tau)\xrightarrow{s_{l}}  \pi_{n}^{c}(cX, s_{l} \tau) \xrightarrow{\theta} \pi_{n}^{L,e}(M_{X},\widetilde{is_{l}\tau})  
\end{align*}
is an isomorphism. 
\end{prop}

Before we begin the proof, let us introduce some auxillary definitions. 

\begin{definition} Let $U$ be a set in $cX$, $S>0$. We define the \textit{$S$-height neighbourhood of $U$} as 
\begin{align*}
N_{S}(U) := \{(hx,h)\in cX\,|\, x\in p'(U), h\in N_{S}(p(U))\}
\end{align*}
where $p:cX\rightarrow [1,\infty)$ is the height projection, $p':cX\rightarrow X$ is the $X$-projection, and $N_{S}(p(U))$ denotes the closed $S$-neighbourhood of the set $p(U)\subset [1,\infty)$. 
\end{definition}

\begin{lemma} \label{shiftweirdness} Let $U$ be a set of diameter $\delta<\frac{1}{4}$ and $u = (tx,t)\in U$. Then for $\varepsilon<\frac{1}{4}$ we have that $iN_{\varepsilon}(U)\subset \langle gst(v)\rangle$ for $v = (z_{x},[t])$, where $[t]$ denotes the closest integer to $t$ and $z_{x}$ a point in the appropriate discretisation $Z_{h}$ of $X$. 
\end{lemma}

\begin{proof} This is a simple application of Lemma \ref{weirdness}. There are two cases to consider since $pN_{\varepsilon}(U)$ has distance $<\frac{1}{2}$ from $t$. Suppose that $pN_{\varepsilon}(U)\subset [2^h,2^{h+1}]$ for some $h\in \mathbb{N}_{0}$. The set $p'N_{\varepsilon}(U)=p'(U)$ has diameter less than $\frac{1}{2^{h+1}}$, so $p'N_{\varepsilon}(U) \subset B(z_{x},\frac{1}{2^h})$ for a point $z_{x}\in Z_{h}$ with $d(x,z_{x})\leq\frac{1}{2^{h+1}}$. Therefore, $iN_{\varepsilon}(U)\subset \langle gst(v)\rangle$ for $v = (q(z_{x}),[t])$, where $q=\phi_{h}$ if $[t]=2^h$ and $q=\id$ otherwise. \\

Alternatively, if the interval containing $pN_{\varepsilon}(U)\subset (2^{h}-1,2^{h}+1)$ intersects a gluing slice then the set $p'(U)$ still has diameter less than $\frac{1}{2^{h+1}}$. In this case, the closest point to $t$ will always be $2^h$, and $iN_{\varepsilon}(U)\subset \langle gst(v)\rangle$ for $v = (\phi_{h}(z_{x}),2^h)$, where $z_{x}\in Z_{h}$ is the closest point to $x$. 
\end{proof}

We now address some technical issues required for the proof of Proposition \ref{isogeneralray}. First, since $\tau$ is an arbitrary coarse base ray, we have to compose with a shrinking map so that the image under $i$ can be simplicially approximated. Additionally, we must take care to ensure that the change of base ray homomorphisms are actually coarse-Lipschitz. \\

For the rest of this subsection, two constants $\delta<\frac{1}{32}$ and $\varepsilon<\frac{1}{8}$. The map $p\circ \tau: [1,\infty)\rightarrow [1,\infty)$ is large-scale Lipschitz, so there exists $L>0$ such that $p\tau(h)< Lh+L$. Now, fix a number $l\in \mathbb{N}$ such that the following conditions are met:

\begin{enumerate}
\item The set $s_{l}\tau [n-1,n+1]$ has diameter less than $\delta$ for all vertices $n\in \mathbb{N}$.  
\item $\frac{Lh+L}{2^l}< h$ for all $h\in [1,\infty)$. 
\end{enumerate}

Let $n\in \mathbb{N}$. We choose $t_n$ as the closest integer to $p s_{l}\tau(n)$, $v_{n}$ the closest vertex to $p's_{l}\tau(n)$ as specified in Lemma \ref{shiftweirdness} and let $z_{n} = (v_{n},t_{n})$. We obtain a geometric simplicial approximation of $is_{l}\tau$, defined by $\widetilde{is_{l}\tau}(n) = z_{n}$, which we fix. For ease of notation, let $\beta = s_{l}\tau$. For any $r\in \mathbb{N}_{0}$ we can fix a geometric simplicial approximation of $is_{r} \beta$ by letting $\widetilde{is_{r}\beta}(n) = q_{[\frac{p\beta(n)}{2^r}]} z_{n}$, the projection of $z_{n}$ to the height which is the closest integer to $\frac{p\beta(n)}{2^r}$. Let
\begin{align*}
H: c[0,1] &\rightarrow cX\\
(h,ht)&\mapsto (\{tp\beta(h) + (1-t)\frac{p\beta(h)}{2^r}\} p' \beta (h), tp\beta(h) + (1-t)\frac{p\beta(h)}{2^r}) := (\gamma_{h,t} p'\beta(h), \gamma_{h,t})
\end{align*}
denote the linear homotopy from  $s_{r} \beta$ to $\beta$. We now prove some technical lemmas about the simplicial approximation of $H$. \\

Divide $[0,1]$ into edges of length $\frac{\varepsilon}{4}$. This simplicial structure induces a simplicial structure on $c[0,1]$, which we fix. Let $(n,nt)$ be a vertex in $c[0,1]$. By Lemma \ref{technicalcone}, we have that $|st(n,nt)| \subset |st(n)| \times I$, where $I \subset [nt-\frac{\varepsilon}{2}, nt + \frac{\varepsilon}{2}]$ is an interval of length at most $\varepsilon$ centred at $nt$. Since $p' H_{t}(|st(n)|) = p'\beta |st(n)|$ for all $t$, we have 
\begin{align*}
iH|st(n,nt)|\subset iH (|st(n)| \times I) \subset iN_{\frac{\varepsilon}{2}}( H_{t} |st(n)|) \subset \langle gst (\mathfrak{q}_{[\gamma_{n,t}]}(z_{n}))\rangle 
\end{align*}
where $\mathfrak{q}_{[\gamma_{n,t}]}$ denotes the deformation retract on $M_{X}$ to the height $[\gamma_{n,t}]$, the closest integer to $\gamma_{n,t}$. The inclusion $iH (|st(n)| \times I) \subset iN_{\frac{\varepsilon}{2}}( H_{t} |st(n)|)$ is guaranteed by condition $2.$ of our choice of $l$: it means that the homotopy is $1$-Lipschitz in the time variable. By letting $\widetilde{iH}(n,nt) = \mathfrak{q}_{[\gamma_{n,t}]}(z_{n})$, we obtain a simplicial approximation of $iH$, which agrees with the fixed simplicial approximation of  $\widetilde{is_{r}\beta}$ and  $\widetilde{i\beta}$  for $t=0,1$. We fix this simplicial approximation and denote it by $\widetilde{iH}$, or $\widetilde{iH^{r}}$. \\

Consider now the $k$-th barycentric subdivision of $c[0,1]$, which we denote by $bs^{k}(c[0,1])$. We denote by $|st^{k}(m,mt)|$ the closed star of a vertex $(m,mt)$ in $bs^{k}(c[0,1])$.

\begin{lemma} \label{twohomotopies}  Let 
\begin{align*}
Q: c[0,1]&\rightarrow M_{X}\\
(h,ht) &\mapsto (p'  \widetilde{i\beta} (h), tp \widetilde{i\beta}(h) + (1-t)\frac{p \widetilde{i\beta}(h)}{2^r}) := (p'\widetilde{i\beta}(h),\overline{\gamma_{h,t}})
\end{align*}
be the linear homotopy from $D_{r} \widetilde{i\beta}$ to $\widetilde{i\beta}$, where $D_{r}: M_{X}\rightarrow M_{X}$ is the map which divides heights by $\frac{1}{2^r}$ and $p'  \widetilde{i\beta} (h)$ is identified with its appropriate image under the bonding maps $\phi$. There is a $k$ sufficiently large and a proper Lipschitz homotopy $c[0,1]^2\rightarrow M_{X}$ between $Q$ and a simplicial approximation of $iH$ with respect to $bs^{k}(c[0,1])$.  
\end{lemma}

To prove this, we show that there is a simplicial map $(\widetilde{iH})_{k}: bs^{k}(c[0,1])\rightarrow M_{X}$ which is a simplicial approximation of $iH$, and is proper Lipschitz homotopic to $Q$. 

\begin{proof} (of Lemma \ref{twohomotopies}) Since the simplicial structure on $c[0,1]$ is uniformly bounded, there exists a $k$ large enough such that for any vertex $(m,mt)$ in $bs^{k}(c[0,1])$, $|st^{k}(m,mt)|$ is contained within $[m-\frac{1}{4}, m+\frac{1}{4}] \times [mt-\frac{\varepsilon}{2}, mt+\frac{\varepsilon}{2}]$. Up to Lipschitz rescaling, we can assume that the image of $p\widetilde{i\beta}[m-\frac{1}{4}, m+\frac{1}{4}]$ has length at most $\frac{2}{3}$. (see Figure \ref{fig:LipschitzApproxBeta})

\begin{figure}[H]
\centering
  \centering
  \includegraphics[width=0.5\linewidth]{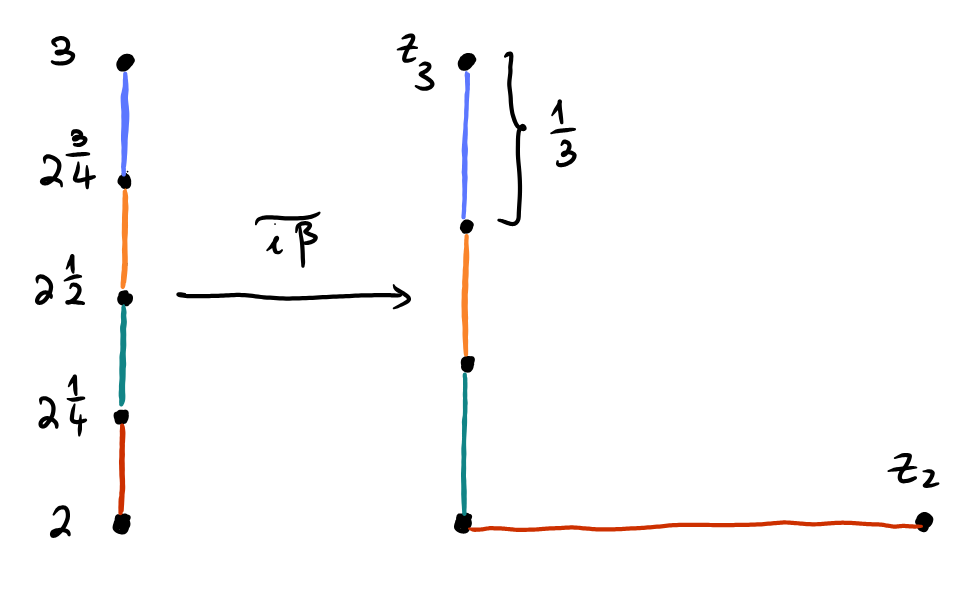}
  \caption{Illustration of the map $\widetilde{i\beta}[2,3]$ if $z_{2}$ and $z_{3}$ lie on different heights.}
  \label{fig:LipschitzApproxBeta}
\end{figure}

Let $(m,mt)$ be a vertex in $bs^{k}(c[0,1])$. There are two cases to consider:\\

Case 1: $m\in (n,n+1)$ for some $n \in \mathbb{N}$. \\

Since $bs^{k}(c[0,1])$ is a subdivision of $c[0,1]$, if $m\in (n,n+1)$ for some $n \in \mathbb{N}$ then $p |st^{k}(m,mt)| \subset [n,n+1]$.\\

If $z_{n}, z_{n+1}$ lie on the same height, it is clear that since $p\beta(n)$ is distance $\leq \frac{1}{2}$ from $p(z_n)$, we have 
\begin{align*}
pi\beta [m-\frac{1}{4}, m+\frac{1}{4}] \subset pi\beta |st(n)| \subset N_{\frac{1}{2}+\delta}(pz_{n}) 
\end{align*}
$\widetilde{i\beta}[m-\frac{1}{4}, m+\frac{1}{4}]$ lies in the simplex with vertices $z_{n},z_{n+1}$, so the set $pi\beta[m-\frac{1}{4}, m+\frac{1}{4}]  \cup p\widetilde{i\beta} [m-\frac{1}{4}, m+\frac{1}{4}] $ has length at most $\frac{1}{2}+\delta$. Therefore the set $p iH |st^{k}(m,mt)| \cup p Q |st^{k}(m,mt)|$ has length at most $\frac{1}{2}+\delta +\varepsilon <1$. Let $I_{m,mt}$ denote the smallest interval containing this set. $I_{m,mt}$ intersects at most one integer height slice on $M_{X}$.  Let $h$ be the integer height slice it intersects. Otherwise, if the set is completely contained between two integer heights, choose the lower one. We can define $(\widetilde{iH})_{k}(m,mt) = \mathfrak{q}_{h}(z_{[m]})$ where $[m]$ is the closest integer to $m$. \\

If $z_{n}, z_{n+1}$ lie on different heights, by construction, $p\beta[m-\frac{1}{4}, m+\frac{1}{4}] \subset p\beta (|st(n)|\cap|st(n+1)|)$ is at most distance $\frac{1}{2}+\delta$ from both $pz_{n}$ and $pz_{n+1}$. So the set $pi\beta [m-\frac{1}{4}, m+\frac{1}{4}]  \cup p\widetilde{i\beta}[m-\frac{1}{4}, m+\frac{1}{4}] $ has length at most $\frac{2}{3}$. Therefore the set  $p iH |st^{k}(m,mt)| \cup p Q |st^{k}(m,mt)|$ has length at most  $\frac{2}{3}+\varepsilon$. Let $I_{m,mt}$ and $h$ be defined as before, and $(\widetilde{iH})_{k}(m,mt) = \mathfrak{q}_{h}(z_{[m]})$.\\

Case 2: $m=n \in \mathbb{N}$.\\

By construction we have $pi\beta [n-\frac{1}{4}, n+\frac{1}{4}]  \subset pi\beta |st(n)| \subset N_{\frac{1}{2}+\delta} (pz_{n})$, and that $p\widetilde{i\beta}|st^{k}(n)| \subset N_{\frac{1}{3}}(pz_{n})$. Therefore the length of the set $piH |st^{k}(n,nt)| \cup pQ|st^{k}(n,nt)|$ is at most $\frac{1}{3}+\frac{1}{2}+\delta +\varepsilon<1$. Again, we define $(\widetilde{iH})_{k}(n,nt) = \mathfrak{q}_{h} (z_n)$, the projection of $z_n$ to the appropriate height.\\

Observe that in all the cases, the geometric star condition is met for the map $iH$ with respect to $bs^{k}(c[0,1])$, and therefore $(\widetilde{iH})_{k}$ can be extended to the interior of simplices to be a simplicial approximation of $iH$. Furthermore, $(\widetilde{iH})_{k}$ is homotopic to $Q$. To see this, let $\sigma$ be a simplex in $bs^{k}(c[0,1])$, $p'\sigma \in [n,n+1]$ for some $n\in \mathbb{N}$. If $(\widetilde{iH})_{k}\sigma$ lies on one height $h$, then $Q\sigma$ lies in the geometric neighbourhood of the simplex $[\mathfrak{q}_{h}z_{n},\mathfrak{q}_{h}z_{n+1}]$. Alternatively, if $(\widetilde{iH})_{k}\sigma$ spans two integer heights $h,h+1$ then by construction $Q\sigma$ lies in the geometric cylinder $[\mathfrak{q}_{h}z_{n},\mathfrak{q}_{h}z_{n+1},\mathfrak{q}_{h+1}z_{n},\mathfrak{q}_{h+1}z_{n+1}]$ (if the interval containing $pQ\sigma$ intersected either $h,h+1$ then by construction $(\widetilde{iH})_{k}\sigma$ would lie on one height). Therefore in both cases, $Q\sigma$ lies in the geometric neighbourhood of the simplex spanned by $\cup_{v\in vert \sigma} (\widetilde{iH})_{k}(v)$.  \\

Now we show that $Q\sigma$ lies in the geometric neighbourhood of the geometric simplex spanned by the vertices $\cup_{v\in vert N\sigma} (\widetilde{iH})_{k}(v)$, where $N\sigma$ consists of all simplices with $\sigma$ as a face. Since $N(\sigma, bs^{k}{c[0,1]}) \subset N(\iota \sigma,c[0,1])$, where $\iota \sigma$ is the unique simplex in $c[0,1]$ containing $\sigma$, we have that $pN(\sigma, bs^{k}{c[0,1]})\subset [n,n+1]$ (Case 1), unless $p\sigma$ is either $n$ or $n+1$ (Case 2). Let $\tau$ be a simplex in $bs^{k}{c[0,1]}$ containing $\sigma$ as a face. \\

Case 1: Suppose that  $(\widetilde{iH})_{k}\sigma$ lies on one height $h$ and $Q\sigma$ lies in the geometric neighbourhood of the simplex $[\mathfrak{q}_{h}z_{n},\mathfrak{q}_{h}z_{n+1}]$. The image of $(\widetilde{iH})_{k}\tau$ may span at most two integer heights $h,h'$ for $h'\in \{h-1,h+1\}$. If this is the case, then $Q\sigma$ must lie in the geometric cylinder of $[\mathfrak{q}_{h}z_{n},\mathfrak{q}_{h} z_{n+1}, \mathfrak{q}_{h'}z_{n},\mathfrak{q}_{h'} z_{n+1}]$, showing that $h'$ must be the same for all possible $\tau$.  Alternatively, if $(\widetilde{iH})_{k}\sigma$ spans two integer heights $h,h+1$, then $(\widetilde{iH})_{k}\tau$ must also lie in the geometric cylinder $[\mathfrak{q}_{h}z_{n},\mathfrak{q}_{h}z_{n+1},\mathfrak{q}_{h+1}z_{n},\mathfrak{q}_{h+1}z_{n+1}]$.\\

Case 2: Let $p\sigma = n\in \mathbb{N}$. If  $(\widetilde{iH})_{k}\tau$ lies on one height, so does $(\widetilde{iH})_{k}\sigma$. If $(\widetilde{iH})_{k}\tau$ spans two integer heights $h,h'$ for $h'\in \{h-1,h+1\}$, then $Q\sigma \subset Q\tau$ lies in the geometric simplex $[\mathfrak{q}_{h}z_{n}, \mathfrak{q}_{h'}z_{n}]$.\\

 By the argument in point $3$. of Proposition \ref{geosa}, $Q$ is properly homotopic to $(\widetilde{iH})_{k}$. Since the images of $(\widetilde{iH})_{k}$ and $Q$ lie in a $2$-dimensional subcomplex of $M_{X}$, the homotopy is Lipschitz. 

\end{proof}. 

This concludes the preliminaries. We now begin the proof of Proposition \ref{isogeneralray}. \\

\begin{proof} (of Proposition \ref{isogeneralray}) This is a generalisation of the proof of Lemma \ref{isoomega}. The map $s_{l}:\pi_{n}^c(cX,\tau)\rightarrow \pi_{n}^c(cX,\beta)$ is induced by a coarse homotopy equivalence $s_{l}: cX\rightarrow cX$, and is therefore an isomorphism. It therefore suffices to define $\theta$ and show that it is an isomorphism. \\

We identify $[0,1]^n$ with $[-1,1]^n$ in the obvious way. Fix a simplicial structure on $[-1,1]^n$. Recall the definition of the set $U = \{(hx_1,\dots,hx_n,h) \subset c[-1,1]^n \,|\, -\frac{1}{2}\leq x_i \leq \frac{1}{2} \, \forall 1\leq i\leq n\}=cW$ and the map $\mathfrak{r}: U\rightarrow c[-1,1]^n$, the height-wise radial multiplication by $2$.\\

Let us define $\theta: \pi_{n}^c(cX,\beta)\rightarrow \pi_{n}^{L,e}(M,\widetilde{i\beta})$. Let $[f]\in \pi_{n}^c(cX,\beta)$. There exists a $r\in \mathbb{N}$ large enough such that $is_{r}f|st(v)|$ has diameter $<\delta$ for all vertices $v\in c[-1,1]^n$ for our fixed simplicial structure on $c[-1,1]^n$. Consider now the simplicial structure on $c[-1,1]^n$, which we denote by $c\overline{[-1,1]^n}$, which restricts to $c[-1,1]^n$ on the set $U$, and on the complement $c[-1,1]^n\setminus U = c(\partial W \times [0,1])$ is the cone of the product simplicial structure of $\partial{[-1,1]^n} \times [0,1]_{\frac{\varepsilon}{4}}$, where the interval $[0,1]$ is divided into segments of length $\frac{\varepsilon}{4}$. Recall that $H: c[0,1]\rightarrow cX$ is the linear homotopy from $s_{r}\beta$ to $\beta$, and $H_{t}: [1,\infty)\rightarrow cX, H_{t}(h) = H(ht,h)$ its restriction to a time slice. We define $\theta[f]$ as follows:
\begin{align*}
\theta: \pi_{n}^c(cX,\beta)&\rightarrow \pi_{n}^{L,e}(M,\widetilde{i\beta})\\
[f] &\mapsto [(ib^H_{s_{r}\beta,\beta}s_{r}f)^{\sim}]
\end{align*}
ie. on the set $U$ we have a simplicial approximation of $is_{r}f\mathfrak{r}$ and on $c[-1,1]^n\setminus U$, we have $\widetilde{i b^H_{s_{r}\beta,\beta}}$ (by abuse of notation, we use this symbol to denote the map as well as the group homomorphism). The simplicial approximation is with respect to $c\overline{[-1,1]^n}$. \\

We now show that $\theta$ is well-defined. Consider a vertex $v$ in $c\overline{[-1,1]^n}$. If $v$ lies in the interior of $U$ then $|st(v, c\overline{[-1,1]^n})| = |st(v,U)|$. We already know that $is_{r}f$ satisfies the geometric star condition. If $v$ lies in the interior of $c\overline{[-1,1]^n}\setminus U$ then by Lemma \ref{technicalcone} we have that $|st(v, c\overline{[-1,1]^n})| = |st(v,c(\partial W \times [0,1]))| \subset |st(qv, c\partial W)| \times I_{\varepsilon}$, where $q: c(\partial W \times [0,1]) \rightarrow c\partial W$ is the projection and $I_{\varepsilon}$ is an interval of length $\varepsilon$. Let $v= (hx,ht,h)$, where $x\in \partial W, t\in [0,1]$. The map $b^H_{s_{r}\beta,\beta}$ is independent of the $\partial W$-variable, and $H_{t}$ is distance-decreasing for all $t\in[0,1]$. This, along with condition $2$ in the choice of $l$ means that $b^H_{s_{r}\beta,\beta}(|st((hx,h), c\partial W)| \times I_{\varepsilon})\subset N_{\frac{\varepsilon}{2}}(H_{t}[h-1,h+1])$, where $H_{t}[h-1,h+1]$ is a set of diameter less than $\delta$. Therefore the conditions of Lemma \ref{shiftweirdness} are met and $v$ satisfies the geometric star condition. Finally, let $v$ be a vertex on the boundary of $U$. $|st(v, c\overline{[-1,1]^n})|= |st(v,U)| \cup |st(v,c\overline{[-1,1]^n}\setminus U)|$. The set $is_{r}f|st(v,U)|$ has diameter less than $\delta$ and $b^H_{s_{r}\beta,\beta}|st(v,c[-1,1]^n\setminus U)|\subset N_{\frac{\varepsilon}{2}}(is_{r}f|st(v,U)|)$. Therefore, in all cases, the geometric star condition is met and the simplicial approximation $(ib^H_{s_{r}\beta,\beta}s_{r}f)^{\sim}$ exists. Additionally, it can be chosen to agree with our fixed simplicial approximation of $\widetilde{i\beta}$ on the boundary $c\partial[-1,1]^n$. \\

We will show that $\theta$ is a group homomorphism later by defining its inverse. Furthermore, $\theta$ is independent of the following choices: 

\begin{enumerate}
\item Indepedence of choice of simplicial approximation: \\

Suppose we choose another simplical approximation $(ib^H_{s_{r}\beta,\beta}s_{r}f)^{\sim'}$ which restricts to $\widetilde{i\beta}$ on $c\partial[-1,1]^n$. By point $1.$ of Proposition \ref{geosa}, $(ib^H_{s_{r}\beta,\beta}s_{r}f)^{\sim'}$ is proper Lipschitz homotopic to $(ib^H_{s_{r}\beta,\beta}s_{r}f)^{\sim}$ relative to the boundary $c\partial[-1,1]^n$. 
\item Indepedence of choice of $r$:\\

 Suppose that $(ib^H_{s_{r}\beta,\beta}s_{r}f)^{\sim}$ already exists. We show that $(ib^H_{s_{r}\beta,\beta}s_{r}f)^{\sim}$ is homotopic to $(ib^{H^{r+r',0}}_{s_{r+r'}\beta,\beta}s_{r+r'}f)^{\sim}$ for all $r'\in \mathbb{N}$. To that end, consider the linear homotopy $\mathbf{H}^{r'}: I_{p}(cX)\rightarrow cX$ between $\id$ and $s_{r'}$. Let $F=b^H_{s_{r}\beta,\beta}s_{r}f$. We define
\begin{align*}
G: c([-1,1]^n \times [0,1])&\rightarrow cX\\
(hx,ht,h)&\mapsto \mathbf{H}^{r'}(pF(hx,h) p'F(hx,h), pF(hx,h)t, pF(hx,h)) 
\end{align*} 
This is a coarse homotopy between $b^H_{s_{r}\beta,\beta}s_{r}f$ and $s_{r'}b^H_{s_{r}\beta,\beta}s_{r}f =b^{s_{r'}H}_{s_{r+r'}\beta,s_{r'}\beta}s_{r+r'}f $.\\

Since the map $F$ is coarse, the composition $p\circ F$ is large-scale Lipschitz with constant $A$. Let $(hx,h)\in c[-1,1]^n$ and $(hy,h)\in c\partial[-1,1]^n$ with $y$ the closest point in $\partial[-1,1]^n$ to $x$. We have
\begin{align*}
|pF(hx,h)| &\leq|pF(hx,h)-pF(hy,h)|+|pF(hy,h)| \\
&< A d_{c[-1,1]^n}((hx,h),(hy,h)) + A + h \\
&<Ah +A+h \leq (2A +1)h
\end{align*}
This means that the homotopy $G$ has speed $< (2A+1)$ in the time variable. Choose $\varepsilon_{F}<\frac{1}{(2A+1)16}$. Consider now the simplicial structure on $c([-1,1]^n \times [0,1])$ which is the cone of the simplicial structure $\overline{[-1,1]^n} \times [0,1]_{{\varepsilon}_{F}}$ where $[0,1]_{\varepsilon_{F}}$ has vertices spaced distance $\varepsilon_{F}$ apart. Let $v=(hx,ht,h)$ be a vertex in $c(\overline{[-1,1]^n} \times [0,1]_{{\varepsilon}_{F}})$. We have that $|st(v)| \subset |st((hx,h),c\overline{[-1,1]^n})|\times I_{4\varepsilon_{F}}$ where $I_{4\varepsilon_{F}}$ is an interval of length $4\varepsilon_{F}$. Therefore, $G|st(v)|\subset N_{2(2A+1)\varepsilon_{F}}(\mathbf{H}^{r'}_{t}F|st(hx,h)|)$. If $(hx,h)$ lies in the interior of  $ U$ then $\mathbf{H}^{r'}_{t}F|st(hx,h)|$ is a set of diameter $<\delta$. If it lies in the interior of $c\overline{[-1,1]^n}\setminus U$ then 
\begin{align*}
N_{2(2A+1)\varepsilon_{F}}(\mathbf{H}^{r'}_{t}F|st(hx,h)|) &\subset N_{2(2A+1)\varepsilon_{F}}(\mathbf{H}^{r'}_{t} N_{\frac{\varepsilon}{2}}(H_{t} [h-1,h+1])) \\
&\subset N_{2(2A+1)\varepsilon_{F}} N_{\frac{\varepsilon}{2}} (\mathbf{H}^{r'}_{t}(H_{t} [h-1,h+1])) \\
&= N_{2(2A+1)\varepsilon_{F} + \frac{\varepsilon}{2}} (\mathbf{H}^{r'}_{t}(H_{t} [h-1,h+1]))
\end{align*}
where $\mathbf{H}^{r'}_{t}(H_{t} [h-1,h+1])$ has diameter $<\delta$. Finally, if $(hx,h)$ lies on the boundary of $U$, then 
\begin{align*}
N_{2(2A+1)\varepsilon_{F}}(\mathbf{H}^{r'}_{t}F|st(hx,h)|) &\subset N_{2(2A+1)\varepsilon_{F}}\mathbf{H}^{r'}_{t}N_{\frac{\varepsilon}{2}} (is_{r}f |st((hx,h),U)|)\\
&\subset N_{2(2A+1)\varepsilon_{F} + \frac{\varepsilon}{2}}(\mathbf{H}^{r'}_{t}is_{r}f |st((hx,h),U)|)
\end{align*}
where $\mathbf{H}^{r'}_{t}is_{r}f |st((hx,h),U)|$ has diameter $<\delta$. Therefore the conditions of Lemma \ref{shiftweirdness} are satisfied and $\widetilde{iG}$ exists.\\

Cosider now the map $\widetilde{iG}$ restricted to the set $c([-1,1]^n \times \{1\}) \cup c(\partial[-1,1]^n \times [0,1])$ (see Figure \ref{fig:geniso1}). This can be identified with $c[-1,1]^n$. Therefore, $\widetilde{iG}$ gives us a homotopy rel boundary between $(ib^{H}_{s_{r}\beta,\beta}s_{r}f)^{\sim}$ and $(i b^{G_{\partial}}_{s_{r'}\beta,\beta}b^{s_{r'}H}_{s_{r+r'}\beta,s_{r'}\beta}s_{r+r'}f)^{\sim}$ where the second simplicial approximation is with respect to the structure $c(\overline{[-1,1]^n} \times \{1\}) \cup c(\partial\overline{[-1,1]^n} \times [0,1]_{\varepsilon_{F}})$ on  $c[-1,1]^n$ and $G_{\partial}$ denotes the homotopy $G$ restricted to $c(\partial[-1,1]^n \times [0,1])$. \\

\begin{figure}
\centering
  \includegraphics[width=0.6\linewidth]{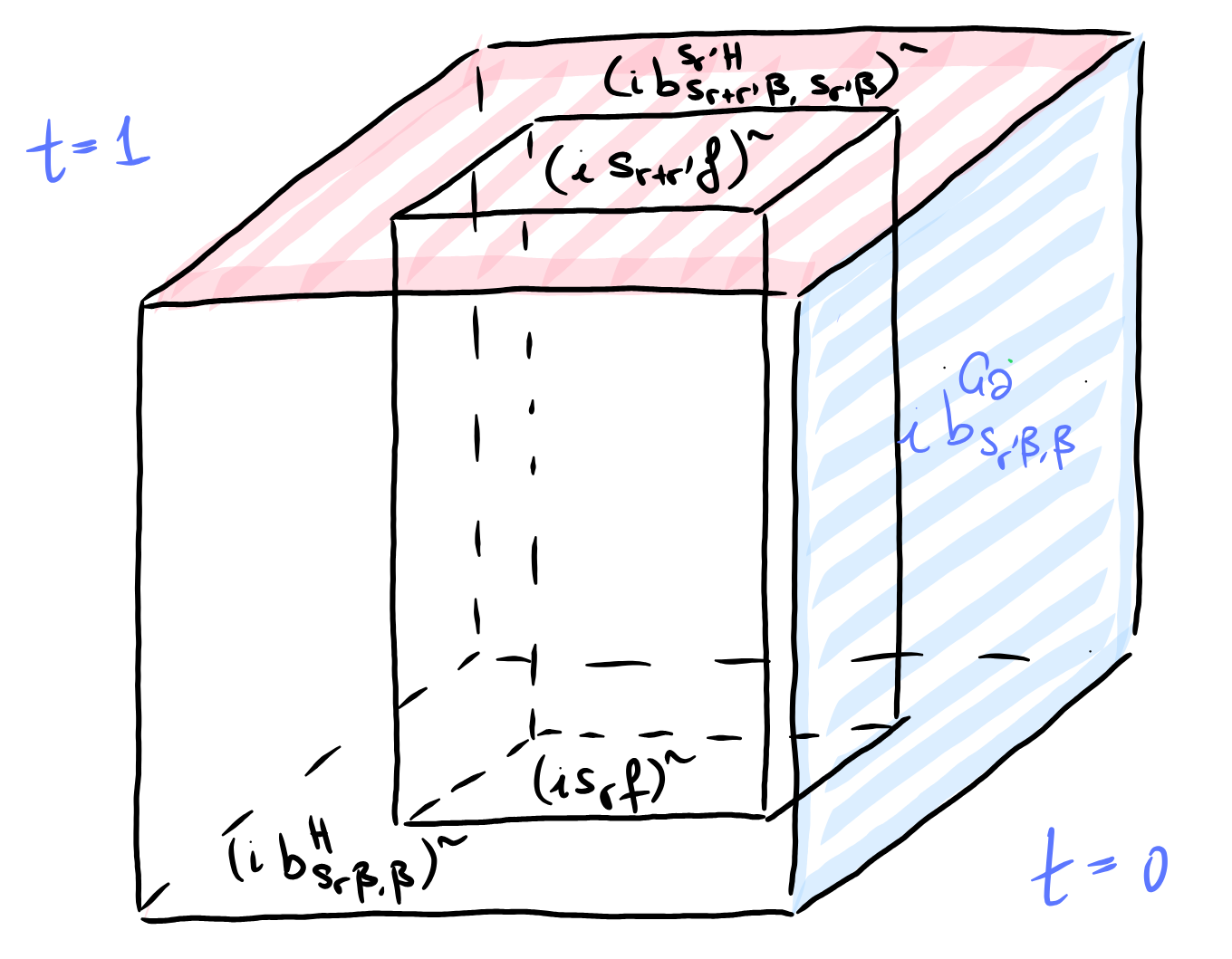}
  \caption{The homotopy $\widetilde{iG}$. The set $c([-1,1]^n \times \{1\})$ is red, and one of the faces of $c(\partial[-1,1]^n \times [0,1])$ is blue.}
  \label{fig:geniso1}
\end{figure}

We show now that $(i b^{G_{\partial}}_{s_{r'}\beta,\beta}b^{s_{r'}H}_{s_{r+r'}\beta,s_{r'}\beta}s_{r+r'}f)^{\sim}$ is homotopic to $(ib^{H^{r+r',0}}_{s_{r+r'}\beta,\beta}s_{r+r'}f)^{\sim}$. To see this, observe that both maps agree when restricted to the set $U$ in $c(\overline{[-1,1]^n} \times \{1\})$ and $c\overline{[-1,1]^n}$ and  respectively. On the complement, both maps are simplicial approximations of the linear homotopy $H^{r+r',0}$ between $s_{r+r'}\beta$ and $\beta$ with respect to different simplicial structures: the complement of $U$ in  $c(\overline{[-1,1]^n} \times \{1\}) \cup c(\partial\overline{[-1,1]^n} \times [0,1]_{\varepsilon_{F}})$ can be identified with $c(\partial W \times ([0,1]_{\frac{\varepsilon}{4}}\cup [0,1]_{\varepsilon_{F}}))$
where the first interval is divided into segments of length $\frac{\varepsilon}{4}$ and the second into segments of length $\varepsilon_{F}$. There is a bi-Lipschitz homeomorphism 
\begin{align*}
\varphi: c(\overline{[-1,1]^n} \times \{1\}) \cup c(\partial\overline{[-1,1]^n} \times [0,1]_{\varepsilon_{F}})\rightarrow c\overline{[-1,1]^n}
\end{align*}
such that $b^{H^{r+r',0}}_{s_{r+r'}\beta,\beta}s_{r+r'}f \circ \varphi = b^{G_{\partial}}_{s_{r'}\beta,\beta}b^{s_{r'}H}_{s_{r+r'}\beta,s_{r'}\beta}s_{r+r'}f$. By point 4. of Proposition $\ref{geosa}$, there exists a proper Lipschitz homotopy between their simplicial approximations which can be chosen to be base ray-preserving. Putting it altogether, this shows that 
\begin{align*}
[(ib^H_{s_{r}\beta,\beta}s_{r}f)^{\sim}]=[(i b^{G_{\partial}}_{s_{r'}\beta,\beta}b^{s_{r'}H}_{s_{r+r'}\beta,s_{r'}\beta}s_{r+r'}f)^{\sim}] = [(ib^{H^{r+r',0}}_{s_{r+r'}\beta,\beta}s_{r+r'}f)^{\sim}] 
\end{align*}
as claimed. 

\item Independence of choice of representative $f$ in its coarse homotopy class.\\

Suppose that $[g]=[f]$ in $\pi_n^c(cX,\beta)$. This means that there is a coarse homotopy $F: (c[-1,1]^{n+1}, c(\partial [-1,1]^{n}\times [0,1]))\rightarrow (cX,\omega)$ which restricts to $f$ and $g$ on $c([-1,1]^n \times \{0\})$ and $c([-1,1]^n \times \{1\})$ respectively (here, $[0,1]$ is the simplicial complex with one edge). There exists a $S>0$ such that $F(|st(v)|)$ has diameter $<S$ for all vertices $v\in c([-1,1]^{n} \times [0,1])$, and a $r$ large enough so that $s_{r} F(|st(v)|)$ has diameter $<\delta$.\\

 Now give $c([-1,1]^n \times [0,1])$ the simpicial structure $c(\overline{[-1,1]^n} \times [0,1])$. Let us define the homotopy $G: c(\overline{[-1,1]^n} \times [0,1]) \rightarrow cX$: on the set $c(W\times [0,1])$ we take the map $F$ and on its complement, we take $b^H_{s_{r}\beta,\beta} \circ q$, where $q: c(\overline{[-1,1]^n} \times [0,1])\rightarrow c\overline{[-1,1]^n}$ is the projection. This is a homotopy between $b^H_{s_{r}\beta,\beta}s_{r}f$ and $b^H_{s_{r}\beta,\beta} s_{r}g$. The simplicial approximation $\widetilde{iG}$ exists, and can be chosen to be fixed as $\widetilde{i\beta}$ on $c(\partial \overline{[-1,1]^n} \times [0,1])$. This shows that $\theta[f] = \theta[g]$. 

\item Independence of barycentric subdivision of $c\overline{[-1,1]^n}$. \\

Fix $r$ sufficiently large such that the simplicial approximation $(ib^H_{s_{r}\beta,\beta}s_{r}f)^{\sim}$ exists. Let $bs^k(c\overline{[-1,1]^n})$ denote the $k$-th barycentric subdivision of $c\overline{[-1,1]^n}$. By proposition \ref{geosa} the simplicial approximation $(ib^H_{s_{r}\beta,\beta}s_{r}f)^{\sim}_{k}$ of $ib^H_{s_{r}\beta,\beta}s_{r}f$ with respect to $bs^k(c\overline{[-1,1]^n})$ exists and is proper Lipschitz homotopic to $(ib^H_{s_{r}\beta,\beta}s_{r}f)^{\sim}$ via a homotopy $G$. Let $J^{k}$ be $G$ restricted to the boundary, ie the homotopy $J^{k}: c(\partial[-1,1]^n) \times [0,1]\rightarrow M$ between $J^{k}_{0}= {(ib^H_{s_{r}\beta,\beta}s_{r}f)^{\sim}_{k}}_{|c\partial[-1,1]^n}$ and  $J^{k}_{1}= {(ib^H_{s_{r}\beta,\beta}s_{r}f)^{\sim}}_{|c\partial[-1,1]^n}$. Let $B^{J^{k}}(ib^H_{s_{r}\beta,\beta}s_{r}f)^{\sim}_{k}$ be the map which is $(ib^H_{s_{r}\beta,\beta}s_{r}f)^{\sim}_{k}$  when restricted to the set $U$, and the homotopy $J^{k}$ (rescaled) on the set $c[-1,1]^n \setminus U$, which can be identified topologically with $c(\partial[-1,1]^n \times [0,1])$. 

We have that
\begin{align*}
[(ib^H_{s_{r}\beta,\beta}s_{r}f)^{\sim}] = [B^{J^{k}}(ib^H_{s_{r}\beta,\beta}s_{r}f)^{\sim}_{k}] 
\end{align*}

For any other choice of simplicial approximation $(ib^H_{s_{r}\beta,\beta}s_{r}f)^{\sim'}_{k}$ of $ib^H_{s_{r}\beta,\beta}s_{r}f$ with respect to $bs^k(c\overline{[-1,1]^n})$ we have 
\begin{align*}
[B^{J}(ib^H_{s_{r}\beta,\beta}s_{r}f)^{\sim}_{k}] =  [B^{J'}(ib^H_{s_{r}\beta,\beta}s_{r}f)^{\sim'}_{k}]
\end{align*}
This follows from the arguments of points $1$. and $2$. of Proposition \ref{geosa}. The details are left as an exercise for the reader. 
\end{enumerate}

Now let us define the inverse $\lambda: \pi_{n}^{L,e}(M,\widetilde{i\beta})\rightarrow \pi_{n}^{c}(cX,\beta)$, which is the composition
\begin{equation*}
	\begin{tikzcd}[row sep=0.5em, column sep=1em]
		\lambda: \pi^{L,e}_{n}(M,\widetilde{i\beta})  \arrow[r, "\id"] & \pi_{n}^c(M,\widetilde{i\beta})   \arrow[r, "R"] & \pi_{n}^c(cX,R\widetilde{i\beta})  \arrow[r] & \pi_{n}^{c}(cX,\beta) \\
		{[g]} \arrow[r,|->] & {[g]} \arrow[r,|->]& {[R\circ g]} \arrow[r,|->] & {[R\circ g]}
	\end{tikzcd}
\end{equation*}

In words, we take a representative $g$ of $[g]\in  \pi_{n}^{L,e}(M,\widetilde{i\beta})$, consider it as a coarse homotopy class, compose with the map $R$, and modify the boundary from $R\widetilde{i\beta}$ to $\beta$. This is well-defined since $R\widetilde{i\beta}$ is close to $\beta$. $\lambda$ is a group homomorphism.\\

Now we show that $\theta$ and $\lambda$ are inverse to each other. First we compute the composition $\lambda \theta: \pi_{n}^{c}(cX,\beta)\rightarrow \pi_{n}^{c}(cX,\beta)$. For a picture of this, see Figure \ref{fig:geniso2}. \\

\begin{figure}
\centering
  \centering
  \includegraphics[width=0.6\linewidth]{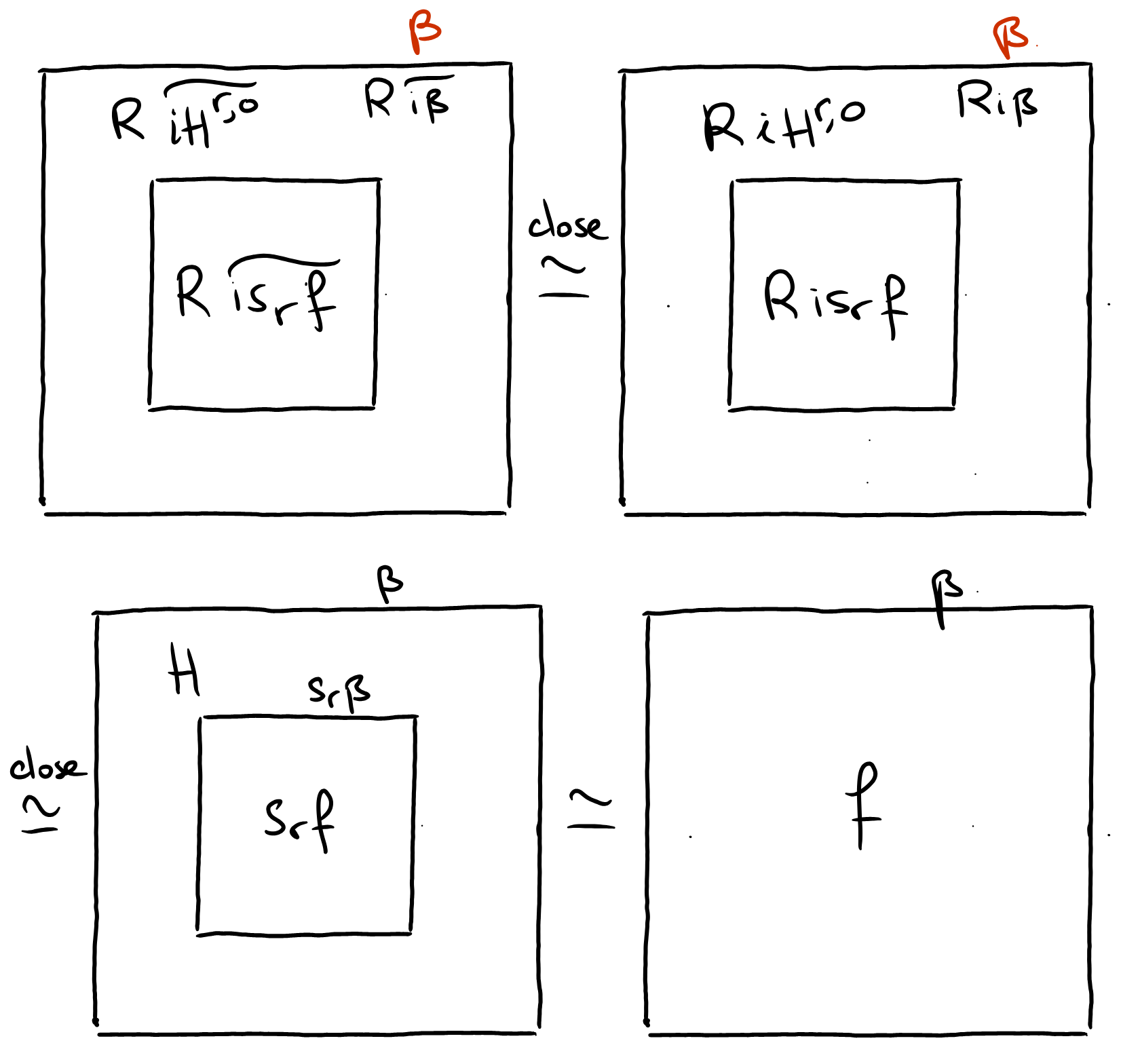}
  \caption{A representative of $\lambda\theta[f]$. The modified map on the boundary is denoted in red.}
  \label{fig:geniso2}
\end{figure}

Let $[f]\in \pi_{n}^{c}(cX,\beta)$. The image $\lambda\theta[f]$ is represented by the homotopy class 
\begin{align*}
 [R(ib^H_{s_{r}\beta,\beta}s_{r}f)^{\sim}]
\end{align*}
Since $(ib^H_{s_{r}\beta,\beta}s_{r}f)^{\sim}$ is close to $ib^H_{s_{r}\beta,\beta}s_{r}f$ and $Ri$ is close to the identity, we have
\begin{align*}
R(ib^H_{s_{r}\beta,\beta}s_{r}f)^{\sim} \simeq Rib^H_{s_{r}\beta,\beta}s_{r}f \simeq b^H_{s_{r}\beta,\beta}s_{r}f \simeq f
\end{align*}
Therefore $\lambda\theta=\id$.\\

Now we compute $\theta \lambda:\pi_{n}^{L,e}(M,\widetilde{i\beta})\rightarrow \pi_{n}^{L,e}(M,\widetilde{i\beta})$. Let $f: (c[-1,1]^n,c\partial[-1,1]^n)\rightarrow (M,\widetilde{i\beta})$ be a proper Lipschitz map representing a class in $\pi_{n}^{L,e}(M, \widetilde{i\beta})$. Since $c[-1,1]^n$ is uniformly bounded, there exists a $k$ large enough such that for every vertex $v\in bs^{k}(c [-1,1]^n)$, we have that  $f(|st(v)|) \subset \langle gst(z_v,t_{v}) \rangle$ where $(z_v,t_{v})$ is the closest vertex to the image of $v$. \\

Choose $r$ sufficiently large such that $(is_{r}Rf)^{\sim}$ is well-defined. Consider now $ib^{\id}_{s_{r}\beta,s_{r}\beta} s_{r}Rf$. We can assume that $r$ is large enough so that this map has a simplicial approximation with respect to the fixed simplicial structure on $c[-1,1]^n$. The underlying identity map on sets is a bi-Lipschitz homeomorphism 
\begin{align*}
\varphi: bs^{k}(c\overline{[-1,1]^n})\rightarrow  bs^{k}(c[-1,1]^n)
\end{align*}
which is a simplicial identity map on the subcomplex $ bs^{k}(c\partial[-1,1]^n)$. Therefore, by point 4. of Proposition $\ref{geosa}$, there exists a proper Lipschitz homotopy between $(ib^{\id}_{s_{r}\beta,s_{r}\beta} s_{r}Rf)^{\sim}_{k}$ with respect to these two different simplicial structures, which can be chosen to be relative to $ bs^{k}(c\partial[-1,1]^n)$.\\

This shows that we can represent $\theta\lambda [f]$ by the map
\begin{align*}
B^{J^{k}} (b^{H^{r}}_{s_{r}\beta,\beta} b^{\id}_{s_{r}\beta,s_{r}\beta} s_{r}Rf)^{\sim}_{k}
\end{align*}
where $(\cdot)^{\sim}_{k}$ is simplicial approximation with respect to the simplicial structure which is $bs^{k}(c\overline{[-1,1]^n})$when restricted to the domain of $ b^{\id}_{s_{r}\beta,s_{r}\beta} s_{r}Rf$, and $bs^{k}(c\overline{[-1,1]^n}\setminus U)$ on the domain of $b^{H^{r}}_{s_{r}\beta,\beta} $. \\

Let us describe the different components of this map (see Figure \ref{fig:geniso3}).

\begin{enumerate}
\item By the same argument as in Theorem  \ref{isoomega} we have that for a sufficiently small subdivision of $c[-1,1]^n$, $is_{r}Rf$ and $D_{r}f$ have the same simplicial approximation, where $D_{r}: M\rightarrow M$ is the map which divides heights by $2^{r}$. Therefore segment $1.$ can be replaced with $(D_{r}f)^{\sim}_{k}$. 
\item This segment is a simplicial approximation of the constant homotopy of the base ray $is_{r}\beta$. The replacement on segment $1$. restricts to a replacement on the boundary $(\id_{is_{r}\beta})^{\sim}_{k} = (\id_{D_{r}\widetilde{i\beta}})^{\sim}_{k}$.   
\item This is a simplicial approximation of the change of base ray $ib^{H}_{s_{r}\beta, \beta}$ with respect to $bs^{k}(c\overline{[-1,1]^n}\setminus U)$. We have the following corollary of Lemma  \ref{twohomotopies}. 
\begin{cor} $(ib^{H^{r}}_{s_{r}\beta, \beta})^{\sim}_{k}$ is proper Lipschitz homotopic to $b^{Q}_{D_{r}\widetilde{i\beta}, \widetilde{i\beta}}$. 
\end{cor}
To see this, we identify the set $c[-1,1]^n\setminus U$ with $c(\partial[-1,1]^n \times [0,1])$ in the obvious way. We can choose a $k$ large enough such that for a vertex $v=(mx,mt,m)$ in $bs^{k}(c\overline{[-1,1]^n}\setminus U)$ the image of $|st^{k}(v)|$ under the map 
\begin{align*}
P: c(\partial[-1,1]^n \times [0,1]) &\rightarrow c[0,1]\\
(hx,ht,h)&\mapsto (h,ht)
\end{align*}
is contained within the set $[m-\frac{1}{4}, m+\frac{1}{4}] \times [mt-\frac{\varepsilon}{2}, mt+\frac{\varepsilon}{2}]$. Observe that $ib^{H^{r}}_{s_{r}\beta, \beta}=iH \circ P$ and $b^{Q}_{D_{r}\widetilde{i\beta}, \widetilde{i\beta}} = Q\circ P$. Then the argument from Lemma \ref{twohomotopies} applies directly. \\
\item This is the change of base ray obtained from the homotopy ${J^{k}}$. By point $4$. (independence of barycentric subdivision) of this proposition, we get the homotopy between $(ib^{H^{r}}_{s_{r}\beta, \beta})^{\sim}_{k}$ and  $b^{Q}_{D_{r}\widetilde{i\beta}, \widetilde{i\beta}}$, restricted to the boundary $c(\partial[-1,1]^n \times \{1\})$, which we denote by $B^{J'}$, can be chosen to agree with $B^{J^{k}}$. 
\end{enumerate}

Putting this altogether, we see that everything in segments $1,2,3$ is a simplicial approximation of the proper Lipschitz map 
\begin{align*}
b^{Q}_{D_{r}\widetilde{i\beta},\widetilde{i\beta}}b^{\id}_{D_{r}\widetilde{i\beta},D_{r}\widetilde{i\beta}}D_{r}f
\end{align*}
Segment $2$ can be deleted, as can the segment outside of the blue boarder, as we run the homotopy and then the inverse. We then obtain the fourth image, which is obviously homotopic to $f$. Therefore $\theta\lambda = \id$. This concludes the proof.

\begin{remark} We defined the simplicial approximations  $(D_{r}f)^{\sim}_{k}$ and $(\widetilde{iH^{r}})_{k}$ so that it is not immediately clear they are compatible on their mutual boundary. We introduced segment $2$ to ensure that they "glue together" in the appropriate way, ie. segment $1$ has no influence on the argument in segment $3$ and vice versa. 
\end{remark}

\begin{figure}[H]
\centering
  \centering
  \includegraphics[width=\linewidth]{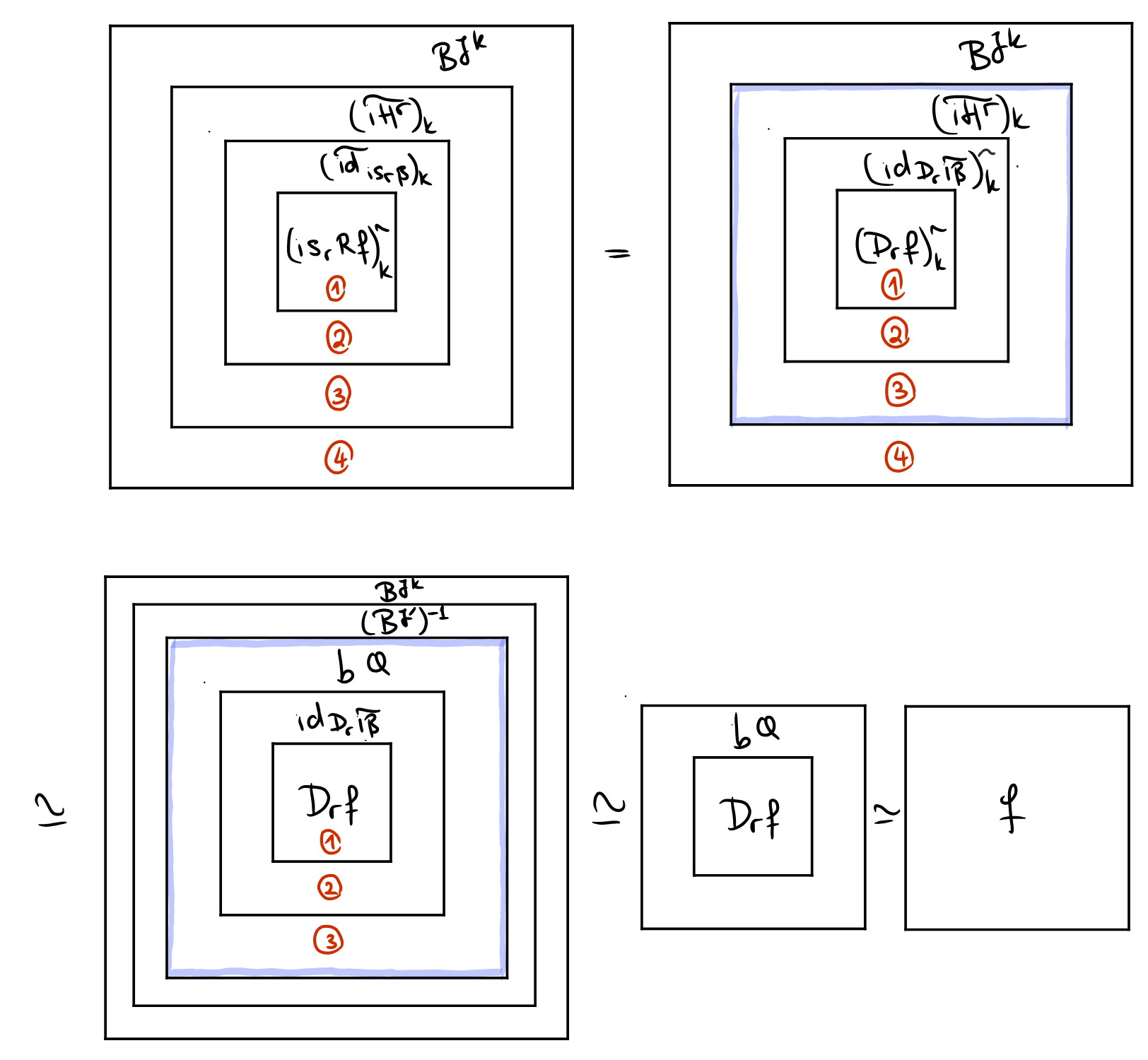}
  \caption{A representative of the composition $\theta \lambda[f]$.}
  \label{fig:geniso3}
\end{figure}
\end{proof}

\begin{prop} \label{psiexistence} Let $\mathfrak{a}: X\rightarrow Y$ be a coarse map with $\mathfrak{a}\omega=\tau$. There exists a $l\in \mathbb{N}$ and a group homomorphism 
\begin{align*}
\psi: \pi_{n}^{e}(M_X,\omega)\rightarrow \pi_{n}^{e}(M_Y,\widetilde{is_{l}\tau})
\end{align*}
induced by a proper, locally Lipschitz, cellular map $\psi: (M_X,\omega)\rightarrow (M_Y, \widetilde{is_{l}\tau})$ such that the diagram
\begin{figure}[H]
\center
\begin{tikzcd}
&\pi_{n}^{c}(cX,\omega) \arrow{r}{\mathfrak{a}}  \arrow{d}{\theta}
&\pi_{n}^{c}(cY,\tau) \arrow{d}{\theta \circ s_{l}} 
  \\ 
& \pi_{n}^{L,e}(M_X,\omega) \arrow{d}{\mathfrak{d}}
& \pi_{n}^{L,e}(M_Y,\widetilde{is_{l}\tau})  \arrow{d}{\mathfrak{d}} 
\\
& \pi_{n}^{e}(M_{X},\omega) \arrow{r}{\psi} 
&\pi_{n}^{e}(M_{Y},\widetilde{is_{l}\tau}) 
\end{tikzcd}
\end{figure}
commutes. If $\mathfrak{a}$ induces an isomorphism $\mathfrak{a}: \pi_{n}^{c}(cX,\omega)\rightarrow \pi_{n}^{c}(cY,\tau)$ for all $n\geq 0$ then $\psi:\pi_{n}^{e}(M_X,\omega)\rightarrow \pi_{n}^{e}(M_Y,\widetilde{is_{l}\tau})$ is an isomorphism for all $n\geq 0$. 
\end{prop}

\begin{proof} Fix an additional constant $\kappa<\frac{1}{8}$. For this proof, we use the same symbol $i$ to denote $i_{X}: cX\rightarrow M_{X}$ and $i_Y: cY\rightarrow M_{Y}$ (resp. $R$ for $R_{X}: M_{X}\rightarrow cX$ and $R_{Y}:M_{Y}\rightarrow cY$). \\

Consider a vertex $v$ of $M_{X}$ and its closed geometric star $|gst(v)|$. Since the map $\mathfrak{a} \circ R$ is coarse, there exists a $l$ large enough, independent of $v$, such that $i s_{l} \mathfrak{a} R(|gst(v)|)\subset \langle gst(z_v) \rangle$ for some vertex $z_v\in vert M_Y$. Additionally, let $l$ be large enough so that for $U\subset cX$, a set of diameter $<\kappa$, we have $s_{l}\mathfrak{a}(U)$ has diameter $<\delta$. Let $\beta = s_{l}\tau$ and let $\widetilde{i\beta}$ denote the fixed simplicial approximation of $is_{l}\tau$, as previously defined. \\

Recall that $M_{X} = \tilde{M}_{X}/_{\sim}$ is constructed as the quotient of 
\begin{align*}
\tilde{M}_{X} = \coprod_{h\in \mathbb{N}_{0}} |\mathcal{U}_{h}| \times [2^h,2^{h+1}] = \coprod_{h\in \mathbb{N}_{0}} \tilde{M}^{h}
\end{align*}
 under the equivalence relation obtained from the simplicial gluing maps $\phi_{h}: |\mathcal{U}_{h}| \times \{2^h\}\rightarrow |\mathcal{U}_{h-1}| \times \{2^{h}\}$. Let $\pi: \tilde{M}_{X}\rightarrow M_{X}$ denote the quotient map.\\

To define $\psi: M_{X}\rightarrow M_{Y}$, we must extend the definition of geometric simplicial approximation to $M_{X}$, which is not a simplicial complex. To do this, we choose a particular CW structure on $M_{X}$ obtained from a choice of product simplicial structure on each $\tilde{M}^{h}$. To this end, choose a total order on the vertices of each $|\mathcal{U}_{h}|$ such that if $v<w$ in $|\mathcal{U}_{h}|$ then $\phi_{h}(v)<\phi_{h}(w)$ in $|\mathcal{U}_{h-1}|$. This can always be accomplished by inducting over $h\in \mathbb{N}_{0}$. This restricts to a total order on each simplex of  $|\mathcal{U}_{h}|$, and induces a canonical simplicial structure on the product $ |\mathcal{U}_{h}| \times [2^h,2^{h+1}]$ with vertices at integer heights. \\

We define the open star of $v$ in $M_{X}$ as  
\begin{align*}
\langle st(v)\rangle = \pi(\bigcup_{w\in \pi^{-1}(v)}\langle st(w)\rangle)
\end{align*}
where $\langle st(w)\rangle$ denotes the star of the vertex $w$ in $\tilde{M}_{X} = \coprod_{h\in \mathbb{N}_0} |\mathcal{U}_{h}| \times [2^h,2^{h+1}]$. Note that since $\pi^{-1} \langle st(v)\rangle = \bigcup_{w\in \pi^{-1}(v)}\langle st(w)\rangle$, $\langle st(v)\rangle$ is an open set in $M_{X}$. We define the closed star $ |st(v)|:= cl\langle st(v)\rangle$ as its closure. Observe that since $|st(v)|\subset |gst(v)|$ we have $i s_{l} \mathfrak{a} R (|st(v)|)\subset \langle gst(z_v) \rangle$ for all vertices $v$. \\

We let $\widetilde{is_{l}\mathfrak{a}R}(v) = z_v$ on vertices and extend skeleta by skeleta in the following way. Assume that $\widetilde{is_{l} \mathfrak{a} R}$ is defined on the $(k-1)$-skeleton of $\tilde{M}_{X}$ and agrees with gluing maps, so that the induced map on the quotient $M_{X}$ is well-defined. Additionally, assume that $\widetilde{is_{l} \mathfrak{a} R}(\Delta^{k-1}) \subset \mathscr{C}_{\Delta^{k-1}}$ where $ \mathscr{C}_{\Delta^{k-1}}$ is the geometric core of the geometric simplex spanned by $\cup_{v\in vert \Delta^{k-1}} z_{v}$, and that $\widetilde{is_{l} \mathfrak{a} R}_{|\Delta^{k-1}}: \Delta^{k-1}\rightarrow \mathscr{C}_{\Delta^{k-1}}$ is Lipschitz, with respect to the product metric on $\Delta^{k-1} \subset \tilde{M}^{h}$ and the path metric on  $\mathscr{C}_{\Delta^{k-1}}$. Consider a $k$-simplex $\Delta^k$ in $\tilde{M}_{X}$. By previous arguments (see Theorem \ref{isoomega}) we have that $\widetilde{is_{l}\mathfrak{a}R}(\partial \Delta^{k})\subset \mathscr{C}_{\Delta^k}$ and that $\widetilde{is_{l}\mathfrak{a}R}_{|\partial \Delta^{k}}: \partial \Delta^{k}\rightarrow \mathscr{C}_{\Delta^{k}}$ is Lipschitz, with respect to the metric on $\Delta^{k} \subset  \tilde{M}^{h}$ and the path metric on $\mathscr{C}_{\Delta^{k}}$ . For the induction step, there are two cases to consider. 

\begin{enumerate}
\item $\pi$ maps the interior of $\Delta^k$ homeomorphically onto its image in $M_{X}$. 
\item $\Delta^k$ lies completely within a gluing component, ie $\Delta^k\subset |\mathcal{U}_{h}| \times \{2^h\}$ for some $h\in \mathbb{N}_{0}$.
\end{enumerate}

In the first scenario, we define the extension to the interior of $\Delta^k$ in the usual way, ie. extend the map on the boundary to the interior by using a Lipschitz deformation retract of the geometric core $\mathscr{C}_{\Delta^k}$ to a point. In the second scenario, recall that the slice $|\mathcal{U}_{h-1}| \times \{2^h\}$ maps homeomorphically onto its image under $\pi$. Let $\phi_{h}: |\mathcal{U}_{h}| \times \{2^h\} \rightarrow |\mathcal{U}_{h-1}| \times \{2^h\}$ be the gluing map. By the first case, or by the induction hypothesis, $\widetilde{is_{l}\mathfrak{a}R}_{|\phi_{h}(\Delta^k)}$ is already defined so we can let $\widetilde{is_{l}\mathfrak{a}R}_{|\Delta^k} = \widetilde{is_{l}\mathfrak{a}R}_{|\phi_{h}(\Delta^k)} \circ \phi_{h}$, and by construction this descends to the quotient. Therefore, we obtain an extension to the $k$-skeleton, with the property that $\widetilde{is_{l}\mathfrak{a}R}(\Delta^{k}) \subset \mathscr{C}_{\Delta^{k}}$, $\widetilde{is_{l}\mathfrak{a}R}_{|\Delta^{k}}: \Delta^{k}\rightarrow \mathscr{C}_{\Delta^{k}}$ is Lipschitz for all $k$-simplices $\Delta^{k}$ (this uses that $\phi_{h}$ is $1$-Lipschitz). By induction, we obtain a globally defined map $\widetilde{is_{l}\mathfrak{a}R}: \tilde{M}_{X}\rightarrow M_{Y}$ which is locally Lipschitz with respect to the set of simplices in $\tilde{M}_{X}$ and descends to a map on the quotient. We denote the map on the quotient also by $\widetilde{is_{l}\mathfrak{a}R}$ and call this the \textit{geometric cellular approximation} of $is_{l}\mathfrak{a}R$.\\

To justify this name, note that for a $k$-simplex $\Delta^k$, $\mathscr{C}_{\Delta^k}$ lies in the $k$-skeleton of the CW structure of $M_{Y}$ given by Example \ref{CW}. Therefore the geometric cellular approximation $\widetilde{is_{l}\mathfrak{a}R}$ is cellular (note that $M_{X}$ and $M_{Y}$ are given different CW structures). Additionally, $\widetilde{is_{l}\mathfrak{a}R}$ is close to $is_{l}\mathfrak{a}R$, and therefore proper. \\

Since $\tilde{M}^h$ is a finite simplicial complex, there exists a Lipschitz constant $m_{h}$ for $\widetilde{is_{l}\mathfrak{a}R}_{|\tilde{M}^h}$. We assume that the sequence $\{m_{h}\}_{h\in \mathbb{N}_{0}}$ is increasing. \\

Consider the finite filtration $\{M_{[1,2^{h}]}\}_{h\in \mathbb{N}}$ of $M_{X}$. We show that $\widetilde{is_{l}\mathfrak{a}R}$ is locally Lipschitz with respect to $\{M_{[1,2^{h}]}\}_{h\in \mathbb{N}}$. Observe that $M_{[1,2^{h}]}\rightarrow M_{X}$ is an isometric embedding. Let  $x$ and $y$ be two points in $M_{[1,2^{h}]}$. A geodesic $\gamma_{xy}$ between them takes the form $\sigma_{x}\ast \eta \ast \sigma_{y}^{-1}$, where $\sigma_{x}, \sigma_{y}$ are straight segments down to a common slice, and $\eta$ is a geodesic path between the projections of $x$ and $y$ onto this slice. Write $\sigma_{x}= \sigma^{1}_{x}\ast \dots \ast \sigma^{j_{x}}_{x}$ as a concatenation of segments that lie entirely in one $\pi(|\mathcal{U}_{k}| \times [2^k,2^{k+1}])$, and analogously for $\sigma_{y}$. We obtain 
\begin{align*}
\gamma_{xy} = \sigma^{1}_{x}\ast \dots \ast \sigma^{j_{x}}_{x} \ast \eta \ast (\sigma^{1}_{y}\ast \dots \ast \sigma^{j_{y}}_{y})^{-1}
\end{align*}
Each of these segments has a lift to $\tilde{M}_{X}$ such that the distance between lifts of its endpoints is equal to the length of the segment in $M_{X}$. Therefore we obtain
\begin{align*}
d_{M_{Y}}(\widetilde{is_{l}\mathfrak{a}R}(x), \widetilde{is_{l}\mathfrak{a}R}(y)) \leq m_{h} d_{M_{X}} (x,y)
\end{align*}

Note that the map $\widetilde{is_{l}\mathfrak{a}R}$ may not be globally Lipschitz because in general the constants $m_{h}$ depend on the maximal dimension of simplices in $\tilde{M}^h$. \\

Now we show that a different choice of geometric cellular approximation $\widetilde{is_{l}\mathfrak{a}R}'$ is locally-Lipschitz homotopic to $\widetilde{is_{l}\mathfrak{a}R}$. Let $M$ be the inverse mapping telescope. Topologically, we have 
\begin{align*}
M \times [0,1] = (\coprod_{h\in \mathbb{N}_{0}} |\mathcal{U}_{h}| \times [2^h,2^{h+1}] \times [0,1])/_{\sim}
\end{align*}
where $\sim$ is the equivalence relation generated by the gluing maps $\phi_{h} \times \id: |\mathcal{U}_{h}| \times \{2^h\} \times [0,1]\rightarrow |\mathcal{U}_{h-1}| \times \{2^h\} \times [0,1]$. \\

The orientation on each simplex of $|\mathcal{U}_{h}|$ induces an orientation on each simplex of $|\mathcal{U}_{h}| \times [2^h,2^{h+1}]$, as well as on each simplex of $|\mathcal{U}_{h}| \times [2^h,2^{h+1}] \times [0,1]$. Let $\langle z_0,\dots, z_{m} \rangle$ be an oriented simplex in $|\mathcal{U}_{h}|$.  An oriented simplex in the product simplicial structure of $|\mathcal{U}_{h}| \times [2^h,2^{h+1}] \times [0,1]$ which lies on the gluing height $\{2^h\}$ is a face of a simplex of the form $\langle z_{0},\dots,z_{j}, w_{j}, \dots, w_{m} \rangle$, where $z's$ live on time $0$ and the $w's$ are copies of the $z's$ at time $1$. The map $\phi_{h} \times \id$ takes this simplex to $\langle  \phi_{h}(z_{0}),\dots, \phi_{h}(z_{j}), \phi_{h}(w_{j}), \dots, \phi_{h}(w_{m})\rangle $ and this is indeed a simplex in $|\mathcal{U}_{h-1}| \times \{2^h\} \times [0,1]$ since $\langle \phi_{h}(z_0), \dots, \phi_{h}(z_{m})\rangle$ is an oriented simplex in $|\mathcal{U}_{h-1}|$. Therefore, the gluing maps $\phi_{h}\times \id$ are simplicial. This defines a CW-structure on $M \times [0,1]$, where $M \times \{0\}$ and $M\times \{1\}$ are subcomplexes. We can define the open star of $v$ in $M \times [0,1]$ as 
\begin{align*}
\langle st(v)\rangle:= \Pi (\bigcup_{w \in \Pi^{-1}(v)} \langle st(w) \rangle)
\end{align*}
where $\Pi: \coprod_{h\in \mathbb{N}_{0}} \tilde{M}^h \times [0,1] \rightarrow M \times [0,1]$ is the projection map obtained from the simplicial gluing maps $\phi_{h}\times \id$, and $\langle st(w)\rangle$ denotes the star of the vertex $w$ in $\coprod_{h\in \mathbb{N}_{0}} \tilde{M}^h \times [0,1]$. We have that for any simplex $\sigma\subset |\mathcal{U}_{h}| \times [2^h,2^{h+1}] \times [0,1]$, 
\begin{align*}
\Pi \sigma \subset \bigcap_{v \in vert \Pi \sigma}|st(v)| \subset \bigcap_{v \in vert \Pi \sigma}|gst(v)|
\end{align*}
where $pr: M_{X}\times [0,1]\rightarrow M_{X}$ is the projection.\\

The argument of point $1$. of Proposition \ref{geosa} applies now. Consider the map $is_{l}\mathfrak{a}R \circ pr: M_{X} \times [0,1] \rightarrow M_{Y}$, which satisfies the geometric star condition since $pr |gst(v)| \subset |gst(pr (v))|$ for all vertices $v\in  M_{X} \times [0,1]$. Therefore, the geometric cellular approximation $(is_{l}\mathfrak{a}R \circ pr)^{\sim}$ with respect to the CW structure on $M_{X}\times [0,1]$ exists, and can be chosen to agree with $\widetilde{is_{l}\mathfrak{a}R}$ and $\widetilde{is_{l}\mathfrak{a}R}'$on $M_{X} \times \{0\}$ and $M_{X} \times \{1\}$ respectively.\\

We define: $\psi: (M_X,\omega) \rightarrow (M_Y,\widetilde{is_{l}\tau})$ as any choice of geometric cellular approximation $\widetilde{is_{l}\mathfrak{a}R}$ for our fixed $l$, which agrees with the fixed simplicial approximation $\widetilde{is_{l}\tau} = \widetilde{i\beta}$ when restricted to the base ray $\omega$. By construction $\psi$ is proper, and so induces a map $\psi: \pi_{n}^{e}(M_X,\omega)\rightarrow \pi_{n}^{e}(M_{Y},\widetilde{is_{l}\tau})$. We use the same symbol $\psi$ for the map itself and the induced homomorphism on end homotopy groups.\\

We show the commutativity of the diagram. Let $[f]\in \pi_{n}^{c}(cX,\omega)$ be a homotopy class. The composition $\theta s_{l} \mathfrak{a}[f]$ is represented by $(ib_{s_{r+l}\tau,s_{l}\tau}s_{l+r}\mathfrak{a}f)^{\sim}$ for sufficiently large $r$. \\

The composition $\psi \theta [f]$ is represented by 
\begin{align*}
(\widetilde{i  s_{l}\mathfrak{a} R})\circ  b^{J}_{\widetilde{s_{k}\omega},\omega} \widetilde{is_{k} f} = b^{\widetilde{i  s_{l}\mathfrak{a} R} \circ J}_{\widetilde{i  s_{l}\mathfrak{a} R} \widetilde{s_{k}\omega},\widetilde{i  s_{l}\mathfrak{a} R} \circ \omega}\widetilde{i  s_{l}\mathfrak{a} R} \circ \widetilde{i s_{k} f}
\end{align*}
for $k$ sufficiently large such that the diameter of $s_{k}f|st(v)|<\kappa$. $J$ is the linear homotopy between $\widetilde{s_{k}\omega}$ and $\omega$. As a reminder: the tildes on top of $is_{l}\mathfrak{a}R$ and $is_{k}f$ are referring to geometric cellular approximation and geometric simplicial approximation respectively.\\

Let us rewrite $\psi \theta[f]$ in a more pleasant form. We first look at the restriction to the subcomplex $U\subset c\overline{[-1,1]^n}$. Let $v$ be a vertex in $U$. From the definitions we have $i s_{k} f(|st(v)|) \subset \langle gst(z_v) \rangle$ for some $z_v\in vert M_{X}$. We have also $is_{l}\mathfrak{a}R |gst(z_v)| \subset \langle gst(w_v) \rangle$ for some $w_v$. Therefore 
\begin{align*}
is_{l}\mathfrak{a}Ri s_{k} f(|st(v)|) \subset is_{l}\mathfrak{a}R \langle gst(z_v)\rangle \subset  \langle gst(w_v) \rangle
\end{align*}
This shows that $is_{l}\mathfrak{a}Ri s_{k} f$ satisfies the geometric star condition and the geometric simplicial approximation $(is_{l}\mathfrak{a}Ri s_{k} f)^{\sim}$ exists. We have 
\begin{align*}
\widetilde{i s_{l}\mathfrak{a}R}\circ \widetilde{i s_{k} f} \simeq (is_{l}\mathfrak{a}R i s_{k} f)^{\sim} = (i s_{l} \mathfrak{a} s_{k} f)^{\sim}
\end{align*}
for sufficiently large $l$. There is subtlety in the first equivalence. Let $\sigma$ be a simplex. By construction $\widetilde{is_{k}f}(v) = z_{v}$ and $\widetilde{is_{l}\mathfrak{a}R}(z_{v}) = w_{v}$ for $v$ the vertices of $\sigma$. $\widetilde{is_{k}f}(\sigma)$ is contained within the geometric core $\mathscr{C}$ of the geometric simplex defined by the vertices $\cup_{v\in vert \sigma} z_{v}$. By construction, the image of $\mathscr{C}$ under $\widetilde{is_{l}\mathfrak{a}R}$ is contained within the geometric core spanned by the vertices $\widetilde{is_{l}\mathfrak{a}R}(\cup_{z\in vert \mathscr{C}} z)$, of which $\cup_{v\in vert \sigma} w_{v}$ is a subset. By the argument in point $1$ of Proposition \ref{geosa}, $\widetilde{i s_{l}\mathfrak{a}R}\circ \widetilde{i s_{k} f}$ is homotopic to $(is_{l}\mathfrak{a}R i s_{k} f)^{\sim}$ .\\

For the second equality, we recall that $Ri$ is $4$-close to the identity. Therefore it is possible to choose $l$ large enough (independent of $f$) such that  $(is_{l}\mathfrak{a}R i s_{k} f)^{\sim}=(i s_{l} \mathfrak{a} s_{k} f)^{\sim}$, since the image $s_{k}f  |st(v)|$, by construction, has diameter $<\kappa$ and lies in the $4$-neighbourhood of $R\langle gst(z_v) \rangle$. \\

Consider now $f$ restricted to the boundary $\partial{U}$. As discussed above, it is generally not true that $\widetilde{i  s_{l}a R} \circ \widetilde{g} = (is_{l}\mathfrak{a}R\circ g)^{\sim}$ for an arbitrary map $g$ since $\widetilde{g}$ only sends a simplex to a geometric simplex. However, in the specific case of $\omega$, $\widetilde{s_{k}\omega}$ sends a simplex to a simplex of $M$. Therefore the homotopy $\widetilde{is_{l}\mathfrak{a}R}\circ \widetilde{is_{k}f} \simeq (is_{l}\mathfrak{a}R is_{k}f)^{\sim}$ can be chosen to be the fixed as the base ray 
\begin{align*}
\widetilde{is_{l}\mathfrak{a}R}\circ \widetilde{is_{k}\omega} = (is_{l} \mathfrak{a} s_{k}\omega)^{\sim} = (is_{l}\tau s_{k})^{\sim}
\end{align*} on the boundary. This shows that we can replace the section represented by $\widetilde{is_{l}\mathfrak{a}R} \circ \widetilde{is_{k}f}$ with $(is_{l} \mathfrak{a} s_{k}f)^{\sim}$. \\

Consider now  $\psi \theta[f]$ restricted to $c\overline{[-1,1]^n}\setminus U = c(\partial W \times [0,1])$. Recall that change of base ray homomorphisms are independent of the $\partial W$-variable. We show that 
\begin{align*}
 b^{\widetilde{i  s_{l}\mathfrak{a} R} \circ J}_{\widetilde{i  s_{l}\mathfrak{a} R} \widetilde{s_{k}\omega},\widetilde{i  s_{l}\mathfrak{a} R} \circ \omega} \simeq   (i b^{s_{l}\mathfrak{a}H^{k}}_{s_{l}\tau s_{k},s_{l}\tau})^{\sim}
\end{align*}
relative to the boundary $c(\partial W \times \partial[0,1])$. Here, $H^{k}: c[0,1]\rightarrow cX$ is the homotopy between $s_{k}\omega$ and $\omega$.\\

First, observe that $ (i b^{s_{l}\mathfrak{a}H^{k}}_{s_{l}\tau s_{k},s_{l}\tau})^{\sim}$ can be chosen to be equal to $b^{\widetilde{is_{l}\mathfrak{a}H^{k}}}_{\widetilde{is_{l}\tau s_{k}},\widetilde{is_{l}\tau}}$. They obviously agree on vertices, and for any simplex $\sigma \in c(\partial W \times [0,1])$ we have that 
\begin{align*}
b^{\widetilde{is_{l}\mathfrak{a}H^{k}}}_{\widetilde{is_{l}\tau s_{k}},\widetilde{is_{l}\tau}}(\sigma) \subset \mathscr{C}_{\sigma}
\end{align*}
where $\mathscr{C}_{\sigma}$ is the geometric core of the simplex spanned by the images of the vertices of $\sigma$ under $b^{\widetilde{is_{l}\mathfrak{a}H^{k}}}_{\widetilde{is_{l}\tau s_{k}},\widetilde{is_{l}\tau}}$. Therefore $b^{\widetilde{is_{l}\mathfrak{a} H^{k}}}_{\widetilde{is_{l}\tau s_{k}},\widetilde{is_{l}\tau}}$ is a geometrically simplicial map such that $i b^{s_{l}\mathfrak{a}H^{k}}_{s_{l}\tau s_{k},s_{l}\tau}$ satisfies the geometric star condition with respect to the restriction $b^{\widetilde{is_{l}\mathfrak{a}H^{k}}}_{\widetilde{is_{l}\tau s_{k}},\widetilde{is_{l}\tau}}$ to vertices. Therefore it is a geometric simplicial approximation.\\
 
Observe that the map $ b^{\widetilde{i  s_{l} \mathfrak{a}R} \circ J}_{\widetilde{i  s_{l} \mathfrak{a}R} \widetilde{s_{k}\omega},\widetilde{i  s_{l}\mathfrak{a} R} \circ \omega}$ is a re-parametrisation of $\widetilde{is_{l}\tau}$. 
This is just because the linear homotopy $J$ has image in $\omega$, and $\widetilde{is_{l}\mathfrak{a}R}\circ \omega$ is just $\widetilde{is_{l}\tau}$. Consider now the homotopy $s_{l}\mathfrak{a}H^{k}$. Give $c[0,1]$ the simplicial structure $c[0,1]_{\frac{\varepsilon}{4}}$, where $\varepsilon<\frac{1}{8}$ is the constant we fixed earlier. $H^{k}$ is $1$-Lipschitz in the time variable. Consider a vertex $v=(h,ht)\in c[0,1]_{\frac{\varepsilon}{4}}$. We have that $|st(v)|\subset [h-1,h+1]\times [ht-\frac{\varepsilon}{2}, ht+\frac{\varepsilon}{2}]$. Therefore we have 
\begin{align*}
 H^{k}|st(v)|\subset N_{\frac{\varepsilon}{2}}(H^{k}_{t} [h-1,h+1])\subset N_{\frac{\varepsilon}{2}}[H^{k}_{t}(h)-1,H^{k}_{t}(h)+1]\\
\subset [[H^{k}_{t}h]-\frac{3}{2}-\frac{\varepsilon}{2},[H^{k}_{t}h]+\frac{3}{2}+\frac{\varepsilon}{2}] \subset [[H^{k}_{t}h]-2,[H^{k}_{t}h]+2]
\end{align*}
where $[H^{k}_{t}h]$ denotes the closest integer to $H^{k}_{t}h$. Recall that we have chosen $l$ large enough so that $s_{l}\tau [h-2,h+2]$ has diameter less than $2\delta <\frac{1}{16}$ for all $h\in \mathbb{N}$. Therefore, the simplicial approximation of $is_{l}\tau H^{k}$ exists by defining 
\begin{align*}
\widetilde{is_{l}\tau H^{k}}(h,ht) := \widetilde{is_{l}\tau} [H^{k}_{t}(h)]
\end{align*}
For $t=0$ we have that $[H^{k}_{0}(h)] = [s_{k}\omega(h)] = \widetilde{is_{k}\omega}(h)$. For $t=1$ we have $[H^{k}_{1}(h)] = [h] = h$. This shows that our homotopy agrees with $\widetilde{is_{l}\mathfrak{a}R}\circ J$ on the boundary.\\

We now show that $\widetilde{is_{l}\tau H^{k}}$ is homotopic to a map $G$ which is a reparametrisation of $\widetilde{is_{l}\tau}$, ie. $G = \widetilde{is_{l}\tau} \circ \overline{G}$ , where $\overline{G}: c[0,1]\rightarrow [1,\infty)$. Since $h\in  |st(h-1)|\cap |st(h)|\cap  |st(h+1)|$, we have that for any three consecutive vertices, the geometric simplex $[\widetilde{is_{l}\tau}(h-1),\widetilde{is_{l}\tau}(h), \widetilde{is_{l}\tau}(h+1)]$ exists. Consider now a simplex $\sigma\in c[0,1]_{\frac{\varepsilon}{4}}$ with $p\sigma\subset [h,h+1]$ and $(h,ht)$ a vertex of $\sigma$. Any other vertex $w=(h', h't')$ of $\sigma$ lies in the set $\{h,h+1\} \times[ht-\frac{\varepsilon}{2}, ht+\frac{\varepsilon}{2}]$. We obtain that $H^k_{t'}(h')\in [H^k_{t}(h)-\frac{\varepsilon}{2}, H^k_{t}(h)+1+\frac{\varepsilon}{2}]$ and hence

\[[H^k_{t'}(h')] \in \begin{cases} 
        [[H^{k}_{t}(h)]-1, [H^{k}_{t}(h)]+1]  & [H^{k}_{t}(h)]\geq H^{k}_{t}(h)\\ 
        [[H^{k}_{t}(h)], [H^{k}_{t}(h)]+2]  &  [H^{k}_{t}(h)]\leq H^{k}_{t}(h) \\
       \end{cases}
    \]

We define the map $\overline{G}: c[0,1]_{\frac{\varepsilon}{4}}\rightarrow [1,\infty)$ by letting $\overline{G}(h,ht) = [H^{k}_{t}(h)]$ on vertices and extending linearly in the interior of simplies. $\overline{G}(\sigma)$ has image in $I_{\sigma}$, the interval of length $\leq 2$ spanned by the vertices $\cup_{v\in \sigma} \overline{G}(v)$. The map $G:= \widetilde{is_{l}\tau}\circ \overline{G}$ is homotopic to $\widetilde{is_{l}\tau H^{k}}$ because for every simplex $\sigma$, we have 
\begin{align*}
\widetilde{is_{l}\tau H^{k}}(\sigma) &\subset \mathscr{C}_{[\cup_{v\in vert \sigma} \widetilde{is_{l}\tau} \overline{G}(v)]} \subset  \mathscr{C}_{[\cup_{v\in vert I_{\sigma}} \widetilde{is_{l}\tau}(v)]} \\
G(\sigma)&\subset \widetilde{is_{l}\tau} (I_{\sigma}) \subset \mathscr{C}_{[\cup_{v\in vert I_{\sigma}} \widetilde{is_{l}\tau}(v)]} 
\end{align*} 
Additionally, since $G$ and $\widetilde{is_{l}\tau H^{k}}$ agree on the boundary $c\partial[0,1]$, this homotopy can be made relative to $c\partial [0,1]$. \\

Any two reparametrisations of $\widetilde{is_{l}\tau}$ are homotopic via a linear homotopy. Therefore, $G$ is homotopic rel $c[0,1]$ to $\widetilde{i  s_{l}\mathfrak{a} R} \circ J$. Putting it together, this means we can replace $b^{\widetilde{i  s_{l}\mathfrak{a} R} \circ J}_{\widetilde{i  s_{l} \mathfrak{a} R} \widetilde{s_{k}\omega},\widetilde{i  s_{l}\mathfrak{a} R} \circ \omega}$ with $(ib^{s_{l}\mathfrak{a}H^{k}}_{s_{l}\tau s_{k},s_{l}\tau})^{\sim} = (b^{is_{l}\mathfrak{a}H^{k}}_{is_{l}\tau s_{k},is_{l}\tau})^{\sim} $. In summary, we have shown that 
\begin{align*}
\psi \theta [f] =[(i b^{s_{l}\mathfrak{a}H^{k}}_{s_{l}\tau s_{k},s_{l}\tau} s_{l}\mathfrak{a} s_{k} f)^{\sim}]
\end{align*}
with respect to the simplicial structure $c\overline{[-1,1]^n}$.\\

Now we construct a homotopy from  $(i b_{s_{l}\tau s_{k},s_{l}\tau} s_{l}\mathfrak{a} s_{k} f)^{\sim}$ to $(ib_{s_{r+l}\tau,s_{l}\tau}s_{l+r}\mathfrak{a}f)^{\sim}$ in two steps:
\begin{align*}
(i b^{s_{l}\mathfrak{a}H^{k}}_{s_{l}\tau s_{k},s_{l}\tau} s_{l}\mathfrak{a} s_{k} f)^{\sim} \overset{\widetilde{iA}}{\simeq} (i b^{s_{l+r}\mathfrak{a}H^{k}}_{s_{l+r}\tau s_{k},s_{l+r}\tau} s_{l+r}\mathfrak{a} s_{k} f)^{\sim} \overset{B}{\simeq} b^{\id}_{\widetilde{is_{l+r}\tau},\widetilde{is_{l+r}\tau}}(i s_{l+r}\mathfrak{a}f)^{\sim}
\end{align*}
$A: c([-1,1]^{n} \times [0,1])\rightarrow cX$ is the homotopy $\mathbf{H}^{r} b_{s_{l}\tau s_{k},s_{l}\tau} s_{l}\mathfrak{a}s_{k}f$, where $\mathbf{H}^{r}: I_{p}(cY)\rightarrow cY$ is the homotopy between $\id$ and $s_{r}$. $B$ restricted to $U$ is the homotopy $(is_{l+r}\mathfrak{a} (\mathbf{H}^{k})^{-1}f)^{\sim}$ where $\mathbf{H}^{k}:I_{p}(cX)\rightarrow cX$ is the homotopy between $\id$ and $s_{k}$. \\

\textbf{Step 1} $(i b^{s_{l}\mathfrak{a}H^{k}}_{s_{l}\tau s_{k},s_{l}\tau} s_{l}\mathfrak{a} s_{k} f)^{\sim} \overset{\widetilde{iA}}{\simeq} (i b^{s_{l+r}\mathfrak{a}H^{k}}_{s_{l+r}\tau s_{k},s_{l+r}\tau} s_{l+r}\mathfrak{a} s_{k} f)^{\sim}$: \\

The homotopy $\widetilde{iA}$ is the same as the one constructed in the proof of Proposition \ref{isogeneralray}, when we show in point $2.$ that $\theta$ is independent of the choice of shrinking map. As before, we let $F= b^{s_{l}\mathfrak{a}H^{k}}_{s_{l}\tau s_{k},s_{l}\tau} s_{l}\mathfrak{a} s_{k} f$ and $L_{F}$ be the large-scale Lipschitz constant for $pF$. For a $\varepsilon_{F}<\frac{1}{(2L_{F}+1)16}$, the simplicial approximation of $iA$ with respect to $c(\overline{[-1,1]^n} \times [0,1]_{\varepsilon_{F}})$ exists; by construction, for any vertex $(hx,ht,h)\in c(\overline{[-1,1]^n} \times [0,1]_{\varepsilon_{F}})$, $F|st(hx,ht,h)|$ lies in 
\begin{align*}
N_{2(2L_{F}+1)\varepsilon_{F} +\frac{\varepsilon}{2}} (U)
\end{align*} 
where $U$ is a set of diameter $<\delta$. Observe that $\widetilde{iA}$ restricted to the boundary $c(\partial\overline{[-1,1]^n} \times [0,1]_{\varepsilon_{F}})$ is a simplicial approximation of the homotopy $\mathbf{H}^{r} s_{l}\tau$, ie. the inverse homotopy to the one used for the change of base ray $(ib_{s_{l+r}\tau,s_{l}\tau})^{\sim}$\\

\textbf{Step 2 $ (i b^{s_{l+r}\mathfrak{a}H^{k}}_{s_{l+r}\tau s_{k},s_{l+r}\tau} s_{l+r}\mathfrak{a} s_{k} f)^{\sim} \overset{B}{\simeq} b^{\id}_{\widetilde{is_{l+r}\tau},\widetilde{is_{l+r}\tau}}(i s_{l+r}\mathfrak{a}f)^{\sim}$}: \\

We first consider the homotopy $B$ restricted to $U$, which we equip with the simplicial structure $c([-1,1]^n \times [0,1]_{\varepsilon})$. Since $s_{l}\mathfrak{a} (\mathbf{H}^{t})^{-1}f$ is a coarse homotopy, there exists a $r$ large enough such that $is_{l+r}\mathfrak{a} (\mathbf{H}^{t})^{-1}f$ has a simplicial approximation with respect to $c([-1,1]^n \times [0,1]_{\varepsilon})$. We let $B_{|U}$ be defined as $(is_{l+r} \mathfrak{a} (\mathbf{H}^{t})^{-1}f)^{\sim}$.\\

There is an obvious identification between the boundary $c(\partial[-1,1]^n \times [0,1]_{\varepsilon})$ and the subcomplex $c\overline{[-1,1]^n}\setminus U = c(\partial W \times [0,1]_{\varepsilon})$. Observe that $B$ restricted to the boundary $c(\partial [-1,1]^n \times [0,1]_{\varepsilon})$ is a simplicial approximation of the change of base ray homormophism $i b^{s_{l+r}\mathfrak{a}H^{k}}_{s_{l+r}\tau s_{k},s_{l+r}\tau}$ under this identification, and by going through the constructions carefully (we define the simplicial approximation by mapping a vertex to the closest vertex to its image under $i b^{s_{l+r}\mathfrak{a}H^{k}}_{s_{l+r}\tau s_{k},s_{l+r}\tau}$), one sees that we can actually choose the simplicial approximations to be the same. Therefore, we can "fill in" the set $c(\overline{[-1,1]^n} \times [0,1]_{\varepsilon})\setminus c(W\times [0,1]_{\varepsilon})$ with a homotopy between $i b^{s_{l+r}\mathfrak{a}H^{k}}_{s_{l+r}\tau s_{k},s_{l+r}\tau}$ and $b^{\id}_{\widetilde{is_{l+r} \tau}, \widetilde{is_{l+r}\tau}}$. We therefore obtain the desired homotopy $B$. Observe that $B$ restricted to $c(\partial \overline{[-1,1]^n} \times [0,1])$ is the fixed base ray $\widetilde{is_{l+r} \tau}$. \\

Now we stack together the homotopies $\widetilde{iA}$ and $B$ to obtain a homotopy 
\begin{align*}
\widetilde{iA} \cup B: c(\overline{[-1,1]^n} \times [0,1]_{\varepsilon_F}) \cup c(\overline{[-1,1]^n} \times [0,1]_{\varepsilon})\rightarrow M_{Y}
\end{align*}
For an image of this, see Figure \ref{fig:natural1}.\\

The set $X = c(\overline{[-1,1]^n} \times \{1\})\cup c(\partial\overline{[-1,1]^n} \times [0,1]_{\varepsilon}) \cup c(\partial \overline{[-1,1]^n} \times [0,1]_{\varepsilon_{F}})$, where 

\begin{align*}
 c(\overline{[-1,1]^n} \times \{1\})\cup c(\partial\overline{[-1,1]^n}\times [0,1]_{\varepsilon}) \subset c(\overline{[-1,1]^n} \times [0,1]_{\varepsilon})\\
c(\partial \overline{[-1,1]^n} \times [0,1]_{\varepsilon_{F}}) \subset c(\overline{[-1,1]^n} \times [0,1]_{\varepsilon_F})
\end{align*} can be identified with $c[-1,1]^n$ in an obvious way: identify $c(W \times \{1\})$ with the set $U$, and rescale the remaining three components to squeeze them into $c[-1,1]^n \setminus U$. Under this identification, we see that $\widetilde{iA}\cup B$ restricted to $X$ can be identified with a rescaled version of 
\begin{align*}
(ib_{s_{l+r}\tau,s_{l}\tau})^{\sim}\circ b^{\id}_{\widetilde{is_{l+r}\tau}, \widetilde{is_{l+r}\tau}}\circ b^{\id}_{\widetilde{is_{l+r}\tau}, \widetilde{is_{l+r}\tau}} (is_{l+r}\mathfrak{a}f)^{\sim} \simeq (ib_{s_{l+r}\tau, s_{l}\tau} s_{l+r} \mathfrak{a} f)^{\sim}
\end{align*}

\begin{figure}[H]
\centering
  \centering
  \includegraphics[width=.9\linewidth]{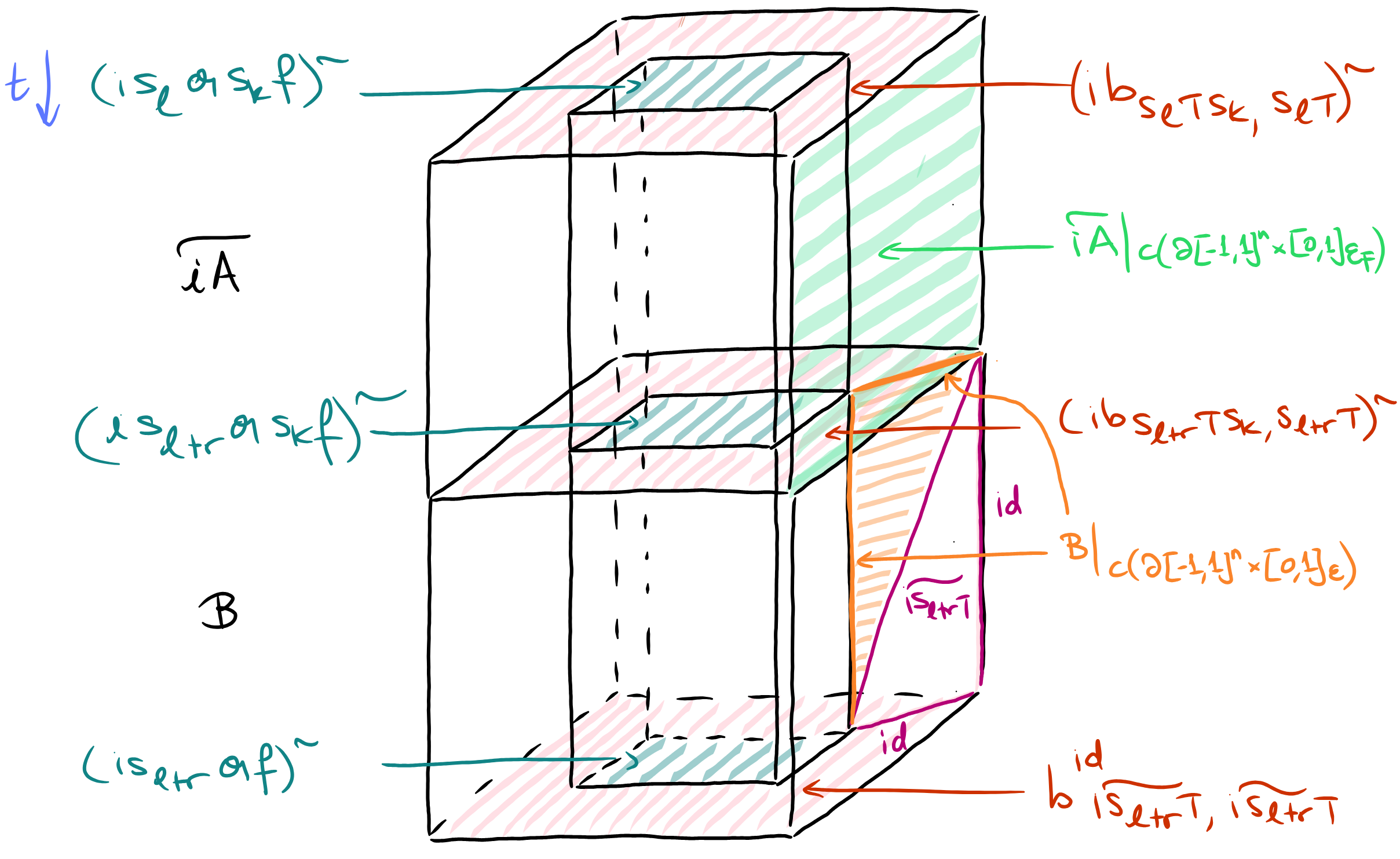}
  \caption{The homotopy $\widetilde{iA}\cup B$. Note that the time is increasing as one goes down the picture.}
  \label{fig:natural1}
\end{figure}

Therefore, we obtain the homotopy 
\begin{align*}
(i b_{s_{l}\tau s_{k},s_{l}\tau} s_{l}\mathfrak{a} s_{k} f)^{\sim} \simeq (ib_{s_{l+r}\tau, s_{l}\tau} s_{l+r} \mathfrak{a} f)^{\sim}
\end{align*}
which shows that $\theta \mathfrak{a}[f]=\psi \theta [f]$. In conclusion, the diagram commutes. The map $\theta \circ s_{l}:\pi_{n}^{c}(cY,\tau)\rightarrow \pi_{n}^{L,e}(M_{Y},\widetilde{is_{l}\tau})$ is an isomorphism by Proposition \ref{isogeneralray}, and 
\begin{align*}
\mathfrak{d}: \pi_{n}^{L,e}(M_{Y}, \widetilde{is_{l}\tau}) \rightarrow \pi_{n}^{e}(M_{Y},\widetilde{is_{l}\tau})
\end{align*}
is also an isomorphism because

\begin{prop}
Let $\beta:[1,\infty)\rightarrow M$ be a simplicial base ray. There is an isomorphism $\pi_{n}^{L,e}(M,\beta)\cong \pi_{n}^{e}(M,\beta)$
\end{prop}
 
\begin{proof} This proof is almost identical to that of Proposition \ref{natiso}, using the same symbols and maps. We record below the small modification required: \\

Let $A$ be the Lipschitz constant for $\beta$. The homotopy 
\begin{align*}
G:I_{p}([1,\infty))&\rightarrow M\\
(h,ht)&\mapsto \beta(th + (1-t)\omega'(h))
\end{align*}
between $\beta \omega'$ and $\beta$ is $2A$-Lipschitz:
\begin{align*}
d_{M} (G(h,ht), G(s,sl)) &\leq A |th+(1-t)\omega'(h) - ls - (1-l)\omega'(s)|  \\
&\leq A(2|ht-sl| + 2|h-s|)
\end{align*}
Therefore $b_{\beta\omega',\beta}[f'] \in \pi_{n}^{L,e}(M,\beta)$, and $\mathfrak{d} b_{\beta\omega',\beta}[f'] = [f]$, which gives us surjectivity of $\mathfrak{d}$. \\

There are no modifications to the rest of the proof. \\
\end{proof}

Putting it all together, the diagram commutes for $n\geq 1$. For $n=0$, let $[\omega'] \in \pi_{0}^{c}(cX,\omega)$. By postcomposing with a shrinking map, we can assume that $\widetilde{i\omega'}$ exists. We have 
\begin{align*}
\psi\theta[\omega'] = [\widetilde{is_{l}\mathfrak{a}R}\widetilde{i\omega'}]= [(is_{l}\mathfrak{a}Ri\omega')^{\sim}] = [(is_{l}\mathfrak{a}\omega')^{\sim}] = \theta s_{l}\mathfrak{a} [\omega']
\end{align*}
which gives us commutativity. \\

Therefore $\psi: \pi_{n}^{e}(M_{X},\omega)\rightarrow \pi_{n}^{e}(M_{Y}, \widetilde{is_{l}\tau})$ is an isomorphism when $\mathfrak{a}$ induces an isomorphism. 
\end{proof}

To summarise, any coarse map $\mathfrak{a}: (cX,\omega)\rightarrow (cY,\tau)$ induces a proper, locally Lipschitz and cellular map $\psi: (M_{X},\omega)\rightarrow (M_{Y}, \widetilde{is_{l}\tau})$ via geometric cellular approximation. If $\mathfrak{a}$ induces an isomorphism on coarse homotopy groups then $\psi$ induces an isomorphism on end homotopy groups. \\

Now we discuss finite dimensionality. Since $X$ and $Y$ have finite shape dimension, by Proposition \ref{finitedim} there exist finite-dimensional shape expansions $\mathcal{X}',\mathcal{Y}'$ and cellular, proper homotopy equivalences $\zeta_{\mathcal{X'}\mathcal{X}}: M_{\mathcal{X}'}\rightarrow M_{X}$ and $\zeta_{\mathcal{Y}\mathcal{Y}'}: M_{Y}\rightarrow M_{\mathcal{Y}'}$. 

\begin{lemma} \label{weakimpliesphe} Let $X$ and $Y$ be compact, connected metric spaces which have finite shape dimension. Let $\psi: M_{X}\rightarrow M_{Y}$ be a proper, cellular map which induces an isomorphism on end homotopy groups. Then there exists a proper homotopy inverse $\xi: M_{Y}\rightarrow M_{X}$. 
\end{lemma}

\begin{proof} The composition
\begin{align*}
M_{\mathcal{X}'}\xrightarrow{\zeta_{\mathcal{X'}\mathcal{X}}} M_{X}\xrightarrow{\psi} M_{Y} \xrightarrow{\zeta_{\mathcal{Y}\mathcal{Y'}}} M_{\mathcal{Y}'}
\end{align*}
is a proper, cellular map between contractible, finite-dimensional, strongly locally finite, path-connected CW complexes, which induces an isomorphism on all end homotopy groups and the set of ends. Therefore Corollary \ref{inverse2} applies, and there exists a proper homotopy inverse $\overline{\xi}:M_{\mathcal{Y}'}\rightarrow M_{\mathcal{X'}}$. We define $\xi = \zeta_{\mathcal{X}' \mathcal{X}} \circ \overline{\xi} \circ \zeta_{\mathcal{Y} \mathcal{Y'}}$. 
\begin{align*}
\xi \psi = \zeta_{\mathcal{X}' \mathcal{X}}  \overline{\xi} \zeta_{\mathcal{Y} \mathcal{Y'}} \psi \simeq \zeta_{\mathcal{X}' \mathcal{X}} ( \overline{\xi}  \zeta_{\mathcal{Y} \mathcal{Y'}} \psi \zeta_{\mathcal{X'}\mathcal{X}}) \zeta_{\mathcal{X}\mathcal{X'}} \simeq \zeta_{\mathcal{X}' \mathcal{X}} \id_{M_{\mathcal{X}'}} \zeta_{\mathcal{X}\mathcal{X}'} \simeq \id_{M_{X}}\\
\psi \xi =  \psi \zeta_{\mathcal{X}' \mathcal{X}}  \overline{\xi}  \zeta_{\mathcal{Y} \mathcal{Y'}} \simeq \zeta_{\mathcal{Y}' \mathcal{Y}}(\zeta_{\mathcal{Y}\mathcal{Y}'} \psi \zeta_{\mathcal {X}' \mathcal{X}} \overline{\xi}) \zeta_{\mathcal{Y} \mathcal{Y}'} \simeq \zeta_{\mathcal{Y}'\mathcal{Y}} \id_{M_{\mathcal{Y}'}} \zeta_{\mathcal{Y}\mathcal{Y}'} \simeq \id_{M_{Y}}
\end{align*}
where the $\simeq$ denote proper homotopy. Hence, $\xi$ is a proper homotopy inverse to $\psi$. 
\end{proof}

\subsection{Constructing a coarse inverse}

The goal of this subsection is to prove:

\begin{theorem} \label{unbasedwhitehead} Let $(X,x_0)$ and $(Y,y_0)$ be connected, compact metric spaces with finite shape dimension. Let $\omega$ be the standard parametrisation of the base ray $c\{x_0\}$. Let $\mathfrak{a}:(cX,\omega)\rightarrow (cY,\mathfrak{a}\omega)$ be a coarse map, which induces an isomorphism on coarse homotopy groups $\mathfrak{a}: \pi_{n}^{c}(cX,\omega)\rightarrow \pi_{n}^{c}(cY,\mathfrak{a}\omega)$ for all $n\geq 0$. Then $\mathfrak{a}$ is a coarse homotopy equivalence, ie. there exists $\mathfrak{b}: cY\rightarrow cX$ such that $\mathfrak{b}\circ \mathfrak{a} \simeq id_{cX}$ and $\mathfrak{a}\circ \mathfrak{b} \simeq id_{cY}$.
\end{theorem} 

We now take the proper homotopy inverse $\xi: M_{Y}\rightarrow M_{X}$ and use it to construct a coarse homotopy inverse to $\mathfrak{a}: cX\rightarrow cY$. The technical issue to overcome is that $i: cY\rightarrow M_{Y}$ is not coarse, so we cannot define the coarse inverse to $\mathfrak{a}$ as simply $R\circ \xi \circ i$. The trick is to precompose with a shrinking map $s: cY\rightarrow cY$. \\

Let $X$ be a metric space and $\mathcal{U}$ a locally finite cover of $X$. We choose a specific partition of unity subordinate to $\mathcal{U}$ as follows. Let $U\in \mathcal{U}$. We define $\varphi_{U}: U\rightarrow [0,1]$ by the formula
\begin{align*}
\varphi_{U}(x) = \frac{d(x,X\setminus U)}{\sum_{V\in \mathcal{U}} d(x,X\setminus V)}
\end{align*}

Recall that the map $\varphi: X\rightarrow |\mathcal{U}|$ induced by the partition of unity $\{\varphi_{U}\}_{U\in \mathcal{U}}$ is defined as $\varphi(x) = \sum_{U\in \mathcal{U}} \varphi_{U}(x)[U]$. 

\begin{lemma} (Lemma $4.3.5$ in \cite{nowak2023large}) Let $\mathcal{U}$ be a uniformly bounded cover of $X$ with multiplicity $\kappa$ and Lebesgue number $2r$. Then for each $U\in \mathcal{U}$ the function $\varphi_{U}$ is $(\frac{1}{r} + \frac{2\kappa}{r^2})$-Lipschitz. 
\end{lemma}

For $x,y\in X$ we have that $d_{std}(\varphi(x),\varphi(y)) = \sum_{U\in \mathcal{U}} |\varphi_{U}(x)-\varphi_{U}(y)| \leq 2\kappa (\frac{1}{r} + \frac{2\kappa}{r^2}) d(x,y)$. Therfore, $\varphi$ is $2\kappa (\frac{1}{r} + \frac{2\kappa}{r^2})$-Lipschitz with respect to the standard metric on $|\mathcal{U}|$. Changing from the standard metric to the spherical metric $d^s$ has a Lipschitz constant dependent on $\kappa = \dim |\mathcal{U}|+1$.\\

Let $X$ be a compact metric space. The open covers $\mathcal{U}_{h}$ are finite, since they are constructed from a $\frac{1}{2^{h+1}}$-separated set. Therefore each $i_{h}: X\rightarrow (|\mathcal{U}_{h}|,d^s)$ has a Lipschitz constant, say $\alpha_{h}$, which only depends on the multiplicity of the cover and Lebesgue number $\frac{1}{2^h}$ of  $\mathcal{U}_{h}$. \\

\begin{lemma} The map $i$ can be made coarse, by precomposing with a shrinking map $s$ on $cX$. \label{shrinkcoarse}
\end{lemma}

\begin{proof} We know that the map 
\begin{align*}
\varphi_{h}: X\rightarrow |\mathcal{U}_{h}|
\end{align*}
is $\alpha_{h}$-Lipschitz. Choose constants $m_h$ such that $2^{m_h} > \alpha_{h}$ and $m_h \leq m_{h+1}$ for all $h\in \mathbb{N}_{0}$.  We can assume that $m_{0}\geq 1$ for later purposes. Consider now the shrinking map $s$, which is defined piecewise as follows:
\begin{align*}
s^{h}: X\times 2^{m_{h}} (2^h,2^{h+1}]&\rightarrow X \times (2^{h},2^{h+1}]\\
(tx,t)&\mapsto (\frac{t}{2^{m_{h}}}x, \frac{t}{2^{m_{h}}})
\end{align*}
and we send points with heights in the interval $(2^{h+m_{h-1}}, 2^{h+m_{h}}]$ to height $2^{h}$ and all points below height $2^{m_0}$ to height $1$. This is a coarse shrinking map, and it is coarsely homotopic to the identity by a linear homotopy. \\

We now check that the composition $i \circ s$ is coarse. Let us first restrict to $X\times 2^{m_{h}} (2^h,2^{h+1}]$. For two points $(tx,t), (ly,l)$ in  $X\times 2^{m_{h}} (2^h,2^{h+1}]$ with distance  $\leq Q$ we have that $\|x-y\|\leq \frac{Q}{2^{h+m_{h}}}$ and $|t-l|\leq Q$. 
\begin{align*}
 d_{M}(i  s^{h}(tx,t), i s^{h}(ly,l)) &= d_{M} ( (\varphi_{h}(x),\frac{t}{2^{m_{h}}}), (\varphi_{h}(y), \frac{l}{2^{m_{h}}}))\\
& \leq d_{|\mathcal{U}_{h}|}(\varphi_{h}(x),\varphi_{h}(y)) + \frac{1}{2^{m_{h}}}|t-l|\\
 &\leq 2^{m_{h}} \|x-y\| +  \frac{1}{2^{m_{h}}}|t-l|< \frac{Q}{2^h} + \frac{Q}{2} < 2Q
\end{align*}
Therefore $i \circ s$ is $2$-Lipschitz when restricted to any $X\times 2^{m_{h}}(2^h,2^{h+1}]$ for any $h\in \mathbb{N}_{0}$ (actually $1$-Lipschitz for $h\geq 1$). For two points with heights in the interval $(2^{h+m_{h-1}}, 2^{h+m_{h}}]$ with distance $\leq Q$ we just get
\begin{align*}
 d_{M}(i s(tx,t), i s(ly,l)) &= d_{M} ( (\varphi_{h-1}(x), 2^h), (\varphi_{h-1}(y), 2^h))\\
 &\leq 2^{m_{h-1}}\|x-y\| < \frac{Q}{2^h} <Q
\end{align*}
The map $i \circ s$ is not continuous at heights $2^{h+m_{h}}$ but since $\phi_{h} \varphi_{h}(x)$ and $\varphi_{h-1}(x)$ lie in a common simplex for any $x$, the discontinuity is uniformly bounded by $1$, the diameter of a simplex in $|\mathcal{U}_{h-1}|$. \\

To see that $i\circ s$ is globally coarse, let $Q>0$. Choose $h\geq 1$ large enough such that $Q<2^{h-1}$, ie. so that a set of diameter $Q$ above height $2^{h-1}$ intersects at most one gluing component.  Consider two points $(tx,t),(ly,l)$ with $t\geq l$ and $\|tx-ly\| + |t-l| < Q$. Observe that for any $l'\in [l,t]$ we have 
\begin{align*}
\|tx-l'y\| + |t-l'|&\leq \|tx-ly\| + |l'-l| + |t-l'| = \|tx-ly\| + |t-l| <Q\\
\|l'y - ly\| + |l'-l| &< 2 |l'-l|< 2Q
\end{align*} If both $t,l > 2^h$ then by possibly choosing $l'$ at a height where $i\circ s$ is discontinuous, we get 
\begin{align*}
d_{M}(is(tx,t), is(ly,l)) &\leq d_{M}(is(tx,t), is(l'y,l'))+ d_{M}(is(l'y,l'), is(ly,l)) + 1 \\
&< Q + 2Q +1 = 3Q+1
\end{align*} 
where the two images of $is(l'y,l')$ are different, to reflect the discontinuity. If only $l\leq 2^h$ then $l>2^{h-1}$ and we can choose $l'= 2^h$.  $i\circ s$ sends $X \times [1,2^h]$ to a bounded set in $M$, say of diameter $<T_{Q}$. Therefore, 
\begin{align*}
d_{M}(is(tx,t), is(ly,l)) &\leq d_{M}(is(tx,t), is(l'y,l'))+ d_{M}(is(l'y,l'), is(ly,l)) + 1 \\
&< Q + T_{Q} +1 
\end{align*} 
Obviously if both $t,l \leq 2^h$ then $d_{M}(is(tx,t), is(ly,l))< T_{Q}$. Therefore
\begin{align*}
d_{cX}((tx,t),(ly,l)) < Q \implies d_{M}(is(tx,t), is(ly,l)) < \max \{3Q+1, Q + T_{Q} +1 \}
\end{align*}
This shows that the map $i \circ s$ is controlled. It is obviously proper.  
\end{proof}

\begin{lemma} \label{chi} $i \circ s$ is a coarse homotopy inverse to $R$. 
\end{lemma}

\begin{proof} We simply compute both compositions. \\

$R \circ i \circ s$: Let $G: c(X\times [0,1]) \rightarrow cX$ be the coarse homotopy with $G_{c(X\times \{0\})} = \id_{cX}$ and $G_{c(X \times \{1\})} = s$. We let $H: c(X\times [0,1])\rightarrow cX$ be defined as $R \circ i \circ G$. Since $R \circ i$ is coarse, $H$ is a coarse homotopy between $R \circ i \circ s$ and $R \circ i$. $R\circ i$ is close and therefore coarsely homotopic to $\id_{cX}$. \\

$i \circ s \circ R$: 
Let $\omega'$ be the linear extension of the function $2^{h+m_{h}}\mapsto 2^{h}, 2^{h+1+m_{h}}\mapsto 2^{h+1}$ to $[1, \infty)\rightarrow [1,\infty)$, so that $s$ is the shrinking map associated to $\omega'$. We show that $i \circ s \circ R$ is close to $\mathfrak{q}\circ i \circ R$ where $\mathfrak{q}$ is a deformation retract on the mapping telescope that takes the height slice $t$ to $\omega'(t)$, for example:

\begin{align*}
\mathfrak{q}: |\mathcal{U}_{h+m_{h}}|\times 2^{m_{h}} (2^h,2^{h+1}]&\rightarrow  |\mathcal{U}_{h}| \times (2^{h},2^{h+1}]\\
(x,t)&\mapsto (\phi_{h+m_{h},h}(x), \omega'(t))
\end{align*}

Let $t\in 2^{m_h}(2^h,2^{h+1}]$. 
\begin{align*}
i s R (x,t) &= i s (t R_{h+m_{h}}(x), t) = i (\omega'(t) R_{h+m_{h}}(x), \omega'(t)) = (i_{h} R_{h+m_{h}}(x), \omega'(t))\\
\mathfrak{q}  i  R (x,t) &= \mathfrak{q} i (tR_{h+m_{h}}(x), t) = \mathfrak{q}(i_{h+m_{h}} R_{h+m_{h}}(x),t) = (\phi_{h+m_h, h}i_{h+m_{h}} R_{h+m_{h}}(x), \omega'(t))
\end{align*}

The points $i_{h}R_{h+m_{h}}(x)$ and $\phi_{h+m_h, h}i_{h+m_{h}}R_{h+m_{h}}(x)$ lie in a common simplex, so differ by a distance $\leq 1$. For a height $t\in(2^{h+m_{h-1}}, 2^{h+m_{h}}]$ we have 
\begin{align*}
i s R (x,t) &= i s (t R_{c(t)}(x), t) = i (2^h R_{c(t)}(x),2^h) = (i_{h-1} R_{c(t)}(x), 2^h)\\
\mathfrak{q}  i  R (x,t) &= \mathfrak{q} i (tR_{c(t)}(x), t) = \mathfrak{q}(i_{c(t)} R_{c(t)}(x),t) = (\phi_{c(t), h-1}i_{c(t)} R_{c(t)}(x), 2^h)
\end{align*}
where $c(t)$ is the largest number such that $2^{c(t)}<t$. As before, the images differ by at most the diameter of a simplex. Therefore $i \circ s\circ R$ is close to $\mathfrak{q} \circ i\circ R$. Since $i\circ R$ is close to the identity and $\mathfrak{q}$ is Lipschitz, $\mathfrak{q}\circ i\circ R$ is close to $\mathfrak{q}$. $\mathfrak{q}$ is coarse-Lipschitz homotopic to the identity, so we are done. \\
\end{proof}

Before we proceed with the technical details, we outline the rest of the proof in steps. 

\begin{enumerate}
\item Find a locally Lipschitz representative $\xi'$ of $\xi$ in its proper homotopy class, such that $\psi \xi'$ and $\xi' \psi$ are locally Lipschitz homotopic to $\id_{M_{Y}}$ and $\id_{M_{X}}$ respectively via homotopies $F$ and $G$. (Lemmas \ref{xi'} and \ref{liphomotopies})
\item Define coarse shrinking maps $S_{\xi'}, S_{\psi\xi'}, S_{\xi' \psi}$ on $M_{Y},I_{p}(M_{Y})$ and $I_{p}(M_{X})$ respectively such that $\xi' S_{\xi'}, FS_{\psi\xi'}$ and $GS_{\xi' \psi}$ are coarse. (Lemma \ref{globallycoarse})
\item Prove that $\xi' S_{\xi'}$ is a coarse homotopy inverse to $is_{l}\mathfrak{a}R$. (Lemma \ref{che})
\item Show that $\mathfrak{b}:= R\circ \xi' S_{\xi'} \circ is$ is a coarse homotopy inverse to $\mathfrak{a}$. (Lemma \ref{defb})
\end{enumerate}

\begin{lemma} \label{xi'} $\xi$ has a locally Lipschitz representative $\xi'$ in its proper homotopy class.
\end{lemma}

\begin{proof} This proof combines the ideas of Proposition \ref{natiso} and Proposition \ref{psiexistence}. Recall that $M_{Y} = \pi(\coprod_{h\in \mathbb{N}_{0}} |\mathcal{U}_{h}| \times [2^h,2^{h+1}])=\pi(\coprod_{h\in \mathbb{N}_{0}}\tilde{M}^h)$. Fix the CW structure on $M_{Y}$ as defined in Proposition \ref{psiexistence}. We inductively define a simplicial map on $\coprod_{h\in \mathbb{N}_{0}}\tilde{M}^h$ with respect to a progressively subdivided simplicial structure, which is locally Lipschitz with respect to $\{\tilde{M}^h\}_{h\in \mathbb{N}_{0}}$ and decends to the quotient.\\

Since $\tilde{M}^0$ is a finite simplicial complex, there exists a barycentric subdivision of $\tilde{M}^0$ such that $\xi_{0}:=\xi \circ \pi_{|\tilde{M}^0}$ satisfies the geometric star condition. Therefore there is a simplicial map $\overline{\xi}'_{0}$, and a homotopy $\overline{H}_{0}: \tilde{M}^0 \times [0,1] \rightarrow M_{X}$ between $\xi_{0}$ and $\overline{\xi}'_{0}$. Extend this homotopy to a homotopy $H_{0}: \tilde{M}_{Y} \times [0,1] \rightarrow M_{X}$ from $\xi \circ \pi$ to a map $\xi'_{0}$ such that $H_{0}$ descends to $M_{Y} \times [0,1]$ and is constant above height $3$.\\

Assume for induction that $H_{0}\ast H_{1}\ast \dots \ast H_{k}$ has already been defined, descends to the quotient, and the endpoint of the homotopy $\xi'_{k}: \tilde{M}_{Y} \rightarrow M_{X}$ is simplicial when restricted to $\coprod_{h\leq k} \tilde{M}^{h}$. Consider now $\tilde{M}^{k+1} = |\mathcal{U}_{k+1}| \times [2^{k+1},2^{k+2}]$. We define $\xi'_{k+1}$ to be $\xi'_{k}\circ \phi_{k+1}$ when restricted to the subcomplex $|\mathcal{U}_{k+1}| \times \{2^{k+1}\}$. $\phi_{k+1}$ is simplicial with respect to iterated barycentric subdivision, so $\xi'_{k}\circ \phi_{k+1}$ is simplicial with respect to the subdivided structure on $\tilde{M}^{k+1}$. By geometric simplicial approximation relative to the subcomplex $|\mathcal{U}_{k+1}| \times \{2^{k+1}\}$ (Theorem \ref{relgeosimpl}), we have a homotopy $\overline{H}_{k+1}: \tilde{M}^{k+1} \times [0,1]\rightarrow M_{X}$ and a simplicial $\overline{\xi}_{k+1}'$. As before, we extend $\overline{H}_{k+1}$ to a homotopy $H_{k+1}: \tilde{M}_{Y}\rightarrow M_{X}$ which descends to the quotient and is constant for heights $\geq 2^{k+2} +1$ and heights $\leq 2^{k+1}$. Call the endpoint of this homotopy ${\xi}_{k+1}'$.\\

The infinite concatenation of homotopies $H:= H_{0}\ast \dots \ast H_{k}\ast H_{k+1}\ast \dots$ is well-defined since any point in $\tilde{M}_{Y}$ is affected by at most $2$ of the $H_{k}$. Call $\xi'$ the endpoint of the homotopy $H$. By construction, it is locally Lipschitz with respect to the filtration $\{\tilde{M}^h\}_{h\in \mathbb{N}_{0}}$ of $\tilde{M}_{Y}$.  Since the metric on $M_{Y}$ is induced from paths in $\tilde{M}_{Y}$, we obtain a locally Lipschitz map $\xi':M_{Y}\rightarrow M_{X}$ with respect to the filtration $\{M_{[1,2^h]}\}_{h\in \mathbb{N}_{0}}$ after descending to the quotient. $\xi'$ is homotopic to $\xi$. (Note that we use the same symbols $\xi',H$ to denote the maps and homotopies on the quotient.) \\

The homotopy $H$ is proper: by the properness of $\xi$ for every $h\in \mathbb{N}$ there exists a $l_{h}$ such that $\xi(M_{[l_{h},\infty)})\subset M_{[h,\infty)}$. We have that $H(M_{[l_{h},\infty)}\times [0,1])\subset M_{[h,\infty)}$. 
\end{proof}

By carefully following the proof above we have:

\begin{lemma} \label{alreadyll} Assume that $f: M_{Y}\rightarrow M_{X}$ is already locally Lipschitz with respect to $\{M_{[1,2^h]}\}_{h\in \mathbb{N}_{0}}$. Then the simplicial approximation $f': M_{Y}\rightarrow M_{X}$ is locally Lipschitz homotopic to $f$ with respect to the filtration $\{M_{[1,2^h]} \times [0,1]\}_{h\in \mathbb{N}_{0}}$. 
\end{lemma}
\begin{proof} It is clear that $f\circ \pi: \tilde{M}_{Y}\rightarrow M_{X}$ is locally Lipschitz with respect to $\{\tilde{M}^{h}\}_{h\in \mathbb{N}_0}$. We first show that the extension $H_{k-1}$ of $\overline{H}_{k-1}$ is Lipschitz when restricted to $\tilde{M}^{k}=|\mathcal{U}_{k}| \times [2^{k}, 2^{k} +1]$. Consider the subcomplex 
\begin{align*}
L_1\cup L_{1}:=(|\mathcal{U}_{k}| \times \{2^k\}) \cup (|\mathcal{U}_{k}| \times \{2^{k}+1\})
\end{align*}
The general procedure which extends homotopies from subcomplexes of CW complexes actually gives us a Lipschitz homotopy $|\mathcal{U}_{k}| \times [2^k, 2^{k}+1] \times [0,1]\rightarrow M_{X}$. This is because there is a retract
\begin{align*}
D: |\mathcal{U}_{k}| \times [2^k, 2^{k}+1] \times [0,1] 
\rightarrow (|\mathcal{U}_{k}| \times [2^k,2^{k}+1] \times \{0\}) \cup (L_1\times [0,1])\cup  (L_2\times [0,1])
\end{align*}
The homotopy extension we want is then the composition $(f\pi \cup \overline{H}_{k-1}\phi_{k} \cup \id_{f\pi})\circ D$. It is clear that this homotopy will be Lipschitz as long as $D$ and $\overline{H}_{k-1}$ are Lipschitz. $\overline{H}_{k-1}$ is a homotopy between a Lipschitz map $f\pi_{||\mathcal{U}_{k-1}| \times \{2^{k}\}}$ from a finite simplicial complex and its geometric simplicial approximation. This is Lipschitz.\\

Consider the unit square $[0,1]^2 \subset \mathbb{R}^2$ and the point $(\frac{1}{2},2)$. For a point $x\in [0,1]^2$, $r(x)$ is defined as the unique point of intersection between the ray from $(\frac{1}{2},2)$ passing through $x$, and the set $([0,1] \times \{0\}) \cup (\partial[0,1] \times [0,1])$. $r:[0,1]^2 \rightarrow ([0,1] \times \{0\}) \cup (\partial[0,1] \times [0,1])$ is Lipschitz, say with constant $A$. We can let $D$ be defined as  $\id_{|\mathcal{U}_{h}|} \times r$ after identifying $[2^{k},2^{k+1}]\cong [0,1]$. Its Lipschitz constant is $A+1$. \\

The homotopy $\overline{H}_{k}: |\mathcal{U}_{k}| \times [2^k, 2^{k+1}]$ is Lipschitz by Theorem \ref{relgeosimpl}. Therefore, the homotopy $H: \tilde{M}_{Y} \times [0,1] \rightarrow M_{X}$, obtained as the concatenation of $2$ Lipschitz homotopies when restricted to any $\tilde{M}^{h} \times [0,1]$, is locally Lipschitz with respect to this filtration.  \\

Now consider $M_{Y} \times [0,1]$, equipped with the product metric. This is a path metric space, where a length-realising geodesic between $(x,t),(y,s)$ is the concatenation $(\gamma_{xy},t)\ast (y,\eta_{ts})$, where $\gamma_{xy}$ is a geodesic in $M_{Y}$ and $\eta_{ts}$ is the straight line segment in $[0,1]$ between $t$ and $s$. This can be broken up into geodesic segments each completely contained in one $M_{[2^{h},2^{h+1}]} \times [0,1]$ for some $h$, whose lifts realise the product metric on  $\tilde{M}^h \times [0,1]$, for which we know that $H$ is Lipschitz. Therefore $H$ is locally Lipschitz with respect to  $\{M_{[1,2^h]} \times [0,1]\}_{h\in \mathbb{N}_{0}}$. 
\end{proof}

\begin{lemma} \label{liphomotopies} $\psi \xi'$ is locally Lipschitz homotopic to $\id_{M_{Y}}$ and $\xi' \psi$ is locally Lipschitz homotopic to $\id_{M_{X}}$.
\end{lemma}

\begin{proof} We know that $\psi \xi'$ is properly homotopic to $\id_{M_{Y}}$. Give $M_{Y} \times [0,1]$ the CW structure as defined in Proposition \ref{psiexistence}. By the same construction as Lemma \ref{xi'} there exists a proper, locally Lipschitz homotopy $\overline{F}: M_{Y} \times [0,1] \rightarrow M_{Y}$ from  $\widetilde{\psi \xi'}:=\overline{F}_{|M_{Y} \times \{0\}}$ to $\widetilde{\id_{M_{X}}}:= \overline{F}_{|M_{Y}\times \{1\}}$. By Lemma \ref{alreadyll}, $\widetilde{\psi \xi'}$ and $\widetilde{\id_{M_{X}}}$ are proper, locally Lipschitz homotopic to $\psi \xi'$ and $\id_{M_{Y}}$ respectively. This gives us a proper, locally Lipschitz homotopy $F: M_{Y} \times [0,1] \rightarrow M_{Y}$ from $\psi \xi'$ to $\id_{M_{Y}}$. Similarly, there exists a proper, locally Lipschitz homotopy $G: M_{X} \times [0,1]\rightarrow M_{X}$ from $\xi' \psi$ to $\id_{M_{X}}$. 
\end{proof}

By stretching out the interval we obtain proper, locally Lipschitz homotopies 
\begin{align*}
F: I_{p}(M_{Y}) \rightarrow M_{Y}\\
G: I_{p}(M_{X})\rightarrow M_{X}
\end{align*}
from $\psi \xi'$ to $\id_{M_{Y}}$ and from $\xi' \psi$ to $\id_{M_{X}}$ respectively. 

\begin{lemma} \begin{enumerate} \label{globallycoarse}
\item There exists a coarse shrinking map $S_{\xi'}$ on $M_{Y}$, defined as the composition $S_{\xi'}:=i' \circ s_{\xi'} \circ R$, such that $\xi' S_{\xi'}$ is coarse. 
\item There exists a coarse shrinking map $S_{\psi \xi'}$ on $I_{p}(M_{Y})$, defined as the composition $S_{\psi \xi'}:= (i' \times \id) \circ s_{\psi \xi'}  \circ (R \times \id)$, where $s_{\psi \xi'}: I_{p}(cY)\rightarrow I_{p}(cY)$ is a shrinking map on $I_{p}(cY)$, such that $F S_{\psi \xi'}$ is coarse. 
\item There exists a coarse shrinking map $S_{\xi'\psi}$ on $I_{p}(M_{X})$, defined as the composition $S_{\xi'\psi}:=(i' \times \id) \circ s_{\xi'\psi}  \circ (R \times \id)$, where $s_{\xi' \psi}: I_{p}(cX)\rightarrow I_{p}(cX)$ is a shrinking map on $I_{p}(cY)$, such that $G S_{\xi' \psi}$ is coarse. 
\item $S_{\xi'}, S_{\psi \xi'}$ and $S_{\xi'\psi}$ are coarsely homotopic to the identity. 
\end{enumerate}
\end{lemma}

\begin{proof} 
\begin{enumerate}
\item We can modify the map $i$ to a close map $i'$ which is continuous. To do this, observe that $\phi_{h}\circ \varphi_{h}$ is homotopic to $\varphi_{h-1}$ by the linear (in particular, Lipschitz) homotopy $\Psi_{h}: Y \times [2^{h}-1,2^{h}] \rightarrow |\mathcal{U}_{h-1}| \times [2^{h}-1,2^{h}]$ within a simplex. We have that $\Psi_{h} (Y\times \{2^{h}\}) = \phi_{h}\circ \varphi_{h}$ and $\Psi_{h}(Y\times \{2^{h}-1\}) = \varphi_{h-1}$. We simply replace $i$ when restricted to height intervals $[2^{h}-1,2^{h}]$ by $\Psi_{h}$ and observe that this glues to a continuous map $i': cY\rightarrow M_{Y}$ which is locally Lipschitz and distance at most the diameter of a simplex away from $i$. Let $a_{h}$ be an increasing sequence of Lipschitz constants of $i'$ restricted to heights $(2^h,2^{h+1}]$, ie
\begin{align*}
i'_{|c_{(2^h,2^{h+1}]} Y}: c_{(2^h,2^{h+1}]} Y\rightarrow \tilde{M}^h
\end{align*}
After descending to the quotient, the same constants $a_{h}$ are appropriate as Lipschitz constants for 
\begin{align*}
i'_{|c_{[2^h,2^{h+1}]} Y}: c_{[2^h,2^{h+1}]} Y\rightarrow  \pi(\tilde{M}^h)
\end{align*}
where $\pi (\tilde{M}^h)$ is equipped with the path metric induced by the product metric on $|\mathcal{U}_{h}| \times [2^h,2^{h+1}] \cup |\mathcal{U}_{h-1}| \times \{2^{h}\}$.  Since the inclusion $\pi (\tilde{M}^h)\rightarrow M_{Y}$ is $1$-Lipschitz, $i': cY\rightarrow M_{Y}$ is also locally Lipschitz with respect to $\{c_{[2^h,2^{h+1}]} Y\}_{h\in \mathbb{N}_{0}}$with the same sequence of constants.  \\

Let $b_{h}$ be the Lipschitz constant of $\xi'$ (thought of as a map $\tilde{M}_{Y}\rightarrow M_{X}$) restricted to the set $\coprod_{k\leq h} \tilde{M}^h$. By previous arguments, after descending to the quotient, $b_h$ is an approxiate Lipschitz constant for $M_{[1,2^{h+1}]} \subset M_{Y}$ with the subspace metric as well as for $\pi (\tilde{M}^h)$ equipped with the path metric. Since $i'$ is height-preserving, the composition $\xi'i'$ has Lipschitz constant $b_{h}a_{h}$. Choose constants $m_{h}$ such that $2^{m_{h}}>b_{h}a_{h}$ and $m_{h}\leq m_{h+1}$ for all $h\in \mathbb{N}_{0}$. By the same construction as in Lemma \ref{shrinkcoarse} we can define a shrinking map $s_{\xi'}$ such that the composition $\xi' i' s_{\xi'}$ is proper and  uniformly locally $1$-Lipschitz with respect to $\{c_{[2^h,2^{h+1}]} Y\}_{h\in \mathbb{N}_{0}}$. Note that we had to make the map $i'$ continuous rather than just coarse, because even a uniformly bounded discontinuity can become arbitrarily large under post-composition with $\xi'$ if the Lipschitz constants $b_{h}$ are unbounded. Precomposing with the coarse map $R$ gives us the desired result.

\item Let $c_{(2^{h},2^{h+1}]}(Y\times [0,1]) = \{(nx,n,nt)\,|\, x\in Y,  n \in (2^{h},2^{h+1}], t\in [0,1]\}$ be the portion of $I_{p}(cY)$ with heights in $(2^{h},2^{h+1}]$. Let $a_{h}$ be the Lipschitz constants for $i'$ as before, and let $b_{h}$ be the Lipschitz constants for $F$ when restricted to heights $[2^{h},2^{h+1}]$. Choose constants $m_{h}$ such that $2^{m_{h}}>b_{h}a_{h}$ and $m_{h}\leq m_{h+1}$ for all $h\in \mathbb{N}_{0}$. We define the shrinking map as follows: 
\begin{align*}
s^{h}_{\psi \xi'}: c_{(2^{h+m_{h}},2^{h+1+m_{h}}]}(Y\times [0,1])&\longrightarrow c_{(2^{h},2^{h+1}]}(Y\times [0,1])\\
(nx,n,nt) &\longmapsto (\frac{n}{2^{m_{h}}}x,\frac{n}{2^{m_{h}}} , \frac{n}{2^{m_{h}}}t)
\end{align*}
and we send points with heights in the interval $(2^{h+m_{h-1}}, 2^{h+m_{h}}]$ to height $2^h$. By the same computation as in Lemma \ref{shrinkcoarse} the composition $F\circ (i'\times \id)\circ s_{\psi \xi'}$ is proper and uniformly locally $1$-Lipschitz with respect to $\{c_{[2^h,2^{h+1}]}(Y \times [0,1])\}_{h\in \mathbb{N}_0}$. Therefore it is globally coarse, by a similar calculation to Lemma \ref{shrinkcoarse} . The map $R\times \id$ is defined as
\begin{align*}
R\times \id: I_{p}(M_{Y})&\rightarrow c(Y\times [0,1])\\
(x,n,nt) &\mapsto (R(x,n),nt)
\end{align*} Since $R\times \id$ is coarse, we get that $FS_{\psi \xi'}$ is coarse. 

\item This is the same proof as the previous point, except with technical caveats required for Lemma  \ref{defb}. First, let $a_{h}$ be the Lipschitz constants of $i'$ as before. Let $(L_{\psi})_{h}$ be the Lipschitz constant of $\psi$ restricted to $M_{[1,2^{h+1}]}\subset M_{X}$. This has bounded image in $M_{[1,2^{l'+1}]}$ for some $l'\in \mathbb{N}$. Let $(L_{\xi'})_{l'}$ be the Lipschitz constant of $\xi'$ restricted to $M_{[1,2^{l'+1}]} \subset M_{Y}$. Note that $(L_{\xi'})_{l'}(L_{\psi})_{h}$ may be strictly larger than just the Lipschitz constant of $\xi'\psi$ restricted to $M_{[1,2^{h+1}]}$. \\

Let $c_{l'} d_{l'}$ be the Lipschitz constant obtained from going between the path metric on $|\mathcal{U}_{l'}|$ induced from the spherical metric, to the standard metric, back to the path metric.
\begin{align*}
(|\mathcal{U}_{l'}|,d^s) \xrightarrow{c_{l'}}(|\mathcal{U}_{l'}|,d_{std}) \xrightarrow{d_{l'}}(|\mathcal{U}_{l'}|,d^s)
\end{align*}
The reason why these metrics are bi-Lipschitz equivalent is because each $|\mathcal{U}_{l'}|$ is a finite simplicial complex. We can assume that all the sequences $\{c_{l'}\}_{l'\in \mathbb{N}_0}, \{d_{l'}\}_{l'\in \mathbb{N}_0}$, $\{(L_{\xi'})_{l'}\}_{l'\in \mathbb{N}_0}$, as well as $\{a_{h}\}_{h\in \mathbb{N}_0}$, $\{(L_{\psi})_{h}\}_{h\in \mathbb{N}_0}$ are increasing and $\geq 2$. Choose an increasing sequence $\{m'_{h}\}_{h\in \mathbb{N}_0}$ such that 
\begin{align*}
2^{m'_{h}} \geq a_{h}(L_{\psi})_{h}(L_{\xi'})_{l'}c_{l'}d_{l'}
\end{align*} 
Let $b_{h}$ be the Lipschitz constants for $G$ when restricted to heights $[2^h,2^{h+1}]$. We can assume that $2^{m'_{h}}\geq a_{h}b_{h}$. The construction of the shrinking map $s_{\xi'\psi}$ and the rest of the proof is identical to point $2$. of this lemma. 
\item We can assume that the $s$ from Lemma \ref{chi} is the same as $s_{\xi'}$. So $S_{\xi'}=i's_{\xi'}R$ is close to $isR$, which by Lemma \ref{chi} is coarsely homotopic to the identity. Similarly, one can show that $S_{\psi \xi'}=(i' \times \id) \circ s_{\psi \xi'}\circ (R\times \id)$ is close to $(i \times \id) \circ s_{\psi \xi'}\circ (R\times \id)$, which is close to $\mathfrak{q} \circ (i R\times \id)$, where $\mathfrak{q}: I_p(M_{X})\rightarrow I_{p}(M_{X})$ is the projection to the height determined by $s_{\psi \xi'}$. The same calculation works for $S_{\xi'\psi}$. 
\end{enumerate}
\end{proof}

\begin{lemma} \label{che} $\xi' S_{\xi'}$ is a coarse homotopy inverse to $is_{l}\mathfrak{a}R: M_{X}\rightarrow M_{Y}$.  
\end{lemma}

\begin{proof} We compute both compositions. Recall that $\psi = \widetilde{is_{l}\mathfrak{a}R}$ for sufficiently large $l$. By Lemma \ref{globallycoarse}, $F S_{\psi \xi'}$ gives us a coarse homotopy between $\psi \xi' S_{\psi \xi'}$ and $\id_{M_{Y}}S_{\psi \xi'}$. (Note that since the shrinking map $S_{\psi \xi'}$ is uniform in time we have used the same symbol for its restriction to $t=0,1$.) We wish to show the homotopy $\psi \xi' S_{\xi'} \simeq \psi \xi' S_{\psi \xi'}$. This is a similar calculation to Lemma \ref{shrinkindep}. Let $\gamma: [1,\infty)\rightarrow [1,\infty)$ be the shrinking map associated to $s_{\xi'}$ and let $\gamma'$ be associated to $s_{\psi\xi'}$. We can assume that $\gamma'(t)\leq \gamma(t)$ for all $t\in [1,\infty)$. Let 
\begin{align*}
H: I_{p}(cY)&\rightarrow cY\\
(nx,n,nt)&\mapsto (\eta(n,t)x,\eta(n,t))
\end{align*}
be the coarse homotopy from $s_{\xi'}$ to $s_{\psi\xi'}$, where $\eta(n,t) = (1-t)\gamma(n)+t\gamma'(n)$.  We show that $\xi' i' H$ is uniformly locally Lipschitz. Let $\{m_{h}\}_{h\in \mathbb{N}_0}$ be the sequence of constants associated to $\gamma$. Consider two points $(nx,n,nt),(sy,s,sl)$ in $c_{2^{m_{h}}(2^{h},2^{h+1}]}(Y \times [0,1])$. 
\begin{align*}
d_{cY}(H(nx,n,&nt),H(sy,s,sl)) \leq \|\eta(n,t)x - \eta(n,t)y\| + 2|\eta(n,t)-\eta(s,l)|\\
&\leq \frac{n}{2^{m_{h}}} \|x-y\| + 2(|\eta(n,t)-\eta(s,t)|+|\eta(s,t)-\eta(s,l)|)\\
&\leq \frac{1}{2^{m_{h}}} \|nx-ny\|+2(\frac{1}{2^{m_{h}}}|n-s|+|t-l||\gamma(s)+\gamma'(s)|)\\
&\leq \frac{1}{2^{m_{h}}} (\|nx-sy\| + |n-s|) + 2(\frac{1}{2^{m_{h}}}|n-s| +\frac{2}{2^{m_{h}}}|st-sl|)\\
&\leq \frac{1}{2^{m_{h}}} (\|nx-sy\| + |n-s|) + 2(\frac{1}{2^{m_{h}}}|n-s| + \frac{2}{2^{m_{h}}}(|nt-sl|+|n-s|))\\
&= \frac{1}{2^{m_{h}}}(\|nx-sy\| + 7|n-s| + 4|nt-sl|)
\end{align*}
The images of the points have heights in $[2^h,2^{h+1}]$. Therefore the restriction of $\xi' i' H$ to $c_{2^{m_{h}}(2^{h},2^{h+1}]}(Y \times [0,1])$ is $7$-Lipschitz. For points in the height interval $(2^{h+m_{h-1}},2^{h+m_{h}}]$ the calculation is analagous. We have  
\begin{align*}
d_{cY}(H(nx,n,nt),H(sy,s,sl)) \leq \frac{1}{2^{m_{h-1}}}(\|nx-sy\| + 7|n-s| + 4|nt-sl|)
\end{align*}
with image at height $2^{h}$. Therefore $\xi' i' H$ restricted to $c_{(2^{h+m_{h-1}},2^{h+m_{h}}]}(Y\times [0,1])$ is also $7$-Lipschitz. This shows that $\xi' i' H$ is uniformly locally $7$-Lipschitz, and therefore coarse. Hence the composition $\psi \circ (\xi' i' H)\circ (R\times \id)$ is a coarse homotopy between $\psi \xi' S_{\xi'}$ and $\psi \xi' S_{\psi \xi'}$. We have the following sequence of coarse homotopies:
\begin{align*}
(is_{l}\mathfrak{a}R) \xi' S_{\xi'} \simeq \psi \xi' S_{\xi'} \simeq \psi \xi' S_{\psi \xi'} \simeq \id_{M_{Y}} S_{\psi \xi'} \simeq \id_{M_{Y}}
\end{align*}
where in the first equivalence we have used that $\psi$ is close to $is_{l}\mathfrak{a} R$.\\

 To show that $ \xi' S_{\xi'} (is_{l}\mathfrak{a}R) \simeq \id_{M_{X}}$ is more complicated. 
We have a coarse homotopy $GS_{\xi' \psi}$ between $\xi' \psi S_{\xi' \psi}$ and $\id_{M_{X}} S_{\xi' \psi}$. We would like to have the following sequence of homotopies:
\begin{align*}
 \xi' S_{\xi'} (is_{l}\mathfrak{a}R) \simeq   \xi' S_{\xi'} \psi \simeq \xi' S_{\xi'} \psi S_{\xi' \psi} \simeq \xi' \psi S_{\xi' \psi} \simeq \id_{M_{X}} S_{\xi' \psi} \simeq \id_{M_{X}}
\end{align*}
The trouble appears in the homotopy $ \xi' S_{\xi'} \psi S_{\xi' \psi} \simeq \xi' \psi S_{\xi' \psi}$. Since $\xi'$ is only locally Lipschitz, we cannot simply run the homotopy $\xi' H \psi S_{\xi' \psi}$ where $H$ is a coarse homotopy between $S_{\xi'}$ and $\id$ without thinking: any bounded distance can become arbitrarily large under $\xi'$. We therefore run the homotopy in two steps:\\

Let $\gamma$ be the shrinking map corresponding to $s_{\xi'}$, $\gamma'$ the shrinking map corresponding to $s_{\xi' \psi}$. We let $\mathbf{D}$ be the deformation retract
\begin{align*}
\mathbf{D}: I_{p}(M_{Y})&\rightarrow M_{Y}\\
(x,n,nt)&\mapsto \mathfrak{q}_{(1-t)n + t\gamma(n)}(x,n)
\end{align*}
where $\mathfrak{q}$ refers to the height projection. $\mathbf{D}_{1}: M_{Y} \rightarrow M_{Y}$ is the map which shrinks heights via $\gamma$. We first compute the Lipschitz constant of $\mathbf{D}$. Let $(x,n,nt),(y,s,sl)$ be two points in $I_{p}(M_{Y})$. Assume that $n\geq s$. 
\begin{align*}
d(\mathbf{D}(x,n,nt),\mathbf{D}(y,s,sl)) \leq d(\mathbf{D}(x,n,nt),\mathbf{D}(\mathfrak{q}_{s}(x,n),sl)) + d(\mathbf{D}(\mathfrak{q}_{s}(x,n),sl),\mathbf{D}(y,s,sl))\\
=d(\mathfrak{q}_{(1-t)n + t\gamma(n)}(x,n), \mathfrak{q}_{(1-l)s + l\gamma(s)}\mathfrak{q}_{s}(x,n)) +d(\mathfrak{q}_{(1-l)s + l\gamma(s)}\mathfrak{q}_{s}(x,n), \mathfrak{q}_{(1-l)s + l\gamma(s)}(y,s))
\end{align*} 
The first distance is just the straight line of length
\begin{align*}
|(1-t)n +t \gamma(n)-(1-l)s - l\gamma(s)| &\leq 2|nt-sl| + 3|n-s|\\
&\leq 2|nt-sl| + 3d((x,n),(y,s))
\end{align*} 
The distance between the second two points is the distance between the projections of $x$ and $y$ to the height $(1-l)s + l\gamma(s)$, measured within the appropriate nerve $|\mathcal{U}|$. Since the bonding maps are $1$-Lipschitz, we have
\begin{align*}
d(\mathfrak{q}_{(1-l)s + l\gamma(s)}\mathfrak{q}_{s}(x,n), \mathfrak{q}_{(1-l)s + l\gamma(s)}(y,s)) \leq d(\mathfrak{q}_{s}(x,n),(y,s)) \leq d((x,n),(y,s))
\end{align*}
Therefore the homotopy $\mathbf{D}$ is $4$-Lipschitz. Consider the composition 
\begin{align*}
\xi'  \mathbf{D}\circ H: I_{p}(cX)&\rightarrow M_{X}\\
(nx,n,nt)&\mapsto \xi' \mathbf{D}(\psi i'(\gamma'(n) x, \gamma'(n)), p\psi i'(\gamma'(n) x, \gamma'(n)) t)
\end{align*}
where $H: I_{p}(cX)\rightarrow I_{p}(M_{Y})$ is the rescaled version of $\psi i' s_{\xi' \psi}\times \id$. $\xi'  \mathbf{D}\circ H$ is a homotopy from $\xi' \psi i' s_{\xi' \psi}$ to  $\xi' \mathbf{D}_{1} \psi i' s_{\xi' \psi}$. We now show that $H$ is Lipschitz. Recall that $(L_{\psi})_{h}$ is the Lipschitz constant of $\psi$ restricted to heights $[1,2^{h+1}]$ and that $\{m'_{h}\}_{h\in \mathbb{N}_0}$ is the sequence of constants associated to $\gamma'$.  We can assume that $\psi = \widetilde{is_{l}\mathfrak{a}R}$ is height-decreasing by choosing $l$ sufficiently large. Let $(nx,n,nt),(sy,s,sl)$ be two points in $c_{ 2^{m'_{h}}[2^h,2^{h+1}]}(X \times [0,1])$. We have that
\begin{align*}
& d(\psi i'(\gamma'(n) x,\gamma'(n)), \psi i'(\gamma'(s) y,\gamma'(s))) \leq a_{h}(L_{\psi})_{h} d((\gamma'(n) x,\gamma'(n)),  (\gamma'(s) y,\gamma'(s)))\\
&= \frac{a_{h}(L_{\psi})_{h}}{2^{m'_{h}}} d((nx,n), (sy,s))\\
& |p\psi i'(\gamma'(n) x, \gamma'(n)) t-p\psi i'(\gamma'(s) y, \gamma'(s)) l)|\\
&\leq |p\psi i'(\gamma'(n) x, \gamma'(n)) t-p\psi i'(\gamma'(n)x, \gamma'(n)) l|+|p\psi i'(\gamma'(n) x, \gamma'(n)) l-p\psi i'(\gamma'(s) y, \gamma'(s)) l|\\
&= |t-l||p\psi i'(\gamma'(n) x, \gamma'(n))| + |p\psi i'(\gamma'(n) x, \gamma'(n)) -p\psi i'(\gamma'(s) y, \gamma'(s))|\\
&\leq |t-l|\frac{n}{2^{m'_{h}}} + \frac{(L_{\psi})_{h}a_{h}}{2^{m'_{h}}} d((nx,n),(sy,s))\\
&\leq \frac{1}{2^{m'_{h}}} (|nt-sl| + |n-s|) + \frac{(L_{\psi})_{h}a_{h}}{2^{m'_{h}}} d((nx,n),(sy,s))
\end{align*}
A similar computation works for heights in the interval in $ [2^{h+m'_{h}+1}, 2^{h+m'_{h+1}+1}]$. This shows that $H$ restricted to $c_{ [2^{h+m'_{h}},2^{h+m'_{h+1}+1}] }(X \times [0,1])$ is $\frac{3a_{h}(L_{\psi})_{h}}{2^{m'_{h}}}$-Lipschitz.
The computation for heights in the interval $[2^{h-1+m'_{h-1}},2^{h+m'_{h}}]$ is analagous, replacing $h$ with $h-1$. \\

Since in the construction of $S_{\xi'\psi}$, we have maximised over the Lipschitz constants of each successive bounded set, and $\mathbf{D}$ is height-decreasing, the homotopy $\xi \mathbf{D} H$ is uniformly locally $12$-Lipschitz with respect to $\{c_{ [2^{h+m'_{h}},2^{h+m'_{h+1}+1}] }(X \times [0,1])\}_{h\in \mathbb{N}_{0}}$ . Hence it is globally coarse. Precomposing with the coarse map $R \times \id$, we obtain a coarse homotopy between $\xi' \psi S_{\xi'\psi}$ and $\xi' \mathbf{D}_{1} \psi S_{\xi' \psi}$. This was step 1.\\

For the second step, we construct a coarse homotopy between $\xi' \mathbf{D}_{1} \psi S_{\xi' \psi}$ and $\xi' S_{\xi'} \psi S_{\xi' \psi}$ as the composition $\xi' \circ \mathbf{G}\circ \mathbf{F} \circ (R\times \id)$. $\mathbf{G}$ is defined as
\begin{align*}
\mathbf{G}: M_{Y} \times [0,1] &\rightarrow M_{Y}\\
(x,n,t) &\mapsto (1-t) \mathbf{D}_{1}(x,n) + t i's_{\xi'} R (x,n) 
\end{align*}
Why does this make sense? Consider a point $(x,n)$ with $x\in int(\sigma)$ in some appropriate nerve $|\mathcal{U}|$. $\mathbf{D}_{1}(\sigma \times \{ n\}) = \mathfrak{q}_{\gamma(n)}(\sigma \times \{n\}) = \phi(\sigma) \times \{\gamma(n)\}$ is a simplex at height $\gamma(n)$ for the appropriate bonding map $\phi$.  The image of $i's_{\xi'} R(int(\sigma) \times \{n\})$ lies in the set 
\begin{align*}
\bigcap_{v \in \phi (\sigma)} \langle st(v) \rangle  \times \{\gamma(n)\}
\end{align*} 
Therefore, the straight line between $ \mathbf{D}_{1}(x,n) $ and $i's_{\xi'} R (x,n)$ makes sense. Note that this homotopy is very much non-continuous. We define $\mathbf{F}$ as the map
\begin{align*}
\mathbf{F}: I_{p}(cX) &\rightarrow M_{Y} \times [0,1]\\
(nx,n,nt) &\mapsto (\psi i'(\gamma'(n)x,\gamma'(n)), t)
\end{align*}
Let $(nx,n,nt),(sy,s,sl)$ be two points heights in the interval  $[2^{h+m'_{h}},2^{h+m'_{h+1}+1}]$. Recall that we have 
\begin{align*}
d(\psi i'(\gamma'(n) x,\gamma'(n)), \psi i'(\gamma'(s) y,\gamma'(s)))
\leq \frac{a_{h}(L_{\psi})_{h}}{2^{m'_{h}}} d((nx,n), (sy,s))\\ 
|t-l| = \frac{1}{n}|nt-nl| \leq \frac{1}{2^{m'_{h}+h}}(|nt-sl| + |n-s|)
\end{align*}

Now we have to show that for a set sufficiently small in $M_{Y}\times [0,1]$, the image under $\mathbf{G}$ is actually sufficiently small. Let $\varepsilon, \delta$ be very small constants, to be determined later. Suppose we have two points $(x,n,t),(y,s,l) \in M_{Y} \times [0,1]$ below height $2^{l'}$, contained within $U_{c_{l'}\varepsilon} \times I_{\varepsilon} \times I_{\delta}$, where $I_{\varepsilon}$ only intersects at most one a gluing slice on the boundary and $U_{c_{l'}\varepsilon}$ is a set of diameter $<c_{l'}\varepsilon$ (with respect to the standard metric) contained within a simplex. The reason why we can restrict to this case is that for an arbitrary set $V$ of diameter $<\varepsilon'$ (with respect to the path metric) in $M_{Y}$ where $c_{l'}\varepsilon'<\frac{1}{2}$,  any two points $(u,n),(w,s)$ in $V$ can be split into at most $5$ segments of the form we have assumed. \\

To see this, let $n\geq s$. First we assume that the shortest path between $(u,n),(w,s)$ consists of $\omega_{(u,n)}\ast \eta_{s}$ where $\omega_{(u,n)}$ is the straight line downards and $\eta_{s} \subset (U_{\varepsilon'} \times \{s\})$ is a geodesic at height $s$ between $\mathfrak{q}_{s}(u)$ and $w$. Split $\omega_{(u,n)}$ into two paths which intersect at most one gluing slice. Now for two points $x,y$ in $U_{\varepsilon'}$, we consider their image in $(|\mathcal{U}|, d_{std})$, which have distance $< c_{l'}\varepsilon'<\frac{1}{2}$. Let $x = \sum_{v \in vert |\mathcal{U}|} x_{v} [v]$ and $y =  \sum_{v \in vert |\mathcal{U}|} y_{v} [v]$. Define
\begin{align*}
Z = \{v\in vert|\mathcal{U}|\,|\, x_{v}\neq 0, y_{v} \neq 0\}
\end{align*}
 and $a = \sum_{v \in Z} x_{v}$ and $b = \sum_{v \in Z} y_{v}$. The distance between $x$ and $y$ is 
\begin{align*}
d_{std}(x,y) = (1-a) + (1-b) + \sum_{v \in Z} |x_{v}-y_{v}| < c_{l'}\varepsilon'
\end{align*} 
Define $z_{x} = \frac{1}{a}\sum_{v\in Z} x_{v} [v]$ and $z_{y} = \frac{1}{b} \sum_{v\in Z} y_{v} [v]$. We compute
\begin{align*}
d_{std}(x,z_{x}) &= (1-a) + \sum_{v\in Z} (\frac{1}{a}-1) x_{v} = 2(1-a) < 2c_{l'}\varepsilon'\\
d_{std}(y,z_{y}) &< 2c_{l'}\varepsilon'\\
a d_{std}(z_{x},z_{y}) &= \sum_{v \in Z} |x_{v}-\frac{a}{b} y_{v}|\leq \sum_{v \in Z} |x_{v}- y_{v}|+ |\frac{a}{b}-1| y_{v}\\
&= (\sum_{v \in Z} |x_{v}- y_{v}|)+ |a-b| \leq (\sum_{v \in Z} |x_{v}- y_{v}|)+ |1-a| + |1-b| < c_{l'}\varepsilon'\\
d_{std}(z_{x},z_{y}) &<\frac{c_{l'}\varepsilon'}{a} <\frac{c_{l'}\varepsilon'}{1-c_{l'}\varepsilon'} <2c_{l'}\varepsilon'
\end{align*}
Therefore $\eta_{s}$ can be split into three straight paths of length less than $2c_{l'}\varepsilon'$, each contained within a simplex. \\

The only other case is if $n,s\in (2^{h},2^{h}+1]$ and the shortest path between $(u,n),(w,s)$ is a concatenation $\omega_{(u,n)} \ast \eta_{2^h} \ast \omega_{(w,s)}^{-1}$, where $\eta_{2^h}$ is a geodesic at height $2^h$. As before, we can split $\eta_{2^h}$ into three paths of the desired form.\\

Let us return to the proof. Recall we have two points $(x,n,t),(y,s,l) \in M_{Y} \times [0,1]$ below height $2^{l'}$, contained within $U_{c_{l'}\varepsilon} \times I_{\varepsilon} \times I_{\delta}$. Assume that $n\geq s$. We consider the points $(x,n,t),(x,s,t),(x,s,l),(y,s,l)$ under the image of $\mathbf{G}$ with respect to the standard metric $d_{std}$:
\begin{enumerate}
\item $d_{std}(\mathbf{G}(x,s,t), \mathbf{G}(x,s,l))$: Since the homotopy is the straight line between two points that lie within a simplex, it is $1$-Lipschitz in the time variable. Therefore $d_{std}(\mathbf{G}(x,s,t), \mathbf{G}(x,s,l))\leq |t-l|<\delta$. 
\item $d_{std}(\mathbf{G}(x,s,l), \mathbf{G}(y,s,l))$: We know that $d_{std}(\mathbf{D}_{1}(x,s), \mathbf{D}_{1}(y,s))\leq d_{std}(x,y)< c_{l}'\varepsilon$. Let $\sigma$ be the simplex that contains $x,y$. We know from previous calculations that $R(\sigma)$ has diameter less than $\frac{69}{4} = C$ (a conservative upper bound). Recall that $i's_{\xi'}$ is $\frac{1}{2^{b_{h'}}}$-Lipschitz for some $h'<l'$ where $2^{b_{h'}}$ is sufficiently large so that $\xi' i's_{\xi'}$ is uniformly locally $1$-Lipschitz. We can assume that $2^{b_{h'}}$ is also greater than the local Lipschitz constant $(L_{\xi'})_{h'}$ multipled by the Lipschitz constant $c_{h'}d_{h'}$ which we obtain from changing between metrics. Therefore $d(i's_{\xi'}R(x,s), i's_{\xi'}R(y,s))\leq \frac{C}{2^{b_{h'}}}$. In standard metric, $d_{std}(\mathbf{G}(x,s,l), \mathbf{G}(y,s,l))\leq \max\{\frac{Cc_{h'}}{2^{b_{h'}}},c_{l}'\varepsilon\}$. To see this, observe that for four points $v_1,v_2,w_1,w_2$ in Euclidean space, 
\begin{align*}
\|(1-t)v_1+tw_1 - (1-t)v_2 + tw_2\| \\
\leq (1-t)\|v_1-v_2\| + t\|w_1-w_2\| \leq \max \{\|v_1-v_2\|, \|w_1-w_2\|\}
\end{align*}
\item $d_{std}(\mathbf{G}(x,n,t), \mathbf{G}(x,s,t))$:  Ignore first the possible discontinuity related to $R$ at the boundary. We have that $R(x,n)$ and $R(x,s)$ get sent to the same $x$-coordinate, so their distance in $cY$ is less than $2|n-s|$. Remember that by construction, $i' s_{\xi'}$ is actually $\frac{1}{2^{b_{h'}}}$-Lipschitz with respect to the product metric on each $\tilde{M}^{h'}$ before gluing.  Because the shrinking maps are constructed using powers of $2$, the image of $s_{\xi'}(x,n)$ and $s_{\xi'} (x,s)$ do not intersect a gluing component except possibly at the boundary - we are allowed to compute distances in one such $\tilde{M}^{h'}$. Therefore, the $|\mathcal{U}_{h'}|$-component of $i's_{\xi'}R(x,n)$ and $i's_{\xi'} R(x,s)$ have distance less than than $2\cdot \frac{1}{2^{b_{h'}}} |n-s|$. Since $\mathbf{D}_1(x,n)$ and $\mathbf{D}_{1}(x,s)$ have the same $x$-coordinate, by the previous case we conclude that 
\begin{align*}
d_{std}(\mathbf{G}(x,n,t), \mathbf{G}(x,s,t))\leq (2\cdot \frac{c_{h'}}{2^{b_{h'}}}+1) |n-s|
\end{align*}
If we do not intersect a gluing slice, this finishes the estimate. If $n$ or $s$ lie on a gluing slice there is a discontinuity to consider. By the previous case, the discontinuity is bounded by $\frac{Cc_{h'}}{2^{b_{h'}}}$, so in total we get
\begin{align*}
d_{std}(\mathbf{G}(x,n,t), \mathbf{G}(x,s,t))< \frac{Cc_{h'}}{2^{b_{h'}}} + (2\cdot \frac{c_{h'}}{2^{b_{h'}}}+1) \varepsilon <\frac{Cc_{h'}}{2^{b_{h'}}} + 3\varepsilon
\end{align*}  
in both cases.
\end{enumerate} 

Therefore, $d_{std}(\mathbf{G}(x,n,t),\mathbf{G}(y,s,l))<  \delta +\max\{\frac{Cc_{h'}}{2^{b_{h'}}},c_{l}'\varepsilon\}+\frac{Cc_{h'}}{2^{b_{h'}}} + 3\varepsilon$. We now let $\varepsilon =  \frac{a_{h}(L_{\psi})_{h}}{2^{m'_{h}}}$ and $\delta = \frac{1}{2^{m'_{h}+h}}<\varepsilon$, which gives us $d_{std}(\mathbf{G}(x,n,t),\mathbf{G}(y,s,l))<  (5c_{l'}\varepsilon + 2\frac{Cc_{h'}}{2^{b_{h'}}})$. Recall that $2^{m_{h}'}> a_{h} (L_{\psi})_{h} (L_{\xi'})_{l'}c_{l'}d_{l'}$ and that $2^{b_{h'}} > c_{h'}d_{h'} (L_{\xi'})_{h'}$.  Postcomposing with the locally Lipschitz map $\xi'$, we obtain 
\begin{align*}
d_{path}(\xi'\mathbf{G}(x,n,t),\xi'\mathbf{G}(y,s,l))&\leq (L_{\xi'})_{h'} d_{h'} (5c_{l'}\varepsilon + 2\frac{Cc_{h'}}{2^{b_{h'}}})\\
&\leq (L_{\xi'})_{l'}d_{l'}\frac{5c_{l'}a_{h}(L_{\psi})_{h}}{2^{m'_{h}}} +2(L_{\xi'})_{h'} d_{h'}\frac{Cc_{h'}}{2^{b_{h'}}}\\
&<5 +2C
\end{align*}
Now we finally show that $\xi\circ \mathbf{G}\circ \mathbf{F}$ is a coarse homotopy. Consider $(nx,n,nt),(sy,s,sl)$, two points in $I_{p}(cX)$ with heights in the interval $ [2^{h+m'_{h}},2^{h+m'_{h+1}+1}]$ and distance less than $Q\in \mathbb{N}$. $\mathbf{F}(nx,n,nt), \mathbf{F}(sy,s,sl)$ live below height $2^{l'}$ and are contained in $V_{Q\varepsilon} \times I_{Q\delta}$. Therefore, the shortest path between them can be split up into $2Q$ segments each lying within $V_{\frac{\varepsilon}{2}} \times I_{\delta}$. By previous discussion we can split this further to obtain $10Q$ path segments each lying within a set of the form $U_{c_{l'}\varepsilon} \times I_{\varepsilon} \times I_{\delta}$ with respect to the standard metric. The endpoints of each of these segments, after composing with $\xi' \mathbf{G}$, are distance $< 5+2C$ apart. We obtain
\begin{align*}
d_{path}(\xi'\mathbf{GF}(nx,n,nt),\xi'\mathbf{GF}(sy,s,sl))&\leq 10Q(5+2C)
\end{align*}
Therefore, a set of diameter $Q\in \mathbb{N}$ in $c_{[2^{h+m'_{h}},2^{h+m'_{h+1}+1}]}(X \times [0,1])$ gets sent to a set of diameter $(50+20C)Q$ under the map $\xi' \mathbf{GF}$. By a similar argument to Lemma \ref{shrinkcoarse}, for arbitrary points $(nx,n,nt),(sy,s,sl)$ in $I_{p}(cX)$ and $Q>0$, there exists a constant $T_{Q}$ such that 
\begin{align*}
d((nx,n,nt),(sy,s,sl))<Q \\
\implies d(\xi' \mathbf{GF}(nx,n,nt),\xi'  \mathbf{GF}(sy,s,sl)) < \max\{3(50+20C)Q, (50+20C)Q+T_{Q}\}
\end{align*} This shows that $\xi' \mathbf{GF}$ is coarse. By precomposing with the coarse map $R \times \id$, we obtain that $\xi' \mathbf{D}_{1} \psi S_{\xi' \psi}$ and $\xi' S_{\xi'} \psi S_{\xi' \psi}$ are coarsely homotopic. \\

The two steps show that $\xi' \psi S_{\xi'\psi}\simeq \xi' \mathbf{D}_{1} \psi S_{\xi' \psi}\simeq \xi' S_{\xi'} \psi S_{\xi' \psi}$ which was the required homotopy to conclude 
\begin{align*}
\xi' S_{\xi'} (is_{l}\mathfrak{a}R) \simeq \id_{M_{X}}
\end{align*}

Therefore, $\xi' S_{\xi'}$ is a coarse homotopy inverse to $is_{l}\mathfrak{a}R$. 
\end{proof}

\begin{lemma} \label{defb} The composition $\mathfrak{b}:=(R \circ \xi' S_{\xi'}\circ is)$ is a coarse homotopy inverse to $\mathfrak{a}$.
\end{lemma}

\begin{proof} We compute
\begin{align*}
R \xi' S_{\xi'} is \mathfrak{a}  \simeq R(\xi' S_{\xi'})(is_{l}\mathfrak{a}R) is \simeq R \id_{M_{X}} is \simeq \id_{cX}\\
\mathfrak{a} R \xi' S_{\xi'} is \simeq R(is_{l}\mathfrak{a}R) (\xi' S_{\xi'}) is \simeq R \id_{M_{Y}}is \simeq \id_{cY}
\end{align*}
\end{proof}

This concludes the proof of Theorem \ref{unbasedwhitehead}. \\

There is also a based version of coarse Whitehead, which can be obtained from observing that all our constructions can be chosen to be maps of pairs, if $\mathfrak{a}$ is a map of pairs.

\begin{cor} \label{basedwhitehead} Let $(X,x_0)$ and $(Y,y_0)$ be connected, compact metric spaces with finite shape dimension. Let $\omega,\tau$ be the standard parametrisations of the base rays $c\{x_0\}, c\{y_0\}$ respectively. Let $\mathfrak{a}:(cX,\omega)\rightarrow (cY,\tau)$ be a coarse map of pairs, which induces an isomorphism on coarse homotopy groups $\mathfrak{a}: \pi_{n}^{c}(cX,\omega)\rightarrow \pi_{n}^{c}(cY,\mathfrak{a}\omega)$ for all $n\geq 0$. Then $\mathfrak{a}$ is a coarse homotopy equivalence of pairs, ie. there a coarse map of pairs $\mathfrak{b}: (cY, \tau)\rightarrow (cX,\omega)$ such that $\mathfrak{b}\circ \mathfrak{a} \sim id_{cX}$ and $\mathfrak{a}\circ \mathfrak{b} \sim id_{cY}$ as coarse homotopies of pairs.
\end{cor}

\newpage

\section{The functor $\mathcal{F}$} 

Recall that $\mathbf{BornCoarse}$ has objects bornological coarse spaces and morphisms proper and controlled maps. Consider now a set $X$ equipped with a coarse structure $\mathcal{C}$. We define a bornology $\mathcal{B}$ on $X$ by letting $B\in \mathcal{B}$ if there exists $x\in X$ such that $B\times \{x\}\in \mathcal{C}$. We say that $\mathcal{B}$ is \textit{induced from the coarse structure}. Recall that a controlled map $f:X\rightarrow Y$, where $Y$ is an arbitrary bornological coarse space, is automatically bornologous. We denote by $\textbf{Coarse}$ the subcategory of $\textbf{BornCoarse}$ consisting of coarse spaces with the induced bornology from the coarse structure, while retaining the morphisms from $\textbf{BornCoarse}$.

\begin{definition} Let $\mathbf{hCoarse}$ be the category with the same objects as $\mathbf{Coarse}$ but where morphisms are coarse homotopy classes of coarse bornologous maps. 
\end{definition}

The reason why we could not define $\mathbf{hBornCoarse}$ directly is because composition of coarse homotopy classes is only well-defined for coarse bornologous maps: see Proposition \ref{properties}. 

\begin{theorem} \label{SSBornCoarse} There is a full functor 
\begin{align*}
\mathcal{F}: \mathbf{StrongShape} \rightarrow \mathbf{hCoarse}
\end{align*}
\end{theorem}

\begin{conjecture} The functor $\mathcal{F}$ from Theorem \ref{SSBornCoarse} is faithful. 
\end{conjecture}

\begin{definition} \label{DefF} (The functor $\mathcal{F}$)
Let $(W,X)$ be an object in $\mathbf{StrongShape}$. We define
\begin{align*}
\mathcal{F}(W,X) := cX 
\end{align*}
the Euclidean cone over the compact metric space $X$. \\

Let $[\overline{\psi}]: (W_{X},X) \rightarrow (W_{Y},Y)$ be a morphism in $\mathbf{StrongShape}$, ie.  a proper homotopy class of proper maps $W_{X}\setminus X\rightarrow W_{Y}\setminus Y$. By Theorem \ref{SSandS}, we know that $W_{X}\setminus X \simeq M_{X}$ and $W_{Y}\setminus Y \simeq M_{Y}$. By pre- and postcomposing with the proper homotopy equivalences $\zeta_{X}:M_{X}\rightarrow W_{X}\setminus X$ and $\zeta^{-1}_{Y}: W_{Y}\setminus Y\rightarrow M_{Y}$, we obtain a proper homotopy class of proper maps $[\zeta^{-1}_{Y}\overline{\psi}\zeta_{X}]=:[\psi]: M_{X}\rightarrow M_{Y}$. We define  
\begin{align*}
\mathcal{F}([\overline{\psi}]) := [R_{Y} \psi' S_{\psi'} i_{X} s_{X}]
\end{align*}
where $\psi'$ is a locally Lipschitz approximation of $\psi$, $S_{\psi'}$ is a coarse shrinking map on $M_{X}$, and $i_{X}s_{X}: cX\rightarrow M_{X}$, $R_{Y}: M_{Y}\rightarrow cY$ are the coarse homotopy equivalences defined previously. 
\end{definition}

We dedicate the rest of this section to showing that $F$ is well-defined, and that it is full. We will outline a proof attempt that it is faithful. \\

\textbf{$F$ is well-defined on morphisms}: Let $\overline{\eta}$ be another representative of $[\overline{\psi}]$. Then $\eta$ is properly homotopic to $\psi$. We have to show that there is a coarse homotopy  
\begin{align*}
R_{Y} \psi' S_{\psi'} i_{X}s_{X} \simeq R_{Y} \eta' S_{\eta' } i_{X} s_{X}
\end{align*}

By Lemma \ref{xi'}, $\psi$ and $\eta$ have locally Lipschitz representatives $\psi', \eta'$. By Lemma \ref{globallycoarse} there are coarse shrinking maps $S_{\psi'}$ and $S_{\eta'}$ such that $\psi' S_{\psi'}$ and $\eta' S_{\eta'}$ are coarse. By Lemma \ref{liphomotopies}  there exists a locally Lipschitz homotopy $F: I_{p}(M_{X})\rightarrow M_{Y}$ from $\psi'$ to $\eta'$. There exists a coarse shrinking map $S_{F}: I_{p}(M_{X})\rightarrow I_{p}(M_{X})$ such that $FS_{F}$ is a coarse homotopy between $\psi' S_{F}$ and $\eta' S_{F}$. By the same argument as in Lemma \ref{che}, we have that $\psi' S_{\psi'} \simeq \psi' S_{F}$ and  $\eta' S_{F} \simeq \eta' S_{\eta'}$. Therefore, we obtain: 
\begin{align*}
 \psi' S_{\psi'} \simeq \psi' S_{F}  \simeq \eta' S_{F} \simeq \eta' S_{\eta'}
\end{align*}
and the claim follows. 

\begin{remark} This proof also tells us that the coarse homotopy class of $R_{Y} \psi' S_{\psi'} i_{X}s_{X}$ is independent of the choice of simplicial approximation $\psi'$ since any two choices of locally Lipschitz simplicial approximation are properly homotopic. 
\end{remark}

\textbf{$\mathcal{F}$ is a functor}: To prove that $\mathcal{F}$ is functorial, we have to show that $\mathcal{F}$ respects composition of morphisms and that $\mathcal{F}([\id_{M_{X}}]) = [\id_{cX}]$. Suppose that $X,Y,Z$ are compact metric spaces and let 
\begin{align*}
W_{X}\setminus X \xrightarrow{\overline{\psi}} W_{Y}\setminus Y \xrightarrow{\overline{\eta}} W_{Z} \setminus Z
\end{align*}
be a composition of morphisms. We have that $\zeta^{-1}_{Z}\overline{\eta} \overline{\psi} \zeta_{X} \simeq \eta \psi$. We show that $\mathcal{F}([\eta\circ \psi]) = \mathcal{F}([\eta]) \circ \mathcal{F}([\psi])$ ie. there is a coarse homotopy 
\begin{align*}
R_{Z}(\eta\psi)' S_{(\eta\psi)'} i_{X} s_{X} \simeq (R_{Z} \eta' S_{\eta'} i_{Y} s_{Y}) \circ (R_{Y} \psi' S_{\psi'} i_{X}s_{X})
\end{align*}
Since $\eta' \psi'$ is a locally Lipschitz map which is properly homotopic to $(\eta\psi)'$, it is also locally Lipschitz homotopic to $(\eta \psi)'$ via a homotopy $F: I_{p}(M_{X})\rightarrow M_{Z}$. As before, there exists a shrinking map $S_{F}$ such that $FS_{F}$ is coarse. Therefore we obtain 
\begin{align*}
\eta' S_{\eta'} \psi' S_{\psi'} \simeq \eta' S_{\eta'} \psi' S_{F} \simeq (\eta' \psi') S_{F} \simeq (\eta \psi)' S_{F} \simeq (\eta\psi)' S_{(\eta\psi)'}
\end{align*}
where the argument from Lemma \ref{che} is required for the homotopy $\eta' S_{\eta'} \psi' S_{F} \simeq (\eta' \psi') S_{F}$. The claim follows. \\

It is clear that $\mathcal{F}([\id_{M_{X}}]) = [\id_{cX}]$ since $R_{X}i_{X}s_{X}$ is coarsely homotopic to the identity. \\

\textbf{$\mathcal{F}$ is full:} To prove that $\mathcal{F}$ is full, we need to show that for a coarse map $\mathfrak{a}: cX\rightarrow cY$, there exists a $\psi: M_{X}\rightarrow M_{Y}$ such that $[\mathfrak{a}] = [R_{Y} \psi' S_{\psi'} i_{X}s_{X}]$. To that end, choose $\psi$ to be the geometric cellular approximation $\widetilde{i_{Y}s_{l}\mathfrak{a}R_{X}}$ for $l$ sufficiently large. It is clear that  $\psi'S_{\psi'} = \psi$ and that 
\begin{align*}
R_{Y} \circ (\widetilde{i_{Y}s_{l}\mathfrak{a}R_{X}})\circ i_{X}s_{X} \simeq (R_{Y}i_{Y})\circ s_{l} \mathfrak{a}\circ (R_{X} i_{X}s_{X}) \simeq s_{l} \mathfrak{a}\simeq \mathfrak{a}
\end{align*}\\
\textbf{$F$ is faithful:} To prove that $F$ is faithful, we need to show that if $[R_{Y} \psi' S_{\psi'} i_{X}s_{X}] = [R_{Y} \eta' S_{\eta'} i_{X}s_{X}]$ then $\psi$ and $\eta$ are properly homotopic. Unfortunately this appears to be difficult. The trouble lies in the fact that even though $S_{\psi'}$ is coarsely homotopic to the identity, it does not follow trivially that $\widetilde{\psi' S_{\psi'}}$ is properly homotopic to $\psi'$. What would need to be proven is the following
\begin{align*}
\widetilde{\psi' S_{\psi'}} \simeq \widetilde{\psi' \widetilde{S_{\psi'}}} \simeq \psi' \widetilde{S_{\psi'}}
\end{align*}
for some suitable definition of simplicial approximation $\widetilde{S_{\psi'}}$. This is technically tricky: recall that by construction $\psi'$ is already cellular, but with respect to a CW-structure on $M_{X}$ that becomes more and more subdivided as we go up the telescope. So we would have to redefine simplicial approximation again for target space this iteratively-subdivided CW complex. We leave this statement for future research.\\

Even without faithfulness, the functor $\mathcal{F}$ respects equivalences. To prove this, we need some auxiliary lemmas.

\begin{lemma} \label{chimpliesph} Let $H:I_p(cX)\rightarrow cY$ be a coarse homotopy between $\mathfrak{a}$ and $\mathfrak{a}'$. Then there exists a proper homotopy $M_{X} \times [0,1] \rightarrow M_{Y}$ from $\widetilde{i_{Y}s_{l}\mathfrak{a}R_{X}}$ to $\widetilde{i_{Y}s_{l}\mathfrak{a}'R_{X}}$ for sufficiently large $l$. 
\end{lemma}

\begin{proof}

Recall that $I_{p}(cX) = c(X\times [0,1])$. We take distances on $X\times [0,1]$ to be measured by the maximum of their distances in $X$ and $[0,1]$ where $[0,1]$ carries the standard metric. 
\begin{align*}
d_{\infty}((x,t),(y,s)) = \max \{d_{X}(x,y), |t-s|\}
\end{align*} Recall that we have a fixed sequence of $\frac{1}{2^{h+1}}$-separated sets $\{Z_{h}\}_{\in \mathbb{N}_0}$ in $X$. We can construct a maximal $\frac{1}{2^{h+1}}$-separated set $Z'_{h}$ on $X\times [0,1]$  by simply letting it be copies of $Z_{h}$ at points separated by $\frac{1}{2^{h}}$ in $[0,1]$:
\begin{align*}
Z'_{h}:= \bigcup_{0\leq k\leq 2^h} Z_{h} \times \{ \frac{k}{2^{h}}\} =: Z_{h} \times  A_{h}
\end{align*}
It is clear that $Z'_{h}$ is a maximal $\frac{1}{2^{h+1}}$-separated set with respect to the metric $d_{\infty}$, and that $Z'_{h}\subset Z'_{h+1}$ for all $h\in \mathbb{N}_{0}$. Take the open cover $\mathcal{U}'_{h}$ of balls of radius $\frac{1}{2^h}$ with respect to $d_{\infty}$ and now we can construct the mapping telescope of $M_{X\times [0,1]}$ exactly as before. There are analogous constructions of $i_{X\times [0,1]}$ and $ R_{X\times [0,1]}$ for $M_{X\times [0,1]}$. It is clear that $c(X\times [0,1])$ and $M_{X\times [0,1]}$ are coarsely homotopy equivalent. \\

Given a coarse homotopy $H: c(X\times [0,1])\rightarrow cY$, there exists a $l$ large enough such that the geometric cellular approximation $(i_{Y}s_{l}HR_{X\times [0,1]})^{\sim}$ exists. This gives us a proper, locally Lipschitz map: 
\begin{align*}
(i_{Y}s_{l}HR_{X\times [0,1]})^{\sim}: M_{X\times [0,1]} \rightarrow M_{Y}
\end{align*}

We define a CW complex $\mathbf{M}_{X \times [0,1]}$, which we will show is a subcomplex of $M_{X \times [0,1]}$. Consider the simplicial complex
\begin{align*}
K_{h} =|\mathcal{U}_{h}| \times [0,2^h]
\end{align*}
where $[0,2^h]$ is the locally ordered simplicial complex with vertices at $\mathbb{N}_{0}$ and total order inherited from $\mathbb{N}_{0}$. Recall that there is a fixed total order on the vertices of $|\mathcal{U}_{h}|$. $|\mathcal{U}_{h}| \times [0,2^h]$ is given the product simplicial structure. There are simplicial maps
\begin{align*}
\psi_{h+1,h}: |\mathcal{U}_{h+1}| \times [0,2^{h+1}] &\rightarrow |\mathcal{U}_{h}| \times [0,2^{h}]
\end{align*}
defined as the linear extension of the vertex map $(z,n)\mapsto (\phi_{h+1}(z),\left \lfloor{\frac{n}{2}}\right\rfloor)$. Let the linear extension of $n\mapsto \left \lfloor{\frac{n}{2}}\right\rfloor$ be the function denoted by $\phi'_{h+1}$. $\psi_{h+1,h}$ can be expressed as $\phi_{h+1} \times \phi'_{h+1}$. \\

Let $\mathbf{M}_{X \times [0,1]}$ be the inverse mapping telescope  
\begin{align*}
\mathbf{M}_{X \times [0,1]}: = K_{0} \times [1,2]\cup_{\psi_{1,0}} K_{1} \times [2,4]\cup_{\psi_{2,1}}\dots 
\end{align*}

A coninuous, proper homotopy $ \mathbf{M}_{X \times [0,1]} \rightarrow M_{Y}$ gives us a homotopy $M_{X}\times [0,1]\rightarrow M_{Y}$. This is because in the time component of $\mathbf{M}_{X \times [0,1]}$  what we have is actually the mapping telescope
\begin{align*}
M_{[0,1]}: = [0,1] \times [1,2] \cup_{\phi'_{1}} [0,2] \times [2,4] \cup_{\phi'_{2}}\dots \cup_{\phi'_{h}} [0,2^h]\times [2^h,2^{h+1}] \cup_{\phi'_{h+1}} \dots
\end{align*}
of a shape expansion of the interval: the simplicial complexes $[0,2^h]$ correspond to the nerves of the open cover
\begin{align*}
\cup_{z\in A_{h}} B(z,\frac{1}{2^{h}})
\end{align*}
and the $\phi'$ are simplicial maps induced by refinement. There is a continuous map 
\begin{align*}
i'_{[0,1]}: [0,1]\times [1,\infty)\rightarrow M_{[0,1]}
\end{align*}
defined piecewise by partitions of unity subordinate to $\cup_{z\in A_{h}} B(z,\frac{1}{2^{h}})$ and modified by coherent homotopies to ensure continuity. Therefore, there is a continuous map
\begin{align*}
 \alpha: M_{X} \times [0,1] &\rightarrow \mathbf{M}_{X \times [0,1]}\\
(x,t) &\mapsto (x, i'_{[0,1]}(t,p(x)))
\end{align*} such that $\alpha$ restricted to $M_{X} \times \{0,1\}$ can be identified with the identity on $M_{X}$.  \\

 A set of vertices in $|\mathcal{U}'_{h}|$ span a simplex exactly when their projections to $|\mathcal{U}_{h}|$ span a simplex, and their projections to $[0,2^h]$ span a simplex. Therefore a simplex $\tau\subset \sigma \times I$, where $I$ denotes an interval of length $1$ in $\mathbf{M}_{X \times [0,1]}$, also exists in $M_{X\times [0,1]}$ but as the face of a high-dimensional simplex that connects all the vertices of $ \sigma \times I$. Moreover, the bonding maps $\psi_{h+1,h}$ agree with those of $M_{X \times [0,1]}$. These facts are sufficient to show that $\mathbf{M}_{X \times [0,1]}$ is a subcomplex of $M_{X\times [0,1]}$.\\

We take the restriction of  $(i_{Y}s_{l}HR_{X\times [0,1]})^{\sim}$ to the subcomplex $\mathbf{M}_{X \times [0,1]}$ and observe that because a star in $\mathbf{M}_{X \times [0,1]}$ is contained inside the star of $M_{X\times [0,1]}$, the restriction of  $(i_{Y}s_{l}HR_{X\times [0,1]})^{\sim}$ to $M_{X} \times \{0\}$ and $M_{X} \times \{1\}$ are simplicial approximations of $i_{Y}s_{l}\mathfrak{a}R_{X}$ and $i_{Y}s_{l}\mathfrak{a}'R_{X}$ respectively (we also require the fact that $R_{X\times [0,1]}$ restricted to the boundaries $M_{X\times \partial [0,1]}$ can be chosen to agree with $R_{X}$). Therefore, we have a proper homotopy from $\widetilde{i_{Y}s_{l}\mathfrak{a}R_{X}}$ and $\widetilde{i_{Y}s_{l}\mathfrak{a}'R_{X}}$ for sufficiently large $l$. 
\end{proof}

\begin{cor} \label{lindep} Let $\mathfrak{a}: cX\rightarrow cY$ be a coarse map. There exists a $l$ sufficiently large such that the geometric cellular approximation $\widetilde{i_{Y}s_{l}\mathfrak{a}R_{X}}$ satisfies 
\begin{align*}
\widetilde{i_{Y}s_{l}\mathfrak{a}R_{X}} \simeq \widetilde{i_{Y}s_{l'}\mathfrak{a}R_{X}}
\end{align*} 
for all $l'\geq l$.
\end{cor}

\begin{proof} Let $\mathbf{H}^1$ be the coarse homotopy from $\id: cY\rightarrow cY$ to $s_1$, the shrinking map which divides heights by $2$. By Proposition \ref{properties} there is a coarse homotopy $I_{p\mathfrak{a}R_{X}}(cX) \rightarrow cY$ between $\mathfrak{a}R_{X}$ and $s_{1}\mathfrak{a}R_{X}$, obtained by composing $\mathbf{H}^1$ with $\mathfrak{a}R_{X}$. By Lemma \ref{p_0} we obtain a coarse homotopy $H: c(X \times [0,1]) \rightarrow cY$ between $\mathfrak{a}R_{X}$ and $s_{1}\mathfrak{a}R_{X}$. Apply now Lemma \ref{chimpliesph} to $H$. There is a $l$ large enough such that $\widetilde{i_{Y}s_{l}\mathfrak{a}R_{X}} \simeq \widetilde{i_{Y}s_{l+1}\mathfrak{a}R_{X}}$. Replace $l$ with $l+1$, and inductively we have our claim. 
\end{proof}

\begin{theorem} \label{SSBornCoarse2} Two objects $(W_{X},X),(W_{Y},Y)$ are equivalent in $\mathbf{StrongShape}$ if and only if their images $\mathcal{F}(W_{X},X),\mathcal{F}(W_{Y},Y)$ are equivalent in $\mathbf{hCoarse}$. 
\end{theorem}
\begin{proof} Suppose we have a coarse homotopy equivalence $\mathfrak{a}:cX\rightarrow cY$ with inverse $\mathfrak{b}: cY\rightarrow cX$. We have to show that $M_{X}$ and $M_{Y}$ are properly homotopy equivalent. To that end, we take the pre-images in the proof of fullness of $\mathcal{F}$ and show that they define a proper homotopy equivalence, ie. that $\psi:=\widetilde{i_{Y}s_{l}\mathfrak{a}R_{X}}: M_{X}\rightarrow M_{Y}$ is a proper homotopy equivalence with inverse $\xi:= \widetilde{i_{X}s_{r}\mathfrak{b} R_{Y}}$. We compute both compositions.
\begin{align*}
\xi\psi=(i_{X}s_{r}\mathfrak{b} R_{Y})^{\sim}\circ ({i_{Y}s_{l}\mathfrak{a}R_{X}})^{\sim}\simeq (i_{X}s_{r}\mathfrak{b} R_{Y} i_{Y}s_{l}\mathfrak{a}R_{X})^{\sim} \simeq (i_{X}s_{r}\mathfrak{b} s_{l}\mathfrak{a}R_{X})^{\sim} \\
\simeq (i_{X}s_{r}\mathfrak{b}\mathfrak{a}R_{X})^{\sim} \simeq (i_{X}s_{r}R_{X})^{\sim} \simeq \id_{M_{X}}
\end{align*}
Note that in this computation, we fix $l$, and can choose $r$ large enough so that all the homotopies exist. \\

The computation for $\psi \xi$ is exactly the same, this time fixing $r$ and choosing a sufficiently large $l$. Therefore $M_{X}$ and $M_{Y}$ are properly homotopy equivalent. So $W_{X}\setminus X$ and $W_{Y}\setminus Y$ are also proper homotopy equivalent. The converse follows from the fact that $\mathcal{F}$ is a functor so sends equivalences to equivalences. 
\end{proof}

For the purpose of studying coarse homology theories, Bunke and Engel define a universal coarsification functor $\mathbf{P}: \mathbf{BornCoarse} \rightarrow \mathbf{Sp}\mathcal{B}$ (Definition $5.1$ in \cite{bunke2020coarse}), where $\mathbf{Sp}\mathcal{B}$ denotes the stable $\infty$-category of closed motivic uniform bornological coarse spectra (Definition $4.14$ in \cite{bunke2020coarse}). $\mathbf{P}$ is u-continuous and invariant under coarse homotopies. Therefore, it is well-defined on $\mathbf{hCoarse}$. The composition
\begin{align*}
\mathbf{StrongShape}\xrightarrow{\mathcal{F}} \mathbf{hCoarse} \xrightarrow{\mathbf{P}} \mathbf{Sp} \mathcal{B}
\end{align*}
is $\mathbf{Sp} \mathcal{B}$-valued local homology theory. \\

We now compare our $\mathcal{F}$ to the universal cone functor $\mathcal{O}:\mathbf{UBC} \rightarrow \mathbf{BornCoarse}$, which takes a uniform, bornological coarse space $X$ to its cone $\mathcal{O}(X)$ and a uniformly continuous map $f: X\rightarrow Y$ to the cone of $f$ which is defined set-wise as $f\times \id$. $\mathcal{O}$ restricts to the subcategory $\mathbf{UC}$, of uniform coarse spaces with bornology induced by the coarse structure, and descends to a well-defined functor $\mathcal{O}: \mathbf{hUC}\rightarrow \mathbf{hCoarse}$. However, this functor is not full. We have discovered that there are many more morphisms between $\mathcal{O}(X)\rightarrow \mathcal{O}(Y)$ (coming from shape maps) than can be obtained from coning a uniformly continuous map $f: X\rightarrow Y$, even when $X,Y$ are compact metric spaces.
\newpage

\part{Proper metric spaces}

We wish to prove and generalise the following theorems in \cite{ashley2025interactions}:

\begin{theorem} \label{surjends} (Theorem $7.4$ in \cite{ashley2025interactions}) Let $X$ be a proper geodesic metric space. There exists a natural surjection $\pi^c_{0}(X)\rightarrow \mathcal{E}nds(X)$. 
\end{theorem}

\begin{theorem}\label{injends} (Theorem $9.9$ in \cite{ashley2025interactions}) Let $X$ be a locally finite geometric tree. There exists a natural isomorphism $\pi^c_{0}(X)\cong \mathcal{E}nds(X)$.
\end{theorem}
\section{$\varprojlim^1$ sequence}

\subsection{End homotopy groups}

Let $X$ be a connected proper metric space. We show $\pi_{n}^{c}(X,\omega)\cong \varinjlim_{h\in \mathbb{N}_{0}}\pi^{e}_{n}(|\mathcal{U}_{h}|, \widetilde{i_{h}\omega})$.\\

Recall that for $C>0$, a subset $Z\subset X$ is $C$-separated if $d(x,y)>C$ for all $x,y\in Z$. By Zorn's lemma, there is a maximal $C$-separated subset for all $C$. Observe that if $X$ is a simplicial complex equipped with the spherical metric, scaled so that edges have length $1$, the set of vertices $Z_{0}$ of $X$ form a maximal $\frac{1}{2}$-separated subset, and the open cover $\cup_{z\in Z_{0}} B(z,1)$ corresponds to the open cover by stars $\cup_{z\in Z_{0}} \langle st(z)\rangle $.\\

Consider $X$ a connected proper metric space. Let $\{Z_{h}\}_{h\in \mathbb{N}_{0}}$ be a maximal $\frac{2^h}{2}$-separated subset, and take the open cover $\mathcal{U}_{h} = \cup_{z\in Z_{h}} B(z,2^h)$. This open cover is locally finite, since $X$ is proper: assume for a contradiction that there was a point $z \in Z_{h}$ such that $B(z,2^{h}) \cap B(w,2^h) \neq \emptyset$ for infinitely many $w$. Since  $B(w,2^h)\subset B(z,2^{h+2})$ for all $w$ and the closure of $B(z,2^{h+2})$ is compact, we can construct a convergent sequence $\{w_{i}\}_{i\in \mathbb{N}}$. But this cannot possibly be Cauchy, since $Z_{h}$ is a $\frac{2^h}{2}$-separated set. Let $|\mathcal{U}_{h}|$ denote the geometric realisation of $|\mathcal{U}_{h}|$.  We equip each simplex in $|\mathcal{U}_{h}|$ with the spherical metric, scaled so that edges have length $1$, and take the induced path metric. $|\mathcal{U}_{h}|$ is a connected, proper geodesic metric space. As a CW complex, it is strongly locally finite (since it is regular and locally finite). There are proper, simplicial maps $\phi_{h,h+1}: |\mathcal{U}_{h}|\rightarrow |\mathcal{U}_{h+1}|$ induced by the refinement of open covers. \\

Let $\rho_{h}= \{\rho_{U}\}_{U\in \mathcal{U}_{h}}$ be a partition of unity subordinate to the cover $\mathcal{U}_{h}$. Define $i_{h}: X\rightarrow |\mathcal{U}_{h}|$ by sending a point $x\in X$ to its image $\sum_{U\in \mathcal{U}_{h}}\rho_{U}(x) [U]$ in barycentric co-ordinates. Define $R_{h}: |\mathcal{U}_{h}| \rightarrow X$ by sending vertices $[B(z,2^{h})]$ to the point $z$, and by sending interiors of simplices $[B(z_0,2^{h}),\dots, B(z_m,2^{h})]$ to a point contained in the intersection $\cap_{i=0}^{m} B(z_i,2^{h})$. 

\begin{notation} From now on we identify points $z \in Z_{h}$ with the vertices $[B(z,2^h)]$ in $|\mathcal{U}_{h}|$. 
\end{notation}

\begin{lemma} Let $X$ be a connected proper metric space and let $h\in \mathbb{N}_{0}$. \begin{enumerate}
\item The map $R_{h}: |\mathcal{U}_{h}| \rightarrow X$ is coarse. 
\item The map $i_{h}: X \rightarrow |\mathcal{U}_{h}|$ is proper. 
\item The compositions $R_{h}\circ i_{h}: X\rightarrow X$ and $i_{h}\circ R_{h}: |\mathcal{U}_{h}|\rightarrow |\mathcal{U}_{h}|$ are close to the identity. 
\item $\phi_{h,h+1} i_{h}$ is close to $i_{h+1}$ and $R_{h+1}\phi_{h,h+1}$ is close to $R_{h}$.
\end{enumerate}
\end{lemma}

\begin{proof} 
\begin{enumerate}
\item The map $R_{h}: |\mathcal{U}_{h}| \rightarrow X$ is controlled: consider the entourage
\begin{align*}
C_{h}:=\{(x,y)\in |\mathcal{U}_{h}| \,|\, x,y \text{ lie in a common simplex }\} \subset D_{1}
\end{align*}
The coarse structure on $|\mathcal{U}_{h}|$ is generated by $C_{h}$. We repeat the argument from the Appendix, for emphasis, on why the spherical metric is important. 

\begin{lemma*} (Lemma \ref{sphericals}) Let $Z$ be a locally finite simplicial complex with a uniform spherical metric. Then the distance from a vertex $v$ of $Z$ to a simplex $\sigma$ in the same component of $Z$ is $\frac{\pi}{2}$ times the length (ie. number of edges) of the shortest simplicial path from $v$ to $\sigma$. In other words there are no shortcuts through the interior of a simplex.
\end{lemma*}

We apply this lemma to the simplicial complex $|\mathcal{U}_{h}|$. Our simplices are scaled so $\pi/2$ is replaced with $1$. Let $x,y$ be two points in $|\mathcal{U}_{h}|$ with distance at most $S\in \mathbb{N}$. Let $x$ be contained in the interior of $\tau_{1}$ and $y$ be contained in the interior of $\tau_2$. Let $x'$ be the closest vertex in $\tau_{1}$ to $x$, and $y'$ the closest vertex in $\tau_2$ to $y$.  We have that 
\begin{align*}
d(x',y')\leq d(y',x) + d(x,x') \leq d(x,y) + d(y,y') + d(x,x') \leq S+1
\end{align*}
$d(x',y')$ is realised by a vertex path of distance at most  $S+1$. Therefore, it travels through at most $S+1$ simplices in $|\mathcal{U}_{h}|$. This shows there is a path between $x$ and $y$ which can be split up into $S+3$ components, each lying within a single simplex. We conclude that $D_{S} \subset C^{\circ (S+3)}_h$. \\

Let $x,y \in |\mathcal{U}_{h}|$ be points lying in a common simplex $\sigma=[z_0,\dots,z_{m}]$. $R_{h}(x), R_{h}(y)$ are contained in the union $\cup_{i=0}^{m}B(z_{i}, 2^h)$, which is a set of diameter $2^{h+2}$, since the intersection $\cap_{i=0}^{m}B(z_{i}, 2^h)$ is non-empty. This shows that $R_{h}$ is controlled.\\

Let $K\subset X$ be a bounded set. Since $X$ is proper, there exist $z_0,\dots, z_{m}\in Z_{h}$ such that $K\subset \cup_{i=0}^{m} B(z_i,2^h)$. We have that $R^{-1}_{h}(\cup_{i=0}^{m} B(z_i,2^h))\subset \cup_{i=0}^{m} |st(z_i)|$: for a point $x\notin \cup_{i=0}^{m} |st(z_i)|$, consider the unique simplex $\tau = [y_0,\dots,y_n]$ such that $x\in int(\tau)$. We have that $R_{h}(x)\subset \cap_{j=0}^{n} B(y_i,2^h)$. If there is an $i$ such that $R_{h}(x)\in B(z_i,2^h)$, then $[y_0,\dots,y_n,z_i]$ is a simplex in $|\mathcal{U}_{h}|$, contradiction. The set $\cup_{i=0}^{m} |st(z_i)|$ has diameter $\leq \max_{i,j=0,\dots,m}\{d(z_i,z_j)\} +2$. This shows that $R_{h}$ is proper. 

\item  The map $i_{h}: X \rightarrow |\mathcal{U}_{h}|$ is proper: the bornology of $|\mathcal{U}_{h}|$ is generated by the metric balls $B(z,S)$, where $z$ is a vertex and $S\in \mathbb{N}$. Suppose that $x \in B(z,S)$. There exists a unique simpliex $\sigma_{x}$ whose interior contains $x$. The shortest simplicial path between $\sigma_{x}$ and $z$ consists of at most $S+1$ edges. Therefore, there is a sequence of open sets $B(z_{0},2^h),\dots, B(z_{m},2^h)$ in $X$ where $z_0 = z$,  $i_{h}^{-1}(x) \subset B(z_{m},2^h)$, $m\leq S+1$, and $B(z_{i},2^h)\cap B(z_{i+1},2^h)\neq \emptyset$ for all $0\leq i< m$. This means that $i_{h}^{-1}(x)$ is at most distance $2^{h+1}(S+1)+ 2^{h}$ from $z\in  i_{h}^{-1}B(z,S)$. This shows that $i_{h}$ is proper. \\

\item The composition $R_{h} \circ i_{h}$ is close to the identity: for a point $x \in B(z,2^h)$ we have $R_{h} \circ i_{h} (x)\subset R_{h} (\langle st(z)\rangle) \subset B(z,2^h)$. The composition $i_{h} \circ R_{h}$ is close to the identity. For a point $x \in \langle st(z) \rangle $ we have $i_{h} \circ R_{h}(x) \subset i_{h}(B(z,2^h)) \subset \langle st(z)\rangle$. 
\item $\phi_{h,h+1} i_{h}$ is close to $i_{h+1}$: let $x\in \cap_{i=0}^m B(z_i,2^h)\subset \cap_{i=0}^m B(\phi_{h,h+1}(z_i),2^{h+1})$. We have that $\phi_{h,h+1} i_{h}(x)\in \phi_{h,h+1} \cap_{i=0}^m \langle st(z_i) \rangle \subset  \cap_{i=0}^m \langle st(\phi_{h,h+1}(z_i)) \rangle$ and $ i_{h+1}(x)\in  \cap_{i=0}^m \langle st(\phi_{h,h+1}(z_i)) \rangle$. \\

$R_{h+1}\phi_{h,h+1}$ is close to $R_{h}$: for $x \in int(\sigma) = [z_0,\dots,z_m]$ we have  $R_{h+1}\phi_{h,h+1}(x)\in \cap_{i=0}^{m}B(\phi_{h,h+1}(z_i),2^{h+1})$ and $R_{h}(x)\in \cap_{i=0}^{m}B(z_i,2^h)\subset \cap_{i=0}^{m}B(\phi_{h,h+1}(z_i),2^{h+1})$. The set $ \cap_{i=0}^{m}B(\phi_{h,h+1}(z_i),2^{h+1})$ has a diameter bounded by $2^{h+2}$. \\

\end{enumerate}

\end{proof}

In the case that $X$ is geodesic, we have the following additional properties. 

\begin{lemma} \label{geoprop} Let $X$ be a connected proper geodesic metric space. Let $h\in \mathbb{N}_{0}$. $i_{h}$ is a coarse equivalence (quasi-isometry) with inverse $R_{h}$. 
$\phi_{h,h+1}$ is a quasi-isometry with inverse $i_{h}\circ R_{h+1}$
\end{lemma}

\begin{proof} It suffices to prove that $i_{h}$ is controlled. Since $X$ is a path metric space, its coarse structure is generated by one entourage $D_{1}:= \{(x,y)\in X\,|\, d(x,y)<\frac{1}{2}\}$. Let $(x,y)\in D_{1}$. There exists a $z\in Z_{h}$ such that $x,y\in B(z,2^h)$. Therefore $i_{h}(x),i_{h}(y) \in \langle st(z) \rangle$, which is a set of diameter $2$. \\

To show that $\phi_{h,h+1}$ is a quasi-isometry, we have
\begin{align*}
i_{h} R_{h+1} \phi_{h,h+1} \simeq_{c} i_{h} R_{h+1} \phi_{h,h+1} i_{h} R_{h}\simeq_{c} i_{h} R_{h+1} i_{h+1} R_{h} \simeq_{c} i_{h} R_{h} \simeq_{c} \id_{X}\\
\phi_{h,h+1} i_{h} R_{h+1}  \simeq_{c} i_{h+1} R_{h+1}\simeq_{c} \id_{|\mathcal{U}_{h+1}|}
\end{align*}
The symbol $\simeq_{c}$ here means that two maps are close. 
\end{proof}

Before proving Theorem \ref{endhomosiso}, let us deal with some technicalities related to the base ray $\omega$. Give $[1,\infty)$ the simplicial structure with edges of length $1$. Since $\omega: [1,\infty)\rightarrow X$ is coarse, there exists a $h_{0}$ large enough such that $\omega[n-1,n+1]$ has diameter less than $\frac{2^{h_0}}{2}$ for all $n\in \mathbb{N}$. Let $h\geq h_0$. We have that $i_{h} \omega [n-1,n+1] \subset \langle st(z_{n,h}) \rangle$ where $z_{n,h}$ denotes the closest vertex in $Z_{h}$ to the point $\omega(n)$ (if there are multiple, choose one arbitrarily). The star condition for $i_{h}\omega$ is met and therefore the vertex map $\widetilde{i_{h}\omega}(n):= z_{n,h}$ can be extended to a simplicial approximation of $i_{h}\omega$. We fix these simplicial approximations and denote them by $\widetilde{i_{h}\omega}$. \\

Assume for the rest of this section that $h\geq h_0$. The proper, simplicial (Lipschitz) maps $\phi_{h,h+1}: |\mathcal{U}_{h}|\rightarrow |\mathcal{U}_{h+1}|$ induce morphisms $(\phi_{h,h+1})_{\ast}: \pi_{n}^e(|\mathcal{U}_{h}|,\widetilde{i_{h}\omega})\rightarrow \pi_{n}^e(|\mathcal{U}_{h+1}|,\phi_{h,h+1}\widetilde{i_{h}\omega})$. Since $\phi_{h,h+1} \widetilde{i_{h}\omega}$ is a simplicial approximation of $i_{h+1} \omega$, there is a proper Lipschitz homotopy $J: c[0,1]\rightarrow |\mathcal{U}_{h+1}|$ between them. We define the bonding maps $\psi_{h,h+1}: \pi^{e}_{n}(|\mathcal{U}_{h}|, \widetilde{i_{h}\omega})\rightarrow \pi_{n}^{e}(|\mathcal{U}_{h+1}|,\widetilde{i_{h+1}\omega})$ as the composition $b^{J}_{\phi_{h,h+1} \widetilde{i_{h}\omega}, \widetilde{i_{h+1}\omega}}(\phi_{h,h+1})_{\ast}$ where $b^{J}_{\phi_{h,h+1} \widetilde{i_{h}\omega}, \widetilde{i_{h+1}\omega}}$ is the change of base ray homomorphism. So $\{\pi_{n}^{e}(|\mathcal{U}_{h}|, \widetilde{i_{h}\omega}), \psi_{h,h+1}; \mathbb{N}_{\geq h_0}\}$ is a directed system with direct limit $\varinjlim_{h\in \mathbb{N}_{\geq h_0}}\pi^{e}_{n}(|\mathcal{U}_{h}|, \widetilde{i_{h}\omega})$. Since all the bonding maps and change of base ray homomorphisms are Lipschitz, $\{\pi_{n}^{L,e}(|\mathcal{U}_{h}|, \widetilde{i_{h}\omega}), \psi_{h,h+1}; \mathbb{N}_{\geq h_0}\}$ is also a directed system of groups with direct limit $\varinjlim_{h\in \mathbb{N}_{\geq h_0}}\pi^{L,e}_{n}(|\mathcal{U}_{h}|, \widetilde{i_{h}\omega})$. 

\begin{theorem} \label{endhomosiso} Let $X$ be a connected proper metric space, $\omega: [1,\infty)\rightarrow X$ a coarse base ray. There is an isomorphism 
\begin{align*}
\pi_{n}^{c}(X,\omega)\cong \varinjlim_{h\in \mathbb{N}_{\geq h_0}}\pi^{L,e}_{n}(|\mathcal{U}_{h}|, \widetilde{i_{h}\omega})
\end{align*} 
\end{theorem}

\begin{proof}

Let $[f]\in \pi_{n}^{c}(X,\omega)$. There exists a $S>0$ such that $f \langle st(v) \rangle$ has diameter $<S$ for all vertices $v \in c[0,1]^n$. We choose $h$ so that $\frac{2^h}{2}>S$. Then for every vertex $v$ there exists a $z_{v} \in Z_{h}$ such that $f \langle st(v) \rangle \subset B(z_{v},2^h)$. Therefore $i_{h}f \langle st(v) \rangle \subset \langle st(z_{v}) \rangle$ and the vertex map $\widetilde{i_{h}f}(v) = z_{v}$ can be extended linearly to obtain a simplicial approximation of $i_{h}f$. We can choose $\widetilde{i_{h}f}$ such that it restricts to $\widetilde{i_{h}\omega}$ on $c\partial[0,1]^n$. We define 
\begin{align*}
 \theta: \pi_{n}^{c}(X,\omega) &\rightarrow \varinjlim_{h\in \mathbb{N}_{\geq h_0}}\pi^{L,e}_{n}(|\mathcal{U}_{h}|, \widetilde{i_{h}\omega})\\
[f] &\mapsto [\widetilde{i_{h} f}]
\end{align*}
Any two choices of simplicial approximation $\widetilde{i_{h} f}$ are properly homotopic by a linear, Lipschitz homotopy. $\theta$ respects the bonding maps $\psi_{h,h+1}$ since $\phi_{h,h+1}\widetilde{i_{h}f}$ and $\widetilde{i_{h+1}f}$ are both simplicial approximations of the map $i_{h+1}f$: there exists a coarse-Lipschitz homotopy $G: c([0,1]^n \times [0,1]) \rightarrow |\mathcal{U}_{h+1}|$ such that $G_{| c([0,1]^n \times \{0\})} = \phi_{h,h+1}\widetilde{i_{h}f}$, $G_{| c([0,1]^n \times \{1\})} =\widetilde{i_{h+1}f}$ and $G_{| c([0,1]^n \times [0,1])} = J$. Therefore $\psi_{h,h+1}[\widetilde{i_{h}f}] = b^{J}_{\phi_{h,h+1} \widetilde{i_{h}\omega}, \widetilde{i_{h+1}\omega}}\phi_{h,h+1}[\widetilde{i_{h}f}] = [\widetilde{i_{h+1}f}]$. \\

Let $[f] = [g]$. There exists a coarse homotopy $H: c([0,1]^n \times [0,1])\rightarrow X$ such that $c([0,1]^n \times \{0\}) = f$, $([0,1]^n \times \{1\}) = g$ and $c(\partial [0,1]^n \times [0,1]) = \omega$. There exists a $h$ large enough such that $i_{h}H$ has a simplicial approximation. The homotopy $\widetilde{i_{h}H}$ is relative to the base ray $\widetilde{i_{h}\omega}$ and restricts to simplicial approximations of $i_{h}f$ and $i_{h}g$ on $c([0,1]^n \times \{0\})$ and $c([0,1]^n \times \{1\})$ respectively. We obtain $[\widetilde{i_{h}f}] = [\widetilde{i_{h}g}]$. This shows that $\theta$ is well-defined.\\

We define the inverse $\lambda:  \varinjlim_{h\in \mathbb{N}_{\geq h_0}}\pi^{L,e}_{n}(|\mathcal{U}_{h}|, \widetilde{i_{h}\omega})  \rightarrow \pi_{n}^{c}(X,\omega)$ as follows. Let $[f]\in \pi^{L,e}_{n}(|\mathcal{U}_{h}|, \widetilde{i_{h}\omega})$. Observe that $R_{h}\widetilde{i_{h}\omega}$ is close to $\omega$. We let $\lambda_{h}[f] = [(R_{h}f)'] \in \pi_{n}^{c}(X,\omega)$ where $(R_{h}f)'$ denotes the map $R_{h}f$ with the boundary replaced by $\omega$. We have that $\lambda_{h+1}\psi_{h,h+1}[f] = \lambda_{h}[f]$ because of the following:
\begin{align*}
(R_{h+1}(b_{\phi_{h,h+1}\widetilde{i_{h}\omega},\widetilde{i_{h+1}\omega}}\phi_{h,h+1}f))'
 \simeq_{c} b^{\id}_{\omega,\omega}(R_{h}f)' \simeq (R_{h}f)' 
\end{align*}
where $\simeq$ denotes coarsely homotopic. The first homotopy comes from the fact that $R_{h+1}\phi_{h,h+1}$ is close to $R_{h}$, $R_{h}i_{h}$ is close to $\id$, and $\phi_{h,h+1} \widetilde{i_{h}\omega}$ is close to $\widetilde{i_{h+1}\omega}$.  Therefore $\lambda$ is well-defined. It is a group homomorphism, since each $\lambda_{h}$ is a composition of group homomorphisms:
\begin{align*}
\pi^{L,e}_{n}(|\mathcal{U}_{h}|, \widetilde{i_{h}\omega}) \xrightarrow{\id_{*}}\pi^{c}_{n}(|\mathcal{U}_{h}|, \widetilde{i_{h}\omega})  \xrightarrow{R_{h}}  \pi_{n}^{c}(X,R_{h}\widetilde{i_{h}\omega})  \xrightarrow{(\cdot)'} \pi_{n}^{c}(X,\omega)
\end{align*}
We have used here that $R_{h}$ is coarse.\\

The composition $\lambda \circ \theta$: let $[f]\in \pi_{n}^{c}(X,\omega)$. We have that $\lambda \theta[f] =  [(R_{h}\widetilde{i_{h}f})']$ for some $h\geq h_0$. The map $R_{h}\widetilde{i_{h}f}$ is close to $f$. Therefore $\lambda \circ \theta = \id$. \\

The composition $\theta \circ \lambda$: let $[f]\in \pi_{n}^{L,e}(|\mathcal{U}_{h}|,\widetilde{i_{h}\omega})$. Since the simplicial structure on $c[0,1]^n$ is uniformly bounded, there exists a small enough barycentric subdivision $bs^k(c[0,1]^n)$ relative to the boundary $c\partial[0,1]^n$, and a Lipschitz map $\mathfrak{h}: c[0,1]^n\rightarrow c[0,1]^n$ such that $f\mathfrak{h}$ has a simplicial approximation, which we denote by $bs^k (\widetilde{f\mathfrak{h}})$. $f$ is coarse-Lipschitz homotopic to $bs^{k}(\widetilde{f\mathfrak{h}})$ relative to $c\partial[0,1]^n$. $\theta \circ \lambda [f]$ is represented by the map 
\begin{align*}
bs^k (i_{h}{(R_{h}bs^k (\widetilde{f\mathfrak{h}}))'})^{\sim}
\end{align*}
where $bs^{k}(\cdot)^{\sim}$ denotes the simplicial approximation with respect to $bs^k(c[0,1]^n)$. The reason why we can use $bs^k(c[0,1]^n)$ instead of the fixed simplicial structure on $c[0,1]^n$ is because the homotopy class of a simplicial approximation is retained under barycentric subdivision relative to a subcomplex. \\

Observe that for a vertex $v$ in $|\mathcal{U}_{h}|$ we have that $i_{h}R_{h} \langle st(v) \rangle\subset \langle st(v) \rangle $. Therefore 
$bs^k (i_{h}{(R_{h}bs^k (\widetilde{f\mathfrak{h}}))'})^{\sim}$ is a simplicial approximation of $bs^k (\widetilde{f\mathfrak{h}})$. This shows that $\theta\circ \lambda = \id$. We conclude that $\theta$ is a group isomorphism with inverse $\lambda$. 

\end{proof}

\begin{lemma} \label{endLipschitz}The identity induces an isomorphism 
\begin{align*}
\varinjlim_{h\in \mathbb{N}_{\geq h_0}}\pi^{L,e}_{n}(|\mathcal{U}_{h}|, \widetilde{i_{h}\omega}) \cong \varinjlim_{h\in \mathbb{N}_{\geq h_0}}\pi^{e}_{n}(|\mathcal{U}_{h}|, \widetilde{i_{h}\omega})
\end{align*}
\end{lemma}
\begin{proof}
The forget map $\mathfrak{d}_{h}: \pi_n^{L,e}(|\mathcal{U}_{h}|,\widetilde{i_{h}\omega})\rightarrow \pi_n^{e}(|\mathcal{U}_{h}|,\widetilde{i_{h}\omega})$ is an isomorphism for every $h\geq h_0$. The proof is analogous to Proposition \ref{natiso}. The only additional points in the argument are as follows. To show surjectivity of $\mathfrak{d}_{h}$, for $k\in \mathbb{N}_{0}$ the image of $c_{[2^k,2^{k+1}]}[0,1]^n$ under $f$ lies in a finite subcomplex of $|\mathcal{U}_{h}|$, so the ordinary relative simplicial approximation theorem (Theorem \ref{relordinarysa}) applies. To show that the homotopy $H$ is proper, we consider a finite filtration $\{L_{j}\}_{j\in \mathbb{N}}$ of $|\mathcal{U}_{h}|$. Consider the double neighbourhood $N^2(L_{j}):=N(N(L_{j},|\mathcal{U}_{h}|), |\mathcal{U}_{h}|)$ of $L_{j}$. Since $f$ is proper, for every $j$, there exists a $l_{j}$ such that $f (c_{2^{l_j}} [0,1]^n) \subset (N^2(L_{j}))^c$. Therefore $H^{-1}(L_{j}) \subset c_{[1,2^{l_j}]}([0,1]^n \times [0,1])$ for all $j\in \mathbb{N}$. For injectivity, we only need the fact that the cover of $|\mathcal{U}_{h}|$ by open stars has Lebesgue number $\frac{1}{4}$. Since $\mathfrak{d}_{h}$ is induced by the identity map, it commutes with the bonding morphisms.
\end{proof} 

\begin{lemma} \label{diff}
There is an isomorphism of directed systems 
\begin{align*}
\{\pi_{n}^{e}(|\mathcal{U}_{h}|,\phi_{h_0,h} \widetilde{i_{h_0}\omega}), (\phi_{h,h+1})_{\ast}; \mathbb{N}_{\geq h_0}\} \cong \{\pi_{n}^{e}(|\mathcal{U}_{h}|,\widetilde{i_{h}\omega}), \psi_{h,h+1}; \mathbb{N}_{\geq h_0}\}
\end{align*}
\end{lemma}

\begin{proof} Let $\phi_{h_0,h}:= \phi_{h-1,h}\circ \dots \circ \phi_{h_0,h_0+1}$. It suffices to show that the diagram 
\[
\begin{tikzcd}
\pi_{n}^{e}(|\mathcal{U}_{h}|, \phi_{h_0,h}\widetilde{i_{h_0}\omega}) \arrow{rrr}{(\phi_{h,h+1})_{\ast}} \arrow{dd}{b^{H^{h}}_{ \phi_{h_0,h}\widetilde{i_{h_0}\omega},\widetilde{i_{h}\omega}}} &&& \pi_{n}^{e}(|\mathcal{U}_{h+1}|, \phi_{h_0,h+1}\widetilde{i_{h_0}\omega})\arrow{dd}{b^{H^{h+1}}_{ \phi_{h_0,h+1}\widetilde{i_{h_0}\omega},\widetilde{i_{h+1}\omega}}} \\
\\
\pi_{n}^{e}(|\mathcal{U}_{h}|, \widetilde{i_{h}\omega})  \arrow{rrr}{b^{J^{h,h+1}}_{ \phi_{h,h+1}\widetilde{i_{h}\omega},\widetilde{i_{h+1}\omega}}(\phi_{h,h+1})_{\ast}} &&&\pi_{n}^{e}(|\mathcal{U}_{h+1}|, \widetilde{i_{h+1}\omega})
\end{tikzcd}
\]
commutes for all $h\geq h_0$. The homotopies $H^{h}, H^{h+1},J^{h,h+1}$ are obtained from the fact that two simplicial approximations are linearly homotopic. \\

To that end, it suffices to show that there exists a homotopy $G: [0,1]^2 \times [1,\infty)\rightarrow |\mathcal{U}_{h+1}|$ such that 
\begin{align*}
G_{|[0,1] \times \{0\} \times [1,\infty)} &= b^{H^{h+1}}_{ \phi_{h_0,h+1}\widetilde{i_{h_0}\omega},\widetilde{i_{h+1}\omega}}\\
G_{|[0,1] \times \{1\} \times [1,\infty)} &= b^{J^{h,h+1}}_{ \phi_{h,h+1}\widetilde{i_{h}\omega},\widetilde{i_{h+1}\omega}} b^{\phi_{h,h+1}H^{h}}_{ \phi_{h,h+1}\phi_{h_0,h}\widetilde{i_{h_0}\omega},\phi_{h,h+1}\widetilde{i_{h}\omega}}
\end{align*}
We give $[0,1]^2 \times [1,\infty)$ the following simplicial structure. The first copy of $[0,1]$ has vertices $\{0,\frac{1}{2},1\}$ (and two edges). The second copy of $[0,1]$ has vertices $\{0,1\}$, and $[1,\infty)$ has vertex set $\mathbb{N}$. Equip $[0,1]^2 \times [1,\infty)$ with the product simplicial structure. Let $\overline{G}: [0,1]^2 \times [1,\infty)$ be defined  on vertices as 
\begin{align*}
\overline{G}(0,0,m) &= \overline{G}(\frac{1}{2},0,m) =\overline{G}(0,1,m) =   \phi_{h_0,h+1} \widetilde{i_{h_0}\omega} (m)\\
\overline{G}(1,0,m) &= \overline{G}(1,1,m)=  \widetilde{i_{h+1}\omega}(m)\\
\overline{G}(\frac{1}{2},1,m) &= \phi_{h,h+1}\widetilde{i_{h}\omega}(m)\\
\end{align*}
for $m\in \mathbb{N}$. See Figure \ref{fig:G}. 

\begin{figure}[H]
\centering
  \centering
  \includegraphics[width=0.5\linewidth]{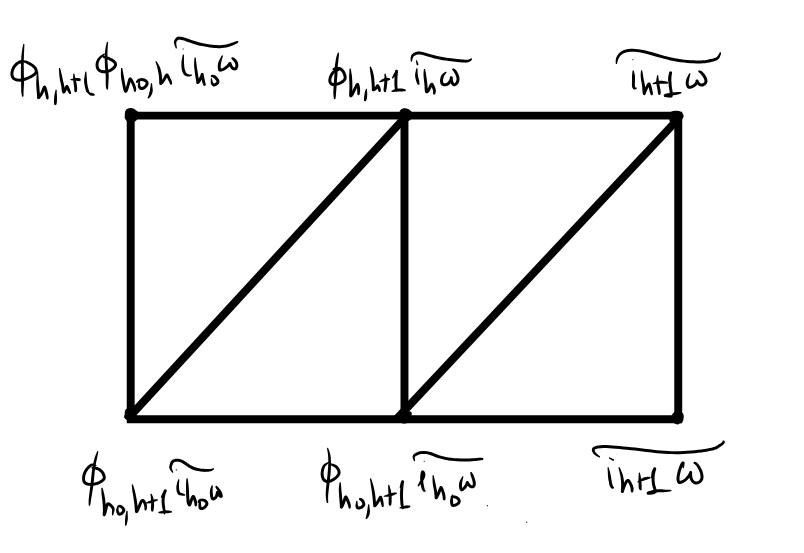}
  \caption{The map $\overline{G}$.}
  \label{fig:G}
\end{figure}
Observe that $\overline{G}$ can be extended linearly to a simplicial map, which is a simplicial approximation of the composition 
\begin{align*}
[0,1]^2 \times [1,\infty) \rightarrow \{\ast\} \times [1,\infty) \xrightarrow{i_{h+1}\omega} |\mathcal{U}_{h+1}|
\end{align*}
Additionally, since the map $\phi_{h,h+1}$ is simplicial, $\overline{G}$ satisfies 
\begin{align*}
\overline{G}_{|[0,\frac{1}{2}] \times \{0\} \times [1,\infty)} &= \id_{\phi_{h_0,h+1}\widetilde{i_{h_0}\omega}}\\
\overline{G}_{|[\frac{1}{2},1] \times \{0\} \times [1,\infty)}&= H^{h+1}\\
\overline{G}_{|[0,\frac{1}{2}] \times \{1\} \times [1,\infty)}&= \phi_{h,h+1} H^{h}\\
\overline{G}_{|[\frac{1}{2},1] \times \{1\} \times [1,\infty)}&= J^{h,h+1}\\
\overline{G}_{|\{0\} \times [0,1] \times [1,\infty)}&= \id_{\phi_{h_0,h+1}\widetilde{i_{h_0}\omega}}\\
\overline{G}_{|\{1\} \times [0,1]  \times [1,\infty)} &= \id_{\widetilde{i_{h+1}\omega}}
\end{align*}
Since $\id_{\phi_{h_0,h+1}\widetilde{i_{h_0}\omega}} \ast H^{h+1} \simeq H^{h+1}$, we can rescale $\overline{G}$ to obtain the desired homotopy $G$. Therefore, the diagram commutes. 
\end{proof}

Assembling together Theorem \ref{endhomosiso}, Lemma \ref{endLipschitz}, and Lemma \ref{diff} we obtain 
\begin{theorem} \label{main1}There is a natural isomorphism 
\begin{align*}
\Theta: \pi_{n}^{c}(X,\omega) &\rightarrow \varinjlim_{h\in \mathbb{N}_{\geq h_0}}\pi^{e}_{n}(|\mathcal{U}_{h}|, \phi_{h_0,h} \widetilde{i_{h_0}\omega})\\
[f] &\mapsto b_{\widetilde{i_{h}\omega},\phi_{h_0,h} \widetilde{i_{h_0}\omega}}[\widetilde{i_{h} f}]
\end{align*}
\end{theorem}

\begin{proof} It only remains to show that this isomorphism is natural. Let $Y$ be a connected, proper metric space with associated sequence of open covers $\{\mathcal{V}_{h}\}_{h\in \mathbb{N}_{0}}$. Let $\mathfrak{a}: (X,\omega)\rightarrow (Y,\mathfrak{a}\omega)$ be a coarse map. Let $h'_{0}\in \mathbb{N}$ be large enough so that $\widetilde{i_{h'_{0}}\mathfrak{a}\omega}$ exists.  We define a group homomorphism $\xi: \varinjlim_{h\in \mathbb{N}_{\geq h_0}} \pi^{e}_{n}(|\mathcal{U}_{h}|, \widetilde{i_{h}\omega}) \rightarrow
\varinjlim_{h'\in \mathbb{N}_{\geq h'_0}} \pi^{e}_{n}(|\mathcal{V}_{h'}|,  \widetilde{i_{h'}\mathfrak{a}\omega})$ which is induced by proper maps $\xi_{h,\rho(h)}: |\mathcal{U}_{h}|\rightarrow |\mathcal{V}_{\rho(h)}|$ for an increasing function $\rho:\mathbb{N}_{\geq h_0}\rightarrow \mathbb{N}_{\geq h'_{0}}$ such that the diagram 

\begin{figure}[H]
\center
\begin{tikzcd}
&\pi_{n}^{c}(X,\omega) \arrow{r}{\mathfrak{a}}  \arrow{d}{\theta}
&\pi_{n}^{c}(Y, \mathfrak{a}\omega) \arrow{d}{\theta} 
\\
&\varinjlim_{h\in \mathbb{N}_{\geq h_0}} \pi^{e}_{n}(|\mathcal{U}_{h}|,  \widetilde{i_{h}\omega}) \arrow{r}{\xi} \arrow{d}{b^{-1}} 
&\varinjlim_{h'\in \mathbb{N}_{\geq h'_0}} \pi^{e}_{n}(|\mathcal{V}_{h'}|, \widetilde{i_{h'}\mathfrak{a}\omega})  \arrow{d}{b'^{-1}} 
  \\ 
&\varinjlim_{h\in \mathbb{N}_{\geq h_0}} \pi^{e}_{n}(|\mathcal{U}_{h}|, \phi_{h_0,h} \widetilde{i_{h_0}\omega}) \arrow{r}{b'^{-1}\xi b}
&\varinjlim_{h\in \mathbb{N}_{\geq h'_0}} \pi^{e}_{n}(|\mathcal{V}_{h'}|, \phi'_{h'_0,h} \widetilde{i_{h'_0}\mathfrak{a}\omega})
\end{tikzcd}
\end{figure}
commutes. \\

Let $h\in \mathbb{N}_{\geq h_0}$. Recall that for a vertex $v\in  |\mathcal{U}_{h}|$, the set $R_{h}\langle st(v) \rangle \subset B(v,2^h)$ has diameter $<2^{h+1}$. Since $\mathfrak{a}$ is coarse we have $d_{X}(x,y)<2^{h+1} \implies d_{Y}(\mathfrak{a}(x),\mathfrak{a}(y)) <S_{h}$. We choose $\rho(h)\geq h'_{0}$ large enough such that $\frac{2^{\rho(h)}}{2}> S_{h}$. Therefore  $i_{\rho(h)} \mathfrak{a} R_{h}$ has a simplicial approximation. Additionally, this simplicial approximation can be chosen to send $\widetilde{i_{h}\omega}$ to $\widetilde{i_{\rho(h)}\mathfrak{a} \omega}$: recall that $\widetilde{i_{h}\omega}(n) = z_{n}$ corresponds to the closest $z_{n}\in Z_{h}$ to $\omega(n)$. By construction, a valid simplicial approximation of $i_{\rho(h)} \mathfrak{a} R_{h}$ sends the vertex $z_n \in |\mathcal{U}_{h}|$ to the closest vertex $w_{y}\in Z'_{\rho(h)}$ to any point $y \in \mathfrak{a} B(z_{n},2^h)$. Obviously the point $y=\mathfrak{a}\omega(n)$ satisfies this property. Therefore, we can define $\xi_{h,\rho(h)}: (|\mathcal{U}_{h}|, \widetilde{i_{h}\omega})\rightarrow (|\mathcal{V}_{\rho(h)}|, \widetilde{i_{\rho(h)}\mathfrak{a}\omega})$ as any choice of base ray-preserving simplicial approximation $(i_{\rho(h)} \mathfrak{a} R_{h})^{\sim}$. We can assume that $\rho: \mathbb{N}\rightarrow \mathbb{N}$ is a strictly increasing sequence.\\

We show now that the sequence of maps $\{\xi_{h,\rho(h)}\}_{h\geq h_0}$ induces a map on the colimit. It suffices to show the commutativity of the diagram

\begin{figure}[H]
\center
\begin{tikzcd}
&\pi^{e}_{n}(|\mathcal{U}_{h}|,  \widetilde{i_{h}\omega}) \arrow{rr}{\xi_{h,\rho(h)}} \arrow{d}{\psi_{h,h+1}} 
&& \pi^{e}_{n}(|\mathcal{V}_{\rho(h)}|, \widetilde{i_{\rho(h)}\mathfrak{a}\omega})  \arrow{d}{\psi'_{\rho(h),\rho(h+1)}} 
  \\ 
& \pi^{e}_{n}(|\mathcal{U}_{h}|,  \widetilde{i_{h+1}\omega}) \arrow{rr}{\xi_{h+1,\rho(h+1)}}
&& \pi^{e}_{n}(|\mathcal{V}_{\rho(h+1)}|,  \widetilde{i_{\rho(h+1)}\mathfrak{a}\omega})
\end{tikzcd}
\end{figure}

We first show that the composition of bonding maps $\psi'_{\rho(h),\rho(h+1)}:= \psi'_{\rho(h+1)-1,\rho(h+1)}\circ \dots \circ  \psi'_{\rho(h),\rho(h)+1}$ is the same as $b_{ \phi'_{\rho(h),\rho(h+1)}\widetilde{i_{\rho(h)}\mathfrak{a}\omega}, \widetilde{i_{\rho(h+1)}\mathfrak{a}\omega}}(\phi'_{\rho(h),\rho(h+1)})_{*}$. We prove this by induction. Assume that for $j\leq k$ the composition of bonding maps $\psi'_{j,k}$ is the same as $b_{ \phi'_{j,k}\widetilde{i_{j}\mathfrak{a}\omega}, \widetilde{i_{k}\mathfrak{a}\omega}}(\phi'_{j,k})_{*}$. Let $[f]\in \pi^{e}_{n}(|\mathcal{V}_{j}|, \widetilde{i_{j}\mathfrak{a}\omega})$.  The class $\psi'_{j,k+1}[f] =  \psi'_{k,k+1} \circ \psi'_{j,k}[f]$ is 
\begin{align*}
b_{ \phi'_{k,k+1} \widetilde{i_{k}\mathfrak{a}\omega},\widetilde{i_{k+1}\mathfrak{a}\omega}} \circ b_{ \phi'_{k,k+1}\phi'_{j,k}\widetilde{i_{j}\mathfrak{a}\omega}, \phi'_{k,k+1} \widetilde{i_{k}\mathfrak{a}\omega}}[\phi'_{k,k+1}\phi'_{j,k}f]
\end{align*} since all maps $\phi'_{k,k+1} \widetilde{i_{k}\mathfrak{a}\omega}, \widetilde{i_{k+1}\mathfrak{a}\omega}, \phi'_{k,k+1}\phi'_{j,k}\widetilde{i_{j}\mathfrak{a}\omega}$ are simplicial approximations of $i_{k+1}\mathfrak{a}\omega$, there is a linear homotopy 
\begin{align*}
b_{ \phi'_{k,k+1} \widetilde{i_{k}a\omega},\widetilde{i_{k+1}\mathfrak{a}\omega}} \circ b_{ \phi'_{k,k+1}\phi'_{j,k}\widetilde{i_{j}\mathfrak{a}\omega}, \phi'_{k,k+1} \widetilde{i_{k}\mathfrak{a}\omega}}\simeq  b_{ \phi'_{j,k+1} \widetilde{i_{j}\mathfrak{a}\omega},\widetilde{i_{k+1}\mathfrak{a}\omega}}
\end{align*}
Therefore $\psi'_{j,k+1}[f] = b_{ \phi'_{j,k+1}\widetilde{i_{j}\mathfrak{a}\omega}, \widetilde{i_{k+1}\mathfrak{a}\omega}}(\phi'_{j,k+1})_{*}[f]$, which proves the claim. \\

Let $[f]\in \pi_{n}^{e}(|\mathcal{U}_{h}|,\widetilde{i_{h}\omega})$. The composition $\psi'_{\rho(h),\rho(h+1)}\xi_{h,\rho(h)}[f]$ is given by 
\begin{align*}
b_{\phi_{\rho(h), \rho(h+1)} (i_{\rho(h)}\mathfrak{a}\omega)^{\sim}, (i_{\rho(h+1)}\mathfrak{a}\omega)^{\sim}}[\phi_{\rho(h), \rho(h+1)}(i_{\rho(h)}\mathfrak{a}R_{h})^{\sim}f]
\end{align*}
The composition $\xi_{h+1,\rho(h+1)}\psi_{h,h+1}[f]$ is given by 
\begin{align*}
b_{(i_{\rho(h+1)}\mathfrak{a}R_{h+1})^{\sim}\phi_{h,h+1}\widetilde{i_{h}\omega},(i_{\rho(h+1)}\mathfrak{a}\omega)^{\sim}}[(i_{\rho(h+1)}\mathfrak{a}R_{h+1})^{\sim}\phi_{h,h+1}f]
\end{align*}

Let $\sigma = [z_0,\dots, z_m]$ be a simplex in $|\mathcal{U}_{h}|$. We have that
\begin{align*}
\phi_{\rho(h), \rho(h+1)}(i_{\rho(h)}\mathfrak{a}R_{h})^{\sim}(z_j) = \phi_{\rho(h), \rho(h+1)}(w_j) = y_j\\
\mathfrak{a}B(z_j,2^h)\subset B(w_j,2^{\rho(h)}) \subset B(y_j,2^{\rho(h+1)})
\end{align*}
 for $0\leq j\leq m$. Similarly, we have 
\begin{align*}
(i_{\rho(h+1)}\mathfrak{a}R_{h+1})^{\sim}\phi_{h,h+1}(z_j) =y'_j\\
\mathfrak{a}B(\phi_{h,h+1}(z_j),2^{h+1})\subset B(y'_j,2^{\rho(h+1)})
\end{align*} 
By definition, the intersection $\cap_{j} B(z_j,2^h)\neq \emptyset$. Since $\cap_{j} B(z_j,2^h) \subset \cap_{j} B(\phi_{h,h+1}(z_j),2^{h+1})$, we have that 
\begin{align*}
\mathfrak{a}(\cap_{j} B(z_j,2^h)) \subset \bigcap_{j} B(y_j,2^{\rho(h+1)}) \cap B(y'_j,2^{\rho(h+1)})
\end{align*} 
is non-empty. Therefore the maps $\phi_{\rho(h), \rho(h+1)}(i_{\rho(h)}\mathfrak{a}R_{h})^{\sim}$ and $(i_{\rho(h+1)}\mathfrak{a}R_{h+1})^{\sim}\phi_{h,h+1}$ are linearly homotopic by a homotopy $H: |\mathcal{U}_{h}| \times [0,1]\rightarrow |\mathcal{V}_{\rho(h+1)}|$. The restriction of $H\circ f$ to the set $c\partial[0,1]^n \times [0,1]$ is a linear homotopy 
\begin{align*}
\phi_{\rho(h), \rho(h+1)} (i_{\rho(h)}\mathfrak{a}\omega)^{\sim}\simeq (i_{\rho(h+1)}\mathfrak{a}R_{h+1})^{\sim}\phi_{h,h+1}\widetilde{i_{h}\omega}  
\end{align*} 
For the base ray, recall that $\widetilde{i_{h}\omega}(n) = z_n$ for $z_n$ the closest vertex to $\omega(n)$. Consider a simplex $[n,n+1]\subset [1,\infty)$. We have that  
\begin{align*}
\mathfrak{a}\omega (n,n+1)\subset \mathfrak{a}(\cap_{j} B(z_j,2^h)) \subset \bigcap_{j} (B(y_j, 2^{\rho(h+1)}) \cap  B(y'_j, 2^{\rho(h+1)}) \cap B(y''_j, 2^{\rho(h+1)}))
\end{align*}
for $j=n,n+1$, where $(i_{\rho(h+1)}\mathfrak{a}\omega)^{\sim}(j) = y''_{j}$. This shows that there is a linear homotopy between the two change of base rays that restricts to $H\circ f$ on the set $\partial U$. This shows that $\psi'_{\rho(h),\rho(h+1)}\xi_{h,\rho(h)}=\xi_{h+1,\rho(h+1)}\psi_{h,h+1}$ so $\xi$ is well-defined. \\

To show that the top square commutes, suppose $[f]\in \pi_{n}^c(X,\omega)$. By construction $\theta[f]$ is represented by $\widetilde{i_{h}f}$ for a $h$ large enough so that $f\langle st(v)\rangle \subset B(z_v,2^h)$ for all vertices $v\in c[0,1]^n$. Since $R_{h}i_{h}B(z_v,2^h)\subset B(z_v,2^h)$, we have that 
$i_{\rho(h)} \mathfrak{a} R_{h}i_{h}f$ has the same simplicial approximation as $i_{\rho(h)} \mathfrak{a} f$. We conclude that 
\begin{align*}
(i_{\rho(h)} \mathfrak{a} R_{h})^{\sim}(i_{h}f)^{\sim} = (i_{\rho(h)} \mathfrak{a} R_{h}i_{h}f)^{\sim} = (i_{\rho(h)} \mathfrak{a} f)^{\sim}
\end{align*}
This shows that $\theta[\mathfrak{a}f]=[\widetilde{i_{\rho(h)}\mathfrak{a}f}] = \xi_{h,\rho(h)}\theta[f] = \xi\theta[f]$. \\

Let us explain what the bottom square means. Recall that we have a sequence of morphisms $b_{h}:=b_{\phi_{h_0,h}\widetilde{i_{h_0}\omega},\widetilde{i_{h}\omega}}: \pi_{n}^{e}(|\mathcal{U}_{h}|, \phi_{h_0,h}\widetilde{i_{h_0}\omega})\rightarrow \pi_{n}^e(|\mathcal{U}_{h}|,\widetilde{i_{h}\omega})$ (resp. $b'_{h}$) which commute with the bonding maps so passes to a map on the colimit, which we denote by $b$ (resp. $b'$). The map $b'^{-1}\xi b$ is defined levelwise as 
\begin{align*}
(b'^{-1}\xi b)_{h}:= b'_{\rho(h)} \xi_{h,\rho(h)} b_{h}=(b_{\phi'_{h'_0,\rho(h)}\widetilde{i_{h'_0}\mathfrak{a}\omega},\widetilde{i_{\rho(h)}\mathfrak{a}\omega}} )^{-1}\xi_{h,\rho(h)}b_{\phi_{h_0,h}\widetilde{i_{h_0}\omega},\widetilde{i_{h}\omega}}
\end{align*}

A small 3D-diagram chase shows that $b'^{-1}\xi b$ passes to the colimit. We have
\begin{align*}
\phi'_{\rho(h),\rho(k)}(b'^{-1}\xi b)_{h} = \phi'_{\rho(h),\rho(k)} (b'_{\rho(h)})^{-1} \xi_{h,\rho(h)} b_{h} = (b'_{\rho(k)})^{-1} \psi'_{\rho(h),\rho(k)}   \xi_{h,\rho(h)} b_{h} \\= (b'_{\rho(k)})^{-1}   \xi_{k,\rho(k)} \psi_{h,k}  b_{h} = (b'_{\rho(k)})^{-1}   \xi_{k,\rho(k)}   b_{k} \phi_{h,k} = (b'^{-1}\xi b)_{k}\phi_{h,k} 
\end{align*}
for $k\geq h$. The bottom square commutes by construction:
\begin{align*}
(b'^{-1}\xi b)_{h} (b_{h})^{-1}= (b'_{\rho(h)})^{-1} \xi_{h,\rho(h)} b_{h}(b_{h})^{-1} =(b'_{\rho(h)})^{-1}\xi_{h,\rho(h)}
\end{align*}

\end{proof}

\subsection{$\pi_0$ case}

Let $\omega: [1,\infty)\rightarrow X$ be the fixed base ray from the previous subsection. We consider the pointed set  $\pi_{0}^c(X,\omega)$ of coarse homotopy classes of coarse base rays with base point $[\omega]$.\\

Let $h\geq h_0$. The proper, simplicial (Lipschitz) maps $\phi_{h,h+1}: |\mathcal{U}_{h}|\rightarrow |\mathcal{U}_{h+1}|$ induce morphisms $(\phi_{h,h+1})_{\ast}: \pi_{0}^e(|\mathcal{U}_{h}|,\widetilde{i_{h}\omega})\rightarrow \pi_{0}^e(|\mathcal{U}_{h+1}|,\phi_{h,h+1}\widetilde{i_{h}\omega}) = \pi_{0}^{e}(|\mathcal{U}_{h+1}|,\widetilde{i_{h+1}\omega})$.  So $\{\pi_{0}^{e}(|\mathcal{U}_{h}|, \widetilde{i_{h}\omega}), (\phi_{h,h+1})_{\ast}; \mathbb{N}_{\geq h_0}\}$ is a directed system with limit $\varinjlim_{h\in \mathbb{N}_{\geq h_0}}\pi^{e}_{0}(|\mathcal{U}_{h}|, \widetilde{i_{h}\omega})$. Since all the bonding maps are Lipschitz and $\phi_{h,h+1}\widetilde{i_{h}\omega}$ is coarse-Lipschitz homotpic to $\widetilde{i_{h+1}\omega}$, $\{\pi_{n}^{L,e}(|\mathcal{U}_{h}|, \widetilde{i_{h}\omega}), (\phi_{h,h+1})_{\ast}; \mathbb{N}_{\geq h_0}\}$ is also a directed system of pointed sets with direct limit $\varinjlim_{h\in \mathbb{N}_{\geq h_0}}\pi^{L,e}_{0}(|\mathcal{U}_{h}|, \widetilde{i_{h}\omega})$.

\begin{theorem} Let $X$ be a connected proper metric space. Let $[\omega]\in \pi_{0}^{c}(X)$ be fixed. There is an isomorphism of pointed sets
\begin{align*}
 \pi_{0}^c(X,\omega)= \varinjlim_{h\in \mathbb{N}_{\geq h_0}} \pi^{L,e}_{0}(|\mathcal{U}_{h}|,\widetilde{i_{h}\omega}) 
\end{align*}
\end{theorem}

\begin{proof}  For $\theta:  \pi_{0}^c(X,\omega)\rightarrow \varinjlim_{h\in \mathbb{N}_{\geq h_0}} \pi^{L,e}_{0}(|\mathcal{U}_{h}|,\widetilde{i_{h}\omega})$, choose the simplicial structure on $[1,\infty)$ with vertices at integer heights. Let $[\tau]\in \pi_{0}^{c}(X)$. As before, there exists a $h\geq h_{0}$ such that $[\tau]\in \pi_0^c(X)$ has a simplicial approximation $\widetilde{i_{h}\tau}\in \pi_0^{L,e}(|\mathcal{U}_{h}|)$. The independence of choice of simplicial approximation, independence of barycentric subdivision, independence of representative in the coarse homotopy class of $\tau$, and commutativity with the bonding maps follow by previous discussion.\\

For the definition of $\lambda$, let $[\tau]\in  \pi^{L,e}_{0}(|\mathcal{U}_{h}|)$. We let $\lambda_{h}[\tau] = [R_{h}\tau]$. Since $R_{h+1} \phi_{h,h+1}$ is close to $R_{h}$, we have 
\begin{align*}
\lambda_{h+1}(\phi_{h,h+1})_{\ast}[\tau] = [R_{h+1}\phi_{h,h+1}\tau] =[R_{h}\tau] = \lambda_{h}[\tau]
\end{align*} Therefore $\lambda:= \varinjlim_{h\in \mathbb{N}_{\geq h_0}} \lambda_{h}$ is well-defined. \\

$\lambda\circ \theta[\tau]$ is represented by $R_{h}\widetilde{i_{h}\tau}$ for some $h\geq h_0$. We have that $R_{h}\widetilde{i_{h}\tau}\simeq_{c} R_{h}i_{h}\tau\simeq_{c}\tau$. This shows that $\lambda \circ \theta = \id$. \\

For $\theta \circ \lambda$, let $[\tau]\in \pi_{0}^{L,e}(|\mathcal{U}_{h}|,\widetilde{i_{h}\omega})$. Choose a small enough barycentric subdivision of $[1,\infty)$ such that $\tau \langle st(v, bs^{k}[1,\infty))\rangle \subset \langle st(z_{v})\rangle$. Recall that $i_{h}R_{h} \langle st(z_v)\rangle\subset \langle st(z_v) \rangle $. We obtain
\begin{align*}
bs^{k}(\widetilde{i_{h}R_{h}\tau})= bs^{k}(\widetilde{\tau}) \simeq \tau
\end{align*}
Hence $\theta\circ \lambda = \id$. \\

\end{proof}

\begin{lemma} The identity induces an isomorphism 
\begin{align*}
\varinjlim_{h\in \mathbb{N}_{\geq h_0}} \pi_0^{L,e}(|\mathcal{U}_{h}|,\widetilde{i_{h}\omega})\cong \varinjlim_{h\in \mathbb{N}_{\geq h_0}}  \pi_0^{e}(|\mathcal{U}_{h}|,\widetilde{i_{h}\omega})
\end{align*}
\end{lemma}

\begin{proof} See Proposition \ref{natiso}. 
\end{proof}

\begin{theorem} \label{main2}There is a natural isomorphism 
\begin{align*}
\theta: \pi_{0}^{c}(X,\omega) &\rightarrow \varinjlim_{h\in \mathbb{N}_{\geq h_0}}\pi^{e}_{0}(|\mathcal{U}_{h}|, \widetilde{i_{h}\omega})\\
[\tau] &\mapsto [\widetilde{i_{h} \tau}]
\end{align*}
\end{theorem}

\begin{proof} Let $Y,\mathfrak{a}, \{\mathcal{V}_{h}\}_{h\in \mathbb{N}_{0}}$ and $\xi_{h,\rho(h)}: |\mathcal{U}_{h}|\rightarrow |\mathcal{V}_{\rho(h)}|$ be as before. For the well-definedness of $\xi$, observe that the diagram 

\begin{figure}[H]
\center
\begin{tikzcd}
&\pi^{e}_{0}(|\mathcal{U}_{h}|,  \widetilde{i_{h}\omega}) \arrow{rr}{\xi_{h,\rho(h)}} \arrow{d}{(\phi_{h,h+1})_{\ast}} 
&& \pi^{e}_{0}(|\mathcal{V}_{\rho(h)}|, \widetilde{i_{\rho(h)}\mathfrak{a}\omega})  \arrow{d}{(\phi'_{\rho(h),\rho(h+1)})_{\ast}} 
  \\ 
& \pi^{e}_{0}(|\mathcal{U}_{h}|,  \widetilde{i_{h+1}\omega}) \arrow{rr}{\xi_{h+1,\rho(h+1)}}
&& \pi^{e}_{n}(|\mathcal{V}_{\rho(h+1)}|,  \widetilde{i_{\rho(h+1)}\mathfrak{a}\omega})
\end{tikzcd}
\end{figure}
commutes. This is because for $[\tau]\in \pi_{0}^{e}(|\mathcal{U}_{h}|,\widetilde{i_{h}\omega})$
\begin{align*}
(\phi_{\rho(h),\rho(h+1)})_{\ast}\xi_{h,\rho(h)}[\tau] = [\phi_{\rho(h),\rho(h+1)}(i_{\rho(h)}\mathfrak{a}R_{h})^{\sim} \tau] = [(i_{\rho(h+1)}\mathfrak{a}R_{h+1})^{\sim}\phi_{h,h+1}\tau] \\
=\xi_{h+1,\rho(h+1)} (\phi_{h,h+1})_{\ast}[\tau] 
\end{align*}\\
Additionally, the diagram 
\begin{figure}[H]
\center
\begin{tikzcd}
&\pi_{0}^{c}(X,\omega) \arrow{r}{\mathfrak{a}}  \arrow{d}{\theta}
&\pi_{0}^{c}(Y, \mathfrak{a}\omega) \arrow{d}{\theta} 
\\
&\varinjlim_{h\in \mathbb{N}_{\geq h_0}} \pi^{e}_{0}(|\mathcal{U}_{h}|,  \widetilde{i_{h}\omega}) \arrow{r}{\xi} 
&\varinjlim_{h'\in \mathbb{N}_{\geq h'_0}} \pi^{e}_{0}(|\mathcal{V}_{h'}|, \widetilde{i_{h'}\mathfrak{a}\omega})  
 
\end{tikzcd}
\end{figure}
commutes: for $[\tau]\in \pi_{n}^{c}(X,\omega)$ we have
\begin{align*}
\xi \theta[\tau] = [(i_{\rho(h)}\mathfrak{a}R_{h})^{\sim}\widetilde{i_{h}\tau}] =[(i_{\rho(h)}\mathfrak{a}\tau)^{\sim}] = \theta \mathfrak{a}[\tau]
\end{align*} 
\end{proof}

\subsection{$\varprojlim^1$ sequence}

We prove the following theorem in this subsection:

\begin{theorem} \label{mainthm}Let $X$ be a connected proper metric space and $[\omega]\in \pi_{0}^{c}(X)$ a fixed base ray such that $\widetilde{i_{h}\omega}$ exists for all $h\geq h_0$.  Let $\mathcal{L}_{h}:=\{L_{j,h}\}_{j\in \mathbb{N}_{0}}$ be a finite filtration of $|\mathcal{U}_{h}|$ such that $\widetilde{i_{h}\omega}$ is well-parametrised with respect to $\mathcal{L}_{h}$. Let $n\geq 1$. There is a short exact sequence of groups
 \begin{align*}
0\rightarrow  \varinjlim_{h} {\varprojlim_{j}}^{1} \pi_{n+1}(L_{j,h}^c,\widetilde{i_{h}\omega}(j))  \rightarrow  \pi_{n}^c(X,\omega) \rightarrow  \varinjlim_{h} \varprojlim_{j} \pi_{n}(L_{j,h}^c,\widetilde{i_{h}\omega}(j)) \rightarrow 0
\end{align*}
For $n=0$, there is a short exact sequence of pointed sets
\begin{align*}
\varinjlim_{h} {\varprojlim_{j}}^{1} \pi_{1}(L_{j,h}^c,\widetilde{i_{h}\omega}(j))  \hookrightarrow \pi^{e}_{0}(Y,\omega) \twoheadrightarrow  \varinjlim_{h} \varprojlim_{j} \pi_{0}(L_{j,h}^c,\widetilde{i_{h}\omega}(j))
\end{align*}
\end{theorem}

We first prove Theorem \ref{mainthm} for $n\geq 1$.  Fix $h\in \mathbb{N}_{\geq h_0}$. Let $\mathcal{L}_{h}:=\{L_{j,h}\}_{j\in \mathbb{N}}$ be a finite filtration of $|\mathcal{U}_{h}|$ such that $\phi_{h_0,h}\widetilde{i_{h_0}\omega}$ is well-parametrised with respect to $\mathcal{L}_{h}$. There is a natural $\varprojlim^{1}$ sequence for each $h$ (Theorem \ref{lim1absolutegeq1}), and we obtain the commutative diagram

\begin{figure}[H]
 \advance\leftskip-3.5cm
\begin{tikzcd}
&
& \dots 
& \dots 
& \dots 
&
\\
& 0 \arrow{r}
&\varprojlim^1 \pi_{n+1}(L_{j,h+1}^c,\phi_{h_0,h+1}\widetilde{i_{h_0}\omega}(j)) \arrow{r} \arrow{u}
&\pi^{e}_{n}(|\mathcal{U}_{h+1}|, \phi_{h_0,h+1}\widetilde{i_{h_0}\omega}) \arrow{u} \arrow{r}
&\varprojlim \pi_{n}(L_{j,h+1}^c,\phi_{h_0,h+1}\widetilde{i_{h_0}\omega}(j)) \arrow{u} \arrow{r}
& 0
  \\
& 0 \arrow{r}
&\varprojlim^1 \pi_{n+1}(L_{j,h}^c,\phi_{h_0,h}\widetilde{i_{h_0}\omega}(j)) \arrow{r} \arrow{u}{(\phi_{h,h+1})_{\ast}}
&\pi^{e}_{n}(|\mathcal{U}_{h}|, \phi_{h_0,h}\widetilde{i_{h_0}\omega}) \arrow{u}{(\phi_{h,h+1})_{\ast}}\arrow{r}
&\varprojlim \pi_{n}(L_{j,h}^c,\phi_{h_0,h}\widetilde{i_{h_0}\omega}(j)) \arrow{u}{(\phi_{h,h+1})_{\ast}}\arrow{r}
& 0
  \\
& 
&\dots \arrow{u}
& \dots \arrow{u}
& \dots \arrow{u}
& 
\\
& 0 \arrow{r}
& \varprojlim^1 \pi_{n+1}(L_{j,h_0}^c,\widetilde{i_{h_0}\omega}(j)) \arrow{r} \arrow{u}{(\phi_{h_0,h_0+1})_{\ast}}
& \pi^{e}_{n}(|\mathcal{U}_{h_0}|, \widetilde{i_{h_0}\omega})\arrow{r} \arrow{u}{(\phi_{h_0,h_0+1})_{\ast}}
&\varprojlim \pi_{n}(L_{j,h_0}^c,\widetilde{i_{h_0}\omega}(j))  \arrow{r} \arrow{u}{(\phi_{h_0,h_0+1})_{\ast}}
& 0
\end{tikzcd}
\end{figure}

since direct limits preserve exactness, we have a short exact sequence 
\begin{align*}
0\rightarrow  \varinjlim_{h} {\varprojlim_{j}}^{1} \pi_{n+1}(L_{j,h}^c,\phi_{h_0,h}\widetilde{i_{h_0}\omega}(j))  \rightarrow  \pi_{n}^c(X,\omega) \rightarrow  \varinjlim_{h} \varprojlim_{j} \pi_{n}(L_{j,h}^c,\phi_{h_0,h}\widetilde{i_{h_0}\omega}(j)) \rightarrow 0
\end{align*}
For $n=0$, Theorem \ref{lim10} give us the commutative diagram

\begin{figure}[H]
\begin{tikzcd}
& \dots 
& \dots 
& \dots 
\\
&\varprojlim^1 \pi_{1}(L_{j,h+1}^c,\widetilde{i_{h+1}\omega}(j)) \arrow{r} \arrow{u}
&\pi^{e}_{0}(|\mathcal{U}_{h+1}|, \widetilde{i_{h+1}\omega}) \arrow{u} \arrow{r}
&\varprojlim \pi_{0}(L_{j,h+1}^c,\widetilde{i_{h+1}\omega}(j)) \arrow{u} 
  \\
&\varprojlim^1 \pi_{1}(L_{j,h}^c,\widetilde{i_{h}\omega}(j)) \arrow{r} \arrow{u}{(\phi_{h,h+1})_{\ast}}
&\pi^{e}_{0}(|\mathcal{U}_{h}|, \widetilde{i_{h}\omega}) \arrow{u}{(\phi_{h,h+1})_{\ast}}\arrow{r}
&\varprojlim \pi_{0}(L_{j,h}^c,\widetilde{i_{h}\omega}(j)) \arrow{u}{(\phi_{h,h+1})_{\ast}}
  \\
&\dots \arrow{u}
& \dots \arrow{u}
& \dots \arrow{u}
\\
& \varprojlim^1 \pi_{1}(L_{j,h_0}^c,\widetilde{i_{h_0}\omega}(j)) \arrow{r} \arrow{u}{(\phi_{h_0,h_0+1})_{\ast}}
& \pi^{e}_{0}(|\mathcal{U}_{h_0}|, \widetilde{i_{h_0}\omega})\arrow{r} \arrow{u}{(\phi_{h_0,h_0+1})_{\ast}}
&\varprojlim \pi_{0}(L_{j,h_0}^c,\widetilde{i_{h_0}\omega}(j))  \arrow{u}{(\phi_{h_0,h_0+1})_{\ast}}
\end{tikzcd}
\end{figure}

Injectivity, surjectivity, and exactness are preserved under taking the direct limit. Therefore, we obtain the short exact sequence 
\begin{align*}
  \varinjlim_{h} {\varprojlim_{j}}^{1} \pi_{1}(L_{j,h}^c,\widetilde{i_{h}\omega}(j))  \hookrightarrow  \pi_{0}^c(X,\omega) \twoheadrightarrow  \varinjlim_{h} \varprojlim_{j} \pi_{0}(L_{j,h}^c,\widetilde{i_{h}\omega}(j)) 
\end{align*}
where $\varinjlim_{h}\varprojlim_{j} \pi_{0}(L_{j,h}^c,\widetilde{i_{h}\omega}(j)) = \varinjlim_{h}\mathcal{E}nds(|\mathcal{U}_{h}|, \widetilde{i_{h}\omega})$. \\

To recover Theorem \ref{surjends} as a corollary, we restrict to the case where $X$ is geodesic.  Since the maps $\phi_{h,h+1}$ are quasi-isometries, and the number of ends of a space is a quasi-isometry invariant, the direct limit on the right collapses to $h=h_{0}$. Moreover, $i_{h_{0}}$ is also a quasi-isometry, so we obtain 
\begin{align*}
 \varinjlim_{h} \varprojlim_{j} \pi_{0}(L_{j,h}^c,\widetilde{i_{h}\omega}(j)) \cong \varprojlim_{j} \pi_{0}(L_{j,h_{0}}^c,\widetilde{i_{h_{0}}\omega}(j)) \cong \mathcal{E}nds(X,\omega)
\end{align*}
where $\mathcal{E}nds(X,\omega)$ denotes the pointed set $\mathcal{E}nds(X)$ with the end defined by $\omega$ as a base point.  This gives us the surjection $ \pi_{0}^c(X,\omega)\rightarrow \mathcal{E}nds(X,\omega)$ as claimed.\\

\newpage
\section{Examples}

We prove Theorem \ref{injends} now. It suffices to show that if $X$ is a locally finite geometric tree, then $ \varinjlim_{h} \varprojlim^1_{j} \pi_{1}(L_{j,h}^c,\widetilde{i_{h}\omega}(j))$ vanishes. In fact, we show quite explicitly that for each $h\geq h_0$ and any finite filtration $\{L_{j,h}\}_{j\in \mathbb{N}}$ of $|\mathcal{U}_{h}|$,  $\pi_{1}(L_{j,h}^c,\widetilde{i_{h}\omega}(j)) = 0$ for all $j\in \mathbb{N}$.\\

\begin{proof} (of Theorem \ref{injends})\\

The CW-complement of $L_{j,h}$ is the subcomplex generated by the vertices not in $L_{j,h}$, ie. a simplex in $L^{c}_{j,h}$ exists if and only if all its vertices are not in $L_{j,h}$. Consider $[\gamma]\in \pi_{1}(L_{j,h}^c,\widetilde{i_{h}\omega}(j))$. By the CW approximation theorem, we can assume that after changing basepoints, $\gamma$ is a closed edge path in $L_{j,h}^c$. Let $(z_0,\dots,z_m)$ be the ordered tuple of vertices traversed by $\gamma$ with $z_0=z_m$. So there exist metric balls $B(z_0,2^h),\dots, B(z_m, 2^h)$ in $X$ with $B(z_i, 2^h)\cap B(z_{i+1},2^h) \neq \emptyset$ for $0\leq i\leq m-1$. \\

We show by induction that $\gamma$ must be contractible in $L_{j,h}^c$. Assume that $m\leq 3$. The cases $m=1,2$ are clear. For $m=3$, we have $B(z_0,2^{h}), B(z_1,2^{h}),B(z_2,2^{h})$ with 
\begin{align*}
B(z_0,2^{h})\cap B(z_1,2^{h}) \neq \emptyset\\
B(z_1,2^{h})\cap B(z_2,2^{h}) \neq \emptyset\\
B(z_2,2^{h})\cap B(z_0,2^{h}) \neq \emptyset
\end{align*}
Choose $x_0 \in B(z_0,2^{h})\cap B(z_1,2^{h}), x_1 \in B(z_1,2^{h})\cap B(z_2,2^{h}), x_2 \in B(z_2,2^{h})\cap B(z_0,2^{h})$. Let $[x_0,x_1],[x_1,x_2],[x_2,x_0]$ denote the unique geodesics in $X$ betweeen the respective points. A geodesic triangle in a tree is a tripod with a unique midpoint 
\begin{align*}
m([x_0,x_1],[x_1,x_2],[x_2,x_0])
\end{align*}    
Since metric balls in a tree are convex, we have that $[x_0,x_1]\subset B(z_1,2^h), [x_1,x_2] \subset B(z_2,2^h), [x_2,x_0]\subset B(z_0,2^h)$. Therefore $m\in B(z_0,2^{h})\cap B(z_1,2^{h})\cap B(z_2,2^{h})$. This shows that the $2$-simplex $[z_0,z_1,z_2]$ exists in $L^c_{j,h}\subset |\mathcal{U}_{h}|$. Therefore $\gamma$ is a contractible loop. \\

Suppose that $m\geq 4$. Assume that all edge paths with length less than $m$ are contractible. Metric balls in trees are subtrees. We apply Gavril's chordal graph theorem (\cite{gavril1974intersection}) to the intersection graph $G$ of $B(z_0,2^h),\dots, B(z_m, 2^h)$. The loop $\gamma$ corresponds to a simple circuit in $G$ with more than three vertices. This means that there is a chord connecting two non-consecutive vertices, ie. there exist $0\leq i<j\leq m-1$ with $B(z_i,2^h)\cap B(z_j,2^h) \neq \emptyset$ and $|i-j|>1$ (we identify $m=0$). $\gamma$ is homotopic to the concatenation of segments $\tau_{1}\ast \eta \ast \tau_{2}$ where
\begin{align*}
\tau_{1} = (z_0,\dots,z_i)\\
\eta = (z_i, z_j, z_{j-1},\dots, z_i)\\
\tau_{2} = (z_i,z_j,z_{j+1}, \dots, z_m)
\end{align*} \\
$\eta$ is a loop of length $<m$, so it is contractible. The remaining segments $\tau_{1}\ast \tau_{2}$ concatenate to form a loop of length $<m$ as well. This shows that all edge paths of length $m$ are contractible. By induction, the statement holds. 
\end{proof}

This recovers the statement that $\pi_{0}^c(X,\omega)\cong \mathcal{E}nds(X,\omega)$ for $X$ a locally finite geometric tree. In fact, something even stronger is true: $L_{j,h}^c$ is homotopy equivalent to a set of points for all $j,h\in \mathbb{N}$. Let $Z^c_{j,h}$ denote the set of vertices that are not in $L_{j,h}$. Observe that $L_{j,h}^c$ is the nerve of the open cover 
\begin{align*}
\bigcup_{z \in Z^c_{j,h}} B(z,2^h)
\end{align*}
We show that this cover is good. Suppose that $\cap_{i=0}^{m} B(z_{i},2^h) \neq \emptyset$ for $z_{0},\dots, z_{m}\in Z^c_{j,h}$. Let $x,y \in \cap_{i=0}^{m} B(z_{i},2^h)$. Since metric balls in a tree are convex, the geodesic $[x,y]\subset \cap_{i=0}^{m} B(z_{i},2^h)$. Thus $\cap_{i=0}^{m} B(z_{i},2^h)$ is a connected subgraph of $X$, and therefore a tree, which is contractible. By the nerve theorem (see the proof of Corollary \ref{good} or Chapter $4.$G of \cite{hatcher2002algebraic}), $L_{j,h}^c$ is homotopy equivalent to $\cup_{z \in Z^c_{j,h}} B(z,2^h)$. \\

$\cup_{z \in Z^c_{j,h}} B(z,2^h)$ is a forest: any connected component is a subgraph of $X$, and therefore a tree. This shows that  $L_{j,h}^c$ is homotopy equivalent to the set of path components of $\cup_{z \in Z^c_{j,h}} B(z,2^h)$. We conclude: 

\begin{lemma} \label{geomtree} Let $X$ be a locally finite geometric tree, $\omega: [1,\infty)\rightarrow X$ a coarse base ray. The coarse homotopy groups of $X$ are

\[ \pi_{n}^{c}(X,\omega) \cong \begin{cases} 
        0 & n\geq 1 \\ 
        \mathcal{E}nds(X,\omega) & n=0\\
       \end{cases}
    \]
\end{lemma}

 There is also an example constructed in \cite{ashley2025interactions} where  $\pi^c_{0}(X)\rightarrow \mathcal{E}nds(X)$ is not injective. We attempt now to connect this example to our framework. 

\begin{example} Let $X$ be the proper geodesic space defined by a pair of geodesic rays constructed by consecutive edges, connected by a distinguished vertex at one end, and $n^2$ edges connecting the $n$-th vertex on each ray from the joined end, equipped with the usual path metric. We write points in $X$ as $(l,x)$, where $l\in [0,\infty)$ denotes the height and $x\in [0,l^2]$ denotes the position along the rungs of the ladder. \\

\begin{figure}
\centering
  \includegraphics[width=.5\linewidth]{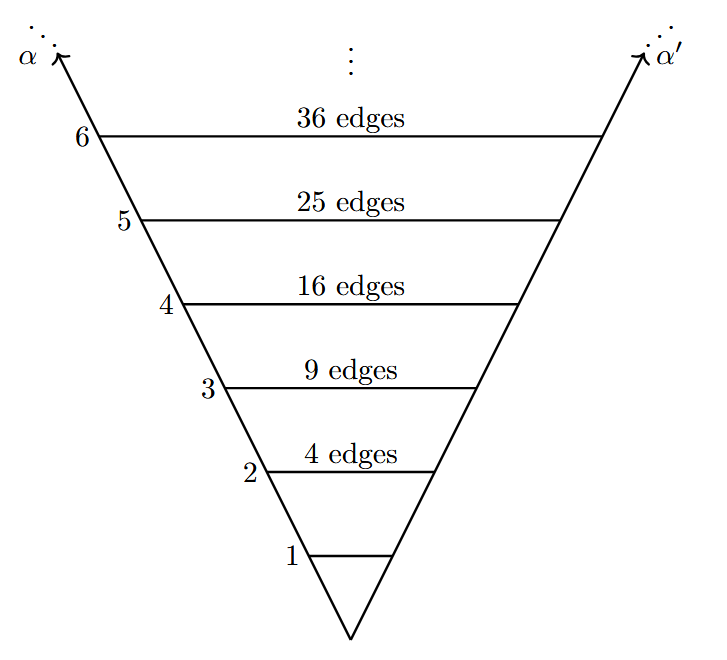}
  \caption{The simplicial complex $X$.}
  \label{fig:notinjexample}
\end{figure}

Let $\alpha$ be the left ray in the graph rooted at the join-point, and similarly for $\alpha'$ the right ray. $X$ has one end, so $\alpha$ and $\alpha'$ define the same end. However, $[\alpha]\neq [\alpha']\in \pi^c_{0}(X)$. 
\end{example}

It is proven in Theorem $7.7$ of \cite{ashley2025interactions} that $[\alpha]\neq [\alpha']\in \pi^c_{0}(X)$ by assuming the existence of a coarse homotopy between $\alpha$ and $\alpha'$ and deriving a contradiction. Using our results, we have an exact algebraic characterisation of the kernel of $\pi^c_{0}(X)\rightarrow \mathcal{E}nds(X)$, namely, it is 
\begin{align*}
\varinjlim_{h} {\varprojlim_{j}}^{1} \pi_{1}(L_{j,h}^c,\widetilde{i_{h}\omega}(j))
\end{align*}
For now, let us focus only on the term $h=0$. In this case, $|\mathcal{U}_{0}| = X$, and a finite filtration is $\{X_{j}\}_{j\in \mathbb{N}_{0}}$, where $X_{j}$ is the subcomplex of $X$ with height $\leq j-1$. $\widetilde{i_{0}\omega}$ is just $\alpha$. We wish to compute
\begin{align*}
{\varprojlim_{j}}^{1} \pi_{1}(X^{c}_{j},\alpha(j))
\end{align*}
where $X^{c}_{j}$ denotes the subcomplex of $X$ with height $\geq j$, which is the CW-complement of $X_{j}$ in $X$ (we shift the index $j$ by $1$ and consider rays $[0,\infty)\rightarrow X$). Let $\gamma_{j}: [0,1]\rightarrow X_{j}$ be the loop based at $\alpha(j)$ which traverses the horizontal steps $j$ and $j+1$ in an anti-clockwise fashion. The pre-image of $[\alpha'] \in \pi^{e}_{0}(X,\alpha)$ under the map  $\varprojlim^1_{j}  \pi_{1}(X^{c}_{j},\alpha(j)) \rightarrow  \pi_{0}^{e}(X,\alpha)$ is the equivalence class of 
\begin{align*}
(\gamma_{j})_{j\in \mathbb{N}_{0}} \in \prod_{j\in \mathbb{N}_{0}} \pi_{1}(X^{c}_{j},\alpha(j))
\end{align*}
This is because the horizontal components cancel out, and we are left with a copy of $\alpha'$, and a copy of $\alpha$ traversed in the opposite direction. It can be checked (Chapter $16.1$ of \cite{geoghegan2007topological}) that this exactly corresponds to $[\alpha'] \in \pi_{0}^{e}(X,\alpha)$. \\

In fact, we can just compute $\varprojlim^1_{j}  \pi_{1}(X^{c}_{j},\alpha(j))$. The fundamental group of $X^{c}_{j}$ is the free group on the generating set $\eta_{k}:= \alpha_{[j,k]} \ast \gamma_{k} \ast (\alpha_{[j,k]})^{-1}$ (traversing the $k$-th circle) which can be identified with the generating set $\mathbb{N}_{0}\setminus \{0,\dots, j-1\} = (\mathbb{N}_{0})_{\geq j}$. The bonding $\psi_{j+1,j}$ maps are induced by the inclusion on spaces which can be identified with the inclusion  
\begin{align*}
F_{(\mathbb{N}_{0})_{\geq j+1}} &\rightarrow F_{(\mathbb{N}_{0})_{\geq j}}
\end{align*}

$\varprojlim^1_{j} \pi_{1}(X^{c}_{j},\alpha(j))$ is the cokernel of the group homomorphism 
\begin{align*}
\varphi: \prod_{j\in \mathbb{N}_{0}} F_{(\mathbb{N}_{0})_{\geq j}} &\rightarrow \prod_{j\in \mathbb{N}_{0}} F_{(\mathbb{N}_{0})_{\geq j}} \\
(a_{j})_{j\in \mathbb{N}_{0}} &\mapsto (a_{j}(\psi_{j+1,j}(a_{j+1}))^{-1})_{j\in \mathbb{N}_{0}}
\end{align*}

$(\gamma_{j})_{j \in \mathbb{N}_{0}}$ is a nontrivial element in the cokernel: suppose not, ie. $\gamma_{j} =  a_{j}(\psi_{j+1,j}(a_{j+1}))^{-1}$ for some $(a_{j})_{j\in \mathbb{N}_{0}}$. To derive a contradiction, we show that $a_{0}$ must traverse each circle exactly once, which cannot happen because the interval is compact. To that end, fix $k\in \mathbb{N}_{0}$. We have that $a_{k}$ must traverse the $k$-th circle anti-clockwise exactly once, since $(\psi_{k+1,k}(a_{k+1}))^{-1}$ is an element of $F_{\mathbb{N}_{\geq k+1}}\subset F_{\mathbb{N}_{\geq k}}$. From the equation $\gamma_{k-1} = a_{k-1} (\psi_{k,k-1}(a_{k}))^{-1}$ we obtain that $a_{k-1}$ traverses the $k$-th circle anti-clockwise exactly once. Going backwards inductively, we obtain that $a_{0}$ traverses the $k$-th circle exactly once. This is true for all $k\in \mathbb{N}_{0}$, so we are done. \\

It remains to show that $(\gamma_{j})_{j \in \mathbb{N}_{0}}$ is also a non-trivial element in the colimit. To do this, we prove that the bonding maps $\phi_{h,h+1}$ induce an isomorphism for all $h\in \mathbb{N}_{0}$.\\

\begin{lemma} Let $X$ be a strongly locally finite, path-connected CW complex, $K$ a compact subcomplex. Let $\omega: [1,\infty)\rightarrow K^c$ be a base ray. The inclusion map induces an isomorphism 
\begin{align*}
\iota: \pi_{n}^{e}(K^c,\omega) \rightarrow  \pi_{n}^{e}(X,\omega)
\end{align*}
with inverse
\begin{align*}
\rho: \pi_{n}^{e}(X,\omega) &\rightarrow   \pi_{n}^{e}(K^c,\omega)\\
[g] &\mapsto b_{\omega_{j},\omega}[g_{j}] 
\end{align*}
where $g_{j}$ is the restriction of $g$ to $[0,1]^n \times [j,\infty)$ for $j$ sufficiently large, and $\omega_{j}(h)= \omega(j+h-1)$. 
\end{lemma}

This lemma further justifies the use of the word "end" to describe $\pi_{n}^{e}$: we can delete any compact subcomplex without affecting the result. 

\begin{proof} We show that $\rho$ is well defined. Suppose we have $k>j$. Consider the homotopy  
\begin{align*}
H:  [0,1]^n \times [1,\infty) \times [0,1] &\rightarrow X \\
H(x,h,t) &\mapsto g(x,(1-t)j+tk+h-1)
\end{align*}
This is a homotopy between $g_{j}$ and $g_{k}$ and restricted to the boundary it is $b_{\omega_{j},\omega_{k}}$. We obtain 
\begin{align*}
b_{\omega_{j},\omega} [g_{j}] = b_{\omega_{j},\omega}b_{\omega_{k},\omega_{j}}[g_{k}] = b_{\omega_{k},\omega} [g_{k}]
\end{align*}
It is obviously independent of the representative of $[g]$ and is a group homomorphism. The composition $\rho \circ \iota$ is the identity on representatives, and for a $[g]\in \pi_{n}^{e}(X,\omega)$ we have $[g] = b_{\omega_{j},\omega} [g_{j}]$. Therefore $\iota \circ \rho = \id$. 
\end{proof}

\begin{lemma} Let $X,Y$ be strongly locally finite, path-connected CW complexes. Let $\psi :X\rightarrow Y$ be a proper map and let $K,L$ be finite subcomplexes of $X,Y$ such that $\psi(K^c)\subset L^c$. Let $\omega:[1,\infty)\rightarrow X$ be a base ray and $j\in \mathbb{N}$ such that $\omega[j,\infty)\subset K^c$. The following diagram commutes

\begin{figure}[H]
\center
\begin{tikzcd}
&\pi_{n}^{e}(K^c,\omega_{j}) \arrow{r}{\psi_{|K^c}}  \arrow{d}{\iota}
&\pi_{n}^{e}(L^c, \psi \omega_{j}) \arrow{d}{\iota} 
  \\ 
& \pi_{n}^{e}(X,\omega_{j}) \arrow{d}{b_{\omega_{j},\omega}} \arrow{r}{\psi}
& \pi_{n}^{e}(Y,\psi\omega_{j})  \arrow{d}{b_{\psi\omega_{j},\psi\omega}} 
\\
& \pi_{n}^{e}(X,\omega) \arrow{r}{\psi} 
&\pi_{n}^{e}(Y,\psi\omega) 
\end{tikzcd}
\end{figure}
\end{lemma}

\begin{proof} The top square obviously commutes. The bottom square commutes since the change of base ray homomorphism is natural. 
\end{proof}

Since the vertical arrows are isomorphisms, this shows that a proper map $\psi: X\rightarrow Y$ induces an isomorphism on end homotopy groups if and only if $\psi_{|K^c}$ is an isomorphism on end homotopy groups for any compact subcomplex $K$. \\

We apply the lemma to the bonding maps $\phi_{h,h+1}: |\mathcal{U}_{h}|\rightarrow |\mathcal{U}_{h+1}|$. It suffices to show that there exist compact subcomplexes $K_{0} \subset |\mathcal{U}_{0}| = X$ and $K_{h} \subset |\mathcal{U}_{h}|$ such that $i_{h}(K^c_{0}) = K^c_{h}$ is a proper homotopy equivalence for all $h\in \mathbb{N}$. This is because $\phi_{h,h+1}i_{h}$ and $i_{h+1}$ have the same simplicial approximation and so are properly homotopic: 
\begin{align*}
\phi_{h,h+1}i_{h} \simeq (\phi_{h,h+1}i_{h})^{\sim}\simeq \widetilde{ i_{h+1}} \simeq  i_{h+1}
\end{align*}
Therefore if $i_{h}$ induces isomorphisms on end homotopy groups for all $h\in \mathbb{N}_{0}$, so do the $\phi_{h,h+1}$. \\

To that end, fix a $h\in \mathbb{N}$. We construct a $\frac{2^h}{2}$-separated subset for $X$. Let $j_{h} = 2^{h+1}$. We have that $5\cdot 2^{h} <(j_{h})^2$. Let
\begin{align*}
c_{l} = \min\{
|l-2^h\lfloor \frac{l}{2^h}\rfloor|,
|l-2^h\lceil \frac{l}{2^h}\rceil|
\}
\end{align*} 
denote the distance to the closest multiple of $2^h$. Consider the vertices 
\begin{align*}
Z'_{h}:= L \bigcup R \bigcup \cup_{l \in \mathbb{N}_{\geq j_{h}}} L_{l} \bigcup \cup_{l \in \mathbb{N}_{\geq j_{h}}}R_{l}\\
:= \{(k2^h,0)\}_{k\in \mathbb{N}}\cup \{(k2^h,(k2^h)^{2})\}_{k\in \mathbb{N}} \cup \{(l,m2^h-c_{l})\}_{l \in \mathbb{N}_{\geq j_{h}}}\cup  \{(l,l^2-m2^h+c_{l})\}_{l \in \mathbb{N}_{\geq j_{h}}}\\
\end{align*}
where $m\in \{1,2\}$. $L,R$ correspond to vertices on the rays $\alpha,\alpha'$ respectively; $L_{l}$ consists of two points on the rung of height $l$ and analogously for $R_{l}$. Observe that $L_{l}, R_{l}$ are separated by a distance greater than $2^h$. Since $c_{l}\leq \frac{2^h}{2}$,  $L_{l}, L_{l+1}$ (resp. $R_{l}, R_{l+1}$) are separated by a distance greater than $2^h$ for all $l\geq j_{h}$. Let $Z_{h}$ be any maximal $\frac{2^h}{2}$-separated subset which contains $Z'_{h}$. Let $K^c_{h}$ be the subcomplex of $|\mathcal{U}_{h}|$ generated by vertices in $Z_{h}$ with height $\geq {j_{h}}$.\\

It can be checked explicitly that $K^c_{h}$ is properly homotopic to $X^c_{j_{h}}$: vertices in rungs at different heights do not share an edge, and a point in the interior of the rung at height $l$ shares an edge with at at most one vertex in $L$ (resp. $R$), unless we are in the specific scenario of $l = (k+\frac{1}{2})2^h$. In this case, the $2$-simplex defined by the three vertices at heights $l, k2^h, (k+1)2^h$ exists. The only other thing to consider is a triple of points in the interior of a rung with pairwise non-empty intersection. Here, a $2$-simplex also exists, because the triple intersection is an interval in the rung.\\

Therefore, $i_{h}: X^c_{j_{h}} \rightarrow K^c_{h}$ induces an isomorphism on end homotopy groups; the colimit collapses to $h=0$ and we obtain
\begin{align*}
\pi^c_{0}(X,\alpha) \cong \pi_{0}^{e}(X,\alpha) \cong  (\prod_{j\in \mathbb{N}_{0}} F_{(\mathbb{N}_{0})_{\geq j}})_{/\image(\varphi)}
\end{align*}
Additionally, is not difficult to show that all higher homotopy groups of $X$ vanish. \\

We now generalise this computation to arbitrary ladders, addressing a question posed in \cite{ashley2025interactions} about how to define the "size" of a hole.

\begin{lemma} \label{laddersize} Let $X$ be a ladder: a $1$-dimensional simplicial complex defined by two rays $\alpha,\alpha'$ joining at a distinguished vertex at height $0$, and rungs of length $l_{m}\in \mathbb{N}$  joining the rays at heights $h_{m} \in \mathbb{N}$. Let $X$ be given the path metric. Assume that $\{h_{m}\}_{m\in \mathbb{N}}$ is a strictly increasing sequence and that $\lim_{m\rightarrow \infty} l_{m}=\infty$. Then the coarse homotopy groups of $X$ are:

\[ \pi_{n}^{c}(X,\alpha) \cong \begin{cases} 
       0 &n\geq 1  \\
        {\varprojlim}^1_{m}  F_{(\mathbb{N}_{0})_{\geq m}} & n=0 \\ 
       \end{cases}
    \]
where the bonding maps $\psi_{m+1,m}: F_{(\mathbb{N}_{0})_{\geq m+1}} \rightarrow F_{(\mathbb{N}_{0})_{\geq m}}$ are the canonical inclusion. 
\end{lemma}

\begin{proof} Fix $h\in \mathbb{N}$. Let $j_{h}$ be a power of $2$ such that $5\cdot 2^h < l_{m}$ for all $h_{m}\geq j_{h}$. Construct the sets $Z_{h}'$ as before and extend to a maximal $\frac{2^h}{2}$-separated subset $Z_{h}$ of $X$. Let $K^c_{h}$ be the subcomplex of $|\mathcal{U}_{h}|$ generated by vertices in $Z_{h}$ with height $\geq j_{h}$. $i_{h}: X^c_{j_h} \rightarrow K^c_{h}$ is a proper homotopy equivalence; the colimit collapses to $h=0$. We have
\begin{align*}
\pi_{n}^c(X,\alpha) \cong \pi_{n}^e(X,\alpha)
\end{align*}
for all $n\geq 0$. We can choose a finite filtration $\{X_{h_{m}}\}_{m\in \mathbb{N}}$ of $X$ corresponding to the rungs of the ladder. We have that 
\begin{align*}
\pi_{1}(X^c_{h_{m}}, \alpha(h_{m})) \cong F_{(\mathbb{N}_{0})_{\geq m}}
\end{align*}
The bonding maps induced by the inclusion of spaces can be identified with the inclusion of groups $F_{(\mathbb{N}_{0})_{\geq m+1}}\rightarrow F_{(\mathbb{N}_{0})_{\geq m}}$. The higher homotopy groups of $X^c_{h_{m}}$ vanish for all $m\in \mathbb{N}$. Since $X^c_{h_{m}}$ is connected, $X$ has one end. Applying the ${\varprojlim}^1$ sequence gives us the claim. 
\end{proof}

This lemma shows that the coarse homotopy groups cannot distinguish between different ladders that have rungs whose lengths grow unboundedly. We now show that there is a ladder $Y$ and a coarse map $\mathfrak{a}: Y\rightarrow X$ which induces an isomorphism on coarse homotopy groups, but is not a coarse homotopy equivalence.  

\begin{example} (Failure of general coarse Whitehead) Let $Y$ be a ladder with rungs of length $2^{m}$ joining the rays at heights $2^{m}$ for $m\geq 10$. Give $Y$ a simplicial structure with edges of length $1$.  Let $X$ be a ladder with rungs of length $m^2$ joining the rays at heights $m$ for $m\in \mathbb{N}$. Let $\mathfrak{a}: Y\rightarrow X$ be the distinguished-vertex-preserving piecewise-linear map which sends rungs at height $2^{m}$ to height $m$. \\

This map is $\frac{1}{4}$-Lipschitz so has a simplicial approximation: let $\tilde{\mathfrak{a}}(v) = x_v$ for the vertices $v \in Y$. To apply Theorems \ref{main1} and \ref{main2} we only need to consider the map $\xi_{0,0}: Y\rightarrow X$ which is the composition 
\begin{align*}
\widetilde{i_{0}\mathfrak{a}R_{0}} 
\end{align*}
Observe that $i_{0}\mathfrak{a}R_{0}$ satisfies the same star condition as $\id_{X}\mathfrak{a} \id_{Y}$: for a vertex $v\in Y$ we have $i_{0}\mathfrak{a}R_{0} \langle st(v) \rangle \subset i_{0}\mathfrak{a} \langle st(v) \rangle \subset i_{0} \langle st(x_{v}) \rangle \subset \langle st(x_{v}) \rangle$. Therefore it suffices to consider the induced map 
\begin{align*}
\mathfrak{a}_{\ast}: \pi_{n}^{e}(Y,\alpha)\rightarrow \pi_{n}^{e}(X,\mathfrak{a}\alpha)
\end{align*}
Since $\mathfrak{a}$ restricted to height $\geq 2^{10}$ is a proper homotopy equivalence, it induces an isomorphism on end homotopy groups for all $n\geq 0$. \\

$\mathfrak{a}$ does not have a coarse homotopy inverse. Suppose for a contradiction there existed a $\mathfrak{b}: X\rightarrow Y$ such that $\mathfrak{b}\circ \mathfrak{a}$ is coarsely homotopic to $\id_{Y}$ via the homotopy $H: I_{p}(Y)\rightarrow Y$ (as usual, we can assume that $p: Y\rightarrow [0,\infty)$ is the projection to the height variable). There exists a $h$ large enough such that $i_{h}\mathfrak{b}R_{0}$ has a simplicial approximation. Then we have the commutative diagram 

\begin{figure}[H]
\center
\begin{tikzcd}
&\pi_{0}^{c}(Y,\alpha) \arrow{r}{\mathfrak{a}_{\ast}}  \arrow{d}{\theta}
&\pi_{0}^{c}(X, \mathfrak{a}\alpha) \arrow{d}{\theta} \arrow{r}{\mathfrak{b}_{\ast}} 
& \pi_{0}^{c}(Y,\alpha) \arrow{d}{\theta}
\\
& \pi^{e}_{0}(Y, \alpha)  \arrow{r}{\mathfrak{a}_{\ast}}
& \pi^{e}_{0}(X,  \mathfrak{a}\alpha) \arrow{r}{\widetilde{i_{h}\mathfrak{b}R_{0}}_{\ast}} 
& \pi^{e}_{0}(|\mathcal{V}_{h}|, \widetilde{i_{h}\alpha})  
\end{tikzcd}
\end{figure}
 By construction we have that $\widetilde{i_{h}\mathfrak{b}R_{0}} = \widetilde{i_{h}\mathfrak{b}}$, which is a proper, simplicial map. $(\mathfrak{b})_{\ast} (\mathfrak{a})_{\ast}$ is the identity: this means that $(\widetilde{i_{h}\mathfrak{b}})_{\ast} \mathfrak{a}_{\ast}$ on end homotopy groups has to be the identity too, after the identification of $\pi^{e}_{0}(|\mathcal{V}_{h}|, \widetilde{i_{h}\alpha})$ with $ \pi_{0}^{e}(Y,\alpha)$. We will show that actually $\widetilde{i_{h}\mathfrak{b}}\mathfrak{a}$ is properly homotopic to $i_{h}: Y\rightarrow |\mathcal{V}_{h}|$. \\

To that end, consider the simplicial complex $I_{p}(Y):=\{(x,p(x)t) \in Y \times [0,\infty) \,|\, t\in [0,1]\}$ with vertex set $Z$ equal to the intersection of the integer lattice $\mathbb{Z}^3$ (in ladder coordinates) with $I_{p}(Y)$. We equip $I_{p}(Y)$ with the metric 
\begin{align*}
d_{\infty}((x,p(x)t), (y,p(y)l)) = \max\{d(x,y), |p(x)t-p(y)l|\}
\end{align*}
This is bi-Lipschitz equivalent to the path metric on $I_{p}(Y)$ induced from the ambient $L^1$-product metric on $Y \times [0,\infty)$. It is clear that $Z$ is a maximal $\frac{1}{2}$-separated set of $I_{p}(Y)$.  We can denote a point in $Z$ by $z_{mlt}:= ((m,l),t)$ for some $m\in \mathbb{N}_{0}, 0 \leq l \leq m$ and $0 \leq t \leq m$, where $(m,l)$ is written in ladder coordinates. $I_{p}(Y)$ is a subcomplex of the nerve $|\mathcal{W}|$ of the open cover $\mathcal{W}:=\cup_{z\in Z} B(z,1)$. Let $i_{\mathcal{W}}: I_{p}(Y)\rightarrow |\mathcal{W}|$ and $R_{\mathcal{W}}: |\mathcal{W}|\rightarrow I_{p}(Y)$ be defined as usual. By construction we have that $\widetilde{{R_{\mathcal{W}}}_{|I_{p}(Y)}} = \id_{I_{p}(Y)}$ \\

We can assume that we have chosen $h$ large enough so that $i_{h}HR_{\mathcal{W}}: |\mathcal{W}|\rightarrow |\mathcal{V}_{h}|$ has a simplicial approximation $(i_{h}HR_{\mathcal{W}})^{\sim}$. We therefore obtain a  proper homotopy 
\begin{align*}
(i_{h}HR_{\mathcal{W}})^{\sim}_{| I_{p}(Y)}: I_{p}(Y)\rightarrow |\mathcal{V}_{h}|
\end{align*}
between $(i_{h}\mathfrak{b} \mathfrak{a}(R_{\mathcal{W}})_{Y\times \{0\}})^{\sim} = (i_{h}\mathfrak{b} \mathfrak{a})^{\sim} = (i_{h}\mathfrak{b})^{\sim} (\mathfrak{a})^{\sim} \simeq (i_{h}\mathfrak{b})^{\sim} \mathfrak{a}$ and $(i_{h}R_{\mathcal{W}})^{\sim} = \widetilde{i_{h}}$. By rescaling, we obtain a proper homotopy $H: Y\times [0,1]\rightarrow |\mathcal{V}_{h}|$ between $(i_{h}\mathfrak{b})^{\sim} \mathfrak{a}$ and $ \widetilde{i_{h}}$. \\

It can be checked explicitly (because of the nicely chosen separation of rungs and rung lengths) that $ \widetilde{i_{h}}: Y\rightarrow |\mathcal{V}_{h}|$ restricted to the subcomplex $Y^c_{h}$ (generated by vertices with height $\geq 2^h$) is a homeomorphism, and the induced map on $\pi_1$ can be canonically identified with the identity. We identify $Y^c_{h}$ with its image $\widetilde{i_{h}}(Y^c_{h})\subset |\mathcal{V}_{h}|$. Choose $j$ large enough such that $H(Y^c_{j}\times [0,1])\subset Y^c_{h}$. Consider the loop $\gamma_{j}$ (traversing the $j$-th circle anti-clockwise) which is a generator of $\pi_1(Y^c_{h},\alpha(2^{j}))$. We have that 
\begin{align*}
(i_{h}\mathfrak{b})^{\sim}_{\ast} \mathfrak{a}_{\ast}[\gamma_{j}] = b\widetilde{i_{h}}_{\ast}[\gamma_{j}] = b[\gamma_{j}]
\end{align*}
where $b$ is the change of base point homomorphism between $\alpha(2^{j})$ and $(i_{h}\mathfrak{b})^{\sim} \mathfrak{a}\alpha(2^{j})$ via the homotopy $H$. We know that $\mathfrak{a}$ sends $\gamma_{j}$ to the $j$-th loop in $X$, which has length $j^2+(j+1)^2+2$. Since $(i_{h}\mathfrak{b})^{\sim}$ is simplicial, it is length-decreasing. But we know that $\gamma_{j}$ has length  $\frac{1}{2^h}(2^j+2^{j+1}+2\cdot 2^j)$ in $|\mathcal{V}_{h}|$, and any element in the homotopy class $ b[\gamma_{j}]$ must have at least this length: contradiction for $j$ sufficiently large.
\end{example}

The example shows us the following:

\begin{theorem} \label{fail} There exist $1$-dimensional simplicial complexes $Y,X$ and a proper Lipschitz map $\mathfrak{a}: Y\rightarrow X$ which induces an isomorphism on all coarse homotopy groups, but is not a coarse homotopy equivalence. 
\end{theorem}

\newpage

\section{Relaxing the restrictions}

We compare the coarse homotopy groups of $X$ to the end homotopy groups of $X$ for $X$ a proper metric space. 

\begin{theorem} \label{endcoarse} Let $X$ be a proper metric space with a continuous coarse base ray $\omega: [1,\infty)\rightarrow X$. Let $n\geq 1$. There exists a group homomorphism 
\begin{align*}
\lambda: \pi_{n}^{e}(X,\omega)&\rightarrow \pi_{n}^{c}(X,\omega)\\
[f] &\mapsto b_{\omega s,\omega} [f\circ s]
\end{align*}
where $s: c[0,1]^n \rightarrow c[0,1]^n$ is a shrinking map.  
\end{theorem}

\begin{proof} Let $h\in \mathbb{N}_{0}$. Consider the set $c_{[2^h,2^{h+1}]}[0,1]^n = [0,1]^n \times [2^h,h^{h+1}]$, the parts of the cone with heights in the interval $[2^{h},2^{h+1}]$. Let $[f]\in \pi_{n}^{e}(X,\omega)$. Since $f$ is continuous, $f_{|c_{[2^h,2^{h+1}]}[0,1]^n}$ is uniformly continuous. There exists a $\varepsilon_{h}>0$ such that 
\begin{align*}
d(x,y)< \varepsilon_{h} \implies d(f(x),f(y))<1
\end{align*}
We can assume that $\{\varepsilon_{h}\}_{h\in \mathbb{N}_{0}}$ is a decreasing sequence. Choose an increasing sequence of constants $\{m_{h}\}_{h\in \mathbb{N}_{0}}$ such that $2^{m_{h}}> \frac{1}{\varepsilon_{h}}$. The projection
\begin{align*}
q_{h}: [0,1]^n \times [2^{h+m_{h}}, 2^{h+1+m_{h}}] &\rightarrow [0,1]^n \times [2^{h}, 2^{h+1}]\\
 (tx,t)&\mapsto (\frac{t}{2^{m_{h}}}x, \frac{t}{2^{m_{h}}})
\end{align*}
has Lipschitz constant $\frac{1}{2^{m_{h}}}$. We have that $q_{h+1}([0,1]^n \times \{2^{h+1+m_{h+1}}\}) = q_{h}([0,1]^n \times \{2^{h+1+m_h}\})$ so we can extend $q:=\{q_{h}\}_{h\in \mathbb{N}_{0}}$ in the obvious way to a continuous map $s_{q}: c[0,1]^n \rightarrow c[0,1]^n$. The composition $f\circ s_{q}$ satisfies
\begin{align*}
d(x,y)< 1 \implies d(fs_{q}(x),fs_{q}(y))<1
\end{align*}
Since $c[0,1]^n$ is a path metric space, $f\circ s_{q}$ is controlled. It is proper because $f$ and $s_{q}$ are proper. Additionally, for any other choice $\{q'_{h}\}_{h\in \mathbb{N}_{0}}$ with $ps_{q'}\leq ps_{q}$, $fs_{q}$ is coarsely homotopic to $fs_{q'}$. Since $\omega$ is coarse, $\omega$ is coarsely homotopic to $\omega s_{q}$ and $\omega s_{q'}$. Furthermore, there exists a coarse homotopy $G: c[0,1]^2\rightarrow X$ such that 
\begin{align*}
G_{| [0,1] \times \{0\}} = b_{\omega s_{q},\omega}\\
G_{| [0,1] \times \{1\}} = b_{\omega s_{q'},\omega}\\
G_{| \{0\} \times [0,1]} = b_{\omega s_{q}, \omega s_{q'}}\\
G_{| [0,1] \times \{1\}} = \id_{\omega}
\end{align*}We obtain 
\begin{align*}
b_{\omega s_{q},\omega} [f\circ s_{q}] = b_{\omega s_{q'},\omega}[f\circ s_{q'}] 
\end{align*}
which shows that $\lambda$ is independent of the choice of $q$. \\

Suppose that $[f] = [g] \in \pi_{n}^{e}(X,\omega)$. Let $H:c ([0,1]^n \times [0,1]) \rightarrow X$ be a continous, proper homotopy between $f$ and $g$. As before, there exist a decreasing sequence of constants $\{\varepsilon_{h}\}_{h\in \mathbb{N}_{0}}$, an increasing sequence of constants $\{m_{h}\}_{h\in \mathbb{N}_{0}}$, and $\frac{1}{2^{m_{h}}}$-Lipschitz maps
\begin{align*}
q_{h}: [0,1]^{n+1} \times [2^{h+m_{h}}, 2^{h+1+m_{h}}] &\rightarrow [0,1]^{n+1} \times [2^{h}, 2^{h+1}]\\
 (tx,t)&\mapsto (\frac{t}{2^{m_{h}}}x, \frac{t}{2^{m_{h}}})
\end{align*} 
such that $H\circ s_{q}$ satisfies $d(x,y)< 1 \implies d(H s_{q}(x),H s_{q}(y))<1$
\end{proof}
This is a coarse homotopy between $fs_{q}$ and $gs_{q}$. After changing base rays, we obtian 
\begin{align*}
b_{\omega s_{q},\omega}[fs_{q}] = b_{\omega s_{q},\omega}[g s_{q}]
\end{align*}
It is clear that $\lambda$ is a group homomorphism. 
\newpage
\part{Conclusion}

\section{Further research}
\subsection{Loose ends}

We discuss some loose ends that have been left by the proofs in this thesis: \\

\textbf{Connectedness of $X$}: \\

We only considered connected spaces $X$ in Theorem \ref{unbasedwhitehead}. An analogous, but more complicated statement may be true for compact metric spaces $X$ such that $M_{X}$ has finitely many ends. \\

Suppose that $M_{X}$ has $t\in \mathbb{N}$ ends. Let $[\tau_{1}],\dots,[\tau_{t}]\in \pi_{0}^{c}(cX)$ be such that $\{\theta[\tau_{1}],\dots, \theta[\tau_{t}]\}$ are proper homotopy classes of proper rays in $M_{X}$ which represent these ends. By post-composing with a shrinking map, we can assume that $\widetilde{i\tau}$ exists for all $\tau\in \{\tau_{1},\dots, \tau_{t}\}$. I expect that by a similar proof to Proposition \ref{psiexistence}, for each $\tau\in \{\tau_{1},\dots, \tau_{t}\}$ there exists a $l \in \mathbb{N}$ such that the diagram

\begin{figure}[H]
\center
\begin{tikzcd}
&\pi_{n}^{c}(cX,\tau) \arrow{r}{\mathfrak{a}}  \arrow{d}{\theta}
&\pi_{n}^{c}(cY,\mathfrak{a}\tau) \arrow{d}{\theta \circ s_{l}} 
  \\ 
& \pi_{n}^{L,e}(M_X,\widetilde{i\tau}) \arrow{d}{\mathfrak{d}}
& \pi_{n}^{L,e}(M_Y,\widetilde{is_{l}\mathfrak{a}\tau})  \arrow{d}{\mathfrak{d}} 
\\
& \pi_{n}^{e}(M_{X},\widetilde{i\tau}) \arrow{r}{b \circ \psi} 
&\pi_{n}^{e}(M_{Y},\widetilde{is_{l}\mathfrak{a}\tau}) 
\end{tikzcd}
\end{figure}
commutes. The map $\psi$ here is the geometric cellular approximation $(is_{l}\mathfrak{a}R)^{\sim}$, and $b$ is the change of base ray homomorphism 
\begin{align*}
b: \pi_{n}^{e}(M_Y,\widetilde{is_{l}\mathfrak{a}R}\widetilde{i\tau}) \rightarrow \pi_{n}^{e}(M_Y,\widetilde{is_{l}\mathfrak{a}\tau}) 
\end{align*}

To show that the diagram commutes would be technically more challenging, since $R$ is not the identity now on $\widetilde{i\tau}$. If this works, we have the following statement:\\

\begin{conjecture} Let $X,Y$ be compact metric spaces with finite shape dimension. Let $\mathfrak{a}: cX\rightarrow cY$ be a coarse map. Assume that $M_{X}$ has $t\in \mathbb{N}$ ends. Let $[\tau_{1}],\dots,[\tau_{t}]\in \pi_{0}^{c}(cX)$ be such that $\{\theta[\tau_{1}],\dots, \theta[\tau_{t}]\}$ are proper homotopy classes of proper rays in $M_{X}$ which represent these ends. Assume that for all $\tau\in \{\tau_{1},\dots, \tau_{t}\}$, there is an isomorphism on coarse homotopy groups $\mathfrak{a}:\pi_{n}^{c}(cX,\tau)\rightarrow \pi_{n}^{c}(cY,\mathfrak{a}\tau)$ for all $n\geq 0$. Then there exists a $l\in \mathbb{N}$ such that $\psi:=(is_{l}\mathfrak{a}R)^{\sim}$ induces an isomorphism 
\begin{align*}
b\circ \psi: \pi_{n}^{e}(M_X,\widetilde{i\tau}) \rightarrow \pi_{n}^{e}(M_Y,\widetilde{is_{l}\mathfrak{a}\tau}) 
\end{align*} 
for all $n\geq 0$ and all $\tau$. \\

If $\psi$ additionally induces an isomorphism on the ends of $M_{X},M_{Y}$, then $\mathfrak{a}$ is a coarse homotopy equivalence. 
\end{conjecture}

Extending the result of coarse Whitehead to inifinitely many ends is significantly more problematic. First, in the proof of Theorem \ref{deforetract}, if $M_{Y}$ has infinitely many ends, the set $K_{j}$ may not be compact. The larger issue in the proof is that in the definition of $j:= \max_{\tau} j_{\tau}$ we have taken the maximum over a finite set. I do not see a way to overcome these problems: the statement itself may be false. \\

\textbf{Arbitrary ladders}:

\begin{conjecture} Let $X$ be a ladder: a $1$-dimensional simplicial complex defined by two rays $\alpha,\alpha'$ joining at a distinguished vertex at height $0$, and rungs of length $l_{m}\in \mathbb{N}$  joining the rays at heights $h_{m} \in \mathbb{N}$, where $\{h_{m}\}_{m\in \mathbb{N}}$ is a strictly increasing sequence. Let $X$ be given the path metric. Assume that the sequence of perimeters of the $m$-th circle $p_{m}=2(h_{m+1}-h_{m})+l_{m}+l_{m+1}$ satisfies $\lim_{m\rightarrow \infty} p_{m} = \infty$. Then the coarse homotopy groups of $X$ are:

\[ \pi_{n}^{c}(X,\alpha) \cong \begin{cases} 
       0 &n\geq 1  \\
        {\varprojlim}^1_{m}  F_{\mathbb{N}_{\geq m}} & n=0 \\ 
       \end{cases}
    \]
where the bonding maps $\psi_{m+1,m}: F_{\mathbb{N}_{\geq m+1}} \rightarrow F_{\mathbb{N}_{\geq m}}$ are the canonical inclusion. If the perimeters $p_{m}$ are bounded, then then $X$ is coarsely homotopy equivalent to a ray.
\end{conjecture}

We proved this statement via an explicit construction for $\lim_{m\rightarrow \infty} l_{m}= \infty$. The proof of the general statement would probably involve the nerve theorem. The conjecture being true would mean that coarse homotopy groups of an arbitrary ladder depend only on whether the perimeters of circles grow to infinity. 

\subsection{Conjectures}

We discuss some conjectures and possible generalisations of  our results.\\

  \textbf{Coarse Hurewicz map}: \\

Let $(X,\mathcal{C},\mathcal{B})$ be a bornological coase space. For $n \in \mathbb{N}$ we consider the abelian group $C\mathcal{X}_{n}$ of functions  $X^{n+1}\rightarrow \mathbb{Z}$ which are locally finite and controlled. An element $c\in C\mathcal{X}_{n}$ can be expressed as an infinite linear combination of $(n+1)$-tuples of points in $X$. $(C\mathcal{X},d)$ is a chain complex with differential
\begin{align*}
d_{n}: C\mathcal{X}_{n}(X)&\rightarrow C\mathcal{X}_{n-1}(X)\\
d_{n}c(x_0,\dots,x_{n-1})&= \sum_{i=0}^{n}(-1)^i \sum_{x\in X} c(x_0,\dots, x_{i-1}, x, x_{i}, \dots, x_n)
\end{align*} 
The coarse homology of $X$, denoted by $H\mathcal{X}(X)$, is the homology of $(C\mathcal{X},d)$. Coarse homology is a functor $H\mathcal{X}: \mathbf{BornCoarse} \rightarrow \mathbf{Ab}^{\mathbb{Z}_{gr}}$ which is excisive, invariant under coarse homotopies, $u$-continuous, strongly additive, and vanishes on flasques. See Chapter $5$ of $\cite{bunke2020lecture}$. \\

I believe there exists a coarse Hurewicz map $\mathcal{H}:\pi^c_{n}(X)\rightarrow H\mathcal{X}_{n+1}(X)$ which is an isomorphism for $n=0$ and is the abelianisation map for $n=1$. The proof I have in mind relies on the following statement:

\begin{conjecture} Let $n\in \mathbb{N}_{0}$. There is an isomorphism between $\pi^c_{n}(X,\omega)$ and the group (or pointed set) of base ray preserving coarse homotopy classes of base ray preserving coarse maps $f: (c'S^{n},c\{x_0\})\rightarrow (X,\omega)$, where 
\begin{align*}
c'S^n=\{(hx,h)\in \mathbb{R}^{n+1} \times [0,\infty)\,|\, x \in S^n, h\in [0,\infty)\}
\end{align*} 
\end{conjecture}

Assuming this is true, I can prove the following:

\begin{conjecture} (The Hurewicz map) \label{hurewicz} Let $(X,d)$ be a proper metric space and $n\in \mathbb{N}_{0}$. There is a natural morphism 
\begin{align*}
\mathcal{H}: \pi_{n}^{c}(X)\rightarrow H\mathcal{X}_{n+1}(X)
\end{align*}
such that the diagram

\begin{figure}[H]
\centering
\begin{tikzcd}
&\pi_{n}^{c}(X) \arrow{r}{\mathcal{H}} \arrow{d}{\cong}
&H\mathcal{X}_{n+1}(X)  \arrow{d}{\mathcal{A}}
  \\
&\varinjlim_{h\in \mathbb{N}} \pi^{e}_{n}(|\mathcal{U}_{h}|) \arrow{r}{\varphi}
&\varinjlim_{h\in \mathbb{N}} H^{lf}_{n+1}(|\mathcal{U}_{h}|) 
\end{tikzcd}
\end{figure}
commutes. 
\end{conjecture}

Let us explain all the maps in the diagram.\\

$\mathcal{H}:$ We can choose a simplicial structure on $S^n$ so that $c'S^n$ is a union of simplices of dimension $n+1$. Let $\Gamma$ be the set of these simplices. For $[f]\in \pi^{c}_{n}(X)$ we let 
\begin{align*}
\mathcal{H}[f] = \sum_{\sigma \in \Gamma} (f(v^{\sigma}_0),\dots, f(v^{\sigma}_{n+1})) \in C\mathcal{X}_{n+1}(X)
\end{align*}
which lies in the kernel of $d_{n+1}$ so defines a coarse homology class. \\

$\varphi:$ The orientation class $\eta: = \sum_{\sigma\in \Gamma} \sigma$ is a generator of $H^{lf}_{n+1}(c'S^n)$. As in ordinary topology, we let
\begin{align*}
\varphi_{h}: \pi^e_{n}(|\mathcal{U}_{h}|)&\rightarrow H^{lf}_{n+1}(|\mathcal{U}_{h}|)\\
\varphi_{h} [g] &= g_{*}[\eta]
\end{align*}
$\{\varphi_{h}\}_{h\in \mathbb{N}}$ commutes with the bonding maps so defines a group homomorphism $\varphi:\varinjlim_{h} \pi^{e}_{n}(|\mathcal{U}_{h}|) \rightarrow \varinjlim_{h} H^{lf}_{n+1}(|\mathcal{U}_{h}|)$.\\

$\mathcal{A}:$ Let $(C\mathcal{X}_{n}(X))_{h}$ be the subgroup of $\frac{2^h}{2}$-controlled chains. We denote by $(H\mathcal{X}(X))_{h}$ the homology of the chain complex $((C\mathcal{X}(X))_{h},d)$. It is a fact that
\begin{align*}
H\mathcal{X}_{n}(X) \cong \varinjlim_{h} (H\mathcal{X}_{n}(X))_{h}
\end{align*}
 For $c = \sum_{\tau \in X^{n+1}} c_{\tau}\tau \in (C\mathcal{X}_{n}(X))_{h}$, we consider each $\tau \in X^{n+1}$ to be a vertex map from the standard simplex $\Delta^n$ to $X$.  Let $r_{h}: |\mathcal{U}_{h}|\rightarrow |\mathcal{U}_{h}|$ to be the map which sends a point to its closest vertex. We define
\begin{align*}
A_{h}(c) = \sum_{\alpha} \sum_{\tau} c_{\tau} [sp(r_{h}i_{h}\tau): e^n_{\alpha}] e^n_{\alpha} \in C^{\infty}_{n}(|\mathcal{U}_{h}|)
\end{align*}
where $C^{\infty}_{n}(|\mathcal{U}_{h}|)$ is the cellular chain complex of locally finite  chains and $sp(r_{h}i_{h}\tau)$ is the simplex spanned by the vertices of $r_{h}i_{h}\tau(\Delta^n)$. The number $ [r_{h}i_{h} \tau: e^n_{\alpha}]$ is $\pm 1$ if and only if this simplex is $e^n_{\alpha}$ (up to orientation). The diagram 
\begin{figure}[H]
\centering
\begin{tikzcd}
&(H\mathcal{X}_{n}(X))_{h}  \arrow{r}{A_{h}} \arrow{d}{\iota}
&H^{lf}_{n}(|\mathcal{U}_{h}|) \arrow{d}{\phi_{h,h+1}}
  \\
&(H\mathcal{X}_{n}(X))_{h+1}  \arrow{r}{A_{h+1}}
&H^{lf}_{n}(|\mathcal{U}_{h+1}|) 
\end{tikzcd}
\end{figure}
 commutes for all $n\in \mathbb{N}_{0}$, and so defines a group homomorphism $\mathcal{A}: H\mathcal{X}_{n}(X) \cong \varinjlim_{h} (H\mathcal{X}_{n}(X))_{h}\rightarrow \varinjlim_{h} H^{lf}_{n}(|\mathcal{U}_{h}|)$. \\

The diagram in the conjecture commutes on the level of chains. \\

There is a conjectural $\varprojlim^1$ sequence for coarse homology: 

\begin{conjecture}  \label{shapeex} Let $Y$ be a compact Hausdroff space. Let $(Y_{i})_{i\in \mathbb{N}}$ be a decreasing family of closed subspaces of $Y$. We consider the set-theoretic intersection $\cap_{i\in \mathbb{N}} Y_{i}$ in $Y$ with the subspace topology. We have an exact sequence
\begin{align*}
0\rightarrow {\lim_{i\in \mathbb{N}}}^{1} H\mathcal{X}_{n+1}\mathcal{O}(Y_{i}) \rightarrow H\mathcal{X}_{n}\mathcal{O} (\cap_{i\in \mathbb{N}} Y_{i}) \rightarrow  \lim_{i\in \mathbb{N}} H\mathcal{X}_{n}\mathcal{O}(Y_{i}) \rightarrow 0
\end{align*}
for $n\in \mathbb{N}_{0}$.
\end{conjecture}

This conjecture is stated as Lemma $7.12$ in \cite{bunke2020lecture}, but the proof is incorrect. I believe I can prove it for $X = \cap_{i\in \mathbb{N}} Y_{i}$ a compact metric space, written as an intersection of ANR neighbourhoods of $X$ in $C_{b}(X)$, using Conjecture \ref{hurewicz}. This would require showing that $\mathcal{A}$ is an isomorphism (which should be true, since $\mathcal{A}$ should coincide with the assembly map in coarse homology), and using some shape-theoretic machinery to prove that 
\begin{align*}
\varinjlim_{h} H^{lf}_{n}(VR_{h}(cX))\cong  H^{lf}_{n}(VR_{1}(cX)) \cong H^{lf}_{n}(M_{X})
\end{align*}
where $VR_{h}(cX)$ is the nerve of the open cover of balls of radius $2^{h}$ centred at a maximal $\frac{2^h}{2}$-separated set of $cX$. \\

 \textbf{Coarse assembly maps}:\\

Cones are an object of interest in studying coarse assembly. 
 
 \begin{prop}  \label{assemblymap} (Proposition $11.23$ in \cite{bunke2020coarse}) Let $E:\mathbf{BornCoarse}\rightarrow \mathbf{C}$ be a coarse homology theory, and let $X$ be a uniform bornological coarse space. Assume:
 \begin{enumerate}
 \item $\mathbf{C}$ is complete. 
 \item $E$ is additive. 
 \item $X$ is homotopy equivalent in $\mathbf{UBC}$ to a locally finite, finite-dimensional simplicial complex equipped with the metric structures.
 \end{enumerate}
 Then we have an equivalence
 \begin{align*}
 (\Sigma E(\ast) \wedge \Sigma_{+}^{\infty})^{lf}(X)\simeq E\mathcal{O}^{\infty}(X)
 \end{align*}
 \end{prop}
 
 The restriction to $X$ homotopy equivalent to a locally finite, finite-dimensional simplicial complex is extremely strong, and is not satisfied for many spaces of interest (eg. arbitrary compact metric spaces or compact Hausdorff spaces). Note that $E$ is applied to the cone of a space, but the locally finite homology theory is applied to the space $X$ itself. This suggests that we are in some sort of degenerate case, and the $X$ on the left hand side should be replaced with something one dimension higher. It is also suspicious that the functor $X\mapsto (\Sigma E(\ast) \wedge \Sigma_{+}^{\infty})^{lf}(X)$ is defined on $\mathbf{TopBorn}$ instead of $\mathbf{UBC}$, and is homotopy invariant with respect to proper homotopies with are not necessarily uniform. This strongly suggests that it should have something to do with (a uniform version of) shape theory.

\begin{definition} An inverse ANR-system $\mathcal{X} = \{X_{\alpha}, \phi_{\alpha,'\alpha}; \mathcal{A}\}$ is said to be uniformly associated with a uniform space $X$ if there exist homotopy classes of uniformly continuous maps $[p_{\alpha}]: X\rightarrow X_{\alpha}$, $\alpha\in \mathcal{A}$, satisfying the following conditions:
\begin{enumerate}
\item $[p_{\alpha}] = [\phi_{\alpha',\alpha}] [p_{\alpha'}]$ for all $\alpha \leq \alpha'$.
\item For any ANR-space $P$ and any homotopy class of uniformly continuous maps $[f]: X\rightarrow P$, there exists an index $\alpha\in \mathcal{A}$ and a homotopy class $[h_{\alpha}]: X_{\alpha}\rightarrow P$ such that $[h_{\alpha}] [p_{\alpha}] = [f]$. 
\item If $[h] [p_{\alpha}] = [h'] [p_{\alpha}]$ for two homotopy classes of uniformly continuous maps $[h],[h']: X_{\alpha}\rightarrow P$, then there exists an index $\alpha'\geq \alpha$ such that $[h][\phi_{\alpha',\alpha}] = [h'] [\phi_{\alpha',\alpha}]$. 
\end{enumerate} 
We call such a system a \textit{uniform shape expansion} of $X$. 
\end{definition}

This definition just replaces all continuous maps and continuous homotopies with uniformly continuous maps and homotopies. In the compact metric case, the uniform condition is vacuous, and for a general proper metric space, the sequence of open covers 
\begin{align*}
\mathcal{U}_{h} := \bigcup_{z\in Z_{h}} B(z,\frac{1}{2^h})
\end{align*} for a maximal $\frac{1}{2^{h+1}}$ separated subset of $X$, is a uniform shape expansion of $X$. \\
 
I have the following idea. Let $\{X_{\alpha},\phi_{\alpha'\alpha}; \mathcal{A}\}$ be a uniform shape expansion $X$, and let $M$ be a model in topological spaces for 
\begin{align*}
M:=\text{holim}_{\alpha\in \mathcal{A}} \{X_{\alpha}, \phi_{\alpha'\alpha}\}
\end{align*}
If $X$ is a compact metric space, $M$ corresponds to the inverse mapping telescope $M_{X}$. There are canonical inclusions 
\begin{align*}
i_{\alpha}: \text{holim}_{\beta \leq \alpha} \{X_{\beta}\} \rightarrow M
\end{align*}
where $M_{\alpha} := \text{holim}_{\beta \leq \alpha} \{X_{\beta}\}$ denotes everything "below" level $\alpha$. The cofibre of the inclusion, denoted $C(i_{\alpha})$, is modelled by the pair $(M,M_{\alpha})$. I believe the following may be true:
 
 \begin{conjecture} Let $E:\mathbf{BornCoarse}\rightarrow \mathbf{C}$ be a coarse homology theory, and let $X$ be a uniform bornological coarse space. Assume:
 \begin{enumerate}
 \item $\mathbf{C}$ is complete. 
 \item $E$ is additive. 
\end{enumerate}
Then we have an equivalence 
 \begin{align*}
 \colim_{\alpha} (E(\ast) \wedge \Sigma_{+}^{\infty})^{lf}(C(i_{\alpha}))\simeq E\mathcal{O}^{\infty}(X)
 \end{align*}
 \end{conjecture}
 
This reduces to the result of Proposition \ref{assemblymap} in the case that $E = H\mathcal{X}$. If $X$ is a locally finite, finite-dimensional simplicial complex, it is in particular a proper metric space. The locally finite homology theory constructed from $H\mathcal{X}$ is ordinary locally finite homology $ (E(\ast) \wedge \Sigma_{+}^{\infty})^{lf} = H^{lf}$. For a simplicial complex, there is a $h'$ large enough such that $|\mathcal{U}_{h}|\simeq X$ for all $h\geq h'$. We have that
 \begin{align*}
 \colim_{h}(E(\ast) \wedge \Sigma_{+}^{\infty})^{lf}_{n}(C(i_{h})) = H_{n}^{lf}(C(i_{h'})) = H_{n}^{lf}(X\times \mathbb{R}) = \Sigma H_{n}^{lf} (X) \\
 = (\Sigma E(\ast) \wedge \Sigma_{+}^{\infty})^{lf}_{n}(X)
 \end{align*}
 
Assume the conjecture is true. If $E = H\mathcal{X}$ and $X$ is a compact metric space, we have
 \begin{align*}
 \colim_{h}(E(\ast) \wedge \Sigma_{+}^{\infty})^{lf}_{n}(C(i_{h})) = H_{n}^{lf}(M_{X}) 
\end{align*}
Therefore, this case reduces to proving Conjecture \ref{shapeex}. \\

 \textbf{Coarse homotopy groups of compact Hausdorff spaces}:\\

Let $(X, \mathcal{T})$ be a compact Hausdorff space. We begin with the theorem: 

\begin{theorem} (Chapter $5.2$, Theorem $7$ in \cite{mardesic1982shape}) Every compact Hausdorff space $X$ is the limit of an inverse system $\mathcal{X} = \{X_{\lambda}, p_{\lambda',\lambda}; \Gamma\}$ of compact polyhedra $X_{\lambda}$ with $PL$-bonding maps. One can always achieve that the cardinal $card(\Gamma) \leq w(X)$, the weight of $X$.
\end{theorem}

One proves this theorem by embedding $X$ as a subset of the Hilbert cube $\prod_{\alpha\in A} I_{\alpha}$, where the cardinality of $A$ is $w(X)$. This is the limit of a system of finite-dimensional cubes with projections as bonding maps, denoted 
\begin{align*}
\prod_{\alpha\in A} I_{\alpha} = \lim_{\alpha} K_{\alpha}
\end{align*}
let $p_{\alpha}: \prod_{\alpha\in A} I_{\alpha} \rightarrow K_{\alpha}$ denote the projection onto $K_{\alpha}$ and $p_{\alpha' \alpha}: K_{\alpha'}\rightarrow K_{\alpha}$ denote the bonding maps. We consider the sets $p_{\alpha}(X) \subset K_{\alpha}$, and take polyhedral neighbourhoods $K_{\alpha i}$ of $p_{\alpha}(X)$. Let $\Gamma = A \times \mathbb{N}$. We let $(\alpha, n)\leq (\alpha', n')$ provided that $\alpha \leq \alpha '$, $p_{\alpha' \alpha}(K_{\alpha' n'}) \subset K_{\alpha n}$ and we define
\begin{align*}
q_{\alpha' n', \alpha n} = {p_{\alpha' \alpha}}_{|K_{\alpha' n'}}: K_{\alpha' n'}\rightarrow K_{\alpha n}\\
q_{\alpha n} = p_{\alpha}|_{X}: X \rightarrow K_{\alpha n}
\end{align*}
It can be checked that $\mathcal{X} =\{K_{\alpha n}, q_{\alpha' n', \alpha n}; \Gamma\}$ is a shape expansion of $X$. \\

It can be shown that any shape expansion $\mathcal{X}$ is associated to a $\mathcal{T}$-admissible function $\phi_{\mathcal{X}}: \mathbb{N}\rightarrow \mathcal{P}(X\times X)$. I expect that $\pi_{n}^{c}(\mathcal{O}_{\phi}(X))$ is related to $\lim_{\alpha} \lim_{n} \pi_{n}(K_{\alpha n})$, possibly via a spectral sequence. Taking the colimit over all possible shape expansions, there should be a relationship between $\pi_{n}^{c}(\mathcal{O}(X))$ and 
\begin{align*}
\colim_{\mathcal{X}} \lim_{\alpha} \lim_{n} \pi_{n}(K_{\alpha n})
\end{align*}

\newpage 
\part{Appendix}
\section{Technical details}\subsection{The metric on $M$}

This subsection proves some elementary statements about the path metric on $M$, which are used in many places in main section. \\

Let $M = \tilde{M}/_{\sim}$ denote the mapping telescope as a quotient space, and let $\pi: \tilde{M}\rightarrow M$ denote the quotient map. Recall 
\begin{align*}
\tilde{M} = \coprod_{h\in \mathbb{N}_{0}} |\mathcal{U}_{h}| \times [2^h,2^{h+1}] = \coprod_{h\in \mathbb{N}_{0}} \tilde{M}^{h}
\end{align*}
and that $\sim$ is the equivalence relation obtained from the simplicial gluing maps $\phi_{h}: |\mathcal{U}_{h}| \times \{2^h\}\rightarrow |\mathcal{U}_{h-1}| \times \{2^{h}\}$. For technical purposes, we now define some projections. Let $h\in \mathbb{N}$. Define the projection to the $2^{h}$-th slice
\begin{align*}
\mathfrak{q}_{2^{h}}: \big\{ \bigcup_{l\geq h} |\mathcal{U}_{l}| \times [2^l,2^{l+1}] \big\} \cup_{\phi_{h}} \big\{ |\mathcal{U}_{h-1}| \times (2^{h-1}, 2^h]\big\} \longrightarrow |\mathcal{U}_{h-1}| \times \{2^h\}
\end{align*}
as $\mathfrak{q}_{2^h}(x,t)=(\phi_{l,h-1}(x),2^h)$ for $t\in [2^l,2^{l+1}], l\geq h$ and as $\mathfrak{q}_{2^h}(x,t) = (x,2^h)$ for $t\in (2^{h-1}, 2^h]$. This is compatible with the gluing maps so passes to the quotient. $\mathfrak{q}_1$ is defined similarly in the obvious way. \\

Let $s\in (2^h,2^{h+1})$. Define the projection to the $s$-th slice 
\begin{align*}
\mathfrak{q}_{s}: \big\{ \bigcup_{l> h} |\mathcal{U}_{l}| \times [2^l,2^{l+1}] \big\} \cup_{\phi_{h+1}} \big\{  |\mathcal{U}_{h}| \times (2^h, 2^{h+1}] \big\}\longrightarrow |\mathcal{U}_{h}| \times \{s\}
\end{align*}
as $\mathfrak{q}_{s}(x,t) = (\phi_{l,h}(x),s)$ for $t\in [2^l,2^{l+1}]$, $l> h$ and as $\mathfrak{q}_{s}(x,t) = (x,s)$ for $t\in (2^h,2^{h+1}]$.

\begin{definition} \label{metricdef} The metric on $M$ is defined as follows: for two points $x,y\in M$ let
\begin{align*}
d(x,y):= \inf_{\gamma_{xy}} \inf_{\tilde{\gamma}_{xy}}  l(\tilde{\gamma}_{xy})
\end{align*}
where $\gamma_{xy}$ is a path in $M$ from $x$ to $y$, $\tilde{\gamma}_{xy}$ is a lift of $\gamma_{xy}$ to $\tilde{M}$ ie. $\pi\circ \tilde{\gamma}_{xy} = \gamma_{xy}$, and $l(\tilde{\gamma}_{xy})$ is calculated as the sum of lengths of path components in each disjoint component $\tilde{M}^{h}$ of $\tilde{M}$ with product metric $d_{\tilde{M}^{h}} = d_{|\mathcal{U}_{h}|} + d_{[2^h,2^{h+1}]}$, where $d_{|\mathcal{U}_{h}|}$ is the path metric induced from the spherical metric on each simplex. 
\end{definition}

Let us describe what some geodesics between two points $x=(u,t)$ and $y = (v,s)$ in $M$ look like. Let $\gamma$ be a path between $x$ and $y$. By "shorter", we mean "length less than or equal to". 

\begin{enumerate}
\item Any path segment $\eta$ of $\gamma$ with endpoints at height $l$ and image at least height $l$ can be replaced with a shorter segment $\eta'$ with the same endpoints. \\

$\eta' = \mathfrak{q}_{l} \circ \eta$ clearly satisfies the desired properties: $\eta'$ is shorter than $\eta$ because the maps $\mathfrak{q}$ are $1$-Lipschitz.

\item Any path segment $\eta$ of $\gamma$ beginning at $x=(u,t)$, ending at height $l<t$ and with image between heights $l,t$ can be replaced with a shorter segment $\eta'= \sigma_x \ast \tau$ with the same endpoints, where $\sigma_x$ travels straight downwards, and $\tau$ has image contained in the slice at height $l$. \\

Let $\tau = \mathfrak{q}_{l} \circ \eta$ and let $\sigma_x$ be the straight downward path from $x$ to $\mathfrak{q}_{l}(x)$. It is clear that $\eta' = \sigma_{x}\ast \tau$ has the same endpoints as $\eta$ and it is shorter, since the metric on $\tilde{M}$ is the product metric, and the $\phi$ used to define $\mathfrak{q}_{l}$ are $1$-Lipschitz.

\item The shortest path $\gamma$ is a concatenation of $\sigma_{x}\ast \eta \ast \sigma_{y}^{-1}$ where $\sigma_{x},\sigma_{y}$ are straight downward paths beginning at $x$ and $y$ respectively, and $\eta$ is a geodesic path with image contained in some height slice $l\leq \min\{t,s\}$. \\

By point $1$ we can assume that the height variable of $\gamma$ never has negative second derivative. Let $l$ be the minimal height reached by $\gamma$. If $l\in \{t,s\}$ we can replace $\gamma$ by $\sigma_{x} \ast \tau$ or $\tau \ast \sigma_{y}^{-1}$ as in point $2$, then replace $\tau$ by a geodesic path $\eta$ with the same endpoints. If $l<\min \{t,s\}$ then $\gamma$ can be written as a concatenation $\gamma_{xz} \ast \gamma_{yz}^{-1}$ where $z$ is a point at height $l$, $\gamma_{xz}$ is a path between $x$ and $z$ and $\gamma_{yz}$ is a path between $y$ and $z$.  By point $2$ we can replace $\gamma$ with $\sigma_{x}\ast \tau_{\mathfrak{q}_{l}(x)z} \ast \tau_{\mathfrak{q}_{l}(y)z}^{-1}\ast \sigma_{y}^{-1}$. Finally we can replace $\tau_{\mathfrak{q}_{l}(x)z} \ast \tau_{\mathfrak{q}_{l}(y)z}^{-1}$ with a geodesic path $\eta$ between $\mathfrak{q}_{l}(x)$ and $\mathfrak{q}_{l}(y)$. 

\item Let $r=\min\{t,s\}$. The distance between $x$ and $y$ is
\begin{align*}
\min \{\min_{2^h<r} l(\tilde{\gamma}^{2^h}), l(\tilde{\gamma}^{r})\} \\
= \min \{\min_{2^h<r} |t-2^h|+|s-2^h| + d_{\tilde{M}^{h-1}}(\widetilde{\mathfrak{q}_{2^h}(x)}, \widetilde{\mathfrak{q}_{2^h}(y)}), |t-s| + d_{r}(\mathfrak{q}_{r}(x), \mathfrak{q}_{r}(y))\}
\end{align*}
where $\gamma^{2^h}=\sigma_{x}\ast \eta_{2^h} \ast (\sigma_{y})^{-1}$ is concatenation as described above with the image of $\eta_{2^h}$ contained in the slice at height $2^h$, and similarly for $\gamma^r$. $d_{r}$ refers to the restriction of the metric of $M$ to height $r$, and the tilde denotes lifts to $\tilde{M}$ with minimal distance.\\

By point $3$, any path can be replaced with one of the form $\sigma_{x}\ast \eta_{l} \ast \sigma_{y}^{-1}$ where $\eta$ lives on height $l$ so it suffices to consider only paths of these form. If $l=r$ then we are in the case of $\gamma^r$. Otherwise, $l<r$. If $l\neq 2^h$ for some $2^h<r$, $\gamma$ can be replaced with a strictly shorter path which is either of the form $\sigma_{x}\ast \eta_{r} \ast \sigma_{y}^{-1}$ or by $\sigma_{x}\ast \eta_{2^h} \ast \sigma_{y}^{-1}$ where $2^h$ is the smallest power of $2$ greater than $l$, by projecting the $\eta$ segment upwards to the appropriate height and shortening $\sigma_{x},\sigma_{y}$.\\

The length of $\gamma^r$ is a sum of the distance $|t-s|$ to the lowest point $r=\min\{t,s\}$ and the minimal distance $d_{r}(\mathfrak{q}_{r}(x),\mathfrak{q}_{r}(y))$ between the endpoints of the path segment with image in slice $r$. Observe that $d_{r}= d_{|\mathcal{U}_{h}| \times \{r\}}$ for $r\in (2^h,2^{h+1})$ because there are unique lifts. \\

The length of $\gamma^h$ is a sum of the lengths of the segments $\sigma_x, \sigma_y$ and the minimal distance $d_{2^h}(\mathfrak{q}_{2^h}(x),\mathfrak{q}_{2^h}(y))$ between the endpoints of $\eta_{2^h}$. Observe that $d_{2^h} = d_{|\mathcal{U}_{h-1}| \times \{2^h\}}$ for $r=2^h$, since lifting to $\tilde{M}^{h-1}$ minimises length due to $\phi_{h}$ being $1$-Lipschitz.
\end{enumerate}

\begin{lemma} $(M,d)$ is a well-defined path metric space with distances bound above by the product metric on each simplex $\sigma \times I \subset \tilde{M}^{h}$.
\end{lemma}

\begin{proof} We need to prove symmetry, triangle inequality, and non-degeneracy. To show that $d(x,y)=d(y,x)$, note that for any  $\gamma_{xy}$ there is an reverse path $\gamma^{-1}_{xy}$ from $y$ to $x$ which satisfies $\inf_{\tilde{\gamma}_{xy}} l(\tilde{\gamma}_{xy}) = \inf_{\tilde{\gamma}^{-1}_{xy}} l(\tilde{\gamma}^{-1}_{xy})$. Therefore
\begin{align*}
d(x,y) =  \inf_{\gamma_{xy}} \inf_{\tilde{\gamma}_{xy}}  l(\tilde{\gamma}_{xy}) =  \inf_{\gamma^{-1}_{xy}} \inf_{\tilde{\gamma}^{-1}_{xy}} l(\tilde{\gamma}^{-1}_{xy}) \leq d(y,x)
\end{align*}
and a similar inequality in the other direction gives us $d(x,y) = d(y,x)$.\\

Triangle inequality: Let $x,y,z\in M$. We show that $d(x,z)\leq d(x,y) + d(y,z)$. Note that for any two paths $\eta_{xy}, \theta_{yz}$ we can concatenate them to obtain a path $\eta_{xy}\ast \theta_{yz}$ between $x$ and $z$. Hence
\begin{align*}
d(x,z) = \inf_{\gamma_{xz}} \inf_{\tilde{\gamma}_{xz}} l(\tilde{\gamma}_{xz}) \leq \inf_{\eta_{xy}\ast \theta_{yz}} \inf_{\widetilde{\eta_{xy}\ast \theta_{yz}}} l(\widetilde{\eta_{xy}\ast \theta_{yz}}) = \inf_{\eta_{xy}} \inf_{\tilde{\eta}_{xy}} l(\tilde{\eta}_{xy}) + \inf_{\theta_{xy}} \inf_{\tilde{\theta}_{xy}} l(\tilde{\theta}_{xy}) \\
= d(x,y) + d(y,z)
\end{align*}

Non-degeneracy: $d(x,x)=0$ since we can take the constant path. Let $x = (u,t)$ and $y=(v,s)$ be different points. If $t\neq s$ then any path has length $\geq |t-s|>0$ so $d(x,y)>0$. Suppose now that $t=s$. From the distance formula we have
\begin{align*}
d(x,y) = \min \{\min_{2^h<t} 2|t-2^h| + d_{\tilde{M}^{h-1}}(\widetilde{\mathfrak{q}_{2^h}(x)}, \widetilde{\mathfrak{q}_{2^h}(y)}), d_{t}(x, y)\}
\end{align*}
Since $d_t$ and $d_{\tilde{M}^{h-1}}$ are metrics, we have that $d(x,y)>0$. This shows non-degeneracy. It is clear by the derivation of the distance formula that $d$ is a path metric and that distances are bound above by the product metric on each $\sigma \times I \subset \tilde{M}^{h}$. 
\end{proof}

\begin{lemma} \label{Lebmetric} The Lebesgue number of the covering of $M$ by open geometric stars is $\frac{1}{4}$. 
\end{lemma}

\begin{proof} Let $U$ be a set of diameter $<\frac{1}{4}$. Let $M_{(2^h-\frac{1}{4}, 2^h+\frac{1}{4})}$ denote the set
\begin{align*}
M_{(2^h-\frac{1}{4}, 2^h+\frac{1}{4})} := \{(u,t)\in M \,|\, 2^h-\frac{1}{4} < t < 2^h + \frac{1}{4} \}
\end{align*}
Case 1: $U\subset M\setminus (\cup_{h\in \mathbb{N}} M_{(2^h-\frac{1}{4}, 2^h+\frac{1}{4})})$. In this case, $U$ has a unique lift to some $\tilde{M}^{h}$. For any $x,y\in U$, the shortest path $\tilde{\gamma}_{xy}$ lies completely in $\tilde{M}^{h}$ (since the vertical components $\sigma_x,\sigma_y$ for any other possibility have length $\geq\frac{1}{4}$). Therefore the metric on $M$ restricted to $U$ is identical to the metric on $\tilde{M}^{h}$. Because $\tilde{M}^{h}$ is a simplicial complex with a nice metric, $U\subset \langle gst(z) \rangle$ for some vertex $z$. To see this, pick a point $x=(u,t) \in U$ and the closest vertex $z=(v,s)$ to $x$. We know that $|t-s|\leq\frac{1}{2}$ and that $u\in \langle st(v)\rangle \subset |\mathcal{U}_{h}|$ with $d(u,v)\leq \frac{1}{2}$. It is clear that $\mathfrak{q}_{s}(U)$ has diameter $<\frac{1}{4}$. We claim that $\mathfrak{q}_{s}(U)  \subset \langle st(v)\rangle$. Since $d_{|\mathcal{U}_{h}|}(\partial \langle st(v)\rangle, v) = 1$ in the spherical metric any point outside of $\langle st(v) \rangle$ has distance at least $1$ from $v$. Furthermore, the heights of points in $U$ are at most $\frac{1}{2}+\frac{1}{4}=\frac{3}{4}$ from $s$. Therefore, $U\subset \langle gst(z)\rangle$. (If $s$ happens to be $2^h$, we identify $z$ with its image in $M$ and the argument still holds.) \\

Case 2: $U\cap M_{(2^h-\frac{1}{4}, 2^h+\frac{1}{4})}\neq \emptyset$ for some $h\in \mathbb{N}$. Consider the image $\mathfrak{q}_{2^h}(U)\subset |\mathcal{U}_{h-1}| \times \{2^h\}$. It is clear that $diam(\mathfrak{q}_{2^h}(U))<\frac{1}{4}$ so by previous arguments there exists $(z,2^h)$ such that $\mathfrak{q}_{2^h}(U)\subset \langle st(z) \rangle \subset |\mathcal{U}_{h-1}|$. We now argue that $U\subset \langle gst(z,2^h)\rangle$. For points $(u,t)$ with $t\leq 2^h$, we just have $\mathfrak{q}_{2^h}(u,t) = (u,2^h)$ and $|t-2^h|<\frac{1}{4}$ so clearly $(u,t)\in \langle gst(z,2^h)\rangle$. For points $(u,t)$ with $t>2^h$ we have that $\mathfrak{q}_{2^h}(u,t) = (\phi_{h}(u),2^h)\in \langle st(z)\rangle$. Since the map $\phi_{h}$ is simplicial, it means that there exists $w\in \phi_{h}^{-1}(z)$ such that $u\in \langle st(w) \rangle\subset |\mathcal{U}_{h}|$. By the definition of geometric stars, this means that $(u,t)\in \langle gst(z,2^h) \rangle $.
\end{proof}

\subsection{Relative (geometric) simplicial approximation}\label{relgeosimplsect}

This subsection contains the details to the proof of Theorem \ref{relgeosimpl}, adapting the proof from Zeeman's classical paper  \cite{zeeman1964relative} on (ordinary) relative simplicial approximation. We also point out a small error in the paper. To avoid confusion, $|K|$ of a simplicial complex $K$ refers to the underlying space, and $K$ refers to $|K|$, along with the data of an abstract simplicial complex. We equip $|K|$ with the metric induced from an embedding into Euclidean space. Throughout this proof, the notation $K-L$ stands for the union of all closed simplices that are not completely contained in the subcomplex $L$. Note that $K-L$ is not a subcomplex of $K$, since a simplex in $K-L$ may have a face in $L$. We use $cl(\cdot)$ to denote the closure of a set. 

\begin{theorem} (\cite{zeeman1964relative}) \label{relordinarysa} Let $K,M$ be finite simplicial complexes, and $L$ a subcomplex of $K$. Let $f: |K|\rightarrow |M|$ be a continous map such that the restriction $f_{|L}$ is a simplicial map from $L$ to $M$. Then there exists an integer $s$ and a simplicial map $g: K_{s}\rightarrow M$ such that $g_{|L}=f_{|L}$ and $g$ is homotopic to $f$ keeping $L$ fixed. \\

Additionally, given an arbitrary neighbourhood $U$ of $|L|$ in $|K|$, we can choose $g$ and the homotopy such that
\begin{enumerate}
\item If $x\notin U$, the homotopy of $x$ is the straight interval in $\Delta(f(x))$ from $f(x)$ to $g(x)$, where $\Delta(f(x))$ is the unique simplex in $M$ which contains $fx$. 
\item If $x\in |L|$, the homotopy leaves $x$ fixed at $f(x)$.
\item If $x\in U - |L|$, the homotopy of $x$ is contained in $N(N(\Delta(f(x)), M), M)$
\end{enumerate}
\end{theorem}

The analagous statement for geometric simplicial approximation is the following. 

\begin{theorem*} (Theorem \ref{relgeosimpl}) Let $K$ be a finite simplicial complex, $L$ a subcomplex of $K$. Let $M$ be the inverse mapping telescope of a compact metric space, and $f: |K|\rightarrow M$ a continuous map such that the restriction $f|_{L}$ is a geometrically simplicial map from $L$ to $M$. Then there exists an integer $s$ and a geometrically simplicial map $g: K_{s}\rightarrow M$ such that $g_{|L} = f_{|L}$ and $g$ is homotopic to $f$ keeping $L$ fixed.\\

Additionally, given an arbitrary neighbourhood $U$ of $|L|$ in $|K|$, we can choose $g$ and the homotopy such that
\begin{enumerate}
\item If $x\notin U$, the homotopy of $x$ is a path from $f(x)$ to $g(x)$ with image contained within the neighbourhood $ \mathcal{N}$ of a geometric simplex.
\item If $x\in |L|$, the homotopy leaves $x$ fixed at $f(x)$.
\item If $x\in U - |L|$, the homotopy of $x$ is contained in $gst^2(\mathscr{C},M)\cup \mathcal{N}$, where $\mathscr{C}$ is the core of a geometric simplex, $gst^2(\mathscr{C},M) := \cup_{v\in vert \mathscr{C}} \langle gst(v) \rangle$, and $ \mathcal{N}$ is the neighbourhood of a geometric simplex.
\end{enumerate}
\end{theorem*}

Before we prove the theorem, we need to define all the terms contained within it. 
\begin{definition} (Faces, stars, neighbourhoods) \begin{enumerate} 
\item Let $K$ be a simplicial complex and $A$ a simplex in $K$. A subset $\{w_{0}, \dots, w_{m}\}\subset vert(A)$ defines a face $[w_0,\dots, w_m]$ of $A$. 
\item Let $A$ be a simplex of $K$. The open star of $A$ in $K$, denoted $\langle st(A,K) \rangle$, is the open subset of $K$ consisting of the union of the interiors of all simplices having $A$ as a face. The closure of $\langle st(A,K) \rangle$ is denoted by $|st(A,K)|$. 
\item The double star of a simplex $A$ in $K$ is
\begin{align*}
\langle st^2(A,K) \rangle = \bigcup_{v \in vert A} \langle st(v,K) \rangle 
\end{align*}
and we denote by  $|st^2(A,K)|$ its closure. 
\item Let $L$ be a subcomplex. We denote by $N(L,K)$ the closed simplicial neighbourhood of $L$ in $K$, which is the subcomplex of $K$ consisting of all interiors of simplices meeting $L$, together with their faces.
\end{enumerate}
\end{definition}

Note that the definition of $\langle st (A,K)\rangle$ agrees the definition in the rest of the thesis if $A$ is a vertex. If $A$ is a simplex we have:
\begin{align*}
N(A,K) = \cup_{v\in vert A} |st(v,K)| = |st^2(A,K)|
\end{align*}

The definition of the barycentric derived complex $K_{r}$ (see Definition \ref{derivedcomplex}) from \cite{zeeman1964relative} is quite nebulously formulated. We explain now what this means more precisely. Let $K_0$ be the original simplicial complex $K$. Vertices in $K_{1}$: for every $[v_0,\dots, v_m] = \sigma \in K-L$ we add a new vertex $b_{\sigma}$ which is the barycentre of $\sigma$. To describe the simplices in $K_{1}$, recall that the simplices in the barycentric subdivision of $\sigma$ are identified with an order $v_{i_{0}}< \dots <v_{i_{m}}$ of the vertices in $\sigma$. Geometrically, this corresponds to the simplex $\tau=[v_{i_{0}}, b_{[v_{i_{0}},v_{i_{1}}]},\dots,b_{\sigma}]$, ie. the vertices are barycentres of an increasing flag of subcomplexes of $\sigma$. A simplex in $K_{1}$ corresponds to a partial order of the vertices in $\sigma$. Suppose that  $[v_{i_0}, \dots, v_{i_{k}}]$ is a simplex in $L$ for some $0\leq k\leq m-1$ and $[v_{i_0}, \dots, v_{i_{k}}, v_{i_{k+1}}] \notin L$ (this can occur even if $v_{i_{k+1}} \in L$). Then the ordering of the first $k+1$ vertices does not matter, but the remaining ones distinguish different simplices in $K_{1}$. Geometrically, the partial order $v_{i_0},\dots,v_{i_{k}} < v_{i_{k+1}}<\dots< v_{i_{m}}$ corresponds to the simplex 
\begin{align*}
[v_{i_0}, \dots, v_{i_{k}}, b_{[v_{i_0}, \dots, v_{i_{k}}, v_{i_{k+1}}]}, \dots, b_{\sigma}]
\end{align*} 
If $v_{i_0}\notin L$ then all the vertices are totally ordered. We call a partial order of the simplices in $\sigma$ described above as admissible. \\

Every simplex in $K_{1}-L$ is uniquely identified with such an admissible partial order, and simplcies in $L$ remain unchanged. Note that it is possible there is some $j>k$ such that $v_{i_{j}}\in L$. The simplices that meet $L$ are exactly those such that $v_{i_0}\in L$, and if  $v_{i_1}\notin L$ we obtain a simplex which coincides with one from the usual barycentric subdivision. We call $K_1$ the first barycentric subdivision of $K$ rel $L$. \\

Let $K$ be a triangulation of $[0,1]^2$ by two triangles. Figure \ref{fig:relsa1} shows an example of the first barycentric subdivision of $K$ rel $L = [0,1] \cup [1,2]$. Observe that the simplices which do not meet $L$ in more than a vertex get barycentrically subdivided, and that $L$ is a subcomplex of $K_{1}$. The definition of $K_{r}$ for $r\geq 2$ proceeds inductively, ie. $K_{r}$ is a subdivision of $K_{r-1}$ relative to $L$. \\

\begin{figure} 
\centering
  \centering
  \includegraphics[width=.5\linewidth]{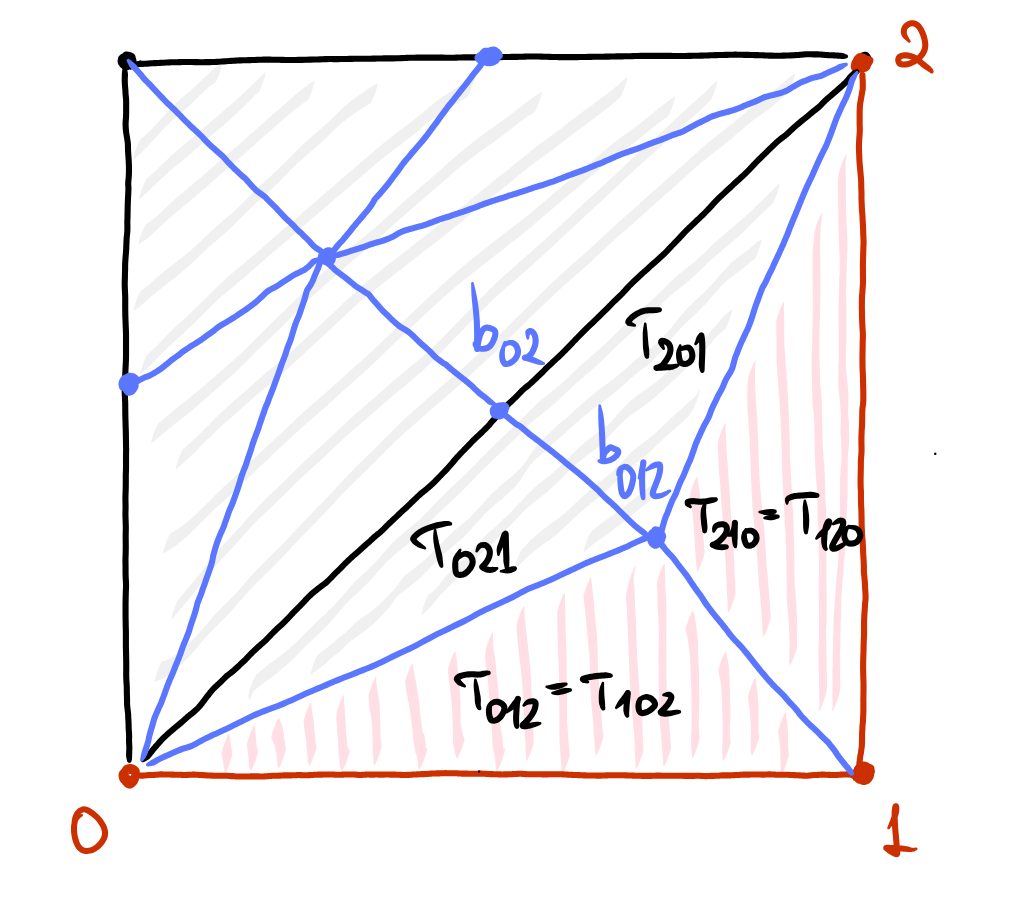}
  \caption{The first barycentric subdivision of $K = [0,1]^2 = \Delta^2 \cup \Delta^2$ relative to $L = [0,1]\cup [1,2]$. Black: vertices and edges in original simplicial complex $K-L$; red: the subcomplex $L$; blue: new vertices and edges in $K_1$; grey: simplices in $K_1$ which coincide with the barycentric subdivision; pink: new types of simplices.}
  \label{fig:relsa1}
\end{figure}

$L$ is full in the barycentric derived complex $K_r$: 
\begin{lemma} \label{fullness} \cite{zeeman1964relative} If $r\geq 1$ then $L$ is full in $K_r$:
\begin{enumerate}
\item No simplex of $K_r - L$ has all its vertices in $L$.
\item Every simplex of $K_r - L$ meets $L$ in a face or the empty set.
\end{enumerate}
\end{lemma}

Let $K$ be a finite simplicial complex.  Let $\mathcal{V}$ denote the open geometric star covering of $M$. The main point is that for an arbitrary continuous map $f: K\rightarrow M$ it may be the case that for no $r$ does the star covering of $K_r$ refine $f^{-1}(\mathcal{V})$. So we first perform a homotopy from $\id_{|K|}$ to a map $\mathfrak{h}:|K|\rightarrow |K| $ which retracts a neighbourhood of $L$ onto $L$, then apply (geometric) simplicial approximation to $f\circ \mathfrak{h}$.\\

Let $B = [v_0,\dots, v_m]$ be a simplex of $L$. Since $f_{|L}$ is simplicial, $f(B)\subset \mathscr{C}_{f(B)}$, where $\mathscr{C}_{f(B)}:= \mathscr{C}_{[f(v_0),\dots, f(v_{m}) ]}$ denotes the core of the geometric simplex $[f(v_0),\dots, f(v_{m}) ]$. Clearly $f(B)\subset gst^2(\mathscr{C}_{f(B)},M) =: \cup_{v\in vert \mathscr{C}_{f(B)}} \langle gst(v,M) \rangle$ so $B$ is contained in the open set $f^{-1}(gst^2(\mathscr{C}_{f(B)},M))$.\\
 
The next definition is a modification of an incorrect definition in the original paper.  $L$ is a full subcomplex of $K_{1}$, so any simplex in $K_1-L$ intersects $L$ in exactly a simplex $A$ in $L$, or the empty set. That means that any simplex $\sigma$ in the subdivision $K_{r} - L$ for $r\geq 1$ lies in a unique simplex $\tau$ in $K_{1}-L$ which intersects $L$ in exactly a simplex $A$ in $L$, or the empty set. 
\begin{definition} (Modified star)
Let $A$ be a simplex of $L$. Denote by $\langle st'(A,K_{r},L) \rangle$ the open subset of $|K_{r}|$ consisting of the union of the interiors of all simplices in $K_{r}-L$ which meet $A$, and are contained within a simplex in $K_{1}$ which intersects $L$ in exactly $A$. Denote by $|st'(A,K_{r},L)|$ the closure of $\langle st'(A,K_{r},L) \rangle$. 
\end{definition} 

The corrected statement from the paper is: let $B$ be a simplex in $L$. As $q\rightarrow \infty$, the compact sets $|st'(B,K_{q},L)|$ converge uniformly in distance to the compact set $B$ (the proof in Lemma $16.3$ of \cite{munkres2018elements} for a special case easily generalises), so for sufficiently large $q$, are contained in the open set $f^{-1}(gst^2(\mathscr{C}_{f(B)},M))$. Since $K$ is finite, we choose $q$ large enough so that $f(|st'(B,K_{q},L)|) \subset f^{-1}(gst^2(\mathscr{C}_{f(B)},M))$ for all simplices $B\in L$. \\

To explain why we needed this modification, the original paper claims the following: \textit{"As $q\rightarrow \infty$, the compact sets $|st(B,K_{q})|$ converge uniformly to the compact set $B$"}, which is false for two reasons. Firstly, the simplices in $L$ which have $B$ as a face obviously do not get smaller, so do not converge to $B$. The larger problem is that if $B$ is not of maximal dimension, then for each subdivision $K_{q}$, there exist simplices which meet $L$ in exactly $B$ which come from a simplex in $K_{q-1}$ which meets $L$ in a larger simplex containing $B$. These simplices cannot be made uniformly close to $B$ for $q\rightarrow \infty$. To illustrate the problem, consider Figure \ref{fig:relsa2}, where $K=[0,1,2]$ is a standard $2$-simplex, $L=[0,1]$ is bottom edge, and $\{1\}=B$ is the bottom right vertex. The problematic simplex in $|st(\{1\},K_{3})|$ is coloured red. Observe that with our modified definition, this simplex belongs to $|st'(L,K_{3},L)|$ but not to $|st'(\{1\},K_{3},L)|$. \\

Short digression: we do not in fact require that $K$ is finite. The statement also holds for uniformly bounded simplicial complexes.

\begin{cor} \label{starconv} Let $K$ be a uniformly bounded simplicial complex. Then for $\varepsilon>0$ there exists a $q$ such that $|st(v,K_{q})| \subset N_{\varepsilon} |st(v,L)|$ for all vertices $v$ in $L$. 
\end{cor} 

\begin{proof} 
Let $\Sigma$ be the set of simplices in $K_1-L$ which meet $L$. For a simplex $\sigma_{B}\in \Sigma$ in $K_1$ which meets $L$ in exactly $B$, there exists a $q$ large enough such that $\sigma_B \cap |st'(B,K_q,L)|$ lies in the $\varepsilon$-neighbourhood of $B$. Since $K_1$ is uniformly bounded, there are only finitely many strong similarity types (with distinguished face) in $\Sigma$. Take the largest-diameter representative in each strong similarity type, then maximise $q$ over the finitely many types. This shows that we can choose a $q$ such that for every simplex $B$ in $L$, $|st'(B,K_q,L)|$ lies in a $\varepsilon$ neighbourhood of $B$. Let $\sigma$ be a simplex in $|st(v,K_{q})|$. Obviously if $\sigma\subset L$ then $\sigma \subset |st(v,L)|$. If $\sigma \subset |st(v,K_q)|\cap K_{q}-L$ there exists a simplex $B\subset L$ such that $\sigma\subset |st'(B,K_{q},L)| \subset N_{\varepsilon} B \subset N_{\varepsilon} |st(v,L)|$.
\end{proof}

Let $K_q'$ denote the subdivision of $K_q$ obtained by barycentrically subdividing all simplices of $K_q - L$ that meet $L$. Note that subdivision-wise, $K_{q}'$ lies between $K_{q}$ and $K_{q+1}$. Let 
\begin{align*}
N:= N(L,K_{q}') = N(L,K_{q+1})
\end{align*}
For a picture of this, see Figure  \ref{fig:relsa3}. \\

\begin{figure}
\centering
  \centering
  \includegraphics[width=.5\linewidth]{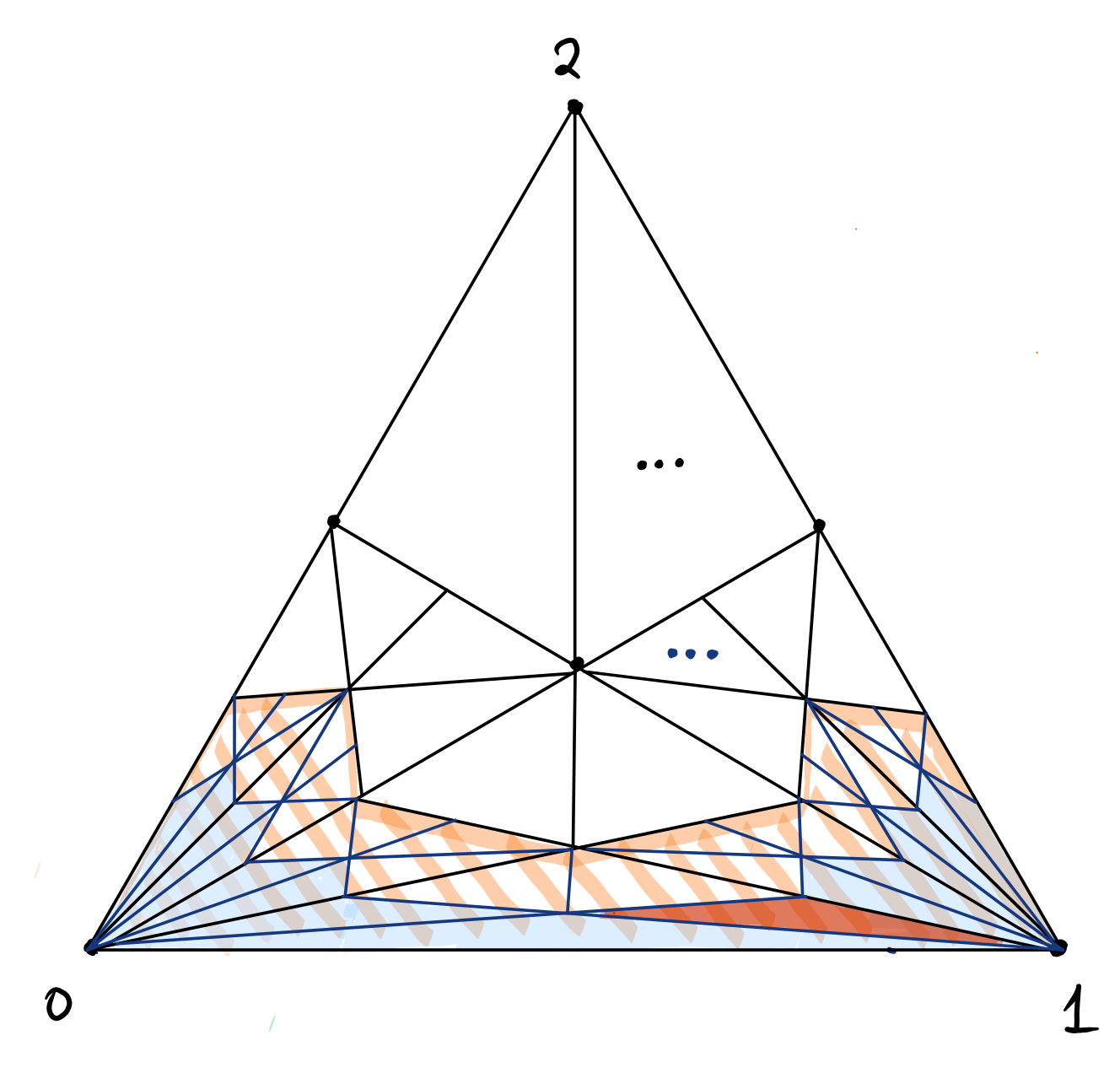}
  \caption{An example of the neighbourhood of the bottom edge $L=[0,1]$ in the second (orange) and third (blue) barycentric subdivision of $K=[0,1,2]$ rel $L$. The "bad" simplex is denoted in red. Dots indicate undrawn simplices.}
  \label{fig:relsa2}
\end{figure}

\begin{figure}
\centering
  \centering
  \includegraphics[width=.8\linewidth]{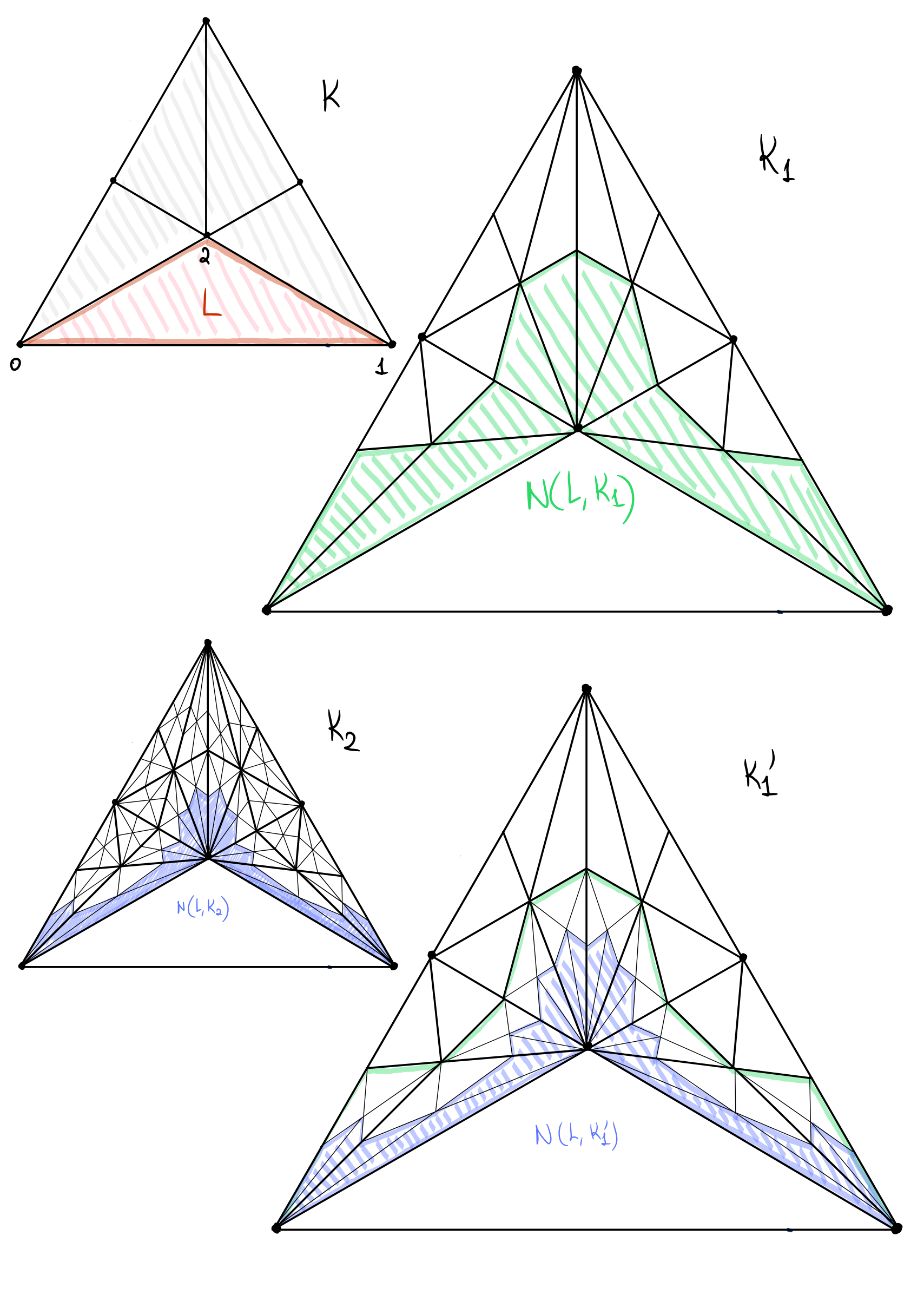}
  \caption{An illustration of $K,K_{1},K_{1}',K_{2}$.}
  \label{fig:relsa3}
\end{figure}

\begin{figure}
\centering
  \centering
  \includegraphics[width=\linewidth]{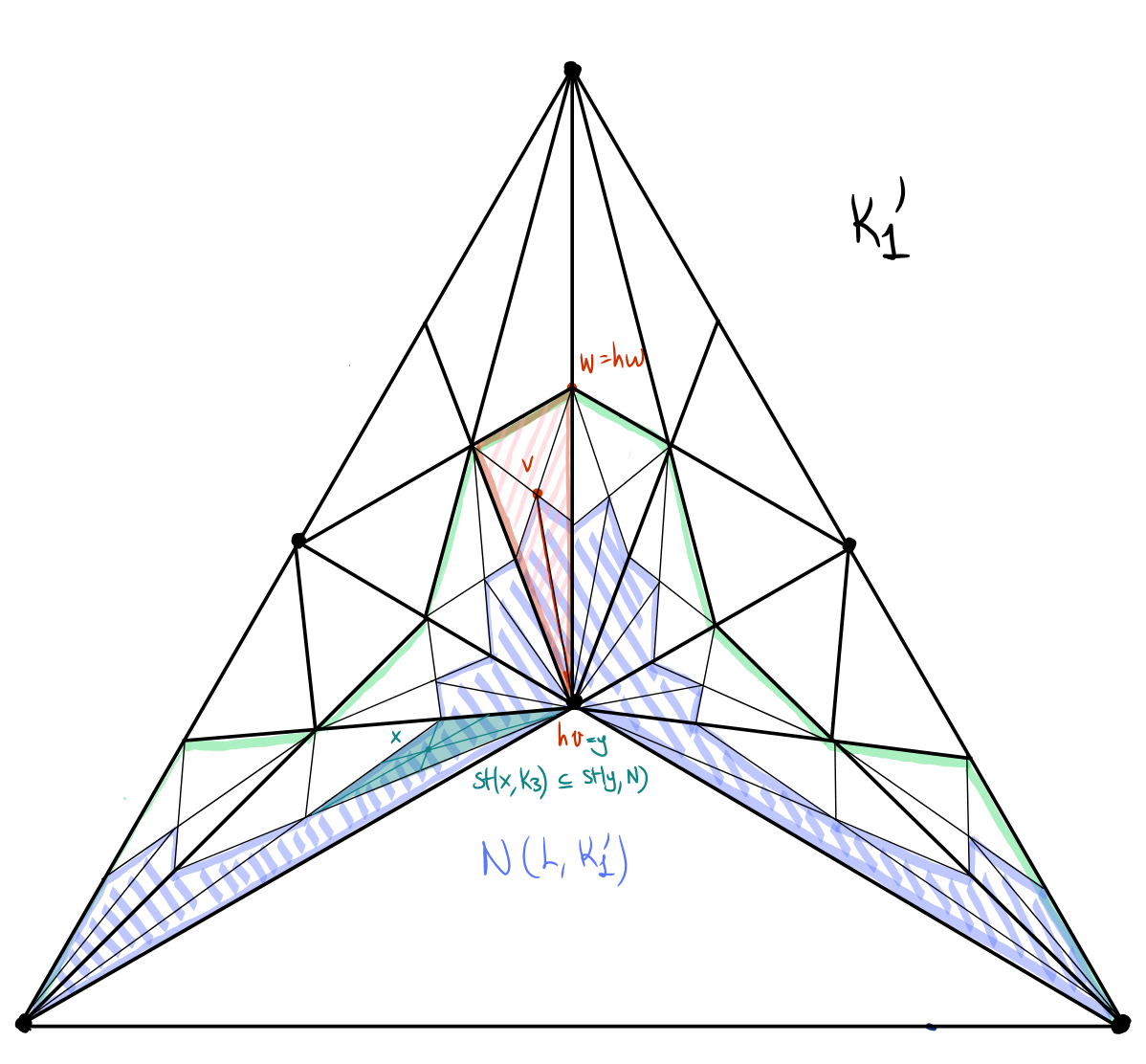}
  \caption{An illustration of the simplicial neighbourhood $N(L,K_{1}') = N(L,K_{2})=: N$ (blue), the neighbourhood $N(L,K_1)$ (outlined in green), and the map $\mathfrak{h}$. $v\in N-L$ and $w\in K'_{1}-N$. $v = b_{\sigma}$ ($\sigma$ shaded red). $x\in N_{3}-V_{3}$ and $x = b_{|st(x,K_3)|}$ ($|st(x,K_3)|$ shaded in turquoise). We have that $st(x,K_{3})\subset st(y,N)$, for both open and closed stars.  }
  \label{fig:relsa4}
\end{figure}

We construct a simplicial map $\mathfrak{h}: K'_{q} \rightarrow K_{q}$ as follows. If $w$ is a vertex in $L \cup (K_{q}' - N)$, define $\mathfrak{h}w= w$. If $v$ is a vertex in $N - L$ then $v$ is the barycentre of some simplex $\sigma$ of $ K_{q}$ meeting $L$: we define $\mathfrak{h}v$ to be a vertex of the face $\sigma \cap L$. See Figure \ref{fig:relsa4}. Therefore $\langle st(v,K'_{q})\rangle \subset \langle st(\mathfrak{h}v,K_q)\rangle$ for all vertices $v\in K'_{q}$, so by the ordinary simplicial approximation theorem, $\mathfrak{h}$ extends linearly to a simplicial map, which is by construction a simplicial approximation to the identity by, say, a homotopy $\mathfrak{h}_t: |K|\rightarrow |K|$ where $\mathfrak{h}_0= \id$ and $\mathfrak{h}_1 = \mathfrak{h}$.  The homotopy leaves the subcomplex $L \cup cl(K_q - N(L,K_q))$ of $K'_q$ fixed. Any other point $x\in |K|$ lies in the interior of some simplex $\sigma \in K_q - L$ which meets $L$, and the homotopy $\mathfrak{h}_t (x)$ of $x$ is the straight interval in $\sigma$ from $x$ to $\mathfrak{h}(x)$. By Lemma \ref{fullness} we have that $N = \mathfrak{h}^{-1}(L)$ and that
\begin{align*}
\mathfrak{h}\langle st(y,N)\rangle = \langle st(y,L) \rangle 
\end{align*}

for all vertices $y\in L$.\\

Since $\id \simeq \mathfrak{h}$ it only remains now to approximate $f\circ \mathfrak{h}$. Let $V$ be the subcomplex of $K_{q+1}$ complementary to $N$, $V:=cl(K_{q+1}-N)$. If $s>q+1$ then the subdivision $K_s$ of $K_{q+1}$ induces subdivisions $N_s, V_s$ of $N$ and $V$. In particular, $V_s$ is the $(s-q-1)$th barycentric derived complex of $V$, because $V$ contains no simplices of $L$.  Let $\alpha_s$ be the star covering of $K_s$ and let $\alpha'_s$ be the subset
\begin{align*}
\alpha'_s:= \{|st(v,K_s)| \,|\, v\in V_s\}
\end{align*}
As $s\rightarrow \infty$, mesh $\alpha'_s\rightarrow 0$. Choose $s$ such that the mesh of $f\alpha'_s$ is less than the Lebesgue number of $\mathcal{V}$, the open cover by geometric stars of $M$ . We can do this because $|K|$ is compact so $f$ is uniformly continuous. We define a vertex map $g: K_{s}\rightarrow M$ as follows: if $v$ is a vertex in $V_{s}$, we let $g(v)$ be a vertex such that $f\mathfrak{h}|st(v,K_{s})| \subset \langle gst(g(v))\rangle$. If $y\in L$ then we define $g(y) = f\mathfrak{h}(y)= f(y)$. \\

If $x$ is a vertex in $N_s - V_s$, we know that  $ st(x,K_s)  =  st(x,N_s)  \subset  st(y,N) $ (for both open and closed stars) for some vertex $y\in L$. Let $g(v): = fy$. To see that this is reasonable, suppose that $\sigma$ is some simplex with $x$ as a vertex. $\sigma$ is in $N_{s}-L$ and lies in a unique simplex $\tau$ in $N$. Recall that $\mathfrak{h} \tau$ becomes a simplex in $L$ (with vertex $y$), denoted $ \tau' $. Therefore
\begin{align*}
f\mathfrak{h}  (\sigma)  \subset f\mathfrak{h}(\tau) = f (\tau') \subset \mathscr{C}_{f(\tau')}  
\end{align*}
since $f_{|L}$ is simplicial. $\mathscr{C}_{f(\tau')}$ denotes the geometric core of the geometric simplex spanned by the vertices $\cup_{v \in \tau'} f(v)$. \\

We now have to check that for all vertices $v$ in $\sigma$, $g(v) \in \mathscr{C}_{f(\tau')}$. Suppose that $v$ is in $\sigma \cap (N_s - V_s)$. $v$ lies in the interior of a unique simplex $\delta$ in $N$ which is a face of $\tau$. Then $st(v,K_{s}) \subset st(z,N)$ for some vertex $z\in \delta \cap L \subset  \tau \cap L $. So we can always choose $g(v) = z \in \mathscr{C}_{f(\tau')}$.\\

Suppose that $v$ is in $\sigma \cap V_{s}$. Then $v$ lies in the interior of a unique simplex $\delta$ which is a face of $\tau$ in $N$. So we have $f\mathfrak{h}(v) \subset f\mathfrak{h}(\delta) \subset f(\tau') \subset \mathscr{C}_{f(\tau')}$. Therefore $f\mathfrak{h}(v)$ is always  $\leq \frac{1}{2}$ in height and $\leq \frac{1}{2}$ in the $|\mathcal{U}_{h}|$-variable away from a vertex in $\mathscr{C}_{f\mathfrak{h}(\delta)}\subset \mathscr{C}_{f(\tau')}$. This shows we can always choose $g(v) = z \in \mathscr{C}_{f(\tau')}$ such that $f\mathfrak{h}|st(v,K_{s})|\subset \langle gst(z) \rangle$. Putting the two cases together yields us $g(vert \sigma)\subset vert \mathscr{C}_{f(\tau')}$.\\

Therefore, the vertex map $g: K_s\rightarrow M$ extends to a Lipschitz map on the interior of simplices. $g$ is not \textit{exactly} the geometric simplicial approximation of $f\mathfrak{h}$ (since $f\mathfrak{h}$ does not satisfy the geometric star condition), but it is homotopic to $f\mathfrak{h}$ by a locally Lipschitz homotopy with respect to the set of simplices in $K_{s}$. The argument for this is a combination of points $1$ and $3$ of Proposition \ref{geosa}: point $1$ applies to simplices $\sigma$ in $K_{s} \times [0,1]$ with the property that $pr (\sigma)$  lies in $V_{s}$, and point $3$ applies if $pr(\sigma)$ lies in $N_{s}$. Since $N_{s}$ and $V_{s}$ are subcomplexes, these two cases do not overlap when increasing dimension of $\sigma$. \\

Both homotopies $\id\simeq \mathfrak{h}, f\mathfrak{h}\simeq g$ keep $L$ fixed, so we have $f\simeq f\mathfrak{h} \simeq g$ keeping $L$ fixed. \\

To prove the second half of the theorem, let $U$ be an arbitrary neighbourhood of $|L|$ in $|K|$. Let 
\begin{align*}
U_q = \cup_{y\in L} \langle st(y,K_q)\rangle 
\end{align*}
As $q\rightarrow \infty$ the sets $cl(U_q) = N(L,K_{q})$ converge uniformly to $|L|$, and so for sufficiently large $q$ are contained in $U$. In the above proof of the theorem we choose $q$ so that in addition $U_q\subset U$. We now verify the three cases:
\begin{enumerate}
\item If $x\notin U$, then $x\notin U_q$, and so $x$ is kept fixed under the homotopy $f\simeq f\mathfrak{h}$, and under the homotopy $f\mathfrak{h} \simeq g$ it moves along a path inside the neighbourhood $\mathcal{N}_{g(\sigma)}$ of the geometric simplex spanned by $\cup_{v \in vert \sigma} g(v)$, for $x \in int (\sigma)$ for $\sigma \in K_{s}$. 
\item If $x\in |L|$, $x$ is kept fixed.
\item If $x\in U - |L|$ then either $x\notin U_q$ or $x\in U_q - |L|$. In the first case condition $1$ is satisfied and therefore also $3$. In the second case $x$ lies in the interior of some simplex $\sigma \in K_q - L$ meeting $L$ in the non-empty face $A$. $\sigma$ comes from a simplex $\tau$ in $K_{1}-L$ meeting $L$ in exactly $B$, which contains $A$. The image of $x$ under the homotopy $\mathfrak{h}_t$ is contained in $\sigma$ and therefore
\begin{align*}
f\mathfrak{h}_t (x)\in f(\sigma) \subset f(|st'(B,K_{q},L)|) \subset gst^2(\mathscr{C}_{f(B)},M)
\end{align*}
The homotopy $f\mathfrak{h}\simeq g$ is contained within the neighbourhood $ \mathcal{N}$ of a geometric simplex.
\end{enumerate}

\subsection{Slightly different definitions of homotopy groups}\label{diffdef}

This subsection will deal with the fact that we use slightly different definitions of coarse homotopy groups and Euclidean cones than Mitchener, Norouzizadeh and Schick in \cite{mitchener2020coarse}. Our cones $cX$ are actually truncated at height $1$ where we delete the tip. We will denote with a prime $'$ all the notation originally introduced in \cite{mitchener2020coarse} (for example, $c'X$ denotes the entire cone). The goal is to prove that our different definitions are equivalent.

\begin{definition} Let $X$ be a bornological coarse space and let $p: X\rightarrow \mathbb[0,\infty)$ be a controlled and bornologous map. We define the cylinder
\begin{align*}
I'_p(X) := \{(x,t)\in X \times [0,\infty)\,|\, t\leq p(x) +1\}
\end{align*}
with inclusions $i_0: X\rightarrow I'_p(X)$ and $i_1: X\rightarrow I'_p(X)$ defined by $i_0(x) = (x,0)$ and $i(x) = (x,p(x)+1)$. 
\end{definition}

Recall that \cite{mitchener2020coarse} actually requires $p$ to be controlled and proper. We have already discussed that properness is not a requirement for any of their results, and that their definition of bornology automatically makes all controlled maps bornologous. 

\begin{definition} Let $X$ and $Y$ be bornological coarse spaces. 
\begin{itemize}
\item A coarse homotopy is a coarse map $H: I'_p(X)\rightarrow Y$ for some controlled and bornologous map $p: X\rightarrow [0,\infty)$. 
\item We call coarse maps $f_0:X\rightarrow Y$ and $f_1: X\rightarrow Y$ coarsely homotopic if there is a coarse homotopy $H: I'_p(X)\rightarrow Y$ such that $f_0 = H\circ i_0$ and $f_1 = H\circ i_1$. 
\end{itemize}
\end{definition}

\begin{definition} Let $X$ be a bornological coarse space with base ray $\omega:[0,\infty)\rightarrow X$. For $n\geq 1$ define the $n$-th coarse homotopy group $\pi_n^{'c}(X,\omega)$ to be the set of relative coarse homotopy classes of maps 
\begin{align*}
f: (c'([0,1]^n), c'(\partial [0,1]^n))\rightarrow (X,\omega[0,\infty))
\end{align*}
such that $f|_{c'(\partial[0,1]^n)} = \omega \circ p$, where $p:c([0,1]^n)\rightarrow [0,\infty); (hx,h)\mapsto h$ denotes the height variable of the cone. 
\end{definition}

\begin{definition} 
For a bornological coarse pair $k_A: A\rightarrow X$ with base ray $\omega:[0,\infty)\rightarrow A$ we define the relative $n$-th coarse homotopy "group" $\pi_{n}^{c'}(X,A, \omega)$ to be the set of relative coarse homotopy classes of maps 
\begin{align*}
f: (c'([0,1]^n), c'(\partial [0,1]^n), c'(\partial_{+}[0,1]^n))\rightarrow (X,A,\omega[0,\infty))
\end{align*}
such that $f|_{c'(\partial_{+}[0,1]^n)}= \omega\circ p$ where $\partial_{+}[0,1]^n := \{(x_1,\dots,x_n)\in \partial [0,1]^n\,|\,x_n>0\}$. 
\end{definition}

\begin{lemma} Let $X$ be a bornological coarse space with $\omega: [0,\infty)\rightarrow X$ a base ray. There is an isomorphism $\pi_{n}^{c}(X,\omega_{|[1,\infty)}) = \pi_{n}^{'c}(X,\omega)$. 
\end{lemma}

Note that any $[1,\infty)$-base ray can be obviously extended to a $[0,\infty)$-base ray. If $X$ is a compact metric space and we let $c'X$ be our bornological coarse space, we obtain $\pi_{n}^{c}(c'X,\omega_{|[1,\infty)}) = \pi_{n}^{'c}(c'X,\omega)$. Since $cX$ is coarsely equivalent to $c'X$, we get $\pi_{n}^{c}(cX,\omega_{|[1,\infty)}) = \pi_{n}^{'c}(c'X,\omega)$.

\begin{proof} Define $\varphi: \pi_{n}^{c}(X,\omega_{|[1,\infty)})\rightarrow \pi_{n}^{'c}(X,\omega)$ as follows. $[f]\in \pi_{n}^{c}(X,\omega_{|[1,\infty)})$ is represented by a map $f: (c[0,1]^n,c\partial [0,1]^n)\rightarrow (X,\omega_{|[1,\infty)})$. We extend $f$ to a map $f':(c'[0,1]^{n},c'\partial [0,1]^n)\rightarrow (X,\omega)$ by letting $f'(hx,h)=\omega(h)$ for $0\leq h<1$. Since $f'$ sends everything below height $1$ to a bounded set, $f'$ is still coarse, and it is obviously base ray preserving. Let $\varphi[f]=[f']$.\\

Independence of homotopy class of $f$: suppose we have a homotopy 
\begin{align*}
H: (c[0,1]^{n+1},c(\partial [0,1]^n \times [0,1]))\rightarrow (X,\omega_{|[1,\infty)})\\
H_{|c([0,1]^n \times \{0\})} = f \quad H_{|c([0,1]^n \times \{1\})} = g
\end{align*} 
This extends to a coarse homotopy $H'$ on $c'[0,1]^{n+1}$ letting ${H'}(hx,h,ht) = \omega(h)$ for $0\leq h<1$ and $t\in [0,1]$ which restricts to $f'$ and $g'$ on $c'([0,1]^n \times \{0\})$ and $c'([0,1]^n \times \{1\})$ respectively. The required homotopy on $I'_{p}(c[0,1]^n)$ is then just defined by extending constantly by $g'$ for $t\in [h,h+1]$. This independence of homotopy class and the fact that $\varphi$ is obviously a group homomorphism, gives us well-definedness of $\varphi$.\\

Define $\psi:\pi_{n}^{'c}(X,\omega) \rightarrow \pi_{n}^{c}(X,\omega_{|[1,\infty)})$ as follows. $[f']\in \pi_{n}^{'c}(X,\omega)$ is represented by a map $(c'[0,1]^n,c'\partial [0,1]^n)\rightarrow (X,\omega)$. We define $f: (c[0,1]^n,c\partial[0,1]^n)\rightarrow(X,\omega_{[1,\infty)})$ to just be the restriction of $f'$. Let $\psi[f']=[f]$. \\

Independence of homotopy class of $f'$: suppose we have a base ray preserving homotopy $H': I'_{p}(c'[0,1]^n) \rightarrow X$ with $H'\circ i_{0} = f', H'\circ i_{1} = g'$. We define a homotopy  $H: c[0,1]^{n+1}\rightarrow X$ by letting $H(hx,h,ht) = H'(hx,h,(h+1)t)$. Observe that $H(hx,h,0)= f$ and $H(hx,h,h)= g$ and that $H$ is obviously base ray preserving. The map which takes $(hx,h,ht)\mapsto (hx,h,(h+1)t)$ is controlled and proper; therefore $H$ is controlled and proper. This gives us well-definedness of $\psi$.\\

$\psi\circ \varphi[f]=[f]$ just by following the definitions. $\varphi\circ \psi[f]$ is represented by a map $g$ which is given by $g(hx,h)= f(hx,h)$ if $h\geq 1$ and by $g(hx,h)= \omega(h)$ otherwise. To show that $f$ and $g$ are close, we only need to consider points for which $g(hx,h) = \omega(h)$. Let $B$ be the bounded set $c'[0,1]^n\setminus c[0,1]^n$. Since $f$ is controlled, we know that $(f\times f )( B \times \{(0,0)\}) =: U$ is an entourage. Similarly, $(\omega \times \omega)(\{0\} \times [0,1])=: V$ is an entourage. We have that for all $(hx,h) \in B$
\begin{align*}
(f \times g)((hx,h),(hx,h)) = (f(hx,h),\omega(h)) \in U \circ V. 
\end{align*}
Therefore $(f\times g)(\Delta) = \Delta \cup (U\circ V)$. 
\end{proof}

Note that the same proof gives us the isomorphism 
\begin{align*}
\pi_{n}^{'c}(X,A,\omega) = \pi_{n}^{c}(X,A,\omega_{[1,\infty)})
\end{align*}

\begin{lemma} \begin{enumerate}
\item Let $(X,\omega),(Y,\tau)$ be bornological coarse spaces with $[0,\infty)$-base rays and $\mathfrak{a}: (X,\omega)\rightarrow (Y,\tau)$ a base ray preserving coarse map. The diagram 

\begin{figure}[H]
\center
\begin{tikzcd}
& \pi_{n}^{c}(X,\omega_{[1,\infty)})  \arrow{r}{\mathfrak{a}_{*}} \arrow{d}{\varphi}
&\pi_{n}^{c}(Y,\tau_{[1,\infty)}) \arrow{d}{\varphi}
  \\
& \pi_{n}^{'c}(X,\omega) \arrow{r}{\mathfrak{a}_{*}}
&\pi_{n}^{'c}(Y,\tau)
\end{tikzcd}
\end{figure}

commutes.
\item If $(X,A,\omega)$ is a based bornological coarse pair with map $k: A\rightarrow X$, the diagram

\begin{figure}[H]
\center
\begin{tikzcd}
& \pi_{n}^{c}(A,\omega_{[1,\infty)})  \arrow{r}{k_{*}} \arrow{d}{\varphi}
&\pi_{n}^{c}(X,\omega_{[1,\infty)}) \arrow{d}{\varphi} \arrow{r}{r}  \arrow{d}{\varphi}
& \pi_{n}^{c}(X,A,\omega_{[1,\infty)}) \arrow{r}{\partial} \arrow{d}{\varphi}
& \pi_{n-1}^{c}(A,\omega_{[1,\infty)}) \arrow{d}{\varphi}
  \\
& \pi_{n}^{'c}(A,\omega) \arrow{r}{k_{*}}
&\pi_{n}^{'c}(X,\omega) \arrow{r}{r}
& \pi_{n}^{'c}(X,A,\omega) \arrow{r}{\partial} 
& \pi_{n-1}^{'c}(A,\omega)
\end{tikzcd}
\end{figure}
commutes 
\end{enumerate}
\end{lemma}

\begin{proof} The diagrams all commute on the level of representatives. 
\end{proof}

\subsection{Spherical simplices} \label{sphs}

We discuss in this section why simplices in $M$ are equipped with the spherical metric. One of the reasons is to show that the map $R: M\rightarrow cX$ is coarse by restricting distance estimates to points lying within a common (geometric) simplex. We begin first with an explanation of why the standard metric is inadequate for our purposes. 

\begin{example} (Example $5.39$ in \cite{bunke2020homotopy}) Let $\Delta^n$ be the standard $n$-simplex. We define the simplicial complex $Z$ as the quotient of $\coprod_{n\in \mathbb{N}} \Delta^n$ by identifying for all $n\in \mathbb{N}$ the last face $\partial_{n}\Delta^n$ of $\Delta^n$ with the first face of the first face $\partial_{0}\partial_{0} \Delta^{n+1}$ of $\Delta^{n+1}$.

\begin{figure}[H]
\centering
  \centering
  \includegraphics[width=0.6 \linewidth]{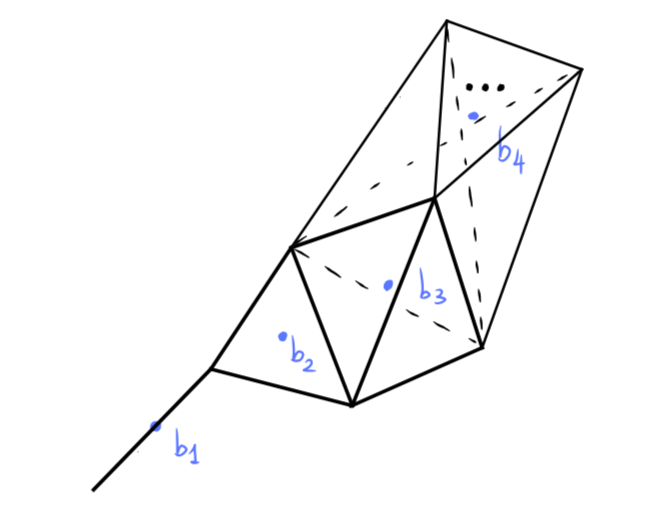}
  \caption{The simplicial complex $Z$. A "bad" path has vertices $b_1,b_2,b_3,b_4,\dots$. }
  \label{fig:Z}
\end{figure}

For the barycentres $b_{n}$ and $b_{n+1}$ of $\Delta^n$ and $\Delta^{n+1}$ we have 
\begin{align*}
d(b_{n},b_{n+1}) = \sqrt{\frac{1}{n(n+1)}}+ \sqrt{\frac{2}{n(n+2)}}
\end{align*}

Define $k: \mathbb{N}\rightarrow \mathbb{N}$ by 
\begin{align*}
k(n):= \max\{m\in \mathbb{N}\,|\, d(b_{n},b_{n+m}\leq 1)\}
\end{align*}
We have that $\lim_{n\rightarrow \infty} k(n) = \infty$. Hence for all $n\in \mathbb{N}$,
\begin{align*}
d(b_{n},b_{n+k(n)})\leq 1
\end{align*}
\end{example}

The example shows that if we take path metric induced from the standard metric on each slice of $M$, there may be points $u,v$ in $M$ such that $d(u,v)\leq 1$ but every path between them passes through an arbitrarily large number of simplices. Therefore, we could not for example have shown that the map $R: M\rightarrow cX$ is coarse by restricting to points lying within a common (geometric) simplex. In the spherical metric, such "bad" paths cannot occur (in the example, we have that in the spherical metric $d^s(b_{n},b_{n+1}) = \frac{\pi}{2}$). In fact, we have the following Lemma:

\begin{lemma} \label{sphericals} (Lemma $A.5$ in \cite{wright2005coarse}) Let $Z$ be a locally finite simplicial complex with a uniform spherical metric. Then the distance from a vertex $v$ of $Z$ to a simplex $\sigma$ in the same component of $Z$ is $\pi/2$ times the length (i.e. number of edges) of the shortest simplicial path from $v$ to $\sigma$. In other words there are no shortcuts through the interior of a simplex. 
\end{lemma}

In our scaled spherical metric, so that the length of edges is $1$, the distance in the Lemma above is just the number of edges. \\

Let $x,y$ be two points in $|\mathcal{U}_{h}|$ with distance at most $Q\in \mathbb{N}$. Let $x$ be contained in the interior of $\tau_{1}$ and $y$ be contained in the interior of $\tau_2$. Let $x'$ be the closest vertex in $\tau_{1}$ to $x$, and $y'$ the closest vertex in $\tau_2$ to $y$.  We have that 
\begin{align*}
d(x',y')\leq d(y',x) + d(x,x') \leq d(x,y) + d(y,y') + d(x,x') \leq Q+1
\end{align*}
$d(x',y')$ is realised by a vertex path of distance at most  $Q+1$. Therefore, it travels through at most $Q+1$ simplices in $|\mathcal{U}_{h}|$. So, there is a path between $x$ and $y$ which can be split up into $Q+3$ components, each lying within a single simplex. \\

Suppose now we have two points $(x,t),(y,s)$ in the mapping telescope, which have distance less than $Q$ from one another. The shortest path between them is a concatenation $\sigma_{x}\ast \eta \ast \sigma_{y}^{-1}$ where $\sigma_{x},\sigma_{y}$ are straight downwards paths, and $\eta$ is a geodesic in $|\mathcal{U}_{h}|$ on a single height slice. $\sigma_{x}$ can be broken up into at most $Q+1$ segments each lying within a geometric cylinder, and similarly for $\sigma_{y}$. The segment $\eta$ can be split up into $Q+3$ components. So the path in total can be split into $3Q+5$ components. Let $C$ be the entourage such that $(u,v)\in C$ if $u$ and $v$ lie in a common geometric cylinder. If we have a coarse map $f: M\rightarrow Y$ such that $(f \times f)(C) \subset W$, where $W$ is an entourage of $Y$, then $(f \times f) (D_{Q} \times D_{Q}) \subset W^{\circ (3Q+5)}$. Or in other words, $C$ generates the coarse structure on $M$. This is why we can restrict to a single geometric cylinder when computing coarse constants. 

\newpage

%

\section{References}

\renewcommand\refname{}
\vspace*{-15pt}
\bibliographystyle{alpha}
\bibliography{bibliography.bib}

\end{document}